\documentclass{article}
\usepackage{cl2e}
\usepackage{amsmath,amssymb,amscd,,amsthm}
\usepackage[english,german]{babel}
\usepackage[utf8]{inputenc}

\usepackage[T1]{fontenc}
\usepackage{imakeidx}
\usepackage{booktabs}
\usepackage{longtable}
\usepackage{rotating}
\usepackage{booktabs}
\usepackage{multirow}
\usepackage{graphics}
\usepackage{tikz}
\usetikzlibrary{arrows.meta,  quotes, positioning,calc, intersections,babel}
\usepackage{hyperref}

\usepackage{listings}

\newcommand{\N}{{\mathbb N}}
\newcommand{\Z}{{\mathbb Z}}
\newcommand{\Q}{{\mathbb Q}}
\newcommand{\R}{{\mathbb R}}
\newcommand{\C}{{\mathbb C}}
\newcommand{\F}{{\mathbb F}}
\newcommand{\cO}{{\mathcal O}}

\newcommand{\Cl}{\operatorname{Cl}\,}
\newcommand{\disc}{\operatorname{disc}\,}
\newcommand{\diff}{\operatorname{diff}\,}
\newcommand{\Tr}{\operatorname{Tr}\,}
\newcommand{\Gal}{\operatorname{Gal}\,}
\newcommand{\ov}{\overline}
\newcommand{\uF}{\underline{F}}
\newcommand{\uK}{\underline{K}}
\newcommand{\cM}{{\mathcal M}}
\newcommand{\ftw}{{\mathfrak 2}}
\newcommand{\fth}{{\mathfrak 3}}
\newcommand{\fsv}{{\mathfrak 7}}
\newcommand{\fel}{{\mathfrak 1}{\mathfrak 1}}
\newcommand{\fp}{{\mathfrak p}}
\newcommand{\fq}{{\mathfrak q}}
\newcommand{\frt}{{\mathfrak 2}}
\newcommand{\eps}{\varepsilon}

\newcommand{\rsp}{\raisebox{0em}[2.7ex][1.3ex]{\rule{0em}{2ex} }}

\renewcommand{\Re}{\operatorname{Re}}
\renewcommand{\Im}{\operatorname{Im}}

\makeindex[name=N,title=Name Index]
\makeindex[name=S,title=Subject Index]
\makeindex[name=dN,title=Namensverzeichnis]
\makeindex[name=dS,title=Sachverzeichnis]
\indexsetup{headers={\indexname}{\indexname}}

\title{Euclidean Rings}
\subtitle{Diploma Thesis 1989}
\author{Franz Lemmermeyer}
\date{}

\begin{document}

\selectlanguage{english}

\maketitle

\newpage

My diploma thesis on Euclidean rings was written in the late 1980s on
an Atari using Signum!. Later I scanned the work, and have now typeset
it in \LaTeX{} with the help of AI, correcting the most serious errors
that arose from the poor quality of the scanner, faulty OCR, and an
occasionally rather bold AI. Apart from minor details I have left
everything in its original state, including inaccurate historical
statements.

I put the German text at the end of the current file; let me know about
any of the errors that have crept into the text or the translation.
For an updated version of the tables of Euclidean minima, see
the web page of
Pierre Lezowski\footnote{\url{https://www.math.u-bordeaux.fr/~plezowsk/index.php}}.

\medskip

This diploma thesis was the basis of my survey article
{\em The Euclidean algorithm in algebraic number fields}.

\vskip 1 cm

Franz Lemmermeyer, Jagstzell 2026

\newpage

\tableofcontents

\chapter*{Overview}
\addcontentsline{toc}{chapter}{Overview}
\markboth{Euclidean Rings}{Overview}

To orient the reader, we begin with a brief overview of the material
that follows. The starting point of our investigations is
Lenstra's\index[N]{Lenstra} criterion~(1.15) \cite{Len74,Len77a},
according to which a number field possessing a sufficiently long
``exceptional sequence'' is norm-Euclidean. Around the same time
(1975--1977) Cooke\index[N]{Cooke} introduced the notion of a
$k$-stage norm-Euclidean number field. This raised the question
whether Lenstra's criterion could be adapted to detect $k$-stage
norm-Euclidean rings. In fact the adaptation is straightforward
(see~(1.16)). The consequent study of continued fractions of length
at most~$k$ with coefficients in the ring under consideration leads
to sets that resemble those introduced by Motzkin\index[N]{Motzkin}
\cite{Mot49}; we present these sets in~\S~0.

The resulting criterion yields, for example, $2$-stage norm-Euclidean
fields of degree~$2$, $3$, $4$, $5$, and higher. One is immediately led
to ask whether any of these fields are already $1$-stage
norm-Euclidean. The most general criterion then available for proving
that a number ring fails to be norm-Euclidean was~(1.6). Van~der
Linden\index[N]{Linden@van der Linden} \cite{Lin83} describes its
principal limitation as follows:

\begin{quote}
  {\em If we generalize the methods for higher degree number fields it
  turns out that the arithmetical methods are applicable for
  extensions of $\Q$ in which at least one prime is totally
  ramified\dots.}
\end{quote}

A careful examination of the proof of~(1.6) shows that the criterion
admits a useful generalization (see~(1.4)): it often suffices that
$K/\Q$ contain an intermediate field~$k$ for which purely ramified
ideals exist in~$K/k$. In~\S~6 we demonstrate the value of this
generalization by determining all norm-Euclidean fields of the form
$\Q(\sqrt{m},\sqrt{n})$ with $m,n\in\Z$ and $m<0$. Moreover,
criterion~(1.4) not only rules out the existence of a Euclidean
algorithm but, in favourable cases, also yields estimates for the
Euclidean minimum of~$K$.

Criterion~(1.8), which was termed ``geometric'' by van~der
Linden,\index[N]{Linden@van der Linden} establishes that a real
quadratic field is not norm-Euclidean. In the form presented here it
is due to Barnes\index[N]{Barnes} and
Swinnerton-Dyer\index[N]{Swinnerton-Dyer} \cite{BS52a}, though the
essential idea already appears in
R\'edei\index[N]{Redei@R\'edei} \cite{Red41a,Red41b} and
Inkeri\index[N]{Inkeri} \cite{Ink49}. Partial generalizations to
cubic fields of unit rank~$1$ were given by
Taylor\index[N]{Taylor} \cite{Tay76} and
Cioffari\index[N]{Cioffari} \cite{Cio79}. Continuing the passage
quoted above, van~der Linden\index[N]{Linden@van der Linden} remarks:

\begin{quote}
  {\em \ldots and the geometrical methods apply for fields with
    $\# S_\infty < 2$}
\end{quote}
(that is, for number fields of unit rank~$1$). In fact~(1.8) extends
readily to arbitrary number fields with infinite unit group
(see~(1.12)). Remarkably, the inequalities employed in that extension
(for instance~(1.9) and its cubic analogue in~\S~4) had already been
used by Cassels\index[N]{Cassels} to improve the bounds on
$|\disc K|$ beyond which quadratic ($r=2$, $s=0$), cubic ($r=1$,
$s=1$), and biquadratic ($r=0$, $s=2$) fields cease to admit a
Euclidean algorithm. This strongly suggests that Cassels' result
itself admits a generalization to fields of arbitrary degree; a
natural first step would be a generalization of~(2.12).

In~\S~2 we describe a method for computing Euclidean minima that
originates essentially with Barnes\index[N]{Barnes} and
Swinnerton-Dyer\index[N]{Swinnerton-Dyer} \cite{BS52a}; important
contributions were also made by Cassels\index[N]{Cassels},
Inkeri\index[N]{Inkeri}, and others. The proofs given in~\S~2 are
substantially clearer and more transparent than the original
arguments. In~\S~3, employing~(1.4) and~(1.8), we give an almost
complete classification of the norm-Euclidean real quadratic fields,
following the classical lines of Behrbohm\index[N]{Behrbohm},
R\'edei\index[N]{Redei@R\'edei}, and Erd\"os\index[N]{Erdos@Erd\"os}
\& Ko\index[N]{Ko}, among others. The results of that section are
used again in~\S~5.

In~\S~4 we turn to cubic fields. To determine all norm-Euclidean
fields of the form $\Q(\sqrt[3]{m})$ we follow Cioffari
\cite{Cio79}\index[N]{Cioffari}. We then compute the Euclidean
minima of cubic fields of small discriminant and improve a bound of
Smith\index[N]{Smith} \cite{Smi69}: there are no cyclic cubic fields
with discriminant~$d$ in the range $157^2\le d\le 10^8$. We extend
the non-existence statement to the larger interval
$157^2\le d\le 2.5\cdot 10^{11}$.

In~\S~5 we prove that only finitely many norm-Euclidean fields of the
form $\Q(\sqrt[4]{m})$ exist, and for $m<0$ we determine them all.
This completes the investigations begun by Cioffari
\cite{Cio79}\index[N]{Cioffari} and Egami
\cite{Ega84}\index[N]{Egami}.

As already indicated, in \S~6 we determine all imaginary bicyclic
biquadratic number fields. Since van~der
Linden\index[N]{Linden@van der Linden} \cite{Lin83} found the two
cyclic imaginary fields of degree~$4$, namely $\Q(\zeta_5)$ and the
unique degree-$4$ subfield of $\Q(\zeta_{13})$, all complex quartic fields
that are Galois and admit a Euclidean algorithm are now known.
The next most challenging class appears to be those fields whose
normal closure has dihedral Galois group.

In~\S~7 we begin the classification of the norm-Euclidean Dirichlet
number fields, determining all such fields for which the relative
discriminant is divisible by~$(1+i)$. A complete solution of the
problem appears to be only a matter of computation.

Section~8 offers suggestions concerning the remaining imaginary
fields of degree~$4$, while~\S~9 provides an outlook on fields of
higher degree. It is not surprising that many questions remain open.
Some of them should be solvable by the methods developed here
(lack of time prevented their full resolution), while others seem to
lie deeper. Section~10 collects a number of these open questions.

Finally, \S~11 describes an algorithm that decides, for a given
number field, whether a Euclidean algorithm exists. Whether the
algorithm always terminates is closely connected with certain still
unproved conjectures of Barnes\index[N]{Barnes} and
Swinnerton-Dyer\index[N]{Swinnerton-Dyer}, as well as of
Lenstra\index[N]{Lenstra}; see Lenstra \cite{Len79b}\index[N]{Lenstra}.

\newpage

We conclude with a list of fields whose norm-Euclidean character was
previously open:
\begin{itemize}
\item $n=3$, $r=3$, $s=0$: $\disc K= 2021, 2024, 2057,
  2101, 2213$;
\item $n=4$, $r=0$, $s=2$: $\Q(\sqrt{2}, \sqrt{5})$,
  $\Q(\sqrt{3}, \sqrt{17})$, $\Q(\sqrt{3},
  \sqrt{19})$, $\Q(\sqrt{7}, \sqrt{5})$,
  $\Q(\sqrt{5+2i})$, $\Q(\sqrt{1+6i})$,
  $\Q(\sqrt{7+2i})$, $\Q(\sqrt{7+4i})$,
  $\Q(\sqrt{1+8i})$, $\Q(\sqrt{3+8i})$,
  $\Q(\sqrt{5+8i})$, $\Q(\sqrt{9+4i})$,
  $\Q(\sqrt{7+8i})$, $\Q(\sqrt{11+4i})$,
  $\Q(\sqrt[4]{12})$;
\item $n=4$, $r=2$, $s=1$: $\Q(\sqrt[4]{2})$, $\Q(\sqrt[4]{5})$;
\item $n=4$, $r=4$, $s=0$: $\Q(\sqrt{3}, \sqrt{5})$,
  $\Q(\sqrt{3}, \sqrt{7})$, $\Q(\sqrt{5}, \sqrt{13})$,
  $\Q(\sqrt{5}, \sqrt{17})$.
\end{itemize}

Many further norm-Euclidean fields of degree~$3$ and~$4$ (and of
higher degree) can certainly be found by the methods of this work.
Besides the fourteen real quadratic $2$-stage norm-Euclidean fields
discovered by Cooke\index[N]{Cooke}, the following examples are also
known: $\Q(\sqrt{m})$ for
$$ m = 47, 59, 62, 67, 71, 109, 149, 157, 161, 173, 193, 201, 213, $$
as well as the cubic fields with discriminant $d=-199$ and $d=-351$.

\newpage

\section*{NOTATION}
\begin{longtable}{ll}
\toprule
$\N$ & the natural numbers \\
$\Z$ & the rational integers \\
$\Q$ & the rational numbers \\
$\R$ & the real numbers \\
$\C$ & the complex numbers \\
$\F_q$ & the finite field with $q$ elements \\
$K$ & an algebraic number field \\
$R$ & an integral domain \\
$\mathcal{O}_K$ & the ring of integers of $K$ \\
$\Phi_K$ & Euler's totient function in $K$ \\
$\|I\|$ & the absolute norm of an ideal $I$ \\
$\Tr_{K/\Q}(\alpha)$ &
   the trace of an element $\alpha \in K$ \\
$N_{K/\Q}(\alpha)$ &
   the norm of an element $\alpha \in K$ \\
$\disc K$ & the discriminant of a field extension \\
$\diff K$ & the different of a field extension \\
$(K:\Q)$ & the degree of a field extension \\
$\Gal(K/\Q)$ & the Galois group \\
$Z(\mathfrak{P}|\mathfrak{p})$ & the decomposition group \\
$T(\mathfrak{P}|\mathfrak{p})$ & the inertia group \\
$V_i(\mathfrak{P}|\mathfrak{p})$ & the ramification groups \\
$\Cl(K)$ & the ideal class group \\
$h(K) = \#\Cl(K)$ & the class number \\
$|\cdot|_i$ & the normalized valuations of $K$ (see pp.\ \pageref{pBew}
   and \pageref{pBew2}) \\
$D(m)$ & the ring of integers of $\Q(\sqrt{m})$ \\
$D(m,n)$ & the ring of integers of $\Q(\sqrt{m},\sqrt{n})$ \\
QRG & quadratic reciprocity law \\
FE & fundamental unit \\
GHB & integral basis \\
PERS & prime residue system of units (p.~\pageref{pPERS}) \\
$E_i$ & see p.\ \pageref{pEi} \\
$E'_i$ & see p.\ \pageref{Ej'} \\
$F_i$ & see p.\ \pageref{pFi} \\
$B_i$ & see p.\ \pageref{DefB} \\
$M_1$ & first Euclidean minimum (see p.\ \pageref{pEuM}) \\
$M_2, \dots$ & second Euclidean minimum, \ldots (see p.\ \pageref{pEuM2}) \\
$M^1$, $M^2$, \ldots
     & minimum for $k$-stage Euclidean functions: see p.\ \pageref{DefMk} \\
$C_1, C_2, \dots$ & see p.\ \pageref{DefC} \\
$M_{r,s}$ & the Minkowski bound (see p.\ \pageref{pEuM}) \\
$\mu_K$ & the Lenstra constants (see p.\ \pageref{Lenmu}) \\
$\lambda_K$ & see p.\ \pageref{Lenmu} \\
$\cM_K$ & the ``general measure'', see p.\ \pageref{GMen} \\
$\uK$ & see p.\ \pageref{puK} \\
\bottomrule
\end{longtable}

\chapter*{\S\ 0 Introduction}
\setcounter{chapter}{0}
\addcontentsline{toc}{chapter}{\S\ 0 Introduction}
\markboth{Euclidean Rings}{\S\ 0 Introduction}

As its name suggests, the history of the Euclidean algorithm is a
long one: it begins with Euclid's ``proof''
that the ring $\Z$ of rational integers is ``Euclidean''
(see the works of Hendy\index[N]{Hendy} \cite{Hen75},
Knorr\index[N]{Knorr} \cite{Kno76}
and Collison\index[N]{Collison} \cite{Col80}).

If $R$ is a ring (here always meaning a commutative
integral domain with identity), we call a function
$f:R\to\N$ a
\textbf{Euclidean function}\index[S]{Euclidean function} or a
\textbf{Euclidean algorithm}\index[S]{Euclidean algorithm}
(EA for short) on $R$ if the following conditions hold:
\begin{enumerate}
\item[](E--1)\ $\forall a\in R: f(a)=0 \Leftrightarrow a=0$
\item[](E--2)\ $\forall a\in R, b\in R\setminus\{0\}
   \exists c\in R: f(a - bc) < f(b).$
\end{enumerate}

If such an $f$ exists, $R$ is called
\textbf{Euclidean with respect to $f$}.\index[S]{Euclidean!with respect to $f$}
Euclid's aim was to prove the so-called
``fundamental theorem of arithmetic'' in $\Z$: every natural number can be
written uniquely, up to order, as a product of primes.

By Euler's time it had been recognized that certain Diophantine
equations, such as $y^2 = x^3 - 2$, can be solved particularly easily by
instead writing $x^3 = (y + \sqrt{-2})(y - \sqrt{-2})$ and assuming
that the ring $\Z[\sqrt{-2}]$ is a UFD, i.e.\ that an analogue of the
fundamental theorem of arithmetic holds in it just as it does in $\Z$.
At the time, however, this assumption was regarded as self-evident.

It was only Gauss\index[N]{Gauss} who, studying the ring $\Z[i]$
(where $i^2 = -1$) for his theory of biquadratic residues, recognized
the need for a proof, and who showed, to this end, that $\Z[i]$ is
Euclidean with respect to the norm.

Afterwards, the theory of quadratic and biquadratic residues
was generalized to $p$-th power residues, and Fermat's last theorem
--- according to which the Diophantine equation $x^p + y^p + z^p = 0$,
for prime $p \ge 3$, has only trivial solutions (that is, solutions with
$xyz = 0$) --- was proved in special cases.

Before Kummer\index[N]{Kummer} developed his theory of ideal
numbers (later expanded by Dedekind\index[N]{Dedekind} into
present-day ideal theory), both tasks required assuming that the
rings $\Z[\zeta_p]$, where $\zeta_p$ is a primitive $p$-th root of
unity, are UFDs. For $p = 3, 5, 7$ this could be proved by showing
that the corresponding rings are norm-Euclidean (see Gauss
\cite{Gau76}\index[N]{Gauss} and Kummer \cite{Kum44}).

A more detailed survey of the history of the Euclidean algorithm, up to
the more recent results of Weinberger\index[N]{Weinberger} and
Lenstra\index[N]{Lenstra}, can be found in Lenstra's\index[N]{Lenstra}
beautifully written
article \cite{Len79b}; further historical
details are given by van der Linden\index[N]{van der Linden}
\cite{Lin85}, as well as in the somewhat older work of
Narkiewicz\index[N]{Narkiewicz} \cite{Nar67}.  We now return to our
definition of a Euclidean function: if $f:R\to\N$ is a function
satisfying (E--1) but not necessarily (E--2), we set
\[
M(f) := \inf \left\{ x \in \R : \forall a,b\in R\setminus\{0\}
\ \exists c \in R: f(a-bc) < x \cdot f(b) \right\}
\]
and $M(f)=\infty$ if no such $x$ exists. We call $M(f)$
the \textbf{Euclidean minimum}\label{pEuM}\index[S]{Euclidean minimum}
of $R$ with respect to $f$.
If $R$ is a number ring and $f$ is the absolute value of the norm, we
also write
$M(R)$ instead of $M(f)$; and if $K$ is a number field with
$R$ the ring of all integers of $K$ (i.e.\ the maximal order), we
usually write $M(K)$ instead of $M(R)$ or $M(f)$.
If $M(f)< 1$, then $R$ is Euclidean with respect to $f$, whereas this
fails for $M(f)>1$. What can happen when
$M(f)=1$ is illustrated by the following examples (which we
shall, however, only be able to justify later):\label{pBei}
\begin{enumerate}
\item Let $R$ be a Euclidean ring and $f$ the ``minimal Euclidean
  function'' on $R$ (which we shall get to know presently);
  then $M(f)=1$ and $R$ is Euclidean with respect to $f$.
\item Let $R$ be the ring of integers of the cubic field of
  discriminant $d=-199$, and $f$ the absolute value of the norm; then
  $M(f) = M(R) = M(K) = 1$, $R$ has class number 1, but $f$ is not
  Euclidean on $R$ (see Taylor\index[N]{Taylor} 1976).
\item Let $R$ be the ring of integers of one of the fields
  $K=\Q(\sqrt{65})$ (see e.g.\ Heinhold\index[N]{Heinhold} 1939),
  $K=\Q(\sqrt{-1},\sqrt{15})$ or
  $K=\Q(\sqrt{-3},\sqrt{13})$ (see §6 below); in these cases
  $M(K)=M(R)=1$, but $K$ has class number 2, so $f$ cannot
  be Euclidean on $R$.
\item Let $R=\Z[\sqrt{-3}]$ (see Cohn\index[N]{Cohn} 1978),
  $R=\Z[\sqrt{5}]$ or $R=\Z[\sqrt{13}]$ and $f$ the
  absolute value of the norm; then in all three cases $M(R)=1$,
  although here $R$ is not even integrally closed, let alone a
  UFD, and certainly not Euclidean. We shall later find further examples
  of such behaviour in number fields of higher degree.
\end{enumerate}
These examples show that, in general, $M(f)=1$ does not let us
conclude whether or not $R$ is Euclidean with respect to $f$.
A lower bound for the Euclidean minimum
$M(f)$ -- sometimes even an exact one -- is given by the following
criterion, which, despite its simplicity,
does not appear to have been stated explicitly before (not even for
the case where $f$ is the absolute value of the norm):

\begin{quote}
  {\bf (0.1)} {\em Let $b \in R \setminus \{0\}$ not be a unit;
    then $M(f) \geq \frac1{f(b)}$.}
\end{quote}

\begin{proof}
  Since $b$ is not a unit, there exists an $a \in R$ with $b \nmid a$;
  if we then write $a = bq + r$ for arbitrary $q, r \in R$, then $r$
  must be $\neq 0$. By (E--1) we therefore have $f(r) \geq 1$ for every choice of
  $q$ and $r$, and hence $\frac{f(r)}{f(b)} \geq \frac{1}{f(b)}$.
  The claim now follows from the definition of $M(f)$.
\end{proof}
If $K$ is a quadratic number field with discriminant $d$, then the
table below gives a $b \in R$ of minimal norm $>1$ and the
resulting bound for $M(K)$:
\begin{center}
$$ \begin{array}{ccc}
\toprule
d & b & M(K) \geq \\
\midrule
\rsp -4 & 1+i & \frac{1}{2} \\
\rsp -3 & \sqrt{-3} & \frac{1}{3} \\
\rsp 5 & 2 & \frac{1}{4} \\
\rsp 8 & \sqrt{2} & \frac{1}{2} \\
\rsp 12 & 1+\sqrt{3} & \frac{1}{2} \\
\rsp 13 & 4+\sqrt{13} & \frac{1}{3} \\
\rsp 17 & \frac{5+\sqrt{17}}{2} & \frac{1}{2} \\
\bottomrule
\end{array} $$
\end{center}

We shall see in §2 that these bounds are in fact exact, and moreover
that these are the only quadratic fields for which the bound
in (0.1) is exact. It is striking that these fields can also
be characterized by $|d|\leq 4$ in the complex case, or by $d\leq 17$ in the real
case -- that is, the fields for which the bound in
(0.1) is exact are precisely those of minimal discriminant!
Similar observations will later be possible for number fields of
degree three and four as well.
We now turn briefly to the general theory of Euclidean
rings. The inclusions
\[
\{\text{Euclidean rings}\} \subset \{\text{principal ideal domains}\}
\subset \{\text{UFDs}\}
\]
are classical. For the number rings considered here the last inclusion
can even be reversed: every UFD here is also a principal ideal
domain. The example $R = \Z[\frac{1+\sqrt{-19}}2]$ shows, however,
that the first inclusion is proper: $R$ is a principal ideal
domain, but, as we shall see, not Euclidean.

To this end we follow an idea of Motzkin\index[N]{Motzkin} \cite{Mot49}
(see also Samuel\index[N]{Samuel} \cite{Sam71}) and try to construct, for a
given ring, a function $f$ that is
Euclidean on $R$. By (E--1) we already know all $a \in R$ with
$f(a)=0$ (only $a=0$), so we ask for which $b \in R$ we may
set $f(b)=1$. By (E--2) this is possible only if,
for every $a \in R$, there exists a $q \in R$ with $f(a-bq) < f(b) = 1$;
this forces $f(a-bq)=0$, and by (E--1) we must then have $a-bq=0$.
This shows that a $b \in R$ with $f(b)=1$ must divide every $a \in R$,
so $b$ must be a unit in $R$. Conversely, if $b$ is
a unit, then (E--2) is certainly satisfied for every given $a \in R$ by taking
$q=ab$, so we may indeed set $f(b)=1$.
If $f(b)=2$, then for every $a \in R$ there should exist a $q \in R$
with $f(a-bq)<2$; thus $a-bq=0$ (if $f(a-bq)=0$) or $a-bq$ is
a unit (if $f(a-bq)=1$) -- in other words, every
$a \in R$ is congruent modulo $b$ to either $0$ or a unit. Continuing
in this way, we are led to the following subsets\label{pEi}
$E_i$ of $R$:
\begin{align*}
  E_0 & = \{0\}, \quad E_1 = E_0 \cup R^\times , \quad
          \text{and in general for all $i \geq 1$} \\
  E_i \setminus E_{i-1} & = \{ b \in R : \text{every residue class} \bmod b
  \text{ contains an element of } E_{i-1} \}.
  \end{align*}
Finally we set
$$ E_\infty = \bigcup_{i \geq 0} E_i. $$
We then have

\begin{quote}
  {\bf (0.2)} {\em A ring $R$ is Euclidean if and only if
    $R = E_\infty$.}
\end{quote}

\begin{proof}
  Let $R = E_\infty$. Define $f_0(a) = \min \{ i \in
  \N : a \in E_i \}$, giving a function $f_0 : R \to
  \N$. Given
  $a, b \in R \setminus \{0\}$ with, say,
  $f_0(b) = i$, we have $i > 1$ since $b \neq 0$, and by
  construction of $f_0$ we have $b \in E_i$. Hence there exists an $r \in
  E_{i-1}$ with $a \equiv r \bmod b$; writing $a = bq + r$, we have
  $q \in R$ and $f_0(r) \leq i-1 < f_0(b)$. Hence $R$ is Euclidean
  with respect to $f_0$.
  Conversely, let $R$ be Euclidean with respect to a function $f$. We
  claim that every $a \in R$ with $f(a) = i$ already lies in $E_i$.
  For $i = 0$ this holds by (E--1).
  Suppose the claim holds for all $i < k$, and let
  $f(b) = k$. Then for every $a \in R$ there exist $q, r \in R$
  with $a = bq + r$ and $f(r) < f(b)$. By the induction hypothesis we have
  $r \in E_{k-1}$, i.e.\ for every $a \in R$ there is an $r \in E_{k-1}$
  with $a \equiv r \bmod b$; this in turn implies $b \in E_k =
  E_{k}(b)$, and consequently $R \subset E_\infty \subset R$.
\end{proof}
Looking again at the second part of the proof, we notice
that we have shown rather more: if $R$ is Euclidean with respect to
$f$, and $f_0$ is the function defined in the first part of the proof,
then the relation $b \in E_{f(b)}$ shows that $f_0(b) \leq f(b)$ for
all $b \in R$, since $f_0(b)$ is by definition the smallest index $i$
for which $b \in E_i$. This leads us to the following definition:
if $R$ is a Euclidean ring, the function $f_m :
R \to \N$ defined by
\[
f_m(a) = \min \{ f(a) : f \text{ is a Euclidean function on } R \}
\]
is called the
\textbf{minimal Euclidean function}\index[S]{Euclidean function!minimal}
on $R$. This name is justified by

\begin{quote}
  {\bf(0.3)} {\em $f_m$ is a Euclidean function on $R$ with the
    property $f_m(a) \leq f(a)$ for all $a \in R$ and all
    Euclidean functions $f$ on $R$. In particular $f_m = f_0$.}
\end{quote}

\begin{proof}
  Let $a, b \in R \setminus \{0\}$. By definition of
  $f_m$ there exists a Euclidean function $f$ with $f_m(b) =
  f(b)$. Hence there exist $q, r \in R$ with $a = bq + r$ and $f(r) <
  f(b)$, so $f_m(r) \leq f(r) < f(b) = f_m(b)$, and
  $f_m$ is indeed a Euclidean function on $R$.
  Since $f_m(a) \leq f(a)$ for every Euclidean function $f$ on
  $R$, and $f_0(a) \leq f_m(a)$ by the remark following (0.2), we must have
  $f_m(a) = f_0(a)$ for all $a \in R$.
\end{proof}

For $R = \Z$, the sets $E_i$ can be described quite simply: we
have $E_0 = \{0\}$, $E_1 = \{0, -1, 1\}$, and $E_2 \setminus E_1$
consists of those $a \in \Z$ for which every residue class modulo $a$
has a representative in $E_1$; since $E_1$ has only three elements,
we need $|a| \le 3$, giving $E_2 = \{0, \pm 1, \pm 2, \pm 3\}$.
By induction it is easy to see that
$E_i = \{ a \in \Z : |a| \leq 2^i - 1 \}$. In particular this gives
$\Z = E_\infty$, i.e.\ $\Z$ is Euclidean.

Matters are somewhat more complicated in imaginary quadratic
number fields. If $m \neq -1, -3$, then $D(m)$ (we shall
from now on write $D(m)$ for the ring of integers of $\Q(\sqrt{m})$),
for $m<0$, is well known to contain only the units $+1$ and $-1$; in
these cases we therefore have $E_1 = \{0, +1, -1\}$. $E_2$ then consists of
all $a \in D(m)$ for which every residue class modulo $a$ contains one of the three
elements $0, 1, -1$. Since $|R/aR| =
N_{K/\Q}(a)$, this is clearly possible only if
$N_{K/\Q}(a) \le 3$. The only imaginary quadratic
rings containing elements of norm 2 or 3 are, however, $D(m)$ for
$m=-1, -2, -3, -7, -11$; for every other $D(m)$ we therefore have $E_1 = E_2 =
\dots = E_\infty = \{0, 1, -1\}$, and we obtain

\begin{quote}
  {\bf (0.4)} {\em Let $m<0$ and $R=D(m)$; then $R$ is Euclidean
    precisely for the values $m=-1, -2, -3, -7, -11$, and in these
    cases the norm is a Euclidean function. For all other values of
    $m$ we have $E_1 = E_2 = \dots = E_\infty = \{0, 1, -1\}$.}
\end{quote}

In particular, the principal ideal domains $D(m)$ with $m=-19, -43,
-67, -163$ are not norm-Euclidean. Gauss\index[N]{Gauss}
showed that $D(m)$ is norm-Euclidean for the values of $m$ noted
above ($m=-1$, $-3$), while the remaining cases are due to
Dedekind\index[N]{Dedekind} \cite{Ded00} and
Dickson\index[N]{Dickson} \cite{Dic27}. Both of the latter authors have
also remarked that there is no other norm-Euclidean $D(m)$ with $m<0$.

The more far-reaching statement (0.4), on the other hand, is due to
Motzkin\index[N]{Motzkin} \cite{Mot49} and, independently, to
Dubois\index[N]{Dubois} and Steger\index[N]{Steger}
\cite{DS58}. Narkiewicz\index[N]{Narkiewicz} writes
\cite[p.\ 175]{Nar67} that Dubois\index[N]{Dubois} and
Steger\index[N]{Steger} had moreover shown that every Euclidean
function $f$ on $D(m)$, $m<0$, must coincide with the norm --- a claim
that is neither historically accurate (they never made it) nor
mathematically true (for $m=-1$, $-2$, $-3$, $-7$, $-11$ the minimal
Euclidean function is also a Euclidean algorithm on $D(m)$, and it
differs from the norm). Further proofs of (0.4) have been given by
Stewart\index[N]{Stewart} and Tall\index[N]{Tall} \cite{ST79} in their
textbook {\em Algebraic Number Theory} and by Castro
Chadid\index[N]{Castro Chadid} \cite{Cha84}. Finally
Campoli\index[N]{Campoli} \cite{Cam88} has once again shown that
$D(-19)$ is a principal ideal domain but not Euclidean. What all these
proofs have in common is the construction of a non-unit of norm $\le 3$ in a
Euclidean $D(m)$.

For the rings $D(-1)$ and $D(-3)$, Lenstra\index[N]{Lenstra}
\cite{Len74} has determined the sets $E_i$ explicitly; beyond these
results, we know only a few elementary properties of the $E_i$:

\begin{quote}
  {\bf (0.5)} {\em Let $R$ be a number ring; then we have
    \begin{enumerate}
    \item[(i)] If $a \in E_\infty$, then every divisor of the ideal $aR$
      is a principal ideal.
    \item[(ii)] If $a = p_1 p_2 \dots p_m$ is the prime factorization of $a \in
      R$, then $a$ does not lie in $E_m$.
    \item[(iii)] For all $i \in \N$ we have $E_i \neq R$.
    \item[(iv)] We have $E_i R^\times = E_i$, i.e.\ with $a \in E_i$ and $u
      \in R^\times$ we also have $au \in E_i$.
    \item[(v)] If $K/\Q$ is Galois, then $E_i^G = E_i$ for all
      $G \in \operatorname{Gal}(K/\Q)$.
    \end{enumerate} }
\end{quote}

\begin{proof}
  (i) Let $a \in E_1$; if then $P$ is a prime ideal in $R$ that
  divides $(a)$, then there exists a $p_0 \in R$ with $P = (a, p_0)$
  (for every ideal in a number field is generated by at most two
  elements, and we may choose $a \in P$ arbitrarily; see
  e.g.\ Cohn\index[N]{Cohn} \cite{Coh78}).

  Since $a \in E_1$ there now exists
  a $p_1 \in E_{k-1}$ with $p_1 \equiv p_0 \mod a$. Thus $P =
  (a, p_0) = (a, p_1)$. Replacing in the above reasoning $a \in E_1$
  by $p_1 \in E_{k-1}$, we obtain $P = (p_1, p_2)$ for a $p_2 \in
  E_{k-2}$ etc. Finally we obtain $P = (p_1, p_2, \dots,
  p_{k-1})$ for a $p_{k-1} \in E_1$ (possibly this already occurs
  earlier). Since $P$ is prime, $p_{k-1}$ cannot be a unit.
  Consequently $p_{k-1} \neq 0$ and $P = (p_{k-1})$ is a
  principal ideal.

  \medskip\noindent
  (ii) The claim is true for $m=0, 1$; suppose it has now been proved
  for all $i \in \N$ with $i<k$ and let $a \in E_k$. We assume that the
  number $m$ of factors of $a$ is greater than or equal to $k$ and choose
  a prime factor $p$ of $a$. Then $a/p \equiv e \mod a$ for some
  $e \in E_{k-1}$, and we have $e \neq 0$ (otherwise $ap \mid a$). Now
  however $e \equiv 0 \mod (a/p)$, i.e.\ $e \in E_{k-1}$ has $k-1$
  prime factors, and this contradicts the induction hypothesis.

  \medskip\noindent
  (iii) follows immediately from (ii) if we observe that there are in $R$
  infinitely many prime ideals (at least one above every
  rational prime).

  \medskip\noindent
  (iv) likewise follows immediately because $R/(a) \cong R/(au)$.

  \medskip\noindent
  (v) The claim holds for $j=0$. If it holds for $j+k-1$ and $a \in
  E_k$, then for every $b \in R$ there exists a $c \in E_{k-1}$ with
  $b \equiv c \mod a$. Hence $b^G \equiv c^G \mod a^G$, and by the
  induction hypothesis we have $c^G \in E_{k-1}^G$. Since $b^G$ runs
  through all of $R$ as $b$ does, it follows that $a^G \in E_k^G$,
  proving the claim.
\end{proof}

We now ask whether there are number rings that are Euclidean with
respect to $f_0$ but not with respect to the norm. To this day not a
single such number ring is known. It would, however, be somewhat hasty
to draw any mathematical conclusion from this: if certain
Riemann hypotheses are true, every number ring of class number 1 and
with infinite unit group is Euclidean with respect to $f_0$; this
result is due to Weinberger\index[N]{Weinberger} \cite{Wei73} and
rests on a generalization of Hooley's ``proof'' of Artin's conjecture
on primitive roots (Hooley had proved this conjecture in 1967 assuming
the Riemann hypothesis).  For a proof of this statement
we refer to Weinberger\index[N]{Weinberger} \cite{Wei73} and
Lenstra\index[N]{Lenstra} \cite{Len74}, and for more details on Artin's
conjecture to Gupta\index[N]{Gupta}, Murty\index[N]{Murty} \& Murty
\cite{GMM87}, Murty \cite{Mur88} and Narkiewicz\index[N]{Narkiewicz}
\cite{Nar88}. Here we wish only to indicate briefly what the minimal
Euclidean function has to do with primitive roots.

To this end, let $a \in R$ be prime; if the unit $u \in R^\times$
is then a primitive root modulo $a$, then $a \in E_2$, because the powers
of $u$, together with $0$, form a complete residue system modulo $a$ and
lie as units in $E_1$. If $u$ is a primitive root for ``many'' $a \in
R$, then $E_2$ becomes ``large'', and we have a good chance of proving
$R = E_\infty$. In fact, Weinberger\index[N]{Weinberger} has shown that
whether or not $R = E_\infty$ can already be read off from $E_3$.

\subsection*{$k$-stage Euclidean rings}
We now come to the definition of $k$-stage Euclidean rings
(see the works of Cooke\index[N]{Cooke} \cite{Coo77}, as well as
Cooke\index[N]{Cooke} and Weinberger\index[N]{Weinberger} \cite{CW75}).
Let $a,b \in R \setminus \{0\}$; a sequence of equations
\begin{align*}
a &= q_1 b + r_1 \\
b &= r_1 q_2 + r_2 \\
r_1 &= r_2 q_3 + r_3 \\
&\vdots \\
r_{k-2} &= r_{k-1} q_k + r_k,
\end{align*}

in which, without loss of generality, at least one of every two
consecutive $q_i$ is different from $0$, and $q_k$ in particular is
different from $0$, is called a
\textbf{division chain}\index[S]{division chain} of length $k$
starting from $(a,b)$; if $r_k = 0$, we call it
\textbf{terminating}.\index[S]{division chain!terminating} A map $f: R
\to \N$ is called a
\textbf{$k$-stage Euclidean function}\index[S]{Euclidean function!$k$-stage}
on $R$ if (E--1) is satisfied and, for all $a,b \in R \setminus \{0\}$,
there exists a division chain of length $k$ starting from $(a,b)$ with
$f(r_k) < f(b)$. If such an $f$ exists, $R$ is called
\textbf{$k$-stage Euclidean}\index[S]{Euclidean!$k$-stage} with
respect to $f$. Finally, we call $R$
\textbf{quasi-Euclidean}\index[S]{quasi-Euclidean} if for all $a,b
\in R \setminus \{0\}$ there exists a terminating division chain
starting from $(a,b)$ (note that being quasi-Euclidean
does not depend on any choice of function). In place of quasi-Euclidean,
Cooke\index[N]{Cooke} used the term $\omega$-stage Euclidean.

We first show that the definitions of quasi-Euclidean rings appearing
in the literature (e.g.\ in Cooke\index[N]{Cooke} \cite{Coo77},
Bougaut\index[N]{Bou76} \cite{Bou76,Bou80} or
Leutbecher\index[N]{Leutbecher} \cite{Leu77}) are all equivalent:

\begin{quote}
  {\bf (0.6)} {\em Let $R$ be a ring; then the following are equivalent:
    \begin{enumerate}
    \item[(i)] $R$ is quasi-Euclidean;
    \item[(ii)] for all $a,b \in R \setminus \{0\}$ there exists a $k \in
      \N$ and a terminating division chain starting from $(a,b)$;
    \item[(iii)] there exists a function $\Phi: R \times R \to \N$
      such that for all $a,b \in R \setminus \{0\}$ there exist numbers $q,r \in R$
      with $a=bq+r$ and $\Phi(b,r) < \Phi(a,b)$.
    \end{enumerate}}
\end{quote}

\medskip\noindent
\textbf{Rem.:} Functions as described in (iii) we call
\textbf{quasi-Euclidean
  functions}\index[S]{quasi-Euclidean function} on $R$.
\begin{proof}
  (i) $\Rightarrow$ (iii): We define
  $\Phi: R \times R \to \N$ by
  \begin{align*}
    \Phi(a,b) & = \min \{k \in \N: \text{there exists a terminating }
                   k\text{-stage division chain} \\
              & \quad \text{starting from } (a,b) .\}
  \end{align*}
  It is easily verified that $\Phi$ is a quasi-Euclidean function.

  (iii) $\Rightarrow$ (i): Let $a,b \in R \setminus \{0\}$ be
  given. Then there exist $q,r \in R$ with $a=bq+r$ and
  $\Phi(b,r) < \Phi(a,b)$. If $r = 0$, we are done; otherwise
  there exist $q_2,r_2 \in R$ for $b,r$ with $b = q_2 r + r_2$ and
  $\Phi(r, r_2) < \Phi(b,r)$. Repeating this step sufficiently often
  yields a terminating division chain starting from $(a,b)$.

  (i) $\Rightarrow$ (ii) is clear.

  (ii) $\Rightarrow$ (i): Let $a,b \in R \setminus \{0\}$ be given.
  Then there exist a $k \in \N$ and a $k$-stage division chain
  starting from $(a,b)$ with $f(r_k) < f(b)$. If already $r_k = 0$, we
  are done. Otherwise we repeat this step with the pair
  $(r_{k-1},r_k)$, obtaining a division chain of length $l$ starting
  from $(r_{k-1},r_k)$ with $f(r_l) < f(r_k)$, and so on. Concatenating
  these division chains yields a terminating division chain
  starting from $(a,b)$.
\end{proof}

The implication (ii) $\Rightarrow$ (i) shows that every $k$-stage
Euclidean ring is also quasi-Euclidean.

If we consider a division chain starting from $(a,b)$ with $r_{k-1}
\neq 0$ and $r_k = 0$, we find $(a,b) = (r_{k-1})$.  By induction
we then conclude that in quasi-Euclidean rings every finitely
generated ideal is a principal ideal.

Quasi-Euclidean rings need not, however, be principal ideal
domains, as the following example shows: let $A$ be the ring of all
algebraic integers; every finitely generated ideal in $A$ is indeed a
principal ideal, yet $A$ is well known not to be a principal ideal
domain. In particular $A$ cannot be Euclidean; on the other hand $A$
is 2-stage Euclidean with respect to every function satisfying
(E--1): if $a,b \in A$ are coprime (i.e.\ $(a,b)=A$), then by
Lenstra\index[N]{Lenstra} \cite{Len73} there exists a $q \in A$ such
that $a-bq$ is a unit in $A$. Such $a,b$ therefore have a terminating
division chain of length $\le 2$. If $a$ and $b$ are not
coprime, choose a $d \in A$ with $(d)=(a,b)$, then set $a'=a/d$,
$b'=b/d$, and multiply the terminating division chain of length 2
starting from $(a',b')$ by $d$.

To study $k$-stage Euclidean functions we introduce an
analogue of $M(f)$, defined as follows:\label{DefMk}
\begin{align*}
M^k(f) := \inf \Bigl\{ \kappa \in \R \;\Big|\;
    &\text{for all } a,b \in R \setminus \{0\}
      \text{ there exists a division chain }  \\
    &\text{ starting from } (a,b) \text{ with }
      f(r_k) < \kappa \cdot f(b) \Bigr\}.
\end{align*}

Since the notations $M_2(f)$, $M_3(f)$, etc., for the second, third,
and further minima of $f$ are already almost standard, we have written
the superscript $M^k(f)$
above. We thus have the chain of inequalities
 $$ M(f) = M^1(f) \ge M^2(f) \ge M^3(f) \ge \dots \ge M^\infty(f) =
           \lim_{k \to \infty} M^k(f) $$
(this limit exists, since $M^k(f)$ is a monotonically
decreasing sequence of real numbers bounded below by $0$). For
describing $k$-stage division chains in a ring $R$, continued
fractions with coefficients from $R$ have proved useful: given a
sequence $q_1, \dots, q_k \in R$ (where of every two
consecutive $q$'s at least one, and in particular $q_k$, is different
from $0$), we define
\[
[q_1] = q_1 = \frac{a_1}{b_1} \quad \text{with } a_1 = q_1, \quad b_1 =
1
\]
\[
[q_1, q_2] = q_1 + \frac{1}{q_2} = \frac{a_2}{b_2} \quad \text{with }
a_2 = q_1 q_2 + 1, \quad b_2 = q_2.
\]
\[
[q_1, \dots, q_k] = \frac{a_k}{b_k} \quad \text{with } a_k = q_k
a_{k-1} + a_{k-2}, \quad b_k = q_k b_{k-1} + b_{k-2}
\]
We readily see that $[q_1, \dots, q_k] = [q_1, \dots, q_{k-1} + 1/q_k]$,
and induction then gives
\[
[q_1, \dots, q_k] = q_1 + \cfrac{1}{q_2 + \cfrac{1}{q_3 + \dots +
    \cfrac{1}{q_{k-1} + \cfrac{1}{q_k}}}}
\]
This last identity in turn immediately gives
$[q_1, \dots, q_k] = q_1 + 1/[q_2, \dots, q_k]$.
Writing $a_k$ for the numerator of the continued fraction and
$b_k$ for its denominator, we find

\begin{quote}
  {\bf (0.7)} {\em The denominator of $[q_1, \dots, q_k]$ is equal to the
    numerator of $[q_2, \dots, q_k]$.}
\end{quote}

By induction we find:

\begin{quote}
  {\bf (0.8)} {\em If $a,b \in R \setminus \{0\}$ and $r_k$ is the
    last remainder of the division chain starting from $(a,b)$ with the
    quotients $q_1, \dots, q_k$, then we have}
    $$ \frac{a}{b} - \frac{a_k}{b_k} = (-1)^k \frac{r_k}{b \cdot b_k}. $$
\end{quote}

If $f$ is multiplicative, we can extend it
multiplicatively from $R$ to $K$, and by (0.8) we find $f(r_k) < f(b)$
if and only if $f(a/b - a_k/b_k) < 1$; thus we obtain

\begin{quote}
  {\bf (0.9)} {\em A ring $R$ is $k$-stage
    norm-Euclidean with respect to a multiplicative function $f$ if and only if
    for every $a/b \in K$ there exists a continued fraction $[q_1, \dots, q_k] =
    a/b$ of length $l \le k$ such that the inequality
    $$ f\left( \frac{a}{b} - \frac{a_k}{b_k} \right) <
       f\left( \frac{1}{b_k} \right) $$
    is satisfied.}
\end{quote}

The difference from the ordinary EA can thus be formulated as
follows: in Euclidean rings an $x \in K$ can be approximated by a $y
\in R$ so that $f(x-y)<1$; in a $k$-stage Euclidean ring,
on the other hand, we may choose $y$ from the larger set, containing $R$, of
continued fractions of length $\le k$, though the approximation must
then be ``more precise''.

In general it is rather difficult to recognize a continued
fraction (even of length 2) as such; for example,
\[
\frac{\sqrt{14}}{2} = -2 + \frac{1}{4 - \sqrt{14}} = [-2, 4 - \sqrt{14}]
\]
is a continued fraction of length 2 with denominator $4 - \sqrt{14}$,
whereas $\frac{1+\sqrt{14}}{2}$ is not a continued fraction of length
2. To prove this, we use

\begin{quote}
  {\bf (0.10)} {\em If $a_k/b_k$ is a continued fraction and we have
    $a/b = a_k/b_k$, then $b_k \mid b$.}
\end{quote}

\begin{proof}
  By induction we show
  $a_k b_{k-1} - a_{k-1} b_k = (-1)^k$.
  Multiplying this identity by $b$ yields
  $b a_k b_{k-1} - b a_{k-1} b_k = (-1)^k b$.
  Since $a b_k = a_k b$, $b_k$ divides the
  left-hand side, and hence also the right-hand side, of this equation.

  If $(1 + \sqrt{14})/2 = [q_1, q_2] = (1 + q_1 q_2)/q_2$ were a
  continued fraction of length 2, we would need $q_2 \mid
  2$; but since $(1 + \sqrt{14})/2$ cannot be reduced further, we
  also need $2 \mid q_2$. Consequently $q_2 = 2e$ for a unit $e \in
  R^\times$. However, $u = 15 + 4\sqrt{14}$ is the fundamental unit
  of $R$, so every unit is $\equiv 1 \bmod 2$. In particular $e
  \equiv 1 \bmod 2$, and we would need $(1 + \sqrt{14})e = 1
  + \sqrt{14} \equiv 1 + q_1 q_2 \bmod 2$, which is evidently false.
\end{proof}

To derive a criterion that lets us write certain numbers
as continued fractions, we introduce the sets\label{Ej'} $E_j'$,
defined as follows: $E_0' = \{0\}$, and
\begin{align*}  
  E_j' \setminus E_{j-1}'
  & = \{ a \in R : \text{every prime residue class } \bmod a \\
  & \qquad \text{ has a representative from } E_{j-1}' \}.
\end{align*}
Then $E_0' = E_0$, $E_1' = E_1$, and
$E_j' \subseteq E_j$ for $j \ge 2$, where the last inclusion is
in general proper: for $D(-5)$ we have
$E_1 = E_2 = \dots = E_\infty = \{0, -1, +1\}$, while
$E_1' = E_2' = \dots = E_\infty' = \{0, \pm 1, \pm 1 \pm \sqrt{-5}\}$.
The criterion promised above now reads

\begin{quote}
  {\bf (0.11)} {\em Let $a, b \in R \setminus \{0\}$, $b \in E_k'$,
    and $(a, b) = 1$; then $a/b$ is a continued fraction of length $\le k$
    with denominator $ub$, where $u \in R^\times$ is a unit in $R$. }
\end{quote}

\begin{proof}
  If $k=1$, then $b$ is a unit and $a/b = q_1 \in R$ is a
  continued fraction of length 1 with denominator $1 = bu$, where $u = \frac1b$
  is a unit in $R$.
 
  Now suppose the claim holds for all $j < k$ and $b \in E_k'$.
  Since $(a, b) = 1$, $a$ lies in a prime residue class mod
  $b$, so there exists an $e \in E_{k-1}'$ with $a \equiv e \pmod{b}$,
  i.e.\ a $q_1 \in R$ with $a = q_1 b + e$. We then see
  that $(b, e) = (b, a) = 1$ and $a/b = q_1 + e/b = q_1 + 1/(b/e)$.
  By the induction hypothesis, $b/e$ is a continued fraction of length
  $\le k-1$ whose denominator (and hence also whose numerator) coincides
  with $e$ (respectively with $b$) up to a factor $u \in R^\times$.
  Writing $b/e = [q_2, \dots, q_k]$, we have $a/b =
  [q_1, \dots, q_k]$, and (0.7) shows that $a/b$ is a continued fraction of
  length $\le k$ with denominator $bu$.
\end{proof}

As an example we show how to apply (0.11) to
write $(1+\sqrt{14})/2$ as a continued fraction. First we note that
$(1+\sqrt{14})/2 = -1 + (3+\sqrt{14})/2$ (this makes the norm of the
denominator small), and then claim that $a = 3+\sqrt{14} \in E'_2$.
\label{p11} For this we must show that every prime
residue class mod $a$ contains a unit. We see this as follows: we have
$u = 15+4\sqrt{14} \equiv 3 \pmod{a}$ (simply note that
$\sqrt{14} \equiv -3 \pmod{a}$), so $u^2 \equiv 4$, $u^3 \equiv 2$,
$u^4 \equiv 1 \pmod{a}$. Since $N_{K/\mathbb{Q}}(a) = 5$, $u$ is thus a
primitive root mod $a$.

In particular we now have $2 \equiv -u \pmod{a}$, and we obtain
\[
\frac{2}{a} = \frac{-u}{a} + 1 + \sqrt{14} = 1 + \sqrt{14} + \frac{1}{-a/u}.
\]
Since $-a/u = 11 - 3\sqrt{14}$ we have
$(1 + \sqrt{14})/2 = [-1, 1 + \sqrt{14}, 11 - 3\sqrt{14}]$,
i.e.\ $(1 + \sqrt{14})/2$ is a
continued fraction of length 2.

We now prove an important property of the sets $E'_j$ (see §1), which
looks more trivial than it is:

\begin{quote}
  {\bf (0.12)} {\em If $c \in E'_k$ and $b \mid c$ for some $b \in R$,
    then also $b \in E'_k$.}
\end{quote}

\begin{proof}
  We shall show a little more: if $B, C$ are integral ideals in $R$ (not
  necessarily principal) and $B \mid C$, then: if every prime
  residue class mod $C$ has a representative in $E'_{k-1}$, the same holds
  for the prime residue classes mod $B$. By induction on the
  number of prime ideal divisors of $CB^{-1}$ we may assume $C=BP$
  for a prime ideal $P$.
  
  Let $r \in R$ be given with $(r) \cdot B = R$; we seek an
  $e \in E'_{k-1}$ with $r \equiv e \pmod{B}$. We
  distinguish two cases: a) $B \cdot P = R$: then, by the Chinese
  remainder theorem, there exists an $s \in R$ satisfying the
  congruences $s \equiv r \pmod{B}$ and $s \equiv 1 \pmod{P}$. With
  this $s$ we then have $(s) \cdot C = (s) \cdot BP = R$, and since
  every prime residue class mod $C$ has, by assumption, a representative
  $e \in E'_{k-1}$, we get $s \equiv e \pmod{C}$, and a fortiori
  $r \equiv s \equiv e \pmod{B}$. b) $B \equiv 0 \pmod{P}$: then also
  $(r) \cdot BP = R$, and we immediately find an $e \in E'_{k-1}$ with
  $r \equiv e \pmod{C}$, whence $r \equiv e \pmod{B}$, as required.
\end{proof}

We next prove an analogue of (0.1):

\begin{quote}
  {\bf (0.13)} {\em Let $b \in R \setminus E'_2$; then
        $M^2(f) \ge \frac{1}{f(b)}$.}
\end{quote}

\begin{proof}
  Let $b \in R \setminus E'_2$. Then there exists an $a \in R$ with
  $(a,b)=1$ such that $a$ is not congruent to a unit mod $b$. We set
  $a = bq_1 + r_1$, $b = r_1 q_2 + r_2$; if $r_2=0$, then $r_1 \mid b$
  and $r_1 \mid a$, hence $r_1 \mid (a,b)=1$. Thus $r_2$ would be a
  unit, contradicting $a \equiv r_1 \pmod{b}$ and the choice of $a$.
  The claim now follows as in (0.1).
\end{proof}

(0.13) is best possible in the sense that the given bound
is sometimes exact. For $R = D(6)$, for instance, we have $M^2(f) = \frac14$,
where $f$ is the absolute value of the norm, and the inequality
$M^2(f) \ge \frac14$ follows from (0.13) because $2 \in E'_2$.

We now notice that we can write (0.1) and (0.13) in the following
form: if $b \in E'_k$, then $M^k(f) \ge 1/f(b)$. The case $k=1$
corresponds to (0.1), the case $k=2$ to (0.13). Whether this
statement also holds for $k=3$ is doubtful.

A further interesting notion for describing $k$-stage
Euclidean rings is due to Cooke and Weinberger \cite{CW75}: let $R$ be
a ring (here an integral domain with identity) and $K$ its field of
fractions; if there exists a division chain of length $\le k$
starting from $(a,b)$ with $f(r_k) < f(b)$ for some $a,b \in R \setminus \{0\}$,
then $K$ is called \textbf{$k$-stage Euclidean} at $x=a/b$
with respect to $f$.

We denote by $F_k$ the set of all $x \in K$ at which $K$ is $k$-stage
Euclidean with respect to $f$. \label{pFi} Thus we have
$F_1 \subset F_2 \subset \dots \subset F_\infty
  = \bigcup_{k=1}^\infty F_k \subset K$.

Evidently $R$ is $k$-stage Euclidean with respect to $f$ if and only
if already $F_k = K$, and quasi-Euclidean if $F_\infty = K$.  If
there exists a $k \in \mathbb{N}$ such that $F_k = F_\infty$
(regardless of whether or not $F_\infty = K$), we call the
smallest such $k$ the \textbf{Euclidean depth} of $K$ with respect to $f$.

In 1975, Cooke and Weinberger were able to prove the following:

\begin{quote}
  {\bf (0.14)} {\em Let $R$ be the ring of integers of an
    algebraic number field $K$ of unit rank $\ge 1$ and let $f$ be
    the absolute value of the norm. If then certain Riemann
    hypotheses are true, then for all $a,b \in R$ with
    $(a,b)=1$ there exists a terminating division chain of
    length $k \le 5$ starting from $(a,b)$. If moreover $K$ has a real embedding, then there exist
    even such of length $k \le 3$.}
\end{quote}

From this theorem they then derive the following corollaries:

\begin{quote}
  {\bf (0.15)} {\em We have $F_5 = F_\infty$: the Euclidean depth of
    $K$ is $\le 5$. If $K$ has a real embedding, then already
    $F_3 = F_\infty$. }
\end{quote}

\begin{proof}
  Assume there exists a terminating division chain of length $k$
  starting from $(a,b)$.  Then $(a,b) = (r_{k-1})$ is a principal
  ideal, and setting $c = a/r_{k-1}$ and $d = b/r_{k-1}$ we have
  $(c,d) = 1$. By (0.14) there then exists a terminating division
  chain of length $\le 5$ (respectively $\le 3$, if $K$ has a real
  embedding) starting from $(c,d)$, and multiplying by $r_{k-1}$
  yields such a chain for $(a,b)$. This shows that $x = a/b
  \in F_k$ already implies $x \in F_5$ (respectively $x \in F_3$), as claimed.
\end{proof}

\begin{quote}
  {\bf (0.16)} {\em If $R$ has class number 1, then $R$ is 4-stage
    norm-Euclidean; if moreover $K$ has a real embedding, then $R$ is
    even 2-stage norm-Euclidean.}
\end{quote}

\begin{proof}
  Let $x \in K$; we may write $x = a/b$ with $a,b \in R$, $(a,b) \neq 1$.
  We seek a division chain of
  length $\le 4$ starting from $(a,b)$ with $|N_{K/\mathbb{Q}}(r_{k-1})| \le
  |N_{K/\mathbb{Q}}(b)|$; for $|N_{K/\mathbb{Q}}(b)| = 1$ there is nothing
  to show. Otherwise (0.14) yields a terminating division chain of length $k \le 5$
  starting from $(a,b)$, in which $r_{k-1}$, as a divisor of $(a,b)$, is a unit, i.e.\
  there exists a division chain of
  length 4 starting from $(a,b)$ with $1 =
  |N_{K/\mathbb{Q}}(r_{k-1})| \le |N_{K/\mathbb{Q}}(b)|$. The claim
  for $K$ with a real embedding follows correspondingly.
\end{proof}

At the end of their paper, Cooke and Weinberger note that, taking into
account a theorem of Lenstra, one obtains a terminating division chain
of length $\le 4$ starting from $(a,b)$ if $R$ has class number $1$.
As in (0.16), it then follows that such rings are already 3-Euclidean.

Cooke seems not to have noticed that (0.15) can be improved quite
simply; in fact we have

\begin{quote}
  {\bf (0.15')} {\em The Euclidean depth is $\le 4$; if $K$ has a
    real embedding, then even $F_2 = F_\infty$. }
\end{quote}

\begin{proof}
  In the proof of (0.15) we saw that $x = a/b$ can be written as
  a continued fraction $a_k/b_k$ of length $k \le 5$. If $b_k
  \in R^\times$ is a unit, then $a_k/b_k$ is a continued fraction of length 1,
  and there is nothing more to show. Otherwise, since
  $$ \frac{a_k}{b_k} - \frac{a_{k-1}}{b_{k-1}} =
     \frac{(-1)^{k-1}}{b_kb_{k-1}}, $$
  $y = a_{k-1}/b_{k-1}$ is a continued fraction of length $k-1
  \le 4$ that approximates $x$ closely enough to show $x \in F_4$.
  The claim for $K$ with a real
  embedding follows correspondingly.
\end{proof}

In particular, quadratic number fields therefore have Euclidean depth
$\le 2$ (if certain Riemann hypotheses are true). In §2 we shall
see, for example, that the fields $\mathbb{Q}(\sqrt{m})$ for $m = 15, 26, 85$
have Euclidean depth 2 and that $M^k(K) = 1$ for all $k\ge2$ there; on
the other hand $\mathbb{Q}(\sqrt{30})$ likewise has Euclidean depth 2,
yet here $M^k(K) = \frac32$ for all $k \ge 2$.

Finally, we wish to give an interpretation of the similarity between the
sets $E_k$ and $E'_k$: to this end we call a map $f:R \to \N$
\textbf{semi-Euclidean} on $R$ if, for all $a,b\in R\setminus\{0\}$, we
have:
\begin{align*}
\text{(E--1)} \quad & f(a)= 0 \Leftrightarrow a = 0; \\ \text{(E'--2)}
\quad & \text{if $(a,b)= 1$ then there exists a $q\in R$ with
  $f(a-bq)<f(b)$.}
\end{align*}
If such an $f$ exists, the ring $R$ is called
\textbf{semi-Euclidean}. As for ordinary Euclidean rings, we now have

\begin{quote}
  {\bf (0.17)} {\em A ring $R$ is semi-Euclidean if and only if
    $R = E'_\infty$. In this case $f_0(a) = \min
    \{j\in\mathbb{N}: a\in E'_j\}$ is a semi-Euclidean function on
    $R$, and this coincides with the minimal semi-Euclidean function defined by
    $$ f_m(a) = \min \{f(a): f
                \text{ is a semi-Euclidean function on } R\}. $$}
\end{quote}

The proof proceeds exactly as in the Euclidean case.

\begin{quote}
  {\bf (0.18)} {\em Let $f$ be multiplicative; then the following are equivalent:
    \begin{enumerate}
    \item[(i)] for all $a,b\in R\setminus\{0\}$ with $(a,b)=1$ there exists
      a $q\in R$ with $f(a-bq)<f(b)$;
    \item[(ii)] for all $a,b\in R\setminus\{0\}$ with $(a,b)=d$ there exists
      a $q\in R$ with $f(a-bq)<f(b)$.
    \end{enumerate} }
\end{quote}

\begin{proof}
  The direction (ii) $\Rightarrow$ (i) is trivial, so assume (i)
  holds. If $(a,b)=(d)$ for some $d\in R$, set $a'=   a/d$,
  $b' = b/d$; by (i) we find a $q\in R$ with $f(a'-qb')< f(b')$.
  Since $f$ is multiplicative, multiplying by $f(d)$ gives
  $f(a-bq)< f(b)$.
\end{proof}

As a corollary we obtain

\begin{quote}
  {\bf (0.19)} {\em A principal ideal domain $R$ is semi-Euclidean
    with respect to a multiplicative function $f$ if and only if
    $R$ is Euclidean with respect to $f$.}
\end{quote}

Thus, to find number rings that are semi-Euclidean but not
Euclidean with respect to the norm, we must search among rings of
class number $> 1$. An elementary, if somewhat lengthy, computation
shows

\begin{quote}
  {\bf (0.20)} {\em If $R=D(m)$, $m\le0$ square-free, then $R$ is
    semi-Euclidean with respect to a function $f$ if and only if $R$ is
    already Euclidean.}
\end{quote}

As an example we consider $D(-5)$: here $E'_1 = \{0, -1, +1\}$,
and for all $a \in E'_2$ we must have $N_{K/\mathbb{Q}}(a) \le 2$ (since
$E'_1$ contains only two non-zero elements). Thus
$N_{K/\mathbb{Q}}(a) \in \{2, 3, 4, 6\}$, and we find $a \in \{+2, +1
\pm \sqrt{-5}\}$; now $\{+1, -1\}$ is a prime residue system mod
$(+1 \pm \sqrt{-5})$, but not mod $2$, so $E'_2 = \{0, +1,
+1 \pm \sqrt{-5}\}$. Considering all $a \in D(-5)$ with
$N_{K/\mathbb{Q}}(a) \le 6$, we find $E'_2 = E'_3 = \dots =
E'_\infty$, and (0.14) shows that $D(-5)$ is not semi-Euclidean.

Incidentally, the proof of (0.20) shows that in the ring
$D(-15)$ we have $E'_2 \subsetneq E'_3 = E'_4$, whereas for every other
$D(m)$ we already have $E'_2 = E'_3$.

The real quadratic case already shows that (0.20) is not typical for
algebraic number fields. For example, Johnson, Queen and Sevilla
\cite{Joh85} showed that $D(10)$ and $D(65)$ are semi-Euclidean rings
of class number 2 with respect to the norm. These are probably even
the only real quadratic rings that are semi-Euclidean with respect to
the norm and have class number different from 1. We shall later see
that the rings of integers of the biquadratic fields
$\mathbb{Q}(\sqrt{-1},\sqrt{15})$ and
$\mathbb{Q}(\sqrt{-3},\sqrt{13})$ also have class number 2 and are
semi-Euclidean.

Analogously, we define $k$-stage semi-Euclidean rings by requiring
that for all $a,b \in R \setminus \{0\}$ with $(a,b) = 1$
there exists a division chain of length $\le k$ starting from $(a,b)$
with $f(r_k) < f(b)$. This also makes clear what is to be
understood by a quasi-semi-Euclidean ring; such rings are actually
called $\mathrm{GE}_2$-rings.

At the end of this section we wish to show how to determine the
Euclidean minima $M(K)$ when $K$ is an imaginary quadratic number field.

\begin{quote}
  {\bf (0.21)} {\em Let $R=D(m)$ be imaginary quadratic; then we have }
    $$ M(K) =
    \begin{cases}
      \dfrac{1+|\,m\,|}{4} & \text{if } m \equiv 2, 3 \pmod{4} \\[1em]
      \dfrac{(1+|\,m\,|)^2}{16|\,m\,|} & \text{if } m \equiv 1 \pmod{4}.
    \end{cases} $$
\end{quote}

In particular we have for the Euclidean rings the following table:
\begin{center}
$$ \begin{array}{c|ccccc}
\toprule
\rsp m & -1 & -2 & -3 & -7 & -11 \\ \midrule
\rsp M(K) & \frac12 & \frac34 & \frac13 & \frac47 & \frac9{11} \\
\bottomrule
\end{array}$$
\end{center}

\begin{proof}
  If $m \equiv 2, 3 \pmod{4}$, then for every $x \in K$ there exists a $y
  \in R$ with $x-y = a+b\sqrt{m}$ and $|a|, |b| \le \frac12$. It follows that
  $N_{K/\mathbb{Q}}(x-y) = a^2 - m b^2 \le (1+|\,m\,|)/4$. To also see
  $M(K) \ge (1+|\,m\,|)/4$, we consider e.g.\ the point
  $x = (1+\sqrt{m})/2$; for all $y \in R$ with $x-y = a+b\sqrt{m}$ we then
  have $|a|, |b| \ge \frac12$, and it follows that
  $N_{K/\mathbb{Q}}(x-y) \ge (1+|\,m\,|)/4$.
 
  
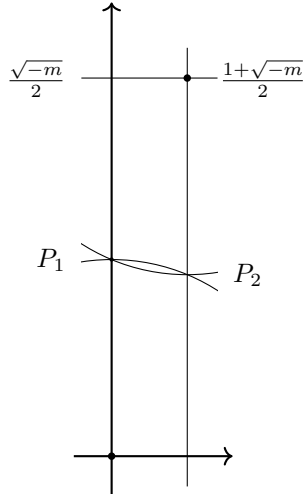
\begin{center}
    \begin{tikzpicture}[scale=2]
      \draw[->,thick] (-0.25,0) -- (0.8,0);
      \draw[->,thick] (0,-0.25) -- (0,3);
      \fill (0,0) circle (0.7pt);
      \fill (0.5,2.5) circle (0.7pt);
      \draw (0.5,-0.2) -- (0.5,2.7);
      \fill (0,1.3) circle (0.4pt);
      \node at (-0.4,1.3) {$P_1$};
      \node at ( 0.9,1.2) {$P_2$};
      \node at ( 1,2.5) {$\frac{1+\sqrt{-m}}2$};
      \node at (-0.5,2.5) {$\frac{\sqrt{-m}}2$};
      \draw (-0.2,2.5) -- (0.7,2.5);
      \begin{scope}
        \clip (-0.2,0) rectangle (0.7,3);
        \draw (0.5,2.5) circle (1.3);
        \draw (0,0) circle (1.3);
      \end{scope}
    \end{tikzpicture}
  \end{center}
   
  In the case $m \equiv 1 \pmod{4}$ we use a geometric method: we
  embed $K$ via $a+b\sqrt{m} \to (a,b\sqrt{|m|})$ into
  $\mathbb{R}^2$. If $\|\cdot\|$ is the ordinary Euclidean norm in
  $\mathbb{R}^2$, then $N_{K/\mathbb{Q}}$ and $\|\cdot\|^2$ coincide
  in the sense that $N_{K/\mathbb{Q}}(a+b\sqrt{m}) =
  \|(a,b\sqrt{|m|})\|^2$. To show that for every $x \in
  K$ there exists a $y \in R$ with $N_{K/\mathbb{Q}}(x-y) < k$, we may
  restrict ourselves to those $x \in K$ with $x=r+s\sqrt{m}$ and $|r|,
  |s| \le \frac12$, since we can reach all other $x$ by adding
  suitable $y \in R$.

  For reasons of symmetry we may even assume $0 \le r, s \le \frac12$.
  We now set $k = \frac{1+|m|}{4\sqrt{|m|}}$ and claim that the
  circles about the points $0$ and $\frac{1+\sqrt{m}}2$ (more precisely,
  about their images $(0,0)$ and $(\frac12, \frac{\sqrt{|m|}}2)$), of
  radius $k$, cover the whole of $F_\varphi$, where $F_\varphi =
  \{x=r+s\sqrt{m}: 0 \le r, s \le \frac12\}$.

  To this end we simply check that these two circles intersect in
  $F_\varphi$ precisely at the points
  $$ P_1 = \Big(0, \frac{(1+|m|)\sqrt{|m|}}{4|m|} \Big)
                 \quad \text{ and } \quad
     P_2 = \Big(\frac12, \frac{(m-1)\sqrt{|m|}}{4|m|} \Big). $$
  Denoting the above embedding by $\varphi$, there therefore exists,
  for every $x \in K$, a $y \in R$ with $\|\varphi(x-y)\|
  \le k$, i.e.\ with $N_{K/\mathbb{Q}}(x-y) \le k^2$ as claimed. The
  points $P_1$ and $P_2$ also show that this bound is best possible,
  which completes the proof.
\end{proof}

Concerning (0.21) we wish to add some further remarks. To this end let $f:R \to
\mathbb{N}$ be a function satisfying (E--1), $a,b \in R \setminus
\{0\}$, and $M(f;a,b) = \inf \{f(a-bq)/f(b): q \in R\}$; we then find
$M(f) = \sup \{M(f;a,b): a,b \in R \setminus \{0\}\}$. If $f$ is
moreover multiplicative, we may write more simply $M(f;x) = \inf
\{f(x-q): q \in R\}$, $x = a/b$. If $a,b \in R
\setminus \{0\}$ and a positive real $\varepsilon > 0$ are
given, then, provided $M(f) < \infty$, we can find a $q \in R$
with $f(a-bq)/f(b) < M(f) + \varepsilon$. On the other hand it is
conceivable that no $q \in R$ satisfies $f(a-qb)/f(b) \le M(f)$
(there might, after all, be a sequence $(q_n)_{n \in \mathbb{N}}$ with
$f(a-q_n b)/f(b) \to M(f)$ as $n \to \infty$). If for all $a,b
\in R \setminus \{0\}$ the inequality $f(a-qb)/f(b) \le M(f)$ can be
satisfied,
we say that $M(f)$ is attained. As the example of the minimal Euclidean
function shows, even when the Euclidean minimum $M(f)$ is attained,
there need not exist $a,b \in R$ with $M(f;a,b) = M(f)$. As we have
seen in (0.21), this is the case in imaginary quadratic number fields
if we let $f$ be the norm and pick $P = (1+\sqrt{m})/2$ for
$m \equiv 2, 3 \pmod{4}$, and $P_1$ and $P_2$ for $m \equiv 1 \pmod{4}$.

We now set\label{DefC}
$C_1 = \{(a,b) \in R \times R : M(f;a,b) = M(f)\}$ and define
\[
M_2(f) = \sup \{ M(f;a,b) : (a,b) \in R \times R \setminus C_1 \}.
\]
If $M_2(f) < M(f)$, then $M(f)$ is called isolated, and $M_2(f)$ the
\textbf{second Euclidean minimum}\label{pEuM2} of $R$ with respect to
$f$. As (0.21) shows, the Euclidean minimum with respect to the norm
is not isolated in imaginary quadratic number fields. Barnes and
Swinnerton-Dyer, on the other hand, have conjectured that it is
always isolated in number fields of unit rank $>1$.

We can continue in this way and ask whether the
second Euclidean minimum is also isolated, and then correspondingly
define $M_3(f)$, and so on. We shall see that there are number fields
possessing an infinite chain $M_1(K), M_2(K), \dots$ of Euclidean minima
(such as $\mathbb{Q}(\sqrt{5})$), while in other
number fields (such as $\mathbb{Q}(\sqrt{23})$) already $M_2(K)$ is not
isolated. Finally, we can introduce corresponding notions
for the $k$-stage minima $M_k(f)$ as well.

\section*{{\sc Remarks on} \S\ 0}

As already mentioned above, Gauss was the first to prove
that an algebraic number field other than $\mathbb{Q}$ (specifically
$\mathbb{Q}(i)$) is norm-Euclidean. It is not entirely correct,
however, that he did this in order to prove the unique prime
factorization theorem in $\Z[i]$. Instead, he proceeded as follows: in
Art.\ 33 (see Werke II, theoria residuorum biquadraticum) he uses
Fermat's insight that every prime $p \equiv 1 \pmod{4}$ can be written
as a sum of two squares, in order to show that such $p$ become the
product of two distinct prime factors in $\Z[i]$.  With the help of
this fact he reduces, in Art.\ 37, the unique prime factorization in
$\Z[i]$ to that in $\Z$.

Only in Art.\ 46 does he then describe the Euclidean algorithm for
determining two numbers from $\Z[i]$.  The proposal to search
number fields for Euclidean functions other than the absolute
value of the norm was probably first made by H.\ Hasse
\cite{Has28}. What such functions might look like has been described
by Lenstra \cite{Len74} and Bedocchi \cite{Bed85}; we shall return to
their proposals in §10. The notion of the Euclidean algorithm that we
defined at the very beginning of this section can be
generalized further: we can, for instance, define left-Euclidean
(respectively right-Euclidean) functions in non-commutative rings and
ask, for example, about Euclidean quaternion algebras.

This has been treated, e.g., by Dickson (\cite{Dic27}), R\'edei
\cite{Red67}, Newman \cite{New72}, Brungs \cite{Bru73} and Lenstra
\cite{Len74,Len78} (see also Felgner \cite{Fel73}).  Further, instead
of maps $f:R \to \N$ we can consider functions mapping $R$
into a well-ordered group $W$: Samuel \cite{Sam71} has shown that the
requirement $W = \mathbb{N}$ is no restriction in rings with finite
residue fields, in the sense that a ring Euclidean with
respect to some $g: R \to W$ is also Euclidean with respect to a
suitable $f: R \to \mathbb{N}$.

Hiblot \cite{Hib75}, on the other hand, has constructed a ring for which
the two definitions are not equivalent. For a thorough investigation
of such and similar questions we refer to Lenstra \cite{Len74}. A
connection between the $\mathrm{GE}_2$-rings, which were only briefly
touched upon in the text, and the Euclidean algorithm can be found in
Hurwitz \cite{Hur19}, P.\ Cohn \cite{Coh66}, Vaserstein \cite{Vas72}
and Cooke \cite{Coo77}.

\chapter*{\S\ 1 Criteria for the Existence of the Euclidean Algorithm}
\setcounter{chapter}{1}
\addcontentsline{toc}{chapter}
                {\S\ 1 Criteria for the Existence of the Euclidean Algorithm}
\markboth{Euclidean Rings}
         {\S\ 1 Criteria for the Existence of the Euclidean Algorithm}

In what follows let $K$ be an algebraic number field and $R$ the ring
of integers of $K$. In order to decide whether $R$ is norm-Euclidean
or not, various criteria have been developed, of which we wish to
describe some.  To this end let $P$ be a prime ideal in $R$; we call a
polynomial
\[
f(x) = x^n + a_{n-1}x^{n-1} + \dots + a_1 x + a_0 \in R[x]
\]
\textbf{Eisenstein}\index[S]{Eisenstein} with respect to $P$ if
$a_i \in P$ for $0 \le i \le n-1$ and $a_0 \in P \setminus P^2$
(note that $f$
must be monic!). The Eisenstein irreducibility criterion
then asserts that such polynomials are irreducible over $K(x)$.
We now wish to consider the following situation: let $L/K$ be a
finite field extension of degree $(L:K) = n$, $S$ the ring of
integers of $L$ and $Q$ a prime ideal in $S$ above $P$:
\begin{center}
  \begin{tikzpicture}[scale=0.6]
    \node (Q) at (0,0) {$\Q$};
    \draw (1,0) node {$\supset$};
    \node (Z) at (2,0) {$\Z$};
    \draw (3,0) node {$\supset$};
    \node (p) at (4,0) {$(p)$};
    \node (K) at (0,2) {$K$};
    \draw (1,2) node {$\supset$};
    \node (R) at (2,2) {$R$};
    \draw (3,2) node {$\supset$};
    \node (P) at (4,2) {$P$};
    \node (L) at (0,4) {$L$};
    \draw (1,4) node {$\supset$};
    \node (S) at (2,4) {$S$};
    \draw (3,4) node {$\supset$};
    \node (Q1) at (4,4) {$Q$};
    \draw (Q) -- (K) -- (L);
    \draw (Z) -- (R) -- (S);
    \draw (p) -- (P) -- (Q1);
  \end{tikzpicture}
\end{center}

The prime ideal $P$ is called \textbf{totally ramified}\index[S]{prime
  ideal!totally ramified} in $L/K$ if $PS = Q^n$. It is well known
that a prime ideal $P$ is totally ramified in $L/K$ if and only if the
extension is generated by a root of a polynomial that is Eisenstein
with respect to $P$; for the sake of completeness we give the proofs
here.

\begin{quote}
  {\bf (1.1)} {\em Let $f(x) = x^n + a_{n-1}x^{n-1} + \dots + a_0 \in
    R[x]$ be Eisenstein with respect to $P$ and $\alpha$ a root of
    $f$; if then $L=K(\alpha)$, one has $(L:K) = n = \deg f$ and $P$
    is totally ramified in $L/K$.}
\end{quote}

\begin{proof}
  By Eisenstein's irreducibility criterion $f$ is
  irreducible, consequently $(L:K) = n$. Moreover we can write $a_i S =
  P A_i$ for certain integral ideals $A_i$ in $S$, where still
  $A_0 \equiv 0 \bmod P$. Since $f(\alpha) = 0$ one has $\alpha^n =
  -(a_{n-1}\alpha^{n-1} + \dots + a_0) \in PS$. Now let $Q$ be any
  prime ideal in $S$ above $P$; then $\alpha^n \in Q$, hence even
  $\alpha \in Q$, because $Q$ is prime, and thus $\alpha^n \in
  Q^n$. Writing $\alpha^n S = (PS)B$ for an integral ideal $B$ in
  $S$, one obtains $Q^n \mid (PS)B$. If we can still show
  that $Q + B = S$, we have $Q^n \mid PS$ and thereby $Q^n = PS$.
  We therefore assume that $Q \mid B$; then however $a_0 \in
  Q(PS)$ because of
  $$ \alpha^n = -(a_{n-1}\alpha^{n-1} + \dots + a_1)\alpha - a_0 \in Q(PS). $$
  Since $a_0 \in R$ and $Q \cap R = P$, this
  implies $a_0 \in P^2$, in contradiction to the hypothesis.
\end{proof}

\begin{quote}
  {\bf (1.2)} {\em Let $P$ be totally ramified in $L/K$, $\pi \in S$ and
    $f$ the characteristic polynomial of $\pi$. If then $\pi \in Q$ and $f(x) =
    x^n + a_{n-1} x^{n-1} + \dots + a_0 \in R[x]$, one has $a_i \in
    P$ for all $0 \le i \le n-1$; if moreover $\pi \in Q \setminus
    Q^2$, then $f$ is even Eisenstein with respect to $P$ and therefore
    equal to the minimal polynomial of $\pi$.}
\end{quote}

\begin{proof}
  Since $f(\pi)=0$ one has $a_0 \equiv 0 \bmod {\pi}$, hence $a_0 \in Q \cap
  R = P$. Considering the congruence $f(\pi) \equiv 0 \bmod {Q^2}$,
  one obtains correspondingly $a_1 \equiv 0 \bmod {P}$ etc. If finally
  $\pi \in Q \setminus Q^2$ and if $a_0 \in P^2$, then
  $\pi^n \equiv 0 \bmod {\pi P}$ would follow, in contradiction to $PS =
  Q^n$. Thus $f$ is Eisenstein with respect to $P$.
\end{proof}
The significance of totally ramified prime ideals for this section rests on

\begin{quote}
  {\bf (1.3)} {\em If $P$ is totally ramified in $L/K$, then for
    every $\alpha \in S$ there exists an $a \in R$ with $\alpha \equiv a \bmod {Q}$
    and one has $N_{L/K}(\alpha) \equiv a^n \bmod {P}$. }
\end{quote}

\begin{proof}
  Let $\alpha \in S$; since $S/Q \cong R/P$ ($P$ is totally ramified,
  and in particular $Q$ has residue degree 1 over $P$) there exists an
  $a \in R$ with $\alpha \equiv a \bmod {Q}$. Thus $\alpha - a \in
  Q$, and the characteristic polynomial of $\alpha - a$ has by (1.2) the form
  $H_{P}(\alpha - a; K)(x) = x^n + a_{n-1} x^{n-1} + \dots + a_0$ with
  $a_i \in P$ for $0 \le i \le n-1$; hence $\alpha$ is a root of
  $(\alpha - a)^n + a_{n-1}(\alpha - a)^{n-1} + \dots + a_0$, and we
  see $N_{L/K}(\alpha) = a^n - a_{n-1}a^{n-1} + \dots + (-1)^n a_0
  \equiv a^n \bmod {P}$.
\end{proof}
All known criteria for the existence of a Euclidean algorithm that rest on totally
ramified prime ideals can be reduced to the following:

\begin{quote}
  {\bf (1.4)} {\em Let $B$ be a product of pairwise distinct prime ideals that are totally
    ramified in $L/K$, $BS=C^n$, and $C$ a
    principal ideal in $S$. If then there exists an $x \in R$ with
    $x^n \equiv e \bmod {B}$ for some $e \in R$, and if there exists no
    $r \in R$ with the properties (a), (b) and (c), then }$M(L) \ge k$.
    \begin{enumerate}
    \item[(a)] $r \equiv e \bmod {B}$ \\
    \item[(b)] $r = N_{L/K}(\alpha)$ {\em for some} $\alpha \in S$ \\
    \item[(c)] $|N_{K/\Q}(r)| < k \cdot \|B\|$.
\end{enumerate}
\end{quote}

\begin{proof}
  We assume that $M(L) < k$; since $C = (c)$ for some $c \in S$
  there then exists a $q \in S$ with $|N_{L/\Q}(x - qc)| < k \cdot
  |N_{L/\Q}(c)|$. Setting $\alpha = x - qc$, one has $\alpha
  \equiv x \bmod {C}$ and $|N_{L/\Q}(\alpha)| < k \cdot
  \|C\|$. With $r = N_{L/K}(\alpha)$ one now has
  (a) $r = N_{L/K}(\alpha) \equiv x^n \equiv e \bmod {B}$; by (1.3)
  this congruence holds for every prime divisor $P$ of $B$.
  Since these $P$ are pairwise distinct by hypothesis,
  the claim follows from the Chinese remainder theorem.
  (b) $r = N_{L/K}(\alpha)$ holds by construction.
  (c) $|N_{K/\Q}(r)| = |N_{K/\Q}(N_{L/K}(\alpha))| =
  |N_{L/\Q}(\alpha)| < k \|C\| = k \|B\|$, for for prime ideals
  $Q$ in $S$ above $P$ of relative residue degree 1 one has $\|Q\|=\|P\|$,
  so that $\|C\|=\|B\|$ follows from the multiplicativity of the ideal norm.
  Since by hypothesis there exists no $r \in R$ with the
  properties (a), (b), (c), one must have $M(L) \ge k$.
\end{proof}
We wish to remark explicitly that $L$ is also not
norm-Euclidean if in (1.4) we can only choose $k=1$: then
$M(N_{L/\Q}; x, c) = 1$, i.e.\ there exists no $q \in
R$ with $|N_{L/\Q}(x-qc)| < |N_{L/\Q}(c)|$.
Hitherto (1.4) was known only in the case $K=\Q$, $k=1$; in order to be able to apply this
criterion, $L/\Q$ therefore had to contain totally
ramified prime ideals. The version given here can
on the other hand also be applied when there exists an intermediate field $K$ with
$\Q \subset K \subset L$ such that $L/K$ possesses totally ramified
ideals (a typical example of this are e.g.\ the
Dirichlet number fields, which are quadratic extensions of
$\Q(i)$, $i^2=-1$).
The (hitherto necessary) condition that $L/\Q$ must contain totally ramified
prime ideals has ultimately led to the fact that much more is known about
the Euclidean algorithm in pure number fields
(Cioffari\index[N]{Cioffari} \cite{Cio79}, Egami\index[N]{Egami}
\cite{Ega79,Ega84}) or in cyclic number fields
(Heilbronn\index[N]{Heilbronn} \cite{Hei50,Hei51}, Egami \cite{Ega84})
than in fields that contain no totally
ramified prime ideals. In the special case $K=\Q$ one obtains from (1.4)

\begin{quote}
  {\bf (1.5)} {\em Let $(L:\Q) = n$, $f \in \N$ a
    product of pairwise distinct rational primes that are totally ramified in $L/\Q$
    and $fS=C^n$ for a
    principal ideal $C$ in $S$. If then $a$ is an $n$-th power residue mod $f$
    and all $r \in \Z$ with $a \equiv r \bmod {f}$ and
    $|r| < k \cdot f$ are not norms from $S$, then $M(L) >
    k$. Moreover $L$ is not norm-Euclidean if $k \ge 1$.}
\end{quote}

\begin{proof}
  (1.4) with $K=\Q$, $B=\{f\}$, $e=a$ and $R=\Z$.
\end{proof}
For $k=1$ one finally obtains from (1.5) the previously already known
criterion (Egami\index[N]{Egami} 1979):

\begin{quote}
  {\bf (1.6)} {\em Let $(L:\Q)=n$, $f \in \N$ a product
    of pairwise distinct rational primes that are totally ramified;
    if then there exist $a, b \in \N$ with $f=ab$ such that
    $a$ is an $n$-th power residue mod $f$ and neither $a$ nor $-b$
    are norms from $S$, then $L$ is not norm-Euclidean. }
\end{quote}

\begin{proof}
  (1.5) with $k=1$; note that $a$ and $-b$ are the only $r \in
  \Z$ with $r \equiv a \bmod {f}$ and $|r| < f$.
\end{proof}
We demonstrate the application of (1.5) by some simple
examples:
\begin{enumerate}
\item[1.] $L = \Q(\sqrt{7})$, $f=14$, $a=9$; then $a$ is a
  quadratic residue mod 14, and the $r \in \Z$ with $r \equiv a \bmod
  {f}$ are $r = \dots, -19, -5, 9, 23, \dots$ etc. If we now choose
  $k = \frac{9}{14}$, then only $r=-5$ satisfies the condition
  $|r| < k \cdot f = 9$;
  now however there exists no $\alpha \in S$ with $N_{L/\Q}(\alpha) = -5$,
  for the prime ideal $(5)$ remains inert in the passage from $\Z$ to $S$
  because $(\frac{28}{5}) = -1$ (Legendre symbol). Thus
  $M(L) \ge \frac{9}{14}$. In order also to find a $z \in L$ with
  $M(z) = \frac{9}{14}$ (here we have written simply $M(z)$ instead of $M(N_{L/\Q}, z)$), we must
  first find a principal ideal $C$ of norm $f=14$: such an ideal is
  e.g.\ $C = (7-3\sqrt{7})$. Then we still need an $x \in \Z$ with
  $x^2 \equiv a \equiv 9 \bmod {14}$, e.g.\ $x=3$, and we have $z-x/c
  = -(21+9\sqrt{7})/14$; since $M(z) = M(z-y)$ for every $y \in R$
  we need only specify $z \bmod R$, i.e.\ we now have
  $z \equiv (7+5\sqrt{7})/14 \bmod {R}$.
\item[2.] $L = \Q(\sqrt{10})$: in order to show $M(L) \ge \frac32$ we use
  $f=10$ and $a=5$; indeed $a$ is a quadratic residue mod $f$ because $5^2 \equiv 5 \bmod
  {10}$. The only $r \in \Z$ with $|r| <
  k \cdot f$ and $r \equiv 5 \bmod {10}$ are $r = \pm 5$; since however
  the prime ideal above $(5)$ is not a principal ideal in $S$, there exist in $S$
  no elements of norm 5, and we have $M(L) \ge \frac32$.
\item[3.] $L = \Q(\sqrt{15})$: we claim $M(L) \ge \frac75$ and
  use for this $f=15$, $a=6$. Thereby $r \in \{-9, 6\}$. Since now
  the prime ideal above $(3)$ is not a principal ideal, all elements of norm
  9 are of the form $3u^k$, where $u=4+\sqrt{15}$ is the fundamental unit of $L$. Because
  $N_{L/\Q}(u) = +1$ one has however $N_{L/\Q}(3u^k) = +9$, i.e.\ there exists in
  $S$ no $\alpha$ with $N_{L/\Q}(\alpha) = -9$. Analogously every
  ideal of norm 6 is generated by an element of the form $\beta = (3+\sqrt{15})u$,
  and because $N_{L/\Q}(\beta) = -6$ there exist in $S$ no
  elements of norm $+6$. This was to be shown.
\end{enumerate}
The following criterion contains (1.4) as a special case (for $m = 0,1$)
and could be generalized still further; we shall however use it only
in the simpler version presented:

\begin{quote}
  {\bf (1.7)} {\em Let $P$ be a prime ideal in $R$ with $PS = Q^m$, $B$
    a product of pairwise distinct prime ideals with $B \cap P
    = R$, and let $B' = P^m B$, as well as $B'S = Q^{mn} C^n$ for a
    principal ideal $Q^m C$. If then $x^n \equiv e \bmod {P}$ for $x, e \in R$,
    and if there exists a $y \in S$ with $y \equiv x \bmod {QC}$ such that
    there exists no $r \in R$ with the properties
    \begin{enumerate}
    \item[(a)] $r \equiv e \bmod {PB}$
    \item[(b)] $r = N_{L/K}(\alpha)$ for
      some $\alpha \in S$ with $\alpha \equiv y \bmod {Q^m C}$
    \item[(c)] $|N_{K/\Q}(r)| < k \cdot \|B\|$
    \end{enumerate}
    then $M(L) \ge k$.}
\end{quote}

\begin{proof}
  Let $M(L) < k$; then there exists an $\alpha \in S$ with $\alpha \equiv y
  \bmod {Q^m C}$ and $|N_{L/\Q}(\alpha)| < k \cdot \|Q^m C\|$,
  as well as an $x \in R$ with $y \equiv x \bmod {QC}$ (for $QC$ is a
  product of pairwise distinct totally ramified
  prime ideals). With $r = N_{L/K}(\alpha)$ the properties (a), (b), (c) can now be verified as in (1.4).
\end{proof}

The difficulty in the formulation of (1.7) rests on the following
fact: from $\alpha \equiv a \bmod {Q}$ it does follow that
$N_{L/K}(\alpha) \equiv a^n \bmod {P}$, but from $\alpha \equiv a
\bmod {Q^m}$ it does not in general follow that $N_{L/K}(\alpha)
\equiv a^n \bmod {P^m}$, but again only the weaker congruence mod~$P$.

We give some examples:
\begin{enumerate}
\item[1.] $L=\Q(\sqrt{14})$, $K=\Q$, $B'=(4)=(2)^2$, so $P=(2)$,
  $m=2$, $B=(1)$; further we have $Q = \frt_1=(4+\sqrt{14})$. We
  now set $k=5/4$ and $y=1+\sqrt{14}$. Thereby $y \equiv 1 \bmod
  {Q}$, hence $x \equiv e = 1$, and we must investigate whether there exists an
  $r \in \Z$ with $|r|<5$ and $r \equiv 1 \bmod {2}$ that is the norm
  of an $\alpha \in S$ with $\alpha \equiv 1+\sqrt{14} \bmod {2}$.
  Since $(3)$ is inert in $L/K$, there exist no elements of norm 3,
  so only $r=1$ comes into question. Now however there exists no unit
  $\alpha \in S$ with $\alpha \equiv 1+\sqrt{14} \bmod {2}$, because the
  fundamental unit $15+4\sqrt{14}$ is already $\equiv 1 \bmod {2}$ (thus every
  unit is $\equiv 1 \bmod {2}$). One may also compare this example
  with the proof on p.\ 11 that $\frac{1+\sqrt{14}}2$ is not a
  continued fraction of length 2 (in particular one should note the role
  that the fundamental unit plays in both proofs).
\item[2.] $L=\Q(\sqrt{26})$, $K=\Q$, $B'=(4)$, $y=26$, $k = \frac52$;
  here $Q =\frt_1 = (2,\sqrt{26})$ is not a principal ideal, but
  $Q^2=(2)$ is. By (1.7) we must investigate the $r \in \Z$ with $r \equiv 0
  \bmod {2}$ and $|r|<10$. Because $y \in Q \setminus Q^2$
  one has however for all $\alpha \equiv y \bmod {Q}$ also $\alpha \in Q
  \setminus Q^2$, and one easily sees that this
  implies $N_{L/\Q}(\alpha) \equiv 2 \bmod {4}$. Thus we need
  only consider the $r \equiv 2 \bmod {4}$. Since however there exist no
  $\alpha \in S$ of norm 2 (for $\frt_1$ is not a principal ideal) or
  of norm 6 (since the ideal $(3)$ remains inert in $L$), it follows in
  fact that $M(L) \ge \frac52$.
\end{enumerate}
The question remains how one should proceed when the criteria (1.4)
and (1.7) yield nothing (e.g.\ when the bounds they furnish
are not good enough or there exist no totally ramified
prime ideals at all). To this end let $K$ be a number field with $M(K)=k$ and $I=(b)$
an integral ideal in $R$ (hence $b \in R$). Then for all $a \in
R$ there exists a $q \in R$ with $|N_{K/\Q}(a-qb)| < k \cdot \|b\|$; because
$a-qb \equiv a \bmod {I}$ we can also express this as follows: every
residue class $\bmod {I}$ contains an $r \in R$ with $|N_{K/\Q}(r)|
< k \cdot \|b\|$. This consideration leads us to the definition
\begin{align*}
  M(K,I) & = \inf \left\{ \kappa \in \R : \forall a \in R \ \exists c
    \in I \text{ with } |N_{K/\Q}(a-c)| < \kappa \cdot \|I\| \right\} \\
         & = \inf \left\{ \kappa \in \R : \forall x \in I^{-1} \ \exists
            y \in R \text{ with } |N_{K/\Q}(x-y)|< \kappa \right\}.
\end{align*}
Since for every $x \in K$ there exists a principal ideal $I$ with $x \in I^{-1}$,
it follows immediately from the definition that
$M(K) = \sup \{ M(K, I) : I \text{ is a principal ideal in } R \}$.
A principal ideal $I$ with $M(K,I)<1$ we call a
\textbf{Euclidean ideal},\index[S]{ideal!Euclidean} one with
$M(K,I)>1$ a \textbf{non-Euclidean
  ideal}.\index[S]{ideal!non-Euclidean}

In (1.4) and (1.7) we have estimated $M(K,I)$ from below for
totally ramified principal ideals $I$. We may expect better
bounds for $M(K)$ when we allow arbitrary principal ideals for $I$,
because $M(K)>M(K,I)$ for every principal ideal $I$.

Consider for example $K = \Q(\sqrt{29})$; here
$\alpha = \frac{1}{2}(3+\sqrt{29})$ generates an ideal $I$ of norm $5$.
In order to determine $M(K,I)$ we first set up a complete
residue system mod $I$; such a system is e.g.\ $\{0, \pm 1, \pm
2\}$. Now we seek in every residue class mod $I$ an element
of minimal norm $> 1$; in the residue classes $0, \pm 1$ mod $I$ these are
the elements $0, \pm 1$ themselves (these are so to speak the
``uninteresting'' residue classes). In order to find in the residue class $2$ mod $I$ an
element of minimal norm, we must ask ourselves whether this
residue class contains a unit. Now however for the fundamental unit $u =
\frac{1}{2}(5+\sqrt{29})$ one has the congruence $u \equiv 1 \bmod {I}$; since by
Dirichlet every unit $e$ can be written in the form $e = u^l$ for some $l \in
\Z$, one therefore has $e \equiv 1 \bmod {I}$ for every
unit $e \in R^\times$, and in particular there exist no units in the
residue classes $2$ mod $I$.
Since there also exist in $R$ no elements of norm $2$ or $3$ (because the
ideals $(2)$ and $(3)$ remain prime in $R$), the residue class $2$ mod
$I$ can contain only elements of norm $\ge 4$, and in fact $r=2$
is such an element. Thereby we have shown $M(K,I) = 4/5$ (this is
somewhat better than the bound $M(K) \ge 23/29 =
0.7931\dots$ obtained with (1.5)). One should further note that in this case $I = (u-1)$,
and that this fact has ensured that $e \equiv 1 \bmod {I}$
for every unit $e \in R^\times$.
The proof of $M(K,I) = 16/9$ for $K =
\Q(\sqrt{85})$ and $I = (u-1) = \frac{1}{2}(7+\sqrt{85})$ runs quite analogously; $I$
has ideal norm 9, and $\{0, \pm 1, \pm 2, \pm 3, \pm 4\}$ is a
complete residue system mod $I$.
It is thus in principle clear how one carries out corresponding computations in
arbitrary number fields for arbitrary principal ideals $I$: one
chooses a system of fundamental units and determines all residue classes
mod $I$ that contain units (for this one simply lets every
fundamental unit run through the powers from $1$ to $\Phi(I)$ and
notes the residue classes in which the products of these powers
lie). Then one seeks a principal ideal $A = (a)$ of minimal norm
and asks in which residue classes mod $I$ the elements $ae$, $e \in
R^\times$, lie etc.
However the necessary computations become rather extensive rather quickly
if one has to consider an $I$ of large norm (e.g.\ for $K =
\Q(\sqrt{31})$ and $I = (u-1)$, $u = 1520 + 273\sqrt{31}$,
here $\|I\| = 3038$) or if the unit rank of $K$ becomes large.
Fortunately the above procedure can be modified in such a way
that it can be programmed rather simply and thus carried out by a
computer. To this end we say that the points $x_1,
\dots, x_t \in K$ are (cyclically) permuted by a unit $u \in R^\times$
if the following congruences hold: $x_1u \equiv x_2
\bmod {R}$, $x_2u \equiv x_3 \bmod {R}$, $\dots$, $x_tu \equiv x_1 \bmod {R}$. In
this case one has incidentally $M(K, x_1) = \dots = M(K, x_t)$, as one
easily sees from the two observations
\begin{enumerate}
\item[1.] for all $y \in R$ one has $M(K, x) = M(K, x-y)$;
\item[2.] for all units $u \in R^\times$ one has $M(K, x) = M(K, ux)$
\end{enumerate}
We now claim:

\begin{quote}
  {\bf (1.8)} {\em Let $R=D(m)$, $m\in\N$ square-free, and let
    $1\le u$ for a unit $u$ in $R$. Further let points $x_1,
    \dots, x_t$ be given that are permuted by $u$. If then
    $M(K,x_j)<k$, there exists a $z=r+s\sqrt{m} \in K$ with the
    properties
    \begin{enumerate}
    \item[(a)] $z \equiv x_j \bmod {R}$ for some $j \in \{1, \dots, t\}$;
    \item[(b)] $|N_{K/\Q}(z)| < k$;
    \item[(c)] $|r| < \mu_1$, $|s| < \mu_2$ with
      $$ \mu_1 = \frac{\sqrt{k}}2 \Big(\sqrt{u} + \frac1{\sqrt{u}} \Big)
      \quad \text{and} \quad
     \mu_2 = \frac{\sqrt{k}}{2\sqrt{m}}\Big(\sqrt{u} + \frac1{\sqrt{u}}\Big). $$
    \end{enumerate} }
\end{quote}

This is essentially ``Theorem B'' of Barnes and Swinnerton-Dyer;
their proof is however completely tailored to quadratic number fields
and cannot be generalized. Before we give
a proof of (1.8), however, we wish to make clear what (1.8) asserts in the
simplest case $t=1$: We then have an $x=x_1$ given with
$xu \equiv x \bmod {R}$; if we then wish to show $M(K,x) \ge k$, we assume
that $M(K,x)<k$. (1.8) then guarantees the existence of a
$z=r+s\sqrt{m} \in K$ with the properties
\begin{enumerate}
\item[(a)] $z \equiv x_j \bmod {R}$ for some $j \in \{1, \dots, t\}$;
\item[(b)] $|N_{K/\Q}(z)| < k$;
\item[(c)] $|r| < \mu_1$, $|s| < \mu_2$
\end{enumerate}
In order to lead $M(K,x)<k$ to a contradiction, we therefore need only show
that among the finitely many $z \in K$ that satisfy the conditions (a) and (c)
none occurs with $|N_{K/\Q}(z)|<k$.
\begin{proof}[Proof of (1.8)]
  Since $M(K,x_1)<k$ there exists a $y \in R$ with
  $|N_{K/\Q}(x_1-y)| < k$; we now set $z_1=x_1-y$ and
  choose $m \in \Z$ so that $\sqrt{k/u} \le |z_1 u^m| <
  \sqrt{k u}$ (this is evidently always possible). We now claim
  that $z=z_1 u^m$ satisfies the conditions (a), (b), (c):
  \begin{enumerate}
  \item[(a)] $z = z_1 u^m = x_1 u^m - y u^m \equiv
    x_1 u^m \equiv x_j \bmod {R}$,
    where $j$ is determined by the congruence $j \equiv 1+m \bmod {t}$.
  \item[(b)]
    $|N_{K/\Q}(z)| = |N_{K/\Q}(z_1)| = |N_{K/\Q}(x_1-y)| < k$.
  \item[(c)] Let $\tilde{z} = r-s\sqrt{m}$ be the conjugate of $z$;
    then we have
    $$ |z'| = |z z'| / |z| = |N_{K/\Q}(z)| / |z| < k/|z| \le \sqrt{k u}, $$
    hence $|z + z'| \le |z| + |z'| < 2\sqrt{k u}$ and
    $|2s\sqrt{m}| = |z - z'| \le |z| + |z'| < 2\sqrt{k u}$.
  \end{enumerate}
\end{proof}
These bounds can however still be improved if, besides the
inequalities $|z|$, $|z'| < \sqrt{k u}$, we also use the inequality $|z z'|
< k$; the claim then follows with $a=\sqrt{k u}$ and
$b=k$ from

\begin{quote}
  {\bf (1.9)} {\em Let $x, y, a, b$ be positive real numbers; from the
    inequalities $x \le a$, $y \le a$, $xy \le b$ it then follows that $x+y \le
    a + \frac{b}{a}$. }
\end{quote}

\begin{proof}
  We have $0 \le (a-x)(a-y) = a^2 - (x+y)a + xy \le a^2 - (x+y)a + b$.
\end{proof}

We remark that the hypothesis $u > 1$ may be replaced by $u \ne 1$;
then we only have to replace $u$ by $|u|$ in (1.8.c). In the case
$N_{K/\Q}(u)=+1$ we moreover recover the bounds given by Barnes
and Swinnerton-Dyer: setting
$u=a+b\sqrt{m}$, one has $\frac{1}{u} = a-b\sqrt{m}$
and thus $(\sqrt{u} + 1/\sqrt{u})^2 = u + 1/u + 2 = 2a+2$, and we find
$\mu_2 = \sqrt{\frac{k(a+1)}{2m}}. $

Example: $K = \Q(\sqrt{19})$; let $x = \frac{20}{57}\sqrt{19}$ and
$u=170+39\sqrt{19}$ the fundamental unit of $K$. We then have $xu =
260+60\sqrt{19} - x \equiv -x \bmod {R}$, so that we would in
principle have to apply (1.8) with $t=2$, $x_1=x$ and $x_2=-x$. If
however we take instead of $u$ the unit $u'=-u$, we have $xu' \equiv x
\bmod {R}$, so that the choice of $u'$ halves the computation time.

We now show that $M(K) \ge M(K,x) = \frac{170}{171}$; the inequality
$M(K,x) \le \frac{170}{171}$ follows from the observation
\[
M(K,x) \le |N_{K/\Q}(x -3-\sqrt{19})| = \frac{170}{171}.
\]
Thus we only need to show $M(K,x) \ge k = \frac{170}{171}$, and
we do this with (1.8): for if $M(K,x)<k$, there would exist a $z
\in K$ with $z \equiv x \bmod {R}$ that satisfies the conditions (a), (b), (c).
Because of $a=170$ we obtain the bound $|s| < 2.12\dots$ for
$z=r+s\sqrt{19}$. Thus there are only four possible values for $s$,
namely $s = i + 20/57$, $i \in \{-2, -1, 0, 1\}$:
$$ \begin{array}{rrrl}
  i & = & -2: & |N_{K/\Q}(z)| \text{ is minimal for $r = \pm 7$ with }
              |N_{K/\Q}(z)| = 457/171; \\
  i & = & -1: & |N_{K/\Q}(z)| \text{ is minimal for $r = \pm 3$ with }
              |N_{K/\Q}(z)| = 170/171; \\
  i & = & 0: & |N_{K/\Q}(z)| \text{ is minimal for $r = \pm 1$ with }
              |N_{K/\Q}(z)| = 229/171; \\
  i & = & +1: & |N_{K/\Q}(z)| \text{ is minimal for $r = \pm 6$ with }
              |N_{K/\Q}(z)| = 227/171.
\end{array} $$
By (1.8) one therefore has $M(K,x) = \frac{170}{171}$, and in fact we shall
see in \S\  2 that even $M(K) = M(K,x) = \frac{170}{171}$
holds. Had we only used the bound $|s| < \sqrt{k u}$, we
would have had to consider eight values for $s$.
Before we generalize (1.8) to arbitrary number fields of unit rank $>1$,
we wish to clarify by means of the cubic case what
difficulties arise. We therefore first show

\begin{quote}
  {\bf (1.10)} {\em Let $K$ be a cubic number field of unit rank 1,
    $u$ a unit with $1<|u|$, and let $x_1, \dots, x_t \in K$ be
    points that are permuted by $u$. Further let $\{\alpha_1,
    \alpha_2, \alpha_3\}$ be a $\Q$-basis of $K$. If then
    $M(K,x_i) < k$, there exists a $z \in K$, $z = r_1\alpha_1 +
    r_2\alpha_2 + r_3\alpha_3$, with the properties
    \begin{enumerate}
    \item[(a)] $z \equiv x_j \bmod {R}$ for some $j \in \{1,\dots,t\}$;
    \item[(b)] $|N_{K/\Q}(z)| < k$;
    \item[(c)] $|r_i| < \mu_i$ for $i=1, 2, 3$, where the $\mu_i$
      are positive real numbers independent of the $x_i$.
    \end{enumerate} }
\end{quote}

\medskip\noindent \textbf{Remark.} The bounds $\mu_i$ will emerge from
the proof and can be determined explicitly. The bounds that we shall
obtain for pure cubic fields are far better than those obtained by
Cioffari\index[N]{Cioffari} \cite{Cio79} in the cases
$K=\Q(\sqrt[3]{m})$, $m=12, 17, 44$. As in the quadratic case the
bounds achieved here can be improved once more; we shall do this in
\S\ 4.

\begin{proof}
  Without loss of generality we may assume that $K$ is real (for since $K$
  has unit rank 1, one must have $r=s=1$, where $r$ denotes the number of real
  embeddings of $K$ into $\C$). Thus we may replace the
  condition $1<|u|$ by $1<u$.
  Now by hypothesis $M(K,x_1)<k$, consequently there exists a
  $y\in R$ with $|N_{K/\Q}(x_1-y)|<k$; we set $z_1=x_1-y$ and
  choose an $m\in\N$ so that
  $$ \sqrt[3]{k/u^2} \le |z_1u^{m-1}| < \sqrt[3]{ku} $$
  holds. With $z=z_1u^m$ the conditions (a) and (b) are then satisfied,
  so that we only still need to take care of (c). We now denote
  the two conjugates of $z$ by $z'$ and $z''$; since the two
  fields $K'$ and $K''$ are complex conjugate, one has $|z'|=|z''|$,
  hence $|N_{K/\Q}(z)| = |z z' z''| = |z| |z'|^2$ and thus
  $|z'|^2 = |N_{K/\Q}(z)|/|z| < k/|z| \le \sqrt[3]{ku^2}$.

  We first consider the somewhat simpler special case of the
  pure cubic number field $K=\Q(\vartheta)$ for $\vartheta^3=m$,
  $m\in\N$; here we choose the basis $\{1, \vartheta, \vartheta^2\}$
  and find $z' = r_1 + r_2 \rho \vartheta + r_3 \rho^2 \vartheta^2$,
  $z'' = r_1 + r_2 \rho^2 \vartheta + r_3 \rho \vartheta^2$, where $\rho$
  is a primitive third root of unity. We can now conclude quite as
  in the quadratic case:
  \begin{align*}
    |3r_1| & = |T_{K/\Q}(z')| = |z+z'+z''| \le |z|+|z'|+|z''| < 3\sqrt[3]{ku}, \\
    |3r_2\vartheta| & = |z+\rho^2 z' + \rho z''| \le |z|+|z'|+|z''|
                        < 3\sqrt[3]{ku}, \\
    |3r_3\vartheta^2| & = |z+\rho z' + \rho^2 z''| \le |z|+|z'|+|z''|
                        < 3\sqrt[3]{ku},
  \end{align*}
and thus we have proved (1.8) in this special case with the
bounds $\mu_1 = \sqrt[3]{ku}$, $\mu_2 = \sqrt[3]{ku/m}$,
$\mu_3 = \sqrt[3]{ku/m^2}$.

Since in the general case we cannot give the conjugates of $z$
explicitly, we must here take a somewhat different route.

To this end let $\{\beta_1, \beta_2, \beta_3\}$ be the dual basis
to $\{\alpha_1, \alpha_2, \alpha_3\}$, which is well known to be
defined by $T_{K/\Q}(\alpha_i \beta_j) = \delta_{ij}$ (Kronecker delta).
Thus we have
\begin{align*}
  T_{K/\Q}(z \beta_j) & = \sum_{i=1}^{3} r_i
  T_{K/\Q}(\alpha_i \beta_j) = r_j, \text{ hence} \\
  |r_j| & = |T_{K/\Q}(z \beta_j)| = |z\beta_j + z'\beta_j + z''\beta_j|
          \le |z\beta_j| + |z'\beta_j| + |z''\beta_j| \\
        & < \sqrt[3]{ku} (|\beta_j| + 2|\beta_j'|).
\end{align*}

We still need to indicate how we can determine the $\beta_j$ from the
$\alpha_i$. To this end we consider the $n \times
n$-matrix $A = (T_{K/\Q}(\alpha_i \alpha_j))$; we know that then
$\det A = \disc_{K/\Q}(\alpha_1, \alpha_2, \dots, \alpha_n) \ne 0$.
Let $B = (b_{ij})$ be the inverse matrix of $A$; then $\beta_i =
b_{i1}\alpha_1 + \dots + b_{in}\alpha_n$ (see e.g.\ Cohn\index[N]{Cohn}
\cite{Coh78} or Marcus\index[N]{Marcus} \cite{Mar77}).

Since we shall deal extensively with Dirichlet number fields in \S\  6,
we will give explicit bounds in this case, too.
To this end let $L = \Q(\sqrt{m}, \sqrt{n})$ be an
imaginary bicyclic biquadratic number field (i.e.\ $m, n \in
\Z, m < 0$). We choose the $\Q$-basis $\alpha_1 = 1,
\alpha_2 = \sqrt{m}, \alpha_3 = \sqrt{n}, \alpha_4 = \sqrt{mn}$ and
assume that $N_{L / \Q}(x_1 - y) < k$ for some $x_1 =
r_1\alpha_1 + r_2\alpha_2 + r_3\alpha_3 + r_4\alpha_4 \in L$, some $y
\in R$ and some $k \in \R$. Since the unit rank is 1,
there exists a fundamental unit $u$ which we can choose so
that $|u| > 1$. With $z_1 = x_1 - y$ we then find an $m \in \Z$ with
$$ \frac{\sqrt[4]{k}}{\sqrt{|u|}} \le |z_1| < \sqrt[4]{k} \cdot \sqrt{|u|}
   = \delta. $$

Now let $\mathrm{Gal}(L/\Q) = \{1, \sigma, \tau, \sigma\tau\}$,
and without loss of generality let $\sigma$ be the automorphism of complex
conjugation (thereby one then has $|z| = |z^\sigma|$; if e.g.\ $m$ and $n$ are
both negative, one can also characterize $\sigma$ by $\sigma: \sqrt{m}
\mapsto -\sqrt{m}, \sqrt{n} \mapsto -\sqrt{n}$). Now
one has $N_{L / \Q}(x) = x \cdot x^\sigma \cdot x^\tau \cdot
x^{\sigma\tau} = |x|^2 |x^\tau|^2$, hence also $|x^\tau|^2 \le N_{L /
  \Q}(x) / |x|^2 < \delta$.

We now set $\tau: \sqrt{m} \mapsto -\sqrt{m}, \sqrt{n} \mapsto
-\sqrt{n}$ and obtain
\begin{align*}
  4 \cdot |r_1| & = |z + z^\sigma + z^\tau + z^{\sigma\tau}| \le |z| +
                    |z^\sigma| + |z^\tau| + |z^{\sigma\tau}| \\
           & = 2 \cdot (|z| + |z^\tau|) < 4\delta, \quad \text{as well as} \\
  4 \cdot |r_2| \cdot \sqrt{|m|} & = |z - z^\sigma + z^\tau - z^{\sigma\tau}|
                  \le 2 \cdot (|z| + |z^\tau|) < 4\delta, \\
  4 \cdot |r_3| \cdot \sqrt{|n|} & = |z - z^\sigma - z^\tau + z^{\sigma\tau}|
                   \le 2 \cdot (|z| + |z^\tau|) < 4\delta, \\
  4 \cdot |r_4| \cdot \sqrt{|mn|}
                 & = |z + z^\sigma - z^\tau -z^{\sigma\tau}|
                   \le 2 \cdot (|z| + |z^\tau|) < 4\delta.
\end{align*}
\end{proof}

If we also use (1.9), then as in the real quadratic case we can 
replace the factor $\sqrt{|u|}$ in $\delta$ by
$(\sqrt{|u|} + 1/\sqrt{|u|})/2$, and we have

\begin{quote}
  {\bf (1.11)} {\em Let $L = \Q(\sqrt{m}, \sqrt{n})$ be an imaginary
    bicyclic number field, $u$ a unit in $L$ with $|u| > 1$, and let
    $x_1, \dots, x_t \in L$ be points that are permuted by $u$.
    
    If then $M(L, x_i) < k$, there exists a
    $z = r_1 + r_2\sqrt{m} + r_3\sqrt{n} + r_4\sqrt{mn} \in L$ with the
    properties
    \begin{enumerate}
    \item[(a)] $z \equiv x_i \bmod {R}$ for some
      $i \in \{1, \dots, t\}$;
    \item[(b)] $N_{L / \Q}(z) < k$;
    \item[(c)] $|r_i| < \mu_i$ for $1 \le i \le 4$ and with
      $$ \mu_1 = \lambda, \quad
         \mu_2 = \frac{\lambda}{\sqrt{|m|}}, \quad
         \mu_3 = \frac{\lambda}{\sqrt{|n|}}, \quad
         \mu_4 = \frac{\lambda}{\sqrt{|mn|}}, $$
      where
      $$ \lambda = \frac{\sqrt[4]{k}}2
             \Big(\sqrt{|u|} + \frac1{\sqrt{|u|}} \Big) $$
      has been set.
    \end{enumerate}}
\end{quote}

It is a curious fact that one can also employ (1.11)
in order to estimate $M(K, x)$ for a real quadratic field: for if $e >
1$ is the fundamental unit of $K$ and $N_{K/\Q}(e) = -1$, then there exists an
imaginary bicyclic number field $L = K(\sqrt{m})$, $m < 0$, in which
$-e = u^2$ becomes a square for a unit $u$ in $S$ (the ring of
integers of $L$; see on this \S\  6 and the literature indicated there).
Instead of then investigating $N_{K/\Q}(x-y) < k$ for some $x \in K$
with $xe \equiv x \bmod {R}$ and $y \in R$,
one considers $N_{L / \Q}(x-y) < k^2$ and
$N_{L/\Q}(xu-y) < k^2$. We shall give an example in \S\  2
that shows that with this idea one can save a lot of work.
We thereby come to the general case:

\begin{quote}
  {\bf (1.12)} {\em Let $K$ be a number field with
    $(K:\Q)=n=r+2s$, $\{u_1,\dots,u_l\}$ a system of
    $l=r+s-1$ independent units, $\{\alpha_1,\dots,\alpha_n\}$
    a $\Q$-basis of $K$ and $x_1,\dots,x_t$ elements of
    $K$ that are somehow permuted by the $u_i$. If then
    $M(K,x_i)<k$, there exists a $z\in K$ with $z = \sum r_i
    \alpha_i$ and the properties
    \begin{enumerate}
    \item[(a)] $z \equiv x_i \bmod {R}$ for some
      $i \in \{1, \dots, t\}$;
    \item[(b)] $|N_{K/\Q}(z)| < k$;
    \item[(c)] $|r_i| < \mu_i$ for $i=1,\dots,n$, where the $\mu_i$
      are positive real numbers that depend only on $k$, the choice of the
      $\alpha_i$ and the $u_i$, but not on the points $x_i$.
    \end{enumerate} }
\end{quote}

Before we prove this, we introduce some further notation.
To this end we think of the $n=r+2s$ embeddings of $K$ into
$\C$ as arranged as follows:
$$ \tau_1, \dots, \tau_r, \tau_{r+1}, \ov{\tau}_{r+1}, \dots,
   \tau_{r+s},\ov{\tau}_{r+s}, $$
where $\tau_1, \dots, \tau_r$ denote the real, $\tau_{r+1}$,
$\ov{\tau}_{r+1}, \dots, \tau_{r+s}, \ov{\tau}_{r+s}$ the pairs of
complex conjugate embeddings. If $|\cdot|$ is the
ordinary absolute value on $\R$ (respectively $\C$), then by
$$ |x_i| =
\begin{cases}
|\tau_i(x)| & \text{for } 1 \le i \le r \\
|\tau_i(x)|^2 & \text{for } r+1 \le i \le r+s
\end{cases} $$
all $r+s$ distinct archimedean valuations\label{pBew} of $K$
are defined. Thereby $N_{K/\Q}(x) = |x_1| \cdots |x_{r+s}|$.
We now claim

\begin{quote}
  {\bf (1.13)} {\em Let $l=r+s-1$ independent units
    $u_1, \dots, u_l$ be given, and let
    $\kappa_i = \big| \log |u_1|_1 + \dots + \log |u_l|_i \big|$ for
    $i=1,\dots,l$; further let $k_i = \exp(\kappa_i)$. If then
    $c_i$ ($1 \le i \le l$) are any positive real numbers,
    then for every $z_1 \in K\setminus\{0\}$ there exists a unit
    $u \in \cO_K^\times$ with $c_i \le |z_1 u|_i < c_i k_i$
    for $i=1, \dots, l$. }
\end{quote}

\begin{proof}
  We define a map $A : K\setminus\{0\} \to \R^l$ by
  \[ A(x) = \begin{pmatrix} \log |x|_1 \\ \vdots \\ \log |x|_l \end{pmatrix}; \]
  $A$ of course depends on which of the $r+s$ valuations of $K$
  we select here. We now claim that the $l$ vectors $v_1 =
  A(u_1), \dots, v_l = A(u_l)$ are independent in $\R^l$. To this end we assume
  that $\sum a_i v_i = 0$ for certain $a_i \in \Z$. Considering
  this equation coordinate-wise, one obtains
  \[
  \sum_{i=1}^{l} a_i \log |u_i|_j = 0 \quad \text{for all } j
  \text{ with } 1 \le j \le l; \]
  adding these $l$ equations and noting
  \[
  \sum_{i=1}^{l} \log |u_i|_j = 0 \quad \text{(because }
  |N_{K/\Q}(u_i)| = 1\text{)},
  \]
  one obtains
  \[ \sum_{i=1}^{l} a_i \log |u_i|_j = 0 \quad \text{also for } j=r+s.
    \text{ Then however } \prod_{i=1}^{l} u_i^{a_i} = 1. \]
  and since the $u_i$ are independent, it follows that $a_i=0$ for all $i$
  ($1 \le i \le l$). Thus the vectors $v_1, \dots, v_l$ are in fact
  linearly independent. Now we can bring an arbitrary vector in $\R^l$
  into the fundamental domain
  \[ F = \left\{ \sum b_i v_i : 0 \le |b_i| \le \frac{1}{2} \right\} \]
  by shifting by suitable integer multiples of the $v_i$.
  For a $w = (w_1, \dots, w_l)^t \in F$ one now has
  \[ w_i = \sum_{j=1}^{l} b_j \cdot \log |u_j|_i = 0, \]
  hence $0 \le |w_i| \le \frac12\kappa_i$. Setting
  \[ w = \begin{pmatrix} \log |z|_1 - \log c_1 - \frac12 \kappa_1 \\
    \vdots \\
    \log |z|_l - \log c_l - \frac12 \kappa_l
  \end{pmatrix}
  \]
  and finding $a_i \in \Z$ so that $w - \sum a_i v_i \in F$,
  we obtain the existence of a unit $u = \prod u_i^{a_i}$
  with
  $$ \begin{array}{rccl}
    \log c_i & \le & \log |z_1 u|_i & < \log c_1 + \kappa_1 \\
                   & & \ldots & \\
    \log c_l & \le & \log |z_1 u|_l & < \log c_l + \kappa_l
    \end{array} $$
  From this our claim follows immediately.
\end{proof}
\begin{proof}[Proof of (1.12)]
  Since $M(K,x_i) < k$ there exists a $y \in R$ with
  $|N_{K/\Q}(x_i - y)| < k$. We set $z_i = x_i - y$ and
  find with (1.13) a unit $u$ with
  $$ \sqrt[n]{k} \le |z_1 u|_i < \sqrt[n]{k} \cdot \kappa_i $$
  for $i=1,\dots,l$. With $z = z_1 u$ (a)
  and (b) are evidently satisfied, and in order also to verify (c) we note
  \[
  |z|_{l+1} = |N_{K/\Q}(z)| / \prod_{j=1}^{l} |z|_j < \sqrt[n]{k},
  \]
  with $x_{l+1} = 1$ hence $|z|_{i} < \kappa_i \cdot \sqrt[n]{k}$ for
  $i=1, \dots, l+1$.
  Now let $\{\beta_1, \dots, \beta_n\}$ be the dual basis to $\{\alpha_1, \dots, \alpha_n\}$.
  As already in the cubic case it now follows also here that
  \[
  |r_j| = |T_{K/\Q}(z \beta_j)| \le \sum_{\tau} |z^\tau| \cdot
  |\beta_j^\tau| < \sqrt[n]{k} \cdot \Bigl( \sum_{\tau} |\lambda_\tau|
  \cdot |\beta_j^\tau| \Bigr) =: \mu_j.
  \]
  Here the $\lambda_\tau$ coincide with the $x_i$ if
  $\tau$ is real ($i=1, \dots, r$), and after suitable
  renumbering with $\sqrt{\kappa_i}$ if $\tau$ is not real
  ($i=r+1, \dots, r+2s$). This annoying renaming stems from the
  definition of the $|\cdot|_i$ before (1.13), where for non-real
  embeddings the square of the ordinary absolute value stands. I have
  however not yet seen how this deficiency can be remedied.
\end{proof}
Undoubtedly the estimates that we have made in the course of the
proof of (1.12) can still be improved. If one therefore
wishes to use (1.12) for the determination of $M(K)$, one will have to give
thought to this.
\subsection*{Tables}
Below we give some tables that are intended to show how one can apply the
criteria 1.5., 1.7.\ and 1.8.\ in real quadratic number fields.
From these tables one can already see that the
rings $D(m)$ with $m \le 101$ can be norm-Euclidean at most for the values
$m = 2, 3, 5, 6, 7, 11, 13, 17, 19, 21, 29, 33, 37, 41, 57, 73$;
in \S\  2 and \S\  3 we shall show that
these rings are the only norm-Euclidean real quadratic number fields.
\begin{table}[ht!]
\centering
\begin{tabular}{ccccc}
\toprule
$m$ & $f$ & $a$ & $k = M(x)$ & $r$ \\
\midrule
7 & 14 & 9 & $9/14$ & $-5$ \\
10 & 10 & 5 & $3/2$ & $5$ \\
11 & 22 & 3 & $19/22$ & $3$ \\
13 & 13 & 4 & $4/13$ & -- \\
15 & 15 & 6 & $7/5$ & $6, -9$ \\
17 & 17 & 8 & $8/17$ & -- \\
19 & 38 & 7 & $31/38$ & $7$ \\
21 & 7 & 2 & $5/7$ & $2$ \\
22 & 22 & 5 & $27/22$ & $5, -17$ \\
23 & 46 & 31 & $77/46$ & $-15, 31, -61$ \\
29 & 29 & 6 & $23/29$ & $6$ \\
31 & 31 & 14 & $45/31$ & $14, -17$ \\
33 & 11 & 5 & $6/11$ & $5$ \\
35 & 35 & 15 & $17/7$ & $15, -20, 50, -55$ \\
37 & 37 & 10 & $27/37$ & $10$ \\
47 & 94 & 65 & $253/94$ & $-29, 65, -123, 159, -217$ \\
51 & 102 & 19 & $287/102$ & (e=23): $-185, -83, 19, 121, 223$ \\
53 & 53 & 15 & $68/53$ & $15, -38$ \\
57 & 19 & 5 & $14/19$ & $5$ \\
59 & 59 & 7 & $125/59$ & $7, -52, 66, -111$ \\
69 & 23 & 2 & $25/23$ & $2, -21$ \\
77 & 11 & 3 & $19/11$ & $3, -8, 14$ \\
79 & 158 & 111 & $585/158$ & $-47, 111, -205, 269, -363, 427, -521$ \\
83 & 166 & 33 & $631/166$ & $33, -133, 199, -299, 365, -465, 531$ \\
85 & 85 & 19 & $151/85$ & $19, -66, 104$ \\
87 & 58 & 53 & $169/58$ & $-5, 53, -63, 111, -121$ \\
93 & 31 & 18 & $44/31$ & $-13, 18$ \\
101 & 101 & 24 & $125/101$ & $24, -77$ \\
\bottomrule
\end{tabular}
\caption{Application of (1.5) to quadratic
  number fields $\Q(\sqrt{m})$}\label{App15}
\end{table}

\medskip\noindent

{\bf Remark}\ on Tab.~\ref{App15}: If $e^2 \equiv a \bmod {f}$ and
$(c)^2 = (f)$ for some $c \in R$, set $x = e/c$; one thereby obtains
an $x \in K$ with $M(K,x)=k$.

In the case $m=51$ it is still to be remarked that 121 is of course a
norm from $S$ because $121 = N_{K/\Q}(11)$; however $11 \not\equiv e
\equiv 23 \bmod {f}$. This difficulty arises here because the
congruence $x^2 \equiv a \equiv 19 \bmod {f}$ has four distinct
solutions, namely $x=11$ and $x=23$.

\begin{table}[h]
\centering
\begin{tabular}{ccccc}
\toprule
$m$ & $B'$ & $y$ & $k=M(x)$ & $r$ \\
\midrule
6 & 4 & $1+\sqrt{6}$ & $3/4$ & $1$ \\
10 & 4 & $\sqrt{10}$ & $3/2$ & $2$ \\
14 & 4 & $1+\sqrt{14}$ & $5/4$ & $1, 3$ \\
15 & 4 & $1+\sqrt{15}$ & $3/2$ & $2$ \\
26 & 4 & $\sqrt{26}$ & $5/2$ & $2, 6$ \\
30 & 4 & $\sqrt{30}$ & $3/2$ & $2$ \\
20 & & $1+\sqrt{30}$ & $29/20$ & $-19, -9, 1, 11, 21$ \\
34 & 4 & $1+\sqrt{34}$ & $9/4$ & $1, 3, 5, 7$ \\
38 & 4 & $1+\sqrt{38}$ & $11/4$ & $1, 3, 5, 7, 11$ \\
39 & 4 & $1+\sqrt{39}$ & $5/2$ & $2, 6$ \\
42 & 4 & $1+\sqrt{42}$ & $7/4$ & $1, 3, 5$ \\
55 & 4 & $\sqrt{55}$ & $9/4$ & $1, 3, 5, 7$ \\
58 & 4 & $\sqrt{58}$ & $3/2$ & $2$ \\
62 & 4 & $1+\sqrt{62}$ & $13/4$ & $1, 3, 5, 7, 11, 13$ \\
66 & 4 & $1+\sqrt{66}$ & $15/4$ & $1, 3, 5, 7, 11, 13$ \\
74 & 4 & $\sqrt{74}$ & $5/2$ & $2, 6$ \\
78 & 4 & $\sqrt{78}$ & $7/2$ & $2, 6, 10$ \\
82 & 4 & $\sqrt{82}$ & $9/2$ & $2, 6, 10, 14$ \\
91 & 4 & $1+\sqrt{91}$ & $5/2$ & $2, 6$ \\
95 & 4 & $1+\sqrt{95}$ & $7/2$ & $2, 6, 10$ \\
\bottomrule
\end{tabular}
\caption{Application of (1.7) to quadratic
  number fields $\Q(\sqrt{m})$}\label{TabB}
\end{table}

\medskip\noindent
{\bf Remark}\ on Tab.~\ref{TabB}: If $B'=(c)^2$ for some $c
\in R$, then $M(x) = k$ for $x=y/c$.
\begin{table}[ht!]
\centering
\begin{tabular}{cccc}
\toprule
$m$ & $x$ & $M(x)$ & $u$ \\
\midrule
19 & $(0, 20/57)$ & $170/171$ & $170 + 39\sqrt{19}$ \\
21 & $(2/5, 1/5)$ & $12/25$ & $(5 + \sqrt{21})/2$ \\
22 & $(0, 9/28)$ & $443/392$ & $197 + 42\sqrt{22}$ \\
29 & $(3/5, 1/5)$ & $4/5$ & $(5 + \sqrt{29})/2$ \\
33 & $(1/2, 2/11)$ & $29/44$ & $23 + 4\sqrt{33}$ \\
34 & $(1/2, 6/17)$ & $135/68$ & $35 + 6\sqrt{34}$ \\
    & $(0, 5/12)$ & $137/72$ & $35 + 6\sqrt{34}$ \\
35 & $(0, 3/7)$ & $17/7$ & $6 + \sqrt{35}$ \\
38 & $(0, 5/12)$ & $173/72$ & $37 + 6\sqrt{38}$ \\
39 & $(0, 5/12)$ & $107/48$ & $25 + 4\sqrt{39}$ \\
41 & $(15/32, 5/32)$ & $23/32$ & $32 + 5\sqrt{41}$ \\
42 & $(0, 7/12)$ & $41/24$ & $13 + 2\sqrt{42}$ \\
43 & $(1/118, 1/2)$ & $11829/6962$ & $3482 + 531\sqrt{43}$ \\
    & $(0, 193/387)$ & $5902/3483$ & $3482 + 531\sqrt{43}$ \\
46 & $(1/2, 311/1058)$ & $76877/48668$ & $24335 + 3588\sqrt{46}$ \\
53 & $(1/14, 3/14)$ & $9/7$ & $(7 + \sqrt{53})/2$ \\
55 & $(1/2, 19/44)$ & $351/176$ & $89 + 12\sqrt{55}$ \\
57 & $(0, 15/76)$ & $219/304$ & $151 + 20\sqrt{57}$ \\
61 & $(1/78, 17/78)$ & $41/39$ & $(39 + 5\sqrt{61})/2$ \\
62 & $(1/2, 48/124)$ & $367/124$ & $63 + 8\sqrt{62}$ \\
    & $(0, 7/16)$ & $367/128$ & $63 + 8\sqrt{62}$ \\
65 & $(1/4, 1/4)$ & $1$ & $8 + \sqrt{65}$ \\
66 & $(0, 7/16)$ & $31/128$ & $65 + 8\sqrt{66}$ \\
67 & $(1/2, 7/18)$ & $341/162$ & $48842 + 5967\sqrt{67}$ \\
70 & $(1/2, 6/25)$ & $891/500$ & $251 + 30\sqrt{70}$ \\
71 & $(0, 216/497)$ & $7393/3479$ & $3480 + 413\sqrt{71}$ \\
73 & $(41, 283)/2136$ & $1541/2136$ & $1068 + 125\sqrt{73}$ \\
85 & $(4/9, 2/9)$ & $16/9$ & $(9 + \sqrt{85})/2$ \\
86 & $(0, 133/473)$ & $10030/5203$ & $10405 + 1122\sqrt{86}$ \\
89 & $(3, 159)/1000$ & $1004/1000$ & $500 + 53\sqrt{89}$ \\
94 & $(0, 3661/14194)$ & $4708623/2143294$ & $2143295 + 221064\sqrt{94}$ \\
97 & $(1496, 2587)/11208$ & $3001/2802$ & $5604 + 569\sqrt{97}$ \\
109 & $(17, 82)/261$ & $289/261$ & $261 + 25\sqrt{109}$ \\
113 & $(3, 219)/1352$ & $967/776$ & $776 + 73\sqrt{113}$ \\
137 & $(3, 447)/3488$ & $2177/1744$ & $1744 + 149\sqrt{137}$ \\
\bottomrule
\end{tabular}
\caption{Application of (1.8) to quadratic number fields
  $\Q(\sqrt{m})$ ($t=1$):}
\end{table}
\begin{table}[h]
\centering
\begin{tabular}{cccc}
\toprule
$m$ & $x_1, x_2$ & $M(x_1)$ \\
\midrule
19 & $(0, 173/494)$ & \multirow{2}{*}{$10579/12844$} \\
    & $(1/2, -115/247)$ & \\
22 & $(1182, 496441/155236)$ & \multirow{2}{*}{$175903/155236$} \\
    & $(-1182, 496441/155236)$ & \\
58 & $(1/2, 197/754)$ & \multirow{2}{*}{$27477/19604$} \\
    & $(1/2, 276/754)$ & \\
61 & $(0, 66/305)$ & \multirow{2}{*}{$1611/1525$} \\
    & $(0, 67/305)$ & \\
74 & $(1/2, 19/43)$ & \multirow{2}{*}{$16255/7396$} \\
    & $(1/86, 1/2)$ & \\
93 & $(0, 44/279)$ & \multirow{2}{*}{$2198/1674$} \\
    & $(1/2, 119/558)$ & \\
\bottomrule
\end{tabular}
\caption{$t=2$}
\end{table}
\medskip
For further examples with $t=2$ see \S\  2.

\subsection*{Sufficient Condition for the Euclidean Algorithm}

After we have so far occupied ourselves with necessary criteria, we
now turn to criteria that guarantee the existence of the Euclidean
algorithm.

In 1974 Lenstra\index[N]{Lenstra} presented a method with which
in the meantime (see e.g.\ Lenstra \cite{Len77a},
Leutbecher\index[N]{Leutbecher} and Martinet\index[N]{Martinet}
\cite{LM81a,LM81b}, Leutbecher \cite{Leu85,Leu86}, Leutbecher and
Niklasch\index[N]{Niklasch} \cite{LN87}) about 400 norm-Euclidean
number fields have been found. A modification of this criterion will
allow us to obtain in a similar manner also $k$-stage
norm-Euclidean number rings (these we shall simply call $k$-Euclidean).

Lenstra's method is based on an idea of Hurwitz\index[N]{Hurwitz}
\cite{Hur19}: the latter has shown that in a number field $K$ there exists a
natural number $m \ge 2$ depending only on $K$ such that:
\begin{quote}
  {\em For all $x \in K$ there exist a $y \in R$ and a $j \in \N$,
    with $|N_{K/\Q}(j(x-y))| < 1$ for all $1 \le j < m$.}
\end{quote}

Note that $K$ is norm-Euclidean precisely if we could choose $m=2$.
Lenstra then observed that in the proof of this
statement it is not essential that the above $j$ be natural numbers,
but that somewhat more generally one has:

\begin{quote}
  {\bf (1.14)} {\em Let $K$ be a number field, $(K:\Q) = n = r+2s$, and
    $d = |\disc K|$. If then $m \in \N$,
    $$ m > M_{r,s}\label{pMrs} =
          \frac1d \cdot \frac{n!}{n^n} \cdot \Big(\frac4\pi\Big)^s $$
    and $\theta_1, \dots, \theta_m$ are
    any pairwise distinct elements of $K$, then for
    all $x \in K$ there exists a $y \in R$ with $|N_{K/\Q}((\theta_i -
    \theta_j)x - y)| < 1$ for certain
    $i, j \in \{1, \dots, m\}$, $1 \le i < j \le m$.}
\end{quote}

Lenstra's idea was the following: if we can choose the $\theta_1,
\dots, \theta_m$ in such a way that the differences $\theta_i -
\theta_j$ for $i < j$ are all units in $R$, then $c =
\frac{y}{\theta_i - \theta_j} \in R$, and (1.14) then states precisely
that for all $x \in K$ there exists a $c \in R$ with $|N_{K/\Q}(x-c)|
< 1$: $R$ would thereby be norm-Euclidean.

The proof of (1.14) that we now wish to present is due to
Lenstra\index[N]{Lenstra} \cite{Len74} and has great similarity with
the classical proof of the finiteness of the class number (in 1977
Lenstra\index[N]{Lenstra} \cite{Len77a} gave a further proof in the
language of packing theory). Some of the statements below
we shall therefore take as known; for proofs
we refer to Marcus\index[N]{Marcus} \cite{Mar77} or
Cassels\index[N]{Cassels} \cite{Cas86}.

\begin{proof}[Proof of (1.14)]
  Now let the embeddings $\tau_1, \dots, \tau_n$ of $K$ into $\C$
  be arranged as before (1.12); we then define an embedding of
  $K$ into $\R^n$ by $\varphi : K \to \R^n$ via
  $$ \alpha \mapsto \varphi(\alpha) =
   (\tau_1(\alpha), \dots, \tau_r(\alpha), \Re
   \tau_{r+1}(\alpha), \Im \tau_{r+1}(\alpha), \dots,
   \Re \tau_{r+s}(\alpha), \Im \tau_{r+s}(\alpha)), $$   
  where $\Re$ respectively $\Im$ denote the real respectively
  imaginary part of a complex number. On $\R^n$ we now define a norm
  for $x = (x_1, \dots, x_n) \in \R^n$ by
\[
N(x) = x_1^2 \cdots x_r^2 (x_{r+1}^2 + x_{r+2}^2) \cdots (x_{n-1}^2 + x_n^2).
\]
Thereby $N_{K/\Q}(\alpha) = N(\varphi(\alpha))$ for all
$\alpha \in K$. Now we set
\[
A = \left\{ x \in \R^n : |x_1| + \cdots + |x_r| +
2 \sqrt{x_{r+1}^2 + x_{r+2}^2} + \cdots + 2 \sqrt{x_{n-1}^2 + x_n^2} < n
\right\}.
\]
It is well known that $A$ is bounded and convex; moreover the
inequality between the geometric and arithmetic means yields the
relation $|N(\alpha)| < 1$ for all $\alpha \in A$. With $U =
\frac{1}{2}A$ one then has $u - v \in A$ for all $u,v \in U$, so that for
such $u,v$ one has $|N(u-v)| < 1$. Finally one still has
$$ \operatorname{vol}(A) = 2^r (\frac{\pi}2)^s \frac{n^n}{n!}
  \quad \text{and} \quad \
  \operatorname{vol}(U) = 2^n \cdot \operatorname{vol}(A), $$
while the fundamental parallelepiped $F$ of the lattice $\Gamma_R := \varphi(\cO_K)$
in $\R^n$ has volume $\operatorname{vol}(F) = 2^{-s} \sqrt{d}$.
In what follows we shall identify the $\alpha \in K$ with their images
$\varphi(\alpha) \in \R^n$; since $\varphi$ is
injective and $N_{K/\Q}(\alpha) = N(\varphi(\alpha))$,
no problems should arise.

Now let $\vartheta_1, \dots, \vartheta_m \in K$ be pairwise
distinct; then for every point $z$ from one of the sets
$\vartheta_i x + U$ ($i=1,\dots,m$) we can find a $y \in R$
such that $z_0 - z = y$ lies within the fundamental parallelepiped
$F$ of the lattice $R$. Now however the total volume of the sets
$\vartheta_i x + U$ equals $m \cdot \operatorname{vol}(U) >
  M_{r,s} \cdot \operatorname{vol}(U) = 2^{-s} \sqrt{d} =
  \operatorname{vol}(F)$. Consequently there exist indices $i,j$ and numbers
$y_1, y_2 \in R$ such that
\[
(\vartheta_i x + U - y_1) \cap (\vartheta_j x + U - y_2) \neq \emptyset
\]
and the pairs $(i,y_1)$ and $(j,y_2)$ are distinct. Thus
we have found $u,v \in U$ with $\vartheta_i x + u - y_1 =
\vartheta_j x + v - y_2$. If $i=j$, then $u - v = y_1 - y_2 \in
R$ would follow; because of $|N(u-v)| < 1$ one must then have $N(y_1 - y_2) = 0$,
and this implies $y_1 = y_2$. This however contradicts the fact
that the pairs $(i,y_1)$ and $(j,y_2)$ are distinct. Therefore $i
\neq j$ and $(\vartheta_i - \vartheta_j)x - y = v - u$ for $y := y_1 -
y_2$. Taking norms now yields
\[
|N((\vartheta_i - \vartheta_j)x - y)| = |N(u - v)| < 1
\]
as claimed.
\end{proof}
We have already remarked that from (1.14) it follows that:

\begin{quote}
  {\bf (1.15)} {\em If there exist numbers $\vartheta_1, \dots, \vartheta_m
    \in K$ whose differences are units in $R$ and if $m >
    M_{r,s}$, then $R$ is norm-Euclidean. }
\end{quote}

A sequence $\vartheta_1, \dots, \vartheta_m$ of elements of $K$
whose differences are all units in $R$ we call an
\textbf{exceptional sequence}.\index[S]{exceptional sequence}
Since only the differences of the $\vartheta_i$ matter, we may without
loss of generality assume $\vartheta_1=0$ (thus all $\vartheta_i \in R$).
Since if $\vartheta_i - \vartheta_j$ is a unit, then so is
$(\vartheta_i - \vartheta_j)/\vartheta_2$,
we see that if $0, \vartheta_2, \dots, \vartheta_m$ isan  exceptional
sequence, then so is
$0, 1, \vartheta_3/\vartheta_2, \dots, \vartheta_m /\vartheta_2$.
Thus we may also assume $\vartheta_2 = 1$.

An example of an exceptional sequence of length $p$ is then in $K =
\Q(\zeta)$, where $\zeta = \zeta_p$ denotes a primitive $p$-th root of
unity:
\[
0, 1, 1+\zeta, 1+\zeta+\zeta^2, \dots, 1+\zeta + \zeta^2 + \dots + \zeta^{p-2};
\]
this shows that the ring $\Z[\zeta]$ is norm-Euclidean if $p >
M_{r,s}$, and this is the case for $p = 3, 5, 7$.  In order to be able
to formulate a result analogous to (1.15) also for $k$-Euclidean
rings,we define: a sequence $\vartheta_1, \dots, \vartheta_m$ is
called a \textbf{$k$-sequence},\index[S]{$k$-sequence}
if the following two conditions are satisfied:
\begin{enumerate}
\item[] (F--1) one has $\vartheta_i - \vartheta_j \in E_k \setminus \{0\}$
  for all $1 \le i < j \le m$;
\item[] (F--2)\label{pF2} every divisor of the ideal
  $(\vartheta_i - \vartheta_j)$ is a principal ideal.
\end{enumerate}
Thus 1-sequences and exceptional sequences are the same, because then the
condition (F--2) is automatically satisfied. We now claim

\begin{quote}
  {\bf (1.16)} {\em If there exists in $R$ a $k$-sequence of length
    $m > M_{r,s}$, then $R$ is $k$-Euclidean.}
\end{quote}

\begin{proof}
  Let $x \in K$ be given; with (1.14) we find a $y \in R$ with
  $|N((\vartheta_i - \vartheta_j)x - y)| < 1$, where $\vartheta_1,
  \dots, \vartheta_m$ is a $k$-sequence of length $m$. By (F--2)
  we may write $y/(\vartheta_i - \vartheta_j) = a/b$ for
  certain $a, b \in R$ with $(a, b)=1$ and $b|(\vartheta_i -
  \vartheta_j)$. By (0.12) one therefore has $b \in E_k$, hence $a/b$ by
  (0.11) is a continued fraction of length $\le k$ with denominator $ub$ for a
  unit $u \in R^\times$, and finally
  \[
  \Big|N_{K/\Q}\Big(x - \frac{a}{b}\Big)\Big| <
  \frac{1}{N_{K/\Q}(\vartheta_i - \vartheta_j)} <
  \frac{1}{N_{K/\Q}(b)}.
  \]
  with (1.14). Now (0.9) shows that $R$ is indeed $k$-Euclidean.
\end{proof}

Note that a sequence $\vartheta_1, \dots, \vartheta_m$ with the
property $\vartheta_i - \vartheta_j \in E_k$ is already a
$k$-sequence: (F--1) is satisfied since $E_k \subset E_k'$, and
(0.5.(i)) shows (F--2). The reason for the introduction of the sets
$E_k'$ lies solely in the fact that these sets are in general strictly
larger than the corresponding $E_k$, which of course facilitates the
finding of (as long as possible) $k$-sequences.

If $I$ is an integral ideal in $R$ and there exists a prime residue
system mod $I$ that consists entirely of units of $R$, then we say
that $I$ possesses a prime unit residue system (for short: a
PERS).\label{pPERS} In particular $I$ has a PERS if there exists a
unit $u \in R$ that is a primitive root mod $I$. On the other hand,
even if $R$ has unit rank 1, there exist ideals that possess a PERS
but no unit as primitive root.

We give some simple examples:

\begin{enumerate}
\item[1.] $K = \Q(\sqrt{5})$: here $\omega = (1+\sqrt{5})/2 \in
  R$, and $0,1, \omega, 1+\omega$ is a 1-sequence, since $1, \omega,
  1+\omega, \omega-1$ are units. On the other hand $0, 1, -1+\omega,
  \omega, 1+\omega$ is not a 1-sequence, because e.g.\ the difference $-2 =
  (\omega-1)-(\omega+1)$ is not a unit. However we have here
  a 2-sequence, because the ideals $I_1=(2)$ and $I_2=(2+\omega)$ possess a
  PERS: $\{1, \omega, 1+\omega\}$ is a prime
  residue system mod $I_1$ and $\{1, \omega, \omega^2, \omega^3,
  \omega^4\}$ one mod $I_2$ (in both cases $\omega$ is a
  primitive root).
\item[2.] $K=\Q(\sqrt{14})$: here $0,1$ is a 1-sequence of length 2,
  and one easily sees that there are no longer ones; for if $u$ and
  $u-1$ are units in $R$, then they are (up to a possible sign) powers
  of the fundamental unit $\varepsilon = 15+4\sqrt{14}$. Now however
  $\varepsilon \equiv 1 \bmod {2}$, so $u$ and $u-1$ would also have
  to be $\equiv 1 \bmod {2}$, which evidently leads to a
  contradiction.
  
  We now claim that $0, 1, 4+\sqrt{14}, 5+\sqrt{14}$ is a
  2-sequence. For this we must show that the ideals
  $I_1=(3+\sqrt{14})$, $I_2=(4+\sqrt{14})$ and $I_3=(5+\sqrt{14})$ possess a
  PERS. In \S\  0 (p.~\pageref{p11}) we have already seen that
  $3+\sqrt{14} \in E_2'$, and that $\varepsilon$ is here a
  primitive root. Because $\|I_2\| = 2$, $\{1\}$ is a PERS of
  $I_2$; finally one checks that $\varepsilon$ is also
  a primitive root mod $I_3$. Since here $M_{r,s} = 3.74\dots$ and
  we have found a 2-sequence of length 4, we can with
  (1.16) conclude that $\Q(\sqrt{14})$ is 2-Euclidean.
\end{enumerate}

In general it is rather tedious to find in a given number field $K$ a
1-sequence (or a $k$-sequence) of sufficiently large length;
before one begins the search for $k$-sequences, one will therefore
want to know whether this effort is worthwhile at all, i.e.\ whether
there can exist in $R$ $k$-sequences of the desired length. In order
to be able to answer this question we set\label{Lenmu}
\begin{align*}
  \mu_k & = \sup \{ m \in \N : \text{there exists a $k$-sequence of length } m
       \text{ in } R \} \qquad \text{and} \\
  \lambda_k & = \inf \{ \|I\| :
       \text{there exists a prime residue class mod } I \\
       & \qquad \text{ that has no representative in } E_0' = \{0\} \}.
\end{align*}
In particular
$\lambda_1 = \min \{ \|I\| \ge 2: \ I \text{ is an integral ideal in } R \}$,
since the prime residue class $1 \bmod I$ has no representative in $E_1 = \{0\}$
(in particular $\lambda_1 \le 2$, since one can choose $I=(2)$),
and $\lambda_2 = \inf \{ \|I\| \ge 2: \ I \text{ possesses no PERS in } R \}$.
Now we have

\begin{quote}
  {\bf (1.17)} {\em If $\lambda_k$ is finite, then one has
    $2 \le \mu_1 \le \mu_2 \le \dots \le \mu_k \le \lambda_k$.}
\end{quote}

\begin{proof}
  Since $0,1$ is a 1-sequence, one has $\mu_1 \ge 2$; further every
  $k$-sequence is a fortiori a $(k+1)$-sequence, hence $\mu_k \le \mu_{k+1}$
  for all $k \in \N$. Thus we only still need to show $\mu_k \le \lambda_k$.
  To this end let $I$ be an integral ideal in $R$ that possesses a prime
  residue class mod $I$ which has no representative in $E_k$
  (such an ideal exists since $\lambda_k$ is finite). If then
  $\vartheta_1, \dots, \vartheta_m$ is a $k$-sequence, one must have
  $\vartheta_i - \vartheta_j \not\equiv 0 \bmod {I}$ for all $i
  \ne j$: otherwise namely $I$ as a divisor of the ideal
  $(\vartheta_i - \vartheta_j)$ would be a principal ideal by (F--2), hence
  $I=(b)$ for some $b \in R$. Because of $b \mid (\vartheta_i - \vartheta_j)$
  and (0.12) one now has $b \in E_k$, while by hypothesis the ideal
  $I=(b)$ possesses a prime residue class that has no representative
  in $E_k$: contradiction! Thus $\vartheta_1, \dots, \vartheta_m$ are
  pairwise incongruent mod $I$ and $m \le \|I\|$. This however implies
  $\mu_k \le \lambda_k$.
\end{proof}

We wish to consider the examples 1.\ and 2.\ given above on page 36
once more in the light of (1.17): in $\Q(\sqrt{5})$ one has
$\lambda_1=4$, because $I=(2)$ is the ideal of minimal norm $>1$. For
every 1-sequence $\vartheta_1, \dots, \vartheta_4$ the $\vartheta_i$
are therefore pairwise incongruent mod $2$, as we have seen in the
proof of (1.17), and for $0,1,\omega,1+\omega$ this is evidently also
the case. Further we can here conclude $\mu_1 = \lambda_1$.  In
$K=\Q(\sqrt{14})$ on the other hand $I_2=(4+\sqrt{14})$ is an ideal of
norm 2, so that here we have $\mu_1 = \lambda_1 = 2$. Further
$\lambda_2 = 4$, since the ideal $I=(2)$ has no PERS. The given
2-sequence of length 4 finally shows $\mu_2 = \lambda_2 = 4$.  While
for $k \ge 3$ I have not succeeded in showing that the $\lambda_k$ are
finite, for $k=2$ one has:

\begin{quote}
  {\bf (1.18)}{\em Let $K$ be a number field and $q$ a rational
    prime. Then there exists an $a \in \N$ such that $I = (q^a)R$ possesses no
    PERS. In particular $\lambda_2 \le (q^a)^n$, where
    $n=(K:\Q)$ is the degree of the field.}
\end{quote}

\begin{proof}
  Let $N$ be the normal closure of $K/\Q$. If $I = (q^a)$ possesses a
  PERS, then for every $r \in \Z$ with $q \nmid r$ there exists a unit
  $u \in R^\times$ with $r \equiv u \bmod {I}$. As $\sigma$ runs through the
  embeddings of $K$ into $\C$, for every $\sigma$ the
  congruence $u^\sigma \equiv r \bmod {I}$ follows in $\cO_N$ (= the
  ring of integers of $N$). Since $r^\sigma = r$ and $1^\sigma = 1$.
  Multiplication over all $\sigma$ then gives $1 = N_{K/\Q}(u) =
  \prod u^\sigma \equiv r^n \bmod {I}$ in $\cO_N$. Because $I
  \cap \Z = q^a \Z$ the congruence $1 \equiv r^n \bmod {q^a}$ also holds
  in $\Z$. By now choosing $r > 1$ and $a$ so large that $q^a >
  r^n + 1$, we obtain the desired contradiction.
\end{proof}
We further remark:
\begin{enumerate}
\item[1.] If $\vartheta_1, \dots, \vartheta_m$ is a 1-sequence in $K$
  and $K \subset L$, then $\vartheta_1, \dots, \vartheta_m$ is also
  a 1-sequence in $L$; for $k$-sequences with $k \ge 2$ this no longer holds
  in general, because the norm of the ideals $(\vartheta_i - \vartheta_j)$
  depends on $L$ if $\vartheta_i - \vartheta_j$ is not a unit.
\item[2.] As Lenstra\index[N]{Lenstra} \cite[p.\ 239]{Len77a} already remarked in the case $k=1$,
  by a small modification of the proof of
  1.13.\ one can obtain an upper bound for the Euclidean minima $M^k(K)$;
  namely under the hypothesis of (1.16) one has
  $M^k(K) \le M_{r,s}/\mu_k$.
\item[3.] Also the definition (1.17) of Lenstra \cite[p.\ 241]{Len77a} can
  be generalized to $k$-sequences with $k > 1$: here one considers
  sequences $\vartheta_1, \dots, \vartheta_m$ such that among $l+1$
  sequence members there are at least two whose difference lies in $E_k$
  (here $l$ is a natural number; for $l=1$ one recovers the
  ordinary $k$-sequences). If one denotes by $\mu_{k,l}$ the
  maximal length of such generalized $k$-sequences and if
  $\mu_{k,l} > l \cdot M_{r,s}$, then one can conclude similarly as in (1.16)
  that $K$ is $k$-Euclidean.
\item[4.] Under the assumption of the extended Riemann hypothesis
  Lenstra has shown \cite[p.\ 242f]{Len77a} that with (1.15) one can find only finitely
  many norm-Euclidean number fields; thereby the
  inequality $\lambda_1 \le 2^n$, $n=(K:\Q)$, which one
  obtains directly from the definition of $\lambda_1$ with $I = (2)$, was used essentially. The
  estimate (1.18) for $\lambda_2$ does not permit a similar conclusion for
  $k$-Euclidean rings with $k \ge 2$.
\end{enumerate}
\subsection*{The Gaussian General Measure}
If we wish to show that a number field $K$ is norm-Euclidean,
we must find for every $x \in K$ a $y \in R$ with
$|N_{K/\Q}(x-y)|<1$. Suppose we know a
function $\cM_K : K \to \R$ with the property
$|N_{K/\Q}(x)| \le \cM_K(x)$ for all $x \in K$, then
it evidently suffices to find for every $x \in K$ a $y \in R$ with
$\cM_K(x-y)<1$. Because $N_{K/\Q}(x) = \prod
\sigma(x)$, where the product runs over all embeddings $\sigma$ of $K$ into
$\C$, e.g.
\[
\cM_K(x) = \cM_{K/\Q}(x) = \frac{1}{n} \Bigl(
\sum |\sigma(x)|^2 \Bigr)
\]
comes into question as such a function; in fact one has\label{GMen}

\begin{quote}
  {\bf (1.19)} {\em Let $K$ be a number field with $(K:\Q) = n$; then for all $x \in K$ one has}
  \[ |N_{K/\Q}(x)| \le \cM_K(x)^{n/2}. \]
\end{quote}

\begin{proof}
  From the inequality between the geometric and arithmetic means
  one obtains
  \[ |N_{K/\Q}(x)|^{2/n} = \Bigl\{ \prod |\sigma(x)|^2
     \Bigr\}^{1/n} \le \Bigl( \sum |\sigma(x)|^2 \Bigr)^{1/n} =
     \cM_K(x). \]
\end{proof}
The function $\cM_K$ was already introduced by Gauss for cyclotomic fields (apart from the
factor $1/n$) (Werke II, p.\ 395) and was investigated extensively by
Cassels\index[N]{Cassels} in 1969.
Lenstra\index[N]{Lenstra} then succeeded in 1974, with the help of
$\cM_K$, in showing that the cyclotomic fields $\Q(\zeta_m)$ for $m =
1, 3, 4, 5, 7, 8, 9, 11, 12, 15, 20$ are norm-Euclidean.
Ojala\index[N]{Ojala} \cite{Oja77} in 1977 also proved
$\Q(\zeta_{16})$ to be norm-Euclidean and thereby (besides
a computer) likewise used $\cM_K$.\footnote{Recently
  R.\ McKenzie has shown with a computer that also $\Q(\zeta_{13})$
  is norm-Euclidean; see on this Leutbecher\index[N]{Leutbecher} \&
  Niklasch\index[N]{Niklasch} 1987.} The advantage of $\cM_K$ lies in the fact that
one can in general compute $\cM_K(x)$ much more simply than $N_{K/\Q}(x)$.
While however norm and trace are transitive in field towers, we cannot
expect this for $\cM_K$: for $x \in K$ one has namely
not necessarily $\cM_K(x) \in \Q$. For if
$\sigma$ is an embedding of $K$ into $\C$, then
$|\sigma(x)|^2 = \sigma(x) \ov{\sigma(x)}$ is a number from the
largest real subfield of the normal closure of $K/\Q$
(with $\ov{\phantom{x}}$ we have as usual denoted complex
conjugation). Matters are however somewhat simpler
if $K/\Q$ is abelian: then complex conjugation is an
element of the Galois group $\mathrm{Gal}(K/\Q)$ and therefore commutes with
every automorphism of $K/\Q$. Thereby it follows
easily that $(\cM_K(x))^\sigma = \cM_K(x)$ for every
$\sigma \in \mathrm{Gal}(K/\Q)$, and we indeed have
$\cM_K(x) \in \Q$ (essential in this proof is the
commutativity of complex conjugation and $\sigma$; therefore
$\cM_K(x) \in \Q$ holds even for all CM-fields,
see Washington\index[N]{Washington} \cite{Was82}).
Directly from the definition it follows for abelian $K$ (respectively for CM-fields
$K$) the relation $\cM_K(x) = \frac1n T_{K/\Q}(x \ov{x})$, again
observing the commutativity of $\sigma$ with the
other automorphisms. With the help of this relation Lenstra
\cite{Len75a,Len75b} proved the next two results for
cyclotomic fields and remarked that they also hold for CM-fields.
If however one dispenses with this formula, one obtains somewhat
more generally:

\begin{quote}
{\bf (1.20)} {\em Let $K \subseteq L$ be number fields with $n = (L:K)$;
  then for all $x \in L$ and all $y \in K$ one has the relation}
  \[ \cM_L(x) - \cM_L(x-y) = \cM_K\Bigl(\frac{1}{n}
     T_{L/K}(x)\Bigr) - \cM_K\Bigl(\frac{1}{n}
     T_{L/K}(x) - y\Bigr). \]
\end{quote}

This will be our substitute for the missing transitivity of
$\cM$; with (1.20) we shall namely be able to make statements about the Euclidean algorithm
in $L$ as soon as we know certain properties of $K$ sufficiently
well.
Before we prove this, we recall some known facts from
field theory (see e.g.\ Marcus\index[N]{Marcus} 1977): Let $\Q
\subset K \subset L \subset N$, $N$ the normal closure of $L/\Q$,
$(K:\Q)=k$, $(L:K)=n$, hence $(L:\Q) = kn$. Then there exist $k$
embeddings $\sigma_1, \dots, \sigma_k$ of $K$ into $\C$, as well as $n$
embeddings $\tau_1, \dots, \tau_n$ of $L$ into $\C$ that leave $K$
elementwise fixed. We can extend the $\sigma_i$ and $\tau_j$ to
automorphisms of $N$ (in general in many different
ways), and we shall denote these (once and for all chosen)
extensions likewise by $\sigma_i$ respectively $\tau_j$. Thereby
we obtain every embedding $\rho$ of $L$ into $\C$ as the restriction
of the $(L:\Q)$ automorphisms $\sigma_i \tau_j$ of $N$ to $L$. Thereby
we have in (1.20)
\begin{align*}
  (L:\Q) & \{\cM_L(x) - \cM_L(x-y)\} =
  \sum |x^\rho|^2 - \sum |(x-y)^\rho|^2 \\
  & = \sum (x^\rho \ov{x^\rho} - (x-y)^\rho \ov{(x-y)^\rho}) \\
  & = \sum (x^\rho \ov{y^\rho} +
    \ov{x^\rho} y^\rho - y^\rho \ov{y^\rho}).
\end{align*}
Now we note $\sigma(y) = y$ (for $\sigma$ leaves $K$ pointwise
fixed):
\[
= \sum_{\sigma, \tau} \{\sigma \tau (x) \ov{\sigma(y)} +
\ov{\sigma \tau (x)} \sigma(y) - \sigma(y)
\ov{\sigma(y)}\}.
\]
By pulling the summation over $\tau$ inwards and noting $T_{L/K}(x)
= \sum \tau(x)$, we obtain further
\begin{align*}
& = \sum_{\sigma} \{\sigma(T_{L/K}(x)) \ov{\sigma(y)} +
    \ov{\sigma(T_{L/K}(x))} \sigma(y) - n \cdot \sigma(y)
    \ov{\sigma(y)}\} \\
& = n \cdot \Bigl\{ \sum_{\sigma} \sigma\Bigl(\frac{1}{n} T_{L/K}(x)\Bigr)
    \ov{\sigma\Bigl(\frac{1}{n} T_{L/K}(x)\Bigr)} \\
& \qquad - \sigma\Bigl(\frac{1}{n} T_{L/K}(x) - y\Bigr)
   \ov{\sigma\Bigl(\frac{1}{n} T_{L/K}(x) - y\Bigr)} \Bigr\} \\
& = (L:\Q) \Bigl\{ \cM_K\Bigl(\frac{1}{n} T_{L/K}(x)\Bigr)
   - \cM_K\Bigl(\frac{1}{n} T_{L/K}(x) - y\Bigr) \Bigr\}
\end{align*}
We now claim

\begin{quote}
  {\bf (1.21)} {\em Let $L = K(\zeta_m)$, where $\zeta_m$ is an $m$-th
    root of unity. Then for all $x \in L$ one has the relation}
  $$ (L:K)\cM_L(x) = \frac{1}{m} \sum_{j=1}^{m}
     \cM_K(T_{L/K}(x \zeta_m^j)). $$
\end{quote}

For the case that $K$ is a cyclotomic field, (1.21) in
this form is due to Lenstra\index[N]{Lenstra} \cite{Len75a,Len75b}, but goes essentially
back to Cassels\index[N]{Cassels} \cite{Cas69}.
\begin{proof}
  In the following computations let $\tau$ and $\tau'$ run through all $(L:K)$
  embeddings of $L$ into $\C$ that leave $K$ elementwise fixed,
  while $\sigma$ runs through the embeddings of $K$ into $\C$. Then
  \begin{align*}
    (K:\Q) & \cdot \sum_{j=1}^{m} \cM_K(T_{L/K}(x \zeta_m^j)) =
    (K:\Q) \sum_{j=1}^{m} \cM_K\Bigl( \sum_{\tau} \tau(x
      \zeta_m^j) \Bigr) \\
     & = \sum_{j=1}^{m} \sum_{\sigma} \sigma\Bigl( \sum_{\tau} \tau(x
      \zeta_m^j) \Bigr) \ov{\sigma\Bigl( \sum_{\tau} \tau(x \zeta_m^j) \Bigr)} \\
     & = \sum_{\sigma} \sum_{\tau,\tau'} \sum_{j=1}^{m} \sigma(\tau(x))
         \ov{\sigma(\tau'(x))} \sigma(\zeta_m^j)
         \ov{\sigma(\zeta_m^j)},
  \end{align*}
  where we again regard $\sigma, \tau, \tau'$ as automorphisms of $N/\Q$.
  The image of a root of unity under such an
  automorphism is again a root of unity, and we find
  \[
  \ov{\sigma(\zeta_m^j)} = \sigma(\zeta_m^{-j}) =
  \sigma(\zeta_m^j)^{-1}.
  \]
  Thus
  \[
  (K:\Q) \sum_{j=1}^{m} \cM_K(T_{L/K}(x \zeta_m^j)) =
  \sum_{\sigma} \sum_{\tau,\tau'} \sum_{j=1}^{m} \sigma(\tau(x))
  \ov{\sigma(\tau'(x))} \sigma(\zeta_m^j)
  \ov{\sigma(\zeta_m^j)}.
  \]
  Now we set $\zeta_{\tau,\tau'} = \sigma(\zeta_m)
  \ov{\sigma(\zeta_m')}$ and find
  \[
  \zeta_{\tau,\tau'} = 1 \Leftrightarrow \tau(\zeta_m) =
  \tau'(\zeta_m) \Leftrightarrow \tau = \tau', \text{ as well as}
  \]
  \[
  \sum_{j=1}^{m} \zeta_{\tau,\tau'}^j =
  \begin{cases}
    m, & \text{if } \tau = \tau'; \\
    0, & \text{otherwise}.
  \end{cases}
  \]
  By pulling the summation over $j$ inwards, it follows that
  \begin{align*}
    (K:\Q) \sum \cM_K(T_{L/K}(x \zeta_m^j))
    & = m \cdot \sum_{\sigma} \sum_{\tau} \sigma(\tau(x))
    \ov{\sigma(\tau(x))}\\
    & = m \cdot \sum_{\sigma} |x^\sigma|^2 = m \cdot (L:\Q)
    \cdot \cM_L(x).
  \end{align*}
  By division by $m \cdot (K:\Q)$ we obtain the claim.
  \end{proof}
Lenstra's idea was now to reduce the existence of the Euclidean algorithm in $L=K(\zeta_m)$ to
properties of $K$. To this end one defines
\begin{align*}
  F & = F_K = \{x \in K: \cM_K(x) \le \cM_K(x-y)
        \text{ for all } y \in R\}, \quad \text{as well as} \\
  c(K) & = \sup \{ \cM_K(x) : x \in F_K \}.
\end{align*}
Directly from the definition it then follows that for every $x \in K$
there exists a $y \in R$ with $\cM_K(x-y) \le c(K)$. By (1.19)
$R$ is certainly norm-Euclidean if $c(K)<1$. In fact often already
$c(K) \le 1$ suffices, as we now wish to show.

\begin{quote}
  {\bf (1.22)} {\em Let $K$ be a number field, $\sigma$ an embedding
    of $K$ into $\C$, and suppose $|\sigma(x)| = |\sigma(x-u)| = 1$
    for some $x \in K$ and a root of unity $u \in R$; then
    $x \in R$.}
\end{quote}

\begin{proof}
  Because of $x^\sigma \ov{x^\sigma} = 1$ and $u^\sigma \ov{u^\sigma} = 1$
  one has
  \[ 1 = |\sigma(x-u)| = (x^\sigma - u^\sigma)(\ov{x^\sigma} -
         \ov{u^\sigma}) = 2 - x^\sigma \ov{u^\sigma} -
         \ov{x^\sigma} u^\sigma. \]
  For $y = \sigma(-x/u)$ one therefore has $y \ov{y} = 1$ and
  \[ y + \ov{y} = -x^\sigma / u^\sigma - \ov{x^\sigma} /
   \ov{u^\sigma} = -(x^\sigma \ov{u^\sigma} +
   \ov{x^\sigma} u^\sigma)/u^\sigma \ov{u^\sigma} = -1. \]
  Thereby $y$ is a root of $z^2 + z + 1 = 0$, i.e.\ $y$ is a
  third root of unity. In particular $y$ is integral, and this
  implies that also $x/u$ and $x$ are integral; this was to be shown.
\end{proof}

A real $c' \in \R$ is called a \textbf{usable bound}
\index[S]{usable bound}
(``usable bound'' in Lenstra) for $K$ if $c' \ge c(K)$
and for all $x \in F_K$ with $\cM_K(x)=c'$ there exists a
root of unity $u \in R$ with $\cM_K(x-u)=c'$. One
notes that every $c' > c(K)$ is usable.

\begin{quote}
  {\bf (1.23)} {\em If $c'=1$ is usable for $K$, then $K$ is
    norm-Euclidean.}
\end{quote}

\begin{proof}
  Let an $x \in K$ be given; we seek a $y \in R$ with
  $|N_{K/\Q}(x-y)| < 1$. Without loss of generality we may assume $x \in F_K$; then
  $\cM_K(x) \le 1$, hence by (1.19)
  $|N_{K/\Q}(x)| \le \cM_K(x) \le 1$. Thus $y=0$ suffices, unless
  precisely $|N_{K/\Q}(x)| = 1$. If this occurs, one must also have
  $\cM_K(x) = 1$, and since $c=1$ is a usable bound
  for $K$, there exists a root of unity $u \in R$ with
  $\cM_K(x) = \cM_K(x-u) = 1$. If $|N_{K/\Q}(x-u)| < 1$,
  we can choose $y=u$; otherwise we have
  $|N_{K/\Q}(x)| = |N_{K/\Q}(x-u)| = \cM_K(x) = \cM_K(x-u) = 1$.
  From equality in (1.19) it then follows that
  $|\sigma(x)| = |\sigma(x-u)| = 1$ for all embeddings $\sigma$, and
  by (1.22) one then has $x \in R$, so that here $y=x$ suffices.
\end{proof}
The central result is now

\begin{quote}
  {\bf (1.24)} {\em Let $\zeta_m$ be an $m$-th root of unity and
    $L = K(\zeta_m)$; then one has
    $$ c(L) \le (L:K)\cdot c(K). $$
    If moreover $c'$ is a usable bound for $K$, then
    $(L:K) \cdot c'$ is one for $L$.}
\end{quote}

\begin{proof}
  Let $x \in F_L$; we must show that $\cM_L(x) \le (L:K)c(K)$.
  Since of $x \in F_L$ we have $\cM_L(x-y) \ge  \cM_L(x)$ for all
  $y \in S$  for all $y \in S$ and thus for all $y \in R$.
  Thus (1.20) guarantees $\frac{1}{n}T_{L/K}(x) \in
  F_K$ for $n = (L:K)$. Now $|\sigma(\zeta_m)| = 1$ for every
  embedding of $L$ into $\C$ and every $j \in \N$; hence
  $\cM_L(x\zeta_m^j) = \cM_L(x)$, and as above
  we now obtain $\frac{1}{n} \cdot T_{L/K}(x\zeta_m^j) \in F_K$ for
  $j=1, \dots, m$. Therefore $\cM_K(T_{L/K}(x\zeta_m^j)) = n^2
  \cdot \cM_K(\frac{1}{n}T_{L/K}(x\zeta_m^j)) \le n^2 \cdot
  c(K)$. With (1.21) we obtain from this
  $$ mn \cdot \cM_L(x) = \sum_{j=1}^{m}
    \cM_K(T_{L/K}(x\zeta_m^j)) \le \sum_{j=1}^{m} n^2 \cdot c(K) =
    mn^2 \cdot c(K), $$
  so that finally we have $\cM_L(x) \le n \cdot c(K)$.
\end{proof}

Now let $c'$ be a usable bound for $K$ and $\cM_L(x) = nc'$ for some
$x \in F_L$. Because of $\cM_L(x) \le n \cdot c(K)$ and $c \ge c(K)$
one must have $c'=c(K)$, and from the above proof it follows that
$\cM_K(\frac{1}{n}T_{L/K}(x\zeta_m^j)) = c(K) = c'$ for all $j \in
\Z$. By choosing $j=0$ we see that $\alpha = \frac{1}{n}T_{L/K}(x)$ is
an element of $F_K$ with $\cM_K(\alpha) = c'$. Since $c'$ is a usable
bound for $K$, there exists a root of unity $u \in R$ with
$\cM_K(\alpha - u) = c'$. With $y=u$ we obtain from (1.20) $\cM_L(x-u)
= \cM_L(x) = nc'$, and this shows that $nc'$ is a usable bound for
$L$.  For $K=\Q$ one has $\cM_\Q(x) = |x|^2$, $F_\Q = [-0.5, +0.5]$
and $c(\Q) = \frac14$; in order to show that $c' = \frac14$ is a
usable bound for $\Q$, we must find for every $x \in F_\Q$ with
$\cM_\Q(x) = \frac14$ a root of unity $u \in \Z$ with $\cM_\Q(x-u) =
\frac14$. Since however $c'$ is attained on $F_\Q$ only at the points
$0.5$ and $-0.5$, $u=1$ suffices.  If therefore $\zeta_m$ is a
primitive $m$-th root of unity and $L=\Q(\zeta_m)$, then from (1.24)
it follows that $c(L) \le (L:\Q)/4$ and $(L:\Q)/4$ is a usable bound.

Thus all cyclotomic fields $L=\Q(\zeta_m)$ with
$(L:\Q) \le 4$ are norm-Euclidean: these are $\Q(\zeta_m)$
for $m = 3, 4, 5, 8, 12$. For cyclotomic fields
$K=\Q(\zeta_p)$, $p$ an odd prime, Lenstra\index[N]{Lenstra} \cite{Len74} has determined the
constants $c(K)$ in an elementary way:

\begin{quote}
  {\bf (1.25)} {\em Let $K = \Q(\zeta_p)$, $p \equiv 1 \bmod {2}$
    prime; then $c(K) = \frac{p+1}{12}$, and $c'=c(K)$ is a
    usable bound for $K$.}
\end{quote}

From this it follows immediately that $\Q(\zeta_7)$ and
$\Q(\zeta_{11})$ are likewise norm-Euclidean; one obtains
however further for
\begin{itemize}
\item $L = \Q(\zeta_9) = K(\zeta_3)$, $K = \Q(\zeta_3)$,
  $c(K) = \frac13$, hence $c(L) \le (L:K) \cdot c(K) = 1$;
\item $L = \Q(\zeta_{15}) = K(\zeta_5)$ for $K = \Q(\zeta_5)$,
  $c(K) = 1/2$, $c(L) \le 2 \cdot 1/2 = 1$;
\item $L = \Q(\zeta_{20}) = K(\zeta_4)$ for $K = \Q(\zeta_5)$,
  $c(K) = 1/2$, $c(L) \le 2 \cdot 1/2 = 1$.
\end{itemize}
Thus also $\Q(\zeta_9)$, $\Q(\zeta_{15})$ and $\Q(\zeta_{20})$
are norm-Euclidean.

If $K$ is an imaginary quadratic number field, then $\cM_K(x) =
N_{K/\Q}(x)$ and thereby $c(K) = M(K)$; (0.18) shows that the
only imaginary quadratic fields with $c(K) \le 1/2$ are the two
cyclotomic fields $\Q(\zeta_3)$ and $\Q(\zeta_4)$.
If on the other hand $K = \Q(\sqrt{m})$ is real quadratic, then for
$x=a+b\sqrt{m}$ one has
\[
\cM_K(x) = \frac{(a+b\sqrt{m})^2 + (a-b\sqrt{m})^2}{2} = a^2 +
m b^2,
\]
and the proof of (0.21) shows

\begin{quote}
  {\bf (1.26)} {\em Let $K = \Q(\sqrt{m})$ be a real quadratic
    number field; then one has}
   \[ c(K) =
   \begin{cases}
     \frac{1+m}{4 } & \text{ {\em for }} m \equiv 2,3 \bmod {4} \\
     \frac{(1+m)^2}{16} & \text{ {\em for }} m \equiv 1 \bmod {4}.
   \end{cases}
   \]
\end{quote}

For $K = \Q(\sqrt{5})$ one therefore has $c(K) = 9/20$, and (1.24)
yields the norm-Euclidean number fields $\Q(\sqrt{5},
\sqrt{-1})$, $\Q(\sqrt{5}, \sqrt{-3})$ and
$\Q(\zeta_5)$.
In order to give a further application of this method, we consider a
quadratic extension of $K = \Q(\sqrt{-3})$: such an extension
can be written in the form $L = K(\sqrt{\mu})$ with square-free $\mu \in
\Z[\sqrt{-3}]$. We set $\mu = a+b\sqrt{-3}$, $\mu' =
a-b\sqrt{-3}$, $\mu\mu' = a^2+3b^2 = p$, and find for an $x =
r+s\rho +t\sqrt{\mu} +u\rho\sqrt{\mu} \in L$ ($\rho =
(-1+\sqrt{-3})/2$) after a little computation
\[
\cM_{L/\Q}(x) = r^2 - rs + s^2 + \sqrt{p}(t^2 - tu + u^2).
\]
In the case $\mu \equiv 1 \bmod {4}$ now $\beta = \frac{1+\sqrt{\mu}}2$ is
integral; by writing $x$ in the form $x = e+f\rho + g\beta + h\rho\beta$,
we obtain in the above notation $r=e+g/2$, $s=f+h/2$,
$t=g/2$, $u=h/2$, hence
\[
\cM_{L/\Q}(x) = (e+g/2)^2 - (e+g/2)(f+h/2) + (f+h/2)^2
+ \sqrt{p}(g^2 - gh + h^2)/4.
\]
Now we know however that $c(K) = \frac13$; consequently we can choose $g$ and
$h \bmod {\Z}$ so that $g^2 - gh + h^2 =
\cM_{K/\Q}(g + h\rho) \le \frac13$. Now we choose
$e+g/2$ and $f+h/2 \bmod \Z$ so that also $(e+g/2)^2 -
(e+g/2)(f+h/2) + (f+h/2)^2 \le \frac13$, and have thereby shown: For
all $x \in L$ we can find a $y \in S$ such that
\[
\cM_{L/\Q}(x-y) \le \frac13 + \frac{\sqrt{p}}{12}
              = \frac{4 + \sqrt{p}}{12}
\]
holds.

Thus  $c(L) \le (4 + \sqrt{p})/12$, and $L$ is
norm-Euclidean if only $\sqrt{p} \le 8$ (and of course
$\mu \equiv 1 \bmod {4}$). Thereby we obtain the following
norm-Euclidean fields: $K(\sqrt{\mu})$ for
$$ \mu = -1+2\sqrt{-3}, \ 3+2\sqrt{-3}, \ 5, \ -5+2\sqrt{-3}, \
         -7, \ -3+4\sqrt{-3}, \ 7-2\sqrt{-3}; $$
the discriminants of these fields are incidentally $\disc L = 117$, $189$,
$225$, $333$, $441$, $513$, $549$. That these fields are norm-Euclidean
was already proved by Lakein\index[N]{Lakein} \cite{Lak72},
though with a completely different method.

\begin{quote}
  {\bf (1.27)} {\em Let $L = \Q(\sqrt{a+b\sqrt{-3}})$, $p=a^2+3b^2$;
    then $c(L) \le \frac{4+\sqrt{p}}{12}$.}
\end{quote}

It may still be remarked that for the function $\cM_K$ an analogue
to (1.14) exists: for if we extend $\cM_K$ to $\R^n$ by
\[
\cM(x) = \frac{x_1^2 + \dots + x_r^2 + 2(x_{r+1}^2 +
  x_{r+2}^2) + \dots + 2(x_{n-1}^2 + x_n^2)}{n}
\]
for $x = (x_1, \dots, x_n) \in \R^n$ and set
$U = \{x \in \R^n : \cM(x) < \frac14\}$, then $\cM(u-v) < 1$ for all
$u,v \in U$. A simple coordinate transformation and
a classical volume computation (``volume of a hypersphere in
$\R^n$'') then shows
$$ \mathrm{vol}(U) = 2^{-s} \Big(\frac{\pi}4\Big)^{n/2}
      \Gamma\Big(1 + \frac{n}2\Big), $$
where $\Gamma$ denotes the gamma function. With a packing argument
(see on this Rogers\index[N]{Rogers} \cite{Rog64} and
Lenstra\index[N]{Lenstra} \cite{Len77a}) one then obtains

\begin{quote}
  {\bf (1.28)} {\em Let $K$ be a number field with $(K:\Q)=n$ and
    $\delta < 1$; then for all $x \in K$ there exists a $y \in R$ and
    a unit $u \in R^\times$ with $\cM_{K/\Q}(u x - y) < \delta$,
    where
    \[ \delta = \sigma_n \Gamma\Big(1 + \frac{n}2\Big)
           \Big(\frac{4}{n\pi}\Big)^{n/2}
           \frac{\sqrt{|\disc K|}}{\mu_1} \]
    holds.}
\end{quote}

Here the $\sigma_n$ are the quantities depending only on $n$ given
e.g.\ by Lenstra\index[N]{Lenstra} \cite{Len77a}, and $\mu_1$ gives
(as on p.\ 33) the length of a maximal 1-sequence in $K$. With (1.19)
it now follows immediately that

\begin{quote}
  {\bf (1.29)} {\em Under the hypotheses of (1.28) one has: for
    all $x \in K$ there exists a $y \in R$ with $|N_{K/\Q}(x-y)| <
    \delta^{n/2}$. In particular $K$ is then norm-Euclidean.}
\end{quote}

Like (1.14), this is a result that goes back to
Lenstra\index[N]{Lenstra} \cite{Len77a}; the following table shows
that (1.29) in the cases $n=2$, $s=1$; $4 \le n \le 7$, $s \ge 2$;
$n=8$, $s \ge 3$ yields better bounds than (1.15) and (1.16):
$$ \begin{array}{ccc}
  \toprule n & \sigma_n \Gamma(1 + \frac n2) / \pi^{n/2} &
               \sigma_n \Gamma(1 + \frac n2) (4/n\pi)^{n/2} \\
  \midrule
  1 & 0.5 & 1 \\
  2 & \sqrt{3}/6 & \sqrt{3}/3 \\
  3 & 0.18613 & 0.286566 \\
  4 & 0.13128 & 0.131280 \\
  5 & 0.09988 & 0.037175 \\
  6 & 0.08113 & 0.024039 \\
  7 & 0.06982 & 0.009848 \\
  8 & 0.06327 & 0.003955 \\ \bottomrule
\end{array} $$
If one replaces $\mu_1$ in (1.28) by $\mu_k$, one of course obtains
corresponding results for $k$-Euclidean rings. Since $\cM_K$ is
not multiplicative, one obtains from (1.28) no bounds for
$c(K)$; if on the other hand one proceeds differently and sets $F(K) = \{x \in K :
\cM(x) \le \cM(u x - y) \text{ for all } u \in R^\times
\text{ and all } y \in R\}$, as well as $c(K) = \sup \{ \cM(x) : x
\in F \}$, then (1.28) does yield estimates for $c(K)$, yet one has
no result like (1.20) or (1.24) available for the $c(K)$.
\section*{{\sc Remarks on} \S\ 1}
\addcontentsline{toc}{section}{{\sc Remarks on} \S\ 1}
                 
The existence of totally ramified prime ideals can also be used for
finding factors of the class number; the best-known result of this
form is probably the genus theory of Gauss, who (in our language) has
shown that the class number of a real quadratic number field with $t$
ramified prime ideals is divisible by $2^{t-2}$ respectively
$2^{t-1}$, according as the norm of the fundamental unit is $+1$ or
$-1$. For similar results in number fields of higher degree see
Ishida\index[N]{Ishida} \cite{Ish76}.  The use of totally ramified
ideals in the sense of (1.5) is found in outline already in
Behrbohm\index[N]{Behrbohm} \& Rédei\index[N]{Redei@R\'edei}
\cite[p.\ 198]{BR36}; their method was then modified by
Erdős\index[N]{Erdos@Erd\"os} and Ko\index[N]{Ko} \cite{EK38} and
further developed by Heilbronn\index[N]{Heilbronn}
\cite{Hei38,Hei50,Hei51}, Cioffari\index[N]{Cioffari} \cite{Cio79} and
Egami\index[N]{Egami} \cite{Ega79}\footnote{At the last writing of the
  work I noticed that for the estimation of $M(K)$ one can also use
  ideals that are not totally ramified; see on this \S\  4.}.

The consideration of the ideals $I = (u-1)$ for units $u$ for the
estimation of $M(K)$ was first proposed by Rédei\index[N]{Redei@R\'edei}
\cite{Red41b}; he was thereby able e.g.\ to show that
$\Q(\sqrt{61})$ and $\Q(\sqrt{109})$ are not norm-Euclidean.
If $0,1,\vartheta_2,\dots,\vartheta_m$ is a 1-sequence, then in any case
$\vartheta_i-1$ for $i \ge 3$ must be a unit. Units
$u$ with the property that also $u-1$ is a unit are called
\textbf{exceptional units}\index[S]{exceptional unit} (exceptional
units).

Exceptional units had already been investigated before Lenstra
recognized their significance for the existence of a Euclidean
algorithm. Julia Robinson conjectured that in an algebraic number
field there are only finitely many exceptional units; this conjecture
was then proved independently by S.~Lang\index[N]{Lang} \cite{Lan60},
S.\ Chowla\index[N]{Chowla} \cite{Cho61} and
T.~Nagell\index[N]{Nagell} \cite{Nag64}.  Nagell occupied himself in a
series of papers with the determination of all exceptional units in
number fields of small unit rank.

Finally Chudnovsky and Chudnovsky\index[N]{Chudnovsky}
\cite{Chu86} have shown that the existence of ``fast'' algorithms for the
multiplication of polynomials over $\Z$ depends on the existence of
number fields with many exceptional units.

\chapter*{\S\ 2 The Determination of Euclidean Minima}
\setcounter{chapter}{2}
\addcontentsline{toc}{chapter}{\S\ 2 The Determination of Euclidean Minima}
\markboth{Euclidean Rings}{\S\ 2 The Determination of Euclidean Minima}

At the beginning of \S\  1 we saw how to find lower bounds for the
Euclidean minimum $M(K)$; here we wish to present methods
that allow one to estimate $M(K)$ from above. If
both bounds coincide, we have determined $M(K)$.

To this end let $\{\alpha_1, \dots, \alpha_n\}$ be a $\Q$-basis of an
algebraic number field $K$ (in general we shall here choose an integral basis;
however it is sometimes advisable, for reasons of symmetry, to take e.g.\
$\{1,\sqrt{m}\}$ for quadratic number fields $\Q(\sqrt{m})$
even in the case $m \equiv 1 \bmod {4}$ -- the reasons for this
will become clear presently). Then we define by
$$ \varphi : K \to \R^n =: \uK\label{puK}:
   (r_1\alpha_1 + \dots + r_n\alpha_n) = (r_1, \dots, r_n) $$
an embedding of $K$ into $\uK = \R^n$. As in \S\  1, we shall
in general identify $\R$ and $\varphi(\R)$. By
\[
|\cdot|_i : \uK \to \R : x = (x_1, \dots, x_n) \mapsto
\Big| \sum_{j=1}^{n} x_j \tau_i(\alpha_j)\Big|^\delta,
\]
where $\delta = 1$ for real and $\delta = 2$ for non-real
embeddings of $K$ into $\C$, we have ``extended'' the valuations of $K$
to $\uK$ in the sense that for all $x \in
K$ we have $|x|_1 = |\varphi(x)|_1$; note however that the $|\cdot|_i$
on $\uK$ are of course no longer valuations\label{pBew2}, if only
because $\uK$ is not a field (we regard $\uK$ only
as an $\R$-vector space). Further it does not in general follow from $|x|_1 = 0$
that $x = 0$! By $N(x) = |x|_1 \cdots |x|_{r+s}$ we can
also extend the norm $|N_{K/\Q}|$ from $K$ to $\uK$ and have
$|N_{K/\Q}(x)| = N(\varphi(x))$ for all $x \in K$. Here too we must
remark that $N$ on $\uK$ is not a norm in the topological sense,
because $N(x)=0$ does not imply $x=0$.

\medskip\noindent
\textbf{Example:}
let $K=\Q(\sqrt{7})$; we choose the $\Q$-basis $\{1, \sqrt{7}\}$;
then $\varphi(a+b\sqrt{7}) = (a,b)$ for all $a,b \in \Q$. Thus
$\varphi(K) = \Q^2$, while $\uK = \R^2$. Further we have
$|x|_1 = |a+b\sqrt{7}|$, $|x|_2 = |a-b\sqrt{7}|$ and $N(x) =
|a^2-7b^2|$ for an $x = (a,b) \in \uK$. In particular $|x|_2 = 0$
and so also $N(x) = 0$ for $x = (\sqrt{7},1) \in \uK \setminus
\{0\}$. Note that $\{x \in \uK: |x|_1 = 0\}$ defines a line
in $\R^2$.
\medskip

By equipping $\uK$ with the ordinary Euclidean metric, the functions
${|\cdot|_1}$, $\ldots$, $|\cdot|_{r+s}$ and $N$ become continuous
functions on $\uK$, and for every $\eps \in \R^+$ the set
$U_\eps = \{x \in \uK: |x|_i < \eps \text{ for } 1 \leq i \leq r+s\}$ is a
bounded neighbourhood of $0$ (the boundedness is easily seen by using
(1.10)). If $\{\alpha_1, \dots, \alpha_n\}$ is an integral basis of
$K$, then
$$ F = \{(x_1, \dots, x_n) \in \uK : - \frac12 \leq x_i < \frac12
     \text{ for } i = 1, \dots, n\} $$
is called a \textbf{fundamental domain}\index[S]{fundamental domain} of
$K$. By $\uF$ we denote the compact closure
$$ \uF = \{(x_1, \dots, x_n) \in \uK: - \frac12 \leq x_i \leq \frac12
  \text{ for } i = 1, \dots, n\} $$
of $F$.

In order to show that for all $x \in K$ there exists a $y \in \cO_K$
with $|N_{K/\Q}|(x-y) < k$, it suffices to find for every $x \in \uK$
a $y \in \cO_K$ with $N(x-y) < k$. This justifies the following
definition:
\[
M(\uK) := \inf \{k \in \R: \forall x \in \uK \ \exists y
\in \cO_K: N(x-y) < k\};
\]
correspondingly we also define the further minima $M_2(\uK)$,
$M_3(\uK)$ etc., insofar as they exist. Because of $K \subset \uK$ the
inequality $M(K) \leq M(\uK)$ is immediately clear. Whether the
isolation of $M(K)$ also implies that of $M(\uK)$ is not known;
however, the isolation of $M(\uK)$ certainly does imply that of $M(K)$,
and we have $M_2(K) \leq M_2(\uK)$. Corresponding statements hold for
the higher minima.

If we set $M(F) = \inf \{k \in \R: \forall x \in F \ \exists y
\in \cO_K \text{ with } N(x-y) < k\}$, then evidently
$M(F) = M(\uF) = M(\uK)$. Since $N$ as a continuous function on the
compact set $\uF$ is bounded, we have
$M(\uF) \le \max \{N(x): x \in \uF\}$, and in particular
$M(\uK)$ (and hence also $M(K)$) is finite.
The method we now present for determining $M(K)$
goes back to Barnes and Swinnerton-Dyer and always proceeds
as follows: we show $M(K) \ge k$ for some $k \in \R$ using
the criteria from \S\  1, and also $M(\uK) \le k$ for the same $k$;
this gives $M(K) \ge M(\uK) \le k$, hence $M(K) = M(\uK) =
k$. This already shows that this method can only be
successful when $M(K) = M(\uK)$ holds; it is conjectured,
however, that this is always the case in number fields. We now
describe this method using the example $K=\Q(\sqrt{14})$.

To this end we choose the integral basis $\{1,\sqrt{14}\}$ and set
$(a,b) \times (c,d) = \{t+s\sqrt{14}: a \le t \le b, c \le s \le d\}$;
so that, for example, $F = (-0.5, 0.5) \times (-0.5, 0.5)$. In order to show that
$M(K) \ge \frac{3}{4}$, we choose a $k \in \R$ that is
somewhat smaller than $\frac{3}{4}$, e.g.\ $k=1.2$ (the exact value of
$k$ plays no role here; in order to determine $M(K)$ it should only
lie between $M_0(K)$ and $M(K)$). We then call a set
$(a,b) \times (c,d) \subset \uK$ covered,\index[S]{covered} if for
all $x \in (a,b) \times (c,d)$ there exists a $y \in \cO_K$ with
$N(x-y) < k$ (whether such a set is covered therefore depends on $k$).
An $x \in \uK$ with $N(x-y) \ge k$ for all $y \in \cO_K$
we call an \textbf{exceptional point}.\index[S]{exceptional point}

The programs that we shall describe more precisely in \S\  11 allow us
to cover the whole of $F$ except for the following regions:
\begin{align*}
  S_1 & = (0.499, 0.5) \times (0.499975, 0.5), \\
  S_2 & = (0.499, 0.5) \times (-0.5, -0.499975), \\
  S_3 & = (-0.5, -0.499) \times (0.499975, 0.5), \\
  S_4 & = (-0.5, -0.499) \times (-0.5, -0.499975).
\end{align*}
We now see that it is somewhat more economical to consider instead
of $F$ the set $F' = (0,1) \times (0,1)$, because then only $S =
(0.499, 0.501) \times (0.499975, 0.500025)$ remains uncovered and we
have to deal with only one set instead of four. Now let $x \in F'$ be
an exceptional point (if such a point exists; in fact we know from
\S\ 1 that $M(K) = \frac{3}{4}$ holds for $x = (1+\sqrt{14})$); then
certainly $x \in S$. With $x$, however, also $xu$ is an exceptional
point, where $u = 15+4\sqrt{14}$ is the fundamental unit of $K$,
i.e.\ there exists an $a \in \R$ with $xu-a \in S$. This consideration
leads us to consider the set $uS = \{ux: x \in S\}$; a simple
computation shows
\[
uS = (35.4836, 35.5164) \times (9.495625, 9.504375), \text{ hence}
\]
\[
uS - (35 + 9\sqrt{14}) = (0.4836, 0.5164) \times (0.495625, 0.504375).
\]
With $a = 35 + 9\sqrt{14}$ we therefore have: for every $x \in S$ the
point $ux-a$ lies either in covered region or again in $S$ (had we
computed less precisely and e.g.\ been able to cover $F'$ only up to
the set $S' = (0.47, 0.53) \times (0.47, 0.53)$, then $uS' = (33.37,
37.63) \times (8.93, 10.07)$ would have resulted and there would exist
no $a \in \R$ with the above property, because e.g.\ $uS' - (35 +
9\sqrt{14}) = (-1.63, 1.63) \times (-0.07, 1.07)$ contains, besides
points from $S'$, also such from $S'-1$ or $S'+1$).  If however we
have a situation as above, we can describe the location of any
exceptional points rather precisely; in fact we have

\begin{quote}
  {\bf (2.1)} {\em Let $K$ be a number field, $u$ a unit in $R$, $S$
    a subset of $F$ and $a \in R$ a point with the
    property that for all $x \in S$ the point $ux-a$ lies in covered region
    or again in $S$. Then for every exceptional point $x_0$
    and for every valuation $| \cdot |_1$ with $1 < |u|_1$ we have $|x_0 -
    z|_1 = 0$; here $z = \frac{a}{u-1}$ is the fixed point of the map
    $x \mapsto ux-a$. If moreover $|u|_1 \neq 1$ for all
    valuations $| \cdot |_1$ of $K$, then by $x_{i+1} = ux_i -
    a$ a sequence $x_0, x_1, x_2, \dots$ of exceptional points is defined
    with $x_i \in S$ and $\lim x_i = z$.}
\end{quote}

\begin{figure}[ht!]
\centering
\begin{tikzpicture}[scale=1.5]
\draw[thin, gray!50] (-2.2,-3) grid [xstep=1, ystep=2] (3.5,2.8);
\draw[thick, ->] (-2.6,0) -- (4,0);
\draw[thick, ->] (0,-2.5) -- (0,3.2);
\draw[line width=1pt] (0.8,1.8) rectangle (1.2,2.2);
\draw[dashed] (0.6,1.6) rectangle (1.4,2.4);
\draw[line width=1pt] (-1,1.8) rectangle (-0.8,2);
\draw (-0.8,1.6) node {$S_3$};
\draw[line width=1pt] (-1,-1.8) rectangle (-0.8,-2);
\draw (-0.8,-1.6) node {$S_4$};
\draw[line width=1pt] ( 1,-1.8) rectangle (1.2,-2);
\draw (1.2,-1.6) node {$S_2$};
\draw (1,2) -- (1.5,3);
\draw (1.5,3.15) node {$z = \frac{1+\sqrt{14}}2$};
\draw (0.1,0.2) node {$0$};
\draw (1.1,0.2) node {$\frac12$};
\draw (2.1,0.2) node {$1$};
\begin{scope}
  \clip (-2.7,-0.3) rectangle (4,3);
  \draw (-3,-2/7) -- ( 4,26/7);
  \draw ( 5,-2/7) -- (-2,26/7);
  \fill (1,2) circle (1pt);
\end{scope}
\draw (-1.5,0.75) node[rotate= 29] {$|x-z|_1 = 0$};
\draw ( 3.5,0.75) node[rotate=-29] {$|x-z|_2 = 0$};
\end{tikzpicture}
\caption{Exceptional sets in $\Q(\sqrt{14}\,)$}
\end{figure}
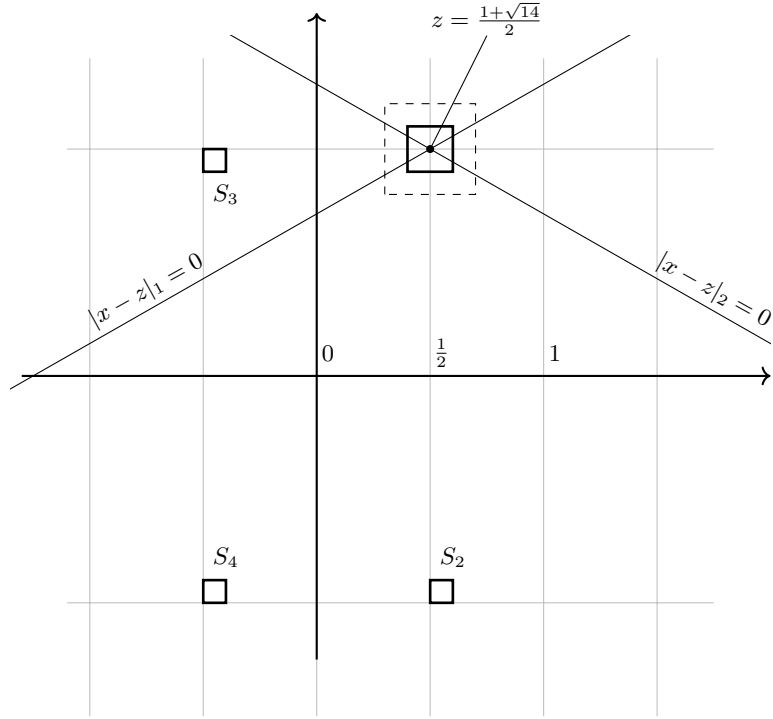

In our example $K = \Q(\sqrt{14})$ we have $u = 15 + 4\sqrt{14}$, $a =
35 + 9\sqrt{14}$, so that $z = a/(u-1) = \frac{35 + 9\sqrt{14}}{14
  + 4\sqrt{14}} = 1 + \sqrt{14}$; the set of all $x \in K$ with
$|x-z|_1 = 0$ is here a line through the point $z$. As the only
exceptional point (for $k=1.2$) we find here $x_0 = z$; from this
$x_1 = ux_0 - a = x_2 = z$, so that the sequence $x_0, x_1,
x_2, \dots$ converges trivially to $z$.

In the case of real quadratic fields (2.1) is due to Barnes and
Swinnerton-Dyer \cite[Theorem C]{BS52a}; the proof that they give for their
Theorem C goes back to Cassels (unpublished; see on this Bambah
\cite{Bam51}) and is much more opaque than the following:

\begin{proof}[Proof of (2.1)]
  Let $| \cdot |_1$ be a valuation of $K$ with $1 < |u|_1$ and $x_0
  \in F$ an exceptional point. Then also $x_1 = ux_0 - a$ must again
  be an exceptional point (otherwise with $x_1$ also $x_0 = (x_1 -
  a)u^{-1}$ would lie in covered region), and by hypothesis we then
  have $x_1 \in S$. By the recursion $x_{i+1} = ux_i - a$ we thus
  obtain a sequence $x_0, x_1, x_2, \dots$ of exceptional points lying
  in $S$. Now however $x_{i+1} - z = u(x_i - z)$ for all $i \in \N_0$;
  since all $x_i$ lie in $S$, $|x_{i+1} - z|_1$ must be bounded for
  all valuations $| \cdot |_1$ of $K$. In the case $|u|_1 > 1$ this is
  however only possible if $|x_i - z|_1 = 0$ for all $i \in \N_0$.  If
  on the other hand $|u|_1 < 1$, then $|x_i - z|_1$ is a null
  sequence; if therefore $|u|_1 \neq 1$ for all valuations
  $| \cdot |_1$ of $K$, we have $\lim |x_i - z|_1 = 0$ for all $j$
  with $1 \le j \le r_0$; this, however, is only possible if $\lim x_i = z$.
\end{proof}

In $K=\Q(\sqrt{14})$ every exceptional point of $S$ therefore lies
on the line $|x-z|_1 = 0$ for that valuation $|\,\cdot\,|_1$ for
which $|15 + 4\sqrt{14}|_1 > 1$. In order to proceed further, we consider
instead of $u$ the unit $u' = u^{-1} = 15 - 4\sqrt{14}$ and find as above
$u'S = -25 + 5\sqrt{14} + (0.4836, 0.5164) \times (0.495625,0.504375)$.
Applying (2.1) to this, it follows because of $1<|u'|_2$
that every exceptional point of $S$ lies on the line $|x-z|_2=0$ through the
point $z = \frac{-21 + 5\sqrt{14}}{14-4\sqrt{14}} = \frac{1+\sqrt{14}}2$.
Since the two lines $|x-z|_1=0$ and $|x-z|_2=0$ intersect precisely at the
point $z$, $z = \frac{1+\sqrt{14}}2$ is the only
exceptional point of $S$. Already in \S\  1 we have seen that
$M(z)=5/4$; we therefore have

\begin{quote}
  {\bf (2.2)} {\em Let $K=\Q(\sqrt{14})$; then we have
  $M(K) = M(\uK) = 5/4$; this Euclidean minimum is
  attained on $F = [0,1] \times [0,1]$ precisely at the point
  $z = (1 + \sqrt{14})/2$; for all $x \in F \setminus \{z\}$
  we have $M(K,x) \le 1.2$, i.e.\ the first Euclidean minimum is isolated.}
\end{quote}

In many cases we will -- as above for $\Q(\sqrt{14})$ --
be able to use both $u$ and $u^{-1}$ for the determination of $M(K)$;
in this case (2.1) takes the following form:

\begin{quote}
{\bf (2.3)} {\em With the designations of (2.1) suppose that:
  \begin{enumerate}
  \item there exists an $a \in \R$ such that $ux-a$ for all $x \in S$
    lies in covered region or again in $S$;
  \item for all $x \in S$ there exists a $b \in \R$ such that
    $u^{-1}x-b$ lies in covered region or again in $S$;
  \item we have $|u|_j \ne 1$ for every valuation $||_j$ of $K$;
  \end{enumerate}
  then $z = \frac{a}{u-1}$ is the only possible exceptional point of $S$. }
\end{quote}

\medskip\noindent
\textbf{Remark.} For real quadratic $K$ this is due to Cassels and was
first published by Bambah \cite{Bam51}; for the original proof the
same remark holds as for that of 2.1. As we easily see, in fields of
unit rank 1 the condition 3.\ can be replaced by the requirement ``$u$
is not a root of unity''.

\begin{proof}
  Let $x \in S$ be an exceptional point; for those valuations with
  $|u|_j > 1$ we have $|x-z|_j = 0$ by (2.4) and condition 1. For
  all other valuations we have $|u|_j < 1$ by 3., so that then $1 <
  |u^{-1}|_j$. We now observe that with $x$ also $x' =
  u^{-1}x - b \in S$ must be an exceptional point; because of $x' \in S$ and
  1.\ we have $ux' - a \in S$, hence $x - ub - a \in S$. Since $S$ is a
  subset of $F$ and $a, ub \in \R$, we must have $a = -ub$: thus
  $b = -au^{-1}$ is independent of $x$. The argument of (2.4) now shows
  that $x$ satisfies the equations $|x-z'|_j = 0$, where
  $$ z' = \frac{b}{u^{-1}-1} = \frac{bu}{1-u} = \frac{a}{u-1} = z $$
  holds. Therefore $|x-z| = 0$ for all $j$ with $1 \le j \le r+s$, and
  consequently $x-z$ is the only possible exceptional point.
\end{proof}

With a little more effort we can also determine the second Euclidean
minimum of $\Q(\sqrt{14})$ (Bedocchi \cite{Bed85} showed in 1985 that
$M_2(K) \le 1$; van der Linden asserts on \cite[p.\ 25]{Lin85} that
$M_2(K) < 1$, but gives no proof for this). To this end we choose
$k=0.96$ and can cover the whole of $F_1 = (-0.1,0.9)\times
(-0.1,0.9)$ except for the following regions:
\begin{align*}
   S^+ & = (-0.0005, 0.0005) \times (0.37475, 0.37525), \\
   S^- & = (-0.0005, 0.0005) \times (-0.37525, -0.37475), \\
   S & = (0.499, 0.501) \times (0.499975, 0.500025);
\end{align*}
with $a = 35 + 9\sqrt{14}$ and $b = -21 + 5\sqrt{14}$ the sets $uS-a$
and $u^{-1}S-b$ lie in covered region or again in $S$ (in particular
both sets have no point in common with $S^+$ and $S^-$); by (2.3)
therefore $z = (1+\sqrt{14})/2$ is the only exceptional point of
$S$. Now we find
\[
uS^+ = 21 + 6\sqrt{14} + (-0.0215,0.0215) \times (-0.38075, -0.36925)
\text{ and}
\]
\[
uS^- = -21 - 6\sqrt{14} + (-0.0215,0.0215) \times (0.36925, 0.38075).
\]
We thus see: for $a=21+6\sqrt{14}$ and $x\in S^+$ the point $ux-a$ lies
either in covered region or in $S^-$, while for $x\in S^-$ the
points $ux-a$ lie in covered region or in $S^+$. This leads us
to the following generalization of (2.1):

\begin{quote}
  {\bf (2.4)} {\em Let $K$ be a number field, $u$ a unit, $S_1, \dots, S_t$
    subsets of $F$ and $a_1, \dots, a_t \in R$ points with the
    property: for all $x\in S_1$ the point $ux-a_1$ lies in covered
    region or in $S_{i+1}$ (here we set $S_{t+1}=S_1$). If
    then $x_1\in S_1$ is an exceptional point, we have $|x_1-z|_i = 0$ for
    every valuation of $K$ with $|u_i| > 1$, where $z = a/(u^t-1)$ with
    $a = u^{t-1}a_1 + u^{t-2}a_2 + \dots + u a_{t-1} + a_t$ is the fixed point
    of the map $x \to ux - a$. }
\end{quote}

\begin{proof}
  Let $x_1\in S_1$ be an exceptional point; then also $x_2 = ux_1 -
  a_1 \in S_2$ is such a point, and we find successively $x_3 = ux_2 -
  a_2 \in S_3, \dots, x_{t+1} = ux_t - a_t \in S_{t+1} = S_1$. By
  computing backwards from here we obtain $x_{t+1} = u^t x_1 - a$ for
  $a = u^{t-1}a_1 + u^{t-2}a_2 + \dots + u a_{t-1} + a_t$; in this way we
  have found an $a \in R$ with the property that for every exceptional
  point $x_1 \in S_1$ also $u^t x_1 - a \in S_1$ holds for $u^t =
  u$. With (2.1) it now follows that every exceptional point $x \in
  S_1$ satisfies the equation $|x-z|_i = 0$, where $z = a/(u^t-1)$ is
  the fixed point of the map $x \to u^t x - a$ and $|u^t|_i$ is a
  valuation with $1 < |u^t|_i$. The claim follows.
\end{proof}

This result too is due to Barnes and Swinnerton-Dyer in the real
quadratic case, and the corresponding generalization of (2.3) to
Cassels. In $\Q(\sqrt{14})$ we have $a_1 = 21 + 6\sqrt{14}$, $a_2 =
-a_1$, hence $a = u a_1 - a_1 = a_1(u-1)$ and consequently $z = a/(u_2-1) =
a_1/(u_1-1) = 3\sqrt{14}/8 = (0, 3/8)$ (from now on we shall write
simply $(a,b)$ instead of $a+b\sqrt{m}$). Using $u^{-1}$ we
easily find that this $z$ is the only possible exceptional point of
$S_1$; from this $z = u z - a_1 = -z$ is the only possible exceptional
point of $S_2$. Now we use (1.8) with $t=2$, $x_1=z$, $x_2 = z$ and
find $M(K,z) = M(K,z') = 31/32 = 0.96875$; this gives the
following sharpening of (2.2):

\begin{quote}
  {\bf (2.2')} {\em Let $K = \Q(\sqrt{14})$: then
    $M(K) = M_1(K) = \frac54$ is the first Euclidean minimum of $K$;
    this is isolated and is attained $\bmod\ R$ precisely at the point
    $(\frac12, \frac12)$. Further
    $M_2(K) = M_2(\uK) = \frac{31}{32}$ is the second
    Euclidean minimum of $K$; this too is isolated and is attained
    $\bmod\ R$ precisely at the points $(0, \frac38)$. }
\end{quote}

The Euclidean minima cannot always be determined as simply as in
$\Q(\sqrt{14})$: a somewhat more complicated example is $K =
\Q(\sqrt{13})$. Here we choose the $\Q$-basis $\{1, \sqrt{13}\}$ and
the ``fundamental domain'' $F = (0, \frac12) \times (0, \frac12)$ (we
may restrict ourselves to this, because for every $x \in \uK$ there
exists a $y \in R$ with $x - y \in F \cup -F \cup F^\sigma \cup
-F^\sigma$; here $\sigma$ stands for the non-trivial automorphism
$\sigma: \sqrt{13} \mapsto -\sqrt{13}$). With $k=0.333$ we can cover
the whole of $F$ except for the following regions
\begin{align*}
  S_1 & = (-0.001, 0.001) \times (0.2128, 0.2129) &
  S_4 & = (0.161, 0.163) \times (0.1678, 0.1682) \\
  S_2 & = (0.116, 0.117) \times (0.1806, 0.1808) &
  S_5 & = (0.165, 0.166) \times (0.1668, 0.1672) \\
  S_3 & = (0.151, 0.152) \times (0.1707, 0.1711) &
  S_6 & = (0.166, 0.167) \times (0.1665, 0.1669),
\end{align*}
as well as the regions arising from these by multiplication with $-1$
and application of the automorphism $\sigma$. With $u =
\frac{3+\sqrt{13}}{2}$ we now find $-uS_1 = -\frac{3+\sqrt{13}}{2}
\times (0.11465, 0.1168) \times (0.18015, 0.1808)$, i.e.\ $-uS_1 + u$
lies in covered region or in $S_2$; further $-uS_2 + u$ lies in
covered region or in $S_3$, $-uS_3 + u$ in covered region or in $S_4$
etc., and $-uS_6 + u$ in covered region or in $S_6$. With $S = S_1
\cup S_2 \cup \dots \cup S_6$ we therefore have: $-uS + u$ lies in
covered region or again in $S$. By (2.1) every exceptional point $x$
of $S$ lies on the line $|x-z|_1 = 0$ through the point $z =
\frac{u}{u+1} = (\frac16, \frac16)$.

Now we use $u' = \frac{3-\sqrt{13}}{2}$ and set $T_1=S_1$, as well as
\begin{align*}
  T_2 & = (-0.117, -0.116) \times (0.1806, 0.1808) \\
  T_6 & = (-0.167, -0.166) \times (0.1665, 0.1669)
\end{align*}
and $T = T_1 \cup T_2 \cup \dots \cup T_6$. As above we now find
that the points $-u'T - u'$ lie in covered region or again in $T$.
Thus every exceptional point $x$ of $T$ lies on the line
$|x-z'|_2 = 0$ through the point $z' = -\frac{u'}{u'+1} = (- \frac16, \frac16)$.

Now $T \cap S = S_1 = T_1$, so that every exceptional point of $S_1$
must lie on the lines $|x-z|_1 = 0$ and $|x-z'|_2 = 0$. These
two lines intersect at the point $x_1 = (0, \frac16 + \frac16\sqrt{13})
\in K \setminus \Q$, so that $x_1$ is the only possible
exceptional point of $S_1$. Further we see immediately that
$x_2 = -u x_1 +u$ is the only possible exceptional point of $S_2$ (for
$uS_2$ lies $\bmod R$ in covered region or in $S_1$). The
recursion $x_{i+1} = u x_i - a$ for $u = -\frac{3+\sqrt{13}}{2}$
and $a=-u$ now yields the sequence of possible exceptional points
\[
x_{k+1} = \Bigl(\frac{1}{6} - \frac{1}{6} e^k,
               \frac{1}{6} + \frac{e^k}{6\sqrt{13}}\Bigr)
\]
for $k = 0, 1, 2, \dots$ and $e = (-3+\sqrt{13})/2$.
In order to show that this sequence gives all possible exceptional points of $S$,
we assume that $x$ is such a point. Since $x$ lies on the
line $|x-z|_2 = 0$, there exists a $k \in \N$ with
$|x_{k+1} - z|_2 < |x - z|_2 < |x_k - z|_2$ (i.e.\ $x$ lies
``between'' two terms of the sequence $\{x_i\}$). Then however
$(-u)^{-k-1} x \bmod R$ lies between $x_1$ and $x_2$, which is evidently not
the case.

Finally we still show that $M(x_1) = \frac{1}{3}$ and that this
minimum is not attained. To this end we note first that $M(x_1) =
M(x_2) = \dots = M(x_k) = \dots$, so that from $\lim N(x_i) = N(\lim
x_i) = N(z) = |N_{K/\Q}(z)| = \frac{1}{3}$ it certainly follows that $M(x_k) \ge
\frac{1}{3}$. If now there existed a $y \in R$ with $N(x_k - y) <
\frac{1}{3}$, then from (1.8) the existence of a $z = x_k \bmod R$
with $z = r + s\sqrt{m}$, $|s| < 0.18957$ and $N(z) <
\frac{1}{3}$ would follow. We easily see, however, that no such $z$ exists.
Thus $M(x_k) = \frac{1}{3}$; since however $N(x_k - y)$ for
all $y \in R$ is irrational, we cannot have $N(x_k - y) = \frac{1}{3}$,
and this shows that $M(x_k)$ is not attained.

We now give the numerical values for the points $x_1, \dots, x_7$, in
order to enable a comparison between the exceptional points and the
sets $S_1, \dots, S_6$:
\begin{align*}
x_1 &= (0, 0.212891683\dots) \\
x_2 &= (0.116204060\dots, 0.180662475\dots) \\
x_3 &= (0.151387818\dots, 0.170904256\dots) \\
x_4 &= (0.162040603\dots, 0.167949705\dots) \\
x_5 &= (0.165266007\dots, 0.167055139\dots) \\
x_6 &= (0.166242581\dots, 0.166784286\dots) \\
x_7 &= (0.166538263\dots, 0.166702279\dots) \quad \text{etc.}
\end{align*}
Still greater difficulties are caused by the determination of a
Euclidean minimum when this is not isolated. That there exist
real quadratic number fields whose second minimum $M_2(K)$
is already no longer isolated was conjectured by Barnes and Swinnerton-Dyer;
Godwin then proved this conjecture in 1963 for $K = \Q(\sqrt{23})$.

It is probably possible with the methods presented here to show that
$M_2(K)$ for $K = \Q(\sqrt{m})$, $m=(2n+1)^2-2$, $n \ge 2$, has a
non-isolated second minimum. Barnes and Swinnerton-Dyer have shown
that in this case we have: $M_1(K) = (8n^3 + 6n^2 - 6n + 1)/2m$; this
minimum is isolated and is attained $\bmod R$ precisely at the points
$\bigl(\frac{1}{2}, \frac{1}{2} - \frac{1}{m}\bigr)$.

Correspondingly Godwin's result can be generalized by
showing that
$$M _2(K) \ge \frac{2n(2n+1)\sqrt{m} - (2n^2 + m)}{2m} $$
holds, and that in the case of equality this minimum is not isolated
(in the case $m=7$ we obtain incidentally $M_3(K) \ge
\frac{6\sqrt{7}-9}{14}$).

In the table below we give the Euclidean minima $M_1(K)$
and $M_2(K)$ for $K = \Q(\sqrt{m})$, $m \le 102$, insofar as
they are known. The first table of this kind is found in
Heinhold \cite{Hei39}; this was then supplemented by Barnes and Swinnerton-Dyer.
The first minima still missing there were finally determined by
Godwin \cite{God55}.

The second minima that were not yet computed in these three papers
are set in italics in the table. The smallest $m$ for which
$M_2(K)$ is still unknown are now $m = 19$ and $m=22$ for $m =
2,3 \bmod 4$, as well as $m = 41$ for $m = 1 \bmod 4$; in these cases
lower bounds for $M_2(K)$ are given in the table.
It may still be remarked that the tables of Barnes and Swinnerton-Dyer
contain some minor misprints; thus the first minimum of
$\Q(\sqrt{19})$ is incorrectly given as $\frac38$; this error
Barnes and Swinnerton-Dyer themselves corrected in the fourth part of their work.
Moreover then ``$C_1 = (\frac12, \frac38)$''
must be changed to ``$C_1 = (0, 5/9)$''.
Further the sets $C_2$ at which $M_2$ is attained are not correct for the
values $m = 17, 35, 37, 65, 101$. For the values
$m = 17$, $37$, $65$, $101$ in each case a sign at the second
point in $C_2$ is to be changed, while for $m = 15$ (respectively $m = 35$)
``$C_2 = (0, \frac12)$'' instead of
``$C_2 = (\frac12, \frac18)$'' (respectively
``$C_2 = (0, \frac12)$'' instead of
``$C_2 = (\frac12, \frac12)$'') must stand.
\begin{table}[h]
\centering
\begin{tabular}{llll}
\toprule
$m$ & $M_1$ & $M_2$ & Rem. \\
\midrule
5 & $1/4$ & $1/5$ & see Davenport 1946 \\
13 & $1/3$ & $4/13$ & \\
17 & $1/2$ & $8/17$ & \\
21 & $5/7$ & $12/23$ & \\
29 & $4/5$ & $23/29$ & \\
33 & $29/44$ & $6/11$ & \\
37 & $3/4$ & $27/37$ & \\
41 & $23/32$ & $\ge 2\,2819/4100$ & \\
53 & $9/7$ & $68/53$ & \\
57 & $14/19$ & $219/304$ & $M_1$: Godwin 1955 \\
61 & $1611/1525$ & $41/39$ & \\
65 & $1$ & $64/65$ & \\
69 & $25/23$ & & \\
73 & $1541/2136$ & irrational; see Godwin 1955 & \\
77 & $19/11$ & $16/11$ & \\
85 & $16/9$ & $151/85$ & \\
89 & $1004287/1000004$ & & \\
93 & $44/31$ & & \\
97 & $33679354/31404817$ & & \\
101 & $5/4$ & $125/101$ & \\
\bottomrule
\end{tabular}
\caption{Euclidean minima of real quadratic number fields
  ($m \equiv 1 \bmod 4$)}\label{App17}
\end{table}
\begin{table}[h]
\centering
\begin{tabular}{llll}
\toprule
$m$ & $M_1$ & $M_2$ & Remarks \\
\midrule
2 & $1/2$ & $1/4$ & see Varnavides \cite{Var48b} \\
3 & $1/2$ & $1/3$ & \\
6 & $3/4$ & $1/2$ & $M_3(K) \ge (6\sqrt{7}-9)/14$ \\
7 & $9/14$ & $1/2$ & \\
10 & $3/2$ & $39/40$ & \\
11 & $19/22$ & $3125/3971$ & \\
14 & $5/4$ & $31/32$ & \\
15 & $3/2$ & $7/5$ & \\
19 & $170/171$ & $\ge 10579/12844$ & \\
22 & $27/22$ & $\ge 175903/155236$ & \\
23 & $77/46$ & $(20\sqrt{23}-31)/46$ & Godwin 1963 \\
26 & $5/2$ & $207/104$ & \\
30 & $3/2$ & $29/20$ & \\
31 & $45/31$ & $3775/3038$ & \\
34 & $9/4$ & $135/68$ & $M_3 = 137/72$ \\
35 & $5/2$ & $17/7$ & \\
38 & $11/4$& $173/72$ & \\
39 & $5/2$ & $107/48$ & \\
42 & $7/4$ & $41/24$ & \\
43 & $11829/6962$ & $5902/3483$ & Godwin 1955 \\
46 & $79877/48668$ & & \\
47 & $253/94$ & $\ge (42\sqrt{47}-65)/94$ & \\
\bottomrule
\end{tabular}
\caption{$m \equiv 2,3 \bmod {4}$}
\end{table}
\begin{table}[h]
\centering
\begin{tabular}{llll}
\toprule
$m$ & $M_1$ & $M_2$ & Remarks \\
\midrule
51 & $287/102$ & & \\
55 & $9/4$ & $351/176$ & \\
58 & $3/2$ & $27477/19604$ & \\
59 & $125/59$ & & $M_3 = 367/128$ \\
62 & $13/4$ & $367/124$ & \\
66 & $15/4$ & $431/128$ & \\
67 & $341/162$ & & Godwin 1955 \\
70 & $891/500$ & & Godwin 1955 \\
71 & $7393/3479$ & & Godwin 1955 \\
74 & $5/2$ & & \\
78 & $7/2$ & & \\
79 & $585/158$ & $\ge (72\sqrt{79}-111)/158$ & \\
82 & $9/2$ & $1311/328$ & \\
83 & $631/166$ & & Godwin 1955 \\
86 & $10030/5203$ & & Godwin 1955 \\
87 & $169/58$ & & \\
91 & $5/2$ & & \\
94 & $4708623/2143294$ & & Godwin 1955 \\
95 & $7/2$ & & \\
102 & $19/4$ & $869/200$ & \\
\bottomrule
\end{tabular}
\caption{$m \equiv 2,3 \bmod {4}$}
\end{table}
The estimate $M(K) \le \sqrt{d}/4$ for the Euclidean minima
$M(K)$ of real quadratic number fields is due to Minkowski
(see Coll.\ Works 1, 320--356 or also Hardy \& Wright, \S\  24) and is
best possible: the latter follows from

\begin{quote}
  {\bf (2.5)} {\em Let $n \in \N$ be odd and $m=n^2+1$
    square-free. Then for the real quadratic number field
    $K = \Q(\sqrt{m})$ we have $M(K) = \frac{n}{2}$, and this minimum
    is attained $\bmod R$ precisely at $z = \frac{\sqrt{m}}{2}$. }
\end{quote}

Because $\disc K = 4m$ in this case we have
$\frac{\sqrt{d}}4 = \frac{\sqrt{n^2+1}}2$, i.e.\ we can make
$\frac{\sqrt{d}}4 - M(K)$ arbitrarily small by choosing $n$
sufficiently large.
For the proof of (2.5) we choose the integral basis $\{1, \vartheta\}$ with
$\vartheta = n+\sqrt{m}$. Then for every $z \in F$ we must find an $a
\in R$ with $|N(z-a)| \le \frac{n}{2}$ and show that
equality can occur only for $z = \sqrt{m}/2$. We claim that
(with this skilful choice of $\vartheta$) even $a=0$ suffices. To this end
we set $z = x+y\vartheta$, regard the norm $N$ as a function of
$x$ and $y$ and ask ourselves where $N$ can attain a maximum or a
minimum on $F$. In the interior of $F$ this is not possible, for
because of
\[
N(x,y) = N(x+y\vartheta) = x^2 + 2nx - y^2,
\]
\[
\frac{\partial N}{\partial x} = 2x + 2ny = 0 \quad \Rightarrow \quad x
= -ny
\]
and
\[
\frac{\partial N}{\partial y} = 2nx - 2y = 0 \quad \Rightarrow \quad y
= nx
\]
$z = 0$ is the only point at which both partial derivatives
vanish (at $z$ neither a maximum nor a minimum is
attained: $z$ is a saddle point of $N$). Consequently the
minimum and maximum of $N$ are attained on the boundary of $F$. Now let
$y=y_0 - \frac{1}{2}$; then $N(x, \frac{1}{2}\vartheta) = x^2 + nx
- \frac{1}{4}$, $\frac{\partial N}{\partial x} = 2x+n$, and consequently
$N$ becomes extremal on the line $y=y_0$ only at $x_0 = -\frac{n}{2}$.
For $n > 2$ however $x_0$ does not lie in $F$, for $n=1$ we have
$x_0 = -\frac{1}{2}$, and we see (after similar computations for
$y_0 = -\frac{1}{2}$, respectively $x_0 = \frac{1}{2}$), that the maximum of
$|N|$ on $F$ is attained at most at the corner points $z_e = (\pm\frac{1}{2},
\pm\frac{1}{2})$. There however we have $N(z_e) = \frac{n}{2}$,
and the claim follows because $(1+\vartheta)/2 \equiv \sqrt{m}/2 \bmod R$.
It remains to show $M(K, x) = \frac{n}{2}$ for $x =
\sqrt{m}/2$. To this end we use (1.8) and assume that $M(K, x)
\le \frac{n}{2}$. Then there exists a $z = x \bmod R$ with $z =
r+s\sqrt{m}$, $|s| < \sqrt{k(n+1)/2m} < \frac{n}{2}$ and $|N(z)| <
k$. The conditions on $s$ allow only $s = \pm\frac{1}{2}$, and for
reasons of symmetry we may assume $s = \frac{1}{2}$. In this case
$N(r+s\sqrt{m}) = |r^2 - \frac{n}{4}|$ and this is evidently minimal
when $|r| \le \frac{n}{2}$. Substituting $r_0 = \frac{n-1}{2}$
and $r_1 = \frac{n+1}{2}$ yields $N(r_0 + s\sqrt{m}) = N(r_1 +
s\sqrt{m}) = \frac{n}{2}$, and that was to be shown.
If $n$ is even and $m=n^2+1$ square-free, then (2.5) no longer holds in this form,
because $\{1, \vartheta\}$ is no longer an integral basis. On the other hand
(2.5) remains correct also in this case if we replace the ring of integers
$D(m)$ by the ring $\Z[\sqrt{m}]$; in particular we have
$M(f)=1$ for $\Z[\sqrt{5}]$ and $f=|N|$ (see the examples on
p.\ \pageref{pBei}).
With a somewhat different method (2.5) can still be sharpened: for if
$z \in F$ is an exceptional point for
$$ k = \frac{n-1}{2} - \frac{1}{4m}, $$
then necessarily
\[
|N(x+y\delta)| = |x^2 + 2nxy - y^2| \ge k.
\]
Because of $|x|, |y| \le \frac{1}{2}$ we have however $|x^2 - y^2| \le
\frac{1}{4}$ and $|2nxy| \le n|y|$, hence
\[
|N(x+y\delta)| \le |x^2 - y^2| + |2nxy| \le n|y| + \frac{1}{4}.
\]
Both inequalities together yield $n|y| \ge k - \frac{1}{4}$ for
every exceptional point $z \in F$, hence
\[
|y| \le \frac{1}{2} - \delta, \quad \delta = \frac{3}{4n} + \frac{1}{4mn}.
\]
This is possible only for $-\frac{1}{2} \le y \le -\frac{1}{2} + \delta$ or
$\frac{1}{2} - \delta \le y \le \frac{1}{2}$.
By replacing $F = (-\frac{1}{2}, \frac{1}{2}) \times (-\frac{1}{2},
\frac{1}{2})$ (as in $\Q(\sqrt{14})$) by $F' = (0,1) \times
(0,1)$, we recognize that every exceptional point must lie in the strip
$|y - \frac{1}{2}| \le \delta$. Correspondingly it now follows from
\[
|N(z)| \le n|x| + \frac{1}{4},
\]
that for an exceptional point $z = x + y\delta \in F'$ also $|x -
\frac{1}{2}| \le \delta$ must hold. Thus every exceptional point $x
+ y\delta$ lies in the rectangle about the point $(1+\delta)/2$ defined by
\[
\Big| x - \frac{1}{2} \Big| \le \delta, \quad
\Big| y - \frac{1}{2} \Big| \le \delta
\]
at which, as we have seen above, $M_1(K)$ is attained. Since now $x$
and $y$ both lie near $\frac{1}{2}$, we can estimate $|x^2 - y^2|$
somewhat better from above than by $\frac{1}{4}$; for $x,y \in F$ we
have $|x - y| \le \delta$, $|x + y| \le 1$ or $|x - y| \le 1$, $|x +
y| \le \delta$. Hence
\[
|x^2 - y^2| = |x - y||x + y| \le \delta,
\]
and as above it now follows that $n|x| \ge k - \delta$ and $n|y| \ge k
- \delta$.  Further we see from this that every exceptional point $z
\in F'$ lies in the rectangle $S$ defined by
\[
\Big|x - \frac{1}{2}\Big| \le \delta_1, \quad
\Big|y - \frac{1}{2}\Big| \le \delta_1
\]
where
\[
\delta_1 = \frac{1}{2n} + \frac{1}{n^2}
\]
holds.
Even if we repeat this step infinitely often, we cannot
improve $\delta_1$ substantially any more. We shall therefore only use the
units $u = \vartheta$ and $u' = \vartheta' = 2n - \vartheta$ in order to
proceed further. We find
\[
z \cdot u = (x + y\vartheta)\vartheta = y \cdot (x + 2ny)\vartheta
\]
and
\[
-z \cdot u' = -(x + y\vartheta)\vartheta' = y - 2nx + x\vartheta.
\]
Since with $z$ also $z\vartheta$ must be an exceptional point, there exists an
$a = r + s\vartheta \in R$ with $z\vartheta - a \in S$, i.e.\ with
\[ \Big|x + 2ny - s - \frac{1}{2}\Big| \le \delta_1; \]
because of $|x - \frac{1}{2}| \le \delta_1$ we must have $|2ny - s| \le 2\delta_1$.
On the other hand it follows from $|y - \frac{1}{2}| \le \delta_1$ that we have
\[ |2ny - n| \le 2n\delta_1 = 1 + \frac{2}{n} \]
and consequently $s \in \{n-1, n, n+1\}$. Dividing
$|2ny - s| \le 2\delta_1$ by $2n$ yields
\[ \Big|y - \frac{1}{2n}\Big| \le \frac{\delta_1}{n} = \eps, \]
and a corresponding computation with $\delta'$ in place of $\delta$
shows that also $|x - \frac{1}{2n}| \le \eps$ holds, where again
$r \in \{n-1, n, n+1\}$.
By combining all possibilities we now obtain the following
nine possible exceptional sets:
\[ S_{i,k}: \quad \Big|x - \frac{n+i}{2n}\Big| \le \eps, \quad
                 \Big|y - \frac{n-k}{2n}\Big| \le \eps, \]
where $i$ and $k$ independently run through the numbers $-1, 0, +1$.
Of these nine sets however $S_{i,k}$ for
$(i,k) = (-1,-1), (-1,+1), (+1,-1), (+1,+1)$ contain no exceptional points
if only $n \ge 5$: for if e.g.\ $|x + \frac{n-1}{2n}| \le \eps$
and $|y - \frac{1}{2n}| \le \eps$, we find $|x^2 - y^2| \le
(n+1)/n^3$ and $2n|x|y| \le (n-1) \cdot 2\delta_1^2/2n$, and for $n
\ge 5$ it follows from this that $|N(z)| < k$. Thus five exceptional sets
remain: $S_{i,k}$ for
$(i,k) = (-1,0), (+1,0), (0,0), (0,-1), (0,+1)$, that is
$$ \begin{array}{lll} \smallskip
  S_1 := S_{-1,0} : & \Big|x - \frac{n-1}{2n}\Big| \le \eps, &
                           \Big|y - \frac{1}{2}\Big| \le \eps, \\ \smallskip
  S_2 := S_{+1,0} : & \Big|x - \frac{n+1}{2n}\Big| \le \eps, &
                           \Big|y - \frac{1}{2}\Big| \le \eps, \\ \smallskip
  S_3 := S_{0,0} : & \Big|x - \frac{1}{2}\Big| \le \eps, &
                           \Big|y - \frac{1}{2}\Big| \le \eps, \\ \smallskip
  S_4 := S_{0,-1} : & \Big|x - \frac{1}{2}\Big| \le \eps, &
                           \Big|y - \frac{n+1}{2n}\Big| \le \eps, \\
  S_5 := S_{0,+1} : & \Big|x - \frac{1}{2}\Big| \le \eps, &
                           \Big|y - \frac{n-1}{2n}\Big| \le \eps.
\end{array} $$
Now we see immediately that $-\vartheta' S_1 \bmod R$ lies in covered region
or in $S_4$ (for $x+y\vartheta \in S_1$ we have $s=x$ with
$r+s\vartheta = -(x+y\vartheta)\vartheta'$); this implies
\[ \Big|y-2nx-(n-1)-\frac{1}{2}\Big| \le \eps, \]
because of $|y-\frac{1}{2}|\le \eps$ it follows that
\[ |2nx-(n-1)|\le \frac{2\eps}{n^2} \]
and thus
\[ \Big|x-\frac{n-1}{2n}\Big| \le \frac{\delta}{n^2}. \]
Quite analogous considerations lead to
$|x-\frac{n-1}{2n}|\le \frac{\delta}{n^2}$ for $z\in S_2$, respectively to
$|y-\frac{n-1}{2n}|\le \frac{\delta}{n^2}$ for $z\in S_4, S_5$.
Now we consider $\vartheta S_1$: this region lies $\bmod R$
either in covered region or in one of the three sets $S_3, S_4,
S_5$. If $(x+y\vartheta) \bmod R$ lies in $S_3$, then we must have
\[ \Big|x+2ny-n-\frac{1}{2}\Big| \le \eps \]
; because of $|x-\frac{n-1}{2n}|\le \eps$ it follows that
$$ \Big|2ny-n-\frac{1}{2n}\Big|\le 2\eps, \quad \text{hence} \quad
   \Big|y-\frac{1}{2}-\frac{1}{4n^2}\Big |\le \frac{\eps}{n}. $$
Correspondingly the possibility that $(x+y\vartheta)\vartheta \bmod R$
lies in $S_4$ leads to the condition $|y-\frac{1}{2}|\le \eta$, and
$(x+y\vartheta)\vartheta \bmod R \in S_5$ finally yields
\[ \Big|y-\frac{1}{2}-\frac{1}{2n^2}\Big|\le \eta. \]
We now claim that the sets
$$ \Big|x-\frac{n-1}{2n}\Big|\le \eta, \quad
   \Big|y-\frac{1}{2}\Big|\le \eta \quad \text{ as well as }
   \Big|x-\frac{n-1}{2n}\Big|\le \eta, \quad
   \Big|y-\frac{1}{2}-\frac{1}{4n^2}\Big| \le \eta $$
can contain no exceptional point. If we have shown this,
it follows that $\vartheta S_1 - n\vartheta$ can only lie in covered region or
in $S_5$. Further $\vartheta S_3 - (n+1)\vartheta$ lies in
covered region or in $S_2$. Now $S_2 = 1+\vartheta - S_1$,
consequently $S_5 - 1 - (n+2)\vartheta$ lies in covered region or in
$-S_1$. Thus we can apply (2.4) with $t=4$, $a_1 = n\vartheta$,
$a_2 = 1+(n+2)\vartheta$, $a_3 = -a_1$, $a_4 = -a_2$, $u=\vartheta$
and find
\[ z_1 = \frac{u^3 a_1 + u^2 a_2 - u a_1 - a_2}{u^4 - 1}
       = \frac{u a_1 + a_2}{u^2 + 1} = \frac{m - n + (m + 1)\vartheta}{2m} \]
as the only possible exceptional point of $S_1$. Further now
\[ z_2 = \frac{m + n + (m - 1)\vartheta}{2m}, \quad
   z_4 = \frac{m - 1 + (m - n)\vartheta}{2m}, \quad
   z_5 = \frac{m + 1 + (m + n)\vartheta}{2m} \]
are the only exceptional points of $S_2, S_4$ and $S_5$.
It now follows immediately that $\vartheta S_3 \bmod R$ can only lie in covered region
or again in $S_3$, and (2.4) yields
\[
z_3 = \frac{1+\vartheta}{2}
\]
as the only possible exceptional point of $S_3$.
We now show that there exists no exceptional point $z=x+y\vartheta$ with
$|x-\frac{n-1}{2n}|\le \eta$ and $|y-\frac{1}{2}|\le \eta$. To this end
we assume without loss of generality $- \frac{1}{2} - \eta \le y \le -\frac{1}{2}$
and find
\begin{align*}
  N(z) & = x^2 + 2nxy - y^2 \le \Big(\frac{n-1}{2n} - \eta\Big)
          + n \Big(\frac{n-1}{2n} + \eta\Big) - \frac{1}{4} \\
     & = \frac{n-1}{2} - \frac{1}{2n} + \frac{1}{4n^2} +
         \eta \cdot \frac{n^3 - n^2 - n + \delta}{n^2}.
\end{align*}
For $n \ge 5$ this is certainly $< k = \frac{n-1}{2} \cdot
\frac{1}{4m}$, so that the above region then contains no exceptional point
(for ``large'' $n$ this follows almost without
computation, because $N(z) \le \frac{n-1}{2} - \frac{1}{4n} + O(n^{-2})$ and
$k = \frac{n-1}{2} \cdot O(n^{-2})$; here we have used the
familiar ``big-O notation''). In much the same way we check
that there also exists no exceptional point with $|x-\frac{n-1}{2n}|\le
\eta$ and $|y-\frac{1}{2}-\frac{1}{4n^2}|\le \eta$ (here we have
$N(z) \le \frac{n-1}{2} - \frac{1}{4n} + O(n^{-2})$).
Now we must still determine $M(z_i)$ for the possible exceptional points $z_i$.
That $M(z_3) = \frac{R}{2}$ we had already
seen. We now use (1.8) and to this end represent the points $z_i$ in
the form $r+s\sqrt{m}$; we find
$$ \begin{cases}
  z_1 & = \frac{m+1}{2m} \sqrt{m}, \quad
  z_2 = \frac{m-1)}{2m} \sqrt{m}, \\
  z_4 & = \frac12 + \frac{m-n}{2m} \sqrt{m}, \quad
  z_5 = \frac12 + \frac{m+n}{2m} \sqrt{m} \end{cases}
  \quad (n \equiv 1 \bmod 2), $$
$$ \begin{cases}
    z_1 & = \frac12 + \frac{m+1}{2m} \sqrt{m}, \quad
    z_2 = \frac12 + \frac{m-1}{2m} \sqrt{m}, \\
    z_4 & = \frac{m-n}{2m} \sqrt{m}, \quad
    z_5 = \frac{m+n)}{2m} \sqrt{m} \end{cases}
  \quad (n \equiv 0 \bmod 2). $$
(1.8) now yields because of $a=n$ the bound
$\mu_2 = \sqrt{\frac{k(n+1)}{2m}} < 0.5$. Now in the case
$n \equiv 1 \bmod 2$ we have
\[\min_{x \in \Z} |N(x + z_2)|
      = \Big|N\Big(\frac{n-1}{2} + z_2\Big)\Big|
      = \frac{n-1}{2} + \frac{1}{4m} \]
and
\[ \min_{x \in \Z} |N(x + z_4)|
   = \Big|N\Big(\frac{n-2}{2} + z_4\Big)\Big|
    = \frac{n-1}{2} - \frac{1}{4m}, \]
and a corresponding computation for the case $n \equiv 0 \bmod 2$ shows

\begin{quote}
  {\bf (2.5')} {\em Let $n \in \N$, $m=n^2+1$,
    $R \equiv \Z[\sqrt{m}]$ and $l$ the absolute value of the
    norm in $\Q(\sqrt{m})$. Then $M_1(f) = \frac n2$ and in the
    case $n \ge 2$ $M_2(f) = \frac{n-1}{2} - \frac{1}{4m}$.
    These minima are attained $\bmod R$ only at $C_1 = \{\sqrt{m}/2\}$ and
    \begin{align*}
      C_2 & = \left\{ \frac{m \pm 1}{2m} \sqrt{m},
                  \frac{1}{2}, \frac{(m \pm n)}{2m} \sqrt{m} \right\}
             \quad \text{for } n \equiv 1 \bmod 2, \text{ and} \\
      C_2 & = \left\{ \frac{1}{2} + \frac{m \pm 1}{2m} \sqrt{m},
                     \frac{m \pm n}{2m} \sqrt{m} \right\}
             \quad \text{for } n \equiv 0 \bmod 2
    \end{align*}
    . If moreover $n$ is odd and $m$ square-free, then
    $R = D[m]$, as well as $M_1(l) = M_1(K)$, $M_2(l) = M_2(K)$. }
\end{quote}

The first minima $M_1(l)$ of the rings $\Z[\sqrt{m}]$ were already
found by Heinhold \cite{Hei39}. Analogous theorems for $m=n^2+r$,
$r=1, 2, 4$, are due to Barnes and Swinnerton-Dyer and to
Varnavides. Both however used rather complicated lemmas for the proof
(e.g.\ the Theorems H, J, K of BSD), which in our proof have turned
out to be superfluous.

It is now not surprising that with the same methods we can also treat
the cases $r \in \{1, 2, 4, n\}$ (quadratic number fields
$\Q(\sqrt{m})$ with $m=n^2+r$, $r|4n$, $-n \le r \le n$ are called
fields of R-D type after C.\ Richaud and G.\ Degert, who explicitly
gave the fundamental units in these fields), even if we cannot always
determine the first two minima so simply (we have already expressed
the conjecture that in the case $m = n^2-2$, $n \ge 5$, $n \equiv 1
\bmod 2$, the second minimum is not isolated). Thus e.g.\ the
following table holds:

\section*{(2.6) Euclidean Minima for Real Quadratic Fields of R-D Type}

Let $R = \Z[\sqrt{m}]$ (for $I = 1$) and $R = \Z[\frac{1+\sqrt{m}}2]$
(for $I = 2$), and $f$ the absolute value of the norm. Then the first
Euclidean minimum is given by the following table:
$$ \begin{array}{llccc}
  \rsp m & n & I & M(K) & C(K) \\ \hline
  \rsp n^2+1 & n \ge 1 & 1 & \frac n2 & (0, \frac12) \\
  \rsp & n \equiv 0\ (2), n \ge 1 & 2 & \frac n8 & \\
  \hline
  \rsp n^2-1 & n \equiv 0\ (2), n \ge 2 & 1 & \frac{n-1}2 & (\frac12, \frac12)
             \\ \hline
  \rsp n^2+2 & n \equiv 0\ (2), n \ge 2 & 1 & \frac{2n-1}4 & \\
  \rsp & n \equiv 1 \bmod 2, n \ge 3 & 1 & \frac{2n^3-3n^2+6n-7}{4m}
             & \\ \hline
  \rsp n^2-2 & n \equiv 0\ (2), n \ge 4 & 1 & \frac{2n-1}4
             & (\frac12, \frac12) \\
  \rsp & n \equiv 1 (2), n \ge 3 & 1 & \frac{2n^3-3n^2-6n+9}{4m}
             & \\ \hline
  \rsp n^2+4 & n \equiv 1 (2), n \ge 3 & 2 & \frac{n^2 - 2n + 1}{4n}
             & \\ \hline
  \rsp n^2-4 & n \equiv 1 (2), n \ge 3 & 2 & \frac{n^2-5}{4n+8} & \\ \hline
  \rsp n^2+n & n \equiv 0\ (2), n \ge 2 & 1 & \frac{n+1}4 & (0, \frac12) \\
  \rsp & n \equiv 1\ (2), n \ge 1 & 1 & \frac{n+1}4 & (\frac12, 0)
\end{array} $$

Moreover the second minimum is known in the following cases:

$$ \begin{array}{llccc}
  \rsp m & n & I & M_2(K) & C_2(K) \\ \hline
  \rsp n^2+1 & n \equiv 0\ (2) & 1 & \frac{n-1}2 - \frac1{4m}
             & (0, \frac{m \pm 1}{2m}), (\frac12, \frac{m \pm n}{2m}) \\
  \rsp & n \equiv 0\ (2) & 1 & \frac{n}8 \big(1 - \frac1{m}\big) & \\
  \rsp n^2-1 & n \equiv 0\ (2) & 1 & \frac{m-1}{2n+2}
             & (0, \pm \frac{n}{2n+2}) \\
  \rsp n^2+2 & n \equiv 0\ (2) & 1 & \frac{2n^3-3n^2+4n-2}{4n} & \\
  \rsp n^2+4 & n \equiv 1\ (2) & 2 & \frac{n^3-2n^2+5n-8}{4m}
\end{array} $$
The work on (2.6) is, as we see, not yet concluded. For
the proof helpful in any case is

\begin{quote}
  {\bf (2.7)} {\em Let $f$ be the absolute value of the norm in
    $\Q(\sqrt{m})$; then }
  \[ M(f) \le \frac{1}{2} \begin{cases} \max \{p, 2n+1-p\}, & \text{where }
    R=\Z[\sqrt{m}], \\
             & \quad m=n^2+p, \quad 1 \le p \le 2n \\
    \max \{p, 2n+2-p\}, & \text{where } R=\Z[\beta_m], \\
             & m=(2n+1)^2+4p, \quad 1 \le p \le 2n+1
              \end{cases} \]
   Here $\beta_m = (1+\sqrt{m})/2$. Further in the case of equality the
   minimum is attained $\bmod R$ at most at the points $z_1 = \sqrt{m}/2$ or
   $z_2 = (1+\sqrt{m})/2$ (if $R = \Z[\sqrt{m}]$), respectively at
   $z_1 = \beta_m/2$ or $z_2 = (1+\beta_m)/2$ (if
   $R = \Z[\beta_m]$).
\end{quote}

The result (2.7) is due to Heinhold \cite{Hei39} (see on this also
Lekkerkerker \cite[p.\ 422]{Lek69}); without difficulty we obtain
from it the Minkowski inequality for quadratic number fields:

\begin{quote}
  {\bf (2.8)} {\em Let $K$ be a quadratic number field with discriminant $d$;
    then $M(K) \le \sqrt{d}/4$. }
\end{quote}

If we choose the basis $\{1, \beta\}$ with $\beta = n+\sqrt{m}$,
where $m=n^2+p$, $1 \le p \le n$, then we quickly find that
$\max\{N(z):z \in F'\}$ equals the bounds of Heinhold given in (2.7);
here $F' = (0, 0.5) \times (0, 0.5)$. Unfortunately this does not
suffice to prove (2.7), because $F' \cup -F' \cup {F'}^\sigma \cup
-{F'}^{\sigma}$ is here not a fundamental domain (as was still the
case with the basis $\{1, \sqrt{m}\}$).

With (2.7) we can estimate some $M(K)$ from (2.6) from above: if
e.g.\ $m = n^2-1$, we write $m = (n-1)^2+2n-2$ and obtain from (2.7)
the bound $M(K) \le \frac{2n-2}{4} = \frac{n-1}2$.  For the estimation
of $M(K)$ from below we use (1.8). It should however be pointed out
that the bounds $\mu_2$ sometimes depend on $n$ (e.g.\ in the case $m
= n^2+2$). The proof that $M(K)$ in (2.6) is then really exact is in
these cases rather tedious (see on this Barnes and Swinnerton-Dyer).

Barnes and Swinnerton-Dyer have in their work proved some interesting
theorems on the Euclidean minima in real quadratic number fields; at
the end of the second part of their work they have pointed out that
these theorems can be generalized to fields of arbitrary unit rank
$>1$. In particular for their Theorem M (our 2.12) this generalization
has to this day only succeeded for fields with a single fundamental
unit (see van der Linden 1985, Theorem 8.7). We therefore wish to
describe more precisely the difficulties that arise here.  To this end
we show

\begin{quote}
  {\bf (2.9)} {\em Let $K$ be a number field and $\{u_1, \dots, u_t\}$ a
    system of independent units. Then there exists a finite set $Z
    \subset R$ with the property: for all $x \in K$ we have
    \[ M(K, x) = \inf_{u} \min_{z \in Z} N(xu + z), \]
    where the infimum is taken over all units $u$ of the subgroup of
    $R^\times$ generated by the $u_i$, and where $x_u \in F$ is
    determined by the congruence $x_u \equiv xu \bmod {R}$. }
\end{quote}

\begin{proof}
  Let $\{\alpha_1, \dots, \alpha_t\}$ be an integral basis and
  $Z = \{z = \sum r_i\alpha_i : |r_i| < \mu_i + 0.5\}$ for the
  $\mu_i$ from (1.12). We claim that then
  $$ M(K, x) = \inf_{u} \min_{z \in Z} N(xu + z) =: k $$
  holds. Otherwise $M(K, x) < k$, and there exists a $y \in
  R$ with $N(x-y) < k$. As in (1.12) there now follows the existence of a
  unit $u$ with $(x-y)u = \sum r_i\alpha_i$ and $|r_i| <
  \mu_i$. Now we have $(x-y)u = xu \equiv x_u \bmod {R}$ for some $x_u
  \in F$; hence $z = (x-y)u - xu \in R$. With $z = \sum s_i\alpha_i$
  and $xu = \sum t_i\alpha_i$ we then have $s_i = r_i - t_i$, hence $|s_i|
  \le |r_i| + |t_i| \le \mu_i + 0.5$ and thus $z \in Z$ in
  contradiction to the assumption $M(K, x) < k$.
\end{proof}

If $x \in K$, then there are evidently only finitely many $x_u \in F$
(for with $x = a/b$ we have $x e = x \bmod {R}$ for every unit $e \in
R^\times$).  If on the other hand $x \notin K$, then we can never
have $x_u = x \bmod {R}$ (for from $x_u = -x = x(e-1) \in R$ it
follows that $x \in K$). In this case the $x_u \in F$ are then
pairwise distinct, and because of the compactness of $F$ there exists
an accumulation point $w$ of the $x_u$. For the estimation of $M(K,
w)$ we use

\begin{quote}
  {\bf (2.10)} {\em With the designations of (2.9) let $w$ be an
    accumulation point of the $x_u$. Then $M(K, x) \le M(K, w)$.}
\end{quote}

\begin{proof}
  If $w$ is an accumulation point of the $x_u$, then there exists a
  subsequence $(x_{u_e})$ of the $(x_u)$ converging to $w$, and since
  the norm is continuous, we obtain
  $$ \inf_{u} N(x_u + z) \le \lim N(x_e + z)
     = N(\lim x_e + z) = N(w + z) \quad
     \text{for all } z \in Z. $$
  Now however with $w$ also every $w_u$ is an accumulation point of the $x_u$,
  because with $x_e \to w$ also $x_{eu} \to w_u$ holds (here of course
  $x_{eu} \in F$ is defined by $xeu \equiv x_{eu} \bmod {R}$ for
  units $e, u \in R^\times$). Thus we can in the
  above reasoning replace $w$ by $w_u$ and have for all $u$
  and for all $z \in Z$ the inequality
  $$ \inf_{u} N(x_u + z) \le N(w_u + z), $$
  and the formation of the minimum over all $z \in Z$ yields
  $$ M(\uK, x) = \inf_{u} \min_{z \in Z} N(x_u + z) \le \min_{z \in Z}
  N(w_u + z). $$

  Here we were allowed to interchange inf and min because $Z$ is
  finite.  Since finally $M(\uK, w) = \inf_{u} \min_{z \in Z} N(w_u +
  z)$, it follows that $M(\uK, x) \le M(\uK, w)$ as claimed.
\end{proof}

We now show the generalization of Theorem L (Barnes and
Swinnerton-Dyer):

\begin{quote}
  {\bf (2.11)} {\em Let $K$ be a number field.}
  \begin{enumerate}
  \item[(i)] {\em For every $x \in \uK$ there exists a $w \in \uK$ with
    $M(\uK,x) = M(\uK,w)$, and $M(\uK,w)$ is attained.}
  \item[(ii)] {\em The set $\{M(\uK,x) : x \in \uK\} \subset \R$
    is closed.}
  \end{enumerate} 
\end{quote}

For a better understanding of (2.11) it is worth keeping the example $D(13)$
in mind: here we had a sequence $x_1, x_2, \dots$ of points with
$M(K,x_i) = \frac{1}{3}$, where this minimum is attained at no $x_i$.
On the other hand $M(K,w) = \frac{1}{3}$ at the accumulation point
$w = (\frac{1}{6}, \frac{1}{6})$ of the $x_i$, and at $w$ the minimum is
attained.

\begin{proof}[Proof of (i)]
  Let $(x_e)$ be a subsequence of the $x_u$ with
  $$ M(\uK,x) = \lim_{e \to \infty} \min_{z \in Z} N(x_e + z). $$

  Without loss of generality we may assume that the $x_e$ are pairwise
  distinct: if $M(\uK,x)$ is not attained, this is certainly
  possible; in the other case there is nothing to show. Now let $w \in
  F$ be an accumulation point of the $x_e$; then there exists a
  subsequence of the $x_e$ converging to $w$, which we shall again
  denote by $(x_e)$. Now we have
  $$ M(\uK,x) = \lim_{e \to \infty} \min_{z \in Z} N(x_e + z) =
     \min_{z \in Z} N(w + z) \ge M(\uK,w), $$
  while from (2.10) $M(\uK,x) \le M(\uK,w)$ follows. Thus
  $M(\uK,x) = M(\uK,w)$, and because of
  $M(\uK,w) = \min_{z \in Z} N(w+z)$ the value $M(\uK,w)$ is
  indeed attained.

  (ii) Let $k \in \R$ be an accumulation point of the set
  $\{M(K,x) : x \in K\}$ and $(x_i)_{i \in \N}$ a sequence of
  points in $F$ with $\lim_{i \to \infty} M(\uK,x_i) = k$. We must then
  find a $w \in F$ with $k = M(\uK,w)$.
  By (i) we may assume that $M(\uK,x_i)$ is attained for every $i \in
  \N$. Thus for every $i \in
  \N$ there exists a unit $u = u(i) \in R^\times$ with the property
  $$ M(\uK,x_i) = \min_{u'} \min_{z \in Z} N(y_{u'} + z) $$
  (where $u'$ need only run through the units $1, u, u^2, \dots, u^{m-1}$).
  
  Now let $w$ be an accumulation point of the $x_{iu}$ and $(x_e)$ a
  subsequence of the $x_{iu}$ converging to $w$.  Then
  \begin{align*}
  k & = \lim_{i \to \infty} M(\uK,x_i)
      = \lim_{i \to \infty} \min_{z \in Z} N(x_{iu} + z) \\
    & = \lim_{e \to \infty} \min_{z \in Z} N(x_e + z)
      = \min_{z \in Z} N(w + z) \ge M(K,w)
  \end{align*}
  By (2.10) we also have $M(\uK,x) \ge \inf_i M(\uK,x_{iu}) = k$,
  consequently $M(\uK,x) = k$.
\end{proof}

An important special case of (2.11.ii) is: there exists an $x \in \uK$ with
$M(\uK,x) = M(\uK)$; this observation is due to Heinhold \cite{Hei39}. Barnes
and Swinnerton-Dyer have conjectured that there even exists an $x \in K$ with
$M(\uK,x) = M(\uK)$, and that in particular $M(\uK) = M(K)$ holds. If
we recall the example $D(13)$, this conjecture seems
at first glance almost trivial: seek an $x \in K$ with
$M(\uK,x) = M(\uK)$ and then choose $w$ as the accumulation point of the $x_u$. The
catch here is of course that $w$ need not lie in $K$
(if the $x_u$ accumulate at only finitely many accumulation points as in $D(13)$, then
there would exist a unit $u$ with $w \equiv w_u \bmod {R}$, and this
would indeed imply $w \in K$; however, it does not seem possible to exclude so easily the possibility
that the $x_u$ have infinitely many accumulation points).
Even the less far-reaching conjecture $M(K) = M(\uK)$ has
hitherto been proved only for fields of unit rank 1
(Theorem M of Barnes and Swinnerton-Dyer (for $n = 2$),
Prop.\ 5.2.\ of van der Linden (for unit rank 1)):

\begin{quote}
  {\bf (2.12)} {\em For all $\eps > 0$ there exists a $y \in K$
    with $M(\uK, y) > M(\uK) - \eps$; in particular $M(K) = M(\uK)$.}
\end{quote}

\begin{proof}
  Let $u$ be a unit of $K$ with $|u|_1 > 1$; from the proof of the
  Dirichlet unit theorem it follows that the powers of $u$
  ``stay away from $1$'' in the sense that there exists a $c > 0$
  with $|u^m - 1|_1 > c$ and $|u^m - 1|_2 > c$ for all $m \in \Z$. Let
  further $x \in K$ be a point with $M(\uK, x) = M(\uK)$ (such an $x$
  exists by 2.11.ii). Without loss of generality we may assume that $x$ is an
  accumulation point of the $x_u$ in $F$ (otherwise we replace $x$ by
  such an accumulation point $x'$, for by (2.9) we have
  $M(\uK, x') \ge M(\uK, x)$, and because of $M(\uK, x) = M(\uK)$ and
  $M(\uK, x') \le M(\uK)$ we must also have $M(\uK, x') = M(\uK)$).
  Thus there exists a power $u^m$ of $u$ such that for an arbitrarily
  prescribed $\delta > 0$ we have $u^m x \equiv x + z \bmod {R}$ and
  $|z|_j < \delta c$ for $j = 1,2$. Therefore there exists an $a \in R$ with
  $u^m x = a + x + z$, i.e.\ with $(u^m - 1)x = a + z$. Now we set
  $y = a/(u^m - 1)$ and have $y \in K$, as well as $x - y = z/(u^m -
  1)$. Thus we have found a $y \in K$ that not only lies near
  $x$, but for which also the powers $uy, u^2 y, \dots, u^m
  y$ lie near $ux, u^2 x, \dots, u^m x$; more precisely: for $0 \le k
  \le m$ we have
  $$ |u^k x - u^k y|_1 = |u^k (x - y)|_1 \le |u^m (x - y)|_1
   = \Big|\frac{z}{1 - u^{-m}}\Big|_1 < \frac{\delta c}{c} = \delta $$
   because of $|u|_1^k \le |u|_1^m$ and $|1 - u^{-m}|_1 < \frac1c$. Analogously
   we obtain
   $$ |u^k x - u^k y|_2 = |u^k (x - y)|_2 \le |(x - y)|_2
   = \Big|\frac{z}{u^m - 1}\Big|_2 < \frac{\delta c}{c} = \delta $$
  because of $|u|_2 < 1$ and $1/|u^m - 1|_2 < \frac1c$. These two
  inequalities can also be written in the form $|x_e - y_e|_1 < \delta$,
  $|x_e - y_e|_2 < \delta$ for $e = u^k$. Our task is
  now to compare $M(\uK, x)$ and $M(\uK, y)$. By (2.9) we have
  $$ M(\uK, x) = \inf_{u} \min_{z \in Z} N(x_u + z), $$
  while because of $u^m y \equiv y \bmod {R} $
  $$ M(\uK, y) = \min_{u'} \min_{z \in Z} N(y_{u'} + z) $$
  holds, where $u'$ need only run through the units $1, u, u^2, \dots, u^{m-1}$.
  Because of $x_u, y_u \in F$ and $z \in Z$ the
  points $x_u + z$ and $y_u + z$ lie in the bounded and closed
  set $Z' = \{ \sum r_i \alpha_i : |r_i| \le \mu_i + 1\}$.
  Consequently there exists a $\delta>0$ with the property that for all
  $x,y,z'$ we have: from $|x-y|_1 < \delta$ and $|x-y|_2 < \delta$ it follows that
  $|N(x)-N(y)| < \eps$ for the prescribed $\eps > 0$.
  As $u'$ runs through the units $1, u, \dots, u^{m-1}$, we have for
  all $z \in Z$ and $j=1,2$ $|(x_u + z)-(y_u + z)|_j < \delta$; therefore
  also $|N(x_u + z)-N(y_u + z)| < \eps$, in particular therefore
  $N(x_u + z) < N(y_u + z) + \eps$. It now follows immediately that
  $$ \min_{u'} \min_{z \in Z} N(x_{u'} + z) <
     \min_{u'} \min_{z \in Z} N(y_{u'} + z) + \eps = M(K, y) + \eps,
     \text{ but because of} $$
   $$ M(\uK, x) = \inf_{u} \min_{z \in Z} N(x_u + z) \le
      \min_{u'} \min_{z \in Z} N(x_{u'} + z) $$
   this is already the claim.
\end{proof}
The difficulty in a generalization of this theorem to fields
of arbitrary unit rank $\ge 1$ lies solely in the
construction of a $y \in K$ with the property that the finitely
many $y_u$ lie near the corresponding $x_u$. I have not, however,
yet succeeded in carrying out such a construction.
For the sake of completeness we still wish to give the generalization of the
``Theorem G'' of Barnes and Swinnerton-Dyer to fields of
unit rank $\ge 1$; it reads

\begin{quote}
  {\bf (2.13)} {\em If $M(\uK, x) < k$ for all $x \in \uF$ except for a
    finite exceptional set $G$, then there exists a $k' < k$ such that
    $M(\uK, x) \le k'$ for all $x \in \uF \setminus G$. }
\end{quote}

From this it follows in particular that $M_1(K)$ is isolated if it is attained only at
finitely many places in $F$. The proof of Barnes
and Swinnerton-Dyer also works in the general case if we
choose the unit $u$ so that $|u|_j \ne 1$ for all $j$ with $1 \le j
\le r+s$ and use (2.3).
We still wish to say a few words on the determination of $M^2(K)$: the
notions and also the theorems from \S\  2 have direct analogues if we
replace $R$ by the set of continued fractions of length $\le k$.
The same also holds for (1.8) and (1.12); however the
analogues of (1.8) and (1.12) are probably worthless, since there exist
e.g.\ infinitely many continued fractions $r+s\sqrt{m}$ of length 2 with
bounded $r,s \in \Q$ and therefore the conditions
(1.8.a,b,c) allow infinitely many possible $z \in K$.
In order therefore to estimate $M^2(K)$ from below, we cannot fall back on (1.8)
or (1.12); besides the estimate (0.13) the following method has proved useful in
rings of class number $>1$, which we wish to demonstrate by the
example $K = \Q(\sqrt{10})$ (this answers a question posed by
Cooke \cite[p.\ 75 bottom]{Coo77}).
We choose a set of continued fractions of length 2, e.g.\
\[
R_2 = \left\{0, \frac{1}{2}, \frac{1+\sqrt{10}}{2}, \frac{4+\sqrt{10}}{6}, \frac{1}{3}, \frac{\sqrt{10}}{3}\right\}
\]
their denominators have the norms 1, 4, 4, 6, 9, 9. Then we try, for
every point $x$ from $F = (0, 0.5) \times (0, 0.5)$, to find a $y$ that
is congruent $\bmod R$ to a continued fraction $a_2/b_2$ from $R_2$ and
satisfies $|N_{K/\Q}(x-y)| < k \cdot |N_{K/\Q}(b_2)|$.
If we choose $k = 0.99$, we can thus cover the whole of $F$ except for the
exceptional set $S = (0, 0.001) \times (0.499, 0.5)$. As in the
norm-Euclidean case it now follows easily that $z = \sqrt{10}/2$ is the
only possible exceptional point. Thus only $M(K,z)$ remains to be
determined. Since $y = a_2/b_2 = \frac{3}{2} = 1 + \frac{1}{2}$ is a
continued fraction of length 2 with $N_{K/\Q}(b_2) = 4$ and
$|N_{K/\Q}(z-y)| = \frac{1}{4}$, we certainly have $M^2(K,z) \le
1$; the corresponding division chain is $\sqrt{10} = 2q_1 + r_1$, $2 =
r_1 q_2 + r_2$ with $q_1 = 1$, $q_2 = 2$, $r_1 = -2 + \sqrt{10}$,
$r_2 = 6 - 2\sqrt{10}$ and $|N_{K/\Q}(r_2)| = 4 = N_{K/\Q}(2)$.
Now however there exists no division chain of length 2 with
$|N_{K/\Q}(r_2)| \le 4$, for if $P = (\sqrt{10}, 2)$ is the
prime ideal of norm 2 above $(2)$, then $r_1, r_2 \in P$, and
in particular $N_{K/\Q}(r_2) \equiv 0 \bmod 2$. Therefore we would have to have
$|N_{K/\Q}(r_2)| = 0$ or $|N_{K/\Q}(r_2)| = 2$;
the latter is impossible since $P$ is not a principal ideal, and for the same
reason we also cannot have $r_2=0$, since otherwise $P = (r_1)$ would follow.
Thus $M^2(K) = M^2(\uK) = M^2(K,z) = 1$, and $K$ has
Euclidean depth 1 (for exactly as above we show that there exists no division chain starting from
$(\sqrt{10},2)$ with $|N_{K/\Q}(r_k)|
\le 4$). If we define the sets\label{DefB}
$B_j$ ($j = 1, 2, \dots, \infty$) by $B_j = (F \cap K) \setminus F_j$
($F$ is a fundamental domain; the $F_j \subset K$ are defined on
p.\ \pageref{pEi}), then
here we have $B_1 = B_2 = \dots = B_\infty = \{\pm \sqrt{10}/2\}$.
In the same manner Table~\ref{Tab6} has come into being:
\begin{table}[h]
\centering
\begin{tabular}{llll}
\toprule
$m$ & $M^2(K)$ & $B_1$ & $B_2 = B_3 = \dots$ \\
\midrule
6 & $1/4$ & $\{\}$ & $\{\}$ \\
10 & $1$ & $\{(0, 1/2)\}$ & $\{(0, 1/2)\}$ \\
14 & $1/4$ & $\{(1/2, 1/2)\}$ & $\{\}$ \\
15 & $1$ & $?$ & $\{(1/2, 1/2)\}$ \\
26 & $1$ & $?$ & $\{(0, 1/2)\}$ \\
30 & $3/2$ & $?$ & $\{(0, 1/2)\}$ \\
34 & $1$ & $?$ & $\{(1/2, 1/2) + (1/2, 1/2)\}$ \\
35 & $7/5$ & $?$ & $\{(0, 1/2), (1/2, 1/2)\}$ \\
39 & $5/2$ & $?$ & $\{(1/2, 1/2)\}$ \\
65 & $1$ & $\{(1/2, 1/2)\}$ & $\{(1/2, 1/2)\}$ \\
85 & $1$ & $?$ & $\{(1/2, 1/2) + (1/2, 1/2)\}$ \\
\bottomrule
\end{tabular}
\caption{Euclidean minima}\label{Tab6}
\end{table}
From this table we can gather that e.g.\ the Euclidean depth
of $D(26)$ and $D(30)$ is equal to 2. Already Cooke has tried to determine the
Euclidean depth of these two rings. For $D(30)$ he was able
to show that the Euclidean depth is $\le 3$, and he conjectured
that it is equal to 2. In the case $D(26)$ he has shown that the
Euclidean depth is finite; that he could not determine it
was due to a computational error. He writes on \cite[p.\ 81]{Coo77} that he can cover the whole of $F$ with
continued fractions of length 2 except for the following points:
$(0, \frac{1}{2})$, as well as
$(\pm \frac{2}{5}, \pm \frac{2}{5})$. However $\sqrt{26}$ is a
continued fraction of length 2 and e.g.\
$$ \Big| N_{K/\Q}\Big(z - \frac{\sqrt{26}}{3}\Big) \Big| = \frac{2}{45}
     = \frac{2}{5} \cdot N_{K/\Q}(3) < N_{K/\Q}(3) $$
for $z = \frac{2 + 2 \sqrt{26}}{5}$.
Thus it also follows from his computations that only the point
$\sqrt{26}/2$ remains uncovered and the Euclidean depth is equal to 2.
Further it remains to be noted that $D(10)$ has Euclidean depth 1
(and not 2, as Cooke writes), since here $B_1 = B_\infty =
\{\sqrt{10}/2\}$; the corresponding holds for $D(65)$.
With this method we can also confirm the results of Cooke \cite{Coo77}
according to which the rings $D(m)$ for
$$ m = 14, 22, 23, 31, 38, 43, 46, 53, 61, 69, 77, 89, 93, 97, 113, 129,
       133, 137, 181, 253 $$
are 2-stage norm-Euclidean; moreover we obtain for
$$ m = 47, 59, 62, 67, 71, 101, 109, 149, 157, 161, 173, 177, 193,
   197, 201, 213 $$
new 2-Euclidean rings.
We must, however, note here that the bounds achieved for $M^2(K)$
as well as the possible exceptional points depend on the number of
continued fractions of length 2 used; thus in $D(14)$ we can
with continued fractions of length 2 whose denominators have norm $\le 11$
cover the whole of $F$ except for the two points $x_1 = (\frac{1}{2},
\frac{9}{14})$ and $x_2 = (\frac{1}{2}, \frac{1}{2})$. By (0.13)
we have $M^2(x_2) = \frac{1}{4}$; because of $|N_{K/\Q}(x_1 -
\frac{2+\sqrt{14}}{7 \cdot 2})| = \frac{1}{28}$ and because
$\frac{2+\sqrt{14}}{14}$ is a continued fraction of length 2 with denominator
$7\cdot 2/14$, we certainly also have $M^2(x_1) \le \frac{1}{4}$ (this
already suffices to prove $M^2(K) = \frac{1}{4}$).
If, however, we also use continued fractions whose denominators have norm 13, and
note that $z = (6+6\sqrt{14})/13$ is one such fraction, we find
that $|N_{K/\Q}(x_1 - z)| < \frac{1}{32}$; this shows that
$x_2$ is the only exceptional point for $k = 1/4$.
As Cooke has remarked, the points from $B_\infty$ stand in a
relation to the class number of $K$; in order to see this we define
a map $\varphi : K \to \mathrm{Cl}(K)$ by
$\varphi(\frac{\alpha}{\beta}) = [(\alpha, \beta)]$ (i.e.\ to an $x =
\frac{\alpha}{\beta}$ is assigned the ideal class generated by the ideal $(\alpha, \beta)$). $\varphi$ is well-defined: if
$\alpha/\beta = \alpha'/\beta'$, then $\alpha\beta' =
\alpha'\beta$ and thus $(\alpha', \beta) \sim (\alpha')(\alpha,
\beta) = (\alpha\alpha', \alpha'\beta) = (\alpha\alpha', \alpha\beta')
\sim (\alpha', \beta')$, and we have $[(\alpha, \beta)] = [(\alpha',
  \beta')]$.
Quite simply there now follows

\begin{quote}
  {\bf (2.14)} {\em $\varphi(B_\infty)$ contains all non-trivial
    ideal classes.}
\end{quote}

\begin{proof}
  Let $I$ be an ideal in $R$, but not a principal ideal. We choose
  $\beta \in I \setminus \{0\}$ so that $|N_{K/\Q}(\beta)|$ becomes minimal,
  and then an $\alpha \in I$ with $I = (\alpha, \beta)$. By
  subtracting from $\alpha$ a suitable multiple of $\beta$,
  we can achieve $\frac{\alpha}{\beta} \in F$ without destroying the
  property $I = (\alpha, \beta)$.
  We now claim that $\frac{\alpha}{\beta} \in B_\infty$. For
  every division chain starting from $(\alpha, \beta)$ we have
  $r_k \in (\alpha, \beta) \setminus \{0\}$ (if $r_k = 0$, then
  $I = (r_{k-1})$ would be a principal ideal), and by the choice of
  $\beta$ we must have $|N_{K/\Q}(r_k)| \ge |N_{K/\Q}(\beta)|$. Thus
  $K$ is not quasi-Euclidean at $x = \frac{\alpha}{\beta}$ and consequently
  $\frac{\alpha}{\beta} \in B_\infty$.
\end{proof}

In this way the knowledge of $B_\infty$ already lets us determine the
class number of $K$: for $D(10)$, $D(15)$, $D(26)$ etc.\ $B_2 =
B_\infty$ contains only a single point, hence $h(K) \mid 2$ in these
cases. On the other hand here $\varphi(B_\infty)$ is never the
principal class, and we can conclude $h=2$ in these cases.  Somewhat
more interesting is $D(35)$: here $B_\infty$ contains three distinct
points, and we first only have $h(K) \mid 4$. However the ideals $5_1
= (2+\sqrt{35}, 5)$ and $2_1 = (1+\sqrt{35}, 2)$ lie in the same ideal
class, so that here also $h=2$ results.  This raises the question when
$\varphi(x) = \varphi(x')$ for two points $x,x'\in K$. To this end let
$SL_2(R)$ be the group of $2\times2$-matrices with entries from $R$
and determinant $+1$; $SL_2(R)$ acts on $K$ by
\[
\begin{pmatrix} a & b \\ c & d \end{pmatrix}
\begin{pmatrix} \alpha \\ \beta \end{pmatrix} =
\frac{a\alpha + b\beta}{c\alpha + d\beta}
\]
(strictly speaking we must replace $K$ by $K \cup \{\infty\}$ if we do
not wish to exclude a priori that $c\alpha + d\beta = 0$;
if, however, we restrict the domain of definition of $SL_2(R)$ to
$B_\infty$, this can never occur). Two points $x =
\alpha/\beta$ and $x' = \alpha'/\beta'$ are now called equivalent ($x
\equiv x'$) if there exists an $A\in SL_2(R)$ with $A(x) = x'$. We then
have

\begin{quote}
  {\bf (2.15)} {\em Let $K$ be an algebraic number field with
    infinite unit group. Then $B_\infty$ contains precisely $h-1$
    equivalence classes of points.}
\end{quote}

For this we must first show that $\varphi(x) = \varphi(x')$ holds
precisely when $x \equiv x'$, and secondly that $B_\infty$ contains no
$x\in K$ such that $\varphi(x)$ is the principal class. The latter
result seems anything but elementary and was proved by Vaserstein
\cite{Vas72}; it seems desirable to me to translate his proof into the
language of algebraic number theory.  In the proof of the first fact
Cooke uses the following observation of Hurwitz \cite{Hur95a}:

\begin{quote}
  {\bf (2.16)} {\em The ideals $(\alpha, \beta)$ and $(\alpha', \beta')$
    are equal if and only if there exists an $A\in SL_2(R)$
    with $A(\frac{\alpha}{\beta}) = \frac{\alpha'}{\beta'}$. }
\end{quote}

The proof that Cooke gives for (2.15) is a pure existence proof
and therefore cannot be satisfactory. The original proof of Hurwitz
on the other hand at least gives us a hint how to construct $A$,
moreover it is simpler than Cooke's proof:

\begin{proof}[Proof of (2.15)]
  Let $A \in \text{SL}_2(\R)$ and $A(\alpha/\beta) =
  \alpha'/\beta'$. Then there exists an $\eta \in \uK$ with $\alpha'\eta =
  a\alpha + b\beta$, $\beta'\eta = c\alpha + d\beta$. This implies
  $(\alpha, \beta) = (\eta)(\alpha', \beta')$. $(\alpha', \beta')$
  hence $\varphi(\frac{\alpha}{\beta}) = \frac{\alpha'}{\beta'}$.
  Now conversely let
  $\varphi(\frac{\alpha}{\beta}) = \frac{\alpha'}{\beta'}$; then
  $(\alpha, \beta) = (\alpha', \beta')$, and by (2.17) there exist
  $\gamma, \delta \in \R$ with $(\gamma, \delta) = (\alpha, \beta)$
  and $\frac{\gamma}{\delta} = \frac{\alpha}{\beta}$. With (2.16)
  we obtain an $A \in \text{SL}_2(\R)$ with
  $A(\frac{\alpha}{\beta}) = (\frac{\alpha'}{\beta'})$, and
  the same $A$ then also satisfies
  $A(\frac{\alpha}{\beta}) = (\frac{\alpha'}{\beta'})$.
\end{proof}
As an example we wish to show that the points $2\sqrt{35}/5$ and
$(1+\sqrt{35})/2$ are equivalent. To this end we write
$(1+\sqrt{35})/2$ with (2.16) in the form $(1+\sqrt{35})/2 =
\alpha/\beta$ with $\alpha' = 20+3\sqrt{35}$ and $\beta' =
5-\sqrt{35}$; hence $(2\sqrt{35}, 5) = (\alpha', \beta')$. Since
$D(35)$ has class number 2 and $(2\sqrt{35}, 5)^2 = (5)$, we can
as in (2.16) write with $\alpha = 2\sqrt{35}$ and $\beta = 5$:
$5 = \alpha\gamma + \beta\delta = \alpha'\gamma' + \beta'\delta'$ with
$\gamma, \delta, \gamma', \delta' \in (2\sqrt{35}, 5)$ for $\gamma =
-\sqrt{35}$, $\delta = 15$, $\gamma' = \sqrt{35}$, $\delta' = -20$ and
obtain the matrix
\[
A = \begin{bmatrix} -41-4\sqrt{35} & 60+17\sqrt{35} \\
    -7-2\sqrt{35} & 29+3\sqrt{35} \end{bmatrix},
\]
which is in fact unimodular and satisfies $A(\alpha/\beta) = \alpha'/\beta'$.
With the help of the map $\varphi: \uK \to
\text{Cl}(\uK)$ defined above we can also give a simple criterion for
semi-Euclidean rings:

\begin{quote}
  \textbf{(2.18)} {\em A number ring $R$ is semi-Euclidean if and only if
    for every $x \in \mathbb{B}_1$ we have $\varphi(x) \neq [(1)]$, i.e.\
    with $x = \frac{\alpha}{\beta}$ the ideal $(\alpha, \beta)$ is not principal.}
\end{quote}

\begin{proof}
  Let $\varphi(x) \neq [(1)]$ for all $x \in B_1$, and let
  $\alpha, \beta \in R$ with $(\alpha, \beta) = R$. Without loss of generality we may
  choose $\alpha \bmod \beta$ so that $\frac{\alpha}{\beta} \in \uF$.
  Because of $\varphi(x) = [(1)]$ for $x = \alpha/\beta$ the point $x$ does
  not lie in $B_1$; this means that there exists a $y \in \R$ with
  $|N_{\uK/\Q}(x-y)| < 1$. Thus $\R$ is semi-Euclidean.
  Conversely let $R$ be semi-Euclidean and
  $x = \frac{\alpha}{\beta} \in \uK$. If $(\alpha, \beta)$ is a
  principal ideal (hence $\varphi(x) = [(1)]$), then, because $R$ is
  semi-Euclidean, there exists a $y \in \R$ with $|N_{\uK/\Q}(x-y)| < 1$.
  Thus $x$ does not lie in $B_1$.
\end{proof}

An analogous result is obtained for $k$-stage semi-Euclidean rings if
we replace $B_1$ by $B_k$. Thus e.g.\ $D(10)$ is semi-Euclidean,
because $B_1 = \{\sqrt{10}/2\}$ and $(\sqrt{10}, 2)$ is not a
principal ideal.

\section*{{\sc Remarks on} \S\ 2}
\addcontentsline{toc}{section}{{\sc Remarks on} \S\ 2}

The notion ``Euclidean minimum'' in number fields is unknown in the
mathematical literature; instead the expression ``inhomogeneous
minimum of the associated norm form'' has hitherto been used. The
first result on Euclidean minima is due to Minkowski, who for real
quadratic number fields of discriminant $d$ has shown that $M(K) =
M_2(K) \le \sqrt{d}/4$. Minkowski has further conjectured that for
totally real number fields of degree $n$ the inequality $M(K) \le
2^{-n}\sqrt{d}$ holds; this has however hitherto been proved only in
special cases ($n = 2, 3, 4, 5$; see on this Remak \cite{Rem34}, Dyson
\cite{Dys48} and Skubenko \cite{Sku72}). According to Chebotarev we
have at least $M(K) \le 2^{-n/2}\sqrt{d}$ (see on this also Hardy and
Wright, \cite[pp.\ 456--458]{HW}. For number fields of arbitrary
signature not even an analogous conjecture is known.

In the opposite direction the result of Davenport \cite{Dav51} is to
be mentioned, according to which for real quadratic fields $M(K) \ge
\frac{\sqrt{d}}{128}$ holds. At the end of this paper Davenport
remarks that Prasad has improved the constant $\frac{1}{128}$ to
$\frac{1}{36}$. As Ennola \cite{Enn58b} writes, Prasad has however
discovered an error in his work (before its publication). The
``Davenport constant'' $\sup \frac{M(K)}{\sqrt{d}} \ge \frac{1}{128}$
was then improved by Cassels \cite{Cas52} to $\frac{1}{36}$, and
finally Ennola succeeded in showing
$$M(K) \ge \frac{\sqrt{d}}{16+6\sqrt{6}}
          \approx \frac{\sqrt{d}}{30.69};$$
this has remained the best result to this day.

The explicit determination of Euclidean minima in real quadratic
number fields was undertaken by Heinhold \cite{Hei39}, Davenport
\cite{Dav46}, Varnavides \cite{Var48a,Var48b}, Bambah
\cite{Bam50,Bam51} and Inkeri \cite{Ink49}. After the fundamental
works of Barnes and Swinnerton-Dyer only Godwin \cite{God62,God65b}
and Varnavides \cite{Var70} have still occupied themselves with the
determination of Euclidean minima in quadratic number fields.  The
first minima of cubic number fields were computed by Davenport
\cite{Dav47}; further results in this direction are due to Prasad
\cite{Pra49}, Clarke \cite{Cla51a}, Samet \cite{Sam54},
Swinnerton-Dyer \cite{Swi54}, Godwin \cite{God55}, Smith
\cite{Smi69,Smi71} and Taylor \cite{Tay75,Tay76} (see on this \S\  4).

For number fields of degree $\ge 4$ no Euclidean
minima have hitherto been determined; only Cohn and Deutsch \cite{CD86} have
conjectured that for $K = \Q(\sqrt{2+\sqrt{2}})$ we have $M(K) =
\frac{1}{2}$ and $M_2(K) = \frac{1}{4}$, while for $K =
\Q(\sqrt{3+\sqrt{2}})$ they have conjectured that $M(K) =
\frac{1}{2}$ and $M_2(K) = \frac{7}{16}$.

\chapter*{\S\ 3 Quadratic Number Fields}
\setcounter{chapter}{3}
\addcontentsline{toc}{chapter}{\S\ 3 Quadratic Number Fields}
\markboth{Euclidean Rings}{\S\ 3 Quadratic Number Fields}

It is already a classical result that quadratic number fields $K =
\Q(\sqrt{m})$ are norm-Euclidean precisely for the values
$$ d = -11, -8, -7, -4, -3, 5, 8, 12, 13, 17, 21, 24, 28, 29, 33,
       37, 41, 44, 57, 73, 76 $$
(where $d = \disc K$ denotes the discriminant of $K$).

We shall see that the proof of this fact can be given rather clearly
and comparatively briefly if we consistently apply (1.5) and (1.8).
Moreover, the ideas presented here make it possible, at least in part,
to classify norm-Euclidean rings of higher degree (see in particular
\S\ 5).

We now collect some well-known properties of quadratic number fields.
Every quadratic number field is of the form
$K = \Q(\sqrt{m})$ for a square-free $m \in \Z$, and, depending on the
case, we have

\begin{align*}
  m \equiv 1 \bmod{4} : & \left\{1, \frac{1+\sqrt{m}}{2}\right\}
                       \text{ is an integral basis, } \disc K = m, \\
  m \equiv 2, 3 \bmod{4} & : \{1, \sqrt{m}\} \text{ is an integral basis, }
                       \disc K = 4m.
\end{align*}

If $(\frac{\cdot}{p})$ denotes the Kronecker symbol, the
decomposition law (by which we mean a description
of how prime ideals behave under field extension) can be
described as follows:

\begin{enumerate}
\item[(I.)] $(\frac{d}{p}) = +1$: $(p) = P_1P_2$, \quad
  $\|P_1\| = \|P_2\| = p$ \quad ($p$ is called \textbf{split})
\item[(II.)] $(\frac{d}{p}) = -1$: $(p) = P$, \quad $\|P\|
  = p^2$ \quad ($p$ is called \textbf{inert})
\item[(III.)] $p \mid d$ : $(p) = P^2$, \quad $\|P\| = p$ \quad ($p$
  is called \textbf{totally ramified}).
\end{enumerate}

The unit group $\cO_K^*$ of quadratic number rings
has a fairly simple structure. If $K$ is imaginary quadratic, then
$\cO_K^* = \{\pm 1\}$, except for $K =
\Q(\sqrt{-1})$ and $K = \Q(\sqrt{-3})$, where the whole of $\cO_K^*$
is formed by the 4th and 6th roots of unity, respectively. In
real quadratic number fields, Dirichlet's
unit theorem tells us that every unit $e \in \cO_K^*$ can be written as $e =
\pm u^k$ for some $k \in \Z$, where $u$ denotes the
fundamental unit of $K$ (which is in general uniquely determined by the
requirement $u>1$).

The parity of the class number $h(K)$ can also be given easily. To this end
let $m \in \N$; we then distinguish

$K = \Q(\sqrt{-m})$: $h(K)$ is odd if and only if $m
\equiv 1, 2$ or $m \equiv 3 \bmod{4}$ is prime, i.e.\ if and only if
$\disc K$ is a prime power.

$K = \Q(\sqrt{m})$: $h(K)$ is odd if and only if:

\begin{enumerate}
\item[a)] $m = p$ is prime;
\item[b)] $m = pq$, $p$ and $q$ prime, $p=2$ or $p\equiv 3 \bmod{4}$,
  $q\equiv 3 \bmod{4}$.
\end{enumerate}

An elementary proof for imaginary quadratic $K$ is found in
Connell\index[N]{Connell} \cite{Con62}, and for real quadratic $K$ in
Redei\index[N]{Redei@R\'edei} \cite{Red60}. Another approach to these
results is provided by Gauss's genus theory (see
Zagier\index[N]{Zagier} \cite{Zag81}). Finally, a
further proof can be obtained from the analytic class number formula (part
of such a proof is found in Hasse\index[N]{Hasse} \cite[Va,Vb]{Has64}).

Since Euclidean rings have class number 1, we need only consider the
possibilities listed above. The simplest case turns out to be
that in which $\disc K \equiv 0 \bmod{4}$. The fact that $K$ can then be
norm-Euclidean only for the values $m = 2, 3, 6, 7, 11, 19$ was
shown independently by J.\ Fox\index[N]{Fox} \cite{Fox35} and
E.\ Berg\index[N]{Berg} \cite{Ber35} (the work of J.\ Fox appears in
the references of almost all authors under the name Fox
Keston; here we follow Bull.\ Am.\ Math.\ Soc.\ 41 (1935), p.\ 186,
where her name is given as Jeanette Fox).

The simplest proofs in the case $\disc K \equiv 0 \bmod{4}$, however,
are probably due to Behrbohm\index[N]{Behrbohm} and
R\'edei\index[N]{Redei@R\'edei} \cite{BR36}, and these we now
wish to present as well.

First we show

\begin{quote}
  \textbf{(3.1)} {\em If $a \in \N$ and $( \frac{m}{a}) = -1$, then
    neither $a$ nor $-a$ is a norm from $\cO_K$.}
\end{quote}

\begin{proof}
  Since $( \frac{m}{a} ) = -1$, $a$ contains a prime factor
  $p$ that divides $a$ to an odd power (i.e.\ $a = p^e b$, $p
  \nmid b$, $e \equiv 1 \bmod{2}$) and for which $( \frac{m}{p}) = -1$.
  Since $p \nmid b$, if $a$ were an ideal norm from
  $\cO_K$, then $p^e$ would also be one; since $e$ is odd and $\|(p)\| = p^2$,
  there would then have to
  exist an ideal of norm $p$ in $\cO_K$, which
  contradicts the decomposition law. Thus $a$ is not an ideal norm,
  and consequently $a$ is not the norm of an element of $\cO_K$.
\end{proof}

Our next result will also play a large role in the investigation of
number fields of the form $\Q(\sqrt[k]{m})$, $k=2^l$:

\begin{quote}
  \textbf{(3.2)} {\em Let $q \equiv 3 \bmod{4}$ be prime and $q \geq 11$;
    then there exist $a,b \in \N$ with $2q = a + b$,
    $a \equiv 5 \bmod{8}$ and $(\frac{a}{q}) = \pm 1$.}
\end{quote}

\begin{proof} We distinguish
  \begin{enumerate}
    \item[1.] $q \equiv 3 \bmod{8}$: for $q=11$ we choose $a=5$; we may
      thus assume $q \geq 19$. There then exists a $c \in
      \N$ with $2q < 16c < 3q$, and we have $q < 16c - q <
      2q$, as well as $0 < 8c - q < q/2$. Thus $8c-q$ and $16c-q$ are
      natural numbers in the interval $(0, 2q)$, and both are
      $\equiv -q \equiv 5 \bmod{8}$. Since $(\frac{2}{q}) = -1$, the
      two numbers have different quadratic character mod $q$
      (note that $16c-q \equiv 2(8c-q) \bmod{q}$), i.e.\ exactly one
      of the two numbers is a quadratic residue, and we may
      choose $a = 8c-q$ or $a = 16c-q$.
    \item[2.] $q \equiv 7 \bmod{8}$: then $q \geq 23$, and we
      choose $c \in \N$ with $5q < 16c < 6q$. As above,
      the two numbers $16c-5q$ and $16c-3q$ lie in the interval
      $(0, 2q)$, are $\equiv 5 \bmod{8}$, and have opposite
      quadratic character mod $q$. The claim now follows as in case 1.
  \end{enumerate}
\end{proof}

It now follows immediately that

\begin{quote}
  \textbf{(3.3)} {\em Let $m \equiv 2 \bmod{4}$, $m \in \N$
    square-free. Then $D(m)$ is norm-Euclidean precisely for $m=2$ and $m=6$.}
\end{quote}

\begin{proof}
  We know from \S\ 2 that $D(2)$ and $D(6)$ are norm-Euclidean.
  We have also already seen in \S\ 1 that $D(14)$ is not
  norm-Euclidean. Let us therefore assume $m \ge 22$ and that $D(m)$ is
  norm-Euclidean. Since $h(K)=1$ we may assume $m=2q$, $q \equiv 3 \bmod
  4$ prime. We now use (1.6) with $f=2q$ and the
  values of $a$, $b$ from (3.2) and show that neither $a$ nor $-b$
  is a norm from $D(m)$. Since $(2/a) = -1$ and $(2/b) = 1$ (note
  here $a \equiv 5 \bmod{8}$, $b = 2q-a \equiv 1 \bmod{8}$)
  we obtain
  $$ (\tfrac{2q}{a}) = -(\tfrac{q}{a}) = -(\tfrac aq) = -1, \text{ as well as }
     (\tfrac{2q}{b}) = (\tfrac qb) = (\tfrac bq) = (\tfrac{2q-a}{q}) =
     (\tfrac{-a}{q}) = -1. $$
  By (3.1) it now follows that neither $a$ nor $-b$ is a norm from $D(m)$,
  and (1.6) yields the claim.
\end{proof}
The case $m \equiv 3 \bmod{4}$ can be treated equally elementarily:

\begin{quote}
  \textbf{(3.4)} {\em Let $q \equiv 3 \bmod{4}$ be prime and $q \ge 23$;
    then there exist $a,b \in \N$ with $q \equiv a \pm b$,
    $(a/q)=+1$ and $a \equiv 5-q \bmod{8}$.}
\end{quote}

\begin{proof}
  For the primes $q$ with $19 < q < 64$ we give the pairs $(q,a)$
  directly: (23,6), (31,14), (43,10), (47,6), (59, 26). We may
  thus assume $q \ge 67$; again we distinguish two cases:
  1.\ $q \equiv 3 \bmod{8}$: we choose $c \in
  \N$ with $2q-16 < 64c < 3q-8$. Then the two numbers $8c+2$ and
  $8(8c+2)-2q$ lie in the interval $(0,q)$, both are $\equiv 2 = 5-q \bmod
  8$, and exactly one of them is a quadratic residue mod $q$.
  2.\ $q \equiv 7 \bmod{8}$: let $c \in \N$ with $q+16 < c <
  2q+16$; then both $8c-2$ and $2q-8(8c-2)$ lie between
  0 and $q$, are $\equiv -2 = 5-q \bmod{8}$ and
  have opposite quadratic character mod $q$.
\end{proof}

\begin{quote}
  \textbf{(3.5)} {\em Let $m \equiv 3 \bmod{4}$ be square-free; then
    $D(m)$ is norm-Euclidean precisely for $m \equiv 3, 7, 11, 19$.}
\end{quote}

\begin{proof}
  For $m \equiv 3, 7, 11, 19$ we have $M(K)<1$, so $D(m)$ is
  norm-Euclidean. Let us therefore take $m=q$, $q \equiv 3 \bmod{4}$ prime and $q \ge
  23$; by (3.4) there exist $a,b \in \N$ with $q \equiv a+b$, $a
  \equiv 5-q \bmod{8}$ and $(a/q)=+1$. With these values of $a$ and
  $-b$ we now use (1.6) with $f=q$, and it remains to show that
  $a$ and $-b$ are not norms from $D(m)$. Since $c = a/2
  \equiv 1 \bmod{2}$ would then also be a norm from $D(m)$ if $a$ were, this
  follows immediately
  from $(q/c) = (q/b) = -1$ and (3.1).
\end{proof}

The case $m \equiv 1 \bmod{4}$ can in part be treated in the same way;
thus, for example, we have

\begin{quote}
  \textbf{(3.6)} (Hofreiter\index[N]{Hofreiter} \cite{Hof34})
    {\em Let $m \equiv 3q$, $q \equiv 7 \bmod{8}$; then $D(m)$ is
    norm-Euclidean precisely for $m = 21$.}
\end{quote}

\begin{proof}
  Since $M(K) = \frac57 < 1$ for $K = \Q(\sqrt{21})$, we may assume $q
  \ge 23$ and $q \equiv 7 \bmod{8}$ prime; we then use
  (1.6) with the following values of $a$ and $b$:
  \begin{center}
    \begin{tabular}{|c|c|c|}
      \hline
      $f=q$ & $a$ & $b$ \\
      \hline
      $q \equiv 7 \bmod{24}$ & 18 & $q-18$ \\
      $q \equiv -1 \bmod{24}$, $q \equiv 5 \bmod{9}$ & 2 & $q-2$ \\
      $q \equiv -1 \bmod{24}$, $q \equiv 2 \bmod{9}$ & 8 & $q-8$ \\
      $q \equiv -1 \bmod{24}$, $q \equiv 8 \bmod{9}$ & 32 & $q-32$ \\
      \hline
    \end{tabular}
  \end{center}
  Since $q \equiv 7 \bmod{8}$, 2 and hence also $a$ is a quadratic
  residue mod $q$; on the other hand, if $a$ were a norm from $D(m)$, then
  so would be 2, which is not the case since $m \equiv 5 \bmod{8}$. Further, in the case
  $q \equiv 7 \bmod{24}$, since $b = q-18 \equiv 13 \bmod{24}$:
  \[ \left(\frac{m}{q}\right) =
  \left(\frac{3}{q}\right)\left(\frac{b}{q}\right) =
  \left(\frac{b}{3}\right)\left(\frac{b}{q}\right) =
  \left(\frac{-18}{q}\right) = \left(\frac{-2}{q}\right) = -1
  \]
  so $-b$ is not a norm from $D(m)$. Finally, $a$ is chosen in the case
  $q \equiv -1 \bmod{24}$ so that $b \equiv 3 \bmod{9}$; if $b$ were a norm
  from $D(m)$, then $c = b/3$ would also have to be one, but
  \[
  \left(\frac{m}{c}\right)
    = \left(\frac{c}{m}\right) = \left(\frac{c}{q}\right)
    =  \left(\frac{3}{q}\right) \left(\frac{3c}{q}\right)
    = \left(\frac{b}{q}\right) = -1
  \]
  shows that this is not the case.
\end{proof}

\begin{quote}
  \textbf{(3.7)} {\em Let $m=3q$, $q \equiv 3 \bmod{8}$; then $D(m)$ is
    norm-Euclidean precisely for $m=33$ and $m=57$.}
\end{quote}

\begin{proof}
  The following proof that only the values $m=33$ and $m=57$ are possible
  goes back to Schuster\index[N]{Schuster} \cite{Sch38}. If there exists
  an $s \in \N$ with $0 < 3s < q$, $s \equiv 1 \bmod{6}$ and $(s/q) =
  -1$, then $D(m)$ is not norm-Euclidean: for if we set $f=q$ and
  $a = 3s$, $b = q-3s$ for $q \equiv 1 \bmod{3}$, or $a = q-3s$, $b
  = 3s$ for $q \equiv 2 \bmod{3}$, then $(a/q) = +1$; further, e.g.,
  in the first case $s = a/3$ would also be a norm if $a$ is, but
  $(\frac{q}{s}) = (\frac{3}{s})(\frac{9}{s}) = -1$ and likewise
  $(\frac{q}{s}) = (\frac{3}{s})(\frac{9}{s}) = -1$. We now show
  the existence of such an $s$ for all primes $q \ge 43$; for
  $q<108$ we give the pairs $(q,s)$ directly: $(43, 7)$, $(59,
  13)$, $(67, 7)$, $(83, 13)$, $(107, 7)$. For the remaining $q$
  we distinguish two cases:
  \begin{enumerate}
  \item[(i)] $q \equiv 11 \bmod{24}$: we choose $c \in \N$ with
    $21q < 108c < 22q$; then the two numbers $6c-q$ and $36c-7q$
    lie in the interval $(0, q/3)$ and have opposite quadratic character
    mod $q$.
  \item[(ii)] $q \equiv 19 \bmod{24}$: let $15q < 108c < 16q$; then one of
    the two numbers $q-6c$ and $36-5q$ satisfies our requirements.
  \end{enumerate}
\end{proof}

This means that, in the case of composite $m \equiv 1 \bmod{4}$, only
those values remain to be investigated for which $m=pq$, $p \equiv q \equiv 3
\bmod{4}$ prime, $7 < p < q$ (in particular $pq \ge 77$).
Hofreiter's\index[N]{Hofreiter} proof \cite{Hof35}
that $D(77)$ is not norm-Euclidean, however, contains an error
that has apparently gone unnoticed until now (see e.g.\ van der
Linden\index[N]{Linden@van der Linden} \cite[p.\ 15 bottom]{Lin85}):
he claims that the equation $77x^2 - Y^2 = -92$ is not
solvable in $\Z$ since $(77/23) = -1$; however $(77/23) = +1$, and
$x=2, Y=20$ is a solution of this equation. That $D(77)$ is indeed
not norm-Euclidean we show with

\begin{quote}
  \textbf{(3.8)} {\em Let $m=pq$, $p \equiv q \equiv 3 \bmod{4}$ prime;
    if then there exists an $r \in \N$ with
    $(\frac{r}{p}) = -(\frac{r}{q}) = (\frac{p-r}{q})$
    and $1 \le r \le p$, then $D(m)$ is not norm-Euclidean.}
\end{quote}

\begin{proof}
  We use (1.6) with $f=p$ and $a=r$, $b=p-r$ if
  $(\frac{r}{p}) = +1$, or $a=p-r$, $b=r$ if
  $(\frac{r}{p}) = -1$. In either case we have $f = a+b$,
  $(\frac{a}{p}) = +1$, and further
  $(\frac{m}{a}) = (\frac{pq}{r}) =
  (\frac{r}{p}) (\frac{r}{q}) = -1$,
  $(\frac{m}{b}) = (\frac{b}{m}) =
  (\frac{p-r}{q})(\frac{-r}{p}) = -1$, i.e.\ by (3.1),
  neither $a$ nor $-b$ is a norm from $D(m)$.
\end{proof}

As a corollary we immediately have

\begin{quote}
  \textbf{(3.9)} {\em If $m=7q$, $q \equiv 11, 19 \bmod{24}$, then
    $D(m)$ is not norm-Euclidean. }
\end{quote}

\begin{proof}
  Set $p=7$, $r=2$ in (3.8).
\end{proof}

In particular, $D(77)$ is therefore not norm-Euclidean. To
complete the classification of norm-Euclidean quadratic number fields,
we now use the result of
Cassels\index[N]{Cassels} \cite{Cas52}, according to which quadratic
number fields $\Q(\sqrt{m})$ with $m \in \N$ can be norm-Euclidean only
if $\disc K \le 2577$. We then have a computer
search, among the 92 values of $m=pq$ with $p \equiv q \equiv 3 \bmod{4}$ prime,
$q \ge 7$, $p < q$, for a solution of (3.8), and obtain

\begin{quote}
  \textbf{(3.10)} {\em Let $m=pq$, $p \equiv q \equiv 3 \bmod{4}$;
    then $D(m)$ is norm-Euclidean precisely for $m=21, 33, 57$. }
\end{quote}

Finally we still need a criterion for prime
$m \equiv 1 \bmod{4}$; such a criterion has been given in the following form by
Erdös\index[N]{Erdos@Erd\"os} and Ko\index[N]{Ko} \cite{EK38}:

\begin{quote}
  \textbf{(3.11)} {\em Let $m = p \equiv 1 \bmod{4}$ be prime; if then there exist
    $r,s,t,u \in \N$ with $p = rs + tu$, $(r,s) = (t,u) = 1$,
    $(\frac{r}{p}) = (\frac{s}{p}) = (\frac{t}{p}) = (\frac{u}{p}) = -1$,
    then $D(m)$ is not norm-Euclidean.}
\end{quote}

\begin{proof}
  (1.6) with $a=rs$, $b=tu$.
\end{proof}

Again we have a computer search, and we find that such a representation
exists for all primes $m \equiv 1 \bmod{4}$ with $m \ge 2577$, with the
following exceptions:
$$ m = 5, 13, 17, 29, 37,
41, 61, 73, 89, 97, 109, 113, 137, 193, 241, 313, 337, 457, 601. $$

With later applications in \S\ 5 in mind, we wish to treat the case
$m \equiv 5 \bmod{24}$ without using the bound of Cassels:

\begin{quote}
  \textbf{(3.12)} {\em If $m \equiv p \equiv 5 \bmod{24}$ is prime,
    then for all $p \ge 29$ there exists a representation $p = rs +
    tu$ as in (3.11). In particular $D(m)$ for such $m$ is
    norm-Euclidean if and only if $m=5$ or $m=29$.}
\end{quote}

\begin{proof}
  We seek an $s \in \N$ with $0 < 3s < p$, $(s,3p)=1$, $s
  \equiv 1 \bmod{4}$ and $(\frac{s}{p}) = -1$. Once we have
  found such an $s$, we have $p-3s \equiv 2 \bmod{4}$, i.e.\
  $p = 3s+2u$ for some $u \in \N$, $u \equiv 1 \bmod{2}$,
  and we are done.
  If $p \equiv 2 \bmod{5}$, we can take $s=5$; for the
  other $p < 432$ we give the pairs $(p,s)$ directly: $(101, 29),
  (149, 13), (269, 29), (389, 29)$. If $p \ge 437$, there exists a $c
  \in \N$ with $2p-36 < 432c < 3p-36$. One of the two numbers
  $12c+1$ and $p-144c-12$ then satisfies the conditions on $s$ above.
\end{proof}

Evidently the proof of (3.12) rests on knowing two
``small'' quadratic non-residues mod $p$ for
$p \equiv 5 \bmod{24}$: $t=2$ and $r=3$. In the case $m \equiv 13 \bmod{24}$
only the non-residue $t=2$ is known, so an analogue of (3.12)
is much harder to prove in this case; Brauer\index[N]{Brauer}
\cite{Bra40} proved in 1940, by a complicated but elementary
method, the existence of a representation (3.11) for primes
$p \equiv 13 \bmod{24}$, $p > 109$.

Since $D(m)$ is norm-Euclidean for $m = 5, 13, 17, 29, 37, 41, 73$, but
not for $m = 61$, $89$, $97$, $109$, $113$, $137$ (see § 1, Table
\ref{App15} and § 2, Table \ref{App17}, etc.), the values
$m = 193, 241, 337, 457$ and $601$ still remain to be investigated. That the Euclidean algorithm does not hold in
$D(m)$ for these $m$ was first shown by Inkeri\index[N]{Inkeri}
\cite{Ink47}, using a method related
to (1.8) that goes back to Redei\index[N]{Redei@R\'edei}. Independently,
Chatland\index[N]{Chatland} and Davenport\index[N]{Davenport} \cite{Cha50}
then gave proofs for this, and in fact shorter ones than Inkeri's; for this
they used
Davenport's method, with which the latter had proved the bound $\disc K < 2^{14}$
for all real quadratic norm-Euclidean number fields.
The computations of Chatland and Davenport can easily be checked with a
computer (in the case $m=601$ an error has crept in; see Ennola\index[N]{Ennola} \cite{Enn58b} on this).

This determines all norm-Euclidean quadratic number fields. It
would, however, be desirable to also have, for the values $m = 193, \dots,
601$, ``simple'' proofs (e.g.\ non-Euclidean ideals
of small norm in these fields). The fact that at least some of these $m$
cannot be excluded with (1.8) is due to the rather large values
of the respective fundamental units; thus e.g.\ $u =
139\,468\,303\,679\,532 + 5\,689\,030\,769\,845\,601$ is the fundamental unit in
$D(601)$.

We now turn to the application of (1.15) and (1.16) in
real quadratic number fields. Using (1.15), only the
norm-Euclidean rings $D(m)$ for $m = 2, 3, 5, 13$ can be found, and
specifically with the trivial 1-sequence $0,1$ (for these rings $M_{r,\beta} \le
\sqrt{3/2} < 2$). It should also be remarked that for all $D(m)$ with
positive $m$ we have $\mu_1 = 2$ (i.e.\ $0, 1$ is a maximal 1-sequence) except
for $D(5)$, where $\mu_1 = 4$.

In finding 2-Euclidean number rings, (1.16) is somewhat
more successful than in the 1-Euclidean case. An example of a
comparatively long 2-sequence is
\begin{align*}
  & 0, 1, 2, 3, 4, \sqrt{2}, 1+\sqrt{2}, 2+\sqrt{2}, 3+\sqrt{2},
  4+\sqrt{2}, 2\sqrt{2}, 1+2\sqrt{2}, 2+2\sqrt{2}, 3+2\sqrt{2}, \\
  & 4+2\sqrt{2}, 3\sqrt{2}, 1+3\sqrt{2}, 2+3\sqrt{2}, 3+3\sqrt{2},
  4+3\sqrt{2}, 5+4\sqrt{2}, 6+4\sqrt{2}, 7+4\sqrt{2}
  \end{align*}
in $D(2)$; this 2-sequence shows $\mu_2 \ge 23$, while by (1.17)
$\lambda_2 \le 25$ (since the ideal $(5)$ possesses no PERS).

The following table gives 2-sequences that prove the
rings $D(m)$ for
\begin{align*}
  m & = 14, 23, 31, 43, 53, 61, 69, 77, 89, 93, 97, 113, 129, 133, 137, \\
    & \qquad 157, 161, 173, 193, 201, 213
\end{align*}
to be 2-Euclidean.

\begin{table}[ht!]
  \centering
  \begin{tabular}{|r|l|} \hline
$m$ & $\theta_1, \dots, \theta_k$  \\ \hline
14 & $0, 1, 4+\sqrt{14}, 5+\sqrt{14}$ \\
23 & $0, 1, 2, 3+\sqrt{23}, 4+\sqrt{23}, 5+\sqrt{23}, 7+2\sqrt{23},
      8+2\sqrt{23}$ \\
31 & $0, 1, 2, 5+\sqrt{31}, 6+\sqrt{31}, 7+\sqrt{31}, 11+2\sqrt{31}$ \\
43 & $0, 1, 2, 8+\sqrt{43}, 14+2\sqrt{43}, 15+2\sqrt{43}, 16+2\sqrt{43}$ \\
53 & $0, 1, 2, 3, 4$ \\
61 & $0, 1, 2, 3, 4, 4+\beta, 6+\beta, 10+\beta$ \\
69 & $0, 1, 2, 2+\beta, 3+\beta, 4+2\beta, 5+2\beta, 6+2\beta$ \\
77 & $0, 1, 2, \beta, 1+\beta, 2+\beta$ \\
89 & $0, 1, 2, 3, 4$ \\
93 & $0, 1, 2, 3+\beta, 4+\beta, 5+2\beta, 6+2\beta$ \\
97 & $0, 1, 2, 3, 7+\beta, 11+2\beta, 12+2\beta, 13+2\beta, 14+2\beta$ \\
113 & $0, 1, 2, 3, 4, 13+2\beta$ \\
129 & $0, 1, 2, 11+2\beta, 12+2\beta, 13+2\beta, 16+3\beta$ \\
133 & $0, 1, 2, 6+\beta, 7+\beta, 8+\beta$ \\
137 & $0, 1, 2, 3, 4, 16+3\beta, 17+3\beta$ \\
157 & $0, 2, 3, 5+\beta, 8+\beta, 13+2\beta, 18+3\beta$ \\
161 & $0, 1, 5+\beta, 6+\beta, 10+2\beta, 11+2\beta, 12+2\beta$ \\
173 & $0, 2, 3, 6+\beta, 12+\beta, 18+2\beta$ \\
193 & $0, 2, 4, 5+\beta, 7+\beta, 9+\beta, 27+4\beta$ \\
201 & $0, 1, 2, 3, 6+\beta, 7+\beta, 8+\beta, 13+2\beta, 14+2\beta$ \\
213 & $0, 1, 2, 3, 5+\beta, 6+\beta, 7+\beta, 8+\beta, 10+2\beta, 11+2\beta$ \\
    & $12+2\beta, 13+2\beta, 17+3\beta, 18+3\beta, 19+3\beta, 20+3\beta$ \\
\hline \end{tabular}
  \caption{2-sequences for 2-Euclidean rings $D(m)$; \qquad
    ($\beta = (1+\sqrt{m})/2$)}
\end{table}

Some of the 2-sequences listed here are longer than needed to prove
that $D(m)$ is 2-Euclidean, while others could still
be lengthened. All these sequences were computed by hand;
it is possible that further quadratic fields could be shown with
(1.16) to be 2-Euclidean.

\section*{{\sc Remarks on} \S\ 3}

A large number of mathematicians have contributed to the classification of all norm-Euclidean quadratic number fields
(for the plural formation, see Mentz\index[N]{Mentz} \cite{Men88}); the following
table is intended to give an overview:

\begin{table}[htbp]
\centering
\begin{tabular}{p{1cm} p{10cm}}
\toprule
\textbf{Year} & \textbf{Event} \\
\midrule
1832 & Gauss\index[N]{Gauss@Gauß} shows, in his investigation of
       biquadratic residues,
       that $D(-1)$ is norm-Euclidean. A similar proof
       for $D(-3)$ is found in his papers. \\
1847 & Wantzel\index[N]{Wantzel} publishes the first proof
       for $D(-3)$. \\
1886 & Legendre\index[N]{Legendre} gives the smallest solution of
       $x^2 - Ny^2 = \pm 1$, $N \le 1003$. \\
1893 & Dedekind\index[N]{Dedekind} gives a proof for $D(-1)$,
       remarking that the cases $D(m)$,
       $m = -11$, $-7$, $-3$, $-2$, $2$, $3$, $5$, $13$,
       can be treated similarly. \\
1907 & Sommer\index[N]{Sommer} gives fundamental units and class numbers of
       $\Q(\sqrt{m})$, $-97 \le m \le 101$. \\
1927 & Dickson\index[N]{Dickson} shows that the rings $D(m)$,
       $m = -11, -7, -3, -2, -1, 2, 3, 5, 13$, are norm-Euclidean and
       claims that there are no other norm-Euclidean quadratic
       rings. \\
1928 & Schaffstein\index[N]{Schaffstein} computes the class numbers of
       real quadratic number fields with prime discriminants
       $\le 12000$. \\
1933 & Perron\index[N]{Perron} establishes the Euclidean algorithm in the rings $D(m)$,
       $m = 6$, $7$, $11$, $17$, $21$, $29$, in doing so uncovering
       Dickson's error. He also suggests that possibly
       every quadratic number field of class number 1 is norm-Euclidean.
       I.\ Schur\index[N]{Schur} then informs Perron that
       $D(47)$ is not norm-Euclidean; this seems to be the first result
       in this direction. \\
1934 & In a letter to Perron, Oppenheim\index[N]{Oppenheim} shows that
       the Euclidean algorithm holds not only in the already known cases but also in $D(33)$,
       $D(37)$ and $D(41)$, while $D(23)$, $D(31)$ and $D(53)$
       are not norm-Euclidean. Independently of Oppenheim,
       Remak\index[N]{Remak} establishes the Euclidean algorithm in $D(33)$, $D(37)$ and $D(41)$. \\
1935 & Fox\index[N]{Fox} and Berg\index[N]{Berg} independently
       show that $D(m)$ in the case $m \equiv 2, 3 \bmod{4}$
       can be norm-Euclidean only for $m \equiv 2, 3, 6, 7, 11, 19$;
       Berg moreover shows that $D(19)$
       is indeed norm-Euclidean. Hofreiter proves that $D(57)$
       is norm-Euclidean and finds all norm-Euclidean $D(m)$
       with $m \equiv 21 \bmod{24}$. \\
1936 & Behrbohm\index[N]{Behrbohm} and Redei\index[N]{Redei@R\'edei}
       find all norm-Euclidean $D(m)$ with $m \equiv 5 \bmod{24}$. \\
1938 & Schuster\index[N]{Schuster} finds all norm-Euclidean $D(m)$ with
       $m \equiv 9 \bmod{24}$; Erdös\index[N]{Erdos@Erd\"os} and
       Ko\index[N]{Ko} show that there are only finitely
       many primes $p$ with $p \equiv 1 \bmod{8}$ such that $D(p)$
       is norm-Euclidean. Heilbronn\index[N]{Heilbronn} shows the same for
       composite $m \equiv 1 \bmod{8}$. \\
1940 & Brauer\index[N]{Brauer} shows that the Euclidean algorithm in
       $D(m)$, $m \equiv 13 \bmod{24}$, can hold only if
       $m \le 109$. \\
\bottomrule
\end{tabular}
\caption{Historical development of the classification of norm-Euclidean
  quadratic number fields}
\end{table}

\begin{table}[htbp]
\centering
\begin{tabular}{p{1cm} p{10cm}}
\toprule
\textbf{Year} & \textbf{Event} \\
\midrule
1942 & Redei\index[N]{Redei@R\'edei} shows that $D(73)$ is norm-Euclidean,
       finds all norm-Euclidean $D(m)$ with $m \equiv 17 \bmod{24}$
       and excludes the Euclidean algorithm in $D(m)$ for $m = 61, 89, 109, 113, 137$.
       This leaves only the residue class $m \equiv 1 \bmod{24}$
       still open. \\
1944 & Hua\index[N]{Hua} shows that for norm-Euclidean fields
       $\Q(\sqrt{p})$, $p \equiv 1 \bmod{4}$ prime, we have
       $\disc K < e^{250}$. Together with Min\index[N]{Min}
       he shows that every
       prime $p \equiv 17 \bmod{24}$ with $p > 137$ has a representation
       $p = rs + tu$, where $(r,s) = (t,u) = 1$ and
       $(r/p) = (s/p) = (t/p) = -1$. \\
1947 & Inkeri\index[N]{Inkeri} shows that there are no norm-Euclidean
       $\Q(\sqrt{m})$ with $\disc K < 5000$ except the already known ones. \\
1948 & Davenport\index[N]{Davenport} shows that $\disc K \le 16384 = 2^{14}$
       for real quadratic norm-Euclidean fields. \\
1949 & Chatland\index[N]{Chatland} finds a representation (3.11) for every
       prime $p \equiv 1 \bmod{4}$,
       $601 < p \le 16384$. Together with
       Inkeri's\index[N]{Inkeri} result (1947), this means that
       all norm-Euclidean quadratic number fields are now known. \\
1950 & Chatland and Davenport\index[N]{Davenport} treat the fields with
       $193 \le \disc K \le 601$ (differently from Inkeri and apparently in
       ignorance of his work). \\
1952 & Barnes\index[N]{Barnes} and Swinnerton-Dyer\index[N]{Swinnerton-Dyer}
       show that $D(97)$ -- contrary to a
       claim made by Redei in 1942 -- is not norm-Euclidean.
       Varnavides\index[N]{Varnavides} gives a uniform proof
       that the fields $\Q(\sqrt{m})$ with
       $m = 2$, $3$, $5$, $6$, $7$, $11$, $13$, $17$, $19$, $21$,
       $29$, $33$, $37$, $41$, $57$, $73$ are norm-Euclidean. \\
1958 & Ennola\index[N]{Ennola} shows once more that $\Q(\sqrt{m})$
       is norm-Euclidean precisely for the values given above. \\
1977 & Cooke\index[N]{Cooke} finds 20 2-stage norm-Euclidean
       quadratic number fields. \\
1985 & Johnson,\index[N]{Johnson} Queen\index[N]{Queen} and
       Sevilla\index[N]{Sevilla} show that $D(10)$ and $D(65)$
       are semi-Euclidean. \\
\bottomrule
\end{tabular}
\caption{Historical development of the classification of norm-Euclidean
  quadratic number fields (continuation)}
\end{table}

\chapter*{\S\ 4 Cubic Number Fields}
\setcounter{chapter}{4}
\addcontentsline{toc}{chapter}{\S\ 4 Cubic Number Fields}
\markboth{Euclidean Rings}{\S\ 4 Cubic Number Fields}

While a quadratic number field is uniquely determined by its
discriminant, this is no longer the case for cubic number fields:
for example, the two fields $\Q(\sqrt[3]{6})$ and $\Q(\sqrt[3]{2})$ have the
same discriminant $d=972$, but are not isomorphic. Nevertheless we
can already read off some properties of the respective field from the
discriminant alone. Since for the sign of the discriminant of a field
with $r$ real and $2s$ non-real embeddings of $K$ into $\C$ the
relation $\operatorname{sign}(\disc K) = (-1)^s$ holds, and in cubic
fields only the two possibilities $r=3$, $s=0$ and $r=1$, $s=1$ occur,
$\disc K$ is positive if and only if $K$ is totally real and has unit
rank $2$. On the other hand $\disc K$ is negative if $r=s=1$ and $K$
has unit rank $1$.

If $K$ is a number field with integral basis $\{\alpha_1, \dots,
\alpha_n\}$, and $\sigma_1, \dots, \sigma_n$ denote the embeddings of
$K$ into $\C$, then as is well known $\disc K = \det
(\sigma_i(\alpha_j))^2$. Thus $\disc K$ is the square of a number that
lies in the Galois closure $L$ of $K/\Q$. This shows that, if $K/\Q$ is
Galois, then $\sqrt{\disc K} \in R$ must hold, where $R$ as usual
denotes the ring of integers in $K$.

\medskip\noindent
{\bf Example:} If $K = \Q(\zeta_p)$ is the field of $p$-th roots of
unity for a prime $p\in\N$, then $\disc K = (-1)^{(p-1)/2}p^{p-2}$; if
we now set $p^* = (-1)^{(p-1)/2}p$, then, since $K/\Q$ is Galois,
we have $\Q(\sqrt{p^*}) = \Q(\sqrt{\disc K}) \subset \Q(\zeta_p)$.

If now $K/\Q$ is Galois and $(K:\Q)$ is odd, then
$K$ can possess no quadratic subfield, and $\sqrt{\disc K} \in
R$ then implies $\sqrt{\disc K} \in \Z$. In the case
$(K:\Q) = 3$ the converse even holds:

\begin{quote}
  {\bf (4.1)} {\em A cubic number field $K$ is Galois if and only if
    $d = \disc K$ is a square in
    $\Z$. If $K$ is not Galois, then $L =
    \Q(\sqrt{d})$ is the Galois closure of $K/\Q$.}
\end{quote}

\begin{proof}
  It suffices to show that $K(\sqrt{d})$ is Galois. With $K =
  \Q(\alpha)$, however (since $\disc
  K/\Q(\alpha) = k^2 d$ for some $k\in\Z$), we have $k\sqrt{d} =
  (\alpha-\alpha')(\alpha' - \alpha'')(\alpha'' - \alpha')$, where
  $\alpha'$, $\alpha''$, $\alpha'''$ denote the conjugates of $\alpha$.
  If $f(x) = x^3+a_2x^2+a_1x+a_0 \in \Z[x]$ is the
  minimal polynomial of $\alpha$, we easily compute that
  $\alpha' - \alpha'' = \pm k\sqrt{d}/f'(\alpha) = \pm
  k\sqrt{d}/(3\alpha^2 + 2a_2\alpha + a_1)$, as well as $\alpha' + \alpha''
  = a_2 - \alpha$. Hence $\alpha'$, $\alpha'' \in
  K(\sqrt{d})$, so $L$ is normal over $\Q$.
\end{proof}

We first wish to investigate the non-Galois case: to this end let
$k_3$ be a cubic field of discriminant $d$, $d$ not a square in $\Z$,
$k_2 = \Q(\sqrt{d})$ a quadratic number field and $K=k_2k_3$ the
normal closure of $k_3/\Q$. Then $K/\Q$ is Galois with Galois group
$S_3$ (the symmetric group of order 6). We now set up the Hilbert
subgroup series for $K/\Q$; from this table and the decomposition law
in quadratic number fields we immediately read off:

\begin{quote}
  {\bf (4.2)} {\em Let $k_3$ be a cubic number field of discriminant $d$;
    then for all primes $p \in \N$ we have:}
  \begin{align*}
    \Big(\frac dp\Big) = +1 &
    \Leftrightarrow p \text{ is completely split or inert in } k_3 \\
    \Big(\frac dp\Big) = -1 &
    \Leftrightarrow p \text{ is a product of two prime ideals in } k_3.
  \end{align*}
\end{quote}

For abelian fields this is trivially correct, because on the one hand
$d$ is then a square and thus never $(d/p) = -1$, and because
on the other hand $p$ never becomes the product of two prime ideals.

\textbf{Hilbert subgroup series for $K = k_2 k_3$}
Here
$k_3$ is a non-abelian cubic number field of discriminant $d$
$d = d_2 f^2$, $d_2 = \disc k_2$, $k_2 =
\Q(\sqrt{d}) = \Q(\sqrt{d_2})$
$K = k_2 k_3$ the normal closure of $k_3$.

\begin{table}[htbp]
\centering
\begin{tabular}{c|c|c|c|c|c|c|c|c|c|c}
\toprule
$(e,f,g)$ & $p$ & $k_2$ & $k_3$ & $\Q$ & $K_{\Z}$
       & $K_T$ & $K_1$ & $K_2$ & $K_3$ & $K$ \\
\midrule
$(1,1,6)$ & $(d/p)=+1$ & $(1,1)$ & $(1,1,1)$ & $\Q$
       & $K$ & $K$ & $K$ & $K$ & $K$ & $K$ \\
$(1,3,2)$ & $(d/p)=+1$ & $(1,1)$ & $(3)$ & $\Q$ & $k_2$
       & $K$ & $K$ & $K$ & $K$ & $K$ \\
$(1,2,3)$ & $(d/p)=-1$ & $(2)$ & $(1,2)$ & $\Q$ & $k_3$
       & $K$ & $K$ & $K$ & $K$ & $K$ \\
$(2,1,3)$ & $p=2$, $p+1$, $p\|d_2$ & $(1^2)$ & $(1,1^2)$ & $\Q$
       & $k_3$ & $k_3$ & $K$ & $K$ & $K$ & $K$ \\
& $p=2$, $2+1$, $4\|d_2$ & & & $\Q$ & $k_3$
       & $k_3$ & $k_3$ & $K$ & $K$ & $K$ \\
& $p=2$, $2+1$, $8\|d_2$ & & & $\Q$ & $k_3$
       & $k_3$ & $k_3$ & $k_3$ & $K$ & $K$ \\
$(3,1,2)$ & $p=3$, $p\|f$, $(d_2/p)=+1$ & $(1,1)$ & $(1^3)$
       & $\Q$ & $k_2$ & $k_2$ & $K$ & $K$ & $K$ & $K$ \\
& $p=3$, $9\|f$, $(d_2/3)=+1$ & & & $\Q$ & $k_2$
       & $k_2$ & $k_2$ & $K$ & $K$ & $K$ \\
$(3,2,1)$ & $p=3$, $p\|f$, $(d_2/p)=-1$ & $(2)$ & $(1^3)$
      & $\Q$ & $\Q$ & $k_2$ & $K$ & $K$ & $K$ & $K$ \\
& $p=3$, $9\|f$, $(d_2/3)=-1$ & & & $\Q$
        & $\Q$ & $k_2$ & $k_2$ & $K$ & $K$ & $K$ \\
$(6,1,1)$ & $p=3$, $3\|f$, $3\|d_2$ & $(1^2)$ & $(1^3)$
     & $\Q$ & $\Q$ & $\Q$ & $k_2$ & $K$ & $K$ & $K$ \\
& $p=3$, $9\|f$, $3\|d_2$ & & & $\Q$ & $\Q$
      & $\Q$ & $k_2$ & $k_2$ & $k_2$ & $K$ \\
\bottomrule
\end{tabular}
\caption{Hilbert subgroup series}
\end{table}

Here e.g.\ $(1,2)$ means that $p$ in $k_3$ is the product of two
prime ideals of residue degree 1 and 2, respectively, while $(1,1^2)$ symbolizes the
decomposition $(p) = P_1 P_2^2$.

A more precise investigation of the factorization of the ideal $(2)$ shows

\begin{quote}
  {\bf (4.3)} {\em Let $k_3$ be a cubic number field of discriminant $d$;
    then only the following possibilities exist:}
\[
\begin{array}{ll}
d \equiv 8 \bmod {16} : & (2) = 2_1 \ftw_2^2 \quad (\text{in } k_3) \\
d \equiv 4 \bmod {16} : & (2) = 2_1^3 \\
d \equiv -4 \bmod {16} : & (2) = 2_1 \ftw_2^2 \\
d \equiv 1 \bmod {8} : & (2) = 2_1 \ftw_2 2_3 \quad \text{{\em or }} (2) = (2) \\
d \equiv 5 \bmod {8} : & (2) = 2_1 \ftw_2
\end{array}
\]
\end{quote}

The proof rests on computing the various
relative discriminants, taking into account the above table.

If we now have a cubic field $k_3$ of discriminant $d_3 = d_2
\cdot f^2$ given and ask ourselves whether it is norm-Euclidean, then
we can work with (1.5) only when there exist in $k_3$ totally
ramified prime ideals. Since precisely the prime divisors of $f$ ramify totally,
(1.5) is useless for fields with $f=1$. We therefore wish to
proceed as follows: for fixed $d_2$, we consider all
cubic fields $k_3$ with $d_3 = d_2 f^2$; we begin with the
simplest case $d_2 = -3$.

Cubic number fields of discriminant $d_3 = -3f^2$ are precisely the
``pure'' number fields $k_3 = \Q(\sqrt[3]{m})$ (see e.g.\ Delone and
Faddeev\index[N]{Faddeev} \cite{Del64} or Cohn\index[N]{Cohn}
\cite[Exercise 18.4]{Coh78}). We note some simple properties:

Let $m = ab^2$ for certain $a,b \in \Z$ with $(a,b)=1$, $a$ and
$b$ square-free, $\theta = \sqrt[3]{ab^2}$, $\theta' =
\sqrt[3]{a^2b}$, $K = \Q(\theta) = \Q(\theta')$;
then in the case
\[
m \equiv 1 \bmod {9} : \{1, \theta, \theta'\} \text{ is an integral basis, }
\disc K = -27a^2b^2;
\]
\[
ab^2 \equiv a^2b \equiv 1 \bmod {9} : \{1, \theta,
(1+\theta+\theta')/3\} \text{ is an integral basis, } \disc K =
-3a^2b^2.
\]

The decomposition law in $k_3$ reads
\[\begin{array}{llll}
\text{a)} & p \nmid \disc K : & (p) = P_1 P_2 \quad
& \Leftrightarrow \quad p \equiv 2 \bmod {3} \\
& & (p) = P_1 P_2 P_3
\quad & \Leftrightarrow \quad p \equiv 1 \bmod {3},\ x^3 \equiv m \text{
  solvable in } \Z; \\
& & (p) = (p) \quad & \Leftrightarrow \quad
p \equiv 1 \bmod {3},\ x^3 \equiv m \text{ not solvable in } \Z; \\
\text{b)} & p \mid \disc K : & (p) = P^3 \quad
& \Leftrightarrow \quad p \neq 3 \text{ or } p=3,\ m \not\equiv \pm 1
\bmod {9} \\
& & (p) = P_1 P_2^2 \quad & \Leftrightarrow \quad p=3 \text{
  and } m \equiv \pm 1 \bmod {9}
\end{array}
\]

Finally $N_{K/\Q}(1+\theta+\theta' +
\theta^2) = r^3 + m s^3 + m^2 t^3 - 3 m r s t$.

Since Cassels\index[N]{Cassels} \cite{Cas52} has shown that a cubic
number field of negative discriminant can be norm-Euclidean only
when $|\disc K| < 420^2$, we need only consider such fields.
The tables of Nakamula\index[N]{Nakamula} \cite{Nak88}
contain all such fields; from them we read off that among
them precisely the following $\Q(\sqrt[3]{m})$ have class number 1:\label{pKh1}
$\Q(\sqrt[3]{m})$ for
\begin{align*}
  m & = 2, 3, 5, 6, 10, 12, 17, 23, 29, 33, 41, 44, 45, 46, 53, 55, 59,
        69, 71, 82, 99, \\
    & \quad 107, 116, 145, 179, 188, 197, 226, 332, 404, 575.
\end{align*}

It will turn out that among these fields precisely those with $m
= 2, 3, 10$ are norm-Euclidean; this result is due to
Cioffari\index[N]{Cioffari} \cite{Cio79}, whose work we shall follow in our
proof. The part of Cioffari's work that deals with
cubic number fields contains the following misprints:
\begin{enumerate}
\item Page 392, table to the corollary of Prop.\ 3: in the row with
  $d = 107$ in the column $f-e$ the number $93$ must stand instead of $91$;
\item Page 392, Prop.\ 5: it must read ``$N(u) \equiv (-25)^3 \equiv +10
  \bmod {53}$'' instead of ``$\dots \equiv -10 \bmod {53}$'';
  correspondingly the signs in the two following
  lines must be changed;
\item Page 393, Prop.\ 7: instead of ``belongs to $\theta (2)$; hence
  $a \equiv \theta \bmod {2}$'' it must read ``belongs to $\theta^2(2)$;
  hence $a \equiv \theta^2 \bmod {2}$'';
\item Page 395, Prop.\ 10: in the row $d=44$ of the table $c =
  5+2\cdot \theta + \theta^2/2$ must stand instead of
  $c = 5+2\theta+\varphi$.
\item Page 396, Prop.\ 12: instead of $(a_1-b_1c_1)(a_2-b_1c_1)$ it should read
  $(a_1 - \sqrt{b_ic_j})(a_2 - \sqrt{b_ic_j})$, and
  correspondingly for the two following lines.
\end{enumerate}

We now have

\begin{quote}
  {\bf (4.4)} {\em Let $(K:\Q) = n\equiv 1 \bmod {2}$, and let the
    rational primes $p_1, \dots, p_t$ be pairwise distinct and
    totally ramified in $K$; moreover let $(p_1-1,n) = \dots = (p_t-1,n) = 1$.

    If then there exists an $e \in \N$ with $1 < e < f = p_1\cdots p_t$,
    such that neither $e$ nor $f-e$ is a norm of an element of $R$,
    then $R$ is not norm-Euclidean.}
\end{quote}

\begin{proof}
  The hypotheses guarantee $(\varphi(f),n) = 1$; thus exponentiation by $n$ is an automorphism of $(\Z/f\Z)^*$ and consequently
  $e$ is an $n$-th power residue mod $f$. The claim now follows with $a = e$,
  $b=f-e$ from (1.6), since for $(K:\Q) \equiv 1 \bmod {2}$
  $b$ is a norm if and only if $-b$ is a norm.
\end{proof}

\begin{quote}
  {\bf (4.5)} {\em For the following values of $m$
    $\Q(\sqrt[3]{m})$ is not norm-Euclidean:} $m = 23, 29, 33,
    41, 46, 59, 69, 71, 82, 107, 188, 197, 226, 332, 404, 575$.
\end{quote}

\begin{proof}
  We use (4.4) with $t=1$ and $t=2$; cf.\ the table.
\end{proof}

\begin{table}[htbp]
\centering
\begin{tabular}{cccc|ccccc}
\toprule
$m$ & $p_1$ & $e$ & $f-e$ & $m$ & $p_1$ & $p_2$ & $e$ & $f-e$ \\
\midrule
59 & 59 & 7 & 52 & 23 & 3 & 23 & 13 & 56 \\
71 & 71 & 19 & 52 & 29 & 3 & 29 & 26 & 61 \\
82 & 41 & 13 & 28 & 33 & 3 & 11 & 7 & 26 \\
107 & 107 & 14 & 93 & 41 & 3 & 41 & 19 & 104 \\
179 & 179 & 7 & 172 & 46 & 2 & 23 & 7 & 39 \\
197 & 197 & 39 & 158 & 69 & 3 & 23 & 26 & 43 \\
226 & 113 & 37 & 76 & 116 & 2 & 29 & 21 & 37 \\
332 & 83 & 7 & 76 & 145 & 5 & 29 & 26 & 119 \\
404 & 101 & 28 & 73 & 188 & 2 & 47 & 37 & 57 \\
 & & & & 575 & 5 & 23 & 37 & 78 \\
\bottomrule
\end{tabular}
\end{table}

For example $a = 14$ for $m = 107 = 14 + 93$ is not a norm from $R$,
since there exists no ideal of norm $7$: the congruence $x^3 \equiv
107 \bmod 7$ is indeed not solvable, and by the decomposition law
$(7)$ remains inert in $R$. Likewise $b = 93$ is not a norm, since in
$R$ there exists no ideal of norm 31.

The possibilities still remaining, $m = 2, 3, 5, 6, 10, 12, 17,
44, 45, 53, 55, 99$, cannot be decided with (4.4) alone.
We shall therefore seek in these fields ideals $I$ with
$M(K,I) > 1$. Particularly suitable for this are divisors of the ideal
$(u-1)$, $u$ a unit in $R$. In this regard we have

\begin{quote}
  {\bf (4.6)} {\em Let $K = \Q(\sqrt[3]{m})$, $m \equiv 1 \bmod 2$
    and $3 \nmid h$, where $h = h(K)$ is the class number of $K$. If there exists
    then a prime $p \neq m$ that is totally ramified in $K$, then
    $u \equiv 1 \bmod 2$ for every unit $u$ in $R$.}
\end{quote}

Rem.: Since $(2) = \ftw_1 \ftw_2$, $\|\ftw_2\| = 4$, the divisor
$\ftw_2$ of $(u-1)$ comes into question as a candidate for $I$.

\begin{proof}
  Let $(p) = P^3$ in $R$; since $(h,3) = 1$, $P$ is a principal ideal in $R$,
  e.g.\ $P = (\pi)$ for some $\pi \in R$. Hence $\pi^3 = \pm p u^k$
  for some $k \in \Z$, where $u$ is the fundamental unit of $R$. By multiplying with
  a suitable power of $u$ we can achieve $k \in \{0,
  1\}$. If $k=0$ were the case, then $\pi = \sqrt[3]{p} \in K$ would follow, which
  is not the case since $p \neq m$. By replacing if necessary $u$ by
  $u^{-1}$ we may finally assume $k=1$. Now
  $(2) = \ftw_1 \ftw_2$ in $R$, i.e.\ $\Phi(2) = \Phi(\ftw_1)\Phi(\ftw_2) = 3$,
  hence $\pi^3 \equiv 1 \bmod 2$ since $(\pi, 2) = 1$. It then follows that $u
  \equiv u^k p \equiv \pi^3 \equiv 1 \bmod 2$.
\end{proof}

Now let $m \equiv 1 \bmod 2$, $K = \Q(\sqrt[3]{m})$; then
$(2) = \ftw_1 \ftw_2$ with $\|\ftw_1\| = 2$, $\|\ftw_2\| = 4$. Setting $\theta
= \sqrt[3]{m}$, then $\{1, \theta, 1+\theta\}$ is a prime
residue system mod $\ftw_2$, and we have $\theta^2 \equiv 1 + \theta
\bmod \ftw_2$. Thus $\ftw_2$ is norm-Euclidean if and only if the
residue classes $\theta, 1+\theta \bmod \ftw_2$ contain elements of norm $< 4$.
By (4.6) these residue classes certainly contain no
units, so that for this only elements of norm 2 or 3 come into
question.
\begin{itemize}
  \item[] $m = 5$: here $\ftw_1 = (3-\theta^2)$, $(3) = Q^3$ with $Q =
    (2-\theta)$, consequently for all elements $\alpha$ of norm 2 we have:
    $\alpha \equiv 3 - \theta^2 \equiv 1 + \theta^2 \equiv \theta
    \bmod \ftw_2$, and for all $\beta$ with $N_{K/\Q}(\beta)=3$: $\beta
    \equiv 2 - \theta \equiv \theta \bmod \ftw_2$. Thus the
    residue class $1+\theta \bmod \ftw_2$ contains only elements of norm $\ge 5$
    (even $\ge 6$, since $\theta_1 = (\theta)$).
\item[] $m = 45$: with $\alpha = \theta^2/\theta$ we have $\{1, \theta,
  \alpha\}$ an integral basis, and we find $\ftw_1 = (22 - 5\theta - \alpha)$,
  $\ftw_2 = (409 + 155\theta + 97\alpha)$, $\theta_1 = (12 - \theta -
  2\alpha)$. Here too the residue class $1+\theta \bmod \ftw_2$ contains
  no elements of norm 2 or 3 (note $\alpha \equiv
  \theta^2 \equiv 1+\theta \bmod \ftw_2$).
\item[] $m = 55$: with $\alpha = (1+\theta+\theta^2)/3$ we have $\{1,
  \theta, \alpha\}$ an integral basis, $\alpha \equiv 0 \bmod \ftw_2$ and $\ftw_1 =
  (341+71+76\alpha)$, $(3) = 3_1 3_2^2$, $3_1 = (5-3\theta + \alpha)$,
  $3_2 = (5+\theta + \alpha)$, so that the residue class $\theta \bmod
  \ftw_2$ contains no elements of norm $<4$.
\item[] $m = 99$: Let $\alpha = \theta^2/3$; then $\{1, \theta,
  \alpha\}$ becomes an integral basis, $\ftw_1^2 = (16-5\theta + \alpha)$, $(3) = 3_1^3$ and
  $16-5\theta + \alpha \equiv 1 \bmod \ftw_2$. If therefore $\pi$ is an
  element of norm $2$ in $R$, then $\pi^2 e = 16 - 5\theta + \alpha
  \equiv 1 \bmod \ftw_2$ for a unit $e$. Since $e \equiv 1 \bmod 2$
  and $\Phi(\ftw_2)=3$, we already have $\pi \equiv 1 \bmod \ftw_2$. Regardless of
  which residue class the generators of $3_1$ lie in, $\ftw_2$ cannot
  be norm-Euclidean.
\item[] $m = 6$: here $u = 1-6\theta + 3\theta^2$ is a fundamental unit of $K$,
  the ideals $(2)$ and $(3)$ are totally ramified, we have $\ftw_1 =
  (2-\theta)$, $3_1 = (3+2\theta + \theta^2)$, and the residue class
  $1+\theta \bmod \ftw_1$ contains no elements of norm $<4$.
\item[] $m = 53$: with $\alpha = (1-\theta + \theta^2)/3$ we have $\{1,
  \theta, \alpha\}$ an integral basis, and we have $(2)=\ftw_1 \ftw_2$, $\|\ftw_1\| =
  2$, $\|\ftw_2\| = 4$, and $(53) = P^3$. We claim that the
  residue class $-25 \bmod \ftw_1 P$ contains no elements of norm $<106$.
  For if $\pi \in R$, $\pi \equiv -25 \bmod P$, then,
  since $P$ is totally ramified, it follows that $N_{K/\Q}(\pi) = (-25)^3 = 10
  \bmod 53$. Thus $N_{K/\Q}(\pi) \in \{-96, -43, 10, 63\}$. Now
  the numbers $96$ and $10$ are divisible by an odd power of $2$,
  and this implies $\pi \equiv 0 \bmod \ftw_1$, in contradiction
  to $\pi \equiv -25 \equiv 1 \bmod \ftw_1$. Since further $(7)$ and $(43)$
  are inert in $K$, $63$ and $-43$ are also not norms from
  $R$. This was to be shown.
\end{itemize}
This leaves from the list on p.\ \pageref{pKh1} only the values
$m = 2, 3, 10, 12, 17, 44$; for $m = 12, 17, 44$ we can exclude the
Euclidean algorithm with the help of (1.10): to this end we give for these $m$ a
(fundamental) unit $u$, the resulting bounds
$\mu_1, \mu_2, \mu_3$ from (1.10) (initially with $k=1$), an $x \in K$
and finally $M(K,x)$:

\begin{table}[h]
\centering
\begin{tabular}{c|ccccc}
\toprule
$m$ & $u^{-1}$ & $|u|$ & $\mu_1$ & $\mu_2$ & $\mu_3$ \\
\midrule
12 & $1 + 3\theta - 3\theta^2$ & $\approx 165$ & 5.5 & 2.4 & 1.1 \\
17 & $18 - 7\theta$ & $\approx 972$ & 9.91 & 3.9 & 1.5 \\
44 & $(113 - 29 - 17\theta^2)/3$ & $\approx 4007$ & 15.89 & 4.5 & 1.28 \\
\bottomrule
\end{tabular}
\end{table}

\begin{table}[h]
\centering
\begin{tabular}{c|cc}
\toprule
$m$ & $x$ & $M(K,x)$ \\
\midrule
12 & $(12 + 15\theta + 8\theta')/18$ & $169/162$ \\
17 & $(-376 + 466\theta - 19\phi)/1028$ & $1115/1028$ \\
44 & $(8 + 12\theta + 36\theta')/59$ & $81/59$ \\
\bottomrule
\end{tabular}
\end{table}

Here $\alpha = (1 - \theta + \theta^2)/3$; moreover for
the $x$ given at $m = 44$ we have the congruence
$x \equiv 12/(5 + 2\cdot \theta + \theta') \bmod R$
(see correction 4.\ before (4.4)). We have thus shown

\begin{quote}
  {\bf (4.7)} {\em $K = \Q(\sqrt[3]{m})$ is norm-Euclidean precisely for
    $m = 2, 3, 10$.}
\end{quote}

That these fields are indeed norm-Euclidean has been shown by Godwin ($m =
2$), Taylor ($m = 3, 10$) and also Cioffari ($m = 2, 3, 10$);
this can easily be checked with a computer. Moreover, for
$K = \Q(\sqrt[3]{2})$ one can determine the Euclidean minimum $M(K) = \frac12$.

I have convinced myself (with the help of class field theory)
that of the fields of the form $d = -4f^2$ at most those with $f =
9, 11$ and $83$ are norm-Euclidean; for the first two Taylor
has established the Euclidean algorithm. For the field with $f=83$, without knowledge
of the generating equation it is rather difficult to make further statements.
I have therefore also dispensed with a precise presentation of the necessary
computations.

In the great majority of cases, in order to decide whether a cubic
field is norm-Euclidean or not, one uses (1.10); without a
table of cubic fields that gives, besides the discriminant, also fundamental unit and
class number, one therefore does not get any further.

I have tried, for cubic fields of small discriminant
($|d|<1300$), to set up such a table (see IV); besides the
generating equation, this contains the class number, as well as a
unit. It must be remarked that the table makes no claim
to completeness, and that the units given are not
necessarily fundamental. I have further determined the Euclidean minima
for the fields in the table below.

From the table one can incidentally see that the cubic number field with
$\disc K = -283$ is semi-Euclidean and has Euclidean depth $1$.
The Euclidean minimum of the field of discriminant $-87$ I have not yet been able
to determine; probably (?), however, $M(K) = 1/3$,
where this minimum, besides at the two rational
points given above, is also attained at infinitely many irrational points
(as in $\Q(\sqrt{13})$).

The ring $R = \Z[\theta]$, where $\theta$ is a root of the polynomial
$f(x) = x^3 + x^2 + 3x -1$, has discriminant $-176$ and is contained in the
field of discriminant $-44$. $R$ is a further
example of a non-integrally closed ring with $M(f) = 1$,
where $M(f)$ mod $R$ is attained at the point $(1+\theta^2)/2$.

\begin{table}[h]
\centering
$$ \begin{array}{c|c|c}
\toprule
\disc K & M(K) & C_1 \\
\midrule
\rsp  -23 & \frac{1}{5} & (1+\theta+2\theta^2)/5,\ (2+2\theta-\theta^2)/5 \\
\rsp  -31 & \frac{1}{3} & (1-\theta-\theta^2)/3 \\
\rsp  -44 & \frac{1}{2} & (1+\theta^2)/2 \\
\rsp  -59 & \frac{1}{2} & (1+\theta+\theta^2)/2 \\
\rsp  -76 & \frac{1}{2} & (1+\theta^2)/2 \\
\rsp  -83 & \frac{1}{2} & (1+\theta+\theta^2)/2 \\
\rsp  -87 & \frac{1}{3} & (1-\theta^2)/3,\ (1+\theta+\theta^2)/3 \\
\rsp -104 & \frac{1}{2} & (\theta+\theta^2)/2 \\
\rsp -107 & \frac{1}{2} & (3+\theta-3\theta^2)/8,\ (1+\theta+\theta^2)/2 \\
\rsp -108 & \frac{1}{2} & \sqrt[3]{4}/2 \\
\rsp -116 & \frac{1}{2} & (\theta+\theta^2)/2,\ (1+\theta^2)/2 \\
\rsp -135 & \frac{3}{5} & (2+2\theta-2\theta^2)/5 \\
\rsp -139 & \frac{1}{2} & (1+\theta+\theta^2)/2 \\
\rsp -140 & \frac{1}{2} & (3+2\theta-3\theta^2)/10,\ (1+\theta^2)/2 \\
\rsp -152 & \frac{1}{2} & (\theta+\theta^2)/2,\ \theta^2/2 \\
\rsp -172 & \frac{3}{4} & (\theta+\theta^2)/2 \\
\rsp -175 & \frac{3}{5} & (2-\theta+2\theta^2)/5 \\
\rsp -199 & 1 & (3+\theta-3\theta^2)/7 \\
\rsp -283 & \frac{3}{2} & (1+\theta+\theta^2)/2 \\
\rsp -307 & \frac{9}{8} & (1+\theta^2)/2 \\
\rsp -528 & \frac{5}{2} & (1+\theta^2)/2 \\
\rsp -891 & \frac{7}{2} & (1+\theta+\theta^2)/2 \\
\bottomrule
\end{array} $$
\end{table}

We have already pointed out that the bounds $\mu_1$ computed in (1.10)
can still be improved in the cubic case;
we now wish to do this. To this end we first consider the pure
cubic case $K = \Q(\theta)$, $\theta^3 = m$, and
choose the $\Q$-basis $\{1, \theta, \theta^2\}$. Then,
with the notation of (1.10) and $v:=|z|$, $w:=|z'|$, we have
$N_{K/\Q}(|z|) = v w^2 \leq k$, moreover $v \leq V$ and $w \leq W$
(where we shall still determine $V$ and $W$ suitably; in (1.10)
we had $V = W = \sqrt[3]{ku}$).

In order to estimate $f(v,w) = v+2w$ from above on the region $0 \leq v \leq V$, $0 \leq w \leq W$,
we remark that $f$ attains its maximum at most
on the boundary. If, however, $v = V$, then $w \leq k/v = k/V$ follows,
and hence $v+2w \leq V + 2\sqrt{k/V}$. In the case $w = W$, on the other hand,
correspondingly $v \leq k/W$ and $v+2w \leq 2W + k/W^2$ follow.

Thus overall we have
$$ v+2w \leq \max \{V + 2\sqrt{k/V}, 2W + k/W^2\}. $$
If we now wish to achieve $V = 2W$, we need only
choose $\sqrt[3]{4k/u^2} \leq |z| \leq \sqrt[3]{4ku}$ and then obtain
$|z'| \leq \sqrt[3]{ku/2}$, hence $V = \sqrt[3]{4ku} = 2W$ and
$$ v + 2w \leq \sqrt[3]{4ku} + \max \{2\sqrt{k/V}, k/W^2\}
   = \sqrt[3]{4ku} + \sqrt[3]{4k/\sqrt{u}} =: S $$
since $2\sqrt{k/V} = \sqrt[3]{4k/\sqrt{u}}$ and $k/W = \sqrt{4k/u}$. Now
the bounds $\mu_1 = S/3$, $\mu_2 = S/3\theta$,
$\mu_3 = S/3\theta^2$ follow (with $\theta^3 = m$).

The asymmetry in the proof leads me to suspect that these bounds are still
not best possible. However, these considerations already yield
far better bounds than (1.10); thus, for example, one finds for $m = 12, 17,
44$:

\begin{table}[h]
\centering
\begin{tabular}{c|ccc}
\toprule $m$ & $\mu_1$ & $\mu_2$ & $\mu_3$ \\ \midrule 12 & 2.92 &
1.27 & 0.59 \\ 17 & 5.25 & 2.05 & 0.80 \\ 44 & 8.41 & 2.38 & 0.68
\\ \bottomrule
\end{tabular}
\end{table}

Noting that the computing time is approximately proportional to the product
$\mu_1\mu_2\mu_3$, and comparing these bounds with those obtained somewhat
further above, the progress over (1.10) becomes
immediately clear.

In the general case we proceed by a different route than in (1.10): we
use the $\Q$-basis $\{1, \alpha, \beta\}$; instead, however, of
using the dual basis, we multiply $z$ by an
element of ``trace'' 0 and then form the ``trace'' (``trace'' stands
here in quotation marks, since e.g.\ $\alpha'-\alpha''$ does not even lie in
$K$, and we consequently may not designate
$(\alpha'-\alpha'')+(\alpha''-\alpha)+(\alpha-\alpha'')$ as
the trace of $(\alpha'-\alpha'')$). Thus, for example, we have
\[
\sum z(\alpha' - \alpha'') = r_3 m, \quad
\sum z(\beta' - \beta'') = r_2 m, \quad
\sum z(\alpha'\beta'' - \alpha''\beta') = r_1 m,
\]
where we regard the summands as elements of the normal closure $L =
K(\sqrt{d})$ of $K$ and sum over the automorphisms of
$L/\Q(\sqrt{d})$ (i.e.\ e.g.\ $\sum z' = z' + z'' +
z$). Hence
\[
m = \alpha'\beta + \alpha''\beta' + \alpha\beta'' - \alpha\beta' - \alpha'\beta'' - \alpha''\beta = \begin{vmatrix}
1 & 1 & 1 \\
\alpha & \alpha' & \alpha'' \\
\beta & \beta' & \beta''
\end{vmatrix} = \pm \sqrt{D},
\]
where as is well known $D = \disc_{K/\Q}(1, \alpha, \beta)$.
With the help of the triangle inequality we now obtain
\begin{align*}
  \|r_3\| & = \Big|\sum z(\alpha' - \alpha'') \Big|
           \le |z||\alpha' - \alpha''| + |z'||\alpha' - \alpha'|
             + |z''||\alpha - \alpha'| \\
          & \le \sqrt[3]{ku}(|\alpha' - \alpha''| + 2|\alpha - \alpha'|),
             \intertext{ and correspondingly } \\
  \|r_2\| & \le \sqrt[3]{ku}(|\beta' - \beta''| + 2|\beta - \beta'|), \quad
             \text{as well as} \\
  \|r_1\| & \le \sqrt[3]{ku}(|\alpha'\beta'' - \alpha''\beta'| +
             2|\alpha\beta' - \alpha'\beta'|).
\end{align*}
Dividing by $|m| = \sqrt{-D}$ now yields the desired bounds
for $r_1, r_2, r_3$. We record:
\begin{quote}
  {\bf (4.8)} {\em Let $K$ be a cubic field of unit rank 1,
    $u$ a unit with $u>1$, and $x_1, \dots, x_t \in K$
    points that are permuted by $u$. Further let $\{1, \alpha,
    \beta\}$ be a $\Q$-basis of $K$. If then $M(K, x_i) <
    k$, there exists a $z = r_1 + r_2\alpha + r_3\beta \in K$ with the
    properties}
    \begin{enumerate}
    \item[(a)] $z \equiv x_i \bmod {R}$ for some $j \in \{1, \dots, t\}$;
    \item[(b)] $N_{K/\Q}(z) < k$;
    \item[(c)] $|r_i| < \mu_i$ {\em for $i=1,2,3$ and}
      \begin{align*}
        \mu_1 & = \frac{\sqrt[3]{ku}(|\alpha'\beta'' - \alpha''\beta'|
          + 2|\alpha\beta' - \alpha'\beta'|)}{\sqrt{-D}}, \\
        \mu_2 & = \frac{\sqrt[3]{ku}(|\beta' - \beta''|
          + 2|\beta - \beta'|)}{\sqrt{-D}}, \\
        \mu_3 & = \frac{\sqrt[3]{ku}(|\alpha' - \alpha''|
          + 2|\alpha - \alpha'|)}{\sqrt{-D}},
      \end{align*}
      {\em where $D = \disc_{K/\Q}(1, \alpha, \beta)$. }
    \end{enumerate}
\end{quote}
With the help of (4.8), some of the cubic fields already investigated by Taylor
can be shown to be not norm-Euclidean: we
let $y$ run through the residue classes mod $(u-1)$ and determine
$M(K,x)$ for $x = y+1$. It seems promising to investigate with
this method those cubic fields for which
it is not yet known whether they are norm-Euclidean or not.
In his work from 1954,
Swinnerton-Dyer\index[N]{Swinnerton-Dyer} \cite{Swi54} has, among other things,
given the following result (without proof):
\begin{quote}
  {\em if $\{1, \theta, \theta^2\}$ is an integral basis of $\cO(\theta)$,
    $\theta^3 + 2a\theta -1 = 0$, and $a\in\N$ is sufficiently
    large, then $M(K) = (a^2 - a +1)/2$, and this minimum is attained mod
    $R$ precisely at the point $P = (\frac12, \frac12, \frac12)$.}
\end{quote}
 
We shall now show that even the following holds:
\begin{quote}
  {\bf (4.9)} {\em Let $\theta^3 + 2a\theta -1 = 0$,
    $R = \Z[\theta]$ and $t$ the absolute value of the
    norm. Then for all $a\in\N$ we have: $M(K) = M(K) = (a^2 - a +1)/2$,
    and this minimum is attained mod $R$ only at $(\frac12, \frac12, \frac12)$.}
\end{quote}
In this case incidentally $\disc_{K/\Q}(1,
\theta, \theta^2) = -(32a^3 + 27)$ and consequently $M(K) \approx
|d|^{2/3}/16^{3/2}$. For $a = 1, 2, 3$ one obtains from (4.9)
$\disc(\theta) = -59, -283, -891$; in these cases
we have already determined the first minimum by computer and
found $M(K) = \frac12$, $\frac32$, $\frac72$ in agreement
with (4.9). In the proof we may therefore restrict ourselves to $a \geq 4$.
The companion matrix of $z = r + s\theta + t\theta^2$ with respect to
the basis $\{1, \theta, \theta^2\}$ is found to be
\[
\begin{vmatrix}
r & s & t \\
t & r - 2a t & s \\
s & t - 2a s & r - 2a t
\end{vmatrix}
\]
Here the second and third rows contain, respectively, the
coordinates of $z\theta$ and $z\theta^2$. Forming the
determinant yields the norm
\[
N(z) = r^3 + s^3 + t^3 - 3rst + 2a r s^2 - 2a s^2 t - 4a t^2 + 4a^2 r t^2.
\]
Since we can shift $z$ mod $R$ so that
$-\frac12 \leq r, s, t \leq \frac12$, we immediately find
\[
|N(z)| \le \frac18 + \frac18 + \frac18 + \frac38
         + \frac a4 + \frac a4 + \frac a2 + 2a^2 t^2.
\]
For every exceptional point $z$ with $|N(z)| \ge k = (a^2 - a +1)/2$ we therefore have
$2a^2 t^2 \ge (2a^2 - 6a -1)/4$, i.e.\ $t^2 \ge (2a^2 - 6a -1)/(8a^2)$;
for all $a \ge 3$ we have the estimate
$(2a^2 - 6a -1)/(8a^2) > (\frac12 - \frac4{5a})^2$, and we see
$|t| \ge \frac12 - \frac4{5a}$.
The coefficient of $\theta^2$ in $z\theta -1$, $z$ and
$z\theta$ is, respectively, $r$, $t$ and $s$; consequently for $a \ge 3$
the estimates $|r| \le \frac12 - \frac58$ and
$|s| \le \frac12 - \frac58$ also hold. By again replacing the fundamental domain
$F = (-0.5,0.5)\times(-0.5, 0.5)\times(-0.5, 0.5)$ by, e.g.,
$F' = (0, 1) \times(-1, 0) \times(0, 1)$ (this choice of $F'$ is
guided by the consideration that $|N(\frac12, \frac12, \frac12)| = k$,
but $|N(\frac12, \frac12, \frac12)| = \frac{a^2}2 > k$), we obtain
the only exceptional set
$ S = (\frac12 - \delta, \frac12 + \delta) \times
     (-\frac12 - \delta, \frac12 + \delta) \times
     (-\frac12 - \delta, \frac12 + \delta)$ with $\delta = \frac4{5a}$.
Now let $z \in S$; we wish to estimate $|N(z)|$ once more from above,
in the hope of obtaining better bounds since $z \in S$.
By multiplying $z$ if necessary by $-1$ and carrying out a
corresponding computation, we may assume $0 \le r \le \frac12$, $0
\le s \le \frac12$, $0 \le t \le \frac12$ and find, after a little
computation, $|N(z)| \le 3a\delta + \delta^2 + \frac12 - \frac{a}2 + 2a^2 t^2$.
As above, this yields the estimate
$t^2 \ge \frac14 - \frac6{5a^2} - \frac2{5a^3}$, and for all
$a \ge 4$ this yields $|t - \frac12| \le \frac1{3a}$. As above, it now follows
also that $|r - \frac12| \le \frac1{3a}$ and $|s - \frac12| \le \frac1{3a}$.

Now we note that $z\theta = t + (r-2at)\theta + s\theta^2$;
since with $z$, $z\theta$ is also an exceptional point, we must have
$|r-2at + b + \frac12| \le \frac1{3a}$ for some $b \in \Z$.
A little computation shows $b = a-1$, and it follows by (2.1) that
$x = - \frac{(a-1)\theta + \theta^2}{\theta-1} = \frac{1- \theta + \theta^2}2$
is the only possible exceptional point of $S$.

It remains to show that $M(K,x) = \frac{a^2-a+1}2$. This is done
as in the work of Swinnerton-Dyer; we need only replace his
asymptotic estimates by explicit ones, which is not difficult.

\subsection*{2-stage norm-Euclidean number fields}

Finally we still give some examples of 2-stage
norm-Euclidean cubic fields:
\begin{enumerate}
\item[1.] $\disc K = -199$: this field is generated by a root
  $\alpha$ of the polynomial $f(x) = x^3+4x^2+x+1$, $\{1, \alpha,
  \alpha^2\}$ is an integral basis, and we have the factorizations $(2) = (2)$,
  $(3) = \fth_1 \fth_2$, $(5) = (5)$ and $(7) = \fsv_1 \fsv_2 \fsv_3$.
  Here $\fth_1 = (\alpha+1)$, $\fsv_1 = (\alpha-1)$,
  $\fsv_2 = (\alpha+2)$, $\fsv_3 = (\alpha+3)$. Since the unit $\alpha$ is
  fundamental, it follows easily that $M(K,\fsv_1) = 1$, since the residue classes
  $2, 3 \bmod \fsv_1$ contain only elements of norm $\ge 7$. According to
  Taylor we even have here $M(K) = 1$; in any case $K$ is not
  norm-Euclidean. In order to show that $K$ is at least 2-stage
  norm-Euclidean, we must find a 2-sequence of length $\ge 4$.
  Since the ideals $(2)$, $7\fsv_2$ and $\fsv_3$ possess a PERS
  and $|N_{K/\Q}(\alpha+4)| = 3$, 
  $0$, $1$, $2$, $\alpha+2$, $\alpha+3$, $\alpha+4$ is a 2-sequence of
  length $6$, and it follows that $K$ is 2-Euclidean.
\item[2.] $\disc K = -351$: $K$ is generated by a root
$\alpha$ of the polynomial $f(x) = x^3+3x+3$; we have the
  factorizations $(2) = (2)$, $(3) = P^3$, $(5) = (5)$,
  $(7) = \fsv_1\fsv_2$, $(11) = \fel_1 \fel_2 \fel_3$, and with
  $N = |N_{K/\Q}|$ we have $N(\alpha) = 3$, $N(\alpha+1) = 1$,
  $N(\alpha+2) = 11$, $N(\alpha-1) = 7$, $N(\alpha-2) = 17$,
  $N(\alpha^2-\alpha+1) = 19$, $N(\alpha^2-\alpha) = 21$.
  One can then check that $0, 1, 2, \alpha$, $\alpha+1$, $\alpha^2+1$
  is a 2-sequence, and this suffices to show that $K$ is 2-Euclidean.
  That $K$ is not norm-Euclidean is easily seen by considering the
  elements in the residue class $5 \bmod \fel_1$, where $\fel_1 = (\alpha+2)$.
\end{enumerate}

We now still wish to say something about the situation in totally real
cubic number fields. Here it is known only that the number of cyclic
norm-Euclidean fields is finite (Heilbronn\index[N]{Heilbronn}
\cite{Hei50}). Smith\index[N]{Smith} \cite{Smi69} investigated in 1969
the fields of discriminant $< 10^8$ with a computer and in doing so
established that the Euclidean algorithm exists at most when
$$ \sqrt{d} = 7, 9, 13, 19, 31, 37, 43, 61, 67, 73, 103, 109, 127, 157$$.
Moreover, he has for the
values $\sqrt{d} \le 43$ and $\sqrt{d} = 73$ determined the Euclidean minimum
$M(K)$ and shown that $K$ for $\sqrt{d} \le 67$ is
norm-Euclidean, but for $\sqrt{d} = 73$ is not. The fields with
$\sqrt{d}= 103, 109, 127, 157$ have not been investigated.

The criterion for the non-existence of the Euclidean algorithm used by
both Heilbronn and Smith\index[N]{Smith} is the following:

\begin{quote}
  {\bf (4.10)} {\em Let $p \equiv 1 \bmod 6$ be prime and $\chi$ a
    non-trivial cubic character mod $p$. If then there exist
    $r, s, t, u \in \N$ with $p = rs-tu$, $(r,s)=(t,u)=1$,
    $\chi(rs) = \chi(tu) = 1$, $\chi(r) \neq 1$, $\chi(s) \neq 1$,
    $\chi(t) \neq 1$, then $K$ is not norm-Euclidean.}
\end{quote}

\begin{proof}
  Since $\chi(rs) = 1$, $a = rs$ is a cubic residue mod $p$, while
  the other conditions imposed on $r, s, t, u$ guarantee that
  neither $rs$ nor $-tu$ is a norm from $R$. (1.6) then yields the
  claim.
\end{proof}

\begin{table}[h]
\centering
\begin{tabular}{c|ccccc}
\toprule
 & $p$ & $r$ & $s$ & $t$ & $u$ \\
\midrule
 & 79 & 2 & 5 & 3 & 23 \\
 & 97 & 2 & 11 & 3 & 25 \\
 & 139 & 2 & 3 & 7 & 19 \\
 & 151 & 2 & 43 & 5 & 13 \\
 & 163 & 2 & 11 & 3 & 47 \\
 & 181 & 2 & 11 & 3 & 53 \\
 & 193 & 2 & 37 & 7 & 17 \\
 & 199 & 2 & 41 & 9 & 13 \\
\bottomrule
\end{tabular}
\caption{Examples}
\end{table}

For the computation of the cubic character $\chi$ Smith\index[N]{Smith}
has drawn on a primitive root. If instead one uses the cubic
reciprocity law (see e.g.\ Ireland and Rosen \cite{IR82}), the
computing time can be substantially shortened, and one obtains, e.g.,

\begin{quote}
  {\bf (4.11)} {\em Let $K$ be a cyclic cubic number field of
    discriminant $d$. If then $197^2 \le d \le 2.5 \cdot 10^{11}$,
    then $K$ is not norm-Euclidean.}
\end{quote}

Since it is improbable that for a field with
$\disc K > 2.5 \cdot 10^{11}$ a representation as in
(4.10) fails to exist, the list given above of cyclic
cubic fields that are possibly norm-Euclidean is probably complete.

In 1971 Smith\index[N]{Smith} \cite{Smi71} determined Euclidean minima of some
real cubic fields and moreover established the Euclidean algorithm in numerous
such fields. In the list of discriminants of
norm-Euclidean fields a misprint has crept in: instead of
disc $K = 1994$ it should be disc $K = 1944$.

The results of Smith can be extended relatively easily: for
example, one can show that the cubic fields of discriminants
$\disc K = 2021$, $2024$, $2037$, $2101$, $2213$ are norm-Euclidean.

At this point we still wish to show how non-totally
ramified ideals can also be used for the estimation of $M(K)$. We have

\begin{quote}
  {\bf (4.12)} {\em Let $K$ be a number field of degree $n$, $R$ the ring
    of integers in $K$, and let $pR =
    \fp^{n-1}\fq$ be the factorization of the ideal
    $(p)$ in $R$. If then $\pi \equiv a \bmod \fp$ and
    $\pi \equiv b \bmod \fq$ for $a, b \in \Z$, then
    $N_{K/\Q}(\pi) = a^{n-1}b \bmod p$.}
\end{quote}

\begin{proof}
  One proceeds as in (1.1), (1.2) and (1.3); one only notes that the
  minimal polynomial $f$ of a $\pi \in P$ has the form $f(x) = x^n +
  a_{n-1}x^{n-1} + \dots + a_0$, where for $i = 0, \dots, n-2$ the
  congruence $a_i \equiv 0 \bmod p$ holds and $a_{n-1} \equiv \pi \bmod \fq$.
\end{proof}

Here one can also prove an analogue of the Eisenstein
irreducibility criterion: if $f(x) = x^n + a_{n-1}x^{n-1} +
\dots + a_0 \in \Z[x]$, $a_i \equiv 0 \bmod p$ for $0 \le i \le n-2$,
$a_0 \not\equiv 0 \bmod p^2$, and if $f$ has no linear factor $x-a$
(with $a \equiv a_{n-1} \bmod p$), then $f$ is irreducible.

Using local methods (see Ishida\index[N]{Ishida} \cite{Ish76}),
one can without difficulty generalize (4.12) to

\begin{quote}
  {\bf (4.13)} {\em Let $pR = \fp_1^{s_1} \dots \fp_g^{s_g}$
    be the factorization of $(p)$ in an algebraic number field. Then
    for all $\alpha \in R$ we have:
    $N_{K/\Q}(\alpha) = a_1^{s_1} \dots a_g^{s_g}$, where
    $a_i \in \Z$ are rational integers.}
\end{quote}

We now give some examples: Let $K$ be a cubic number field with $5
\parallel$ disc $K$; then $(5) = 5_1 5_2^2$. If now $\alpha \in
R$, $\alpha \equiv \pm 2 \bmod 5_1$ and $\alpha \equiv b \bmod 5_2$, then
from (4.12) it follows that $N_{K/\Q}(\alpha) = \pm 2 b^2 \equiv 0,
\pm 2 \bmod 5$. Since $\alpha \neq 0$, $5_1$ is not a Euclidean
ideal if there exist no elements of norm 2 or 3 in $R$. One
thus obtains
\begin{table}[h]
\centering
\begin{tabular}{c|c}
\toprule
$\disc K$ & $M(K, 5_1)$ \\
\midrule
985 & 1 \\
1345 & $7/5$ \\
3305 & $\ge 1$ \\
4345 & $\ge 1$ \\
6185 & $\ge 1$ \\
\bottomrule
\end{tabular}
\end{table}

\noindent
\textbf{Rem.}: I have not got round to determining $M(K, 5_1)$ in all
cases, or to investigating whether this consideration can also be
successfully carried out with prime ideals above other $p \equiv 1 \bmod 4$.

\section*{{\sc Remarks on} \S\ 4}
\addcontentsline{toc}{section}{{\sc Remarks on} \S\ 4}

\subsection*{1. Cubic number fields of negative discriminant}

\begin{tabular}{p{1cm} p{10cm}}
\toprule
\textbf{Year} & \textbf{Event} \\
\midrule
1892 & Markov\index[N]{Markov} gives units of $\Q(\sqrt[3]{m})$,
  $m \le 70$ \\
1940 & Delone\index[N]{Delone} and Faddeev\index[N]{Faddeev}
  publish a table of cubic fields with $-999 \le \disc K \le -23$ \\
1949 & $M(K) = 1/5$ for the field with
  $\disc K = -23$ (Prasad)\index[N]{Prasad} \\
1950 & There exist only finitely many norm-Euclidean cubic fields
  of negative discriminant (Davenport)\index[N]{Davenport} \\
1952 & Cassels\index[N]{Cassels} improves Davenport's bound: we have
  $|d| < 420^2$ \\
1954 & Swinnerton-Dyer\index[N]{Swinnerton-Dyer} proves for complex
  cubic fields with $d \le -1237$ the inequality
  $M(K) \le |d|^{2/3}/16^{3/2}$ and shows that this is best possible. \\
1957 & Godwin\index[N]{Godwin} gives a table of complex cubic fields \\
1967 & The fields with $-23 \le \disc K \le -152$
  are norm-Euclidean (Godwin) \\
1973 & Angell\index[N]{Angell} tabulates the cubic fields with $-23 \le
  \disc K \le -20000$ and gives in each case fundamental unit and class number \\
1976 & Taylor\index[N]{Taylor} determines all norm-Euclidean fields with
  $|\disc K| < 680$ \\
1979 & Cioffari\index[N]{Cioffari} finds all pure cubic
  number fields with Euclidean algorithm \\
1988 & Nakamula\index[N]{Nakamula} publishes a table of pure
  cubic fields and gives fundamental units together with class number \\
\bottomrule
\end{tabular}

\subsection*{2. Cubic number fields of positive discriminant}

\begin{tabular}{p{1cm} p{10cm}}
\toprule
\textbf{Year} & \textbf{Event} \\
\midrule
1923 & Remak\index[N]{Remak} shows $M(K) \le \sqrt{d}/8$;
  this implies in particular that the cubic field of $\disc K = 49$
  is norm-Euclidean \\
1947 & For the fields of $\disc K = 49, 81$ $M(K)$ is determined
  (Davenport)\index[N]{Davenport} \\
1950 & There exist only finitely many cyclic cubic fields with
  Euclidean algorithm; an upper bound is not given (Heilbronn)\index[N]{Heilbronn} \\
1951 & For the field of $\disc K = 148$
  $M(K)$ is determined (Clarke) \index[N]{Clarke}\\
1954 & Samet\index[N]{Samet} determines $M(K)$ for a series of cubic fields \\
1956 & Billevic\index[N]{Billevic} publishes a table of cubic
  fields with $\disc K \le 1296$ and gives two fundamental units \\
1959 & Godwin\index[N]{Godwin} and Samet\index[N]{Samet} tabulate
  real cubic fields with $\disc K \le 20.000$ \\
1969 & Smith\index[N]{Smith} investigates the cyclic fields with
  $\disc K \le 10^8$ \\
1971 & Smith determines the Euclidean minima of numerous
  cubic fields and finds all norm-Euclidean fields with $\disc
  K \le 1957$ \\
1975 & Gras\index[N]{Gras} gives a table of cyclic fields together with units and
  class number \\
1985 & Ennola\index[N]{Ennola} and Turunen\index[N]{Turunen} determine
  the real cubic fields with $\disc K < 5 \cdot 10^5$ \\
1987 & Cusick\index[N]{Cusick} and Schoenfeld\index[N]{Schoenfeld}
  tabulate the real cubic fields with $\disc K \le 6885$ and
  give class number and fundamental units \\
1988 & Llorente\index[N]{Llorente} and Quer\index[N]{Quer} compute all
  cubic fields with $0 < \disc K < 10^7$. \\
\bottomrule
\end{tabular}

\chapter*{\S\ 5 Pure Number Fields of $2$-Power Degree}
\setcounter{chapter}{5}
\addcontentsline{toc}{chapter}{\S\ 5 Pure Number Fields of $2$-Power Degree}
\markboth{Euclidean Rings}{\S\ 5 Pure Number Fields of $2$-Power Degree}

We begin with pure biquadratic number fields and to this end use
the following notation:
\[
m = ab^2 c^2 \text{ with } (a,b) = (b,c) = (c,a) = 1, \quad a,b,c
\text{ square-free,}
\]
\[
K_1 = \Q(\sqrt[4]{m}), \quad K_2 = \Q((1+i)\sqrt{m}) =
\Q(\sqrt{-4m}), \quad K_3 = \Q(i,\sqrt{ac}),
\]
\[
k_1 = \Q(\sqrt{ac}), \quad k_2 = \Q(\sqrt{-ac}), \quad
k_3 = \Q(i), \quad L = K_1(i) = K_2(i) = \Q(i,
\sqrt[4]{m}).
\]

First we ask when the polynomial $f(x) = x^4 - m$ is irreducible over
$\Q$; since by hypothesis $m$ is not a
4th power, $f$ has no roots in $\Q$, so
$f$ can factor at most into a product of two quadratic factors.
If we write
\[
f(x) = (x^2 + rx + s)(x^2 + tx + u) \text{ for } r,s,t,u \in
\Z, \text{ then it follows that}
\]
\begin{enumerate}
\item[(a)] $r=t=0, s=-u$, hence $m=s^2$: $m$ is a square in $\Z$; or
\item[(b)] $s=u$, $r=-t$, $r^2=2s$, $r^4=4s^2=-4m$: $-4m$ is a
  4th power in $\Z$.
\end{enumerate}

Thus $f$ is reducible over $\Q$ if and only if $ac=1$
or $ac=-1$, $b \equiv 0 \mod 2$. We wish
to exclude these two cases from now on.

$K_1$ then arises by adjoining a root $\alpha$ of the
polynomial $f$; the other roots of $f$ are $-\alpha$ and
$\pm i\alpha$. If $K_1/\Q$ is Galois, then $i \in
K_1$ must hold. In that case $f$ factors over $\Q(i)$ into the product
of two quadratic polynomials with coefficients from $\Z[i]$,
and comparing coefficients yields $m=-b^2$ for some $b \in
\Z$. Indeed $f(x) = (x^2 + ib)(x^2 - ib)$,
and $K_1=\Q(i,\sqrt{2b})$. So if $\pm m$ is not a square
in $\Z$, then $f$ is irreducible, $L$ is the
splitting field of $f$, and $[L:K_j] = 2$ for $j=1,2,3$.

The Galois group of $f$ respectively of $L/\Q$ is the
dihedral group
\[
D_4 = \{1, s, s^2, s^3, t, ts, ts^2, ts^3\}
\]
with the automorphisms $s: i \mapsto i$, $\alpha \mapsto i\alpha$
and $t: i \mapsto -i$, $\alpha \mapsto \alpha$.

$D_4$ has the following subgroups:
\begin{align*}
  G_1 & = \{1, t\}, \quad G_1' = s^{-1}G_1 s = \{1, ts^2\}, \\
  G_2 & = \{1, ts\}, \quad G_2' = s^{-1}G_2 s = \{1, ts^3\}, \\
  G_3 & = \{1, s^2\}, \quad G_3' = s^{-1}G_3 s = \{1, ts^2\}, \\
  G_4 & = \{1, ts\}, \quad G_4' = s^{-1}G_4 s = \{1, ts^3\},
\intertext{as well as the centre $Z(D_4)$:}
  g_1 & = \{1, s^2, t, ts^2\}, \quad g_2 = \{1, s^3, ts, ts^3\},
  \quad g_3 = \{1, s^2, s^3\}.
\end{align*}

Here $K_j$ is the fixed field of $G_j$, $k_j$ the fixed field of
$g_j$, respectively for $j=1,2,3$. Further
\begin{align*}
  K_1' & = \Q(i\sqrt[4]{m}) \text{ is the fixed field of $G_1'$ and} \\
  K_2' & = \Q(i\sqrt[4]{-4m}) = \Q((1-i)\sqrt[4]{m})
             \text{ that of $G_2'$.}
\end{align*}

We thus have the following field diagram:

\begin{center}
  \begin{tikzpicture}
    \node (Q) at (0,0) {$\Q$};
    \node (k1) at (-1, 1) {$k_1$};
    \node (k3) at ( 0, 1) {$k_3$};
    \node (k2) at ( 1, 1) {$k_2$};
    \node (K1) at (-2, 2) {$K_1$};
    \node (K1') at (-1, 2) {$K_1'$};
    \node (K3) at ( 0, 2) {$K_3$};
    \node (K2) at ( 1, 2) {$K_2$};
    \node (K2') at ( 2, 2) {$K_2'$};
    \node (L) at ( 0, 3) {$L$};
    \draw (Q) -- (k1) -- (K1) -- (L) -- (K2) -- (k2) -- (Q);
    \draw (Q) -- (k3) -- (K3) -- (L);
    \draw (k1) -- (K1') -- (L);
    \draw (k2) -- (K2') -- (L);
    \draw (k1) -- (K3) -- (k2);
  \end{tikzpicture}
\end{center}

\section*{Hilbert Subgroup Series}

$$ \footnotesize \begin{array}{ll|ccc|cccccccccccc} \hline
  p & & e & f & g & K_1 & K_2 & Z & T &
       V_1 & V_2 & V_3 & V_4 & V_5 & V_6 & V_7 & V_8 \\ \hline
  1 \ (4) & (\frac mp) = +1, (\frac mp)_4 = +1 &
    1 & 1 & 8 & (1,1,1,1) & (1,1,1,1)
       & 1 & 1 & 1 & 1 & 1 & 1 & 1 & 1 & 1 & 1 \\
                   & (\frac mp) = +1, (\frac mp)_4 = -1 &
    1 & 2 & 4 & (2,2) & (2,2)
       & G_3 & 1 & 1 & 1 & 1 & 1 & 1 & 1 & 1 & 1 \\
                   & (\frac mp) = -1 &
    1 & 4 & 2 & (4) & (4)
       & g_3 & 1 & 1 & 1 & 1 & 1 & 1 & 1 & 1 & 1 \\
    & (\frac{ac}{p}) = +1, p \mod b & 2 & 1 & 4 & (1^2,1^2) & (1^2,1^2)
       & g_3 & G_3 & & 1 & 1 & 1 & 1 & 1 & 1 & 1 \\
    & (\frac{ac}{p}) = -1, p \mod b & 2 & 2 & 2 & (2^2) & (2^2)
       & g_3 & G_3 & & 1 & 1 & 1 & 1 & 1 & 1 & 1 \\
    & p \mid ac & 4 & 1 & 2 & (1^4) & (1^4)
       & g_3 & g_3 & & 1 & 1 & 1 & 1 & 1 & 1 & 1 \\ \hline
    3 \ (4) & (\frac mp) = +1 & 1 & 2 & 4 & (1,1,2) & (2,2)
         & G_1 & 1 & 1 & 1 & 1 & 1 & 1 & 1 & 1 & 1 \\
       & (\frac mp) = +1 & 1 & 2 & 4 & (2,2) & (1,1,2)
         & G_2 & 1 & 1 & 1 & 1 & 1 & 1 & 1 & 1 & 1 \\
       & (\frac{ac}{p}) = +1, p \mod b & 2 & 2 & 2 & (1^2,1^2) & (2^2)
         & g_1 & G_3 & 1 & 1 & 1 & 1 & 1 & 1 & 1 & 1 \\
       & (\frac{ac}{p}) = -1, p \mod b & 2 & 2 & 2 & (1^2) & (1^2,1^2)
         & g_2 & G_3 & 1 & 1 & 1 & 1 & 1 & 1 & 1 & 1 \\
      & p \mid ac & 4 & 1 & 2 & (1^4) & (1^4)
         & G & g_3 & 1 & 1 & 1 & 1 & 1 & 1 & 1 & 1 \\ \hline
    p=2 & m \equiv 1 \bmod 16 & 2 & 1 & 4 & (1,1,1^2) & (1^2,1^2)
         & G_1 & G_1 & G_1 & 1 & 1 & 1 & 1 & 1 & 1 & 1\\
       & m \equiv 9 \mod 16 & 2 & 2 & 2 & (2,1^2) & (1^2,1^2)
         & g_1 & G_1 & G_1 & 1 & 1 & 1 & 1 & 1 & 1 & 1 \\
       & \frac m4 \equiv -1 \mod 16, 2| b & 2 & 1 & 4 & (1^2,1^2) & (1,1,1^2)
         & G_2 & G_2 & G_2 & 1 & 1 & 1 & 1 & 1 & 1 & 1 \\
       & \frac m4 \equiv -9 \mod 16, 2| b & 2 & 2 & 2 & (1^2,1^2) & (2,1^2)
         & g_2 & G_2 & G_2 & 1 & 1 & 1 & 1 & 1 & 1 & 1\\
       & ac \equiv 1 \mod 8, 2| b & 4 & 1 & 2 & (1^2,1^2) & (1^4)
         & g_1 & g_1 & g_1 & G_3 & G_3 & & 1 & 1 & 1 & 1\\
       & ac \equiv 5 \mod 8, 2 \nmid b & 4 & 2 & 1 & (1^2) & (1^4)
         & G & g_1 & g_1 & 1 & 1 & 1 & 1 & 1 & 1 & 1\\
       & ac \equiv 5 \mod 8, 2 \mid b & 4 & 2 & 1 & (1^2) & (1^4)
         & G & g_1 & g_1 & G_3 & G_3 & 1 & 1 & 1 & 1 & 1\\
      & ac \equiv 7 \mod 8, 2\nmid b & 2 & 1 & (1^2) & (1^4)
          & g_2 & g_2 & g_2 & G_3 & G_3 & 1 & 1 & 1 & 1 & 1 & 1\\
      & ac \equiv 3 \mod 8, 2 \nmid b & 4 & 2 & 1 & (1^4) & (2^2)
        & G & g_2 & g_2 & G_3 & G_3 & 1 & 1 & 1 & 1 & 1\\
      & ac \equiv 3 \mod 8, 2 \mid b & 4 & 2 & 1 & (1^4) & (2^2)
        & G & g_2 & g_2 & 1 & 1 & 1 & 1 & 1 & 1 & 1\\
      & 2 \mid ac & 8 & 1 & 1 & (1^4) & (1^4)
        & G & G & G & g_3 & g_3 & G_3 & G_3 & G_3 & G_3 & 1
  \end{array} $$

\section*{Integral Basis of $\Q(\sqrt[4]{m})$}

\begin{table}[h]
\centering
\begin{tabular}{llll}
\toprule
\textbf{m} & \textbf{disc $K_1$} & \textbf{disc $L$} & \textbf{integral basis}  \\
\midrule
$2 \nmid ac$ & $2^8 a^3 b^2 c^3$ & $2^{18} a^6 b^4 c^6$ & $\{1, \alpha_1, \alpha_2, \alpha_3\}$   \\
$ac \equiv 3 \mod 4$, $2 \nmid b$ & $2^8 a^3 b^2 c^3$ & $2^{18} a^6 b^4 c^6$ & $\{1, \alpha_1, \alpha_2, \alpha_3\}$  \\
$ac \equiv 3 \mod 8$, $2 \nmid b$ & $2^4 a^3 b^2 c^3$ & $2^8 a^6 b^4 c^6$ & $\{1, \alpha_1, \beta_2, \beta_3\}$  \\
$ac \equiv 7 \mod 8$, $2 \nmid b$ & $2^2 a^3 b^2 c^3$ & $2^4 a^6 b^4 c^6$ & $\{1, \alpha_1, \beta_2, \beta_3\}$ \\
$ac \equiv 1 \mod 8$, $2 \nmid b$ & $2^2 a^3 b^2 c^3$ & $2^8 a^6 b^4 c^6$ & $\{1, \alpha_1, \beta_2, \beta_3\}$ \\
$ac \equiv 5 \mod 8$, $2 \nmid b$ & $2^4 a^3 b^2 c^3$ & $2^{12} a^6 b^4 c^6$ & $\{1, \alpha_1, \beta_2, \beta_3\}$ \\
$ac \equiv 1 \mod 4$, $2 \nmid b$ & $2^4 a^3 b^2 c^3$ & $2^{12} a^6 b^4 c^6$ & $\{1, \alpha_1, \beta_2, \beta_3\}$  \\
\bottomrule
\end{tabular}

$$ \begin{array}{lll}
\toprule
 m & \beta_2 & \beta_3 \\
ac \equiv 3 \mod 8, 2 \nmid b &
  (1 + \alpha_1 + \alpha_2)/2 &  (\alpha_1 + \alpha_3)/2 \\
ac \equiv 7 \mod 8, 2 \nmid b & 
  (1 + \alpha_1 + \alpha_2)/2 & (\alpha_1 + 2\alpha_2 + \alpha_3)/2 \\
ac \equiv 1 \mod 8, 2 \nmid b &
  (1 + \alpha_2)/2 & (1 + b\alpha_1 + ab\alpha_2 + \alpha_3)/4 \\
ac \equiv 5 \mod 8, 2 \nmid b &
  (1 + \alpha_2)/2 & (\alpha_1 + \alpha_3)/2 \\
ac \equiv 1 \mod 4, 2 \nmid b &
  (1 + \alpha_2)/2 & (\alpha_1 + \alpha_3)/2 
\end{array} $$
\caption{Integral basis of $\Q(\sqrt[4]{m})$}
\end{table}

\medskip

The Galois group $D_4$ was investigated in 1987 by C.~E.\ van der Ploeg;
his work, however, contains some errors: for
instance, it is not correct that every non-normal number field of degree 4 contains a
quadratic number field, or that $Z_4$, $V_4$ and $D_4$ are the only Galois groups
of polynomials of degree 4 that can occur ($A_4$ and $S_4$ are missing;
number fields of degree 4 whose
Galois closure has one of these groups as its Galois group have,
for instance, no quadratic subfields); moreover, in the field diagrams given there,
the field we have denoted $k_3$ is missing.

To find an integral basis, the discriminant, and the decomposition law of
$\Q(\sqrt[4]{m})$, we proceed as follows: we
set up the Hilbert subgroup series for all prime ideals $P$ in
$L$, read off from it the decomposition law in $L$, $K_1$, $K_2$,
determine the contributions of the individual $P$ to the different, and finally compute
the discriminant and the integral basis of $K_1$. The
results of these computations are recorded in the tables above;
here
\begin{itemize}
\item $p \in \N$ denotes a rational prime
\item $P$ a prime ideal in $L$ above $(p)$
\item $Z = Z(P|p)$ the decomposition group
\item $T = T(P|p)$ the inertia group
\item $V_j = V_j(P|p)$ the $j$-th ramification group, $V_0 = T$
\item $e = e(P|p)$ the ramification index
\item $f = f(P|p)$ the residue degree
\item $\alpha_1 = \alpha = \sqrt[4]{m}$, \quad $m = ab^2c^2$, \quad $\alpha_2 = \sqrt{ac}$, \quad $\alpha_3 = \sqrt[4]{a^3b^2c^3}$
\item $( \cdot / p )$ the Legendre symbol
\item $( \cdot / p )_4$ the biquadratic residue symbol (if $( \cdot / p ) = +1$).
\end{itemize}

Further, e.g.\ $(1,1,2)$ means that $(p)$ factors in the corresponding field
into a product of three prime ideals, one of which has
residue degree 2 and the other two residue degree 1,
i.e.\ $(p) = P_1 P_2 P_3$, $\|P_1\| = \|P_2\| = p$, $\|P_3\| =
p^2$. Correspondingly, $(1,1,1^2)$ symbolizes the decomposition $(p) = P_1 P_2
P_3^2$, and so on.

An integral basis for $\Q(\sqrt[4]{m})$ was first given by
Ljunggren\index[N]{Ljunggren} \cite{Lju36},
though without proof. By a different route from the one described
here, Funakura\index[N]{Funakura} \cite{Fun84} has computed an
integral basis for $\Q(\sqrt[4]{m})$. The parity of the class numbers
of $K_1$, $K_2$ and $L$ has been determined by Parry\index[N]{Parry}
\cite{Par75a}; the decompositions of the ideal $(2)$ in the
intermediate fields of $L/\Q$ used in his proof are incorrect in the
two cases $m \equiv 1 \mod 8$ and $\sqrt[4]{m} \equiv -1 \mod 8$,
because these -- as a comparison with our table shows -- depend on the
residue classes mod $16$. One can, however, easily convince oneself
that this error does not affect Parry's results; these read:

\begin{enumerate}
\item Let $K = \Q(\sqrt[4]{m})$, $m \in \N$; then
$h(K)$ is odd if and only if $m$ has the following form:
\begin{align*}
  m & = 2,\ 2p^2\ (p \equiv 3 \bmod {8}), \\
  m & = p \ (p \equiv \pm 3 \bmod {8}), \\
  m & = 4p\ (p \equiv \pm 3, 7 \bmod {8}), \\
  m & = 2p\ (p \equiv 3 \bmod {8}),\ \\
  m & = 8p\ (p \equiv 3 \bmod {8}).
\end{align*}
\item Let $K = \Q(\sqrt[4]{-m})$, $m \in \N$; then
$h(K)$ is odd if and only if $m$ has the following form:
$$ m = 2,\ p (p \equiv 3 \bmod {4}),   \
   m = 4p    (p \equiv 3 \bmod {8}. $$
\item Let $L = \Q(i, \sqrt[4]{m})$; then $h(L)$ is odd if and only if
$m$ has the following form:
$$ m = 2,\ p (p \equiv 3 \bmod {8}\,),\
   m = 4p    (p \equiv 3 \bmod {4}. $$
\end{enumerate}
Here $p$ always denotes a rational positive prime.

Taking these results of Parry into account, it follows from Theorem B
of Cioffari\index[N]{Cioffari} \cite{Cio79}: if $m$ is not a square
and $K = \Q(\sqrt[4]{-m})$ is norm-Euclidean, then $m \in \{2, 3, 7,
12, 44, 67\}$. We shall show here that precisely the
fields $\Q(\sqrt[4]{-m})$ with $m=2, 3, 7, 12$ are norm-Euclidean; in doing
so we shall take a somewhat different path from Cioffari, which will
allow us to dispense with the deep result of Stark (specifically the
determination of all imaginary quadratic number fields of class number 1).

The second part of Cioffari's work contains two further small
misprints:
\begin{enumerate}
\item In Prop.\ 14 the first $K$ must be replaced by $k$;
  \item In the remark to Prop.\ 16 Cioffari writes: ``Cassels
    proved that $K$ cannot be Euclidean if $D_{K/\Q} > 5300^2$;
    hence $\Q(\sqrt[4]{-163})$ and
    $\Q(\sqrt[4]{652})$ are not Euclidean.''
\end{enumerate}

This should of course read $\Q(\sqrt[4]{-652})$.
In addition, van der Linden remarked in 1983 that Cassels
had miscalculated in his work and that the correct bound
is $D_{K/\Q} > 15170^2$. As a result,
$\Q(\sqrt[4]{-163})$ can no longer be excluded with the help of Cassels's bound.

In what follows we shall consider fields of the form $K = \Q(\sqrt[k]{-m})$,
$k = 2^l$, $l \ge 1$; here we have

\begin{quote}
  \textbf{(5.1)} {\em Let $m \in \Z$ and $\pm m$ not a square;
    with $K = \Q(\sqrt[k]{m})$ and $L = \Q(\sqrt[2k]{m})$ for some
    $k \in \N$, we then have $h(K) \mid h(L)$. }
\end{quote}

Before we prove (5.1) we recall some facts from
class field theory: if $K$ is an algebraic number field of
class number $h$ (in the wider sense), then there exists a field $L$ abelian over $K$
with the properties
\begin{itemize}
\item[] (CF--1) $L/K$ is unramified at all places (including the
  infinite ones);  
\item[] (CF--2) The Galois group $\mathrm{Gal}(L/K)$ is isomorphic to the
  ideal class group $\mathrm{Cl}(K)$; in particular $(L:K) = h =
  |\mathrm{Cl}(K)|$.
\end{itemize}

$L$ is uniquely determined by $K$, is called the Hilbert class field
of $K$ and is denoted by $\mathrm{CF}(K)$. Conversely if $L/K$ is
abelian and everywhere unramified, then $(L:K) \mid h(K)$.

\begin{proof}[Proof of 5.1]
  Because of the inclusion $K \subset L \cap \mathrm{CF}(K) \subset L$
  and $(L:K)=2$, we have either $K = L \cap \mathrm{CF}(K)$ or $L
  \subset \mathrm{CF}(K)$; in the latter case, however, $L/K$ would be
  unramified, which is not the case: since $m$ is not a square,
  there exists a rational prime $p$ dividing $m$ to an odd
  power. This prime is then totally ramified in $L/\Q$ and in particular also in
  $L/K$. Thus $L \cap \mathrm{CF}(K) = K$,
  hence $L\cdot\mathrm{CF}(K)/L$ is an abelian unramified
  extension of degree $(L\cdot\mathrm{CF}(K):L) = (\mathrm{CF}(K):L
  \cap \mathrm{CF}(K)) = (\mathrm{CF}(K):K) = h(K)$, which shows
  $h(K) \mid h(L)$.
\end{proof}

In the case $k=2$ this statement is probably due to
Chevalley\index[N]{Chevalley} \cite{Che31}. Now let $m=qs^2$ for a
square-free $q>1$: with $L = \Q(\sqrt[k]{-m})$, $k = 2^l$, we have
$K = \Q(\sqrt{-q})$ as the quadratic subfield of $L$. If
$L$ is to be Euclidean, then $h(K)=1$ must hold because of (5.1), and
as is well known this happens at most when $q=2$ or $q
\equiv 3 \bmod {4}$ is prime. If, however, $q \ge 19$ and $h(K)=1$,
then the ideals $(2)$ and $(3)$ must remain prime under the passage $K/\Q$,
since $K$ can contain no integral elements of norm $2$ or
$3$. By the decomposition law in quadratic
number fields we then have $(-q/2) = (-q/3) = -1$, and by the quadratic reciprocity law
$(2/q) = (3/q) = -1$. Since $(q)$ is totally ramified in $L/\Q$,
we can apply (1.6) with $f=q$, $a=6$, $b=q-6$ and
claim that $a=6$ is a $2^l$-th power residue mod $q$.
Because of $q \equiv 3 \bmod {4}$ and $(a/q) = (2/q)(3/q) = +1$
this is indeed the case (it suffices that $a$ be a quadratic residue
mod $q$). We now need only show that $a$ and
$-b$ are not norms from $S$. Since $K \subset L$, by the
norm tower formula it suffices to show that
$a$ and $-b$ are not norms from $\Q(\sqrt{-q})$, and this is
clear.

\begin{quote}
  \textbf{(5.2)} {\em Let $L = \Q(\sqrt[k]{-m})$ with $k = 2^l$,
    $l > 1$ and $m = qs^2$ for a square-free $q \in \N$,
    $q > 1$. Then $L$ is norm-Euclidean at most when $q \in
    \{2, 3, 7, 11\}$. If even $k \ge 4$, then in the case $q =
    11$ we must also have $s \equiv \pm 2 \bmod {3}$.}
\end{quote}

Only the statement about $q=11$ for $k \ge 4$ remains to be proved. To this end
we set $i=11$, $a=5$, $b=6$ in (1.6); since $(\frac{5}{11}) = +1$
and $11 \equiv 3 \bmod {4}$, 5 is a $2^1$-th power residue $\bmod
11$. If $s \equiv +1 \bmod {5}$, then we have
\[
\left(\frac{-m}{5}\right)^4 = \left(-\frac{11}{5}\right)^4
\left(\frac{5}{11}\right)^4 = -1,
\]
and by the decomposition law in $K = \Q(\sqrt[4]{-m})$ there exist
in $K$ (and hence a fortiori in $L$) no ideals of norm
$5$. Because $(\frac{19}{7}) = -1$, 14 is not even a
norm from $\Q(\sqrt{19})$, and by (1.6) it follows that $L$ also
cannot be norm-Euclidean for $p=19$. We now easily see
that $\Q(\sqrt[4]{11})$ is likewise not norm-Euclidean: since there exist
in $D(11)$, and hence a fortiori in $\cO_L$, no ideals of
norm 3, it suffices to show that the fundamental units of $L$
are both $\equiv 1 \bmod {(3+\sqrt{11})}$, since the ideal $I =
(3+\sqrt{11})$ is then not Euclidean because $\Phi(11) = 2$. Now
$u_1 = 10 + 3\sqrt{11}$ and
$u_2 = 881 + 477 \theta + 264 \theta^2 + 147 \theta^3$
are independent, not squares, and both $\equiv 1 \bmod I$, so the
claim follows.

\begin{quote}
  \textbf{(5.3)} {\em Let $L = \Q(\sqrt[4]{-m})$, $m \in \N$ not a
    square and free of 4th powers; then $L$ is norm-Euclidean
    precisely for $m=2,3,7,12$.}
\end{quote}

The proof that these fields are norm-Euclidean is carried out as in
\SS\ 2, 3, 4 with a computer (see \S\ 11); for $m=3$
Lakein\index[N]{Lakein} \cite{Lak72} had already established the
Euclidean algorithm, and for $m = 2, 7$ this has been done by
Cioffari\index[N]{Cioffari} \cite{Cio79}. It remains to show that $L =
\Q(\sqrt[4]{-44})$ is not norm-Euclidean.

To this end we use (1.4) with $K = \Q(\sqrt{-11})$, $B = (2)$
and $e = (1-\sqrt{-11})/2$ (since $\Phi_K(2) = 3$, all numbers in $K$
are quadratic residues $\bmod 2$). To show $M(L) \ge
\frac{5}{4}$ we must consider the elements $r \in D(-11)$
that satisfy $r \equiv e \bmod {2}$ and $N_{K/\Q}(r) <
\frac{5}{4}$. The only such $r$, however, is $r = e$,
and it suffices to show that $e$ is not a norm from $L$.
If $[\,\frac{\cdot}{\cdot}\,]$ denotes the quadratic
residue symbol in $D(-11)$, this follows immediately from
\[
[\sqrt{-44}/e] = [2\sqrt{-11}/e] = [-1/e] = (-1/3) = -1
\]
and the decomposition law in relatively quadratic number fields.

We can, however, also deduce this from the decomposition law in $L$:
we have $(3) = 3_1 3_2 3_3$ with $\|3_1\| = \|3_2\| = 3$,
$\|3_3\| = 9$; thus precisely one
of the two ideals $P_1 = \frac{1+\sqrt{-11}}2$ and $P_2 = (e)$ remains
prime under the passage to $L$, and this is $P_2$, because
$P_1 \subset L = 3_1 3_2$ for
$3_1 = (\frac{3 - \sqrt[4]{-44} + \sqrt{-11}}2)$ and
$3_2 = (\frac{3 + \sqrt[4]{-44} + \sqrt{-11}}2)$.
With (1.4) it now follows that $M(L) \ge \frac{5}{4}$, and (5.3) is proved.

Now we turn to the fields $\Q(\sqrt[k]{-m})$ with $k=2^l$,
$l>1$, $m \in \N$. To my knowledge, the only known results for
$l>2$ concerning the Euclidean algorithm in these fields are
due to Egami\index[N]{Egami} \cite{Ega79}, who,
using the results of Parry and deep estimates of
character sums, has shown:
\begin{quote}
  {\em the number of norm-Euclidean number fields of the form
    $\Q(\sqrt[4]{m})$, where $m$ is free of 4th powers and $m \ne
    2p^2$ for primes $p \equiv 3 \bmod {8}$, is finite.}
\end{quote}

We shall now show that, restricting ourselves to purely
elementary methods, we can first free ourselves of the restriction $m \ne 2p^2$,
and second can even give explicitly this finite list of fields
$\Q(\sqrt[4]{m})$ that are possibly norm-Euclidean.

\begin{quote}
  \textbf{(5.4)} {\em Let $L = \Q(\sqrt[k]{2^m p})$ for
    $k=2^l$, $l \ge 1$, $m \equiv 1 \bmod {2}$ and $p \equiv 3
    \bmod {4}$ prime. Then $L$ is norm-Euclidean at most for $p=3$.}
\end{quote}

\begin{proof}
  For $l=1$ we already know this from \S\ 3. So let $l \ge 2$; by
  Parry we may then assume $p \equiv 3 \bmod {8}$, $p \ge 11$.
  By (3.2) there then exist $a,b \in \N$ with $2p=a \cdot b$,
  $a \equiv 5 \bmod {8}$ and $(\frac{a}{p}) = +1$. Since $p \equiv 3
  \bmod {4}$ and $a \equiv 1 \bmod {2}$, $a$ is thus a $2^l$-th
  power residue $\bmod {2p}$, so we can apply (1.6) with $f=2p$
  (since the ideals $(2)$ and $(p)$ are totally ramified in $L/\Q$). With
  $K=\Q(\sqrt{2p})$, however, $a$ and $-b$ are not norms from
  $\mathcal{O}_K$, and since $K \subset L$ they are therefore a fortiori not norms
  from $\mathcal{O}_L$. Thus $L$ is not norm-Euclidean.
\end{proof}

For $k=4$ there thus remain $L=\Q(\sqrt[4]{6})$ and
$L=\Q(\sqrt[4]{24})$ to be investigated; both fields are, however,
not norm-Euclidean: for this we need only show that $M(L,I) = 9/8$
for the ideal $I$ of norm $8$ in both fields $L$. If
$\alpha = 1+\sqrt{6} \bmod I$ in $\cO_L$, it follows quickly that
$N_{L/K}(\alpha) = 7 + 2\sqrt{6} \equiv 3 \bmod {4 + 2 \sqrt{6}}$ and
$N_{L/\Q}(\alpha) \equiv 1 \bmod {8}$, where $K = \Q(\sqrt{6})$
is the quadratic subfield of $L$. The congruence in $\cO_K$
shows that $\alpha$ is not a unit (the fundamental unit of $\cO_K$ is
$5 + 2 \sqrt{6} \equiv 1\bmod {4 + 2 \sqrt{6}}$). Since $(7)$ remains
prime in $K$, we thus have $|N_{L/\Q}(\alpha)| \ge 9$, and since
$|N_{L/\Q}(3 + \sqrt{6})| = 9$, it indeed follows that $M(L,I) = 9/8$.

\begin{quote}
  \textbf{(5.5)} {\em Let $L = \Q(\sqrt[4]{2^m p})$, $k=2^l$,
    $l \ge 2$, $m \equiv 0 \bmod {2}$, $p \equiv 3 \bmod {8}$ prime. Then
    $L$ is norm-Euclidean at most for $p \le 19$.}
\end{quote}

\begin{proof}
  Let $p \ge 43$; by (3.4) there exist $a,b \in \N$ with $p=a + b$,
  $a \equiv 2 \bmod {8}$ and $(a/p) = +1$. Since $p \equiv 3
  \bmod {4}$, $a$ is thus a $k$-th power residue $\bmod {p}$. From the
  proof of (3.5) it follows that neither $a$ nor $-b$ are norms from
  $\cO_K$ for $K = \Q(\sqrt{p})$. Since $K \subset
  L$, $L$ is therefore not norm-Euclidean.
  In the case $p=19$, $m=0 \bmod {4}$, we set $a=5$, $b=14$; then
  $5 \equiv 3^4 \bmod {19}$, and $(\frac{19}{5})_4 = -1$
  shows that in $\Q(\sqrt[4]{p})$ there exist no ideals of norm 5.
  Because $(\frac{19}{7}) = -1$, 14 is not even a
  norm from $\Q(\sqrt{19})$, and by (1.6) it follows that $L$ also
  cannot be norm-Euclidean for $p=19$. We now easily see
  that $\Q(\sqrt[4]{11})$ is likewise not norm-Euclidean: since there exist
  in $D(11)$, and hence a fortiori in $\cO_L$, no ideals of
  norm 3, it suffices to show that the fundamental units of $L$
  are both $\equiv 1 \bmod {(3+\sqrt{11})}$, since the ideal $I =
  (3+\sqrt{11})$ is then not Euclidean because $\Phi(11) = 2$. Now
  $u_1 = 10 + 3\sqrt{11}$ and
  $u_2 = 881 + 477 \theta + 264 \theta^2 + 147 \theta^3$
  are independent, not squares, and both $\equiv 1 \bmod I$, so the
  claim follows.
\end{proof}

Whether $\Q(\sqrt[4]{3})$ is norm-Euclidean or not can
probably not be decided so easily. In any case $M(L) \ge
\frac{11}{12}$, as is easily established by considering the ideal $I$ of norm 12.

The cases $m=2 \mod 4$, $p=11$ and $p=19$, can be decided easily
for $k \ge 4$: we simply use (1.6) with $f = 11$, $a=5$,
$b=6$, respectively $f = 38$, $a=17$, $b=21$, and note $(\frac{44}{3})_4 =
-1$, $(\frac{11}{3})_4 = -1$, respectively $(\frac{76}{17})_4 = -1$,
$(\frac{12}{7})_4 = -1$.

\begin{quote}
  \textbf{(5.6)} {\em Let $L = \Q(\sqrt[k]{p})$, $k=2^l$, $l \ge 2$,
    $p \equiv 5 \bmod {8}$ prime. Then $L$ is norm-Euclidean at most
    for the values $p = 5, 13, 37, 61$.}
\end{quote}

\begin{proof}
  From the works of Behrbohm\index[N]{Behrbohm}
  and R\'edei\index[N]{Redei@R\'edei} \cite{BR36}; see also \S\ 3) and
  Brauer\index[N]{Brauer} \cite{Bra40}, it follows that for primes $p \equiv 5
  \bmod {8}$, $p > 109$, a representation $p = rs+tu$ exists with
  $r,s,t,u \in \N$, $(r,s) = (t,u) = 1$ and $(\frac{r}{p}) =
  (\frac{s}{p}) = (\frac{t}{p}) = -1$ (see 3.12.). Since $-1$ is a
  biquadratic non-residue mod $p$, either $rs$ or $tu$ is a
  biquadratic residue mod $p$ (note $rs = -tu \mod p$). So let
  $rs$ be, without loss of generality, a biquadratic residue mod $p$.
  Since neither $a=rs$ nor  $-b=-tu$ are norms from $\Q(\sqrt[4]{p})$,
  $L$ is by (1.6) not norm-Euclidean.
  \begin{enumerate}
    \item[] $p=29$: here $24 \equiv 4^4 \bmod {p}$, $29 = 24 +
      5$. Since $(3)$ remains prime in $\Q(\sqrt[4]{29})$, $24$ is certainly
      not a norm from $\mathcal{O}_L$. Further, $5$ is not an ideal norm,
      since $(\frac{29}{5})_4 = -1$;
  \item[] $p=53$: $53 = 15 + 38$, $8^4 = 15 \bmod {53}$;
  \item[] $p=109$: $109 = 104 + 5$, $(\frac{104}{109})_4 =
    (\frac{5}{109})_4 = +1$.
  \end{enumerate}
\end{proof}

In the case $l=2$ we can also exclude $p=37$: to this end we note
that the ideal $I = (\frac{5 - \sqrt{37}}2)$ of norm $3$ remains prime in $L$,
since $[\frac{\sqrt{37}}{I}] = [\frac 5I] = (\frac{5}{3}) =
-1$. We then set $B = (2)$, $e = (\frac{5 - \sqrt{37}}2)$, $k =
\frac74$ in (1.4) and consider all $r \in \mathcal{O}_K$,
$K=\Q(\sqrt[4]{37})$, with $r \equiv e \bmod {2}$ and $|N_{K/\Q}(r)| <
7$. Since $u = 6+\sqrt[4]{37}$ is the fundamental unit of $K$ and $u \equiv 1 \bmod
{2}$, we certainly have $|N_{K/\Q}(r)| \in \{3, 5\}$, and since $(5)$ is
prime in $K$, we must even have $|N_{K/\Q}(r)| = 3$. Thus $r =
u^{e'}$ for $e' = (5+\sqrt[4]{37})/2$. This contradicts the
congruence $r \equiv e \bmod {2}$, while the possibility $(r) = (e)$
fails because $I$ remains prime in $L$.

\begin{quote}
  \textbf{(5.7)} {\em Let $L = \Q(\sqrt[k]{2^mp})$,
    $k=2^l$, $m \equiv 0 \bmod {2}$, $p \equiv 5 \bmod {8}$ prime. Then
    $L$ is norm-Euclidean at most for $p = 5, 13, 29, 37, 61, 109$.}
\end{quote}

\begin{proof}
  As (5.6).
\end{proof}

In the case $k=4$ the fields $\Q(\sqrt[4]{4 \cdot 37})$
and $\Q(\sqrt[4]{4 \cdot 61})$ can easily be excluded: in the
first case we use the intermediate field $K = \Q(\sqrt{37})$
and proceed exactly as above for
$\Q(\sqrt[4]{37})$. Using (1.6) with $f=61$, $a=56$
and $b=5$ settles the second case.

Finally we must still treat the case where no rational
prime except $p=2$ is totally ramified in $L/\Q$, so that
the criterion (1.6) cannot be applied. Surprisingly, however,
we find

\begin{quote}
  \textbf{(5.8)} {\em Let $L = \Q(\sqrt[4]{2p^2})$,
    $p \equiv 3 \bmod {8}$ prime. Then $L$ is not norm-Euclidean.}
\end{quote}

\begin{proof}
  $L$ has two fundamental units, $u=1+\sqrt{2}$ and $v$;
  we show that $v \equiv 1 \bmod {\sqrt{2}}$, so that, since
  $\Phi_L(1) = 2$ and $\|1\| = 4$, the ideal $I = (\sqrt{2})$ is not
  Euclidean (since there exist no elements of norm 3).
  Since the ideal $(p)$ is ramified in $L/K$, there exists a $\pi \in \cO_L$ with
  $p = \pm \pi^2 u^m v^n$. If $n \equiv 0 \bmod {2}$ were the case,
  $\sqrt{\pm p} \in L$ or $\sqrt{\pm u} \in L$ would follow,
  which is not the case. Now, $p \equiv 1 \bmod {2}$ and
  $u \equiv 1 \bmod {\sqrt{2}}$, so $v^n \equiv 1
  \bmod {\sqrt{2}}$ results. Since $\Phi_L(1) = 2$, we have $v^2 \equiv 1
  \bmod {I}$, and $v \equiv 1 \bmod {I}$ then follows immediately.
\end{proof}

Summarizing, we can now establish that $L =
\Q(\sqrt[4]{m})$ is norm-Euclidean at most for the values $m = 2, 3, 5, 12, 13,
20, 28, 37, 52, 61, 116, 436$. It is to be assumed
that some non-norm-Euclidean fields still remain among these.
With the help of a computer we can further confirm that
the two fields $\Q(\sqrt[4]{2})$ and
$\Q(\sqrt[4]{5})$ are norm-Euclidean.

\section*{{\sc Remarks on} \S\ 5}
\addcontentsline{toc}{section}{{\sc Remarks on} \S\ 5}
\medskip

\begin{tabular}{p{1cm} p{10cm}}
\toprule
\textbf{Year} & \textbf{Event} \\
\midrule
1936 & Ljunggren gives integral bases and investigates units
  of pure biquadratic number fields \\
1975 & Parry gives all pure biquadratic number fields of odd
  class number \\
1979 & Egami shows that there exist only finitely many norm-Euclidean
  fields of the form $\Q(\sqrt[4]{m})$, provided that not $m =
  2p^2$, $p \equiv 3 \bmod {8}$ prime. The case $m = 2p^2$ remains
  unfinished. \\
     & Cioffari shows $m \in \{2, 3, 7, 12, 44, 67\}$ or $m = 2p^2$ for
  norm-Euclidean fields $\Q(\sqrt[4]{m})$. The case $m = 2p^2$ can
  however be settled with the results of Parry (which Cioffari
  apparently did not know). \\
1980 & Parry develops a genus theory for pure
  biquadratic fields \\
1984 & Funakura gives an integral basis for pure biquadratic fields \\
1987 & Buchmann computes a fundamental unit of rings of the form
  $\Z[\sqrt[4]{-m}]$ \\
\bottomrule
\end{tabular}

\chapter*{\S\ 6 Bicyclic Biquadratic Number Fields}
\setcounter{chapter}{6}
\addcontentsline{toc}{chapter}{\S\ 6 Bicyclic Biquadratic Number Fields}
\markboth{Euclidean Rings}{\S\ 6 Bicyclic Biquadratic Number Fields}

Let $m$ and $n$ be square-free integers and $K = \Q(\sqrt{m},\sqrt{n})$.
Then $K/\Q$ is a Galois extension of $\Q$ of degree 4 with abelian Galois
group $G = \mathrm{Gal}(K/\Q) = V_4$ (Klein four-group). $K$ contains
precisely three non-trivial intermediate fields: $k_1 = \Q(\sqrt{m})$,
$k_2 = \Q(\sqrt{n})$, and $k_3 = \Q(\sqrt{mn})$.

If $l = (m,n)$, we write $m = lm'$, $n = ln'$, so that
$k_3 = \Q(\sqrt{m'n'})$ with square-free $m'n'$. Every $\alpha \in K$
can be represented in the form $\alpha = r + s\sqrt{m} + t\sqrt{n} +
u\sqrt{mn}$, $r,s,t,u \in \Q$; if we denote the
relative norm from $K$ to $k_i$ by $N_i$, we easily find
\[
N_1(\alpha) = a + b\sqrt{m}; \quad N_2(\alpha) = c + d\sqrt{n}; \quad
N_3(\alpha) = e + f\sqrt{mn} \quad \text{with}
\]
\[
a = r^2 + m s^2 - n t^2 - m n u^2, \quad b = 2(rs - n t u),
\]
\[
c = r^2 - m s^2 + n t^2 - m n u^2, \quad d = 2(r t - m s u),
\]
\[
e = r^2 - m s^2 - n t^2 + m n u^2, \quad f = 2(r u - s t).
\]

The norm tower formula immediately yields
\[
N_{K/\Q}(\alpha) = a^2 - m b^2 = c^2 - n d^2 = e^2 - m n f^2.
\]

If we denote by $T_i$ the relative trace from $K$ to $k_i$, then
an $\alpha \in K$ is integral if and only if $N_i(\alpha)$ and
$T_i(\alpha)$ are integral for some $i \in \{1,2,3\}$. After a somewhat
lengthy but elementary computation (see e.g.\ Williams \cite{Wil70}),
we obtain the following: $\alpha = (r + s\sqrt{m} + t\sqrt{n} + u\sqrt{mn})/g
\in K$ is integral if and only if $r, s, t, u, g \in \Z$
satisfy the conditions (x):

\begin{center}
\begin{tabular}{|c|c|c|c|}
\hline
$m \bmod {4}$ & $n \bmod {4}$ & $d$ & (x) \\
\hline
1 & 1 & 4 & $r \equiv s \equiv t \equiv u \bmod {2}$,
$r+s+t+u \equiv 0 \bmod {4}$ \\
1 & 2 & 2 & $r \equiv s \bmod {2}$, $t \equiv u \bmod {2}$ \\
1 & 3 & 2 & $r \equiv u \bmod {2}$, $s \equiv t \bmod {2}$ \\
2 & 3 & 2 & $r \equiv t \equiv 0 \bmod {2}$, $s \equiv u \bmod {2}$ \\
\hline
\end{tabular}
\end{center}

Here $m$ and $n$ are given only modulo 4; it is readily seen that,
because of $\Q(\sqrt{m},\sqrt{n}) =
\Q(\sqrt{m},\sqrt{mn}) = \Q(\sqrt{n},\sqrt{m})$
etc.\, all possible cases are covered.

If we set $d = \mathrm{disc}\, K$, $d_i = \mathrm{disc}\, k_i$
($i=1,2,3$), a straightforward computation shows that in every case the
relation $d = d_1 d_2 d_3$ holds (a proof that does not require knowledge
of the integral basis proceeds via the ``conductor-discriminant
formula''). The decomposition law in $K$ then reads
\begin{align*}
(p) & = \fp_1 \fp_2 \fp_3 \fp_4,
\quad \text{if } \big(\frac{d_1}{p}\big) =
\big(\frac{d_2}{p}\big) = \big(\frac{d_3}{p}\big) = +1; \\
(p) & = \fp_1 \fp_2, \quad \text{if }
\left(\frac{d_1}{p}\right) = \left(\frac{d_2}{p}\right) = +1,
\left(\frac{d_3}{p}\right) = -1; \\
(p) & = \fp_1^2 \fp_2^2, \quad \text{if } p \nmid d_1,
   p \nmid d_2 \text{ and } \left(\frac{d_3}{p}\right) = +1; \\
   (p) & = \fp^2, \quad \text{if } p \nmid d_1, p \nmid d_2 \text{ and }
   \left(\frac{d_3}{p}\right) = -1, \\
(p) & = \fp^4, \quad \text{if } p=2, \, 2\mid d_1, \, 2\mid
d_2, \, 2\mid d_3.
\end{align*}

Here $(\cdot / p)$ denotes the Kronecker symbol, which coincides with
the Legendre symbol for odd $p$, and for $p=2$ and $d \equiv 1
\bmod {4}$ is defined by $(d/2) = (2/d)$.

The table of higher ramification groups for $K$ given by
Cohn\index[N]{Cohn} \cite[p.\ 202]{Coh78} contains two errors. To
correct these, we adopt the notation used there and distinguish the
following cases:

\begin{enumerate}
\item[(I)] $(d_1/p) = (d_2/p) = (d_3/p) = +1$
\item[(II)] $(d_1/p) = (d_2/p) = +1$, $(d_3/p) = -1$
\item[(III)] $p \equiv 1 \bmod {2}$, $p \nmid d_1$, $p \nmid d_2$, $(d_3/p) = +1$
\item[(IV)] $p \equiv 1 \bmod {2}$, $p \nmid d_1$, $p \nmid d_2$, $(d_3/p) = -1$
\item[(V)] $p=2$, $d_1 \equiv d_2 \equiv 12 \bmod {16}$, $d_3 \equiv 1 \bmod {8}$
\item[(VI)] $p=2$, $d_1 \equiv d_2 \equiv 12 \bmod {16}$, $d_3 \equiv 5 \bmod {8}$
\item[(VII)] $p=2$, $d_1 \equiv d_2 \equiv 8 \bmod {16}$, $d_3 \equiv 1 \bmod {8}$
\item[(VIII)] $p=2$, $d_1 \equiv d_2 \equiv 8 \bmod {16}$, $d_3 \equiv 5 \bmod {8}$
\item[(IX)] $p=2$, $d_1 \equiv d_2 \equiv 8 \bmod {16}$, $d_3 \equiv 12 \bmod {8}$
\end{enumerate}

Cohn gives the same subgroup series for cases V.\ and VII.\ (respectively
VI.\ and VIII.); this would, however, imply that $d = \mathrm{disc}\, K$
is divisible by the same power of 2 in the cases V.\ and VII.\
(respectively VI.\ and VIII.), which is not the case (recall $d = d_1 d_2 d_3$).

The correct subgroup series are given by the following table:

\begin{center}
$$ \begin{array}{c|cccccccc}
\text{Type} & \Q & K_{\mathfrak{z}} & K_{\mathrm{T}} & K_1 & K_2 & K_3 & K_4 & K \\
\midrule
\text{I}   & \Q & K & K & K & K & K & K & K \\
\text{II}  & \Q & k_3 & K & K & K & K & K & K \\
\text{III} & \Q & k_3 & k_3 & K & K & K & K & K \\
\text{IV}  & \Q & \Q & k_3 & K & K & K & K & K \\
\text{V}   & \Q & k_3 & k_3 & k_3 & K & K & K & K \\
\text{VI}  & \Q & \Q & k_3 & k_3 & K & K & K & K \\
\text{VII} & \Q & k_3 & k_3 & k_3 & k_3 & K & K & K \\
\text{VIII}& \Q & \Q & k_3 & k_3 & k_3 & k_3 & K & K \\
\text{IX}  & \Q & \Q & \Q & \Q & k_3 & k_3 & K & K \\
\end{array} $$
\end{center}

From now on let $m$ and $n$ be negative; then $K$ is totally imaginary,
and its unit group has rank 1. If $u$ denotes the fundamental unit of
$k_3$ (the largest real subfield of $K$), there are two possibilities:

\begin{enumerate}
\item[1.] $u$ is also the fundamental unit of $K$; in this case we set $q(K)=1$;
\item[2.] $u$ is not the fundamental unit of $K$; then there exists a root of unity
  $\zeta_K$ such that $u = \zeta \cdot e^2$. In this case $e$ is
  the fundamental unit of $K$, and we write $q(K)=2$.
\end{enumerate}

The constant $q(K)$ is called the \textbf{unit index of $K$} and
can be determined fairly easily:

\begin{enumerate}
\item[(a)] $K$ contains $\Q(i)$: then $q(K)=2$ if and only if the
  ideal $(2)$ in $k_3$ is the square of a principal ideal; if we write
  $(2) = (\alpha)^2$, then $e = \alpha^{-1}$ is the fundamental unit of $K$.
\item[(b)] $K$ does not contain $\Q(i)$: then $q(K)=2$ if and only if
  the ideal $(m_1)$ in $k_3$ is the square of a principal ideal; with
  $(m_1) = (\alpha)^2$ we have that $u = \alpha / \sqrt{m_1}$ is the
  fundamental unit of $K$.
\end{enumerate}

\noindent
\textbf{Rem.:} Part (a) also appears in Cohn\index[N]{Cohn} \cite{Coh78}
as part of Theorem 19.8; there, however, the condition that $(2)$
be the square of a principal ideal is missing, although it is used in the proof.

From the analytic class number formula there follows fairly easily the
elegant relation $H = q(K) h_1 h_2 h_3$, where $H, h_1, h_2, h_3$
respectively denote the class numbers of $K, k_1, k_2, k_3$ and $K \ne
\Q(\sqrt{-1},\sqrt{2})$. From this we in turn obtain $h_3 \mid H$
for the class number $h_3$ of the real quadratic subfield (this is
a special case of the more general relation $h' \cdot h$ that holds for
CM-fields $K$ with maximal real subfield $K^*$; see
Washington\index[N]{Washington} 1982 on this). If $K$ is Euclidean, then
in particular $H=1$, so $h_3 =1$ must also hold.

Proofs of these and further results can be found in
Kuroda\index[N]{Kuroda} \cite{Kur43,Kur50},
Kubota\index[N]{Kubota} \cite{Kub56},
Wada\index[N]{Wada} \cite{Wad66},
Fröhlich\index[N]{Frohlich@Fröhlich} \cite{Fro83}, as well as in the
classical work of Hasse\index[N]{Hasse} \cite{Has85}.

In her paper \cite{Sau73} from 1972/73,
J.~Sauvageot\index[N]{Sauvageot} writes:

\begin{quote}
  {\em On sait quels corps quadratiques sont euclidiens. La question
    reste ouverte pour les corps biquadratiques. Elle devrait être
    bientôt (?) résolue pour ceux d'entre eux qui sont bicyclique
    imaginaires, \dots''.}
\end{quote}

In fact, we have the following (by $D(m,n)$ we denote the ring of
integers of $\Q(\sqrt{m},\sqrt{n})$):

\begin{quote}
  \textbf{(6.1)} {\em Let $m$ be negative; then precisely the following
    rings $D(m,n)$ are norm-Euclidean:
    \begin{itemize}
    \item[] $D(-1,n)$ for $n = 2, 3, 5, 7$;
    \item[] $D(-2,n)$ for $n = -3, 5$;
    \item[] $D(-3,n)$ for $n = 2, 5, -7, -11, 17, -19$, and
    \item[] $D(-7,5)$.
      \end{itemize} }
\end{quote}

For fields $\Q(\sqrt{m}, \sqrt{n})$ with $m \equiv n \equiv 1
\bmod {2}$, $m < 0$, J.\ Sauvageot had already attempted a classification.
Despite an error (in the case $m = n \equiv 1 \bmod {4}$),
the conjectures she states at the end of her paper have been confirmed;
there she writes:

\begin{quote}
  {\em \dots les seuls survivants de ces éliminatoires sont
    \begin{itemize}
      \item $\Q(i\sqrt{3}, i\sqrt{7}); \Q(i\sqrt{3}, i\sqrt{11});
        \Q(i\sqrt{3}, \sqrt{5}) $ dont LAKEIN a montré qu'ils le
        sont.
      \item[] $\Q(i\sqrt{3}, i\sqrt{19})$; $\Q(i\sqrt{3}, \sqrt{17})$
        dont j'espère montrer qu'ils le sont.
      \item[] $\Q(i\sqrt{3}, i\sqrt{43})$, que je crois non-euclidien
        et $\Q(i\sqrt{7}, \sqrt{5})$ dont je ne sais rien.
    \end{itemize}   
    Les mots ``j'espère'' and \"je crois'' sont conséquence
    d'explorations du problème sur ordinateurs que j'exposerai \dots
    si elles aboutissent.}
\end{quote}

Before proving (6.1), we give a table of the first minima of certain
$D(m,n)$, together with a set $C_1$ of points at which these minima are
attained.

\begin{table}[ht!]
$$ \begin{array}{cr|c|l}
\rsp m & n & M_1(L) & C_1 \\
\midrule
\multirow{10}{*}{$-1$}
\rsp & 2 & \frac{1}{2} & (\frac{1}{2}, \frac{1}{2}, \frac{1}{2}, 0) \\
\rsp & 3 & \frac{1}{4} & (\frac{1}{4}, \frac{1}{4}, \frac{1}{4}, \frac{1}{4}) \\
\rsp & 5 & \frac{5}{16} & (\frac{1}{4}, \frac{1}{4}, \frac{1}{4}, \frac{1}{4}), (\frac{1}{2}, \frac{1}{2}, \frac{1}{2}, 0), (\frac{1}{2}, \frac{1}{2}, \frac{1}{2}, \frac{1}{2}) \\
\rsp & 6 & \frac{3}{2} & (\frac{1}{2}, \frac{1}{2}, \frac{1}{2}, 0) \\
\rsp & 7 & \frac{1}{2} & (\frac{1}{4}, \frac{1}{4}, \frac{1}{4}, \frac{1}{4}) \\
\rsp & 10 & \frac{5}{2} & (\frac{1}{2}, \frac{1}{2}, \frac{1}{2}, 0) \\
\rsp & 11 & \frac{5}{4} & (\frac{1}{4}, \frac{1}{4}, \frac{1}{4}, \frac{1}{4}) \\
\rsp & 13 & \ge 1 & (0, 0, \frac{1}{13}, \frac{1}{13}) \\
\rsp & 14 & \frac{9}{2} & (\frac{1}{2}, \frac{1}{2}, \frac{1}{2}, 0) \\
\rsp & 15 & 1 & (\frac{1}{4}, \frac{1}{4}, \frac{1}{4}, \frac{1}{4}) \\
\midrule
\multirow{4}{*}{$-2$}
\rsp & -3 & \frac{1}{3} & (\frac{1}{2}, \frac{1}{2}, \frac{1}{2}, \frac{1}{2}) \\
\rsp & 5 & \frac{11}{16}& (\frac{1}{4}, \frac{1}{4}, \frac{1}{4}, \frac{1}{4}), (\frac{1}{4}, \frac{1}{4}, \frac{1}{4}, 0), (\frac{1}{2}, \frac{1}{2}, \frac{1}{2}, \frac{1}{2}) \\
\rsp & -7 & \frac{9}{8} & (\frac{1}{4}, \frac{1}{4}, \frac{1}{4}, \frac{1}{4}), M_2 < 0.999 \\
\rsp & -11& \ge \frac{6323}{5808} & (0, 0, \frac{1}{13}, \frac{1}{13}) \\
\midrule
\multirow{7}{*}{$-3$}
\rsp & 2 & \ge \frac{1}{4} & (0, 0, \frac{1}{2}, \frac{1}{2}) \\
\rsp & 5 & \frac{1}{4} & (0, \frac{1}{4}, \frac{1}{4}, 0), (\frac{3}{8}, \frac{1}{8}, \frac{1}{8}, \frac{1}{8}), (\frac{1}{2}, \frac{1}{2}, \frac{1}{2}, 0) \\
\rsp & -7 & \frac{4}{9} & (\frac{1}{4}, \frac{1}{12}, \frac{1}{4}, \frac{1}{12}) \\
\rsp & -11& <0.46 & \\
\rsp & 13 & 1 & (\frac{1}{4}, \frac{1}{12}, \frac{1}{4}, \frac{1}{12}) \\
\rsp & 17 & \frac{13}{16} & (0, \frac{1}{4}, \frac{1}{4}, 0), (\frac{3}{8}, \frac{1}{8}, \frac{1}{8}, \frac{1}{8}), (\frac{1}{2}, \frac{1}{4}, \frac{1}{4}, 0) \\
\rsp & -19& <0.95 & \\
\midrule
\rsp -7 & 5 & \frac{9}{16} & (0, \frac{1}{4}, \frac{1}{4}, 0), (\frac{3}{8}, \frac{1}{8}, \frac{1}{8}, \frac{1}{8}), (\frac{1}{2}, \frac{1}{4}, \frac{1}{4}, 0) \\
\bottomrule
\end{array} $$
  \caption{Euclidean minima for rings $D(m,n)$}
\end{table}

Here the points in $C_1$ are given with respect to the basis $\{1,
\sqrt{m}, \sqrt{n}, \sqrt{mn}\}$. We obtain all further points at which
$M(L)$ is attained by multiplying by $-1$ and applying the three
non-trivial automorphisms of $L/\Q$. We also note that the rings
$D(-1,15)$ and $D(-3,13)$ are semi-Euclidean; possibly $D(-1,13)$
and $D(-1,17)$ are semi-Euclidean as well.

In order to prove (6.1) we proceed as follows: we consider a
residue class $\alpha \bmod {2}$ in $D(m,n)$. If $\sigma$ is that
automorphism of $\Q(\sqrt{m},\sqrt{n})$ that leaves precisely
$\Q(\sqrt{m})$ elementwise fixed, then
$N_1(\beta) = \beta \beta^\sigma$. From $\beta \equiv \alpha \bmod {2}$
it follows that $\beta^\sigma \equiv \alpha^\sigma \bmod {2}$ (because of
$(2)^\sigma = (2)$), hence $N_1(\alpha) \equiv N_1(\beta) \bmod {2}$
for all $\alpha \equiv \beta \bmod {2}$. The same of course holds
for the two other relative norms $N_2$ and $N_3$.

This observation will allow us, in the two imaginary quadratic
subfields of $\Q(\sqrt{m},\sqrt{n})$, to establish the existence of
integral elements of small norm. If we have, for example, shown that
$D(m)$, with $m$ negative, contains an element of norm 2, then it
follows easily that $m \in \{-1, -2, -7\}$ (we need only investigate
the solvability of the Diophantine equation $x^2 - m y^2 = 2$ for
$m \equiv 2,3 \bmod {4}$, respectively of $x^2 - m y^2 = 8$,
$x \equiv y \bmod {2}$, for $m \equiv 1 \bmod {4}$).

We divide the proof of (6.1) into two cases:

\subsubsection*{I. There exists an ideal of norm 2 in $D(m,n)$.}

Since $D(m,n)$ is Euclidean, this ideal is principal, and taking
relative norms shows that both imaginary quadratic subfields contain
elements of norm 2. Thus two of the three fields $\Q(\sqrt{-1})$,
$\Q(\sqrt{-2})$, $\Q(\sqrt{-7})$ are contained in $\Q(\sqrt{m},\sqrt{n})$,
leaving the possibilities $D(-1,2)$, $D(-1,7)$ and $D(-2,-7)$. The table
above shows that of these precisely the rings $D(-1,2)$ and $D(-1,7)$
are norm-Euclidean.

\subsubsection*{II. There exists no ideal of norm 2 in $D(m,n)$.}

Then only the following possibilities arise:
\begin{enumerate}
\item[(a)] $m \equiv 2, 3 \bmod {4}$, $n \equiv 5 \bmod {8}$
\item[(b)] $m \equiv 1 \bmod {8}$, $n \equiv 5 \bmod {8}$
\end{enumerate}

We first consider case (a) and distinguish whether $n > 0$
or not:

\begin{enumerate}
\item[(a1)] $m \equiv 2, 3 \bmod {4}$, $n \equiv 5 \bmod {8}$, $-m, n
  \in \N$. Here $(2)$ ramifies in both imaginary quadratic subfields
  $\Q(\sqrt{-m})$ and $\Q(\sqrt{-mn})$; if the prime ideals above
  $(2)$ in both fields were not principal, then by the class
  number formula we could not have $H=1$. Thus at least one of the two
  imaginary quadratic fields contains an element of norm 2, and
  since $m \equiv 2, 3 \bmod {4}$ we have $m \in \{-1, -2\}$. If
  $D(-1,n)$ is norm-Euclidean, then there exists an $\alpha \in
  D(-1,n)$ with $\alpha \equiv (1+\sqrt{-n})/(1+i) \bmod {2}$ and
  $N_{K/\Q}(\alpha)<16 = N_{K/\Q}(2)$. If $N_2$ denotes the norm from
  $K$ to $\Q(\sqrt{-n})$, we find $N_2(\alpha) \equiv \sqrt{-n}
  \bmod {2}$. Thus $k_2 = \Q(\sqrt{-n})$ contains an element $\beta =
  N_2(\alpha) \equiv \sqrt{-n} \bmod {2}$ with $N(\beta) < 16$ ($N$
  here of course denotes the norm from $k_2$ to $\Q$). Since
  $N(\beta) \ge N(\sqrt{-n}) = n$, we must therefore have $n<16$, and
  since $n \equiv 5 \bmod {8}$, only the two possibilities $n=5$ and
  $n=13$ remain.

That $D(-1,13)$ is not norm-Euclidean can be seen either from the above
table or by using (1.4) with $L = \Q(\sqrt{-13})$, $K =
\Q(\sqrt{-7})$, $B = (3+2i)$, $x=2$ and $e=4$. It then only remains
to show that there exists no $r \in \Z[i]$ with
$N(r)<13$ and $r \equiv 4 \bmod {3+2i}$ such that $r$ is a norm from $D(-1,13)$.
The first two conditions are satisfied only by $r = 1-2i$ and
$r = -1+i$; since $\left(\frac{13}{2}\right) =
\left(\frac{13}{3}\right) = -1$, both values of $r$ remain prime in $L$.

If $D(-2,n)$ is norm-Euclidean, we correspondingly consider the
residue class $\alpha \equiv \sqrt{-2n} \cdot (1 + \sqrt{n})/2 \bmod {2}$
and find, as above, $\beta = N_2(\alpha) \equiv 1 + \sqrt{-2n}
\bmod {2}$ and $N(\beta) < 16$. From $n \equiv 5 \bmod {8}$ and $1 + 2n <
16$ it then follows that $n=5$.

\item[(a2)] $m \equiv 2, 3 \bmod {4}$, $n \equiv 5 \bmod {8}$, $-m, -n
  \in \N$.
Here (1.4), with $K = \Q(\sqrt{n})$, $L = K(\sqrt{m})$, $B
= \fp_1$ and $\epsilon = (1 + \sqrt{n})/2$ ($\epsilon$ is a
quadratic residue $\bmod {\fp_1}$ since
$\Phi_K(\fp_1) = 3$), shows that $\Q(\sqrt{n})$ contains an element
of odd norm $< 4$; thus $n=-3$ or $n=-11$.

In the case $n=-11$, moreover, $(1 + \sqrt{-11})/2$ must still be a norm
from $D(m,-11)$, so $D(m)$ contains an element of norm 3, and this
implies $m = -2$.

If $D(m,-3)$ is norm-Euclidean, we consider the residue class
$\alpha \equiv \rho + \sqrt{m} \bmod {2}$, where $\rho$ is a primitive 3rd
root of unity. This gives $N_1(\alpha) \equiv m+1 + \sqrt{m} \bmod {2}$.
In the case $m \equiv 2 \bmod {4}$, $\Q(\sqrt{m})$ thus contains a
$\beta \equiv 1 + \sqrt{m} \bmod {2}$ with $N(\beta) < 16$, and this
yields $m \in \{-2, -6, -10, -14\}$. Among these, of the rings $D(m,-3)$
only $D(-2,-3)$ and $D(-2,-3) = D(6,-3)$ have class number 1.

Correspondingly, in the case $m \equiv 3 \bmod {4}$ we obtain only the
possibilities $m \in \{-1, -3, -13\}$, of which however only $D(-1,-3)$
has class number 1.

\item[(b)] $m \equiv 1 \bmod {8}$, $n \equiv 5 \bmod {8}$.
\end{enumerate}

Let $\alpha \equiv (1 + \sqrt{m})/2 \bmod {2}$; then $N_1(\alpha)
\equiv (\pm 1 + \sqrt{m})/2 \bmod {2}$. If, therefore, the residue class
$\alpha \bmod {2}$ is to contain an element of norm $< 16$, then
$D(m)$ must also contain such an element; this yields $m < 64$. Thus
each of the two imaginary quadratic number fields contains an element
of norm $< 64$, and among the rings still remaining only the following
have class number 1: $D(-3,5)$, $D(-3,-7)$, $D(-3,-11)$, $D(-3,17)$,
$D(-3,-19)$, $D(-3,-43)$, $D(-7,5)$, $D(-7,-11)$, $D(-7,-19)$,
$D(-7,-43)$ and $D(-11,-19)$. J.\ Sauvageot\index[N]{Sauvageot}, in
1972/73, attempted, by considering the residue class $\bmod {2}$, to
exclude some further rings, but made an error in doing so. In fact,
here we can only proceed further with (1.4) (or 1.12):

$$ \begin{array}{cc|c|c|c|c|c}
    \toprule
    m & n & x & K & B & r \bmod {B} & M(K)\ge \\ \midrule
    -3 & -43 & \frac{1 + \sqrt{-43}}2 & \Q(\sqrt{-43}) & (3)
         & \sqrt{-43} & \frac{13}{9} \\
    -7 & -11 & 3 & \Q(\sqrt{-7}) & 2 + \sqrt{-7} & -2 & \frac{16}{11} \\
    -7 & -19 & 2 & \Q(\sqrt{-19}) & \frac{3 + \sqrt{-19}}2 & -3 & 1 \\
    -7 & -43 & \frac{3 + \sqrt{-43}}2 & \Q(\sqrt{-43}) & (7)
         & \frac{-3 + 3\sqrt{-43}}2 & \frac{99}{49} \\
   -11 & -19 & 4 & \Q(\sqrt{-19}) & \frac{5 + \sqrt{-19}}2
         & 5 & 1 \\ \bottomrule
\end{array} $$

This proves (6.1).

Some of the fields listed under (6.1) can be recognised fairly easily
as norm-Euclidean. To this end we note the following: if $R$ is a
number ring, $\alpha, \beta \in R$ and if also $(\alpha + \beta)/2$ is
integral, then we can shift $x\alpha + y\beta$ (with $x, y \in \Q$)
$\bmod {R}$ so that $|x|, |y| \le \frac{1}{2}$ and
$|x \cdot y| \le \frac{1}{4}$. We can clearly achieve $|x|, |y| \le
\frac{1}{2}$; if then, say, $x \cdot y \ge \frac{1}{2}$ (so that
$x, y \ge 0$), we set $x' = x - \frac{1}{2}$, $y' = y - \frac{1}{2}$
and have $|x'|, |y'| \le \frac{1}{2}$ and $|x' \cdot y'| \le \frac{1}{4}$.

We further note that $K$ is totally complex, so $N_{K/\Q}(\alpha) \ge 0$
for all $\alpha \in K$. To show $|N_{K/\Q}(\alpha)| < 1$ it therefore
suffices to show $N_{K/\Q}(\alpha) < 1$. We shall usually do this by
writing $N_{K/\Q}(\alpha) = X^2 - n Y^2$ and then establishing $|X| < 1$;
for then $N_{K/\Q}(\alpha) = X^2 - n Y^2 \le X^2 < 1$ (better
estimates are of course obtained if we can also show $|Y| \ge c > 0$ for
certain $\alpha$, since then $N_{K/\Q}(\alpha) = X^2 - n c^2$).

We shall use these considerations constantly in what follows; thus,
e.g., we choose, in the case $D(-1,2)$, the basis
$\{1, i, \sqrt{2}, \sqrt{-2}\}$, note that $\frac{\sqrt{2}+\sqrt{-2}}2$
is integral, and shift $x + y i + z \sqrt{2} + w \sqrt{-2} \bmod {R}$ so
that $|x|, |y|, |z|, |w| \le \frac{1}{2}$ and $|z| + |w| \le \frac{1}{2}$.
The last inequality yields $z^2 + w^2 \le \frac{1}{4}$ and hence
$x^2 + y^2 + 2 z^2 + 2 w^2 \le 1$ (with equality holding precisely then
when $|x| = |y| = \frac{1}{2}$ and $|z| = \frac{1}{2}$, $w = 0$ or $z = 0$,
$|w| = \frac{1}{2}$).

Now $N_{K/\Q}(x + y i + z \sqrt{2} + w \sqrt{-2}) =
X^2 - 2 Y^2$ with $X = x^2 + y^2 + 2 z^2 + 2 w^2$ and $Y = 2(x z + y
w)$, so $N_{K/\Q}(x + y i + z \sqrt{2} + w \sqrt{-2})
\le X^2 \le 1$. Since equality can occur only in the two cases indicated
above, where $|Y| = 2 |x z| = \frac{1}{2}$ and hence
$X^2 - 2 Y^2 = \frac{1}{2}$, we have in every case obtained
$N_{K/\Q}(x + y i + z \sqrt{2} + w \sqrt{-2}) < 1$.

In $D(-1,n)$, $n \in \{-3, 5\}$, we proceed somewhat differently: here
we choose the basis $\{1, i, \sqrt{n}, \sqrt{-n}\}$ and note that
$(1+\sqrt{n})/2$ and $(i+\sqrt{-n})/2$ are integral. We can now shift
$x + y i + z \sqrt{n} + w \sqrt{-n} \bmod {R}$ so that, in the case
\begin{align*}
n & = -3: & x^2 + 3 w^2 &\le \frac{1}{3}, & y^2 + 3 z^2 &\le \frac{1}{3} \\
n & = 5: & x^2 + 5 z^2 &\le \frac{9}{20}, & y^2 + 5 w^2 &\le \frac{9}{20}
\end{align*}
(this follows from the proof of (0.21)). It again follows that
\begin{align*}
  N_{K/\Q}& (x + y i + z \sqrt{n} + w \sqrt{-n})
        = X^2 - |n| Y^2 \le X^2 \quad \text{with} \\
  |X| & = x^2 + y^2 + 3 z^2 + 3 w^2 \le \frac{4}{3} \quad
         \text{in the case } D(-1,-3) \quad \text{and} \\
  |X| &= x^2 + y^2 + 5 z^2 + 5 w^2 \le \frac{9}{10} \quad
         \text{in the case } D(-1,5).
\end{align*}

Quite analogously we can now show that the rings $D(-3,2)$, $D(-3,-2)$,
$D(-3,5)$, and $D(-3,-7)$ are norm-Euclidean. That $D(-1,2)$ and
$D(-1,3)$ are norm-Euclidean was already known to Eisenstein (1850),
whose proof is rather similar to the one above. Masley\index[N]{Masley}
\cite{Mas72}, in 1972 (unaware of Eisenstein's work), gave further
proofs of this, though far more complicated than those given above.
Lakein \cite{Lak72}, likewise in 1972, then established that the
Euclidean algorithm also holds in the rings $D(-1,5)$, $D(-1,7)$,
$D(-3,2)$, $D(-3,-2)$, $D(-3,5)$, $D(-3,-7)$, and $D(-3,-11)$. This can
be checked by computer; moreover, we obtain the new norm-Euclidean
rings $D(-2,5)$, $D(-3,17)$, $D(-3,-19)$ and $D(-7,5)$.

Finally, as in the real quadratic case, we can give the first minima
$M(K)$ for a whole class of number fields:

\begin{quote}
  \textbf{(6.2)} {\em Let $m = n^2 + 1$, $n \equiv 1 \bmod {2}$, $R =
    \Z[i, \sqrt{m}, (\sqrt{m} + \sqrt{-m})/2]$, and $f$ the
    absolute value of the norm. Then $M(f) = \frac{m}{4}$, and this
    minimum is attained mod $R$ only at the point $(1 + i + \sqrt{m})/2$.
    If $m$ is square-free, then $M(f) = M(K)$.}
\end{quote}

We proceed as in $\Q(\sqrt{m})$ (see \S\ 2) and set
$\theta = n + \sqrt{m}$. Then $\{1, i, \theta,
\frac{1+i+\theta+i\theta}2 \}$ is an integral basis, and with
$x = a+bi$, $y = c+di$ we have $N_0(x+y\theta) = x^2 + 2nxy - y^2$, as
well as $N_{0/\Q}(x+y\theta) = |x^2 + 2nxy - y^2|^2$ (here we regard
$x^2 + 2nxy - y^2$ as an element of $\mathbb{C}$; $|\cdot|$ is the
ordinary absolute value on $\mathbb{C}$). For every exceptional point
$z = x+y\theta$ we therefore have $\frac{m}{4} \le N_{0/\Q}(x+y\theta)$
and thus $\frac{n}{2} < |x^2 + 2nxy - y^2|$.

We can now choose $z \bmod {R}$ so that first $|c|, |d| \le
\frac{1}{2}$ and $|c|+|d| \le \frac{1}{2}$, and then $|a|,|b| \le
\frac{1}{2}$. This gives $|x|^2 = a^2 + b^2 \le \frac{1}{2}$
and $|y|^2 = c^2 + d^2 \le \frac{1}{4}$, hence
\[
\frac{n}{2} \le |N_0(z)| = |x^2 + 2nxy - y^2|
            \le |x|^2 + 2n |xy| + |y|^2 \le n |x| \cdot \frac{3}{4}.
\]
This implies $|x| \ge \frac{1}{2} - \frac{3}{4n}$ for every
exceptional point $z = x+y\theta$. By instead choosing $z \bmod {R}$ so
that first $|a|, |b| \le \frac{1}{2}$, $|a|+|b| \le \frac{1}{2}$ and
then $|c|,|d| \le \frac{1}{2}$, it follows correspondingly that
$|y| \ge \frac{1}{2} - \frac{3}{4n}$. If we now choose the fundamental
domain $F$ suitably, every exceptional point $z = x+y\theta$ lies in the
set defined by the inequalities
$\bigl||x| - \frac{1}{2}\bigr| \le \frac{3}{4n}$,
$\bigl||y| - \frac{1}{2}\bigr| \le \frac{3}{4n}$.

Now let again $|c|,|d| \le \frac{1}{2}$ and $|c|+|d| \le
\frac{1}{2}$. Together with $|y| \ge \frac{1}{2} - \frac{3}{4n}$,
it follows, since
$c^2 + d^2 = |y|^2 \ge \frac{1}{4} - \frac{3}{4n} + (\frac{3}{4n})^2$
and $c^2 \le (\frac{1}{2} - |d|)^2$, that
\[
\frac{1}{4} - \frac{3}{4n} \cdot \Big(\frac{3}{4n}\Big)^2
   \le c^2 + d^2 \le 2d^2 - |d| + \frac{1}{4}
\]
holds. For $n \ge 7$ we then find
\[
\Big(|d| - \frac{1}{4}\Big)^2
   \ge \frac{1}{16} - \frac{3}{8n} + \frac{9}{32n^2}
   \ge \Big(\frac{1}{4} - \frac{1}{n}\Big)^2,
\]
and, taking into account that $|x|$ and $|y|$ lie close to
$\frac{1}{2}$, we can also establish this for $n \ge 5$.
Thus $|d| < \frac{1}{n}$, $|c| - \frac{1}{2} < \frac{1}{n}$ or
$|c| < \frac{1}{n}$, $|d| - \frac{1}{2} < \frac{1}{n}$. Quite analogously,
the corresponding inequalities for $a$ and $b$ follow. This leaves
the following possible exceptional sets ($\delta = \frac{1}{n}$):
\begin{align*}
  S_1 & = \Big(0, \frac{1}{2}, \frac{1}{2}, 0\Big) +
          (-\delta, \delta) \times (-\delta, \delta) \times
          (-\delta, \delta) \times (-\delta, \delta), \\
   \quad \text{and} \\
   S_2 & = \Big(\frac{1}{2}, 0, \frac{1}{2}, 0\Big) +
          (-\delta, \delta) \times (-\delta, \delta) \times
          (-\delta, \delta) \times (-\delta, \delta), \\
\intertext{ since, for instance, $S_1$ and }
S_3 & = \Big(\frac{1}{2}, 0, 0, \frac{1}{2}\Big) +
          (-\delta, \delta) \times (-\delta, \delta) \times
          (-\delta, \delta) \times (-\delta, \delta),
\end{align*}
are congruent modulo $R$. It can now easily be computed that $S_2$
contains no exceptional point, so $S_1$ is the only exceptional set.

A further estimate of the norm on $S$ now shows that we may replace
$\delta = \frac{1}{n}$ by $\epsilon = \frac{1}{3n}$. We then note that
$\theta$ is a unit and that $(a+bi+c\theta + di\theta)\theta =
c+di+(a+2nc)\theta + (b+2nd)i\theta$ holds, and establish that the
bounds allow the application of (2.3). Thus $z = (i+\theta)/2 \bmod {R}$
is the only possible exceptional point, and since $N_{0/\Q}(z) =
\frac{m}{4}$, it follows that $M(K) \le \frac{m}{4}$.

In order to estimate $M(K)$ from below we proceed as follows: let,
e.g.\, $a = d = 0$, $b \equiv c \equiv \frac{1}{2} \bmod {1}$. Then, with
$x=a+bi$, $y=c+di$, the congruences $2xy = 2(ac-bd) + 2(ad+bc)i =
\frac{1}{2} \bmod {\Z[i]}$ and $x^2 - y^2 = a^2 - b^2 - c^2 +
d^2 + 2(ab-cd)i = \frac{1}{2} \bmod {\Z[i]}$ hold. Thus
$|\operatorname{Re}(x^2 + 2nxy - y^2)| \ge \frac{1}{2}$ and
$|\operatorname{Im}(x^2 + 2nxy - y^2)| \ge \frac{1}{2}$, and consequently
$N_{0/\Q}(z-a) \ge (1+n^2)/4 = \frac{m}{4}$ for all $a \in
R$. This proves all the claims.

\section*{{\sc Remarks on} \S\ 6}
\addcontentsline{toc}{section}{{\sc Remarks on} \S\ 6}
\medskip
         
\begin{tabular}{p{1cm} p{10cm}}
\toprule
\textbf{Year} & \textbf{Event} \\
\midrule
1850 & In the course of his investigations on eighth-power residues,
  Eisenstein\index[N]{Eisenstein} shows that $D(-1,2)$ is norm-Euclidean,
  and remarks that $D(-1,3)$ can be treated quite analogously. \\
1894 & Hilbert\index[N]{Hilbert} writes a treatise on Dirichlet number
  fields, in which he investigates, in particular, the quadratic
  reciprocity law in these fields. \\
1943 & Kuroda\index[N]{Kuroda} begins to study Dirichlet number fields;
  these investigations are continued by Kubota\index[N]{Kubota}
  \cite{Kub56} and Wada\index[N]{Wada} \cite{Wad66}. \\
1958 & Hasse\index[N]{Hasse} \cite{Has85} determines all imaginary bicyclic
  fields of odd class number. \\
1970 & Williams\index[N]{Williams} publishes an elementary
  computation of the integral basis of special Dirichlet number fields. \\
1972 & Lakein\index[N]{Lakein} determines some norm-Euclidean Dirichlet
  number fields; independently, Masley shows that $D(-1,2)$ and
  $D(-1,3)$ are norm-Euclidean. Sauvageot gives, for the first time,
  some Dirichlet number fields that have class number 1 but are not
  norm-Euclidean. \\
1974 & Brown\index[N]{Brown} and Parry\index[N]{Parry} determine all
  imaginary bicyclic number fields of class number 1. This was
  made possible by the solution of the ``class number 2'' problem
  for imaginary quadratic fields. \\
\bottomrule
\end{tabular}

\chapter*{\S\ 7 Dirichlet Number Fields}
\setcounter{chapter}{7}
\addcontentsline{toc}{chapter}{\S\ 7 Dirichlet Number Fields}
\markboth{Euclidean Rings}{\S\ 7 Dirichlet Number Fields}

Let $K = \Q(i)$ and $R = \Z[i]$; $R$ shares many properties with the
ring $\Z$, some of which we list here. We denote elements of $R$ by
lowercase letters ($a$, $b$, $p$, $q$), and rational integers by
capital letters ($A$, $B$, $P$, $Q$, etc.). By ``rational prime'' we
always mean a natural number. Further, let $N$ always denote the norm
from $K$ to $\Q$.

\begin{enumerate}
\item The unit group $R^\times$ is finite:
  $R^\times = \{ \pm 1, \pm i \}$;
\item $R$ is a UFD;
\item every prime $p \in R$ is associated with one of the following numbers:
  \begin{enumerate}
  \item $q = 1+i$, $N(q) = 2$,
  \item $p = A+Bi$, $A-1 \equiv B \equiv 0 \bmod {2}$,
    $N(p) = A^2+B^2=P$, where $P \equiv 1 \bmod {4}$ is a rational prime,
  \item $q = Q$, where $Q \equiv 3 \bmod {4}$ is a rational prime;
  \end{enumerate}
\item If $a \in R$, then the congruences
  $$ \begin{aligned}
  a^2 &\equiv 0 \bmod {4} &&\Leftrightarrow && a \equiv 0 \bmod {2} \\
  a^2 &\equiv 2i \bmod {4} &&\Leftrightarrow && a \equiv 1+i \bmod {2} \\
  a^2 &\equiv 1 \bmod {4} &&\Leftrightarrow && a \equiv 1 \bmod {2} \\
  a^2 &\equiv 3 \bmod {4} &&\Leftrightarrow && a \equiv i \bmod {2};
  \end{aligned} $$
\item If $p \in R$ is prime and $(a,p)=1$ for some $a \in R$, then
  $a^{Np-1} \equiv 1 \bmod {p}$.
\end{enumerate}

The ``little Fermat'' theorem in 5.\ allows us to define a quadratic
residue symbol in $R$: for this, let $a,p \in R$, with $p$ prime,
$Np \equiv 1 \bmod {2}$, and $(a,p)=1$. We then set $[a/p] = \pm 1$,
choosing the sign to match that in $a^{(Np-1)/2} \equiv
\pm 1 \bmod {p}$. In analogy with the Jacobi symbol, this residue
symbol can also be generalized to composite denominators; we have

\begin{enumerate}
\item[6.] The quadratic reciprocity law in $\Z[i]$: for $a,b \in R$
  with $a \equiv b \equiv 1 \bmod {2}$, we have $[a/b] = [b/a]$.
  Moreover, writing $a = A+Bi$ and $Q = A+B$, we have the
  supplementary laws $[\frac ia] = (-1)^{B/2}$ and $[\frac{1+i}a] =
  (\frac{2}{Q})$. If, in addition, $b = C+Di$ is prime and $D \neq
  0$, then via $[\frac ab] = (\frac{AC + BD}p)$, $P = N(b) =
  C^2+D^2$, the residue symbol $[\frac ab]$ can be reduced to a
  Legendre symbol; this identity also allows the QRL in $R$ to be
  derived from that in $\Z$.
\item[7.] Let $L = K(\sqrt{m})$ with $m \in R$ square-free, and let
  $S = \cO_L$ be the ring of integers in $L$; then $S$ has integral
  basis $\{1, \beta\}$ over $R$, and relative discriminant $d$, as
  follows:
  \begin{align*}
    \text{I.} \quad & m \equiv 1 \bmod {4} : \quad
      \beta = \frac{1+\sqrt{m}}{2}, \quad d = m; \\
    \text{II.} \quad & m \equiv \pm 1 + 2i \bmod {4} : \quad
      \beta = \frac{1+\sqrt{m}}{1+i}, \quad d = 2im; \\
    \text{III.} \quad & m \equiv i \bmod {2} \text{ or }
      m \equiv 0 \bmod {(1+i)} : \quad
      \beta = \sqrt{m}, \quad d = 4m.
  \end{align*}
\end{enumerate}

If $\alpha = r + s\sqrt{m}$ for $r,s \in K$, then
$\alpha' = r - s\sqrt{m}$ is called the conjugate of $\alpha$. The
relative discriminant is then given by $d = (\alpha - \alpha')^2$,
where it should be noted (here and in what follows) that $d$ is
defined only up to a factor of $\pm 1$. This factor, however, has no
effect on the residue symbol $[\frac{d}{p}]$ or on the extension
$K(\sqrt{d})$, since $\pm 1$ is a square. Finally, the absolute
discriminant of $L$ is obtained from $d$ by
$\disc (L/\Q) = 16N(d)$.

By 7.\ we have $d \equiv 0 \bmod {2}$ or $d \equiv 1 \bmod {4}$ for
the relative discriminant $d$ of quadratic extensions of $K$; this in
fact holds for extensions of arbitrary degree: if, for instance, the
relative discriminant $d$ of an extension $L/K$ satisfied $d \equiv
\pm 1 \bmod {4}$, then $L/K$ would certainly not be Galois (since $d$
is not a square in $K$). Furthermore, $(1+i)$ is unramified in $L/K$,
and hence the same holds for the normal closure $M/K$. This normal
closure $M$, however, contains the quadratic subfield $L(\sqrt{d})$,
so $(1+i)$ is also unramified in $L(\sqrt{d})$. By the decomposition
law for relatively quadratic extensions (see Hilbert), the congruence
$x^2 \equiv d \bmod {4}$ would then have to be solvable in $K$, which
is not the case, given $d \equiv \pm 1 \bmod {4}$.

The same argument clearly also applies to the case
$d \equiv \pm 1 + 2i \bmod {4}$. Finally, if $d \equiv 0 \bmod
{(1+i)}$, then $(1+i)$ is ramified in $L/K$; writing $e = e(L/K)$
for the ramification index, we then have either $e=2$ or $e \ge 3$,
and in both cases known theorems (on wild ramification) give
$(1+i)^2 \mid d$.

To describe the decomposition law in $L$ in terms of the residue
symbol $[\frac{d}{p}]$, we extend it, in analogy with the Kronecker
symbol, by setting for $p=1+i$
\[
\left( \frac{d}{p} \right) =
\begin{cases}
+1, & \text{if } d \equiv 1 \bmod {p^5} \\
-1, & \text{if } d \equiv 5 \bmod {p^5}
\end{cases}
\]

Thus we have

\begin{enumerate}
\item[8.] Let $p \in R$ be prime and $L/K$ a quadratic extension of
  $K$ with relative discriminant $d$; then there are the following cases:
  $$ \begin{array}{rclrcl}
    {} [\tfrac{d}{p}] & = & +1: & (p) & = & P_1 P_2
         \text{ for two distinct prime ideals } P_1, P_2 \text{ in } S; \\
    {} [\tfrac{d}{p}] & = & -1: & (p) & = & (p) \text{ remains prime;} \\
    {} [\tfrac{d}{p}] & = &  0: & (p) & = & P^2 \text{ ramifies in } L/K.
  \end{array} $$
\end{enumerate}

Since $L$ is totally complex, it follows first that
$N_{L/\Q}(\alpha) \ge 0$ for all $\alpha \in L$, and second, by
Dirichlet's unit theorem,
\begin{enumerate}
\item[9.] Every unit $e \in S^\times$ can be written as $e =
  \zeta u^k$, where $\zeta$ is a root of unity and $u$ is the
  fundamental unit of $L$. Here $\zeta = \zeta_8$ if
  $L = K(\sqrt{i}) = \Q(i,\sqrt{2})$; $\zeta = \zeta_{12}$ if
  $L = K(\sqrt{-3})$; and $\zeta = \zeta_4 = i$ otherwise.
\end{enumerate}

For real quadratic number fields, the question of interest is
whether the norm of the fundamental unit equals $\pm 1$ or $-1$;
here the corresponding question is whether $N_{L/K}(u) = \pm 1$ or
$\pm i$, and the answers turn out to be largely analogous.

The properties of Dirichlet number fields (that is, quadratic
extensions of $K$) listed so far are classical and can already be
found in Dirichlet and Hilbert. With quite elementary methods,
however, further-reaching results can be obtained here as well.

Our first result concerns the parity of the class number of $L$; the
only prior work on this is due to Popović\index[N]{Popović}
\cite{Pop38}, who showed that the class number $h=h(L)$ can be odd
only if the relative discriminant $d$ is prime, or $d \equiv \pm 1
\bmod {4}$ and $d$ has exactly two prime divisors (incidentally, the
claims made in the review of this paper (Reviews in Number Theory
1940--1972, R 16-33) are incorrect). This result can, however, be
sharpened further, and for this we need

\begin{quote}
  \textbf{(7.1)} {\em Let $p \equiv \pm 1 \bmod {4}$ be prime in $R$, then
    there exist $a, b \in R$ with $p = a^2 + i b^2$,
    $a \equiv 1 \bmod {(1+i)}$. }
\end{quote}

\begin{proof}
  Let $L = K(\sqrt{i}) = \Q(i,\sqrt{2})$; by \S\ 1 (or also \S\ 5),
  $L$ is norm-Euclidean and hence has class number 1. Since
  $[\frac{d}{p}] = +1$, $p$ splits in $S$, so there exists a $\pi \in
  S$ with $N_{L/K}(\pi) = \zeta p$, where $\zeta$ is a unit in $R$
  (so $\zeta = \pm i, \pm 1$). Since $\{1, \sqrt{i}\}$ is an integral
  basis of $S$ over $R$, we may write $\pi = a + b\sqrt{i}$ for some
  $a, b \in R$, giving $\zeta p = a^2 + i b^2$. Dividing by $\zeta$
  now yields the claim.
\end{proof}

\textbf{Rem.:} The analogue in $\Z$ is the representation of primes
$P \equiv 1 \bmod {4}$ as a sum of two squares; in the proof of
(7.2) as well, which in the real quadratic case goes back to
Rédei\index[N]{Redei@R\'edei} \cite{Red60}, $p = a^2 + i b^2$ takes
over the role played by $P = a^2 + b^2$ in Rédei's proof.

\begin{quote}
  \textbf{(7.2)} {\em Let $L = K(\sqrt{m})$, $m \in R$ square-free. If
    $h(L) \equiv 1 \bmod {2}$, then only the following possibilities
    can occur:}
  \begin{enumerate}
  \item $m = i$
  \item $m = p$, $p \in R$ {\em prime}, $p \not\equiv \pm i \bmod {4}$
  \item $m = pq \equiv \pm 1 \bmod {4}$ {\em for primes} $p, q \in R$,
    {\em where either $p, q \equiv \pm 1 + 2i \bmod {4}$ or
      $p \equiv 1 \pm i$, $q \equiv \pm 1 + 2i \bmod {4}$.}
  \end{enumerate}
\end{quote}

\begin{proof}
  Let $d$ be the relative discriminant of $L/K$, as defined in 7.
  Since $d$ (by 7.) is not a unit, there exist primes $p \in R$ with
  $p \mid d$. By the decomposition law, we have $(p) = \fp^2$ for
  every such $p$, where $\fp$ is a prime ideal in $S$ above $p$.
  Since $L$ has odd class number by hypothesis, $\fp$ must be a
  principal ideal whenever $\fp^2$ is, so there exist $\pi \in S$
  and a unit $e \in S^\times$ with $p e = \pi^2$, i.e.\ with
  $\sqrt{p e} \in S$ (this unit $e$ of course depends on the
  particular $p$ under consideration).

  By multiplying $e$ by a suitable power of the fundamental unit $u$
  if necessary, we may restrict ourselves to the two possibilities
  $e=\zeta$ and $e=\zeta u$, where $\zeta$ is a 4th root of unity.
  If $e=\zeta$ for some $p$ with $p \mid d$, then
  $\sqrt{\zeta p} \in S \setminus R$, and we must have
  $L = K(\sqrt{e})$, hence $m = e p$: so case i) or case ii) occurs.
  So let $m \ne i$ and $e = \zeta u$ for all $p \mid d$ (where, as
  already noted, $\zeta$ depends on $p$ and need not be the same
  for every $p$). It then follows that
  \begin{enumerate}
  \item[I.] $N_{L/K}(u) = 1$ (and not $= i$): since $pe$ is a
    square in $S$, so is $pe'$, i.e.\ we have $\sqrt{p e} \cdot
    \sqrt{p e'} = p \sqrt{e e'} \in S$. But $m \ne i$, hence
    $\sqrt{i}$ does not lie in $S$, so we must have
    $e e' = \pm u u' = \pm 1$.
  \item[II.] $d$ has at most two prime divisors: if $p \mid d$
    and $q \mid d$ for non-associated primes $p,q \in R$, then,
    with $e_1 = \zeta_1 u$ and $e_2 = \zeta_2 u$, we get
    $\sqrt{p e_1} \cdot \sqrt{q e_2} = u \sqrt{p q \zeta_1 \zeta_2}
    \in S$, so $m \sim p q$. If, moreover, $p$ and $q$ are both
    different from $1+i$, then the congruence $m \equiv \pm 1 \bmod
    {4}$ must hold, since otherwise $1+i$ would, by 7., also be a
    divisor of $d$ besides $p$ and $q$. We must still exclude the
    following possibilities (here we go beyond Popović's results):
  \begin{enumerate}
  \item[(a)] $m = p q$, $p = 1 \pm i$, $q \equiv \pm 1, \pm i \bmod {4}$:
    (in the case $q \equiv 2 \pm i \bmod {4}$, we replace $p = 1 \pm i$
    by $p = 1-i$, so that this case need not be treated separately).
    By (7.1), there exist $a,b \in R$ with $m = p q = a^2 + i
    b^2$ and $a \equiv b \equiv 1 \bmod {1+i}$. We now set $\pi = a
    \cdot \sqrt{m}$ and $\rho = (b, \pi)$; here $\pi$ is of course
    determined only up to a unit in $S$. We then see:
    \begin{itemize}
    \item $(\pi, \pi') = 1$: a common divisor $\sigma$ of $\pi$
      and $\pi'$ would also divide $\pi \cdot \pi' = 2a$ and $\pi -
      \pi' = 2\sqrt{m}$; since $(a,m)=1$, a fortiori
      $(a,\sqrt{m})=1$, so $\sigma$ would have to divide $(2)$.
      But this would give $(1+i) N_{L/K}(\pi) = a^2 - m = i b^2$,
      contradicting $b \equiv 1 \bmod {1+i}$.
    \item $\rho^2 / \pi$ is a unit; this follows from
      $(\rho^2) = (b,\pi)^2 = (b^2, b\pi, \pi^2) = (\pi \pi', b\pi,
      \pi^2) = \pi (\pi', b, \pi) = (\pi)$.
    \end{itemize}
    Now $N_{L/K}(\rho) = \zeta b$ for some 4th root of unity
    $\zeta$, hence $N_{L/K}(\rho)^2 = \pm b^2$; and we have already
    seen above that $N_{L/K}(\pi) = i b^2$. Thus the unit
    $\rho^2 / \pi$ has norm $N_{L/K}(\rho^2 / \pi) = \pm b^2 /
    i b^2 = \pm i$, contradicting I.
  \item[(b)] $m = p q \equiv 1 \bmod {4}$, $p \equiv q \equiv 1 \bmod {4}$:
    as in (a), we again have $m = a^2 + i b^2$, but now $a-1
    \equiv b \equiv 0 \bmod {2}$. We set $\pi = (a \cdot
    \sqrt{m})/2$, $\rho = (\pi', b/2)$, and conclude as above that
    $\rho^2 / \pi$ is a unit of norm $\pm i$.
  \item[(c)] $m = p$, $p \equiv \pm i \bmod {4}$: we show that
    $M = L(\sqrt{i})$ is an unramified quadratic (hence abelian)
    extension of $L$; by class field theory (or
    Hilbert's Theorem 94), the class number of $L$ must then be even.
    Since $L$ is totally complex, it suffices to consider finite places.
    Because $i$ is a unit, only prime ideals above $(2)$ can ramify
    in $M/L$. If this were the case, $(2)$
    would have to be totally ramified in $M$. However,
    $M = F(\sqrt{m})$ can be written with $F = \Q(\zeta_8)$, and by
    the decomposition law for relatively quadratic extensions
    (Hilbert), the prime ideals above $(2)$ in $M/F$ are unramified
    (note the congruence $m \equiv i \equiv (\sqrt{i})^2 \bmod {4}$
    in $F$).
  \end{enumerate}  
  \end{enumerate}
  
\end{proof}

\textbf{Rem.:} It would also be desirable to have a proof of case c)
that dispenses with class field theory.

It is interesting to note that the analogy between Dirichlet fields
and quadratic number fields extends further: exactly as
Cohn\index[N]{Cohn} \cite{Coh62}, we can show:

\begin{itemize}
\item Let $m = p \equiv \pm 1 \bmod {4}$ be prime and $L = K(\sqrt{m})$;
  then there exists a unit $u \in S$ with $N_{L/K}(u) = i$, and $L$
  has odd class number.
\item Let $p,q \equiv \pm 1 + 2i \bmod {4}$ be non-associated primes,
  $m = pq$, and $L = K(\sqrt{m})$; writing $pS = P^2$ and $qS = Q^2$
  for prime ideals $P,Q$ in $S$, then $P$ and $Q$ are principal
  ideals. Moreover, there exist $x,y \in R$ with $4\zeta = p x^2 + q
  y^2$, where $\zeta$ is a 4th root of unity.
\end{itemize}

Exactly as in Cohn, we can derive from these two facts a proof
of the QRL in $R$.

We now turn to the question of which Dirichlet number fields are
norm-Euclidean. As in the real quadratic case, this question can be
treated most simply when the relative discriminant is divisible by
a high power of 2:

\begin{quote}
  \textbf{(7.3)} {\em If the relative discriminant $d$ of $L/K$ is divisible by
    $4$, then $L$ is norm-Euclidean if and only if $m=i$ or
    $m=1+i$.}
\end{quote}

For the proof we use

\begin{quote}
  \textbf{(7.4)} {\em If the relative discriminant $d$ of $L/K$ is divisible by
    $4$, the class number $h$ of $L$ is odd and
    $m \not\equiv i, 1+i$, then for the fundamental unit $u$ the
    congruence $u \equiv 1 \bmod {(1+i)}$ holds.}
\end{quote}

\begin{proof}
  Since the ideal $(1+i)$ is ramified in $L/K$, there exists a prime
  ideal $\fp_1$ of absolute norm $2$ in $L$. Since $L$ has odd class
  number, $\fp_1$ is a principal ideal, i.e.\ there exists $\pi \in
  S$ with $(\pi)^2 = (1+i)$. Since $d \equiv 0 \bmod {4}$,
  $\{1, \sqrt{m}\}$ is an integral basis over $R$, so we can write
  $\pi = a + b\sqrt{m}$ for some $a, b \in R$. Thus $(a +
  b\sqrt{m})^2 = (1+i)e$ for some unit $e \in S$. By the unit
  theorem, $e = \zeta u^k$ (here $\zeta$ is a 4th root of unity,
  since $L$ contains no 8th roots of unity because $m \ne i$, and no
  12th roots of unity because $m \ne -3$). Moreover $k \equiv 1 \bmod
  {2}$, since otherwise $\sqrt{1 \pm i} \in L$ would follow. By
  multiplying $a + b\sqrt{m}$ by a suitable power of $u$ if
  necessary, we may assume $k=1$. Now $\pi \equiv \pi' \bmod {2}$, so
  $$ \zeta u = \frac{\pi^2}{1+i} = \frac{\pi \pi'}{1+i}
             = \frac{1 \pm i}{1+i} \equiv 1 \bmod (1+i), $$
  and since $\zeta$ is a 4th root of unity, we also have
  $\zeta \equiv 1 \bmod {(1+i)}$.
\end{proof}

\begin{proof}[Proof of (7.3)]
  We show that $M(L,I) \ge 5/4$ for the ideal $I = (1+i)$; since
  $\Phi_L(I) = 2$, there are two prime residue classes mod $I$. By
  (7.2), the residue class different from $1 \bmod {I}$ contains no
  units; furthermore, it contains no elements of norm $3$ (since
  $(3)$ is prime in $K$). The claim follows.
\end{proof}

We now give a short table of some Euclidean minima:

\begin{center}
$$ \begin{array}{|l|c|c|c|}
\hline
\rsp & \disc K & M(K) & C_1 \\ \hline
\rsp i    & 144 & \frac{1}{2} & \text{see \S\, 6, } D(-1,2) \\
\rsp 1+i  & 512 & \frac{1}{2} & \\ 
\rsp 2+i  & 1280 & \frac{5}{4} & M_2 \le 0.999 \\
\rsp 3i   & 2304 & \frac{5}{2} & \\ 
\rsp 3+ i & 2560 & \frac{5}{4} & M_2 \le 0.999 \\
\rsp 1+3i & 2560 & \frac{5}{4} & M_2 \le 0.999 \\ \hline
\end{array} $$
\end{center}

The case $L = K(\sqrt{p})$, $p \equiv \pm 1+2i \bmod {4}$, is already
somewhat more difficult; here, however, we have

\begin{quote}
  \textbf{(7.5)} {\em Let $L = K(\sqrt{p})$, $p \equiv \pm 1+2i \bmod 4$
    prime; if then there exist $s, t \in R$ with $s \equiv t \equiv 1 \bmod 2$,
    $(\frac{s}{p}) = (\frac{t}{p}) = -1$, $(s,t)=1$ and
    $|st| \le k \cdot |p|$ for $k = \sqrt{2}/(1+\sqrt{2})$,
    then $L$ is not norm-Euclidean.}
\end{quote}

Here $|\cdot|$ is the ordinary absolute value on $\mathbb{C}$.

\begin{proof}
   We use (1.4) with $B = (1+i)p$, $e=st$, and $k=1$; we must show
   that there is no $r \equiv e \bmod {(1+i)p}$ in $R$ with $N(r) <
   2P$ and $r = N_{L/K}(\alpha)$ for some $\alpha \in S$ (here $P =
   N(p) = |p|^2$). Now $N(r) < 2P$, so $|r| < \sqrt{2}|p|$. We
   set $r = st + l(1+i)p$ for some $l \in R$, and since $|st|
   \le k|p|$ and $|r| < \sqrt{2}|p|$, we find:
   \[ |l| \cdot \sqrt{2} |p| = |r - st| \le |r| + |st|
        \le (\sqrt{2} + k)|p| = 2|p|.\]
   Thus $|l| < \sqrt{2}$, so $l \in \{0, \pm 1, \pm i\}$.
   If $l \equiv 1 \bmod (1+i)$, then $r \equiv i \bmod 2$ would follow;
   but such an $r$ cannot be a norm from $S$, by (7.6).
   Thus $l=0$; but then $r=st$ is also not a norm from $S$, because
   $(s,t)=1$ and $(\frac{s}{p}) = (\frac{t}{p}) = -1$ (this follows as
   in the real quadratic case).
\end{proof}

\begin{quote}
  \textbf{(7.6)} {\em Let $L = K(\sqrt{m})$,
    $m \equiv \pm 1 \pm 2i \bmod {4}$ and $\alpha \in S$. Then
    $N_{L/K}(\alpha) \equiv \pm i \bmod {2}$.}
\end{quote}

\begin{proof}
  It remains to prove (7.6): let
  $\alpha = \frac{a + b\sqrt{m}}{1+i}$, $a \equiv b \bmod {1+i}$, and
  $N_{L/K}(\alpha) = c + di = \beta$. Then $a^2 - m b^2 = 2i \beta$,
  and considering this equation mod $m$ and mod $\beta$ gives
  $[\ \frac{\beta}{m} ] = [ \frac{m}{\beta}] = \pm 1$.
  If $\beta \equiv \pm i \bmod {2}$, the QRL in $R$ would give
  \[ \left[ \frac{\beta}{m} \right] = \left[ \frac{i}{m} \right]
     \left[ \frac{m}{\beta} \right] = \left[ \frac{i}{m} \right]
     \left[ \frac{m}{i} \right] = \left[ \frac{i}{m} \right]
     \left[ \frac{m}{i} \right], \]
  which, since $[ \frac{i}{m}] = -1$, contradicts the observations
  above.
\end{proof}

\begin{quote}
  \textbf{(7.7)} {\em Let $L = K(\sqrt{p})$, $p \equiv \pm 1 + 2i \bmod 4$;
    then $L$ is norm-Euclidean if and only if $p \in \{1+2i, 3+2i,
    5+2i, 1+6i, 7+2i\}$.}
\end{quote}

\begin{proof}
  That $L$ is norm-Euclidean for $p = 1+2i$ and $p = 3+2i$ has
  already been shown by Lakein\index[N]{Lakein} \cite{Lak72}. If $L$
  is norm-Euclidean, van der Linden\index[N]{Linden@van der Linden}
  \cite[p.\ 163]{Lin85} has shown that $\disc L < 16 \cdot 899\,225$.
  Since in our case $\disc L = 64 \cdot N(p)$,
  we may assume $N(p) < 224\,805$. For these $p$, a computer search
  finds a representation as in (7.4), except when
  $$ N(p) \in \{5, 13, 29, 37, 53, 61, 101, 157, 181, 229, 349, 541\}.$$
  For example:

\begin{center}
\begin{tabular}{|c|c|c|c||c|c|c|c|}
\hline
$p$ & $p$ & $s$ & $t$ & $P$ & $P$ & $s$ & $t$ \\
\hline
109 & $3+10i$ & $1-2i$ & $1+2i$ & 293 & $17+2i$ & $1-2i$ & 3 \\
149 & $7+10i$ & $1-2i$ & $1+2i$ & 317 & $11+14i$ & $1-2i$ & 3 \\
173 & $13+2i$ & $1+2i$ & 3 & 373 & $7+18i$ & $1+2i$ & $3-2i$ \\
197 & $1+14i$ & $1-2i$ & 3 & 389 & $17+10i$ & $1-2i$ & $1+2i$ \\
269 & $13+10i$ & $1-2i$ & $1+2i$ & 397 & $19+6i$ & $1-2i$ & $3-2i$ \\
277 & $9+14i$ & $1+2i$ & $1+4i$ & 421 & $15+14i$ & $1-2i$ & $1+2i$ \\
\hline
\end{tabular}
\end{center}

We now exclude the $p \ge 101$ directly using (1.4):

\begin{center}
\begin{tabular}{|c|c|c|c|}
\hline
$P$ & $p$ & $s$ & $B$ \\
\hline
61 & $5+6i$ & 5 & $(1+i)p$ \\
101 & $1+10i$ & $4-i$ & $(p)$ \\
157 & $11+6i$ & $-4+5i$ & $(p)$ \\
181 & $9+10i$ & $-4+5i$ & $(p)$ \\
229 & $15+2i$ & $4+i$ & $(p)$ \\
349 & $5+18i$ & $7i$ & $(p)$ \\
541 & $21+10i$& $-5+14i$& $(1+i)p$ \\
\hline
\end{tabular}
\end{center}
\end{proof}

The computing time required for (7.7) can be shortened somewhat by
using the following corollary of (7.6):

\begin{quote}
  \textbf{(7.8)} {\em Let $p \equiv \pm 1, 2i \bmod {4}$ be prime, $N(p)
    \ge 109$, and $a \equiv \pm 2$, $b \equiv 0 \bmod {5}$ or $a
    \equiv 0$, $b \equiv \pm 1 \bmod {5}$. Then $L=K(\sqrt{p})$
    is not norm-Euclidean.}
\end{quote}

\begin{proof}
  (7.5) with $s=1-2i$, $t=1+2i$.
\end{proof}

Here too we give some Euclidean minima:

\begin{center}
\begin{tabular}{|l|c|c|c|}
\hline
$m$ & $\operatorname{disc} K$ & $M(K)$ & $C_1$ \\
\hline
$1+2i$ & 320 & $\frac{1}{2}$ & \\
$3+2i$ & 832 & $\frac{1}{2}$ & \\
$5+2i$ & 1856 & & \\
$1+6i$ & 2368 & & \\
$3+6i$ & 2880 & & \\
$7+2i$ & 3392 & $\ge 50/53$ & \\
\hline
\end{tabular}
\end{center}

\textbf{Rem.:} The next most difficult case is $m \equiv 5
\bmod {(1+i)^6}$; here we have

\begin{quote}
  \textbf{(7.9)} {\em Let $m \equiv 5 \bmod {(1+i)^6}$, $m \ne \pm 3$,
    and $u$ the fundamental unit of $L$. If then the form $u =
    \frac{x + y\sqrt{m}}{1+i}$ with $x \equiv y \equiv 1 \bmod {1+i}$ holds,
    then $L$ is not norm-Euclidean.}
\end{quote}

\begin{proof}
   We show that the ideal $(1+i)$ is not Euclidean. Because
   $\Phi_L(1+i) = 3$, $\lVert(1+i)\rVert = 4$, and since $L$ contains
   no ideals of norm 2 or 3, it suffices to show that $u
   \equiv 1 \bmod {1+i}$. This follows immediately from
   $u-1 = \frac{(x-1) + (y-1)\sqrt{m}}{1+i}$, since
   $\frac{(x-1) + (y-1)\sqrt{m}}{2}$ is integral.
\end{proof}

With this criterion, however, only few $m$ can be excluded;
at least we have

\begin{quote}
  \textbf{(7.10)} {\em Let $m \equiv \pm 3 \bmod {(1+i)^6}$ and $m = x^2
    \pm 2$ for some $x \in R$. Then
    $\frac{2 \pm i + \sqrt{m}}{1-i} = (\frac{1 + \sqrt{m}}{2})^3$ for
    $m = 1 + 4i$, while for all other $m$ the fundamental unit
    $u \equiv \frac{x + \sqrt{m}}{1 + i}$.}
\end{quote}

It follows, for instance, that $L$ cannot be norm-Euclidean for the
following values of $m$:

\begin{align*}
m &= -11 &&= (3i)^2 - 2 &&\quad (\text{see also §\,6}) \\
m &= 7 + 12i &&= (3 + 2i)^2 + 2, &&\quad P = 193 \\
m &= 13 + 8i &&= (4 + i)^2 - 2, &&\quad P = 233 \\
m &= 5 + 24i &&= (4 + 3i)^2 - 2, &&\quad P = 601 \\
m &= 23 + 20i &&= (5 + 2i)^2 + 2, &&\quad P = 929.
\end{align*}

A somewhat more useful criterion is

\begin{quote}
  \textbf{(7.11)} {\em Let $m \equiv \pm 3 \bmod {(1+i)^6}$,
    $s \equiv t \equiv 1 \bmod {2}$, $(s,t) = 1$,
    $[ \frac{s}{m} ] = [ \frac{t}{m} ] = -1$ and
    $|st| \le (\sqrt{2}-1)|m|$, $|st+m| \ge |m|$. If then
    $|st-m| \ge |m|$ or $|st-m|$ is not a norm from $S$, then
    $L = K(\sqrt{m})$ is not norm-Euclidean.}
\end{quote}

\begin{proof}
  We use (1.4) with $e = s \cdot t$, and note that for
  $r = e + l m$, $l \in R$, we have
  $|l m| = |r - e| \le |r| + |e| < |m| + (\sqrt{2}-1)|m| = \sqrt{2}|m|$,
  hence $|l| < \sqrt{2}$, so $l \in \{0, \pm 1, \pm i\}$.

  If $l=0$, then $r=e$ is not a norm; if $l=\pm i$, then $r = e +
  l m \equiv 1 + i \bmod {2}$ would follow, and since $(1+i)$ is inert
  in $L/K$, $r$ cannot then be a norm from $S$. This leaves
  only the possibilities $l = \pm 1$, and the hypotheses guarantee
  that then either $|r| \ge |m|$ or $r$ is not a norm from $S$.
\end{proof}

\chapter*{\S\ 8 Other Number Fields of Degree 4}
\setcounter{chapter}{8}
\addcontentsline{toc}{chapter}{\S\ 8 Other Number Fields of Degree 4}
\markboth{Euclidean Rings}{\S\ 8 Other Number Fields of Degree 4}

We first consider totally complex fields $K$ of degree 4; we classify these
according to the Galois group of their normal
closure $N$:

\begin{enumerate}
\item $\Gal(N/\Q) = \Gal(K/\Q) = Z_4$ (cyclic
  group of order 4);
\item $\Gal(N/\Q) = \Gal(K/\Q) = V_4$ (Klein
  four-group);
\item $\Gal(N/\Q) = D_4$ (dihedral group of order 8);
\item $\Gal(N/\Q) = A_4$ (alternating group of order
  12);
\item $\Gal(N/\Q) = S_4$ (symmetric group of order 24).
\end{enumerate}

If $K$ has a quadratic subfield, then its Galois group is $Z_4$, $V_4$
or $D_4$; in the first and third cases $K$ has exactly one such subfield,
and in the second case it has three.

In the first case, according to van der Linden, there exist exactly
two norm-Euclidean fields: the field $\Q(\zeta_5)$ of 5th roots
of unity and the degree-4 subfield of $\Q(\zeta_{13})$.

In the second case we determined all norm-Euclidean fields in \S\ 6. So let
$\Gal(N/\Q) = D_4$ be the dihedral group. Then $K$ has exactly
one quadratic subfield $k = \Q(\sqrt{m})$, and we
distinguish the following cases:

\begin{enumerate}
  \item[3.1] $\Q(\sqrt{m})$ is imaginary quadratic. For $K$
    to be Euclidean, it must have class number 1; this
    implies that the quadratic subfield $\Q(\sqrt{m})$
    has class number 1 or 2. This is actually an elementary
    result (see e.g.\ Narkiewicz \cite{Nar74}), though it is usually
    proved using Hilbert's Theorem 94 or class field theory. If
    $k$ had class number 2, then $K$ would be the Hilbert class field of $k$,
    and hence $K/\Q$ would be Galois; but this is not the case.
    
    So $k = \Q(\sqrt{m})$ has class number 1, and we may restrict
    ourselves to the possibilities $m = -1, -2, -3, -7, -11, -43, -67,
    -163$. I do not expect that investigating these values of $k$
    will raise problems beyond those already encountered in \S\ 7 (where
    $m = -1$); rather, we should be able to proceed exactly as in the
    case of the Dirichlet fields.
    
    The two most difficult cases are naturally
    $m=-1$ and $m=-3$, because of the existence of non-trivial units; most
    of the norm-Euclidean fields considered here will probably
    contain one of these two fields as well.
\item[3.2] $\Q(\sqrt{m})$ is real quadratic. Exactly as in 3.1, we find
  that $\Q(\sqrt{m})$ must have class number 1. If we wish to apply
  criterion (1.4) here, we quickly find that the presence of infinitely
  many units in $\Q(\sqrt{m})$ has a rather disturbing effect. This
  drawback can, however, be compensated for by the following
  considerations: if $u > 1$ is the fundamental unit of
  $k$, then either $u$ is also a fundamental unit of $K$, or we have
  $K = k(\sqrt{-u})$, with $\sqrt{-u}$ a fundamental unit of $K$ (since
  $+1$ and $-1$ are the only roots of unity in $K$). For $K$ to be
  totally complex in the latter case, $u$ must be
  totally positive (i.e.\ $u>0$ and $u'>0$), and then $N_{K/\Q}(u) =
  1$. This, in turn, implies (see e.g.\ Kubota\index[N]{Kubota}
  \cite{Kub56}) that $\Gal(k(\sqrt{-u})/\Q) = V_4$. Thus $u$ is also a
  fundamental unit of $K$, so that for a given $k$ we already know the
  unit group of $K$. We use the example $k = \Q(\sqrt{2})$ to illustrate
  how to exploit this.

  We distinguish among the possible factorizations of $(\sqrt{2})$:
  \begin{enumerate}
  \item[(a)] $(\sqrt{2}) = P^2$: then $\Phi(P^2) = 2$, $u=1+\sqrt{2}
    \equiv 1 \pmod{P^2}$, and consequently every unit in $K$ is
    $\equiv 1 \pmod{P^2}$. If $K$ were norm-Euclidean, the prime
    residue class $1 + \alpha \bmod P^2$ (where $\alpha$ is a prime
    element for $P$) would have to contain an element of norm $<
    \|P^2\| = 4$.  This element cannot be a unit, so it would have to
    have norm 3. But this cannot be, because $(3)$ is inert in
    $k$. Since $(5)$ is also inert in $k$, it follows that
    $M(K) \ge \frac74$ for all such $K$.
  \item[(b)] $(\sqrt{2}) = P$ is inert: then $\Phi(P) = 3$, and
    as above it follows that the ideal $P$ -- and hence $K$ itself -- is not
    norm-Euclidean.
  \item[(c)] $(\sqrt{2}) = 2_1 2_2$ splits: by the decomposition law for
    relative quadratic extensions (see Hilbert), we have $K = k(\sqrt{\mu})$
    for some $\mu \equiv 1, 3+2\sqrt{2} \bmod 4$, i.e.\ with
    $\mu = a+b\sqrt{2}$ we have
    $a-1 \equiv b \equiv 0 \bmod 2$. Now let $P = (p)$ be a prime ideal in
    $R = \Z[\sqrt{2}]$ (i.e.\ $p \in R$), and suppose $P$ ramifies in $K/k$:
    $(p) = Q^2$. If $K$ has odd
    class number, then $Q$ must be a principal ideal, so there exists a
    $\pi \in S$ with $(\pi^2) = (p)$. Since all units lie in $k$,
    we can absorb them into $p$ and find $\pi^2 =
    p$. In this case, therefore, $\mu = \pi$. This also shows
    immediately that at most one prime ideal of $K/k$ may ramify
    if $K$ is to have odd class number: in this case 
    $a^2 - 2b^2 = q$ is a rational prime.
  \end{enumerate}  
\end{enumerate}

For $K$ to also have class number 1, the ideals $2_1$ and $2_2$
must be principal ideals; this means that the equation
$4 \cdot p = \xi^2 - \mu \eta^2$ must be solvable for some $p \in R$
of norm $2$ (in $R=\Z[\sqrt{2}]$; the factor $4$ appears because
$(1+\sqrt{\mu})/2$ or $(1+\sqrt{2} + \sqrt{\mu})/2$ is an integer). Since,
with $-\mu$, the whole right-hand side of the equation is totally positive,
we may assume $p = 2 + \sqrt{2}$ (absorbing any factors $u^{2k}$
into $\xi$ and $\eta$) and obtain
$2 + \sqrt{2} = \xi^2 - \mu \eta^2$. Now let the two Archimedean
valuations of $k$ be defined by
$|a+b\sqrt{2}|_1 := |a+b\sqrt{2}|$
and $|a+b\sqrt{2}|_2 := |a - b\sqrt{2}|$. Then, since $\xi^2$ and $-\mu
\eta^2$ are both totally positive,
\begin{align*}
  4 \cdot (2 + \sqrt{2}) & = |2+\sqrt{2}|_1
    = |\xi^2 - \mu \eta^2|_1 \ge |\mu \eta^2|_1
     \quad \text{and correspondingly} \\
  4 \cdot (2 - \sqrt{2}) & = |2 - \sqrt{2}|_2
  = |\xi^2 - \mu \eta^2|_2 \ge |\mu \eta^2|_2.
\end{align*}
Multiplying both inequalities and using
$|\alpha_1| \cdot |\alpha_2| = |N_{k/\Q}(\alpha)|$ for $\alpha \in k$
yields $32 \ge |N_{k/\Q}(\mu \eta^2)|$. Since $\eta \ne 0$ we certainly have
$|N_{k/\Q}(\eta)| \ge 1$, and we finally find
$|N_{k/\Q}(\mu)| \le 32$.

Together with the conditions on $\mu$ already derived above, this leaves
open only the possibility $\mu = \pm 5 + 2\sqrt{2}$, with
$N_{k/\Q}(\mu) = 17$ and $\disc K = 8^2 \cdot 17 = 1088$.
We have shown:

\begin{quote}
  \textbf{(8.1)} {\em If $K$ is a totally complex field of degree 4
    with $\Q(\sqrt{2}) \subset K$, then the norm-Euclidean fields
    $K$ are, at most, the following:}
    $$ \Q(\sqrt{2}, \sqrt{-1}), \quad
       \Q(\sqrt{2}, \sqrt{-3}) \quad \text{{\em and}} \quad
       \Q\Big( \sqrt{-5 - 2\sqrt{2}}\Big). $$
\end{quote}

Corresponding computations can also be carried out for other real
quadratic fields. It should be noted that the above
considerations were motivated by the class numbers of imaginary
quadratic extensions of $\Q(\sqrt{5})$ computed by Cohn\index[N]{Cohn}
\cite{Coh58}. The book by Pólya\index[N]{Polya} \cite{Pol69} also
deserves explicit mention here.

In the case $\Gal(N/\Q) = A_4$, not a single
example of a norm-Euclidean field of degree 4 is known to this day; this
is probably because no such fields exist with $\disc K <
3136 = 2^6 7^2$.

Finally, the great majority of totally complex fields of degree 4 have
the group $S_4$ as the Galois group of their normal closure. Here we
will probably be able to say little without a table of all such fields,
together with their unit and ideal class groups.

Even less can be said about fields of degree 4 with $r=2$, $s=1$ or
$r=4$, $s=0$. In isolated cases I have verified the Euclidean algorithm
by computer; however, a good deal of further computation will be
necessary before we find suitable starting points.

\section*{{\sc Remarks on} \S\ 8}
\addcontentsline{toc}{section}{{\sc Remarks on} \S\ 8}

\begin{tabular}{p{1cm} p{10cm}}
\toprule
\textbf{Year} & \textbf{Event} \\
\midrule
1937 & Dribin\index[N]{Dribin} determines the Hilbert subgroup series for
  number fields with Galois group $S_4$. \\
1956 & Godwin\index[N]{Godwin} gives a table of totally real fields of degree 4 with
  small discriminant. \\
1957 & Godwin publishes corresponding tables for number fields
  of degree 4 with $r=2$ and $r=0$. \\
1958 & H.\ Cohn\index[N]{Cohn} computes the class number of
  imaginary quadratic extensions of $\Q(\sqrt{5})$, using an
  extremely unusual method: Gauss had already discovered that the
  number of representations of a prime $p \equiv 3 \pmod{8}$ as
  a sum of three squares is closely connected with the class number of $\Q(\sqrt{-p})$;
  Cohn made use of a corresponding formula with $\Q(\sqrt{5})$ as
  base field. \\
1980 & Edgar\index[N]{Edgar} and Peterson\index[N]{Peterson}
  describe cyclic fields of degree 4. \\
1985 & Godwin extends his 1957 table for number fields
  of degree 4 with $r=2$, $s=1$. \\
1986 & H.\ Cohn and J.\ Deutsch\index[N]{ Deutsch} show by computer that
  $\Q\big(\sqrt{2 + \sqrt{2}}\big)$ and $\Q\big(\sqrt{3 + \sqrt{2}}\big)$
  are norm-Euclidean. \\
1989 & Buchmann\index[N]{Buchmann} and Pohst\index[N]{Pohst} give all
  totally real number fields of degree 4 with $\disc K < 10^6$, together
  with their units and class numbers. \\
\bottomrule
\end{tabular}

\chapter*{\S\ 9 Number Fields of Higher Degree}
\setcounter{chapter}{9}
\addcontentsline{toc}{chapter}{\S\ 9 Number Fields of Higher Degree}
\markboth{Euclidean Rings}{\S\ 9 Number Fields of Higher Degree}

We begin with a quotation from H.\ Stark\index[N]{Stark} \cite{Sta88}:

\begin{quote}
  {\em Heilbronn seems to suspect \dots that there may be infinitely
    many totally real Euclidean cubic fields. \dots he goes even
    further and names two families of quartic and sextic fields that
    should be investigated.}
\end{quote}

Here Stark refers to the following two remarks of
Heilbronn:\index[N]{Heilbronn}
\begin{itemize}
\item \cite{Hei50}: {\em Theorem 1: E.A. holds only in a finite number
  of cyclic cubic fields. The question of E.A. in real non-cyclic
  fields is thus left open. Though I cannot prove any result in the
  opposite direction I should be surprised to learn that the analogue
  of theorem 1 is true in that case.}
\item \cite{Hei51}: {\em Finally I should like to mention two types of
  cyclic fields for which E.A. may possibly hold in an infinity of
  cases.  
  \begin{enumerate}
  \item[(a)] The real quartic fields $\Q(\sqrt{\omega p})$,
    $\omega = \frac{5+\sqrt{5}}{2}$, of discriminant $125p^2$, where
    $p \equiv 3 \bmod 20$ is a prime;
  \item[(b)] The complex sextic field $\Q(\zeta_9 + \zeta_9^{-1},\sqrt{-p})$
    of discriminant $-3^8p^3$, where $p \equiv 3 \bmod 4$ is a prime.
  \end{enumerate} }
\end{itemize}

Part (b) of Heilbronn's 1951 conjecture is rather
mysterious: $K = \Q(\zeta_9 + \zeta_9^{-1}, \sqrt{-p})$ is
a cubic extension of $k=\Q(\sqrt{-p})$, and indeed $K/k$ is
ramified at the prime ideals above $(3)$. So $K$ does not lie in the
Hilbert class field of $k$, and we have $h(k)\mid h(K)$. However,
Heilbronn himself proved that there exist only finitely many $k$ of
class number 1; the above divisibility relation then implies
that there also exist only finitely many $K$ of class number 1. What is
almost more astonishing still is that this was also overlooked by Stark
(who, as is well known, found all imaginary quadratic number fields of
class number 1).

In fact, it is not difficult to show that $K$ can be norm-Euclidean
only for $p=3$ or $p=11$ (for $p=3$ we have $K=\Q(\zeta_9)$, so that
$K$ is norm-Euclidean). If $p>11$ and $k = \Q(\sqrt{-p})$
has class number 1, then we must have $p \equiv 1 \bmod 3$ (otherwise
$k$ would have to possess elements of norm 3, which is not the
case). So $(1 - \sqrt{-p})^2 = 1 - p - 2\sqrt{-p} \equiv \sqrt{-p}
\bmod 3$; however, the element of minimal norm in the residue class $\sqrt{-p}
\bmod 3$ is $\frac{3 - \sqrt{-p}}2$, with norm
$\frac{p+9}4$, and for $p>43$ this norm is $>9$. By (1.4) $K$ then
cannot be norm-Euclidean. In the case $p=19$, $\frac{3 - \sqrt{-p}}2$
is not a norm from $O_K$, because $(7)$ is inert in $\Q(\zeta_9 +
\zeta_9^{-1})$; so $K$ is not norm-Euclidean for $p>11$. Finally, if
$p=7$, then neither $\sqrt{-p}$ nor $\frac{3 - \sqrt{-p}}2$ is a norm
from $O_K$, which proves the claim made above.

Further classes of fields that can be investigated by simple means
are, e.g.\ $\Q(\rho, \sqrt[3]{m}\,)$ or $\Q(i, \sqrt[4]{m})$ for
$m \in \mathbb{Z}$. Here new results can be obtained without much
difficulty. The same holds for number fields of the form $\Q(\sqrt[n]{m})$
with, e.g.\ $n=3, 6, 8$, etc., or more generally for
$\Q(\sqrt[n]{m}, \sqrt{-l})$, $l \in \N$.

On the other hand, comparatively little is known about the Euclidean
algorithm in cyclotomic fields. It is known that $\Q(\zeta_m)$ is
norm-Euclidean for the values
$$ m = 1, 3, 4, 5, 7, 8, 9, 11, 12, 13, 15, 16, 20, 24. $$
Since the investigations of K.\ Uchida, H.L.\ Montgomery,
and J.M.\ Masley (see e.g.\ Masley and Montgomery \cite{MM76}),
it is known that $\Q(\zeta_m)$ has class number 1 exactly for
\begin{align*}
  m & = 1, 3, 4, 5, 7, 8, 9, 11, 12, 13, 15, 16, 17, 19, 20, 21, 24, 25, 27,
  28, 32, 33, \\
  & \quad 35, 36, 40, 44, 45, 48, 60, 84.
\end{align*}
For $m = 32$ we can show without much effort that the
residue class $1+\lambda^5 \bmod \lambda^6$ contains no units, where
$\lambda = 1+\zeta_{32}$ is an element of norm 2 in $\Q(\zeta_{32})$
(see Lenstra \cite[p.\ 14]{Len79c} on this); since $\mathbb{Z}[\zeta_{32}]$
contains no elements of norm 33 or 65, the residue class above contains
only elements of norm $\ge 97$, and it follows that $M(K) \ge
\frac{97}{64}$ (note that $N_{K/\Q}(\alpha) \equiv 1
\bmod 32$ for $\alpha \in \mathbb{Z}[\zeta_{32}]$; further details
needed for a complete proof, e.g.\ information about the units in
$\mathbb{Z}[\zeta_{32}]$, can be found in Washington \cite{Was82}).

\chapter*{\S\ 10 Open Questions}
\setcounter{chapter}{10}
\addcontentsline{toc}{chapter}{\S\ 10 Open Questions}
\markboth{Euclidean Rings}{\S\ 10 Open Questions}

\begin{quote}
  {\em There are a lot of things I don't understand myself} (Linus van Pelt)
\end{quote}

The following questions are numbered consecutively; the first
digit refers to the paragraph to which I have assigned the
question. Questions that I believe are solvable with the methods
presented in \SS\ 1--9 are formulated as problems. Finally, particularly
important questions (or ones I consider important) are printed in bold.

\subsection*{1.}

\begin{enumerate}
\item[1.1.] The bounds found in (1.10) and (1.12) depend strongly on
  the choice of the integral basis; thus the question arises whether,
  and if yes, how to choose the integral basis so that the bounds
  $\mu_i$ become as small as possible.
\item[1.2.] Are the Lenstra constants $\mu_i$ defined before (1.17)
  and the associated bounds $\lambda_i$ for $n \ge 3$ still
  finite (if $R$ is not a principal ideal domain, this certainly holds;
  for this the condition (F-2), p.~\pageref{pF2}, takes care).
\item[1.3.] In (1.21) a statement is proved for $L=Kk$, where
  $k/\Q$ is abelian. Does something similar hold if we only
  require that $k/\Q(\sqrt{-m})$ is abelian (i.e.\ can
  we replace the $m$-th roots of unity by elliptic units?)
\item[1.4.] Can we  dispense with the hypothesis ``$n$ is
  usable'' in (1.24), or can it at least be replaced by a
  weaker one?
\item[1.5.] Determine $c(K)$ also for other classes of cyclotomic fields
  than for full cyclotomic fields of prime-power discriminant.
\item[1.6.] Determine $c(K)$ for $K = \Q(\sqrt{a+b\sqrt{-3}})$,
  in particular in the case $a=-1$, $b=2$; here we know that $c(K) \le
  (4 + \sqrt{13})/12 < 0.6338$. If we could show $c(K) \le \frac12$,
  it would follow that the field of degree 8 $L = K(i)$
  is norm-Euclidean. However, recent computations of mine
  lead me to suspect that here $c(K) > \frac12$.
\end{enumerate}

\subsection*{2.}

\begin{enumerate}
\item[2.1.] Determine the missing second minima in the tables on
  pages 47 and 48.
\item[2.2.] In the theorems (2.5) to (2.9), can we also determine the
  higher minima ($M_3$, $M_4$, etc.)?
\item[2.3.] \textbf{Does (2.12) also hold in fields of unit rank $\ge 2$?}
\item[2.4.] Let $C(K)$ be the set of all $x \in K$ with $M(x) =
  M(K)$. Are then the following statements correct?
  \begin{enumerate}
  \item[(a)] The accumulation points of $C(K)$ lie in $K$;
  \item[(b)] $C(K)$ contains at least one $x$ from $K$.
  \end{enumerate}
\end{enumerate}

\subsection*{3.}

\begin{enumerate}
\item[3.1.] Find simple proofs that $D(m)$ for the values
  $$ m = 193, 241, 337, 457, 601 $$
  is not norm-Euclidean.
\item[3.2.] Show that for all primes $p \equiv 1 \pmod{24}$ with
  $p > 601$ there exist natural numbers $r, s, t, u$ with $(r,s) = (t,u) = 1$
  and $(\frac rp) = (\frac sp) = (\frac tp) = -1$.
\item[3.3.] Show that $D(10)$ and $D(65)$ are the only
  semi-Euclidean quadratic rings of class number different from 1.
\item[3.4.] \textbf{Is the ring $R = D(14)$ Euclidean?} We
  have already seen that $R$
  is not norm-Euclidean. The question of whether
  there exist other functions on $R$ with respect to which $R$ is Euclidean
  was first explicitly posed, in this special case, by Samuel\index[N]{Samuel}
  \cite{Sam71}. Lenstra\index[N]{Lenstra}
  then proposed in \cite{Len74} that we consider, in this regard, ``weighted
  norms'': to this end, choose an ideal $P$ in $R$, set
  $N(P)=c$ for a positive real $c$, and set $N(Q) = || Q ||$ for
  all other prime ideals $Q$. We then continue $N$ multiplicatively to
  all ideals in $R$ and define $N(a) := N(aR)$ for all $a \in K$.
  $N$ is called the norm weighted by $N(P)=c$. Instead of
  maps $f:R \to \mathbb{N}$, we now also wish to
  admit maps $f:R \to
  \mathbb{R}$ as Euclidean functions, provided that for every
  $a \in R$ there exist only finitely many function values
  $<f(a)$ (so we can enumerate the set $f(R)$ in order of
  increasing function values, giving a bijection $g:
  f(R) \to \mathbb{N}$; then $g \circ f: R \to \mathbb{N}$
  is a possible Euclidean function in the previous sense). We
  now consider the ideal $P = (4 + \sqrt{14})$ of norm 2 in $R$ and
  ask whether we can choose $N(P) = c$ so that the
  weighted norm becomes a Euclidean function on $R$. For the
  ideal $I = (7 + 2\sqrt{14})$ of norm 7 to become Euclidean with respect to $N$,
  every prime residue class mod $I$ must contain an element $a$
  with $N(a) < 7 = N(7)$. Since $u = 15 + 4\sqrt{14}
  \equiv 1 \mod I$, only the residue classes $\pm 1 \mod I$ contain
  units. Further, the elements $\pm 3 + \sqrt{14}$ have norm 5
  and are $\equiv 3 \mod I$. Since $P = (4 + \sqrt{14})$ and $4 +
  \sqrt{14} \equiv -3 \mod I$, $P$ also contains no element $\equiv 2
  \mod I$ of norm $<7$. Since the only ideals of norm $<7$ are the
  ideals $R$, $(\pm 3 + \sqrt{14})$, and possibly powers of
  $P$, we have $N(a) > N(P)^2$ for every $a \equiv 2 \mod I$. For
  $N(a) < 7$ to hold we must therefore choose $c^2 = N(P)^2 < 7$.

  Now we consider $I = (2) = P^2$; if the residue class $1 +
  \sqrt{14} \mod I$ is to contain an element with norm $< c^2$, we must
  have $\delta < c^2$, because $J = R$ is the only ideal
  with $I + J = R$ and $N(J) < 5$.

  Thus the norm weighted by $N(4 + \sqrt{14}) = c$ can at most
  be a Euclidean function on $R$ if $\sqrt{5} < c < \sqrt{7}$
  holds. The simplest possible choice of $c$ is therefore
  probably $c = \sqrt{6}$. Indeed, by a quite
  different route, Bedocchi found this same function in 1985 and conjectured that it
  is Euclidean on $R$.
\end{enumerate}

\subsection*{4.}

\begin{enumerate}
\item[4.1.] In quadratic number rings $K$ with $d = \disc K$ we have,
  as is well known,
  $$ \frac{\sqrt{d}}{16 + 6\sqrt{6}} \le M(K) \le \frac{\sqrt{d}}4. $$
  In cubic number fields of unit rank 1, on the other hand, we know
  $$ \frac{\sqrt{d}}{420} \le M(K) \le \frac{|d|^{2/3}}{16^{3/2}}, $$
  where the exponent $2/3$ and the factor $1/16^{3/2}$ are best possible
  (though this last inequality is proved only for
  $\disc K < -1236$). In analogy with the
  quadratic case, we could now conjecture that the lower
  bound $\sqrt{d}/420$ can be improved to $k_1 |d|^{2/3}$ for some
  real $k_1$. Van der Linden writes in his proof of the
  ``Cassels bound'' $\sqrt{d}/420$ that there is a certain asymmetry
  between the proof in the
  cubic and in the quadratic case, and he conjectures that this is due to
  some of his
  estimates in the cubic case not being best possible.
  If the exponent $\frac12$ in the Cassels bound
  can indeed be improved to $\frac23$, this would probably also
  have consequences for the totally real cubic case (respectively for
  number fields of arbitrary degree and arbitrary signature in general):
  here we know that $M(K) \le \sqrt{d}/8$ holds and that this
  bound is best possible, so in this case too we conjecture
  that there exists a real $k_2$ with $k_2 \sqrt{d} \le M(K)$.
  Finally, we are led to the conjecture that only
  finitely many norm-Euclidean number fields exist for a given degree
  $n$.
\item[4.2.] Show (e.g.\ with a computer) that the inequality $M(K)
  \le |d|^{2/3}/16^{3/2}$ also holds for the fields with $-300 \ge
  \disc K \ge -1236$. For the fields with $-23 \ge
  \disc K \ge -199$ this follows from the values for
  $M(K)$ that we determined in \S\ 4. For the Euclidean
  fields with $-200 \ge \disc K \ge -1236$, too, the
  inequality follows immediately.
\item[4.3.] Determine $M(K)$ for the cubic field with
  $\disc K = -87$.
\item[4.4.] Are the cubic cyclic fields with $\sqrt{d} = 103,
  109, 127$ and $157$ norm-Euclidean?
\item[4.5.] Do there exist cubic cyclic fields with Euclidean algorithm and $\sqrt{d} \ge
  5 \cdot 10^5$?
\end{enumerate}

\subsection*{5.}

\begin{enumerate}
\item[5.1.] Find the best upper bound for $M(K)$, where $K$ is a totally
  complex number field of degree 4. According to Davenport and
  Swinnerton-Dyer $M(K) \le c \cdot d^{3/4}$ is best possible; however,
  they gave no proof of this.
\item[5.2.] Does the estimate
  $M(K) \le \sqrt{d}/16$ hold for arbitrary Dirichlet fields?
\item[5.3.] Does $M(K) \le c \sqrt{d}$ hold for quadratic extensions of
  imaginary quadratic
  number fields?
\item[5.4.] Find a class of totally complex biquadratic
  number fields for which $M(K) \le c \cdot d^{3/4}$ is best possible.
\item[5.5.] Finish the classification of the pure biquadratic
  number fields with Euclidean algorithm.
\item[5.6.] Do there exist only finitely many pure number fields of
  2-power degree?
\end{enumerate}

\subsection*{6.}

\begin{enumerate}
\item[6.1.] Determine $M(K)$ for $D(-1,13)$, $D(-1,17)$, $D(-3,2)$,
  $D(-3,-11)$ and $D(-3,-19)$.
\item[6.2.] Determine $M(K)$ for further classes of bicyclic
  number fields, for example for $D(q,m)$, $q = -1, -2, -3$, $m = n^2
  \pm r$, $r|4n$.
\end{enumerate}

\subsection*{7.}

\begin{enumerate}
\item[7.1.] Finish the classification of the norm-Euclidean
  Dirichlet number fields.
\item[7.2.] Determine $M(K)$ for further classes of Dirichlet
  number fields.
\end{enumerate}

\subsection*{8.}

\begin{enumerate}
\item[8.1.] Find all norm-Euclidean number fields of degree 4 that
  possess a quadratic subfield.
\item[8.2.] Show that there exist only finitely many norm-Euclidean number fields
  $K$ with $(K:\mathbb{Q}) = 4$, $K/\mathbb{Q}$ cyclic,
  $\disc K = 125p^3$ (this would refute a conjecture of
  Heilbronn).
\end{enumerate}

\subsection*{9.}

\begin{enumerate}
\item[9.1.] Are the fields of $p$-th roots of unity ($p = 17, 19$)
  norm-Euclidean?
\end{enumerate}

\chapter*{\S\ 11 Description of the Computer Programs}
\setcounter{chapter}{11}
\addcontentsline{toc}{chapter}{\S\ 11 Description of the Computer Programs}
\markboth{Euclidean Rings}{\S\ 11 Description of the Computer Programs}

Let an algebraic number field $K$ of degree $n$ be given with an order
$R$ that is generated over $\mathbb{Z}$ by the elements $\alpha_1,
\dots, \alpha_n$ (if $R$ is equal to the maximal order, then
$\{\alpha_1, \dots, \alpha_n\}$ forms an integral basis). Further let
independent units $u_1, \dots, u_l$ be given, as well as a $k \in R$,
and asked is for the exceptional sets with respect to $k$.

To this end we prescribe by 
\[
S = \{\alpha \in R : \alpha = \sum a_i \alpha_i \, , \, -a \leq a_i \leq a \}
\]
a subset $S$ of $R$ (e.g.\ with $a=10$, $10^2$, \dots). Then
we subdivide the fundamental domain $F = (-\frac{1}{2},
\frac{1}{2}) \times \dots \times (-\frac{1}{2}, \frac{1}{2})$ (we
adopt the notation of § 2) into small $n$-dimensional cubes of
edge length 0.1 (or also 0.25). For every such cube $W$ we now
proceed as follows:

\begin{enumerate}
\item Check whether for all $x \in W$ there exists a $y \in S$ with $N(x-y)
  \leq k$; if yes, then check the next cube, otherwise go
  to 2.
\item Subdivide $W$ into $2^n$ subcubes by halving the
  edge length. If the edge length is now smaller than a prescribed
  bound (e.g.\ $10^{-2}$, $10^{-3}$, etc.), then we count $W$ as belonging to the
  exceptional set and with the next cube go back to
  1.\ (with the last cube we go to 3.); otherwise we go
  with every individual subcube back to 1.  
\item At this point we have found all subcubes $W$ that
  possibly contain an exceptional point (if no cube
  has remained, one has $M(K) < k$). We now check for which
  cubes $W$ the transforms $u_i W \bmod R$ lie in covered region
  (by testing whether $u_i W$ has a point in common with one of the other $W$).
  Every cube $W$ for which some $u_i W$ lies entirely in
  covered region we can forget.
    
  Now it is conceivable that some $u_i W_1 \bmod R$ intersects a cube $W_2$
  which at a later time is ``forgotten'',
  because $u_j W_2 \bmod R$ lies in covered region. We therefore repeat
  step 3.\ until no further change results.
\end{enumerate}

Now we check (by hand) whether the theorems from § 2 can be applied.
If this is not the case, then there are the following possibilities:

\begin{enumerate}
\item[a)] Decrease the bound for the edge length of the subcubes in 2.
\item[b)] Enlarge $S$.
\item[c)] If the $u_i$ do not form a system of fundamental units, one can
  try to find such a system, thereby replace the $u_i$ and
  repeat step 3.
\item[d)] Enlarge $k$ to $k_1$; it is indeed quite conceivable
  that in the first attempt one has chosen a much too small $k$.
\end{enumerate}

We must still show how one estimates the norm on a given
subcube $W$ from above. To this end let $W$ (or also a $W$ shifted by a
$y \in S$) be given by the inequalities $r_1 \leq x_1 \leq
s_1, \dots, r_n \leq x_n \leq s_n$ (here $-\frac{1}{2}
\leq r_i < s_i \leq \frac{1}{2}$). We must then find an upper bound
for $N(\sum x_i \alpha_i)$. As is well known however $N(z) = |z_1|
\cdot \dots \cdot |z_r| \cdot \dots \cdot |z_r+s|$ (for the definition
of the $| \cdot |_j$ see § 2), so that it suffices to estimate each factor $|z_j|$ from
above.

If $| \cdot |_j$ is defined by $|z_j| = |\sum x_i \sigma_j(\alpha_i)|$
and $\sigma_j$ is a real embedding, then one has
\[
\sup \{|z_j| : z \in W\} = \max \{ |\sum x_i \sigma_j(\alpha_i)| : x_i \in \{r_i, s_i\} \},
\]
i.e.\ one can find for $|z_j|$ with $2^n$ evaluations an upper bound
if one only knows the real numbers $\sigma_j(\alpha_i)$ (sufficiently
accurately). If $\sigma_j$ on the other hand is a non-real embedding,
then one can quite analogously estimate the real and imaginary parts
of $\sum x_i \sigma_j(\alpha_i)$, or one computes the real
coefficients of $x_i x_k$ in the product $\sum x_i \sigma_j(\alpha_i)
\sum x_i \overline{\sigma_j(\alpha_i)}$ and estimates this product
from above.

As an example we indicate how this procedure looks in the case $n=3$,
$r=s=1$; to this end we assume that $K$ is given by a polynomial $f$
of index $1$ (if the index is $\neq 1$, then one must modify the following
accordingly): Input for the function NORMAX, which estimates the
absolute value of the norm on a cube $W$ from above, are here $r_1,
s_1, t_1$, the ``step size'' $sw$, the real root $x$ of
$f$, as well as the complex roots $w_1 \pm i w_2$. We then form
$r_2 = r_1 + sw$, $s_2 = s_1 + sw$, $t_2 = t_1 + sw$ and $y = w_1^2 -
w_2^2$ (since $1, x, x^2$ form an integral basis, one has $\sigma(x^2) =
w_1^2 - w_2^2 + 2 w_1 w_2 i$). If finally $x_1$ denotes the upper
bound for $|z_1|$, $y_2$ the upper bound for the real part, $y_3$
that for the imaginary part of $z'$ (respectively $z''$), then we have in
GFA-Basic:

\begin{verbatim}
Function NORMAX(r1,s1,t1,sw)
    r2=r1+sw
    s2=s1+sw
    t2=t1+sw
    x2=x^2
    y=w1^2 - w2^2
    x1=MAX(ABS(r1*s1*x+t1*x2),ABS(r1*s1*x+t2*x2),ABS(r1*s2*x+t1*x2))
    x1=MAX(ABS(r1*s2*x+t2*x2),ABS(r2*s1*x+t1*x2),ABS(r2*s1*x+t2*x2),x1)
    x1=MAX(ABS(r2*s2*x+t1*x2),ABS(r2*s2*x+t2*x2), x1)
    y2=MAX(ABS(r1*s1*w1+t1*y),ABS(r1*s1*w1+t2*y),ABS(r1*s2*w1+t1*y))
    y2=MAX(ABS(r1*s2*w1+t2*y),ABS(r2*s1*w1+t1*y),ABS(r2*s1*w1+t2*y),y2)
    y2=MAX(ABS(r2*s2*w1+t1*y),ABS(r2*s2*w1+t2*y), y2)
    y3=MAX(ABS(s1+2*t1*w1),ABS(s1+2*t2*w1),ABS(s2+2*t1*w1),ABS(s2+2*t2*w2))
    y3=y3*abs(w2)
    y1=y2^2 + y3^2
    RETURN x1*y1
ENDFUNC
\end{verbatim}

$y_1 = y_2^2 + y_3^2$ is thus an upper bound for $|z_2|$, and $x_1
y_1$ such a one for the absolute value of the norm on $W$.

In practice the program runs as follows: after input of a generating
polynomial (with coefficients as small as possible; this is of importance for the
error estimate) one computes with the help of the
Newton method the real root $x$ of $f$ with an
accuracy of at least 10 decimal places, and after
dividing out this root obtains a quadratic polynomial of negative
discriminant, from which one can easily determine the values $w_1$ and $w_2$.
We describe the set $S$ by three arrays $a(n)$, $b(n)$
and $c(n)$, which denote the (integral) coordinates with respect to $1$, $x$ and $x^2$;
at the beginning we start with only one $y \in S$, namely
$y=0$, i.e.\ we set $n=0$, $a(0)=b(0)=c(0)=0$. Then we choose
e.g.\ the step sizes $sw_1=0.1$ and $sw_2=0.05$ and divide $F=(0,
0.5)\times(-0.5, 0.5)\times(-0.5, 0.5)$ into small subcubes of
edge length $sw_1$. Then we estimate the norm on $y+W$ for every
subcube $W$ and all $y \in S$ with NORMAX; if the estimate
is larger than the entered bound, then we halve the edge lengths
of $W$ and repeat the step for each of the eight subcubes $V$
of $W$. If the bound for one of these subcubes for every $y \in
S$ is larger than the bound, then we try to find a $y \in R$
such that the norm on $y+V$ becomes smaller than the bound. This we do
with the procedure SEARCH, by trying every $y \in R$
whose coordinates have absolute value $\leq 8$:

\begin{verbatim}
PROCEDURE SEARCH(r1,s1,t1,sw)
    FOR e = -8 TO 8
        FOR f = -8 TO 8
            FOR g = -8 TO 8
                norm=NORMAX(r1+e,s1+f,t1+g,sw)
                IF norm<schranke
                    n=n+1
                    a(n)=e
                    b(n)=f
                    c(n)=g
                    e=9
                    f=9
                    g=9
                ENDIF
            NEXT g
        NEXT f
    NEXT e
    RETURN
ENDPROC
\end{verbatim}

\section*{The Main Program}

The main program thereby looks as follows:
\lstset{
  basicstyle=\ttfamily\small,
  breaklines=true,
  columns=fullflexible
}
\begin{lstlisting}
FOR a1=0 TO 0.49 STEP sw1 
FOR b1=0 TO 0.49 STEP sw1 
FOR c1=0 TO 0.49 STEP sw1 
index=-1 
REPEAT 
  index=index+1 
  norm1=NORMAX(a1+a(index),b1+b(index),c1+c(index),sw1) 
UNTIL index=n or norm1<schranke 
IF norm1>schranke 
  FOR a2=a1 TO a1+0.049 STEP sw2 
  FOR b2=b1 TO b1+0.049 STEP sw2 
  FOR c2=c1 TO c1+0.049 STEP sw2 
  index=-1 
  REPEAT 
    index=index+1 
    norm2=NORMAX(a2+a(index),b2+b(index),c2+c(index),sw2) 
  UNTIL index=n or norm2<schranke 
  IF norm2>schranke 
    SEARCH(a2,b2,c2,sw2) 
    IF NORMAX(a2+a(n),b2+b(n),c2+c(n),sw2)>schranke 
      SCHREIB(a2,b2,c2,sw2) 
    ENDIF 
  ENDIF 
  NEXT c2 
  NEXT b2 
  NEXT a2 
ENDIF 
NEXT c1 
NEXT b1 
NEXT a1
\end{lstlisting}

Here \textbf{SCHREIB} is a procedure that outputs the considered
subcube $V$ together with the transformed cubes $uV$ for given
units $u$ (e.g.\ the fundamental unit and its inverse) to screen, printer
or disk. How one has to supplement this program if
one wishes to subdivide the subcubes further is clear.

\chapter*{Tables}
\setcounter{chapter}{12}
\addcontentsline{toc}{chapter}{Tables}
\markboth{Euclidean Rings}{Tables}

\section*{Cubic Number Fields of Negative Discriminant}
\addcontentsline{toc}{section}
                     {Cubic Number Fields of Negative Discriminant}

We use the following notation:
\begin{enumerate}
\item[] $\disc K$: the discriminant of the number field
\item[] E: $K$ is norm-Euclidean
\item[] NE: $K$ is not norm-Euclidean
\item[] H: $K$ has class number different from $1$
\item[] $a_2, a_1, a_0$: the coefficients of a generating polynomial
  \[ f(x) = x^3 + a_2 x^2 + a_1 x + a_0 \]
\item[] integral basis: gives an integral basis, if not
  $\{1, \vartheta, \vartheta^2\}$ is such,
  where $\vartheta$ is a root of $f$
\item[] $h$: the class number of the field
\item[] $u$: a unit
\end{enumerate}

$$ 
\begin{array}{rrrrrrrrcc}
\toprule
Nr. & \disc K & & a_2 & a_1 & a_0 & GHB & h & u & M(K) \\
\midrule
1. & -23 & E & 0 & -1 & & 1 & 1 & \vartheta & \frac{1}{5} \\
2. & -31 & E & 0 & 1 & & 1 & 1 & \vartheta & \frac{1}{3} \\
3. & -44 & E & 1 & -1 & & 1 & 1 & \vartheta & \frac{1}{2} \\
4. & -59 & E & 0 & 2 & & 1 & 1 & \vartheta & \frac{1}{2} \\
5. & -76 & E & 3 & 1 & & 1 & 1 & \vartheta & \frac{1}{2} \\
6. & -83 & E & 2 & -2 & & 1 & 1 & \vartheta & \frac{1}{2} \\
7. & -87 & E & 2 & -1 & & 1 & 1 & \vartheta & \geq \frac{1}{3} \\
8. & -104 & E & 0 & -1 & & 2 & 1 & 3+2\vartheta & \frac{1}{2} \\
9. & -107 & E & 4 & 2 & & 1 & 1 & \vartheta & \frac{1}{2} \\
10. & -108 & E & 3 & -3 & & 1 & 1 & \vartheta & \frac{1}{2} \\
\midrule
11. & -116 & E & 1 & 0 & & 2 & 1 & 1 - \vartheta + \vartheta^2 & \frac{1}{2} \\
12. & -135 & E & 3 & 0 & & 1 & 1 & \vartheta & \frac{3}{5} \\
13. & -139 & E & 4 & 6 & & 1 & 1 & \vartheta & \frac{1}{2} \\
14. & -140 & E & 5 & 3 & & 1 & 1 & \vartheta & \frac{1}{2} \\
15. & -152 & E & 1 & -2 & & 2 & 1 & 3 - \vartheta - \vartheta^2 & \frac{1}{2} \\
16. & -172 & E & 2 & 0 & & 2 & 1 & 1 + \vartheta^2 & \frac{3}{4} \\
17. & -175 & E & -2 & 3 & & 1 & 1 & \vartheta & \frac{3}{5} \\
18. & -199 & NE& 4 & 1 & & 1 & 1 & \vartheta & 1 \\
19. & -200 & E & 1 & 2 & & -2& 1 & 1 - \vartheta - \vartheta^2 & <0.6 \\
20. & -204 & E & 1 & 1 & & 3 & 1 & 1 - \vartheta - \vartheta^2 & <0.64 \\
\bottomrule
\end{array} $$

\clearpage

\small
$$ \begin{array}{rrrrrrrrcc}
\toprule
Nr. & \disc K & & a_2 & a_1 & a_0 & GHB & h & u & M(K) \\
\midrule
21. & -211 & E & 3 & 1 & 2  & & 1 & 3+\vartheta & \\
22. & -212 & E & 1 & 4 & 2  & &  1 & 1+2\vartheta & \frac{5}{8} \\
23. & -216 & E & 0 & 3 & 2  & &  1 & 1+\vartheta-\vartheta^2 & \\
24. & -231 & E & -4 & 5 & 1 & &  1 & \vartheta & \geq \frac{7}{9} \\
25. & -239 & E & 0 & -1 & 3 & &  1 & 5+3\vartheta & \frac{8}{9} \\
26. & -243 & E & 0 & 0 & 3 & &  1 & 2-\vartheta^2 & \\
27. & -244 & E & 2 & -3 & 2 &  & 1 & 1-2\vartheta+2\vartheta^2 & \\
\bottomrule
\end{array} $$

\clearpage

$$
\begin{array}{rrrrrrrlrc}
\toprule
Nr. & \disc K & a_2 & a_1 & a_0 & GHB & h & u & M(K) & Typ \\
\midrule
61. & -516 & 4 & 1 & 6 & \frac{\vartheta+\vartheta^2}{2} & 1 & -16+(21\vartheta+7\vartheta^2)/2 & & E \\
62. & -519 & 2 & -3 & 3 &  & 1 & 1-\vartheta+\vartheta^2 & \geq \frac{44712}{45747} & E \\
63. & -524 & -2 & 4 & 2 &  & 1 & 1+2\vartheta-\vartheta^2 & & NE \\
64. & -527 & 0 & 5 & 1 &  & 1 & \vartheta & & NE \\
65. & -543 & -2 & 3 & 3 &  & 1 & 1+\vartheta-\vartheta^2 & \geq \frac{7593}{8343} & E \\
66. & -547 & 1 & -3 & 4 &  & 1 & 9-2\vartheta-2\vartheta^2 & & NE \\
67. & -563 & 2 & 6 & 1 &  & 2 & \vartheta & & H \\
\bottomrule
\end{array} $$

\clearpage

$$ \begin{array}{rrrrrlrrrl}
\toprule
Nr. & \disc K & a_2 & a_1 & a_0 & \text{GHB} & h & u & M(K) & Typ \\
\midrule
101 & -823 & 3 & -2 & 3 & & 1 & 1-\vartheta+\vartheta^2 & & NE \\
102 & -835 & 2 & 0 & 5 & & 1 & 18+4\vartheta-\vartheta^2 & & \\
103 & -839 & 4 & 3 & 5 & & 1 & 3+\vartheta+2\vartheta^2 & & NE \\
104 & -843 & 7 & -7 & 2 & & 1 & 3-11\vartheta+12\vartheta^2 & & \\
105 & -856 & -1 & 1 & 11 & frac{1+\vartheta^2}{2} & 1 & 5+10\vartheta+4\vartheta^2 & & \\
106 & -863 & 7 & -4 & 1 & & 1 & \vartheta & & NE \\
107 & -867 & 7 & 5 & 2 & & 1 & 9+21\vartheta+28\vartheta^2 & & NE \\
\bottomrule
\end{array} $$

\clearpage

$$ \begin{array}{rrlrrrlrrr}
\toprule
\text{Nr.} & \disc K & & a_2 & a_1 & a_0 & \text{GHB} & h & u & M(K) \\
\midrule
141 & -1099 & & 1 & -1 & 6 & & & & \\
142 & -1107 & & 0 & 6 & 3 & & & & \\
143 & -1108 & & 5 & 0 & 2 & & 1 & & \\
144 & -1135 & NE & 5 & 2 & 3 & & 1 & & \\
145 & -1144 & & -1 & 6 & 2 & & 1 & & \\
146 & -1147 & NE & 2 & -2 & 5 & & 1 & & \\
147 & -1164 & NE & -1 & -1 & 7 & & 1 & & \\
148 & -1172 & NE & 5 & -4 & 2 & & 1 & & \\
149 & -1175 & NE & 3 & 8 & 1 & & 1 & & \\
150 & -1176 & & 1 & -2 & 6 & & & & \\
\midrule
151 & -1187 & NE & 1 & 7 & 2 & & 1 & & \\
152 & -1188 & & 6 & -3 & 4 & \frac{\vartheta+\vartheta^2}{2} & 1 & & \\
153 & -1191 & NE & 4 & -5 & 3 & & 1 & & \\
154 & -1192 & & 7 & -8 & 4 & \frac{\vartheta+\vartheta^2}{2} & 1 & & \\
155 & -1203 & & 4 & 2 & 5 & & 1 & & \\
156 & -1207 & NE & 4 & -1 & 3 & & 1 & & \\
157 & -1208 & NE & 4 & 13 & 2 & \frac{\vartheta+\vartheta^2}{2} & 1 & & \\
158 & -1219 & NE & 3 & -5 & 4 & & 1 & & \\
159 & -1231 & NE & 2 & -3 & 5 & & 1 & & \\
160 & -1235 & NE & 2 & -4 & 5 & & 1 & & \\
\midrule
161 & -1236 & & 7 & -3 & 3 & \frac{1+\vartheta^2}{2} & 1 & & \\
162 & -1255 & & 4 & 9 & 1 & & & & \\
163 & -1259 & NE & 5 & -1 & 2 & & 1 & & \\
164 & -1267 & NE & 1 & 7 & 4 & & 1 & & \\
165 & -1272 & & -7 & 13 & 9 & \frac{1+\vartheta^2}{2} & 1 & & \\
166 & -1291 & NE & 3 & 1 & 6 & & 1 & & \\
167 & -1292 & NE & 6 & 2 & 2 & & 1 & & \\
168 & -1295 & NE & -4 & 5 & 5 & & 1 & & \\
169 & -1300 & & 9 & -8 & 2 & & 1 & & \\
170 & -1315 & NE & 5 & -3 & 2 & & 1 & & \\
\midrule
171 & -1316 & NE & 7 & 4 & 2 & & 1 & & \\
172 & -1319 & NE & 11 & -6 & 1 & & 1 & & \\
173 & -1323 & & 12 & 6 & 1 & & & & \\
174 & -1327 & NE & 0 & 1 & 7 & & 1 & & \\
175 & -1336 & & 5 & -2 & 2 & & 1 & & \\
176 & -1347 & & 4 & -2 & 3 & & 1 & & \\
177 & -1351 & NE & 1 & 0 & 7 & & 1 & & \\
178 & -1355 & NE & 0 & 2 & 7 & & 1 & & \\
179 & -1363 & NE & 4 & -4 & 3 & & 1 & & \\
180 & -1371 & & -4 & 8 & 1 & & & & \\
\bottomrule
\end{array} $$

\clearpage

\section*{Totally real cubic fields}
\addcontentsline{toc}{section}
                {Totally real cubic fields}

\begin{table}[h]
\centering
\small
\begin{tabular}{rrrrrrccc}
\toprule
Nr. & disc K & $a_2$ & $a_1$ & $a_0$ & GHB & h & M(K) \\
\midrule
1.  & 49   & 1  & -2 & -1 &  & 1 & $1/7$   & E \\
2.  & 81   & 0  & -3 & 1  &  & 1 & $1/3$   & E \\
3.  & 148  & 3  & -1 & -1 &  & 1 & $1/2$   & E \\
4.  & 169  & 1  & -4 & 1  &  & 1 & $5/13$  & E \\
5.  & 229  & 0  & -4 & 1  &  & 1 & $1/2$   & E \\
6.  & 257  & 3  & -2 & -1 &  & 1 & $1/3$   & E \\
7.  & 316  & 1  & -4 & -2 &  & 1 & $1/2$   & E \\
8.  & 321  & 1  & -4 & -1 &  & 1 & $1/3$   & E \\
9.  & 361  & -2 & -5 & -1 &  & 1 & $8/19$  & E \\
10. & 404  & 1  & -5 & 1  &  & 1 & $1/2$   & E \\
\midrule
11. & 469  & 4  & -2 & -1 &  & 1 & $1/2$   & E \\
12. & 473  & 0  & -5 & 1  &  & 1 & $1/3$   & E \\
13. & 564  & -7 & 11 & 1  &  & 1 & $1/2$   & E \\
14. & 568  & 1  & -6 & 2  &  & 1 & $1/2$   & E \\
15. & 621  & -6 & 6  & 1  &  & 1 & $1/2$   & E \\
16. & 697  & 3  & -4 & -1 &  & 1 & $13/31$ & E \\
17. & 733  & -2 & -6 & -1 &  & 1 & $1/2$   & E \\
18. & 756  & 0  & -6 & 2  &  & 1 & $1/2$   & E \\
19. & 761  & 1  & -6 & 1  &  & 1 & $1/2$   & E \\
20. & 785  & -2 & -5 & 1  &  & 1 &         & E \\
\midrule
21. & 788  & 9  & 7  & 1  &  & 1 &         & E \\
22. & 837  & 0  & -6 & 1  &  & 1 & $1/2$   & E \\
23. & 892  & -5 & 0  & 2  &  & 1 &         & E \\
24. & 940  & 3  & -4 & -2 &  & 1 &         & E \\
25. & 961  & -2 & -9 & -2 & $\frac{\theta+\theta^2}{2}$ & 1 & $16/31$ & E \\
26. & 985  & 1  & -6 & -1 &  & 1 & 1       & NE \\
27. & 993  & 1  & -6 & -3 &  & 1 &         & E \\
28. & 1016 & -1 & -6 & 2  &  & 1 &         & E \\
29. & 1076 & -5 & -3 & 1  &  & 1 &         & E \\
30. & 1101 & -5 & -1 & 2  &  & 1 &         & E \\
\bottomrule
\end{tabular}
\end{table}

$$ \begin{array}{rrrrrrrrl}
\toprule
Nr. & \disc K & a_2 & a_1 & a_0 & GHB & h & M(K) & \text{type} \\
\midrule
31. & 1129 & -2 & -6 & 1 &  & 1 &        & E \\
32. & 1229 & -6 & -2 & 1 &  & 1 &        & E \\
33. & 1257 & 2  & -7 & 1 &  & 1 &        & E \\
34. & 1300 & 3  & -7 & 1 &  & 1 &        & E \\
35. & 1304 & 11 & 10 & 2 &  & 1 &        & E \\
36. & 1345 & 0  & -7 & 1 &  & 1 & 7/5    & NE \\
37. & 1369 & 4  & -7 & 1 &  & 1 &        & E \\
\bottomrule
\end{array} $$

$$ \begin{array}{rrrrrrrl}
\toprule
Nr. & disc K & a_2 & a_1 & a_0 & \text{GHB} & h & M(K) \\
\midrule
71. & 2241 & 0  & -9 & 5  &  & 1 & E \\
72. & 2292 & 10 & -10& 2  &  & 1 & E \\
73. & 2296 & -7 & 2  & 2  &  & 1 & E \\
74. & 2300 & -1 & -8 & 2  &  & 1 & E \\
75. & 2349 & -6 & 0  & 3  &  & 1 & E \\
76. & 2429 & -7 & 1  & 2  &  & 1 & E \\
77. & 2505 & -4 & -5 & 3  &  & 1 & E \\
78. & 2557 & 1  & -9 & 2  &  & 1 & E \\
79. & 2589 & 10 & -12& 3  &  & 1 & E \\
80. & 2597 & -8 & -2 & 1  &  & 3 & E \\
\midrule
81. & 2636 & -7 & 0  & 2  &  & 1 & E \\
82. & 2673 & 0  & -9 & 3  &  & 1 & E \\
83. & 2677 & -3 & -7 & 2  &  & 1 & E \\
84. & 2700 & -6 & -3 & 2  &  & 1 & E \\
85. & 2708 & -4 & -6 & 2  &  & 1 & E \\
86. & 2713 & -6 & -1 & 3  &  & 1 & E \\
87. & 2777 & -5 & -6 & 1  &  & 2 & E \\
88. & 2804 & 1  & -9 & 1  &  & 1 & E \\
89. & 2808 & 0  & -9 & 2  &  & 1 & E \\
90. & 2836 & -2 & -8 & 2  &  & 1 & E \\
\bottomrule
\end{array} $$

\clearpage

\section*{Totally complex biquadratic fields with quadratic subfield}
\addcontentsline{toc}{section}
    {Totally complex biquadratic fields with quadratic subfield}

$$ \begin{array}{rrclll}
\toprule
Nr. & \disc K & \disc k & \mu & M(K) & \text{type} \\
\midrule
1.  & 117 & -3     & -1+2\sqrt{-3}          & \geq 1/7 & E \\
2.  & 125 & 5      & -5+2\sqrt{5}           & \geq 1/5 & E \\
3.  & 144 & -3,-4,12 &                        & 1/4      & E \\
4.  & 189 & -3     & 3+2\sqrt{-3}           & \geq 1/3 & E \\
5.  & 225 & -3,-15,5 &                      & 1/4      & E \\
6.  & 256 & -4,-8,8 &                        & 1/2      & E \\
7.  & 272 & -4     & 1+4i                   & \geq 1/4 & E \\
8.  & 320 & -4     & 1+2i                   & 1/2      & E \\
9.  & 333 & -3     & -5+2\sqrt{-3}          &          & E \\
10. & 392 & -7     & (-1-\sqrt{-7})/2       &          & E \\
\midrule
11. & 400 & -4,-20,5 &                      & 5/16     & E \\
12. & 432 & -3     & \sqrt{-3}              &          & E \\
13. & 441 & -3,-7,21 &                      & 4/9      & E \\
14. & 512 & -4     &                        & 1/2      & E \\
15. & 513 & -3     & -3+4\sqrt{-3}          &          & E \\
16. & 549 & -3     & 7-2\sqrt{-3}           &          & E \\
17. & 576 & -3,24,-8 &                      & 1/3      & E \\
18. & 576 & -3,-24,8 &                      &          & E \\
19. & 605 & -11    & (3+\sqrt{-11})/2       &          & E \\
20. & 656 & -4     & 5+4i                   & 1/2      & E \\
\midrule
21. & 657 & -3     & 5-4\sqrt{-3}           & 1/2      & E \\
22. & 784 & -4,-7,28 &                      & 1/2      & E \\
23. & 832 & -4     & 3+2i                   & 1/2      & E \\
24. & 837 & -3     & -9+2\sqrt{-3}          &          & E \\
25. & 873 & -3     & -7+4\sqrt{-3}          &          & E \\
26. & 981 & -3     & -1+6\sqrt{-3}          &          & E \\
27. & 1008& -3     & 2-\sqrt{-3}            &          & E \\
28. & 1008& -3     & -2+\sqrt{-3}           &          & E \\
29. & 1025& 5      & -11+4\sqrt{5}          &          & E \\
30. & 1040& -4     & 7+4i                   &          & E \\
\bottomrule
\end{array} $$

\clearpage

$$ \begin{array}{rrclll}
\toprule
Nr. & \disc K & \disc k & \mu & M(K) & \text{type} \\
\midrule
31. & 1040 & -4     & 1+8i                    &      & E \\
32. & 1088 & -8     & 3+2\sqrt{-2}            &      & E \\
33. & 1088 & 8      & -5+2\sqrt{2}            &      & E \\
34. & 1089 & -3,-11,33 &                      &      & E \\
35. & 1161 & -3     & 9+4\sqrt{-3}            &      & E \\
36. & 1168 & -4     & 3+8i                    &      & E \\
37. & 1197 & -3     & 11+2\sqrt{-3}           &      & E \\
\bottomrule
\end{array} $$

\clearpage

$$ \begin{array}{rrccll}
\toprule
Nr. & \disc K & \disc k & \mu & M(K) & \text{type} \\
\midrule
71. & 2304 & -4,-24,24 &         & 5/2  & NE \\
72. & 2304 & -4        & 3i      & 5/2 & NE \\
73. & 2304 & -8,12,-24 &         &     & E   \\
74. & 2320 & -4        & 1+12i   & & \\
75. & 2320 & -4        & 9+8i    & & \\
76. & 2368 & -4        & 1+6i    & & E \\
77. & 2448 & -4        & 3+12i   & & NE \\
78. & 2457 & -3        & 9+8\sqrt{-3} & & \\
79. & 2493 & -3        & -13+6\sqrt{-3} & & \\
80. & 2560 & -4        & 1+3i    & 5/4 & NE \\
\midrule
81. & 2560 & -4        & 3+i     & 5/4 & NE \\
82. & 2597 & -7        & -5+2\sqrt{-7} & & \\
83. & 2601 & -3,17,-51 &               & 13/16   & E \\
84. & 2624 & -8        & -3+4\sqrt{-2} & & \\
85. & 2624 & 8         & -7+2\sqrt{-2} & & \\
86. & 2704 & -4,13,-52 &               & \geq 1  & NE \\
87. & 2709 & -3        & 1+10\sqrt{-3} & & \\
88. & 2709 & -3        & 17+2\sqrt{-3} & & \\
89. & 2736 & -3        & 4+\sqrt{-3} & & \\
90. & 2744 & -7        & (7+\sqrt{-7})/2 & &\\
\midrule
91. & 2781 & -3        & 3+10\sqrt{-3} & & \\
92. & 2817 & -3        & -11+8\sqrt{-3} & & \\
93. & 2873 & 13        & -15+4\sqrt{13} & & \\
94. & 2880 & -4        & 3+6i    & & NE \\
95. & 2925 & -3        & -5+10\sqrt{-3} & & \\
96. & 2960 & -4        & 13+4i   & & \\
97. & 2960 & -4        & 11+8i   & & \\
98. & 3024 & -3        & -3+2\sqrt{-3} & &  \\
99. & 3033 & -3        & 17+4\sqrt{-3} & & \\
100.& 3072 & -8        & 1+\sqrt{-2} & & \\
\midrule
101.& 3072 & -8        & -1-\sqrt{-2} & & \\
102.& 3088 & -4        & 7+12i   & & NE \\
103.& 3136 & -7,-8,56  &         & 9/8     & NE \\
104.& 3141 & -3        & 7+10\sqrt{-3} & & \\
105.& 3249 & -3,-19,57 &         & & E \\
106.& 3328 & -4        & 2+3i    & & \\
107.& 3357 & -3        & 19+2\sqrt{-3} & & \\
108.& 3392 & -4        & 7+2i    & 50/53 & E \\
109.& 3429 & -3        & -9+10\sqrt{-3} & &  \\
110.& 3573 & -3        & 17+6\sqrt{-3} & & \\
\bottomrule
\end{array} $$

\newpage

$$ \begin{array}{rccll}
\toprule
\text{Nr.} & \disc K & \disc k & \mu & \\
\midrule
111. & 3600 & -4,-15,60 &  & NE \\
112. & 3648 & -8 & 7 + 2\sqrt{-2} & \\
113. & 3648 & -8 & 5 + 4\sqrt{-2} & \\
114. & 3681 & -3 & -19 + 4\sqrt{-3} & \\
115. & 3725 & 5 & -13 - 2\sqrt{5} & NE \\
116. & 3728 & -4 & 13 + 8i & NE \\
117. & 3753 & -3 & -15 + 8\sqrt{-3} & \\
118. & 3773 & -7 & 7 + 2\sqrt{-7} & \\
119. & 3789 & -3 & 11 + 10\sqrt{-3} & \\
120. & 3856 & -4 & 15 + 4i & \\
\midrule
121. & 3897 & -3 & 1 + 12\sqrt{-3} & \\
122. & 3904 & -4 & 5 + 6i & NE \\
123. & 4032 & -3 & 4 + 2\sqrt{-3} & \\
124. & 4032 & -3 & -4 - 2\sqrt{-3} & \\
125. & 4077 & -3 & -21 + 2\sqrt{-3} & \\
126. & 4112 & -4 & 1 + 16i & \\
127. & 4113 & -3 & 5 - 12\sqrt{-3} & \\
128. & 4221 & -3 & -13 + 10\sqrt{-3} & \\
129. & 4221 & -3 & 19 + 6\sqrt{-3} & \\
130. & 4240 & -4 & 3 + 16i & \\
\midrule
131. & 4240 & -4 & 11 + 12i & \\
132. & 4329 & -3 & -7 + 12\sqrt{-3} & \\
133. & 4352 & -4 & 4 + i & NE \\
134. & 4352 & -8 & -3 + 2\sqrt{-2} & \\
135. & 4400 & -5 & -8 + 2\sqrt{5} & \\
136. & 4401 & -3 & 21 + 4\sqrt{-3} & \\
137. & 4464 & -3 & 2 + 3\sqrt{-3} & \\
138. & 4464 & -3 & -2 - 3\sqrt{-3} & \\
139. & 4496 & -4 & 5 + 16i & NE \\
140. & 4508 & -7 & 4 + \sqrt{-7} & \\
\midrule
141. & 4624 & -4,17,-68 &  & NE \\
142. & 4672 & -8 & -1 + 6\sqrt{-2} & \\
143. & 4725 & -3 & 15 + 10\sqrt{-3} & \\
144. & 4761 & -3,-23,69 &  & NE \\
145. & 4869 & -3 & 23 + 2\sqrt{-3} & \\
146. & 4880 & -4 & 17 + 4i & \\
147. & 4880 & -4 & 7 + 16i & \\
148. & 4977 & -3 & -19 + 8\sqrt{-3} & \\
149. & 4977 & -3 & -11 + 12\sqrt{-3} & \\
150. & 5008 & -4 & 13 + 12i & \\
\bottomrule
\end{array} $$

\newpage

$$ \begin{array}{rccll}
\toprule
\text{Nr.} & \disc K & \disc k & \mu & \\
\midrule
151. & 5193 & -3 & -23 + 4\sqrt{-3} & \\
152. & 5301 & -3 & -17 + 10\sqrt{-3} & \\
153. & 5301 & -3 & -1 + 14\sqrt{-3} & \\
154. & 5328 & -3 & 5 + 2\sqrt{-3} & \\
155. & 5341 & -7 & -9 + 2\sqrt{-7} & \\
156. & 5392 & -4 & 9 + 16i & \\
157. & 5440 & -4 & 9 + 2i & \\
158. & 5488 & -7 & 2\sqrt{-7} & \\
159. & 5537 & -7 & 1 + 4\sqrt{-7} & \\
160. & 5648 & -4 & 17 + 8i & \\
\midrule
161. & 5696 & -8 & -9 + 2\sqrt{-2} & NE \\
162. & 5776 & -4,-19,76 & & NE \\
163. & 5929 & -7,-11,77 & & NE \\
164. & 6032 & -4 & 19 + 4i & \\
165. & 6032 & -4 & 11 + 16i & \\
166. & 6144 & -8 & 2 + \sqrt{-2} & \\
167. & 6144 & -8 & -2 + \sqrt{-2} & \\
168. & 6208 & -8 & -5 + 6\sqrt{-2} & \\
169. & 6400 & -4 & 5i & \\
170. & 6400 & -8,-20,40 & & NE \\
\midrule
171. & 6416 & -4 & 1 + 20i & \\
172. & 6464 & -4 & 1 + 10i & \\
173. & 6544 & -4 & 3 + 20i & \\
174. & 6656 & -4 & 5 + i & NE \\
175. & 6656 & -4 & -1 + 5i & NE \\
176. & 6713 & -7 & 5 + 4\sqrt{-7} & \\
177. & 6725 & 5 & -17 - 2\sqrt{5} & \\
178. & 6800 & -4 & 5 + 20i & NE \\
179. & 6928 & -4 & 17 + 12i & \\
180. & 6976 & -4 & 3 + 10i & NE \\
\midrule
181. & 7056 & -4,21,-84 & & NE \\
182. & 7184 & -4 & 7 + 20i & \\
183. & 7312 & -4 & 21 + 4i & \\
184. & 7424 & -4 & 2 + 5i & NE \\
185. & 7488 & -4 & 9 + 6i & NE \\
186. & 7696 & -4 & 9 + 20i & \\
187. & 7696 & -4 & 15 + 16i & NE \\
188. & 7744 & -8,-11,44 & & NE \\
189. & 8448 & -8 & 5 + 2\sqrt{-2} & \\
\bottomrule
\end{array} $$

\clearpage

\section*{Real biquadratic fields with quadratic subfields}

\medskip
Notation as above; in addition:
\begin{tabbing}
\hspace*{1.5cm}\= \kill
$V_4$\> the Klein four-group \\
$Z_4$\> the cyclic group of order 4 \\
$D_4$\> the dihedral group of order 8 \\
$A_4$\> the alternating group of order 12 \\
$S_4$\> the symmetric group of order 24
\end{tabbing}

\medskip
\[
\begin{array}{r@{\quad}r@{\quad}c@{\quad}l@{\quad}l@{\quad}c}
\hline
\text{Nr.} & \text{disc}\,K &  & \text{disc}\,k & \mu & \mathrm{Gal}(K/\Q) \\
\hline
1.  & 725  & \mathrm{E} & 5           & 7+2\sqrt{5}     & D_4 \\
2.  & 1125 & \mathrm{E} & 5           & 15+6\sqrt{5}    & Z_4 \\
3.  & 1600 & \mathrm{E} & 5,8,40      &                 & V_4 \\
4.  & 2000 &           & 5           & 5+2\sqrt{5}     & Z_4 \\
5.  & 2048 & \mathrm{E} & 8           & 2+\sqrt{2}      & D_4 \\
6.  & 2225 & \mathrm{E} & 5           & 13+4\sqrt{5}    & D_4 \\
7.  & 2304 & \mathrm{E} & 8,12,24     &                 & V_4 \\
\hline
\end{array}
\]

\clearpage

\section*{Survey of the Distribution of Norm-Euclidean Number Fields}
\addcontentsline{toc}{section}
    {Survey of the Distribution of Norm-Euclidean Number Fields}

The number of known norm-Euclidean number fields has in the last
twenty years strongly increased. The following tables, which show the
distribution of the known norm-Euclidean number rings with respect to
$n$ (degree of the field) and $r+s$ (unit rank + 1), are to make this
clear:

\bigskip

\textbf{I. 1967 (Kummer (1844), Eisenstein (1850), Godwin (1965b, 1967))}

\begin{center}
\begin{tabular}{c|ccccc|c}
\toprule
$r+s$ & 1 & 2 & 3 & 4 & 5 & $\Sigma$ \\
\midrule
1 & 1 & 5 & & & & 6 \\
2 & & 16 & 15 & 3 & & 34 \\
3 & & & 9 & 0 & 0 & 9 \\
4 & & & & 8 & 0 & 8 \\
5 & & & & & 1 & 1 \\
\midrule
$\Sigma$ & 1 & 21 & 24 & 11 & 1 & 58 \\
\bottomrule
\end{tabular}
\end{center}

\textbf{II. 1974 (Lenstra (1974), Smith (1969, 1971), Lakein (1972))}

\begin{center}
\begin{tabular}{c|cccccccccc|c}
\toprule
$r+s$ & 1 & 2 & 3 & 4 & 5 & 6 & 7 & 8 & 9 & 10 & $\Sigma$ \\
\midrule
1 & 1 & 5 & & & & & & & & & 6 \\
2 & & 16 & 15 & 31 & & & & & & & 62 \\
3 & & & 57 & 2 & 1 & 2 & & & & & 62 \\
4 & & & & 9 & 0 & 0 & 0 & 2 & & & 11 \\
5 & & & & & 1 & 1 & 0 & 0 & 0 & 1 & 3 \\
6 & & & & & & 0 & 0 & 0 & 0 & 0 & 0 \\
\midrule
$\Sigma$ & 1 & 21 & 72 & 42 & 2 & 3 & 0 & 2 & 0 & 1 & 144 \\
\bottomrule
\end{tabular}
\end{center}

\clearpage

\textbf{III. Sept. 1976 (Lenstra (1977), Taylor (1976))}

\begin{center}
\begin{tabular}{c|cccccccccc|c}
\toprule
$r+s$ & 1 & 2 & 3 & 4 & 5 & 6 & 7 & 8 & 9 & 10 & $\Sigma$ \\
\midrule
1 & 1 & 5 & & & & & & & & & 6 \\
2 & & 16 & 52 & 32 & & & & & & & 100 \\
3 & & & 57 & 11 & 12 & 28 & & & & & 108 \\
4 & & & & 9 & 10 & 25 & 23 & 24 & & & 91 \\
5 & & & & & 1 & 2 & 0 & 0 & 0 & 1 & 3 \\
6 & & & & & & 2 & 0 & 0 & 0 & 0 & 0 \\
\midrule
$\Sigma$ & 1 & 21 & 109 & 52 & 23 & 57 & 23 & 24 & 0 & 1 & 311 \\
\bottomrule
\end{tabular}
\end{center}

\textbf{VII. August 1989:} Besides the norm-Euclidean fields found in this work,
$\mathbb{Q}(\zeta_{13})$ is also listed in the
following table (see on this Leutbecher and Niklasch 1987)

\begin{center}
\begin{tabular}{c|cccccccccccc|c}
\toprule
$r+s$ & 1 & 2 & 3 & 4 & 5 & 6 & 7 & 8 & 9 & 10 & 11 & 12 & $\Sigma$ \\
\midrule
1 & 1 & 5 & & & & & & & & & & & 6 \\
2 & & 16 & 52 & 50 & & & & & & & & & 118 \\
3 & & & 62 & 13 & 13 & 28 & & & & & & & 116 \\
4 & & & & 15 & 10 & 37 & 39 & 45 & & & & & 146 \\
5 & & & & & 1 & 12 & 26 & 65 & 92 & 50 & & & 246 \\
6 & & & & & & 3 & 3 & 3 & 0 & 0 & 0 & 1 & 10 \\
\midrule
$\Sigma$ & 1 & 21 & 114 & 78 & 24 & 80 & 68 & 113 & 92 & 50 & 0 & 1 & 642 \\
\bottomrule
\end{tabular}
\end{center}

\printindex[N]
\printindex[S]

\selectlanguage{german}

\chapter*{Überblick}
\addcontentsline{toc}{chapter}{Überblick}
\markboth{Euklidische Ringe}{Überblick}

Zum besseren Verständnis der Arbeit sei im Folgenden kurz beschrieben,
was den Leser in den nächsten Paragraphen erwartet: Ausgangspunkt aller
Überlegungen war das Kriterium (1.15) von Lenstra\index[dN]{Lenstra}
\cite{Len74,Len77a}, wonach ein Zahlkörper mit einer hinreichend
langen \glqq{}Ausnahmefolge\grqq{} normeuklidisch ist. Etwa zur gleichen Zeit
(nämlich in den Jahren 1975 bis 1977) hat Cooke\index[dN]{Cooke} den
Begriff eines $k$-stufig normeuklidischen Zahlkörpers eingeführt. Dies
legte die Frage nahe, ob sich das Lenstra-Kriterium so modifizieren
lässt, dass man damit auch $k$-stufig normeuklidische Ringe finden
kann. Tatsächlich gelingt dies recht einfach (s. (1.16)). Das hierzu
notwendige Studium von Kettenbrüchen der Länge $\leq k$ mit
Koeffizienten aus dem betrachteten Ring führt auf Mengen, die denen
ähneln, die Motzkin\index[dN]{Motzkin} \cite{Mot49} eingeführt hat, und
die wir im § 0 vorstellen.

Damit hat man nun ein Kriterium, mit dem man (beispielsweise) 2-stufig
normeuklidische Zahlkörper vom Grad 2, 3, 4, 5 usw. finden kann, und
dies wirft die Frage auf, ob diese Körper nicht vielleicht schon
(1-stufig) normeuklidisch sind. Das allgemeinste Kriterium, mit dessen
Hilfe man zeigen konnte, dass ein gegebener Zahlring nicht
normeuklidisch ist, war (1.6); den Nachteil dieses Kriteriums
beschreibt van der Linden \index[dN]{Linden@van der Linden} \cite{Lin83} so:

\begin{quote}
  {\em If we generalize the methods for higher degree number fields it
  turns out that the arithmetical methods are applicable for
  extensions of $\Q$ in which at least one prime is totally
  ramified\dots.}
\end{quote}

Eine genaue Analyse des Beweises von (1.6) zeigt jedoch, dass sich
dieses Kriterium zumindest etwas verallgemeinern lässt (s. (1.4)),
dieweil es nämlich manchmal genügt, dass $K/\Q$ einen Zwischenkörper
$k$ enthält, sodass es in $K/k$ rein verzweigte Ideale gibt. Dass
diese Verallgemeinerung tatsächlich von Nutzen ist, zeigt sich
u. a. in § 6, wo wir alle normeuklidischen Zahlkörper der Form
$\Q(\sqrt{m},\sqrt{n})$ mit $m, n \in \Z$, $m < 0$ finden. Weiter kann
man mit (1.4) ggf. nicht nur die Existenz eines EA ausschließen,
sondern sogar Abschätzungen für das euklidische Minimum von $K$ geben.
Ein (von van der Linden\index[dN]{Linden@van der Linden} mit dem
Attribut ``geometrisch'' belegtes)
Kriterium, einen reellquadratischen Zahlkörper als nicht
normeuklidisch nachzuweisen, ist (1.8). In der hier vorgestellten Form
stammt es von Barnes\index[dN]{Barnes} und
Swinnerton-Dyer\index[dN]{Swinnerton-Dyer} \cite{BS52a}, jedoch lässt es
sich mehr oder weniger explizit bereits bei
Rédei\index[dN]{Redei@R\'edei} \cite{Red41a,Red41b} und
Inkeri\index[dN]{Inkeri}  \cite{Ink49} finden.
Zaghafte Verallgemeinerungen auf den kubischen Fall mit Einheitenrang
1 findet man bei Taylor\index[dN]{Taylor} \cite{Tay76} und bei
Cioffari\index[dN]{Cioffari} \cite{Cio79};
van der Linden\index[dN]{Linden@van der Linden} schreibt hierzu
(in Fortsetzung von oben):
\begin{quote}
  {\em \ldots and the geometrical methods apply for fields with
    $\# S_\infty < 2$}
\end{quote}
(damit sind Zahlkörper mit Einheitenrang 1 gemeint). In der Tat lässt
sich (1.8) jedoch ganz einfach auf beliebige Zahlkörper mit
unendlicher Einheitengruppe verallgemeinern (s. (1.12)); es ist
bemerkenswert, dass die dabei benutzten Ungleichungen (wie z.B. (1.9),
oder die entsprechende kubische in § 4) von
Cassels\index[dN]{Cassels} benutzt wurden, um
die Schranken für $|\disc K|$ zu verbessern, oberhalb derer es keine
quadratischen ($r=2$, $s=0$), kubischen ($r=1$, $s=1$) und
biquadratischen ($r=0$, $s=2$) Körper mit Euklidischem Algorithmus
mehr gibt. Dies ist eigentlich ein mehr als deutlicher Hinweis darauf,
dass sich auch das Casselsche Ergebnis auf Zahlkörper beliebigen
Grades verallgemeinern lässt. Ein erster Schritt in diese Richtung
wäre eine Verallgemeinerung von (2.12).

In § 2 beschreiben wir eine Methode zur Bestimmung euklidischer
Minima; diese geht im Wesentlichen auf Barnes\index[dN]{Barnes} und
Swinnerton-Dyer \index[dN]{Swinnerton-Dyer} zurück (\cite{BS52a};
wichtige Beiträge hierzu wurden aber u. a. von
Cassels\index[dN]{Cassels} und Inkeri\index[dN]{Inkeri} geliefert).
Die in § 2 gegebenen Beweise dieser Theoreme erscheinen mir wesentlich
klarer und durchsichtiger als die originalen.

In § 3 geben wir (fast) eine Klassifikation aller normeuklidischen
reellquadratischen Zahlkörper mittels (1.4) und (1.8); Vorbild sind
die klassischen Beweise von Behrbohm\index[dN]{Behrbohm},
Rédei\index[dN]{Redei@R\'edei}, Erdős\index[dN]{Erdos@Erd\"os} \&
Ko\index[dN]{Ko} usw. Die hier erzielten Ergebnisse verwenden wir auch in § 5.

In § 4 untersuchen wir kubische Körper; zur Auffindung aller
normeuklidischen Körper der Form $\Q(\sqrt[3]{m})$ folgen wir
der Arbeit Cioffaris \cite{Cio79}\index[dN]{Cioffari}. Dann bestimmen
wir die euklidischen Minima kubischer Körper mit kleiner Diskriminante,
und schließlich verbessern wir eine von Smith\index[dN]{Smith} \cite{Smi69}
erhaltene Schranke, wonach es keine kubischen zyklischen Körper mit
Diskriminante $d$, $157^2 \le d \le 10^8$ gibt, indem wir dieses Ergebnis
auf alle solchen Körper mit $157^2 \le d \le 2.5 \cdot 10^{11}$ ausweiten.

In § 5 zeigen wir, dass es nur endlich viele normeuklidische Zahlkörper
der Form $\Q(\sqrt[4]{m})$ gibt, und für $m<0$ bestimmen wir
sie alle. Damit werden die diesbezüglichen Untersuchungen von Cioffari
\cite{Cio79}\index[dN]{Cioffari} und Egami \cite{Ega84}\index[dN]{Egami}
abgeschlossen.

In § 6 bestimmen wir, wie schon erwähnt, alle imaginären, bizyklischen
biquadratischen Zahlkörper. Nachdem van der
Linden\index[dN]{Linden@van der Linden} \cite{Lin83} alle
zyklischen imaginären Zahlkörper 4. Grades gefunden hat (es gibt nur
zwei: $\Q(\zeta_5)$ und den Teilkörper 4. Grades von
$\Q(\zeta_{13})$), sind jetzt alle galoisschen imaginären
Körper 4. Grades mit EA bekannt. Es scheint, als ob das
nächstschwierigere Problem die Körper sind, deren normaler Abschluss
die Diedergruppe als Galoisgruppe besitzt.

In § 7 beschäftigen wir uns mit der Klassifikation der
normeuklidischen Dirichletschen Zahlkörper. Wir werden dabei alle
solchen finden, für welche die Relativdiskriminante durch $(1+i)$
teilbar ist. Die vollständige Lösung dieses Problems scheint nur noch
eine Frage der (Rechen-)Zeit zu sein.

In § 8 machen wir einige Andeutungen, wie man die restlichen
imaginären Zahlkörper 4. Grades behandeln kann, und in § 9 geben wir
einen Ausblick auf die Untersuchungen von Körpern höheren Grades.

Nach dem bisher Gesagten ist es nicht verwunderlich, dass auch viele
Fragen offen bleiben mussten. Einige dieser Probleme sollten sich mit
den hier vorgestellten Methoden lösen lassen (``je n'ai pas le
temps''), andere scheinen tiefer zu liegen. Jedenfalls sind in § 10
einige dieser Fragen aufgelistet.

§ 11 schließlich enthält die Beschreibung eines Algorithmus, mit dem
man den EA in einem gegebenen Zahlkörper nachweisen oder widerlegen
kann (die Frage, ob dieser Algorithmus immer terminiert, also eine
Antwort liefert, hängt eng mit gewissen -- noch nicht bewiesenen --
Vermutungen von Barnes\index[dN]{Barnes} und
Swinnerton-Dyer\index[dN]{Swinnerton-Dyer}, sowie von
Lenstra\index[dN]{Lenstra} zusammen
(s. dazu Lenstra \cite{Len79b}\index[dN]{Lenstra}).

An dieser Stelle sei noch eine Liste der Körper angegeben, von denen
bisher nicht bekannt war, ob sie normeuklidisch sind oder nicht:

\begin{itemize}
\item $n=3$, $r=3$, $s=0$: $\disc K= 2021, 2024, 2057,
  2101, 2213$;
\item $n=4$, $r=0$, $s=2$: $\Q(\sqrt{2}, \sqrt{5})$,
  $\Q(\sqrt{3}, \sqrt{17})$, $\Q(\sqrt{3},
  \sqrt{19})$, $\Q(\sqrt{7}, \sqrt{5})$,
  $\Q(\sqrt{5+2i})$, $\Q(\sqrt{1+6i})$,
  $\Q(\sqrt{7+2i})$, $\Q(\sqrt{7+4i})$,
  $\Q(\sqrt{1+8i})$, $\Q(\sqrt{3+8i})$,
  $\Q(\sqrt{5+8i})$, $\Q(\sqrt{9+4i})$,
  $\Q(\sqrt{7+8i})$, $\Q(\sqrt{11+4i})$,
  $\Q(\sqrt[4]{12})$;
\item $n=4$, $r=2$, $s=1$: $\Q(\sqrt[4]{2})$, $\Q(\sqrt[4]{5})$;
\item $n=4$, $r=4$, $s=0$: $\Q(\sqrt{3}, \sqrt{5})$,
  $\Q(\sqrt{3}, \sqrt{7})$, $\Q(\sqrt{5}, \sqrt{13})$,
  $\Q(\sqrt{5}, \sqrt{17})$.
\end{itemize}

Zweifellos lassen sich mit den vorgestellten Methoden noch viele
weitere normeuklidische Zahlkörper 3. und 4. (oder auch höheren)
Grades finden.

Weiter sind mir neben den 14 von Cooke\index[dN]{Cooke} gefundenen
reellquadratischen, 2-stufig normeuklidischen Körpern noch die
folgenden bekannt:
$\Q(\sqrt{m})$ für $m = 47, 59, 62, 67, 71, 109, 149, 157,
161, 173, 193, 201, 213$, sowie die kubischen Körper mit Diskriminante
$d = -199, -351$.

\newpage
\section*{NOTATION}
\begin{longtable}{ll}
\toprule
$\N$ & die natürlichen Zahlen \\
$\Z$ & die ganzen rationalen Zahlen \\
$\Q$ & die rationalen Zahlen \\
$\R$ & die reellen Zahlen \\
$\C$ & die komplexen Zahlen \\
$\F_q$ & der endliche Körper mit $q$ Elementen \\
$K$ & ein algebraischer Zahlkörper \\
$R$ & ein nullteilerfreier, kommutativer Ring mit 1 \\
$\mathcal{O}_K$ & der Ring ganzer Zahlen in $K$ \\
$\Phi_K$ & die Eulersche Phi-Funktion in $K$ \\
$\|\mathfrak{I}\|$ & die Absolutnorm eines Ideals $\mathfrak{I}$ \\
$\Tr_{K/\Q}(\alpha)$ &
   die Spur eines Elements $\alpha \in K$ \\
$N_{K/\Q}(\alpha)$ &
   die Norm eines Elements $\alpha \in K$ \\
$\disc K$ & die Diskriminante einer Körpererweiterung \\
$\diff K$ & die Differente einer Körpererweiterung \\
$(K:\Q)$ & der Grad einer Körpererweiterung \\
$\Gal(K/\Q)$ & die Galoisgruppe \\
$Z(\mathfrak{P}|\mathfrak{p})$ & die Zerlegungsgruppe \\
$T(\mathfrak{P}|\mathfrak{p})$ & die Trägheitsgruppe \\
$V_i(\mathfrak{P}|\mathfrak{p})$ & die Verzweigungsgruppen \\
$\Cl(K)$ & die Idealklassengruppe \\
$h(K) = \#\Cl(K)$ & die Klassenzahl \\
$|\cdot|_i$ & die normierten Bewertungen von $K$ (s. S. \pageref{pBew}
   und  \pageref{pBew2}) \\
$D(m)$ & der Ring ganzer Zahlen in $\Q(\sqrt{m})$ \\
$D(m,n)$ & der Ring ganzer Zahlen in $\Q(\sqrt{m},\sqrt{n})$ \\
QRG & quadratisches Reziprozitätsgesetz \\
FE & Fundamentaleinheit \\
GHB & Ganzheitsbasis \\
PERS & primes Einheitenrestsystem (S.~\pageref{pPERS}) \\
$E_i$ & s. S. \pageref{pEi} \\
$E'_i$ & s. S. \pageref{Ej'} \\
$F_i$ & s. S. \pageref{pFi} \\
$B_i$ & s. S. \pageref{DefB} \\
$M_1$ & erstes euklidisches Minimum (s. S. \pageref{pEuM}) \\
$M_2, \dots$ & zweites euklidisches Minimum, \dots (s. S. \pageref{pEuM2}) \\
$M^1$, $M^2$, \ldots
     & Minimum für $k$-stufige euklidische Funktionen: s. S. \pageref{DefMk} \\
$C_1, C_2, \dots$ & s. S. \pageref{DefC} \\
$M_{r,s}$ & die Minkowski-Schranke (s. S. \pageref{pEuM}) \\
$\mu_K$ & die Lenstra-Konstanten (s. S. \pageref{Lenmu}) \\
$\lambda_K$ & s. S. \pageref{Lenmu} \\
$\cM_K$ & die ``Generalmensur'', s. S. \pageref{GMen} \\
$\uK$ & s. S. \pageref{puK} \\
\bottomrule
\end{longtable}

\chapter*{\S\ 0 Einführung}
\setcounter{chapter}{0}
\addcontentsline{toc}{chapter}{\S\ 0 Einführung}
\markboth{Euklidische Ringe}{\S\ 0 Einführung}

Die Geschichte des Euklidischen Algorithmus ist, wie der Name Euklid
schon vermuten lässt, recht lang: sie beginnt mit Euklids ``Beweis'',
dass der Ring $\Z$ der ganzen rationalen Zahlen ``euklidisch''
ist (s. dazu die Arbeiten von Hendy\index[dN]{Hendy} \cite{Hen75},
Knorr\index[dN]{Knorr} \cite{Kno76}
und Collison\index[dN]{Collison} \cite{Col80}).

Ist $R$ ein Ring (darunter wollen wir hier nur einen kommutativen,
nullteilerfreien Ring mit Eins verstehen), so nennen wir eine Funktion
$f:R\to\N$ eine
\textbf{euklidische Funktion}\index[dS]{euklidische Funktion} oder einen
\textbf{euklidischen Algorithmus}\index[dS]{euklidischer Algorithmus}
(kurz EA) auf $R$, wenn gilt:

\begin{enumerate}
\item[](E--1)\ $\forall a\in R: f(a)=0 \Leftrightarrow a=0$
\item[](E--2)\ $\forall a\in R, b\in R\setminus\{0\} 
   \exists c\in R: f(a - bc) < f(b).$  
\end{enumerate}

Falls es ein solches $f$ gibt, heißt $R$
\textbf{euklidisch bezüglich $f$}.\index[dS]{euklidisch!bezüglich $f$}
Ziel der Euklidschen Überlegungen war es, in $\Z$ den
sogenannten ``Fundamentalsatz der Arithmetik'' zu zeigen: dass nämlich
jede natürliche Zahl sich bis auf die Reihenfolge eindeutig als
Produkt von Primzahlen schreiben lässt.

Zu Eulers Zeiten wurde dann erkannt, dass sich gewisse diophantische
Gleichungen wie z.B. $y^2 = x^3 - 2$ besonders einfach lösen lassen,
wenn man stattdessen $x^3 = (y + \sqrt{-2})(y - \sqrt{-2})$ schreibt und
annimmt, dass der Ring $\Z[\sqrt{-2}]$ ein ZPE-Ring ist,
d.h. dass in ihm gewissermaßen ein Analogon zum Fundamentalsatz der
Arithmetik in $\Z$ gilt. Jedoch galt diese Annahme damals als
selbstverständlich.

Erst Gauß\index[dN]{Gauß}, der für seine Theorie der biquadratischen
Reste den Ring $\Z[i]$ studierte (wo $i^2 = -1$ ist), sah die
Notwendigkeit eines Beweises ein und zeigte dazu, dass $\Z[i]$
euklidisch bezüglich der Norm ist.

Danach versuchte man einerseits, die Theorie der quadratischen und
biquadratischen Reste zu verallgemeinern auf $p$-te Potenzreste,
andererseits wollte man Fermat's letzte Vermutung beweisen, wonach die
diophantische Gleichung $x^p + y^p + z^p = 0$ für prime $p \ge 3$ nur
triviale Lösungen besitzt (das sind solche mit $xyz = 0$).

Bevor aber Kummer\index[dN]{Kummer} seine Theorie der idealen Zahlen
entwickelt hatte (die dann von Dedekind\index[dN]{Dedekind} zur
heutigen Idealtheorie ausgebaut wurde), musste man für beide Aufgaben
voraussetzen, dass die Ringe $\Z[\zeta_p]$, wo $\zeta_p$ eine
primitive $p$-te Einheitswurzel ist, ZPE-Ringe sind. Für $p = 3, 5, 7$
konnte man dies beweisen, indem man zeigte, dass die entsprechenden
Ringe normeuklidisch sind (s. Gauß \cite{Gau76}\index[dN]{Gauß} und
Kummer \cite{Kum44}).

Einen genaueren Überblick über die Geschichte des euklidischen
Algorithmus bis hin zu den neueren Ergebnissen von
Weinberger\index[dN]{Weinberger} und Lenstra\index[dN]{Lenstra} findet
man in dem sehr nett geschriebenen Artikel von
Lenstra\index[dN]{Lenstra} \cite{Len79b}; weitere historische Details
stehen bei van der Linden\index[dN]{van der Linden} \cite{Lin85}, sowie
in der bereits etwas älteren Arbeit von
Narkiewicz\index[dN]{Narkiewicz} \cite{Nar67}.

Wir wollen nun wieder zu unserer Definition einer euklidischen
Funktion zurückkehren; ist $f:R\to\N$ eine Funktion, die zwar
(E--1), aber nicht notwendig (E--2) erfüllt, so setzen wir
\[
M(f) := \inf \left\{ x \in \R : \forall a,b\in R\setminus\{0\}
\ \exists c \in R: f(a-bc) < x \cdot f(b) \right\}
\]
und $M(f)=\infty$, falls es kein solches $x$ gibt. Wir nennen $M(f)$
das \textbf{euklidische Minimum}\label{dpEuM}\index[dS]{euklidisches Minimum}
von $R$ bezüglich $f$.

Ist $R$ ein Zahlring und $f$ der Absolutbetrag der Norm, so schreiben
wir auch $M(R)$ statt $M(f)$; ist $K$ schließlich ein Zahlkörper und
$R$ der Ring aller ganzen Zahlen in $K$ (also die Hauptordnung), so
schreiben wir i.A. $M(K)$ statt $M(R)$ oder $M(f)$.

Ist $M(f)< 1$, so ist $R$ euklidisch bezüglich $f$, während dies
für $M(f)>1$ nicht der Fall ist. Was alles passieren kann, wenn
$M(f)=1$ ist, sollen die folgenden Beispiele zeigen (die wir
allerdings erst später werden rechtfertigen können):\label{dpBei}

\begin{enumerate}
\item $R$ sei ein euklidischer Ring und $f$ die ``minimale euklidische
  Funktion'' auf $R$ (die wir gleich nachher kennenlernen werden);
  dann ist $M(f)=1$ und $R$ euklidisch bezüglich $f$.
\item Sei $R$ der Ring ganzer Zahlen in dem kubischen Körper der
  Diskriminante $d=-199$, und $f$ der Absolutbetrag der Norm; dann ist
  $M(f) = M(R) = M(K) = 1$, $R$ hat Klassenzahl 1, aber $f$ ist nicht
  euklidisch auf $R$ (s. Taylor\index[dN]{Taylor} 1976).
\item Sei $R$ der Ring ganzer Zahlen in einem der Körper
  $K=\Q(\sqrt{65})$ (s. z.B. Heinhold\index[dN]{Heinhold} 1939),
  $K=\Q(\sqrt{-1},\sqrt{15})$ oder
  $K=\Q(\sqrt{-3},\sqrt{13})$ (s. dazu § 6); in diesen Fällen
  ist $M(K)=M(R)=1$, aber $K$ hat Klassenzahl 2, somit kann $f$ nicht
  euklidisch sein auf $R$.
\item Sei $R=\Z[\sqrt{-3}]$ (s. Cohn\index[dN]{Cohn} \cite{Coh78}),
  $R=\Z[\sqrt{5}]$ oder $R=\Z[\sqrt{13}]$ und $f$ der
  Absolutbetrag der Norm; dann ist in allen drei Fällen $M(R)=1$,
  obwohl $R$ hier nicht einmal ganz abgeschlossen, geschweige denn
  ZPE-Ring oder gar euklidisch ist. Wir werden später in Zahlkörpern
  höheren Grades noch weitere Beispiele für ein solches Verhalten
  finden.
\end{enumerate}

Diese Beispiele zeigen, dass man i.A. aus $M(f)=1$ nicht folgern kann,
dass $R$ euklidisch (oder nicht euklidisch) bezüglich $f$ ist.

Eine manchmal sogar exakte untere Schranke für das euklidische Minimum
$M(f)$ gibt das folgende Kriterium, das trotz seiner Einfachheit
bisher nirgendwo explizit zu finden ist (nicht einmal in dem Falle, wo
$f$ der Absolutbetrag der Norm ist):

\begin{quote}
  {\bf (0.1)} {\em Sei $b \in R \setminus \{0\}$ keine Einheit;
    dann ist $M(f) \geq \frac1{f(b)}$.}
\end{quote}

\begin{proof}
  Da $b$ keine Einheit ist, gibt es ein $a \in R$ mit $b \nmid a$;
  schreibt man dann $a = bq + r$ für irgendwelche $q, r \in R$, so muss $r
  \neq 0$ sein. Wegen (E--1) ist also $f(r) \geq 1$ für jede Wahl von
  $q$ und $r$, und damit $\frac{f(r)}{f(b)} \geq \frac{1}{f(b)}$.
  Die Behauptung folgt nun nach Definition von $M(f)$.
\end{proof}

Ist $K$ ein quadratischer Zahlkörper mit Diskriminante $d$, so gibt
die untenstehende Tafel ein $b \in R$ mit minimaler Norm $>1$ und die
daraus resultierende Schranke für $M(K)$:

\begin{center}
$$ \begin{array}{ccc}
\toprule
d & b & M(K) \geq \\
\midrule
\rsp -4 & 1+i & \frac{1}{2} \\
\rsp -3 & \sqrt{-3} & \frac{1}{3} \\
\rsp  5 & 2 & \frac{1}{4} \\
\rsp  8 & \sqrt{2} & \frac{1}{2} \\
\rsp 12 & 1+\sqrt{3} & \frac{1}{2} \\
\rsp 13 & 4+\sqrt{13} & \frac{1}{3} \\
\rsp 17 & \frac{5+\sqrt{17}}{2} & \frac{1}{2} \\
\bottomrule
\end{array} $$
\end{center}

Wir werden in § 2 sehen, dass diese Schranken sogar exakt sind; weiter
sind dies die einzigen quadratischen Körper, für welche die Schranke
aus (0.1) exakt ist. Auffallend hierbei ist, dass man diese Körper auch
noch durch $|d|\leq 4$ im komplexen, bzw. durch $d\leq 17$ im reellen
Fall charakterisieren kann, d.h.: die Körper, für die die Schranke aus
(0.1) exakt ist, sind gerade diejenigen mit minimaler Diskriminante!
Ähnliche Beobachtungen werden wir später auch für Zahlkörper dritten
und vierten Grades machen können.

Wir wollen uns nun etwas mit der allgemeinen Theorie euklidischer
Ringe befassen. Die Inklusionen
\[
\{\text{Euklidische Ringe}\} \subset \{\text{Hauptidealringe}\}
\subset \{\text{ZPE-Ringe}\}
\]
sind klassisch. Für die hier betrachteten Zahlringe lässt sich die
letzte Inklusion sogar umkehren, d.h. hier ist jeder ZPE-Ring auch
Hauptidealring. Das Beispiel $R = \Z[\frac{1+\sqrt{-19}}2]$ zeigt
allerdings, dass die erste Inklusion echt ist: $R$ ist nämlich
Hauptidealring, aber nicht euklidisch, wie wir gleich sehen werden.

Dazu folgen wir einer Idee Motzkins\index[dN]{Motzkin} \cite{Mot49}
(s. auch Samuel\index[dN]{Samuel} \cite{Sam71}) und versuchen, zu einem
gegebenen Ring eine Funktion $f$ so zu konstruieren, dass $f$
euklidisch auf $R$ wird. Wegen (E--1) kennen wir alle $a \in R$ mit
$f(a)=0$ (nämlich $a=0$) und fragen uns, für welche $b \in R$ wir
$f(b)=1$ setzen dürfen. Wegen (E--2) ist dies nur dann möglich, wenn
für alle $a \in R$ ein $q \in R$ existiert mit $f(a-bq) < f(b) = 1$;
dies geht nur für $f(a-bq)=0$, und wegen (E--1) muss $a-bq=0$
sein. Dies zeigt, dass ein $b \in R$ mit $f(b)=1$ jedes $a \in R$
teilen muss: $b$ muss also Einheit in $R$ sein. Ist andererseits $b$
eine Einheit, so lässt sich (E--2) für jedes gegebene $a \in R$ mit
$q=ab$ sicher erfüllen, und wir dürfen in der Tat $f(b)=1$ setzen.

Ist nun $f(b)=2$, so soll es zu jedem $a \in R$ ein $q \in R$ geben
mit $f(a-bq)<2$; also ist $a-bq=0$ (im Falle $f(a-bq)=0$) oder gleich
einer Einheit (falls $f(a-bq)=1$ ist), anders ausgedrückt: jedes
$a \in R$ ist modulo $b$ entweder $0$ oder einer Einheit kongruent. Indem
wir so fortfahren, werden wir auf die folgenden Teilmengen\label{dpEi}
$E_i$ von $R$ geführt:

\begin{align*}
  E_0 & = \{0\}, \quad E_1 = E_0 \cup R^\times , \quad
          \text{und allgemein für alle $i \geq 1$} \\
  E_i \setminus E_{i-1} & = \{ b \in R : \text{jede Restklasse} \bmod b
  \text{ enthält ein Element aus } E_{i-1} \}.
  \end{align*}

Schließlich setzen wir noch
$$ E_\infty = \bigcup_{i \geq 0} E_i. $$

Damit gilt

\begin{quote}
  {\bf (0.2)} {\em Ein Ring $R$ ist genau dann euklidisch, wenn
    $R = E_\infty$ gilt.}
\end{quote}

\begin{proof}
  Sei $R = E_\infty$. Wir definieren durch $f_0(a) = \min \{ i \in
  \N : a \in E_i \}$ eine Funktion $f_0 : R \to
  \N$. Sind dann $a, b \in R \setminus \{0\}$ gegeben und ist
  z.B. $f_0(b) = i$, so ist $i > 1$ wegen $b \neq 0$, und nach
  Konstruktion von $f_0$ ist $b \in E_i$. Also gibt es ein $r \in
  E_{i-1}$ mit $a \equiv r \bmod b$; setzt man $a = bq + r$, so ist
  $q \in R$ und $f_0(r) \leq i-1 < f_0(b)$. Also ist $R$ euklidisch
  bezüglich $f_0$.

  Sei andererseits $R$ euklidisch bezüglich einer Funktion $f$. Wir
  behaupten, dass dann alle $a \in R$ mit $f(a) = i$ bereits in $E_i$
  enthalten sind. Für $i = 0$ ist dies wegen (E--1) sicher
  richtig. Haben wir die Behauptung für alle $i < k$ gezeigt und ist
  $f(b) = k$, dann existieren zu jedem $a \in R$ Elemente $q, r \in R$
  mit $a = bq + r$ und $f(r) < f(b)$. Nach Induktionsvoraussetzung ist
  $r \in E_{k-1}$, d.h. zu jedem $a \in R$ gibt es ein $r \in E_{k-1}$
  mit $a \equiv r \bmod b$; dies impliziert aber $b \in E_k =
  E_{k}(b)$, und folglich ist $R \subset E_\infty \subset R$.
\end{proof}

Wenn wir uns den zweiten Teil des Beweises anschauen, so bemerken wir,
dass wir viel mehr gezeigt haben: ist nämlich $R$ euklidisch bezüglich
$f$ und $f_0$ die im ersten Teil des Beweises definierte Funktion,
dann zeigt die Relation $b \in E_{f(b)}$, dass $f_0(b) \leq f(b)$ für
alle $b \in R$ ist, denn $f_0(b)$ gibt ja den kleinsten Index $i$ an,
für den $b \in E_i$ ist. Wir werden daher auf folgende Definition
geführt: ist $R$ ein euklidischer Ring, dann heißt die Funktion $f_m :
R \to \N$, die durch
\[
f_m(a) = \min \{ f(a) : f \text{ ist euklidische Funktion auf } R \}
\]
definiert wird, die
\textbf{minimale euklidische Funktion}\index[dS]{euklidische Funktion!minimale}
auf $R$. Diese Bezeichnung wird gerechtfertigt durch

\begin{quote}
  {\bf(0.3)} {\em $f_m$ ist eine euklidische Funktion auf $R$ mit der
    Eigenschaft $f_m(a) \leq f(a)$ für alle $a \in R$ und alle
    euklidischen Funktionen $f$ auf $R$. Insbesondere gilt
    $f_m = f_0$.}
\end{quote}

\begin{proof}
  Seien $a, b \in R \setminus \{0\}$ gegeben. Nach Definition von
  $f_m$ gibt es eine euklidische Funktion $f$ mit $f_m(b) =
  f(b)$. Damit existieren $q, r \in R$ mit $a = bq + r$ und $f(r) <
  f(b)$. Daraus ergibt sich $f_m(r) \leq f(r) < f(b) = f_m(b)$, und
  $f_m$ ist in der Tat eine euklidische Funktion auf $R$.

  Wegen $f_m(a) \leq f(a)$ für alle euklidischen Funktionen $f$ auf
  $R$ und $f_0(a) \leq f_m(a)$ wegen der Bemerkung nach (0.2) muss nun
  $f_m(a) = f_0(a)$ für alle $a \in R$ gelten.
\end{proof}

Im Falle $R = \Z$ lassen sich die Mengen $E_i$ ganz einfach
beschreiben: man hat nämlich $E_0 = \{0\}$, $E_1 = \{0, -1, 1\}$, und
$E_2 \setminus E_1$ enthält diejenigen $a \in \Z$, die in
jeder Restklasse modulo $a$ einen Vertreter aus $E_1$ haben; da $E_1$ nur
drei Elemente enthält, muss $|a| \le 3$ gelten, und in der Tat findet man
$E_2 = \{0, \pm1, \pm2, \pm3\}$. Mittels vollständiger Induktion sieht
man nun leicht ein, dass $E_i = \{ a \in \Z : |a| \leq 2^i - 1
\}$ gilt. Insbesondere ergibt sich hieraus $\Z = E_\infty$,
d.h. $\Z$ ist euklidisch.

Etwas komplizierter liegen die Dinge in imaginärquadratischen
Zahlkörpern. Ist $m \neq -1, -3$, so enthält $D(m)$ (wir werden
künftig den Ring ganzer Zahlen in $\Q(\sqrt{m})$ mit $D(m)$
bezeichnen) für $m<0$ bekanntlich nur die Einheiten $+1$ und $-1$; in
diesen Fällen ist also $E_1 = \{0, +1, -1\}$. $E_2$ besteht nun aus
allen $a \in D(m)$, für die jede Restklasse modulo $a$ eines der drei
Elemente $0, 1$ oder $-1$ enthält. Wegen $|R/aR| =
N_{K/\Q}(a)$ ist dies offenbar nur dann möglich, wenn
$N_{K/\Q}(a) \le 3$ ist. Die einzigen imaginärquadratischen
Ringe, die Elemente der Norm 2 oder 3 enthalten, sind aber $D(m)$ für
$m=-1, -2, -3, -7, -11$; für alle anderen $D(m)$ ist also $E_1 = E_2 =
\dots = E_\infty = \{0, 1, -1\}$, und wir haben

\begin{quote}
  {\bf (0.4)} {\em Sei $m<0$ und $R=D(m)$; dann ist $R$ genau für die
      Werte $m=-1, -2, -3, -7, -11$ euklidisch, und in diesen Fällen
      ist die Norm eine euklidische Funktion. Für alle anderen Werte
      von $m$ ist $E_1 = E_2 = \dots = E_\infty = \{0, 1, -1\}$.}
\end{quote}

Insbesondere sind also die Hauptidealringe $D(m)$ mit $m=-19, -43,
-67, -163$ nicht normeuklidisch. Dass $D(m)$ für die oben angegebenen
Werte von $m$ normeuklidisch ist, stammt von Gauß\index[dN]{Gauß}
($m=-1, -3$), bzw. von Dedekind\index[dN]{Dedekind} \cite{Ded00} und
Dickson\index[dN]{Dickson} \cite{Dic27}. Die beiden letzteren haben
auch bemerkt, dass es keine anderen normeuklidischen $D(m)$ mit $m<0$
gibt. Die weitergehende Aussage (0.4) dagegen stammt von
Motzkin\index[dN]{Motzkin} \cite{Mot49} und unabhängig davon von
Dubois\index[dN]{Dubois} und Steger\index[dN]{Steger}
\cite{DS58}. Narkiewicz\index[dN]{Narkiewicz} schreibt
\cite[S. 175]{Nar67}, dass Dubois\index[dN]{Dubois} und
Steger\index[dN]{Steger} darüberhinaus gezeigt hätten, dass jede
euklidische Funktion $f$ auf $D(m)$, $m<0$, mit der Norm
übereinstimmen muss, was aber weder historisch (die beiden haben das
nicht behauptet) noch mathematisch wahr ist (denn für $m=-1, -2, -3,
-7, -11$ ist ja auch die minimale euklidische Funktion ein
euklidischer Algorithmus auf $D(m)$, und dieser unterscheidet sich von
der Norm). Weitere Beweise von (0.4) haben Stewart\index[dN]{Stewart}
und Tall\index[dN]{Tall} \cite{ST79} in ihrem Lehrbuch {\em Algebraic
  number theory} und Castro Chadid\index[dN]{Castro Chadid}
\cite{Cha84} gegeben. Schließlich hat Campoli\index[dN]{Campoli}
\cite{Cam88} noch einmal gezeigt, dass $D(-19)$ zwar Hauptidealring,
aber nicht euklidisch ist. Allen Beweisen gemeinsam ist die
Konstruktion einer Nichteinheit der Norm $\le 3$ in einem euklidischen
$D(m)$.

Für die Ringe $D(-1)$ und $D(-3)$ hat Lenstra\index[dN]{Lenstra}
\cite{Len74} die Mengen $E_i$ explizit bestimmt; ansonsten kennt man
nur einige elementare Eigenschaften der $E_i$:

\begin{quote}
  {\bf (0.5)} {\em Sei $R$ ein Zahlring; dann gilt
    \begin{enumerate}
    \item[(i)] Ist $a \in E_\infty$, dann ist jeder Teiler des Ideals $aR$
      ein Hauptideal.
    \item[(ii)] Ist $a = p_1 p_2 \dots p_m$ die Primzerlegung von $a \in
      R$, dann liegt $a$ nicht in $E_m$.
    \item[(iii)] Für alle $i \in \N$ ist $E_i \neq R$.
    \item[(iv)] Es gilt $E_i R^\times = E_i$ d.h. mit $a \in E_i$ und $u
      \in R^\times$ ist auch $au \in E_i$.
    \item[(v)] Ist $K/\Q$ galoisch, dann ist $E_i^G = E_i$ für alle
      $G \in \operatorname{Gal}(K/\Q)$.
    \end{enumerate} }
\end{quote}

\begin{proof}
  (i) Sei $a \in E_1$; ist dann $P$ ein Primideal in $R$, welches
  $(a)$ teilt, dann gibt es ein $p_0 \in R$ mit $P = (a, p_0)$ (denn
  jedes Ideal in einem Zahlkörper wird von maximal zwei Elementen
  erzeugt, wobei man $a \in P$ beliebig vorgeben darf;
  s. z.B. Cohn\index[dN]{Cohn} 1978). Wegen $a \in E_1$ existiert nun
  ein $p_1 \in E_{k-1}$ mit $p_1 \equiv p_0 \mod a$. Damit ist $P =
  (a, p_0) = (a, p_1)$. Ersetzt man in obigem Gedankengang $a \in E_1$
  durch $p_1 \in E_{k-1}$, so folgt $P = (p_1, p_2)$ für ein $p_2 \in
  E_{k-2}$ usw. Schließlich erhalten wir $P = (p_1, p_2, \dots,
  p_{k-1})$ für ein $p_{k-1} \in E_1$ (möglicherweise tritt dies auch
  schon früher ein). Da $P$ prim ist, kann $p_{k-1}$ keine Einheit
  sein. Folglich ist $p_{k-1} \neq 0$ und $P = (p_{k-1})$ ein
  Hauptideal.

  (ii) Die Behauptung ist richtig für $m=0, 1$; sei sie nun bewiesen
  für alle $i \in \N$ mit $i<k$ und sei $a \in E_k$. Wir nehmen an, es
  sei die Anzahl $m$ der Faktoren von $a$ größer gleich $k$ und wählen
  einen Primfaktor $p$ von $a$. Dann ist $a/p \equiv e \mod a$ für ein
  $e \in E_{k-1}$, und es ist $e \neq 0$ (sonst wäre $ap \mid a$). Nun
  gilt aber $e \equiv 0 \mod (a/p)$, d.h. $e \in E_{k-1}$ hat $k-1$
  Primfaktoren, und dies widerspricht der Induktionsvoraussetzung.

  (iii) folgt sofort aus ii), wenn man beachtet, dass es in $R$
  unendlich viele Primideale gibt (nämlich mindestens eines über jeder
  rationalen Primzahl).

  (iv) folgt ebenfalls sofort, und zwar wegen $R/(a) \cong R/(au)$.

  (v) Die Behauptung ist richtig für $j=0$. Gilt sie für $j+k-1$ und
  ist $a \in E_k$ dann gibt es zu jedem $b \in R$ ein $c \in E_{k-1}$
  mit $b \equiv c \mod a$. Also ist $b^G \equiv c^G \mod a^G$, und
  nach Induktionsvoraussetzung ist $c^G \in E_{k-1}^G$. Da mit $b$
  auch $b^G$ ganz $R$ durchläuft, folgt $a^G \in E_k^G$ und damit die
  Behauptung.
\end{proof}

Man wird sich nun fragen, ob es Zahlringe gibt, die zwar bezüglich
$f_0$, nicht aber bezüglich der Norm euklidisch sind. Bis heute ist
kein einziger solcher Zahlring bekannt. Es wäre allerdings etwas
voreilig, daraus auf irgendeine mathematische Wahrheit schließen zu
wollen: falls nämlich gewisse Riemannsche Vermutungen richtig sind,
ist jeder Zahlring mit Klassenzahl 1 und unendlicher Einheitengruppe
euklidisch bezüglich $f_0$; dieses Resultat stammt von
Weinberger\index[dN]{Weinberger}
\cite{Wei73} und beruht auf einer Verallgemeinerung von Hooleys
``Beweis'' der Artinschen Vermutung über Primitivwurzeln (Hooley hatte
diese Vermutung 1967 unter Annahme der Riemannschen Vermutung
bewiesen). Für einen Beweis dieser Aussage verweisen wir auf
Weinberger\index[dN]{Weinberger} \cite{Wei73} und Lenstra\index[dN]{Lenstra}
\cite{Len74}, für Näheres über
die Artinsche Vermutung auf Gupta\index[dN]{Gupta},
Murty\index[dN]{Murty} \& Murty \cite{GMM87}, Murty \cite{Mur88} und
Narkiewicz\index[dN]{Narkiewicz} \cite{Nar88}. Hier wollen wir nur kurz
andeuten, was die minimale euklidische Funktion mit Primitivwurzeln
zu tun hat.

Dazu sei $a \in R$ prim; ist dann die Einheit $u \in R^\times$ eine
Primitivwurzel modulo $a$, dann gilt $a \in E_2$, weil die Potenzen von
$u$ zusammen mit $0$ ein vollständiges Restsystem modulo $a$ bilden und
als Einheiten in $E_1$ liegen. Ist $u$ Primitivwurzel für ``viele'' $a
\in R$, so wird $E_2$ ``groß'' und man hat gute Chancen, $R =
E_\infty$ nachzuweisen. In der Tat hat Weinberger\index[dN]{Weinberger}
gezeigt, dass man
bereits an $E_3$ ablesen kann, ob $R = E_\infty$ ist oder nicht.

Wir kommen nun zur Definition von $k$-stufig euklidischen Ringen
(s. dazu die Arbeiten von Cooke\index[dN]{Cooke} \cite{Coo77}, sowie
Cooke\index[dN]{Cooke} und Weinberger\index[dN]{Weinberger} \cite{CW75}.
Seien nun $a,b \in R \setminus \{0\}$; eine Folge von Gleichungen
\begin{align*}
a &= q_1 b + r_1 \\
b &= r_1 q_2 + r_2 \\
r_1 &= r_2 q_3 + r_3 \\
&\vdots \\
r_{k-2} &= r_{k-1} q_k + r_k,
\end{align*}
wobei oBdA von zwei aufeinanderfolgenden $q_i$ mindestens eines und
außerdem $q_k$ von $0$ verschieden sind, heißt eine von $(a,b)$
ausgehende \textbf{Teilerkette}\index[dS]{Teilerkette} der Länge $k$;
ist $r_k = 0$, so nennt man sie
\textbf{abbrechend}.\index[dS]{Teilerkette!abbrechend} Eine
Abbildung $f: R \to \N$ heißt eine \textbf{$k$-stufige euklidische
  Funktion}\index[dS]{euklidische Funktion!$k$-stufig} auf $R$,
wenn (E--1) erfüllt ist und es für alle $a,b \in R \setminus \{0\}$
eine von $(a,b)$ ausgehende Teilerkette der Länge $k$ gibt mit $f(r_k)
< f(b)$. Gibt es ein solches $f$, dann heißt $R$ \textbf{$k$-stufig
  euklidisch}\index[dS]{euklidisch!$k$-stufig} bezüglich $f$. Schließlich
nennen wir $R$ \textbf{quasi-euklidisch},\index[dS]{quasi-euklidisch}
wenn es für alle $a,b \in R
\setminus \{0\}$ eine von $(a,b)$ ausgehende, abbrechende Teilerkette
gibt (man beachte, dass die Eigenschaft quasi-euklidisch nicht von
irgendwelchen Funktionen abhängt). Statt quasi-euklidisch hat
Cooke\index[dN]{Cooke}
den Begriff $\omega$-stufig euklidisch gewählt. Wir zeigen zuerst,
dass die in der Literatur auftauchenden Definitionen
quasi-euklidischer Ringe (z.B. in Cooke\index[dN]{Cooke} \cite{Coo77},
Bougaut\index[dN]{Bou76} \cite{Bou76,Bou80} 
oder Leutbecher\index[dN]{Leutbecher} \cite{Leu77} alle gleichbedeutend sind:

\begin{quote}
  {\bf (0.6)} {\em Sei $R$ ein Ring; dann sind äquivalent:
    \begin{enumerate}
    \item[(i)] $R$ ist quasi-euklidisch;
    \item[(ii)] für alle $a,b \in R \setminus \{0\}$ existiert ein $k \in
      \N$ und eine von $(a,b)$ ausgehende, abbrechende
      Teilerkette;
    \item[(iii)] es gibt eine Funktion $\Phi: R \times R \to \N$,
      sodass für alle $a,b \in R \setminus \{0\}$ Zahlen $q,r \in R$
      existieren mit $a=bq+r$ und $\Phi(b,r) < \Phi(a,b)$.
    \end{enumerate}}
\end{quote}

\medskip\noindent
\textbf{Bem.:} Funktionen, wie sie in (iii) beschrieben sind, nennen
wir \textbf{quasi-euklidische
  Funktionen}\index[dS]{quasi-euklidische Funktion} auf $R$.

\begin{proof}
  (i) $\Rightarrow$ (iii): Wir definieren
  $\Phi: R \times R \to \N$ durch
  \begin{align*}
    \Phi(a,b) & = \min \{k \in \N: \text{es gibt eine von }
                   (a,b) \text{ ausgehende, } k\text{-stufige} \\
              & \quad \text{abbrechende Teilerkette} .\}
  \end{align*}
  Man rechnet leicht nach, dass $\Phi$ eine quasi-euklidische Funktion ist.

  (iii) $\Rightarrow$ (i) Seien $a,b \in R \setminus \{0\}$
  gegeben. Dann existieren $q,r \in R$ mit $a=bq+r$ und
  $\Phi(b,r) < \Phi(a,b)$. Ist $r = 0$, so sind wir fertig; andernfalls
  gibt es zu $b,r$ Zahlen $q_2,r_2 \in R$ mit $b = q_2 r + r_2$ und
  $\Phi(r, r_2) < \Phi(b,r)$. Indem wir diesen Schritt genügend oft
  wiederholen, erhalten wir eine von $(a,b)$ ausgehende, abbrechende
  Teilerkette.

  (i) $\Rightarrow$ (ii) ist klar.

  (ii) $\Rightarrow$ (i): Seien $a,b \in R \setminus \{0\}$ gegeben.
  Dann gibt es ein $k \in \N$ und eine von $(a,b)$ ausgehende
  $k$-stufige Teilerkette mit $f(r_k) < f(b)$. Ist bereits $r_k = 0$,
  so sind wir fertig. Andernfalls wiederholen wir diesen Schritt mit
  dem Paar $(r_{k-1},r_k)$ und erhalten eine von $(r_{k-1},r_k)$
  ausgehende Teilerkette der Länge $l$ mit $f(r_l) < f(r_k)$
  usw. Indem wir diese Teilerketten aneinanderhängen, bekommen wir
  eine von $(a,b)$ ausgehende, abbrechende Teilerkette.
\end{proof}

Die Implikation ii) $\Rightarrow$ i) zeigt, dass jeder $k$-stufig
euklidische Ring auch quasi-euklidisch ist.

Betrachtet man eine von $(a,b)$ ausgehende Teilerkette mit $r_{k-1}
\neq 0$ und $r_k = 0$, so stellt man fest, dass $(a,b) = (r_{k-1})$
ist. Mittels vollständiger Induktion kann man dann folgern, dass in
quasi-euklidischen Ringen jedes endlich erzeugte Ideal ein Hauptideal
ist. Jedoch sind quasi-euklidische Ringe nicht notwendig
Hauptidealringe, wie das folgende Beispiel belegt:

Sei $A$ der Ring aller ganzen algebraischen Zahlen; jedes endlich
erzeugte Ideal in $A$ ist zwar Hauptideal, dennoch ist $A$ bekanntlich
kein Hauptidealring. Insbesondere kann $A$ nicht euklidisch sein;
dagegen ist $A$ 2-stufig euklidisch bezüglich jeder Funktion, die
(E--1) erfüllt: sind nämlich $a,b \in A$ teilerfremd (d.h. ist
$(a,b)=A$), so gibt es nach Lenstra\index[dN]{Lenstra} \cite{Len73} ein
$q \in A$, sodass $a-bq$
eine Einheit in $A$ ist. Daher haben solche $a,b$ eine abbrechende
Teilerkette der Länge $\le 2$. Sind $a$ und $b$ aber nicht
teilerfremd, so wähle man ein $d \in A$ mit $(d)=(a,b)$, setze dann
$a'=a/d$, $b'=b/d$ und multipliziere die von $(a',b')$ ausgehende,
abbrechende Teilerkette der Länge 2 mit $d$.

Zur Untersuchung $k$-stufiger euklidischer Funktionen führen wir ein
Analogon zu $M(f)$ ein, indem wir definieren:\label{dDefMk}
\begin{align*}
M^k(f) := \inf \Bigl\{ \kappa \in \R \;\Big|\;
    &\text{für alle } a,b \in R \setminus \{0\}
      \text{ existiert eine von } (a,b) \\
    &\text{ ausgehende Teilerkette mit }
      f(r_k) < \kappa \cdot f(b) \Bigr\}.
\end{align*}

Da die Bezeichnungen $M_2(f)$, $M_3(f)$ usw. für das zweite, dritte
usw. Minimum von $f$ schon fast klassisch sind, haben wir $M^k(f)$
oben indiziert. Wir haben damit die Ungleichungskette
 $$ M(f) = M^1(f) \ge M^2(f) \ge M^3(f) \ge \dots \ge M^\infty(f) =
           \lim_{k \to \infty} M^k(f) $$
(dieser Grenzwert existiert, da die Folge $M^k(f)$ eine monoton
fallende, durch 0 nach unten beschränkte Folge reeller Zahlen ist).

Zur Beschreibung $k$-stufiger Teilerketten in einem Ring $R$ haben
sich Kettenbrüche mit Koeffizienten aus $R$ bewährt: ist eine Folge
$q_1, \dots, q_k \in R$ gegeben (wobei von zwei aufeinanderfolgenden
$q$ mindestens eines und außerdem noch $q_k$ von $0$ verschieden ist),
so definiert man
\[
[q_1] = q_1 = \frac{a_1}{b_1} \quad \text{mit } a_1 = q_1, \quad b_1 =
1
\]
\[
[q_1, q_2] = q_1 + \frac{1}{q_2} = \frac{a_2}{b_2} \quad \text{mit }
a_2 = q_1 q_2 + 1, \quad b_2 = q_2.
\]
\[
[q_1, \dots, q_k] = \frac{a_k}{b_k} \quad \text{mit } a_k = q_k
a_{k-1} + a_{k-2}, \quad b_k = q_k b_{k-1} + b_{k-2}
\]
Man sieht dann leicht ein, dass $[q_1, \dots, q_k] = [q_1, \dots,
  q_{k-1} + 1/q_k]$ gilt, und daraus folgt durch vollständige
Induktion
\[
[q_1, \dots, q_k] = q_1 + \cfrac{1}{q_2 + \cfrac{1}{q_3 + \dots +
    \cfrac{1}{q_{k-1} + \cfrac{1}{q_k}}}}
\]
Diese letzte Identität wiederum impliziert sofort $[q_1, \dots, q_k] =
q_1 + 1/[q_2, \dots, q_k]$. Nennt man $a_k$ den Zähler und $b_k$ den
Nenner des Kettenbruchs, so hat man also
\begin{quote}
  {\bf (0.7)} {\em Der Nenner von $[q_1, \dots, q_k]$ ist gleich dem
    Zähler von $[q_2, \dots, q_k]$.}
\end{quote}

Mittels vollständiger Induktion folgt jetzt
\begin{quote}
  {\bf (0.8)} {\em Sind $a,b \in R \setminus \{0\}$ und ist $r_k$ der
    letzte Rest der von $(a,b)$ ausgehenden Teilerkette mit den
    Quotienten $q_1, \dots, q_k$, dann gilt}
    $$ \frac{a}{b} - \frac{a_k}{b_k} = (-1)^k \frac{r_k}{b \cdot b_k}. $$
\end{quote}

Ist $f$ eine multiplikative Funktion, so kann man $f$ multiplikativ
von $R$ nach $K$ fortsetzen, und wegen (0.8) ist genau dann $f(r_k) <
f(b)$, wenn $f(a/b - a_k/b_k) < 1$ gilt, damit erhält man
\begin{quote}
  {\bf (0.9)} {\em Ein Ring $R$ ist genau dann $k$-stufig
    normeuklidisch bezüglich einer multiplikativen Funktion $f$, wenn
    es zu jedem $a/b \in K$ einen Kettenbruch $[q_1, \dots, q_k] =
    a/b$ der Länge $l \le k$ gibt, sodass die Ungleichung
    $$ f\left( \frac{a}{b} - \frac{a_k}{b_k} \right) <
       f\left( \frac{1}{b_k} \right) $$
    erfüllt ist.}
\end{quote}

Der Unterschied zum gewöhnlichen EA lässt sich also folgendermaßen in
Worte fassen: in euklidischen Ringen lässt sich ein $x \in K$ so durch
ein $y \in R$ approximieren, dass $f(x-y)<1$ wird; in einem $k$-stufig
euklidischen Ring dagegen darf man $y$ aus der $R$ umfassenden Menge
der Kettenbrüche der Länge $\le k$ wählen, allerdings muss die
Approximation ``genauer'' sein.

Nun ist es i.A. recht schwer, einen Kettenbruch (auch nur der Länge 2)
als solchen zu erkennen; so ist z.B.
\[
\frac{\sqrt{14}}{2} = -2 + \frac{1}{4 - \sqrt{14}} = [-2, 4 - \sqrt{14}]
\]
ein Kettenbruch der Länge 2 mit Nenner $4 - \sqrt{14}$, während
$\frac{1+\sqrt{14}}{2}$ kein Kettenbruch der Länge 2 ist. Um dies zu
beweisen, verwenden wir

\begin{quote}
  {\bf (0.10)} {\em Ist $a_k/b_k$ ein Kettenbruch und gilt
    $a/b = a_k/b_k$, dann ist $b_k \mid b$.}
\end{quote}

\begin{proof}
  Mittels vollständiger Induktion zeigt man
  $a_k b_{k-1} - a_{k-1} b_k = (-1)^k$.
  Multipliziert man diese Identität mit $b$, so folgt
  $b a_k b_{k-1} - b a_{k-1} b_k = (-1)^k b$.
  Wegen $a b_k = a_k b$ teilt $b_k$ die
  linke und damit auch die rechte Seite der letzten Gleichung.

  Wäre nun $(1 + \sqrt{14})/2 = [q_1, q_2] = (1 + q_1 q_2)/q_2$ ein
  Kettenbruch der Länge 2, so müsste $q_2 \mid 2$ gelten; da sich aber
  $(1 + \sqrt{14})/2$ nicht mehr kürzen lässt, gilt auch $2 \mid
  q_2$. Folglich ist $q_2 = 2e$ für eine Einheit $e \in R^\times$. Nun
  ist aber $u = 15 + 4\sqrt{14}$ die FE von $R$, sodass jede Einheit
  $\equiv 1 \bmod 2$ ist. Insbesondere ist $e \equiv 1 \bmod 2$, und
  es müsste $(1 + \sqrt{14})e = 1 + \sqrt{14} \equiv 1 + q_1 q_2 \bmod
  2$ gelten, was offenbar falsch ist.
\end{proof}

Um nun ein Kriterium herzuleiten, das es uns erlaubt, gewisse Zahlen
als Kettenbrüche zu schreiben, führen wir die Mengen\label{dEj'} $E_j'$
ein, die wie folgt definiert sind: $E_0' = \{0\}$, und
\[
E_j' \setminus E_{j-1}' = \{ a \in R : \text{jede prime Restklasse }
   \bmod a \text{ hat einen Vertreter aus } E_{j-1}' \}.
\]
Man stellt dann fest, dass $E_0' = E_0$, $E_1' = E_1$, und
$E_j' \subseteq E_j$ für $j \ge 2$ gilt, wobei die letzte Inklusion
i.A. echt ist: für $D(-5)$ ist nämlich
$E_1 = E_2 = \dots = E_\infty = \{0, -1, +1\}$, während
$E_1' = E_2' = \dots = E_\infty' = \{0, \pm 1, \pm 1 \pm \sqrt{-5}\}$
gilt. Das bereits angekündigte Kriterium lautet nun

\begin{quote}
  {\bf (0.11)} {\em Seien $a, b \in R \setminus \{0\}$, $b \in E_k'$,
    und $(a, b) = 1$; dann ist $a/b$ ein Kettenbruch der Länge $\le k$
    mit Nenner $ub$, wo $u \in R^\times$ eine Einheit in $R$ ist. }
\end{quote}

\begin{proof}
  Ist $k=1$, so ist $b$ eine Einheit und $a/b = q_1 \in R$ ein
  Kettenbruch der Länge 1 mit Nenner $1 = bu$, wo $u = \frac1b$
  Einheit in $R$ ist.
  
  Sei nun die Behauptung bewiesen für alle $j < k$ und $b \in E_k'$.
  Wegen $(a, b) = 1$ liegt $a$ in einer primen Restklasse mod
  $b$, folglich gibt es ein $e \in E_{k-1}'$ mit $a \equiv e \pmod{b}$,
  d.h. es existiert ein $q_1 \in R$ mit $a = q_1 b + e$. Jetzt sehen
  wir, dass $(b, e) = (b, a) = 1$ und $a/b = q_1 + e/b = q_1 + 1/(b/e)$
  ist.

  Nach Induktionsvoraussetzung ist nun $b/e$ ein Kettenbruch der Länge
  $\le k-1$, dessen Nenner (und damit auch dessen Zähler) mit $e$
  (bzw. mit $b$) bis auf einen Faktor $u \in R^\times$
  übereinstimmt. Schreibt man $b/e = [q_2, \dots, q_k]$, so ist $a/b =
  [q_1, \dots, q_k]$, und (0.7) zeigt, dass $a/b$ ein Kettenbruch der
  Länge $\le k$ mit Nenner $bu$ ist.
\end{proof}

Als Beispiel zeigen wir, wie man (0.11) anwenden kann, um
$(1+\sqrt{14})/2$ als Kettenbruch zu schreiben. Zuerst beachten wir
$(1+\sqrt{14})/2 = -1 + (3+\sqrt{14})/2$ (damit haben wir die Norm des
Nenners klein gemacht) und behaupten dann, dass $a = 3+\sqrt{14} \in E'_2$
ist.\label{dp11} Dazu müssen wir zeigen, dass jede prime
Restklasse mod $a$ eine Einheit enthält. Dies sehen wir so ein: es
gilt $u = 15+4\sqrt{14} \equiv 3 \pmod{a}$ (man beachte einfach
$\sqrt{14} \equiv -3 \pmod{a}$), also ist $u^2 \equiv 4$, $u^3 \equiv
2$, $u^4 \equiv 1 \pmod{a}$. Wegen $N_{K/\mathbb{Q}}(a) = 5$ ist $u$
damit Primitivwurzel mod $a$.

Insbesondere ist jetzt $2 \equiv -u \pmod{a}$, und wir erhalten
\[
\frac{2}{a} = \frac{-u}{a} + 1 + \sqrt{14} = 1 + \sqrt{14} + \frac{1}{-a/u}.
\]
Wegen $-a/u = 11 - 3\sqrt{14}$ haben wir
$(1 + \sqrt{14})/2 = [-1, 1 + \sqrt{14}, 11 - 3\sqrt{14}]$,
d.h. $(1 + \sqrt{14})/2$ ist ein
Kettenbruch der Länge 2.

Wir beweisen nun eine wichtige Eigenschaft der Mengen $E'_j$ (sh. §
1), die trivialer aussieht als sie ist, nämlich
\begin{quote}
  {\bf (0.12)} {\em Ist $c \in E'_k$ und $b \mid c$ für ein $b \in R$,
    dann ist auch $b \in E'_k$.}
\end{quote}

\begin{proof}
  Wir werden etwas mehr zeigen: sind $B, C$ ganze Ideale in $R$ (nicht
  notwendig Hauptideale) und ist $B \mid C$, dann gilt: hat jede prime
  Restklasse mod $C$ einen Vertreter in $E'_{k-1}$, dann gilt dasselbe
  für die primen Restklassen mod $B$. Indem wir Induktion über die
  Anzahl der Primidealteiler von $CB^{-1}$ machen, dürfen wir $C=BP$
  für ein Primideal $P$ annehmen.

  Sei nun ein $r \in R$ gegeben mit $(r) \cdot B = R$; wir suchen ein
  $e \in E'_{k-1}$ mit der Eigenschaft $r \equiv e \pmod{B}$. Dazu
  unterscheiden wir: a) $B \cdot P = R$: dann gibt es nach dem
  chinesischen Restsatz ein $s \in R$ das den Kongruenzen $s \equiv r
  \pmod{B}$ und $s \equiv 1 \pmod{P}$ genügt. Mit diesem $s$ gilt dann
  $(s) \cdot C = (s) \cdot BP = R$, und weil jede prime Restklasse mod
  $C$ nach Voraussetzung einen Repräsentanten $e \in E'_{k-1}$ hat,
  ist $s \equiv e \pmod{C}$, also erst recht $r \equiv s \equiv e
  \pmod{B}$.  b) $B \equiv 0 \pmod{P}$: dann ist auch $(r) \cdot BP =
  R$, und wir finden sofort ein $e \in E'_{k-1}$ mit $r \equiv e
  \pmod{C}$, womit erst recht $r \equiv e \pmod{B}$ ist.

  Das war zu zeigen.
\end{proof}

Als nächstes zeigen wir ein Analogon zu (0.1):

\begin{quote}
  {\bf (0.13)} {\em Sei $b \in R \setminus E'_2$; dann ist
        $M^2(f) \ge \frac{1}{f(b)}$.}
\end{quote}

\begin{proof}
  Sei $b \in R \setminus E'_2$. Dann gibt es ein $a \in R$ mit $(a,b)=1$,
  sodass $a$ keiner Einheit kongruent mod $b$ ist. Nun setzen wir
  $a = bq_1 + r_1$, $b = r_1 q_2 + r_2$; wäre $r_2=0$, so folgt $r_1 \mid b$
  und $r_1 \mid a$, also $r_1 \mid (a,b)=1$. Damit wäre $r_2$ Einheit im
  Widerspruch zu $a \equiv r_1 \pmod{b}$ und der Wahl von $a$. Die
  Behauptung folgt nun wie in (0.1).
\end{proof}

Auch (0.13) ist bestmöglich in dem Sinne, dass die angegebene Schranke
manchmal exakt ist. So ist z.B. für $R = D(6)$ $M^2(f) = \frac14$, wenn
$f$ der Absolutbetrag der Norm ist, wobei die Ungleichung
$M^2(f) \ge \frac14$ aus (0.13) folgt, weil $2 \in E'_2$ ist.

Man bemerkt nun, dass man (0.1) und (0.13) in folgender Form
aufschreiben kann: Sei $b \in E'_k$; dann ist $M^k(f) \ge 1/f(b)$. Der
Fall $k=1$ entspricht dabei (0.1), der Fall $k=2$ dagegen (0.13). Ob
diese Aussage allerdings auch für $k=3$ gilt, ist zweifelhaft.

Ein weiterer interessanter Begriff zur Beschreibung $k$-stufig
euklidischer Ringe stammt von Cooke und Weinberger \cite{CW75}: Sei $R$ ein
Ring (hier also Integritätsbereich mit Eins) und $K$ sein
Quotientenkörper; existiert dann eine von $(a,b)$ ausgehende
Teilerkette der Länge $\le k$ mit $f(r_k) < f(b)$ für $a,b \in R
\setminus \{0\}$, so heißt $K$ in $x=a/b$ \textbf{$k$-stufig
  euklidisch} bezüglich $f$. Die Menge aller $x \in K$, in denen $K$
$k$-stufig euklidisch bezüglich $f$ ist, wollen wir mit $F_k$
bezeichnen.\label{dpFi} Damit ist
$F_1 \subset F_2 \subset \dots \subset F_\infty
  = \bigcup_{k=1}^\infty F_k \subset K$.

Offenbar ist $R$ genau dann $k$-stufig euklidisch bezüglich $f$, wenn
bereits $F_k = K$ ist, und quasi-euklidisch, wenn $F_\infty = K$
ist. Gibt es nun ein $k \in \mathbb{N}$ derart, dass $F_k = F_\infty$
ist (unabhängig davon, ob $F_\infty = K$ ist oder nicht), so nennt man
das kleinste solche $k$ die \textbf{euklidische Tiefe} von $K$
bezüglich $f$.

Cooke und Weinberger konnten 1975 folgende Aussagen beweisen:
\begin{quote}
  {\bf (0.14)} {\em Sei $R$ der Ring ganzer Zahlen in einem
    algebraischen Zahlkörper $K$ mit Einheitenrang $\ge 1$ und sei $f$
    der Absolutbetrag der Norm. Sind dann gewisse Riemannsche
    Vermutungen richtig, dann gibt es für alle $a,b \in R$ mit
    $(a,b)=1$ eine von $(a,b)$ ausgehende, abbrechende Teilerkette der
    Länge $k \le 5$. Hat $K$ außerdem eine reelle Einbettung, so gibt
    es sogar solche der Länge $k \le 3$.}
\end{quote}

Aus diesem Satz leiten sie dann folgende Korollare ab:
\begin{quote}
  {\bf (0.15)} {\em Es ist $F_5 = F_\infty$: die euklidische Tiefe von
    $K$ ist $\le 5$. Hat $K$ eine reelle Einbettung, so ist bereits
    $F_3 = F_\infty$. }
\end{quote}
\begin{proof}
  Wir nehmen an, es gibt eine von $(a,b)$ ausgehende, abbrechende
  Teilerkette der Länge $k$. Dann ist $(a,b) = (r_{k-1})$ ein
  Hauptideal, und indem wir $c = a/r_{k-1}$ und $d = b/r_{k-1}$
  setzen, wird $(c,d) = 1$. Nach (0.14) gibt es dann eine von $(c,d)$
  ausgehende, abbrechende Teilerkette der Länge $\le 5$ (bzw. $\le 3$,
  falls $K$ eine reelle Einbettung hat), und Multiplikation mit
  $r_{k-1}$ liefert eine solche für $(a,b)$. Wir haben damit gezeigt,
  dass $x = a/b \in F_k$ sogar $x \in F_5$ (bzw. $x \in F_3$)
  impliziert, qed.
\end{proof}

\begin{quote}
  {\bf (0.16)} {\em Hat $R$ Klassenzahl 1, dann ist $R$ 4-stufig
    normeuklidisch; hat außerdem $K$ eine reelle Einbettung, so ist $R$
    sogar 2-stufig normeuklidisch.}
\end{quote}
\begin{proof}
  Sei $x \in K$; wir dürfen $x = a/b$ mit $a,b \in R$, $(a,b) \neq 1$
  schreiben. Wir suchen eine von $(a,b)$ ausgehende Teilerkette der
  Länge $\le 4$ mit $|N_{K/\mathbb{Q}}(r_{k-1})| \le
  |N_{K/\mathbb{Q}}(b)|$; für $|N_{K/\mathbb{Q}}(b)| = 1$ ist nichts
  zu zeigen. Ansonsten liefert (0.14) eine von $(a,b)$ ausgehende,
  abbrechende Teilerkette der Länge $k \le 5$; für diese ist $r_{k-1}$
  als Teiler von $(a,b)$ eine Einheit, d.h. es existiert eine von
  $(a,b)$ ausgehende Teilerkette der Länge 4 mit $1 =
  |N_{K/\mathbb{Q}}(r_{k-1})| \le |N_{K/\mathbb{Q}}(b)|$. Entsprechend
  folgt die Behauptung, falls $K$ eine reelle Einbettung besitzt.
\end{proof}

Am Schluss ihrer Arbeit schreiben Cooke und Weinberger, dass man bei
Berücksichtigung eines Theorems von Lenstra unter den Voraussetzungen
von (0.14) eine von $(a,b)$ ausgehende, abbrechende Teilerkette der
Länge $\le 4$ erhält, falls $R$ Klassenzahl 1 hat. Wie in (0.16) folgt
dann, dass solche Ringe schon 3-euklidisch sind.

Wie es scheint, hat Cooke nicht bemerkt, dass man (0.15) ganz einfach
verbessern kann; es gilt nämlich

\begin{quote}
  {\bf (0.15')} {\em Die euklidische Tiefe ist $\le 4$; hat $K$ eine
    reelle Einbettung, so gilt sogar $F_2 = F_\infty$. }
\end{quote}

\begin{proof}
  Im Beweis von (0.15) haben wir gesehen, dass sich $x = a/b$ als
  Kettenbruch $a_k/b_k$ der Länge $k \le 5$ schreiben lässt. Ist $b_k
  \in R^\times$ eine Einheit, so ist $a_k/b_k$ Kettenbruch der Länge 1
  und somit nichts zu zeigen. Andernfalls ist wegen
  $$ \frac{a_k}{b_k} - \frac{a_{k-1}}{b_{k-1}} =
     \frac{(-1)^{k-1}}{b_kb_{k-1}} $$
  $y = a_{k-1}/b_{k-1}$ ein Kettenbruch der Länge $k-1
  \le 4$, der $x$ genügend genau approximiert, um $x \in F_4$ zu
  zeigen. Entsprechend folgt die Behauptung, wenn $K$ eine reelle
  Einbettung besitzt.
\end{proof}

Insbesondere haben also quadratische Zahlkörper euklidische Tiefe $\le 2$
(falls gewisse Riemannsche Vermutungen richtig sind). In § 2 werden
wir z.B. sehen, dass die Körper $\mathbb{Q}(\sqrt{m})$ für $m = 15, 26,
85$ euklidische Tiefe 2 haben und dass hier $M^k(K) = 1$ für alle
$k\ge2$ gilt; dagegen hat $\mathbb{Q}(\sqrt{30})$ zwar ebenfalls
euklidische Tiefe 2, jedoch ist hier $M^k(K) = \frac32$ für alle $k \ge 2$.

Schließlich wollen wir der Ähnlichkeit der Mengen $E_k$ und $E'_k$
eine Deutung geben: dazu nennen wir eine Abbildung $f:R \to
\mathbb{N}$ \textbf{semi-euklidisch} auf $R$, wenn für alle $a,b\in
R\setminus\{0\}$ gilt:
\begin{align*}
\text{(E--1)} \quad & f(a)= 0 \Leftrightarrow a = 0; \\ \text{(E'--2)}
\quad & \text{ist $(a,b)= 1$ so gibt es ein $q\in R$ mit
  $f(a-bq)<f(b)$.}
\end{align*}
Falls ein solches $f$ existiert, heißt der Ring $R$
\textbf{semi-euklidisch}. Wie für gewöhnlich euklidische Ringe gilt
nun
\begin{quote}
  {\bf (0.17)} {\em Ein Ring $R$ ist genau dann semi-euklidisch, wenn
    $R = E'_\infty$ gilt. In diesem Fall ist $f_0(a) = \min
    \{j\in\mathbb{N}: a\in E'_j\}$ eine semi-euklidische Funktion auf
    $R$, und diese stimmt mit der durch
    $$ f_m(a) = \min \{f(a): f
                \text{ ist semi-euklidische Funktion auf } R\} $$
    definierten minimalen semi-euklidischen Funktion überein.}
\end{quote}
Der Beweis verläuft exakt wie im euklidischen Fall.

\begin{quote}
  {\bf (0.18)} {\em Sei $f$ multiplikativ; dann sind äquivalent:
    \begin{enumerate}
    \item[(i)] für alle $a,b\in R\setminus\{0\}$ mit $(a,b)=1$ existiert
      ein $q\in R$ mit $f(a-bq)<f(b)$;
    \item[(ii)] für alle $a,b\in R\setminus\{0\}$ mit $(a,b)=d$ existiert
      ein $q\in R$ mit $f(a-bq)<f(b)$.
    \end{enumerate} }
\end{quote}
\begin{proof}
  Die Richtung ii) $\Rightarrow$ i) ist trivial; sei also i)
  richtig. Ist dann $(a,b)=(d)$ für ein $d\in R$, so setzen wir $a'=
  a/d$, $b' = b/d$, und finden wegen i) ein $q\in R$ mit $f(a'-qb')<
  f(b')$. Weil $f$ multiplikativ ist, folgt dann $f(a-bq)< f(b)$ durch
  Multiplikation mit $f(d)$.
\end{proof}

Als Korollar hiervon hat man sofort
\begin{quote}
  {\bf (0.19)} {\em Ein Hauptidealring $R$ ist genau dann
    semi-euklidisch bezüglich einer multiplikativen Funktion $f$, wenn
    $R$ euklidisch bezüglich $f$ ist.}
\end{quote}

Will man daher Zahlringe finden, die semi-euklidisch, aber nicht
euklidisch bezüglich der Norm sind, so muss man unter Ringen mit
Klassenzahl $> 1$ suchen. Eine elementare, aber etwas langwierige
Rechnung zeigt
\begin{quote}
  {\bf (0.20)} {\em Ist $R=D(m)$, $m\le0$ quadratfrei, dann ist $R$
    genau dann semi-euklidisch bezüglich einer Funktion $f$, wenn $R$
    schon euklidisch ist.}
\end{quote}

Als Beispiel betrachten wir $D(-5)$: hier ist $E'_1 = \{0, -1, +1\}$
und für alle $a \in E'_2$ muss $N_{K/\mathbb{Q}}(a) \le 2$ gelten (weil
$E'_1$ nur zwei vom 0 verschiedene Elemente enthält). Also ist
$N_{K/\mathbb{Q}}(a) \in \{2, 3, 4, 6\}$ und wir sehen $a \in \{+2, +1
\pm \sqrt{-5}\}$; nun ist $\{+1, -1\}$ zwar ein primes Restsystem mod
$(+1 \pm \sqrt{-5})$, aber nicht mod 2, und wir haben $E'_2 = \{0, +1,
+1 \pm \sqrt{-5}\}$. Jetzt betrachtet man alle $a \in D(-5)$ mit
$N_{K/\mathbb{Q}}(a) \le 6$ und findet $E'_2 = E'_3 = \dots =
E'_\infty$, und aus (0.14) folgt, dass $D(-5)$ nicht semi-euklidisch
ist.

Aus dem Beweis von (0.20) ergibt sich übrigens, dass im Ring $D(-15)$
$E'_2 \subsetneq E'_3 = E'_4$ gilt, während für alle anderen $D(m)$
bereits $E'_2 = E'_3$ ist.

Schon der reellquadratische Fall zeigt, dass (0.20) für algebraische
Zahlkörper nicht typisch ist. So zeigten beispielsweise Johnson, Queen
und Sevilla \cite{Joh85}, dass $D(10)$ und $D(65)$ semi-euklidische Ringe
bezüglich der Norm mit Klassenzahl 2 sind. Wahrscheinlich sind dies
sogar die einzigen reellquadratischen, bezüglich der Norm
semi-euklidischen Ringe mit von 1 verschiedener Klassenzahl. Wir
werden später sehen, dass auch die Ringe ganzer Zahlen in den
biquadratischen Körpern $\mathbb{Q}(\sqrt{-1},\sqrt{15})$ und
$\mathbb{Q}(\sqrt{-3},\sqrt{13})$ Klassenzahl 2 haben und
semi-euklidisch sind.

In Analogie hierzu definieren wir $k$-stufig semi-euklidische Ringe
durch die Forderung, dass es für alle $a,b \in R \setminus \{0\}$ mit
$(a,b) = 1$ eine von $(a,b)$ ausgehende Teilerkette der Länge $\le k$
gibt mit $f(r_k) < f(b)$. Damit wäre auch klar, was man unter einem
quasi-semi-euklidischen Ring zu verstehen hat; glücklicherweise
existiert jedoch in der mathematischen Literatur bereits ein Begriff
mit derselben Bedeutung: solche Ringe heißen $\mathrm{GE}_2$-Ringe.

Zum Schluss dieses Paragraphen wollen wir zeigen, wie man die
euklidischen Minima $M(K)$ bestimmen kann, wenn $K$ ein
imaginärquadratischer Zahlkörper ist.
\begin{quote}
  {\bf (0.21)} {\em Sei $R=D(m)$ imaginärquadratisch; dann gilt }
    $$ M(K) =
    \begin{cases}
      \dfrac{1+|\,m\,|}{4} & \text{falls } m \equiv 2, 3 \pmod{4} \\[1em]
      \dfrac{(1+|\,m\,|)^2}{16|\,m\,|} & \text{falls } m \equiv 1 \pmod{4}.
    \end{cases} $$
\end{quote}
Insbesondere haben wir für die euklidischen Ringe die folgende Tabelle:

\begin{center}
$$ \begin{array}{c|ccccc}
\toprule
\rsp m    & -1 & -2 & -3 & -7 & -11 \\  \midrule
\rsp M(K) & \frac12 & \frac34 & \frac13 & \frac47 & \frac9{11} \\
\bottomrule
\end{array}$$
\end{center}

\begin{proof}
  Ist $m \equiv 2, 3 \pmod{4}$, dann gibt es für alle $x \in K$ ein $y
  \in R$ mit $x-y = a+b\sqrt{m}$ und $|a|, |b| \le \frac12$. Damit wird
  $N_{K/\mathbb{Q}}(x-y) = a^2 - m b^2 \le (1+|\,m\,|)/4$. Um auch
  $M(K) \ge (1+|\,m\,|)/4$ einzusehen, betrachten wir z.B. den Punkt
  $x = (1+\sqrt{m})/2$; für alle $y \in R$ mit $x-y = a+b\sqrt{m}$ ist
  dann offenbar $|a|, |b| \ge \frac12$, und es folgt
  $N_{K/\mathbb{Q}}(x-y) \ge (1+|\,m\,|)/4$.

  
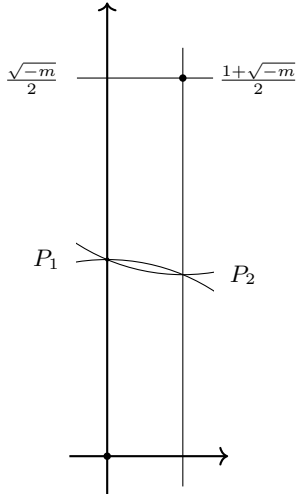
\begin{center}
    \begin{tikzpicture}[scale=2]
      \draw[->,thick] (-0.25,0) -- (0.8,0);
      \draw[->,thick] (0,-0.25) -- (0,3);
      \fill (0,0) circle (0.7pt);
      \fill (0.5,2.5) circle (0.7pt);
      \draw (0.5,-0.2) -- (0.5,2.7);
      \fill (0,1.3) circle (0.4pt);
      \node at (-0.4,1.3) {$P_1$};
      \node at ( 0.9,1.2) {$P_2$};
      \node at (   1,2.5) {$\frac{1+\sqrt{-m}}2$};
      \node at (-0.5,2.5) {$\frac{\sqrt{-m}}2$};
      \draw (-0.2,2.5) -- (0.7,2.5);
      \begin{scope}
        \clip (-0.2,0) rectangle (0.7,3);
        \draw (0.5,2.5) circle (1.3);
        \draw (0,0) circle (1.3);
      \end{scope}
    \end{tikzpicture}    
  \end{center}

  Im Falle $m \equiv 1 \pmod{4}$ verwenden wir eine geometrische
  Methode: wir betten $K$ via $a+b\sqrt{m} \to (a,b\sqrt{|m|})$ in den
  $\mathbb{R}^2$ ein. Ist $\|\cdot\|$ die gewöhnliche euklidische Norm
  im $\mathbb{R}^2$, so stimmen $N_{K/\mathbb{Q}}$ und $\|\cdot\|^2$
  überein in dem Sinne, dass $N_{K/\mathbb{Q}}(a+b\sqrt{m}) =
  \|(a,b\sqrt{|m|})\|^2$ ist. Um nun zu zeigen, dass es zu jedem $x \in
  K$ ein $y \in R$ gibt mit $N_{K/\mathbb{Q}}(x-y) < k$, dürfen wir uns auf
  diejenigen $x \in K$ mit $x=r+s\sqrt{m}$ und $|r|, |s| \le \frac12$
  beschränken, weil wir alle anderen $x$ durch Addition gewisser $y \in R$
  erreichen können. Aus Symmetriegründen dürfen wir sogar $0 \le r, s \le
  \frac12$ annehmen. Wir setzen nun $k = (1+|m|)/(4\sqrt{|m|})$ und
  behaupten, dass die Kreise um die Punkte $0$ und $\frac{1+\sqrt{m}}2$
  (genauer um deren Bilder $(0,0)$ und $(\frac12, \frac{\sqrt{|m|}}2)$)
  mit Radius $k$ ganz $F_\varphi$ bedecken, wo
  $F_\varphi = \{x=r+s\sqrt{m}: 0 \le r, s \le \frac12\}$ ist. Dazu
  rechnen wir einfach nach, dass sich diese beiden Kreise in
  $F_\varphi$ genau in den Punkten
  $$ P_1 = \Big(0, \frac{(1+|m|)\sqrt{|m|}}{(4|m|} \Big)
                 \quad \text{ und } \quad
     P_2 = \Big(\frac12, \frac{(m-1)\sqrt{|m|}}{4|m|} \Big) $$
  schneiden. Bezeichnen wir die obige Einbettung mit $\varphi$, so
  gibt es also zu jedem $x \in K$ ein $y \in R$ mit $\|\varphi(x-y)\|
  \le k$, d.h. mit $N_{K/\mathbb{Q}}(x-y) \le k^2$ wie behauptet. Die
  Punkte $P_1$ und $P_2$ zeigen auch, dass diese Schranke bestmöglich
  ist, womit wir alles bewiesen haben.
\end{proof}

Zu (0.21) wollen wir noch einige Bemerkungen machen. Dazu sei $f:R \to
\mathbb{N}$ eine Funktion, die (E--1) erfüllt, $a,b \in R \setminus
\{0\}$ und $M(f;a,b) = \inf \{f(a-bq)/f(b): q \in R\}$; wir finden
dann $M(f) = \sup \{M(f;a,b): a,b \in R \setminus \{0\}\}$. Ist $f$
darüberhinaus multiplikativ, so dürfen wir einfacher $M(f;x) = \inf
\{f(x-q): q \in R\}$, $x = a/b$ schreiben. Sind nun $a,b \in R
\setminus \{0\}$ und ist ein positives, reelles $\varepsilon > 0$
gegeben, so können wir, falls $M(f) < \infty$ ist, ein $q \in R$
finden mit $f(a-bq)/f(b) < M(f) + \varepsilon$. Dagegen ist es
denkbar, dass wir kein $q \in R$ finden, sodass $f(a-qb)/f(b) \le M(f)$
wird (es könnte ja eine Folge $(q_n)_{n \in \mathbb{N}}$ geben mit
$f(a-q_n b)/f(b) \to M(f)$, $n \to \infty$). Falls wir für alle $a,b
\in R \setminus \{0\}$ die Ungleichung $f(a-qb)/f(b) \le M(f)$
erfüllen können, so sagen wir, dass $M(f)$ erreicht wird. Wie das
Beispiel der minimalen euklidischen Funktion zeigt, braucht es im
Falle eines erreichten euklidischen Minimums $M(f)$ keine $a,b \in R$
zu geben mit $M(f;a,b) = M(f)$. Wie wir in (0.21) gesehen haben, ist
dies aber in imaginärquadratischen Zahlkörpern der Fall, wenn man für
$f$ die Norm nimmt, nämlich für $P = (1+\sqrt{m})/2$, falls $m \equiv
2, 3 \pmod{4}$ ist, und für die Punkte $P_1$ und $P_2$, falls $m
\equiv 1 \pmod{4}$ ist.

Wir setzen\label{dDefC} jetzt
$C_1 = \{(a,b) \in R \times R : M(f;a,b) = M(f)\}$ und definieren
\[
M_2(f) = \sup \{ M(f;a,b) : (a,b) \in R \times R \setminus C_1 \}.
\]
Ist $M_2(f) < M(f)$, so heißt $M(f)$ isoliert und $M_2(f)$ das
\textbf{zweite euklidische Minimum}\label{dpEuM2} von $R$ bezüglich
$f$. Wie (0.21) zeigt, ist das euklidische Minimum bezüglich der Norm
in imaginärquadratischen Zahlkörpern nicht isoliert. Dagegen haben
Barnes und Swinnerton-Dyer vermutet, dass dies in Zahlkörpern mit
Einheitenrang $>1$ immer der Fall ist.

Man kann nun in dieser Art weitermachen und sich fragen, ob auch das
zweite euklidische Minimum isoliert ist und dann dementsprechend
$M_3(f)$ definieren usw. Wir werden sehen, dass es Zahlkörper gibt, die
eine unendliche Kette $M_1(K), M_2(K), \dots$ von euklidischen Minima
besitzen (wie z.B. $\mathbb{Q}(\sqrt{5})$), während in anderen
Zahlkörpern (wie $\mathbb{Q}(\sqrt{23})$) bereits $M_2(K)$ nicht
isoliert ist. Schließlich kann man entsprechende Begriffsbildungen
auch für die $k$-stufigen Minima $M_k(f)$ einführen.

\section*{{\sc Anmerkungen zu} \S\ 0}

Wie wir schon eingangs erwähnten, war Gauß der erste, der einen
algebraischen Zahlkörper $\neq \mathbb{Q}$ (nämlich $\mathbb{Q}(i)$)
als normeuklidisch nachwies. Allerdings ist es nicht ganz richtig, dass
er dies tat, um in $\Z[i]$ den Satz von der eindeutigen
Primfaktorzerlegung zu beweisen. Dazu ging er nämlich wie folgt vor:
in Art. 33 (sh. Werke II, theoria residuorum biquadraticum) benutzt er
die Fermat'sche Erkenntnis, dass sich jedes prime $p \equiv 1 \pmod{4}$
als Summe zweier Quadrate schreiben lässt, um zu zeigen, dass solche $p$
in $\Z[i]$ zum Produkt zweier verschiedener Primfaktoren
werden. Mit Hilfe dieser Tatsache führt er in Art. 37 die eindeutige
Primfaktorzerlegung in $\Z[i]$ auf diejenige in $\Z$
zurück. Erst in Art. 46 beschreibt er dann den Euklidischen
Algorithmus zur Bestimmung zweier Zahlen aus $\Z[i]$.

Den Vorschlag, in Zahlkörpern nach euklidischen Funktionen zu suchen,
die vom Absolutbetrag der Norm verschieden sind, hat wohl als erster
H. Hasse \cite{Has28} gemacht. Wie solche Funktionen aussehen könnten,
haben Lenstra \cite{Len74} und Bedocchi \cite{Bed85} beschrieben; wir
werden in §10 auf deren Vorschläge zurück kommen. Der Begriff des
Euklidischen Algorithmus, den wir ganz am Anfang dieses Paragraphen
definiert haben, lässt sich noch verallgemeinern: so kann man in
nichtkommutativen Ringen linkseuklidische (bzw. rechtseuklidische)
Funktionen definieren und z.B. nach euklidischen Quaternionenalgebren
fragen. Damit haben sich z.B. Dickson (\cite{Dic27}, Redei \cite{Red67},
Newman \cite{New72}, Brungs \cite{Bru73} und Lenstra \cite {Len74,Len78}
befasst (sh. auch Felgner \cite{Fel73}).

Weiter kann man statt Abbildungen $f:R \to \mathbb{N}$ Funktionen
betrachten, die $R$ in eine wohlgeordnete Gruppe $W$ abbilden: Samuel
\cite{Sam71} hat gezeigt, dass die Forderung $W = \mathbb{N}$ in
Ringen mit endlichen Restklassenkörpern keine Einschränkung ist in dem
Sinne, dass ein Ring, der bezüglich einem $g: R \to W$ euklidisch ist,
auch bezüglich einem geeigneten $f: R \to \mathbb{N}$ euklidisch
ist. Dagegen hat Hiblot \cite{Hib75} einen Ring konstruiert, für den
beide Definitionen nicht äquivalent sind. Für eine gründliche
Untersuchung solcher und ähnlicher Fragen verweisen wir auf Lenstra
\cite{Len74}.  Einen Zusammenhang der im Text nur kurz gestreiften
$\mathrm{GE}_2$-Ringe mit dem Euklidischen Algorithmus findet man bei
Hurwitz \cite{Hur19}, P. Cohn \cite{Coh66}, Vaserstein \cite{Vas72}
und Cooke \cite{Coo77}.

\chapter*{\S\ 1 Kriterien für die Existenz des Euklidischen Algorithmus}
\setcounter{chapter}{1}
\addcontentsline{toc}{chapter}
                {\S\ 1 Kriterien für die Existenz des Euklidischen Algorithmus}
\markboth{Euklidische Ringe}
         {\S\ 1 Kriterien für die Existenz des Euklidischen Algorithmus}                 

Im Folgenden sei $K$ algebraischer Zahlkörper und $R$ der Ring ganzer
Zahlen in $K$. Um zu entscheiden, ob $R$ normeuklidisch ist oder
nicht, hat man verschiedene Kriterien entwickelt, von denen wir einige
beschreiben wollen. Dazu sei $P$ ein Primideal in $R$; wir nennen ein
Polynom
\[
f(x) = x^n + a_{n-1}x^{n-1} + \dots + a_1 x + a_0 \in R[x]
\]
\textbf{eisensteinsch}\index[dS]{eisensteinsch} bezüglich $P$, wenn
$a_i \in P$ für $0 \le i \le n-1$ und $a_0 \in P \setminus P^2$ ist
(man beachte, dass $f$
normiert sein muss!). Das Eisensteinsche Irreduzibilitätskriterium
besagt dann, dass solche Polynome über $K(x)$ irreduzibel sind.
Wir wollen nun folgende Situation betrachten: es sei $L/K$ eine
endliche Körpererweiterung vom Grade $(L:K) = n$, $S$ der Ring der
ganzen Zahlen in $L$ und $Q$ ein Primideal in $S$ über $P$:

\begin{center}
  \begin{tikzpicture}[scale=0.6]
    \node (Q) at (0,0) {$\Q$};
    \draw (1,0) node {$\supset$};
    \node (Z) at (2,0) {$\Z$};
    \draw (3,0) node {$\supset$};
    \node (p) at (4,0) {$(p)$};
    \node (K) at (0,2) {$K$};
    \draw (1,2) node {$\supset$};
    \node (R) at (2,2) {$R$};
    \draw (3,2) node {$\supset$};
    \node (P) at (4,2) {$P$};
    \node (L) at (0,4) {$L$};
    \draw (1,4) node {$\supset$};
    \node (S) at (2,4) {$S$};
    \draw (3,4) node {$\supset$};
    \node (Q1) at (4,4) {$Q$};
    \draw (Q) -- (K) -- (L);
    \draw (Z) -- (R) -- (S);
    \draw (p) -- (P) -- (Q1);
  \end{tikzpicture}
\end{center}

Das Primideal $P$ heißt
\textbf{rein verzweigt}\index[dS]{Primideal!rein verzweigt} in $L/K$,
wenn $PS = Q^n$ gilt. Bekanntlich ist ein Primideal $P$ in $L/K$ genau
dann rein verzweigt, wenn die Erweiterung von der Wurzel eines bezüglich
$P$ eisensteinschen Polynoms erzeugt wird; der Vollständigkeit halber
wollen wir die Beweise hier geben.

\begin{quote}
  {\bf (1.1)} {\em Sei $f(x) = x^n + a_{n-1}x^{n-1} + \dots + a_0 \in
    R[x]$ eisensteinsch bezüglich $P$ und $\alpha$ eine Wurzel von
    $f$; ist dann $L=K(\alpha)$, so gilt $(L:K) = n = \deg f$ und $P$
    ist rein verzweigt in $L/K$.}
\end{quote}

\begin{proof}
  $f$ ist wegen des Eisensteinschen Irreduzibilitätskriteriums
  irreduzibel, folglich gilt $(L:K) = n$. Außerdem können wir $a_i S =
  P A_i$ schreiben für gewisse ganze Ideale $A_i$ in $S$, wobei noch
  $A_0 = 0 \bmod P$ gilt. Wegen $f(\alpha) = 0$ gilt $\alpha^n =
  -(a_{n-1}\alpha^{n-1} + \dots + a_0) \in PS$. Sei nun $Q$ irgendein
  Primideal in $S$ über $P$; dann ist $\alpha^n \in Q$, also sogar
  $\alpha \in Q$, weil $Q$ prim ist, und damit $\alpha^n \in
  Q^n$. Schreibt man $\alpha^n S = (PS)B$ für ein ganzes Ideal $B$ in
  $S$, so folgt $Q^n \alpha^n S = (PS)B$. Wenn wir noch zeigen können,
  dass $Q + B = S$ ist, haben wir $Q^n \mid PS$ und damit $Q^n = PS$.
  Wir nehmen daher an, es sei $Q \mid B$; dann folgt aber $a_0 \in
  Q(PS)$ wegen
  $$ \alpha^n = -(a_{n-1}\alpha^{n-1} + \dots + a_1)\alpha - a_0 \in Q(PS). $$
  Wegen $a_0 \in R$ und $Q \cap R = P$ impliziert
  dies aber $a_0 \in P^2$ im Widerspruch zur Voraussetzung.
\end{proof}

\begin{quote}
  {\bf (1.2)} {\em Sei $P$ in $L/K$ rein verzweigt, $\pi \in S$ und
    $f$ das Hauptpolynom von $\pi$. Ist dann $\pi \in Q$ und $f(x) =
    x^n + a_{n-1} x^{n-1} + \dots + a_0 \in R[x]$, dann gilt $a_i \in
    P$ für alle $0 \le i \le n-1$; ist außerdem $\pi \in Q \setminus
    Q^2$, so ist $f$ sogar eisensteinsch bezüglich $P$ und damit
    gleich dem Minimalpolynom von $\pi$.}
\end{quote}

\begin{proof}
  Wegen $f(\pi)=0$ ist $a_0 \equiv 0 \bmod {\pi}$, also $a_0 \in Q \cap
  R = P$. Betrachtet man die Gleichung $f(\pi) \equiv 0 \bmod {Q^2}$,
  so erhält man entsprechend $a_1 \equiv 0 \bmod {P}$ usw. Ist
  schließlich $\pi \in Q \setminus Q^2$ und wäre $a_0 \in P^2$, so
  folgte $\pi^n \equiv 0 \bmod {\pi P}$ im Widerspruch zu $PS =
  Q^n$. Also ist $f$ eisensteinsch bezüglich $P$.
\end{proof}

Die Bedeutung rein verzweigter Primideale für diesen Paragraphen beruht auf

\begin{quote}
  {\bf (1.3)} {\em Ist $P$ rein verzweigt in $L/K$, dann gibt es zu
    jedem $\alpha \in S$ ein $a \in R$ mit $\alpha \equiv a \bmod {Q}$
    und es gilt $N_{L/K}(\alpha) \equiv a^n \bmod {P}$. }
\end{quote}

\begin{proof}
  Sei $\alpha \in S$; wegen $S/Q \cong R/P$ ($P$ ist rein verzweigt,
  und insbesondere hat $Q$ den Trägheitsgrad 1 über $P$) existiert ein
  $a \in R$ mit $\alpha \equiv a \bmod {Q}$. Damit ist $\alpha - a \in
  Q$, und das Hauptpolynom von $\alpha - a$ hat nach (1.2) die Gestalt
  $H_{P}(\alpha - a; K)(x) = x^n + a_{n-1} x^{n-1} + \dots + a_0$ mit
  $a_i \in P$ für $0 \le i \le n-1$; also ist $\alpha$ Nullstelle von
  $(\alpha - a)^n + a_{n-1}(\alpha - a)^{n-1} + \dots + a_0$, und wir
  sehen $N_{L/K}(\alpha) = a^n - a_{n-1}a^{n-1} + \dots + (-1)^n a_0
  \equiv a^n \bmod {P}$.
\end{proof}

Alle bekannten Kriterien für die Existenz eines EA, die auf rein
verzweigten Primidealen beruhen, lassen sich auf das folgende
zurückführen:

\begin{quote}
  {\bf (1.4)} {\em Sei $B$ Produkt von paarweise verschiedenen, in
    $L/K$ rein verzweigten Primidealen, $BS=C^n$, und $C$ ein
    Hauptideal in $S$. Gibt es dann ein $x \in R$ mit
    $x^n \equiv e \bmod {B}$ für ein $e \in R$, und existiert kein
    $r \in R$ mit den Eigenschaften (a), (b) und (c), dann ist }$M(L) \ge k$.
    \begin{enumerate}
    \item[(a)] $r \equiv e \bmod {B}$ \\
    \item[(b)] $r = N_{L/K}(\alpha)$ {\em für ein} $\alpha \in S$ \\
    \item[(c)] $|N_{K/\Q}(r)| < k \cdot \|B\|$.
\end{enumerate}
\end{quote}

\begin{proof}
  Wir nehmen an, es sei $M(L) < k$; wegen $C = (c)$ für ein $c \in S$
  gibt es dann ein $q \in S$ mit $|N_{L/\Q}(x - qc)| < k \cdot
  |N_{L/\Q}(c)|$. Setzt man $\alpha = x - qc$, so wird $\alpha
  \equiv x \bmod {C}$ und $|N_{L/\Q}(\alpha)| < k \cdot
  \|C\|$. Mit $r = N_{L/K}(\alpha)$ wird jetzt
  (a) $r = N_{L/K}(\alpha) \equiv x^n \equiv e \bmod {B}$; nach (1.3)
  ist diese Kongruenz nämlich für jeden Primteiler $P$ von $B$
  richtig. Da diese $P$ nach Voraussetzung paarweise verschieden sind,
  ergibt sich die Behauptung aus dem chinesischen Restsatz.
  (b) $r = N_{L/K}(\alpha)$ gilt nach Konstruktion.
  (c) $|N_{K/\Q}(r)| = |N_{K/\Q}(N_{L/K}(\alpha))| =
  |N_{L/\Q}(\alpha)| < k \|C\| = k \|B\|$, denn für Primideale
  $Q$ in $S$ über $P$ mit relativem Trägheitsgrad 1 ist $\|Q\|=\|P\|$,
  sodass $\|C\|=\|B\|$ aus der Multiplikativität der Idealnorm
  folgt. Da es nach Voraussetzung aber kein $r \in R$ mit den
  Eigenschaften a), b), c) gibt, muss $M(L) \ge k$ sein.
\end{proof}

Wir wollen ausdrücklich bemerken, dass $L$ auch dann nicht
normeuklidisch ist, wenn wir in (1.4) nur $k=1$ wählen können: dann
ist nämlich $M(N_{L/\Q}; x, c) = 1$, d.h. es gibt kein $q \in
R$ mit $|N_{L/\Q}(x-qc)| < |N_{L/\Q}(c)|$.

Bisher war (1.4) nur im Fall $K=\Q$, $k=1$ bekannt; um dieses
Kriterium anwenden zu können, musste $L/\Q$ also rein
verzweigte Primideale enthalten. Die hier gegebene Version lässt sich
dagegen auch dann anwenden, wenn es einen Zwischenkörper $K$ mit
$\Q \subset K \subset L$ gibt, sodass $L/K$ rein verzweigte
Ideale besitzt (ein typisches Beispiel hierfür sind z.B. die
Dirichletschen Zahlkörper, das sind quadratische Erweiterungen von
$\Q(i)$, $i^2=-1$).

Die (bisher notwendige) Bedingung, dass $L/\Q$ rein verzweigte
Primideale enthalten muss, hat letztendlich dazu geführt, dass über
den Euklidischen Algorithmus in reinen Zahlkörpern
(Cioffari\index[dN]{Cioffari} \cite{Cio79}, Egami\index[dN]{Egami}
\cite{Ega79,Ega84}) oder in zyklischen Zahlkörpern
(Heilbronn\index[dN]{Heilbronn} \cite{Hei50,Hei51}, Egami \cite{Ega84})
viel mehr bekannt ist als in Körpern, die keine rein
verzweigten Primideale enthalten. Im Spezialfall $K=\Q$ wird aus (1.4)

\begin{quote}
  {\bf (1.5)} {\em Sei $(L:\Q) = n$, $f \in \N$ ein
    Produkt von in $L/\Q$ rein verzweigten, paarweise
    verschiedenen rationalen Primzahlen und $fS=C^n$ für ein
    Hauptideal $C$ in $S$. Ist dann $a$ ein $n$-ter Potenzrest mod $f$
    und sind alle $r \in \Z$ mit $a \equiv r \bmod {f}$ und
    $|r| < k \cdot f$ keine Normen aus $S$, dann ist $M(L) >
    k$. Darüberhinaus ist $L$ nicht normeuklidisch, falls $k \ge 1$
    gilt.}
\end{quote}

\begin{proof}
  (1.4) mit $K=\Q$, $B=\{f\}$, $e=a$ und $R=\Z$.
\end{proof}

Für $k=1$ schließlich erhält man aus (1.5) das bisher bereits bekannte
Kriterium (Egami\index[dN]{Egami}  1979):

\begin{quote}
  {\bf (1.6)} {\em Sei $(L:\Q)=n$, $f \in \N$ Produkt
    von rein verzweigten, paarweise verschiedenen rationalen
    Primzahlen; gibt es dann $a, b \in \N$ mit $f=ab$, sodass
    $a$ ein $n$-ter Potenzrest mod $f$ ist und weder $a$ noch $-b$
    Normen aus $S$ sind, dann ist $L$ nicht normeuklidisch. }
\end{quote}

\begin{proof}
  (1.5) mit $k=1$; man beachte, dass $a$ und $-b$ die einzigen $r \in
  \Z$ sind mit $r \equiv a \bmod {f}$ und $|r| < f$.
\end{proof}

Wir demonstrieren die Anwendung von (1.5) an einigen einfachen
Beispielen:
\begin{enumerate}
\item[1.] $L = \Q(\sqrt{7})$, $f=14$, $a=9$; dann ist $a$
  quadratischer Rest mod 14, und die $r \in \Z$ mit $r \equiv a \bmod
  {f}$ sind $r = \dots, -19, -5, 9, 23, \dots$ usw. Wählen wir nun
  $k = \frac{9}{14}$, so erfüllt nur $r=-5$ die Bedingung
  $|r| < k \cdot f = 9$;
  nun gibt es aber kein $\alpha \in S$ mit $N_{L/\Q}(\alpha) = -5$,
  denn das Primideal $(5)$ bleibt beim Übergang von $\Z$ nach $S$
  träge wegen $(\frac{28}{5}) = -1$ (Legendre-Symbol). Also ist
  $M(L) \ge \frac{9}{14}$. Um auch ein $z \in L$ zu finden mit
  $M(z) = \frac{9}{14}$ (wir haben hier statt $M(N_{L/\Q}, z)$
  einfach $M(z)$ geschrieben), müssen wir
  zuerst ein Hauptideal $C$ mit Norm $f=14$ finden: ein solches ist
  z.B. $C = (7-3\sqrt{7})$. Dann brauchen wir noch ein $x \in \Z$ mit
  $x^2 \equiv a \equiv 9 \bmod {14}$, z.B. $x=3$, und wir haben $z-x/c
  = -(21+9\sqrt{7})/14$; wegen $M(z) = M(z-y)$ für jedes $y \in R$
  brauchen wir $z \bmod R$ anzugeben, d.h. wir haben jetzt
  $z \equiv (7+5\sqrt{7})/14 \bmod {R}$.
\item[2.] $L = \Q(\sqrt{10})$: um $M(L) \ge \frac32$ zu zeigen, verwenden
  wir $f=10$ und $a=5$; tatsächlich ist $a$ wegen $5^2 \equiv 5 \bmod
  {10}$ quadratischer Rest mod $f$. Die einzigen $r \in \Z$ mit $|r| <
  k \cdot f$ und $r \equiv 5 \bmod {10}$ sind $r = \pm 5$; weil jedoch
  das Primideal über $(5)$ kein Hauptideal in $S$ ist, gibt es in $S$
  keine Elemente der Norm 5, und wir haben $M(L) \ge \frac32$.
\item[3.] $L = \Q(\sqrt{15})$: wir behaupten $M(L) \ge \frac75$ und
  verwenden dazu $f=15$, $a=6$. Damit wird $r \in \{-9, 6\}$. Weil nun
  das Primideal über $(3)$ kein Hauptideal ist, sind alle Elemente der Norm
  9 von der Form $3u^k$, wo $u=4+\sqrt{15}$ die FE von $L$ ist. Wegen
  $N_{L/\Q}(u) = +1$ ist aber $N_{L/\Q}(3u^k) = +9$, d.h. es gibt in
  $S$ kein $\alpha$ mit $N_{L/\Q}(\alpha) = -9$. Analog wird jedes
  Ideal der Norm 6 von einem Element der Form $\beta = (3+\sqrt{15})u$
  erzeugt, und wegen $N_{L/\Q}(\beta) = -6$ gibt es in $S$ keine
  Elemente der Norm $+6$. Dies war zu zeigen.
\end{enumerate}
Das folgende Kriterium enthält (1.4) als Spezialfall (für $m = 0,1$)
und ließe sich noch weiter verallgemeinern; wir werden es jedoch nur
in der vorgestellten einfacheren Version benutzen:

\begin{quote}
  {\bf (1.7)} {\em Es sei $P$ ein Primideal in $R$ mit $PS = Q^m$, $B$
    ein Produkt von paarweise verschiedenen Primidealen mit $B \cap P
    = R$, und es sei $B' = P^m B$, sowie $B'S = Q^{mn} C^n$ für ein
    Hauptideal $Q^m C$. Ist dann $x^n \equiv e \bmod {P}$ für $x, e \in R$,
    und existiert ein $y \in S$ mit $y \equiv x \bmod {QC}$, sodass
    es kein $r \in R$ gibt mit den Eigenschaften
    \begin{enumerate}
    \item[(a)] $r \equiv e \bmod {PB}$
    \item[(b)] $r = N_{L/K}(\alpha)$ für
      ein $\alpha \in S$ mit $\alpha \equiv y \bmod {Q^m C}$
    \item[(c)] $|N_{K/\Q}(r)| < k \cdot \|B\|$
    \end{enumerate}
    dann ist $M(L) \ge k$.}
\end{quote}

\begin{proof}
  Sei $M(L) < k$; dann gibt es ein $\alpha \in S$ mit $\alpha \equiv y
  \bmod {Q^m C}$ und $|N_{L/\Q}(\alpha)| < k \cdot \|Q^m C\|$,
  sowie ein $x \in R$ mit $y \equiv x \bmod {QC}$ (denn $QC$ ist
  Produkt von paarweise verschiedenen, rein verzweigten
  Primidealen). Mit $r = N_{L/K}(\alpha)$ lassen sich nun wie in (1.4)
  die Eigenschaften (a), (b), (c) nachweisen.
\end{proof}

Der Schwierigkeit bei der Formulierung von (1.7) liegt die folgende
Tatsache zugrunde: aus $\alpha \equiv a \bmod {Q}$ folgt zwar
$N_{L/K}(\alpha) \equiv a^n \bmod {P}$, aber aus $\alpha \equiv a
\bmod {Q^m}$ folgt eben i.A. nicht $N_{L/K}(\alpha) \equiv a^n
\bmod {P^m}$, sondern wieder nur die schwächere Kongruenz $\bmod {P}$.

Wir geben einige Beispiele:
\begin{enumerate}
\item[1.] $L=\Q(\sqrt{14})$, $K=\Q$, $B'=(4)=(2)^2$, also $P=(2)$,
  $m=2$, $B=(1)$; weiter haben wir $Q = \frt_1=(4+\sqrt{14})$.  Wir
  setzen nun $k=5/4$ und $y=1+\sqrt{14}$. Damit ist $y \equiv 1 \bmod
  {Q}$, also $x \equiv e = 1$, und wir müssen untersuchen, ob es ein
  $r \in \Z$ mit $|r|<5$ und $r \equiv 1 \bmod {2}$ gibt, das Norm
  eines $\alpha \in S$ mit $\alpha \equiv 1+\sqrt{14} \bmod {2}$
  ist. Da (3) in $L/K$ träge ist, gibt es keine Elemente der Norm 3,
  also kommt nur $r=1$ in Frage. Nun gibt es aber keine Einheit
  $\alpha \in S$ mit $\alpha \equiv 1+\sqrt{14} \bmod {2}$, weil die
  FE $15+4\sqrt{14}$ bereits $\equiv 1 \bmod {2}$ ist (damit ist jede
  Einheit $\equiv 1 \bmod {2}$). Man vergleiche dieses Beispiel auch
  mit dem Beweis auf S. 11, dass $\frac{1+\sqrt{14}}2$ kein
  Kettenbruch der Länge 2 ist (insbesondere achte man auf die Rolle,
  die die FE in beiden Beweisen spielt).
\item[2.] $L=\Q(\sqrt{26})$, $K=\Q$, $B'=(4)$, $y=26$, $k = \frac52$;
  hier ist zwar $Q =\frt_1 = (2,\sqrt{26})$ kein Hauptideal, wohl aber
  $Q^2=(2)$. Nach (1.7) müssen wir die $r \in \Z$ mit $r \equiv 0
  \bmod {2}$ und $|r|<10$ untersuchen. Wegen $y \in Q \setminus Q^2$
  ist aber für alle $\alpha \equiv y \bmod {Q}$ auch $\alpha \in Q
  \setminus Q^2$, und man sieht leicht ein, dass dies
  $N_{L/\Q}(\alpha) \equiv 2 \bmod {4}$ impliziert. Also brauchen wir
  nur die $r \equiv 2 \bmod {4}$ zu betrachten. Weil es aber keine
  $\alpha \in S$ der Norm 2 (denn $\frt_1$ ist kein Hauptideal) oder
  der Norm 6 gibt (da das Ideal (3) in $L$ träge bleibt), folgt in der
  Tat $M(L) \ge \frac52$.
\end{enumerate}
Es bleibt die Frage, wie man vorgehen soll, wenn die Kriterien (1.4)
und (1.7) nichts bringen (z.B. wenn die von ihnen gelieferten
Schranken nicht gut genug sind oder es gar keine rein verzweigten
Primideale gibt). Dazu sei $K$ ein Zahlkörper mit $M(K)=k$ und $I=(b)$
ein ganzes Ideal in $R$ (also $b \in R$). Dann gibt es für alle $a \in
R$ ein $q \in R$ mit $|N_{K/\Q}(a-qb)| < k \cdot \|b\|$ wegen
$a-qb \equiv a \bmod {I}$ können wir dies auch so ausdrücken: jede
Restklasse $\bmod {I}$ enthält ein $r \in R$ mit $|N_{K/\Q}(r)|
< k \cdot \|b\|$. Diese Überlegung führt uns auf die Definition
\begin{align*}
  M(K,I) & = \inf \left\{ \kappa \in \R : \forall a \in R \ \exists c
    \in I \text{ mit } |N_{K/\Q}(a-c)| < \kappa \cdot \|I\| \right\} \\
         & = \inf \left\{ \kappa \in \R : \forall x \in I^{-1} \ \exists
            y \in R \text{ mit } |N_{K/\Q}(x-y)|< \kappa \right\}.
\end{align*}
Weil es zu jedem $x \in K$ ein Hauptideal $I$ mit $x \in I^{-1}$ gibt,
folgt sofort aus der Definition
$M(K) = \sup \{ M(K, I) : I \text{ ist Hauptideal in } R \}$.

Ein Hauptideal $I$ mit $M(K,I)<1$ nennen wir ein
\textbf{euklidisches},\index[dS]{Ideal!euklidisches} ein solches mit
$M(K,I)>1$ ein \textbf{nicht euklidisches
  Ideal}.\index[dS]{Ideal!nicht euklidisches}

In (1.4) und (1.7) haben wir nichts anderes gemacht, als $M(K,I)$ für
rein verzweigte Hauptideale $I$ nach unten abzuschätzen. Wir dürfen
daher wegen $M(K)>M(K,I)$ für jedes Hauptideal $I$ i.A. bessere
Schranken für $M(K)$ erwarten, wenn wir für $I$ beliebige Hauptideale
zulassen.

Betrachten wir beispielsweise $K = \Q(\sqrt{29})$; hier
erzeugt $\alpha = \frac{1}{2}(3+\sqrt{29})$ ein Ideal $I$ der Norm
5. Um $M(K,I)$ zu bestimmen, stellen wir zuerst ein vollständiges
Restsystem mod $I$ auf; ein solches ist z.B. $\{0, \pm 1, \pm
2\}$. Jetzt suchen wir in jeder Restklasse mod $I$ ein Element
minimaler Norm $> 1$; in den Restklassen $0, \pm 1$ mod $I$ sind dies
die Elemente $0, \pm 1$ selbst (dies sind sozusagen die
``uninteressanten'' Restklassen). Um in der Restklasse $2$ mod $I$ ein
Element minimaler Norm zu finden, müssen wir uns fragen, ob diese
Restklasse eine Einheit enthält. Nun gilt aber für die FE $u =
\frac{1}{2}(5+\sqrt{29})$ die Kongruenz $u \equiv 1 \bmod {I}$; da sich
nach Dirichlet jede Einheit $e$ in der Form $e = u^l$ für ein $l \in
\Z$ schreiben lässt, ist also $e \equiv 1 \bmod {I}$ für jede
Einheit $e \in R^\times$, und insbesondere gibt es keine Einheiten in den
Restklassen $2$ mod $I$.

Da es in $R$ auch keine Elemente der Norm $2$ oder $3$ gibt (weil die
Ideale (2) und (3) in $R$ prim bleiben), kann die Restklasse $2$ mod
$I$ nur Elemente der Norm $\ge 4$ enthalten, und in der Tat ist $r=2$
ein solches Element. Damit haben wir $M(K,I) = 4/5$ gezeigt (dies ist
etwas besser als die mit (1.5) erhaltene Schranke $M(K) \ge 23/29 =
0.7931\dots$). Man beachte noch, dass in diesem Fall $I = (u-1)$ gilt,
und dass diese Tatsache dafür gesorgt hat, dass $e \equiv 1 \bmod {I}$
für jede Einheit $e \in R^\times$ ist.

Ganz analog läuft der Beweis von $M(K,I) = 16/9$ für $K =
\Q(\sqrt{85})$ und $I = (u-1) = \frac{1}{2}(7+\sqrt{85})$; $I$
hat Idealnorm 9, und $\{0, \pm 1, \pm 2, \pm 3, \pm 4\}$ ist
vollständiges Restsystem mod $I$.

Es ist damit prinzipiell klar, wie man entsprechende Rechnungen in
beliebigen Zahlkörpern für beliebige Hauptideale $I$ durchführt: man
wählt sich ein System von Grundeinheiten und bestimmt alle Restklassen
mod $I$, die Einheiten enthalten (dazu lässt man einfach jede
Grundeinheit die Potenzen von $1$ bis $\Phi(I)$ durchlaufen und
notiert sich die Restklassen, in denen die Produkte dieser Potenzen
liegen). Dann sucht man sich ein Hauptideal $A = (a)$ minimaler Norm
und fragt, in welchen Restklassen mod $I$ die Elemente $ae$, $e \in
R^\times$, liegen usw.

Allerdings werden die nötigen Rechnungen recht schnell umfangreich,
wenn man ein $I$ mit großer Norm zu betrachten hat (z.B. für $K =
\Q(\sqrt{31})$ und $I = (u-1)$, $u = 1520 + 273\sqrt{31}$,
hier ist $\|I\| = 3038$) oder wenn der Einheitenrang von $K$ groß
wird. Glücklicherweise lässt sich das obige Verfahren so modifizieren,
dass es sich recht einfach programmieren lässt und damit von einem
Rechner durchgeführt werden kann. Dazu sagen wir, die Punkte $x_1,
\dots, x_t \in K$ werden von einer Einheit $u \in R^\times$ (zyklisch)
permutiert, wenn die folgenden Kongruenzen gelten: $x_1u = x_2
\bmod {R}$, $x_2u = x_3 \bmod {R}$, $\dots$, $x_tu = x_1 \bmod {R}$. In
diesem Falle gilt übrigens $M(K, x_1) = \dots = M(K, x_t)$, wie man
leicht aus den beiden Beobachtungen
\begin{enumerate}
\item[1.] für alle $y \in R$ ist $M(K, x) = M(K, x-y)$;
\item[2.] für alle Einheiten $u \in R^\times$ ist $M(K, x) = M(K, ux)$
\end{enumerate}
folgt. Wir behaupten nun:

\begin{quote}
  {\bf (1.8)} {\em Sei $R=D(m)$, $m\in\N$ quadratfrei, und sei
    $1\le u$ für eine Einheit $u$ in $R$. Weiter seien Punkte $x_1,
    \dots, x_t$ gegeben, die von $u$ permutiert werden. Ist dann
    $M(K,x_j)<k$, so gibt es ein $z=r+s\sqrt{m} \in K$ mit den
    Eigenschaften
    \begin{enumerate}
    \item[(a)] $z = x_j \bmod {R}$ für ein $j \in \{1, \dots, t\}$;
    \item[(b)] $|N_{K/\Q}(z)| < k$;
    \item[(c)] $|r| < \mu_1$, $|s| < \mu_2$ mit
      $$ \mu_1 = \frac{\sqrt{k}}2 \Big(\sqrt{u} + \frac1{\sqrt{u}} \Big)
      \quad \text{und} \quad
     \mu_2 = \frac{\sqrt{k}}{2\sqrt{m}}\Big(\sqrt{u} + \frac1{\sqrt{u}}\Big). $$
    \end{enumerate} }
\end{quote}

Dies ist im wesentlichen ``Theorem B'' von Barnes und Swinnerton-Dyer;
deren Beweis ist jedoch völlig auf quadratische Zahlkörper
zugeschnitten und lässt sich nicht verallgemeinern. Bevor wir aber
einen Beweis für (1.8) geben, wollen wir uns klarmachen, was (1.8) im
einfachsten Fall $t=1$ aussagt: Wir haben dann ein $x=x_1$ gegeben mit
$xu=x \bmod {R}$; wenn wir dann $M(K,x) \ge k$ zeigen wollen, nehmen
wir an, es sei $M(K,x)<k$. (1.8) garantiert dann die Existenz eines
$z=r+s\sqrt{m} \in K$ mit den Eigenschaften

\begin{enumerate}
\item[(a)] $z = x_j \bmod {R}$ für ein $j \in \{1, \dots, t\}$;
\item[(b)] $|N_{K/\Q}(z)| < k$;
\item[(c)] $|r| < \mu_1$, $|s| < \mu_2$
\end{enumerate}

Um $M(K,x)<k$ zum Widerspruch zu führen, müssen wir also nur zeigen,
dass unter den endlich vielen $z \in K$, die den Bedingungen (a) und (c)
genügen, keines mit $|N_{K/\Q}(z)|<k$ vorkommt.

\begin{proof}[Beweis von (1.8)]
  Wegen $M(K,x_1)<k$ gibt es ein $y \in R$ mit
  $|N_{K/\Q}(x_1-y)| < k$; wir setzen nun $z_1=x_1-y$ und
  wählen $m \in \Z$ so, dass $\sqrt{k/u} \le |z_1 u^m| <
  \sqrt{k u}$ wird (dies ist offenbar immer möglich). Jetzt behaupten
  wir, dass $z=z_1 u^m$ den Bedingungen (a), (b), (c) genügt:
  \begin{enumerate}
  \item[(a)] $z = z_1 u^m = x_1 u^m - y u^m \equiv
    x_1 u^m \equiv x_1 \bmod {R}$,
    wobei $j$ durch die Kongruenz $j \equiv 1+m \bmod {t}$ bestimmt ist.
  \item[(b)]
    $|N_{K/\Q}(z)| = |N_{K/\Q}(z_1)| = |N_{K/\Q}(x_1-y)| < k$.
  \item[(c)] Sei $\tilde{z} = r-s\sqrt{m}$ die Konjugierte von $z$;
    dann haben wir
    $|z'| = |z z'| / |z| = |N_{K/\Q}(z)| / |z| < k/|z| \le \sqrt{k u}$,
    also $|z + z'| \le |z| + |z'| < 2\sqrt{k u}$ und
    $|2s\sqrt{m}| = |z - z'| \le |z| + |z'| < 2\sqrt{k u}$.
  \end{enumerate}
\end{proof}

Diese Schranken lassen sich aber noch verbessern, wenn wir neben den
Ungleichungen $|z|$, $|z'| < \sqrt{k u}$ auch die Ungleichung $|z z'|
< k$ verwenden; die Behauptung folgt nämlich mit $a=\sqrt{k u}$ und
$b=k$ aus

\begin{quote}
  {\bf (1.9)} {\em Seien $x, y, a, b$ positive reelle Zahlen; aus den
    Ungleichungen $x \le a$, $y \le a$, $xy \le b$ folgt dann $x+y \le
    a + \frac{b}{a}$. }
\end{quote}

\begin{proof}
  Es ist $0 \le (a-x)(a-y) = a^2 - (x+y)a + xy \le a^2 - (x+y)a + b$.
\end{proof}

Wir bemerken noch, dass die Voraussetzung $u > 1$ durch $u \ne 1$ ersetzt
werden darf; man muss dann nur in (1.8.c) $u$ durch $|u|$ ersetzen. Im
Falle $N_{K/\Q}(u)=+1$ erhalten wir übrigens die von Barnes
und Swinnerton-Dyer angegebenen Schranken zurück: setzt man
$u=a+b\sqrt{m}$, so ist $\frac{1}{u} = a-b\sqrt{m}$
und damit $(\sqrt{u} + 1/\sqrt{u})^2 = u + 1/u + 2 = 2a+2$, und wir finden
$\mu_2 = \sqrt{\frac{k(a+1)}{2m}}. $

Beispiel: $K = \Q(\sqrt{19})$; sei $x = \frac{20}{57}\sqrt{19}$ und
$u=170+39\sqrt{19}$ die FE von $K$. Wir haben dann $xu =
260+60\sqrt{19} - x = -x \bmod {R}$, sodass wir (1.8) an sich mit $t=2$,
$x_1=x$ und $x_2=-x$ anwenden müssten. Wenn wir aber statt $u$ die
Einheit $u'=-u$ nehmen, haben wir $xu' = x \bmod {R}$, sodass die Wahl
von $u'$ die Rechenzeit halbiert.

Wir zeigen nun, dass $M(K) \ge M(K,x) = \frac{170}{171}$ gilt; die
eine Ungleichung $M(K,x) \le \frac{170}{171}$ folgt dabei aus der
Beobachtung
\[
M(K,x) \le |N_{K/\Q}(x -3-\sqrt{19})| = \frac{170}{171}.
\]
Wir brauchen daher nur noch $M(K,x) \ge k = \frac{170}{171}$ zu zeigen, und
dies machen wir mit (1.8): wäre nämlich $M(K,x)<k$, so gäbe es ein $z
\in K$ mit $z \equiv x \bmod {R}$, das den Bedingungen a), b), c)
genügt. Wegen $a=170$ erhalten wir die Schranke $|s| < 2.12\dots$, für
$z=r+s\sqrt{19}$. Damit kommen für $s$ nur vier Werte in Frage,
nämlich $s = i + 20/57$, $i \in \{-2, -1, 0, 1\}$:

$$ \begin{array}{rrrl}
  i & = & -2: & |N_{K/\Q}(z)| \text{ wird minimal für $r = \pm 7$ mit }
              |N_{K/\Q}(z)| = 457/171; \\
  i & = & -1: & |N_{K/\Q}(z)| \text{ wird minimal für $r = \pm 3$ mit }
              |N_{K/\Q}(z)| = 170/171; \\
  i & = &  0: & |N_{K/\Q}(z)| \text{ wird minimal für $r = \pm 1$ mit }
              |N_{K/\Q}(z)| = 229/171; \\
  i & = & +1: & |N_{K/\Q}(z)| \text{ wird minimal für $r = \pm 6$ mit }
              |N_{K/\Q}(z)| = 227/171.
\end{array} $$

Nach (1.8) ist damit $M(K,x) = \frac{170}{171}$, und in der Tat werden
wir in § 2 sehen, dass sogar $M(K) = M(K,x) = \frac{170}{171}$
gilt. Hätten wir nur die Schranke $|s| < \sqrt{k u}$ verwendet, so
hätten wir acht Werte für $s$ betrachten müssen.

Bevor wir (1.8) auf beliebige Zahlkörper mit Einheitenrang $>1$
verallgemeinern, wollen wir uns am kubischen Fall klar machen, welche
Schwierigkeiten dabei auftreten. Wir zeigen daher zuerst

\begin{quote}
  {\bf (1.10)} {\em Sei $K$ kubischer Zahlkörper mit Einheitenrang 1,
    $u$ eine Einheit mit $1<|u|$, und $x_1, \dots, x_t \in K$ seien
    Punkte, die von $u$ permutiert werden. Weiter sei $\{\alpha_1,
    \alpha_2, \alpha_3\}$ eine $\Q$-Basis von $K$. Ist dann
    $M(K,x_i) < k$, so gibt es ein $z \in K$, $z = r_1\alpha_1 +
    r_2\alpha_2 + r_3\alpha_3$, mit den Eigenschaften
    \begin{enumerate}
    \item[(a)] $z \equiv x_j \bmod {R}$ für ein $j \in \{1,\dots,t\}$;
    \item[(b)] $|N_{K/\Q}(z)| < k$;
    \item[(c)] $|r_i| < \mu_i$ für $i=1, 2, 3$, wobei die $\mu_i$
      positive reelle Zahlen sind, die von den $x_i$ unabhängig sind.
    \end{enumerate} }
\end{quote}

\medskip\noindent

\textbf{Bem.:} Die Schranken $\mu_i$ werden sich aus dem Beweis
ergeben und lassen sich explizit bestimmen. Die Schranken, die wir
speziell in rein kubischen Körpern erhalten werden, sind weitaus
besser als diejenigen, die Cioffari\index[dN]{Cioffari} \cite{Cio79} in
den Fällen $K=\Q(\sqrt[3]{m})$, $m=12, 17, 44$, erhalten hat. Ähnlich
wie im quadratischen Fall lassen sich die hier erreichten Schranken
noch einmal verbessern; wir werden dies in § 4 tun.

\begin{proof}
  OBdA dürfen wir annehmen, dass $K$ reell ist (denn da $K $
  Einheitenrang 1 hat, muss $r=s=1$ sein, wo $r$ die Anzahl der reellen
  Einbettungen von $K$ in $\C$ angibt). Damit dürfen wir die
  Bedingung $1<|u|$ durch $1<u$ ersetzen.
  Nun ist nach Voraussetzung $M(K,x_1)<k$, folglich existiert ein
  $y\in R$ mit $|N_{K/\Q}(x_1-y)|<k$; wir setzen $z_1=x_1-y$ und
  wählen ein $m\in\N$ so, dass
  $$ \sqrt[3]{k/u^2} \le |z_1u^{m-1}| < \sqrt[3]{ku} $$
  wird. Mit $z=z_1u^m$ sind dann die Bedingungen a) und b) erfüllt,
  sodass wir uns nur noch um c) zu kümmern brauchen. Wir bezeichnen nun
  die beiden Konjugierten von $z$ mit $z'$ und $z''$; da die beiden
  Körper $K'$ und $K''$ konjugiert-komplex sind, gilt $|z'|=|z''|$,
  also $|N_{K/\Q}(z)| = |z z' z''| = |z| |z'|^2$ und damit
  $|z'|^2 = |N_{K/\Q}(z)|/|z| < k/|z| \le \sqrt[3]{ku^2}$.

  Wir betrachten nun zuerst den etwas einfacheren Spezialfall des
  reinen kubischen Zahlkörpers $K=\Q(\vartheta)$ für $\vartheta^3=m$,
  $m\in\N$; hier wählen wir die Basis $\{1, \vartheta, \vartheta^2\}$
  und finden $z' = r_1 + r_2 \rho \vartheta + r_3 \rho^2 \vartheta^2$,
  $z'' = r_1 + r_2 \rho^2 \vartheta + r_3 \rho \vartheta^2$, wo $\rho$
  eine primitive dritte Einheitswurzel ist. Jetzt können wir ganz wie
  im quadratischen Fall schließen:
  \begin{align*}
    |3r_1| & = |T_{K/\Q}(z')| = |z+z'+z''| \le |z|+|z'|+|z''| < 3\sqrt[3]{ku}, \\
    |3r_2\vartheta| & = |z+\rho^2 z' + \rho z''| \le |z|+|z'|+|z''|
                        < 3\sqrt[3]{ku}, \\
    |3r_3\vartheta^2| & = |z+\rho z' + \rho^2 z''| \le |z|+|z'|+|z''|
                        < 3\sqrt[3]{ku},
  \end{align*}
und damit haben wir (1.8) in diesem Spezialfall bewiesen mit den
Schranken $\mu_1 = \sqrt[3]{ku}$, $\mu_2 = \sqrt[3]{ku/m}$,
$\mu_3 = \sqrt[3]{ku/m^2}$. Da wir im allgemeinen Fall die Konjugierten
von $z$ nicht explizit angeben können, müssen wir hier einen etwas
anderen Weg einschlagen.

Sei dazu $\{\beta_1, \beta_2, \beta_3\}$ die zu $\{\alpha_1, \alpha_2,
\alpha_3\}$ duale Basis, die bekanntlich durch
$T_{K/\Q}(\alpha_i \beta_j) = \delta_{ij}$ (Kroneckerdelta)
definiert ist. Damit wird
\begin{align*}
  T_{K/\Q}(z \beta_j) & = \sum_{i=1}^{3} r_i
  T_{K/\Q}(\alpha_i \beta_j) = r_j, \text{ also} \\
  |r_j| & = |T_{K/\Q}(z \beta_j)| = |z\beta_j + z'\beta_j + z''\beta_j|
          \le |z\beta_j| + |z'\beta_j| + |z''\beta_j| \\
        & < \sqrt[3]{ku} (|\beta_j| + 2|\beta_j'|).
\end{align*}
Daher müssen wir nur noch angeben, wie wir die $\beta_j$ aus den
$\alpha_i$ bestimmen können. Dazu betrachten wir die $n \times
n$-Matrix $A = (T_{K/\Q}(\alpha_i \alpha_j))$; wir wissen, dass dann
$\det A = \disc_{K/\Q}(\alpha_1, \alpha_2, \dots, \alpha_n) \ne 0$
ist. Sei $B = (b_{ij})$ die Umkehrmatrix von $A$; dann ist $\beta_i =
b_{i1}\alpha_1 + \dots + b_{in}\alpha_n$ (s. z.B. Cohn\index[dN]{Cohn}
\cite{Coh78} oder Marcus\index[dN]{Marcus} \cite{Mar77}).

Da wir uns in § 6 ausführlich mit Dirichletschen Zahlkörpern
beschäftigen werden, sei auch dieser Fall noch explizit
ausgeführt. Sei dazu $L = \Q(\sqrt{m}, \sqrt{n})$ ein
imaginärer, bizyklischer, biquadratischer Zahlkörper (d.h. $m, n \in
\Z, m < 0$). Wir wählen die $\Q$-Basis $\alpha_1 = 1,
\alpha_2 = \sqrt{m}, \alpha_3 = \sqrt{n}, \alpha_4 = \sqrt{mn}$ und
nehmen an, es sei $N_{L / \Q}(x_1 - y) < k$ für ein $x_1 =
r_1\alpha_1 + r_2\alpha_2 + r_3\alpha_3 + r_4\alpha_4 \in L$, ein $y
\in R$ und ein $k \in \R$. Da der Einheitenrang 1
hat, gibt es eine Fundamentaleinheit $u$, die wir so wählen können,
dass $|u| > 1$ ist. Mit $z_1 = x_1 - y$ finden wir dann ein $m \in \Z$ mit
$$ \frac{\sqrt[4]{k}}{\sqrt{|u|}} \le |z_1| < \sqrt[4]{k} \cdot \sqrt{|u|}
   = \delta. $$
Sei nun $\mathrm{Gal}(L/\Q) = \{1, \sigma, \tau, \sigma\tau\}$,
und oBdA sei $\sigma$ der Automorphismus der komplexen
Konjugation (damit ist dann $|z| = |z^\sigma|$; sind z.B. $m$ und $n$
beide negativ, so kann man $\sigma$ auch durch $\sigma: \sqrt{m}
\mapsto -\sqrt{m}, \sqrt{n} \mapsto -\sqrt{n}$ charakterisieren). Nun
ist $N_{L / \Q}(x) = x \cdot x^\sigma \cdot x^\tau \cdot
x^{\sigma\tau} = |x|^2 |x^\tau|^2$, also auch $|x^\tau|^2 \le N_{L /
  \Q}(x) / |x|^2 < \delta$.
Wir setzen jetzt $\tau: \sqrt{m} \mapsto -\sqrt{m}, \sqrt{n} \mapsto
-\sqrt{n}$ und erhalten
\begin{align*}
  4 \cdot |r_1| & = |z + z^\sigma + z^\tau + z^{\sigma\tau}| \le |z| +
                    |z^\sigma| + |z^\tau| + |z^{\sigma\tau}| \\
                & = 2 \cdot (|z| + |z^\tau|) < 4\delta, \quad \text{sowie} \\
  4 \cdot |r_2| \cdot \sqrt{|m|} & = |z - z^\sigma + z^\tau - z^{\sigma\tau}|
                  \le 2 \cdot (|z| + |z^\tau|) < 4\delta, \\
  4 \cdot |r_3| \cdot \sqrt{|n|} & = |z - z^\sigma - z^\tau + z^{\sigma\tau}|
                   \le 2 \cdot (|z| + |z^\tau|) < 4\delta, \\
  4 \cdot |r_4| \cdot \sqrt{|mn|}
                 & = |z + z^\sigma - z^\tau -z^{\sigma\tau}|
                   \le 2 \cdot (|z| + |z^\tau|) < 4\delta.
\end{align*}
\end{proof}

Benutzt man auch (1.9), so kann man wie im reellquadratischen Fall den
Faktor $\sqrt{|u|}$ in $\delta$ ersetzen durch $(\sqrt{|u|} +
1/\sqrt{|u|})/2$, und wir haben

\begin{quote}
  {\bf (1.11)} {\em Sei $L = \Q(\sqrt{m}, \sqrt{n})$ ein
    imaginärer, bizyklischer Zahlkörper, $u$ eine Einheit in $L$ mit
    $|u| > 1$, und $x_1, \dots, x_t \in L$ seien Punkte, die von $u$
    permutiert werden. Ist dann $M(L, x_i) < k$, so gibt es ein
    $z = r_1 + r_2\sqrt{m} + r_3\sqrt{n} + r_4\sqrt{mn} \in L$ mit den
    Eigenschaften
    \begin{enumerate}
    \item[(a)] $z \equiv x_i \bmod {R}$ für ein
      $i \in \{1, \dots, t\}$;
    \item[(b)] $N_{L / \Q}(z) < k$;
    \item[(c)] $|r_i| < \mu_i$ für $1 \le i \le 4$ und mit
      $$ \mu_1 = \lambda, \quad
         \mu_2 = \frac{\lambda}{\sqrt{|m|}}, \quad 
         \mu_3 = \frac{\lambda}{\sqrt{|n|}}, \quad
         \mu_4 = \frac{\lambda}{\sqrt{|mn|}}, $$  
      wobei
      $$ \lambda = \frac{\sqrt[4]{k}}2
             \Big(\sqrt{|u|} + \frac1{\sqrt{|u|}} \Big) $$
      gesetzt wurde.
    \end{enumerate}}
\end{quote}

Es ist eine merkwürdige Tatsache, dass man (1.11) auch einsetzen kann,
um $M(K, x)$ für ein reellquadratisches abzuschätzen: ist nämlich $e >
1$ die FE von $K$ und $N_{K/\Q}(e) = -1$, so gibt es einen
imaginären, bizyklischen Zahlkörper $L = K(\sqrt{m})$, $m < 0$, indem
$-e = u^2$ zum Quadrat wird für eine Einheit $u$ in $S$ (dem Ring der
ganzen Zahlen in $L$; s. hierzu § 6 und die dort angegebene
Literatur). Anstatt dann $N_{K/\Q}(x-y) < k$ für ein $x \in K$
mit $xe \equiv x \bmod {R}$ und $y \in R$ zu
untersuchen, betrachtet man $N_{L / \Q}(x-y) < k^2$ und
$N_{L/\Q}(xu-y) < k^2$. Wir werden in § 2 ein Beispiel geben,
das zeigt, dass man mit dieser Idee eine Menge Arbeit sparen kann.

Wir kommen damit zum allgemeinen Fall:

\begin{quote}
  {\bf (1.12)} {\em Sei $K$ ein Zahlkörper mit
    $(K:\Q)=n=r+2s$, $\{u_1,\dots,u_l\}$ ein System von
    $l=r+s-1$ unabhängigen Einheiten, $\{\alpha_1,\dots,\alpha_n\}$
    eine $\Q$-Basis von $K$ und $x_1,\dots,x_t$ Elemente von
    $K$, die von den $u_i$ irgendwie permutiert werden. Ist dann
    $M(K,x_i)<k$, dann gibt es ein $z\in K$ mit $z = \sum r_i
    \alpha_i$ und den Eigenschaften
    \begin{enumerate}
    \item[(a)] $z \equiv x_i \bmod {R}$ für ein
      $i \in \{1, \dots, t\}$;
    \item[(b)] $|N_{K/\Q}(z)| < k$;
    \item[(c)] $|r_i| < \mu_i$ für $i=1,\dots,n$, wo die $\mu_i$
      positive reelle Zahlen sind, die nur von $k$, der Wahl der
      $\alpha_i$ und der $u_i$, nicht aber von den Punkten $x_i$
      abhängen.
    \end{enumerate} }
\end{quote}

Bevor wir dies beweisen, führen wir noch einige Bezeichnungen
ein. Dazu denken wir uns die $n=r+2s$ Einbettungen von $K$ in
$\C$ folgendermaßen angeordnet:
$$ \tau_1, \dots, \tau_r, \tau_{r+1}, \ov{\tau}_{r+1}, \dots,
   \tau_{r+s},\ov{\tau}_{r+s}, $$
wo $\tau_1, \dots, \tau_r$ die reellen, $\tau_{r+1}$,
$\ov{\tau}_{r+1}, \dots, \tau_{r+s}, \ov{\tau}_{r+s}$ die Paare
konjugiert-komplexer Einbettungen bezeichnen. Ist $|\cdot|$ der
gewöhnliche Absolutbetrag auf $\R$ (bzw. $\C$), so sind durch
$$ |x_i| =
\begin{cases}
|\tau_i(x)| & \text{für } 1 \le i \le r \\
|\tau_i(x)|^2 & \text{für } r+1 \le i \le r+s
\end{cases} $$
alle $r+s$ verschiedenen archimedischen Bewertungen\label{dpBew} von $K$
definiert. Damit ist $N_{K/\Q}(x) = |x_1| \cdots |x_{r+s}|$.
Wir behaupten nun

\begin{quote}
  {\bf (1.13)} {\em Es seien $l=r+s-1$ unabhängige Einheiten
    $u_1, \dots, u_l$ gegeben, und es sei
    $\kappa_i = \big| \log |u_1|_1 + \dots + \log |u_l|_i \big|$ für
    $i=1,\dots,l$; weiter sei $k_i = \exp(\kappa_i)$. Sind dann
    $c_i$ ($1 \le i \le l$) irgendwelche positiven reellen Zahlen,
    dann gibt es zu jedem $z_1 \in K\setminus\{0\}$ eine Einheit
    $u \in \cO_K^\times$ mit $c_i \le |z_1 u|_i < c_i k_i$
    für $i=1, \dots, l$. }
\end{quote}

\begin{proof}
  Wir definieren eine Abbildung $A : K\setminus\{0\} \to \R^l$ durch
  \[ A(x) = \begin{pmatrix} \log |x|_1 \\ \vdots \\ \log |x|_l \end{pmatrix}; \]
  $A$ hängt natürlich davon ab, welche der $r+s$ Bewertungen von $K$
  wir hier auswählen. Wir behaupten nun, dass die $l$ Vektoren $v_1 =
  A(u_1), \dots, v_l = A(u_l)$ im $\R^l$ unabhängig sind. Dazu nehmen
  wir an, es sei $\sum a_i v_i = 0$ für gewisse $a_i \in \Z$. Indem
  wir diese Gleichung koordinatenweise betrachten, folgt
  \[
  \sum_{i=1}^{l} a_i \log |u_i|_j = 0 \quad \text{für alle } j
  \text{  mit } 1 \le j \le l; \]
  addiert man diese $l$ Gleichungen auf und beachtet
  \[
  \sum_{i=1}^{l} \log |u_i|_j = 0 \quad \text{(wegen }
  |N_{K/\Q}(u_i)| = 1\text{)},
  \]
  so folgt
  \[ \sum_{i=1}^{l} a_i \log |u_i|_j = 0 \quad \text{auch für } j=r+s.
    \text{ Dann ist aber } \prod_{i=1}^{l} u_i^{a_i} = 1. \]
  und da die $u_i$ unabhängig sind, folgt $a_i=0$ für alle $i$
  ($1 \le i \le l$). Also sind die Vektoren $v_1, \dots, v_l$ in der Tat
  linear unabhängig. Nun können wir einen beliebigen Vektor im $\R^l$
  durch Verschieben um geeignete ganzzahlige Vielfache der $v_i$ in
  den Fundamentalbereich
  \[ F = \left\{ \sum b_i v_i : 0 \le |b_i| \le \frac{1}{2} \right\} \]
  bringen. Für ein $w = (w_1, \dots, w_l)^t \in F$ gilt nun
  \[  w_i = \sum_{j=1}^{l} b_j \cdot \log |u_j|_i = 0, \]
  also $0 \le |w_i| \le \frac12\kappa_i$.   Indem wir
  \[ w = \begin{pmatrix} \log |z|_1 - \log c_1 - \frac12 \kappa_1 \\
    \vdots \\
    \log |z|_l - \log c_l - \frac12 \kappa_l
  \end{pmatrix}
  \]
  setzen und $a_i \in \Z$ so finden, dass $w - \sum a_i v_i \in F$
  wird, erhalten wir die Existenz einer Einheit $u = \prod u_i^{a_i}$
  mit
  $$ \begin{array}{rccl}
    \log c_i & \le & \log |z_1 u|_i & < \log c_1 + \kappa_1 \\
                   &     & \ldots & \\
    \log c_l & \le & \log |z_1 u|_l & < \log c_l + \kappa_l 
    \end{array} $$
  Daraus folgt aber sofort unsere Behauptung.
\end{proof}

\begin{proof}[Beweis von (1.12)]
  Wegen $M(K,x_i) < k$ existiert ein $y \in R$ mit
  $|N_{K/\Q}(x_i - y)| < k$. Wir setzen $z_i = x_i - y$ und
  finden mit (1.13) eine Einheit $u$ mit
  $$ \sqrt[n]{k} \le |z_1 u|_i < \sqrt[n]{k} \cdot \kappa_i $$
  für $i=1,\dots,l$. Mit $z = z_1 u$ sind a)
  und b) offenbar erfüllt, und um auch c) nachzuweisen, beachten wir
  \[
  |z|_{l+1} = |N_{K/\Q}(z)| / \prod_{j=1}^{l} |z|_j < \sqrt[n]{k},
  \]
  mit $x_{l+1} = 1$ also $|z|_{i} < \kappa_i \cdot \sqrt[n]{k}$ für
  $i=1, \dots, l+1$.

  Sei nun $\{\beta_1, \dots, \beta_n\}$ die zu $\{\alpha_1, \dots, \alpha_n\}$
  duale Basis. Wie schon im kubischen Fall folgt nun auch   hier
  \[
  |r_j| = |T_{K/\Q}(z \beta_j)| \le \sum_{\tau} |z^\tau| \cdot
  |\beta_j^\tau| < \sqrt[n]{k} \cdot \Bigl( \sum_{\tau} |\lambda_\tau|
  \cdot |\beta_j^\tau| \Bigr) =: \mu_j. 
  \]
  Hierbei stimmen die $\lambda_\tau$ mit den $x_i$ überein, falls
  $\tau$ reell ist ($i=1, \dots, r$), und nach geeigneter
  Umnumerierung mit $\sqrt{\kappa_i}$, falls $\tau$ nicht reell ist
  ($i=r+1, \dots, r+2s$). Diese lästige Umbenennung rührt von der
  Definition der $|\cdot|_i$ vor (1.13) her, wo bei nichtreellen
  Einbettungen das Quadrat des gewöhnlichen Betrages steht. Ich habe
  aber bisher nicht gesehen, wie sich dieser Mangel beheben lässt.
\end{proof}

Zweifellos lassen sich die Abschätzungen, die wir im Verlauf des
Beweises von (1.12) gemacht haben, noch verbessern. Wenn man daher
(1.12) zur Bestimmung von $M(K)$ benutzen will, wird man sich hierüber
Gedanken machen müssen.

\subsection*{Tabellen}

Nachstehend geben wir einige Tabellen, die zeigen sollen, wie man die
Kriterien 1.5., 1.7. und 1.8. in reellquadratischen Zahlkörpern
anwenden kann. Aus diesen Tabellen kann man bereits ersehen, dass die
Ringe $D(m)$ mit $m \le 101$ höchstens für die Werte
$m = 2, 3, 5, 6, 7, 11, 13, 17, 19, 21, 29, 33, 37, 41, 57, 73$
normeuklidisch sein können; in § 2 und § 3 werden wir zeigen, dass
diese Ringe die einzigen normeuklidischen reellquadratischen Zahlkörper
sind.

\begin{table}[ht!]
\centering
\begin{tabular}{ccccc}
\toprule
$m$ & $f$ & $a$ & $k = M(x)$ & $r$ \\
\midrule
7 & 14 & 9 & $9/14$ & $-5$ \\
10 & 10 & 5 & $3/2$ & $5$ \\
11 & 22 & 3 & $19/22$ & $3$ \\
13 & 13 & 4 & $4/13$ & -- \\
15 & 15 & 6 & $7/5$ & $6, -9$ \\
17 & 17 & 8 & $8/17$ & -- \\
19 & 38 & 7 & $31/38$ & $7$ \\
21 & 7 & 2 & $5/7$ & $2$ \\
22 & 22 & 5 & $27/22$ & $5, -17$ \\
23 & 46 & 31 & $77/46$ & $-15, 31, -61$ \\
29 & 29 & 6 & $23/29$ & $6$ \\
31 & 31 & 14 & $45/31$ & $14, -17$ \\
33 & 11 & 5 & $6/11$ & $5$ \\
35 & 35 & 15 & $17/7$ & $15, -20, 50, -55$ \\
37 & 37 & 10 & $27/37$ & $10$ \\
47 & 94 & 65 & $253/94$ & $-29, 65, -123, 159, -217$ \\
51 & 102 & 19 & $287/102$ & (e=23): $-185, -83, 19, 121, 223$ \\
53 & 53 & 15 & $68/53$ & $15, -38$ \\
57 & 19 & 5 & $14/19$ & $5$ \\
59 & 59 & 7 & $125/59$ & $7, -52, 66, -111$ \\
69 & 23 & 2 & $25/23$ & $2, -21$ \\
77 & 11 & 3 & $19/11$ & $3, -8, 14$ \\
79 & 158 & 111 & $585/158$ & $-47, 111, -205, 269, -363, 427, -521$ \\
83 & 166 & 33 & $631/166$ & $33, -133, 199, -299, 365, -465, 531$ \\
85 & 85 & 19 & $151/85$ & $19, -66, 104$ \\
87 & 58 & 53 & $169/58$ & $-5, 53, -63, 111, -121$ \\
93 & 31 & 18 & $44/31$ & $-13, 18$ \\
101 & 101 & 24 & $125/101$ & $24, -77$ \\
\bottomrule
\end{tabular}
\caption{Anwendung von (1.5) auf quadratische
  Zahlkörper $\Q(\sqrt{m})$}\label{dApp15}
\end{table}

\medskip\noindent

Bem. zu Tab.~\ref{dApp15}: Ist $e^2 \equiv a \bmod {f}$ und $(c)^2 = (f)$
für ein $c \in R$, so setze man $x = e/c$; man erhält damit ein
$x \in K$ mit $M(K,x)=k$.

Im Falle $m=51$ ist noch zu bemerken, dass 121 natürlich Norm aus $S$
ist wegen $121 = N_{K/\Q}(11)$; jedoch ist $11 \not\equiv e
\equiv 23 \bmod {f}$. Diese Schwierigkeit tritt hier auf, weil die
Kongruenz $x^2 \equiv a \equiv 19 \bmod {f}$ vier verschiedene Lösungen
hat, nämlich $x=11$ und $x=23$.

\begin{table}[h]
\centering
\begin{tabular}{ccccc}
\toprule
$m$ & $B'$ & $y$ & $k=M(x)$ & $r$ \\
\midrule
6 & 4 & $1+\sqrt{6}$ & $3/4$ & $1$ \\
10 & 4 & $\sqrt{10}$ & $3/2$ & $2$ \\
14 & 4 & $1+\sqrt{14}$ & $5/4$ & $1, 3$ \\
15 & 4 & $1+\sqrt{15}$ & $3/2$ & $2$ \\
26 & 4 & $\sqrt{26}$ & $5/2$ & $2, 6$ \\
30 & 4 & $\sqrt{30}$ & $3/2$ & $2$ \\
20 & & $1+\sqrt{30}$ & $29/20$ & $-19, -9, 1, 11, 21$ \\
34 & 4 & $1+\sqrt{34}$ & $9/4$ & $1, 3, 5, 7$ \\
38 & 4 & $1+\sqrt{38}$ & $11/4$ & $1, 3, 5, 7, 11$ \\
39 & 4 & $1+\sqrt{39}$ & $5/2$ & $2, 6$ \\
42 & 4 & $1+\sqrt{42}$ & $7/4$ & $1, 3, 5$ \\
55 & 4 & $\sqrt{55}$ & $9/4$ & $1, 3, 5, 7$ \\
58 & 4 & $\sqrt{58}$ & $3/2$ & $2$ \\
62 & 4 & $1+\sqrt{62}$ & $13/4$ & $1, 3, 5, 7, 11, 13$ \\
66 & 4 & $1+\sqrt{66}$ & $15/4$ & $1, 3, 5, 7, 11, 13$ \\
74 & 4 & $\sqrt{74}$ & $5/2$ & $2, 6$ \\
78 & 4 & $\sqrt{78}$ & $7/2$ & $2, 6, 10$ \\
82 & 4 & $\sqrt{82}$ & $9/2$ & $2, 6, 10, 14$ \\
91 & 4 & $1+\sqrt{91}$ & $5/2$ & $2, 6$ \\
95 & 4 & $1+\sqrt{95}$ & $7/2$ & $2, 6, 10$ \\
\bottomrule
\end{tabular}
\caption{Anwendung von (1.7) auf quadratische
  Zahlkörper $\Q(\sqrt{m})$}\label{dTabB}
\end{table}

\medskip\noindent Bem. zu Tab.~\ref{dTabB}: Ist $B'=(c)^2$ für ein $c
\in R$, so ist $M(x) = k$ für $x=y/c$.

\begin{table}[ht!]
\centering
\begin{tabular}{cccc}
\toprule
$m$ & $x$ & $M(x)$ & $u$ \\
\midrule
19 & $(0, 20/57)$ & $170/171$ & $170 + 39\sqrt{19}$ \\
21 & $(2/5, 1/5)$ & $12/25$ & $(5 + \sqrt{21})/2$ \\
22 & $(0, 9/28)$ & $443/392$ & $197 + 42\sqrt{22}$ \\
29 & $(3/5, 1/5)$ & $4/5$ & $(5 + \sqrt{29})/2$ \\
33 & $(1/2, 2/11)$ & $29/44$ & $23 + 4\sqrt{33}$ \\
34 & $(1/2, 6/17)$ & $135/68$ & $35 + 6\sqrt{34}$ \\
    & $(0, 5/12)$ & $137/72$ & $35 + 6\sqrt{34}$ \\
35 & $(0, 3/7)$ & $17/7$ & $6 + \sqrt{35}$ \\
38 & $(0, 5/12)$ & $173/72$ & $37 + 6\sqrt{38}$ \\
39 & $(0, 5/12)$ & $107/48$ & $25 + 4\sqrt{39}$ \\
41 & $(15/32, 5/32)$ & $23/32$ & $32 + 5\sqrt{41}$ \\
42 & $(0, 7/12)$ & $41/24$ & $13 + 2\sqrt{42}$ \\
43 & $(1/118, 1/2)$ & $11829/6962$ & $3482 + 531\sqrt{43}$ \\
    & $(0, 193/387)$ & $5902/3483$ & $3482 + 531\sqrt{43}$ \\
46 & $(1/2, 311/1058)$ & $76877/48668$ & $24335 + 3588\sqrt{46}$ \\
53 & $(1/14, 3/14)$ & $9/7$ & $(7 + \sqrt{53})/2$ \\
55 & $(1/2, 19/44)$ & $351/176$ & $89 + 12\sqrt{55}$ \\
57 & $(0, 15/76)$ & $219/304$ & $151 + 20\sqrt{57}$ \\
61 & $(1/78, 17/78)$ & $41/39$ & $(39 + 5\sqrt{61})/2$ \\
62 & $(1/2, 48/124)$ & $367/124$ & $63 + 8\sqrt{62}$ \\
    & $(0, 7/16)$ & $367/128$ & $63 + 8\sqrt{62}$ \\
65 & $(1/4, 1/4)$ & $1$ & $8 + \sqrt{65}$ \\
66 & $(0, 7/16)$ & $31/128$ & $65 + 8\sqrt{66}$ \\
67 & $(1/2, 7/18)$ & $341/162$ & $48842 + 5967\sqrt{67}$ \\
70 & $(1/2, 6/25)$ & $891/500$ & $251 + 30\sqrt{70}$ \\
71 & $(0, 216/497)$ & $7393/3479$ & $3480 + 413\sqrt{71}$ \\
73 & $(41, 283)/2136$ & $1541/2136$ & $1068 + 125\sqrt{73}$ \\
85 & $(4/9, 2/9)$ & $16/9$ & $(9 + \sqrt{85})/2$ \\
86 & $(0, 133/473)$ & $10030/5203$ & $10405 + 1122\sqrt{86}$ \\
89 & $(3, 159)/1000$ & $1004/1000$ & $500 + 53\sqrt{89}$ \\
94 & $(0, 3661/14194)$ & $4708623/2143294$ & $2143295 + 221064\sqrt{94}$ \\
97 & $(1496, 2587)/11208$ & $3001/2802$ & $5604 + 569\sqrt{97}$ \\
109 & $(17, 82)/261$ & $289/261$ & $261 + 25\sqrt{109}$ \\
113 & $(3, 219)/1352$ & $967/776$ & $776 + 73\sqrt{113}$ \\
137 & $(3, 447)/3488$ & $2177/1744$ & $1744 + 149\sqrt{137}$ \\
\bottomrule
\end{tabular}
\caption{Anwendung von (1.8) auf quadratische Zahlkörper
  $\Q(\sqrt{m})$ ($t=1$):}
\end{table}

\begin{table}[h]
\centering
\begin{tabular}{cccc}
\toprule
$m$ & $x_1, x_2$ & $M(x_1)$ \\
\midrule
19 & $(0, 173/494)$ & \multirow{2}{*}{$10579/12844$} \\
    & $(1/2, -115/247)$ & \\
22 & $(1182, 496441/155236)$ & \multirow{2}{*}{$175903/155236$} \\
    & $(-1182, 496441/155236)$ & \\
58 & $(1/2, 197/754)$ & \multirow{2}{*}{$27477/19604$} \\
    & $(1/2, 276/754)$ & \\
61 & $(0, 66/305)$ & \multirow{2}{*}{$1611/1525$} \\
    & $(0, 67/305)$ & \\
74 & $(1/2, 19/43)$ & \multirow{2}{*}{$16255/7396$} \\
    & $(1/86, 1/2)$ & \\
93 & $(0, 44/279)$ & \multirow{2}{*}{$2198/1674$} \\
    & $(1/2, 119/558)$ & \\
\bottomrule
\end{tabular}
\caption{$t=2$}
\end{table}

\medskip

Für weitere Beispiele mit $t=2$ s. § 2.

\subsection*{Hinreichende Bedingung für den Euklidischen Algorithmus}

Nachdem wir uns bisher mit notwendigen Kriterien beschäftigt haben,
wenden wir uns nun solchen zu, die die Existenz des EA garantieren.

1974 hat Lenstra\index[dN]{Lenstra} eine Methode vorgestellt, mit der
inzwischen (s. z.B. Lenstra \cite{Len77a},
Leutbecher\index[dN]{Leutbecher} u. Martinet\index[dN]{Martinet}
\cite{LM81a,LM81b}, Leutbecher \cite{Leu85,Leu86}, Leutbecher und
Niklasch\index[dN]{Niklasch} \cite{LN87} ca. 400 normeuklidische
Zahlkörper gefunden wurden. Eine Modifikation dieses Kriteriums wird
es uns erlauben, auf ähnliche Art und Weise auch $k$-stufig
normeuklidische Zahlringe zu erhalten (diese werden wir in Zukunft
einfach $k$-euklidisch nennen).

Lenstras Methode basiert auf einer Idee von Hurwitz\index[dN]{Hurwitz}
\cite{Hur19}: dieser hat gezeigt, dass es in einem Zahlkörper $K$ eine
nur von $K$ abhängige natürliche Zahl $m \ge 2$ gibt, sodass gilt:

\begin{quote}
  {\em Für alle $x \in K$ gibt es ein $y \in R$ und ein $j \in \N$,
    $1 \le j < m$ mit $|N_{K/\Q}(j(x-y))| < 1$.}
  \end{quote}

Man beachte, dass $K$ genau dann normeuklidisch ist, wenn wir $m=2$
wählen würden. Lenstra hat dann bemerkt, dass es im Beweis dieser
Aussage nicht darauf ankommt, dass die obigen $j$ natürliche Zahlen
sind, sondern dass etwas allgemeiner gilt:

\begin{quote}
  {\bf (1.14)} {\em Sei $K$ ein Zahlkörper, $(K:\Q) = n = r+2s$, und
    $d = |\disc K|$. Ist dann $m \in \N$,
    $$ m > M_{r,s}\label{dpMrs} =
          \frac1d \cdot \frac{n!}{n^n} \cdot \Big(\frac4\pi\Big)^s $$
    und sind $\theta_1, \dots, \theta_m$
    irgendwelche paarweise verschiedenen Elemente von $K$, dann gibt
    es für alle $x \in K$ ein $y \in R$ mit $|N_{K/\Q}((\theta_i -
    \theta_j)x - y)| < 1$ für gewisse
    $i, j \in \{1, \dots, m\}$, $1 \le i < j \le m$.}
\end{quote}

Lenstras Idee war nun die folgende: wenn wir die $\theta_1, \dots,
\theta_m$ so wählen können, dass die Differenzen $\theta_i - \theta_j$
für $i < j$ lauter Einheiten in $R$ sind, dann ist
$c = \frac{y}{\theta_i - \theta_j} \in R$, und (1.14) besagt dann
gerade, dass es für alle $x \in K$ ein $c \in R$ mit $|N_{K/\Q}(x-c)| < 1$
gibt: $R$ wäre damit normeuklidisch.

Der Beweis von (1.14), den wir nun vorstellen wollen, stammt von
Lenstra\index[dN]{Lenstra} \cite{Len74} und hat große Ähnlichkeit mit
dem klassischen Beweis für die Endlichkeit der Klassenzahl (1977 hat
Lenstra\index[dN]{Lenstra} \cite{Len77a} einen weiteren Beweis in der
Sprache der Packungstheorie gegeben). Einige der untenstehenden
Aussagen werden wir daher als bekannt voraussetzen; für Beweise
verweisen wir auf MarcusMarcus\index[dN]{Marcus} \cite{Mar77} oder
Cassels\index[dN]{Cassels} \cite{Cas86}.

\begin{proof}[Beweis von (1.14)]
  Seien nun die Einbettungen $\tau_1, \dots, \tau_n$ von $K$ in $\C$
  wie vor (1.12) angeordnet; wir definieren dann eine Einbettung von
  $K$ in den $\R^n$ durch $\varphi : K \to \R^n$ durch
  $$ \alpha \mapsto \varphi(\alpha) =  
   (\tau_1(\alpha), \dots, \tau_r(\alpha), \Re
   \tau_{r+1}(\alpha), \Im \tau_{r+1}(\alpha), \dots,
   \Re \tau_{r+s}(\alpha), \Im \tau_{r+s}(\alpha)), $$   

  wobei $\Re$ bzw. $\Im$ den Real- bzw. Imaginärteil einer komplexen
  Zahl bezeichnen. Auf $\R^n$ definieren wir nun eine Norm für
  $x = (x_1, \dots, x_n) \in \R^n$ durch
\[
N(x) = x_1^2 \cdots x_r^2 (x_{r+1}^2 + x_{r+2}^2) \cdots (x_{n-1}^2 + x_n^2).
\]
Damit wird $N_{K/\Q}(\alpha) = N(\varphi(\alpha))$ für alle
$\alpha \in K$. Nun setzen wir
\[
A = \left\{ x \in \R^n : |x_1| + \cdots + |x_r| +
2 \sqrt{x_{r+1}^2 + x_{r+2}^2} + \cdots + 2 \sqrt{x_{n-1}^2 + x_n^2} < n
\right\}.
\]
Bekanntlich ist $A$ beschränkt und konvex, außerdem liefert die
Ungleichung zwischen geometrischem und arithmetischem Mittel die
Beziehung $|N(\alpha)| < 1$ für alle $\alpha \in A$. Mit $U =
\frac{1}{2}A$ ist dann $u - v \in A$ für alle $u,v \in U$, sodass für
solche $u,v$ $|N(u-v)| < 1$ gilt. Schließlich ist noch
$$ \operatorname{vol}(A) = 2^r (\frac{\pi}2)^s \frac{n^n}{n!}
  \quad \text{und} \quad \
  \operatorname{vol}(U) = 2^n \cdot \operatorname{vol}(A), $$
während die Fundamentalmasche $F$ des Gitters $\Gamma_R := \varphi(\cO_K)$
im $\R^n$ das Volumen $\operatorname{vol}(F) = 2^{-s} \sqrt{d}$ besitzt.

Im folgenden wollen wir die $\alpha \in K$ mit ihren Bildern
$\varphi(\alpha) \in \R^n$ identifizieren; da $\varphi$
injektiv ist und $N_{K/\Q}(\alpha) = N(\varphi(\alpha))$ gilt,
sollten deswegen keine Probleme auftreten.

Seien nun $\vartheta_1, \dots, \vartheta_m \in K$ paarweise
verschieden; dann können wir zu jedem Punkt $z$ aus einer der Mengen
$\vartheta_i x + U$ ($i=1,\dots,m$) ein $y \in R$ finden,
sodass $z_0 - z = y$ innerhalb der Fundamentalmasche $F$ des Gitters
$R$ liegt. Nun ist aber das Gesamtvolumen der Mengen $\vartheta_i x +
U$ gleich $m \cdot \operatorname{vol}(U) > M_{r,s} \cdot
\operatorname{vol}(U) = 2^{-s} \sqrt{d} =
\operatorname{vol}(F)$. Folglich gibt es Indizes $i,j$ und Zahlen
$y_1, y_2 \in R$, sodass
\[
(\vartheta_i x + U - y_1) \cap (\vartheta_j x + U - y_2) \neq \emptyset
\]
ist und die Paare $(i,y_1)$ und $(j,y_2)$ verschieden sind. Damit
haben wir $u,v \in U$ gefunden mit $\vartheta_i x + u - y_1 =
\vartheta_j x + v - y_2$. Wäre $i=j$, so folgte $u - v = y_1 - y_2 \in
R$; wegen $|N(u-v)| < 1$ muss dann $N(y_1 - y_2) = 0$ sein,
und dies impliziert $y_1 = y_2$. Dies widerspricht aber der Tatsache,
dass die Paare $(i,y_1)$ und $(j,y_2)$ verschieden sind. Daher ist $i
\neq j$ und $(\vartheta_i - \vartheta_j)x - y = v - u$ für $y := y_1 -
y_2$. Normbildung liefert nun
\[
|N((\vartheta_i - \vartheta_j)x - y)| = |N(u - v)| < 1
\]
wie behauptet.
\end{proof}

Wir haben bereits bemerkt, dass aus (1.14) folgt:

\begin{quote}
  {\bf (1.15)} {\em Gibt es Zahlen $\vartheta_1, \dots, \vartheta_m
    \in K$, deren Differenzen Einheiten in $R$ sind und ist $m >
    M_{r,s}$, dann ist $R$ normeuklidisch. }
\end{quote}

Eine Folge $\vartheta_1, \dots, \vartheta_m$ von Elementen aus $K$,
deren Differenzen lauter Einheiten in $R$ sind, nennen wir eine
\textbf{Ausnahmefolge}\index[dS]{Ausnahmefolge} (``exceptional sequence'').
Da es nur auf die
Differenzen der $\vartheta_i$ ankommt, dürfen wir oBdA $\vartheta_1=0$
annehmen (damit sind alle $\vartheta_i \in R$). Weil weiter mit
$\vartheta_i - \vartheta_j$ auch $(\vartheta_i -
\vartheta_j)/\vartheta_2$ Einheit ist, ist mit $0, \vartheta_2, \dots,
\vartheta_m$ auch $0, 1, \vartheta_3/\vartheta_2, \dots, \vartheta_m /
\vartheta_2$ eine Ausnahmefolge. Also dürfen wir auch noch
$\vartheta_2=1$ annehmen. Ein Beispiel für eine Ausnahmefolge der
Länge $p$ ist dann in $K = \Q(\zeta)$, wo $\zeta = \zeta_p$
eine primitive $p$-te Einheitswurzel bezeichnet:
\[
0, 1, 1+\zeta, 1+\zeta+\zeta^2, \dots, 1+\zeta + \zeta^2 + \dots + \zeta^{p-2};
\]
dies zeigt, dass der Ring $\Z[\zeta]$ normeuklidisch ist, falls $p >
M_{r,s}$ gilt, und dies ist für $p = 3, 5, 7$ der Fall.

Um ein zu (1.15) analoges Ergebnis auch für $k$-euklidische Ringe
formulieren zu können, definieren wir: eine Folge $\vartheta_1, \dots,
\vartheta_m$ heiße \textbf{$k$-Folge},\index[dS]{$k$-Folge}
wenn die beiden folgenden Bedingungen erfüllt sind:
\begin{enumerate}
\item[] (F--1) es ist $\vartheta_i - \vartheta_j \in E_k \setminus \{0\}$
  für alle $1 \le i < j \le m$;  
\item[] (F--2)\label{dpF2} jeder Teiler des Ideals $(\vartheta_i - \vartheta_j)$
  ist Hauptideal.  
\end{enumerate}
Damit sind 1-Folgen und Ausnahmefolgen dasselbe, weil dann die
Bedingung (F--2) automatisch erfüllt ist. Wir behaupten nun

\begin{quote}
  {\bf (1.16)} {\em Gibt es in $R$ eine $k$-Folge der Länge
    $m > M_{r,s}$, dann ist $R$ $k$-euklidisch.}
\end{quote}

\begin{proof}
  Sei $x \in K$ gegeben; mit (1.14) finden wir ein $y \in R$ mit
  $|N((\vartheta_i - \vartheta_j)x - y)| < 1$, wobei $\vartheta_1,
  \dots, \vartheta_m$ eine $k$-Folge der Länge $m$ ist. Wegen (F--2)
  dürfen wir $y/(\vartheta_i - \vartheta_j) = a/b$ schreiben für
  gewisse $a, b \in R$ mit $(a, b)=1$ und $b|(\vartheta_i -
  \vartheta_j)$. Wegen (0.12) ist also $b \in E_k$, daher $a/b$ nach
  (0.11) ein Kettenbruch der Länge $\le k$ mit Nenner $ub$ für eine
  Einheit $u \in R^\times$, und schließlich
  \[
  \Big|N_{K/\Q}\Big(x - \frac{a}{b}\Big)\Big| <
  \frac{1}{N_{K/\Q}(\vartheta_i - \vartheta_j)} <
  \frac{1}{N_{K/\Q}(b)}.
  \]
  mit (1.14). Nun zeigt (0.9), dass $R$ wirklich $k$-euklidisch ist.
\end{proof}

Man beachte, dass eine Folge $\vartheta_1, \dots, \vartheta_m$ mit der
Eigenschaft $\vartheta_i - \vartheta_j \in E_k$ schon eine $k$-Folge
ist: wegen $E_k \subset E_k'$ ist nämlich (F--1) erfüllt, und (0.5.(i))
zeigt (F--2). Der Grund für die Einführung der Mengen $E_k'$ liegt
allein darin, dass diese Mengen i.A. echt größer als die entsprechenden
$E_k$ sind, was das Auffinden von (möglichst langen) $k$-Folgen
natürlich erleichtert.

Ist $I$ ein ganzes Ideal in $R$ und gibt es ein primes Restsystem mod
$I$, welches ganz aus Einheiten von $R$ besteht, so sagen wir, $I$
besitze ein primes Einheitenrestsystem (kurz: ein PERS).\label{dpPERS}
Insbesondere hat $I$ ein PERS, wenn eine Einheit $u \in R$ existiert,
die Primitivwurzel mod $I$ ist. Dagegen gibt es selbst dann, wenn $R$
Einheitenrang 1 hat, Ideale, die zwar ein PERS, aber keine Einheit als
Primitivwurzel besitzen.

Wir geben einige einfache Beispiele:
\begin{enumerate}
\item[1.] $K = \Q(\sqrt{5})$: hier ist $\omega = (1+\sqrt{5})/2 \in
  R$, und $0,1, \omega, 1+\omega$ ist eine 1-Folge, da $1, \omega,
  1+\omega, \omega-1$ Einheiten sind. Dagegen ist $0, 1, -1+\omega,
  \omega, 1+\omega$ keine 1-Folge, weil z.B. die Differenz $-2 =
  (\omega-1)-(\omega+1)$ keine Einheit ist. Allerdings haben wir hier
  eine 2-Folge, weil die Ideale $I_1=(2)$ und $I_2=(2+\omega)$ ein
  PERS besitzen: es ist $\{1, \omega, 1+\omega\}$ ein primes
  Restsystem mod $I_1$ und $\{1, \omega, \omega^2, \omega^3,
  \omega^4\}$ eines mod $I_2$ (in beiden Fällen ist $\omega$
  Primitivwurzel).
\item[2.] $K=\Q(\sqrt{14})$: hier ist $0,1$ eine 1-Folge der Länge 2,
  und man sieht leicht ein, dass es keine längeren gibt; sind nämlich
  $u$ und $u-1$ Einheiten in $R$, so sind sie (bis auf ein etwaiges
  Vorzeichen) Potenzen der FE $\varepsilon = 15+4\sqrt{14}$. Nun ist
  aber $\varepsilon \equiv 1 \bmod {2}$, also müssten auch $u$ und
  $u-1 \equiv 1 \bmod {2}$ sein, was offenbar zum Widerspruch
  führt. Wir behaupten nun, dass $0, 1, 4+\sqrt{14}, 5+\sqrt{14}$ eine
  2-Folge ist. Dazu müssen wir zeigen, dass die Ideale
  $I_1=(3+\sqrt{14})$, $I_2=(4+\sqrt{14})$ und $I_3=(5+\sqrt{14})$ ein
  PERS besitzen. In § 0 (S.~\pageref{p11}) haben wir bereits gesehen, dass
  $3+\sqrt{14} \in E_2'$ ist, und dass $\varepsilon$ hier
  Primitivwurzel ist. Wegen $\|I_2\| = 2$ ist $\{1\}$ ein PERS von
  $I_2$; schließlich rechnet man noch nach, dass $\varepsilon$ auch
  mod $I_3$ eine Primitivwurzel ist. Da hier $M_{r,s} = 3,74\dots$ ist
  und wir eine 2-Folge der Länge 4 gefunden haben, können wir mit
  (1.16) folgern, dass $\Q(\sqrt{14})$ 2-euklidisch ist.
\end{enumerate}

Im allgemeinen ist es recht mühsam, in einem gegebenen Zahlkörper $K$
eine 1-Folge (oder gar eine $k$-Folge) genügend großer Länge zu
finden; bevor man mit der Suche nach $k$-Folgen beginnt, wird man
daher wissen wollen, ob sich diese Mühe überhaupt lohnt, d.h. ob es in
$R$ $k$-Folgen der gewünschten Länge geben kann. Um diese Frage
beantworten zu können, setzen wir\label{dLenmu}
\begin{align*}
  \mu_k & = \sup \{ m \in \N : \text{es gibt eine $k$-Folge der Länge } m
       \text{ in } R \} \qquad \text{und} \\
  \lambda_k & = \inf \{ \|I\| : \text{es gibt eine prime Restklasse mod } I, \\
       & \qquad \text{ die keinen Vertreter in } E_0' = \{0\} \text{ hat} \}.
\end{align*}
Insbesondere ist
$\lambda_1 = \min \{ \|I\| \ge 2: \ I \text{ ist ganzes Ideal in } R \}$,
da die prime Restklasse $1 \bmod I$ keinen Vertreter in $E_1 = \{0\}$
hat (insbesondere ist $\lambda_1 \le 2$, da man $I=(2)$ wählen kann),
und $\lambda_2 = \inf \{ \|I\| \ge 2: \ I \text{ besitzt kein PERS in } R \}$.

Mit diesen Bezeichnungen gilt dann

\begin{quote}
  {\bf (1.17)} {\em Ist $\lambda_k$ endlich, dann gilt
    $2 \le \mu_1 \le \mu_2 \le \dots \le \mu_k \le \lambda_k$.}
\end{quote}

\begin{proof}
  Da $0,1$ eine 1-Folge ist, gilt $\mu_1 \ge 2$; weiter ist jede
  $k$-Folge erst recht eine $(k+1)$-Folge, also $\mu_k \le \mu_{k+1}$
  für alle $k \in \N$. Also müssen wir nur noch $\mu_k \le \lambda_k$
  zeigen. Sei dazu $I$ ein ganzes Ideal in $R$, das eine prime
  Restklasse mod $I$ besitzt, welche keinen Vertreter in $E_k$ hat
  (ein solches existiert, da $\lambda_k$ endlich ist). Ist dann
  $\vartheta_1, \dots, \vartheta_m$ eine $k$-Folge, dann muss
  $\vartheta_i - \vartheta_j \not\equiv 0 \bmod {I}$ sein für alle $i
  \ne j$: ansonsten wäre nämlich $I$ als Teiler des Ideals
  $(\vartheta_i - \vartheta_j)$ ein Hauptideal wegen (F--2), also
  $I=(b)$ für ein $b \in R$. Wegen $b \mid (\vartheta_i - \vartheta_j)$
  und (0.12) ist nun $b \in E_k$, während nach Voraussetzung das Ideal
  $I=(b)$ eine prime Restklasse besitzt, welche keinen Vertreter
  in $E_k$ hat: Widerspruch! Also sind $\vartheta_1, \dots, \vartheta_m$
  paarweise inkongruent mod $I$ und $m \le \|I\|$. Dies aber impliziert
  $\mu_k \le \lambda_k$.
\end{proof}

Wir wollen die oben gegebenen Beispiele 1. und 2. von Seite 36 noch
einmal im Lichte von (1.17) betrachten: in $\Q(\sqrt{5})$ ist
$\lambda_1=4$, weil $I=(2)$ das Ideal mit minimaler Norm $>1$ ist. Für
jede 1-Folge $\vartheta_1, \dots, \vartheta_4$ sind also die
$\vartheta_i$ paarweise inkongruent mod $2$, wie wir im Beweis von
(1.17) gesehen haben, und für $0,1,\omega,1+\omega$ ist dies offenbar
auch der Fall. Weiter können wir hier $\mu_1 = \lambda_1$ schließen.

In $K=\Q(\sqrt{14})$ dagegen ist $I_2=(4+\sqrt{14})$ ein Ideal
der Norm 2, sodass wir hier $\mu_1 = \lambda_1 = 2$ haben. Weiter ist
$\lambda_2 = 4$, da das Ideal $I=(2)$ kein PERS hat. Die angegebene
2-Folge der Länge 4 zeigt schließlich $\mu_2 = \lambda_2 = 4$.

Während es mir für $k \ge 3$ nicht gelungen ist zu zeigen, dass die
$\lambda_k$ endlich sind, gilt für $k=2$:

\begin{quote}
  {\bf (1.18)}{\em Sei $K$ ein Zahlkörper und $q$ eine rationale
    Primzahl. Dann gibt es ein $a \in \N$, sodass $I = (q^a)R$ kein
    PERS besitzt. Insbesondere ist $\lambda_2 \le (q^a)^n$, wo
    $n=(K:\Q)$ der Körpergrad ist.}
\end{quote}

\begin{proof}
  Sei $N$ der normale Abschluss von $K/\Q$. Besitzt $I = (q^a)$ ein
  PERS, dann gibt es zu jedem $r \in \Z$ mit $q \nmid r$ eine Einheit
  $u \in R^\times$ mit $r \equiv u \bmod {I}$. Durchläuft $\sigma$ die
  Einbettungen von $K$ in $\C$, so folgt für jedes $\sigma$ die
  Kongruenz $u^\sigma \equiv r \bmod {I}$ in $\cO_N$ (= der
  Ring der ganzen Zahlen in $N$). Da $r^\sigma = r$ und $1^\sigma = 1$
  ist. Multiplikation über alle $\sigma$ gibt dann $1 = N_{K/\Q}(u) =
  \prod u^\sigma \equiv r^n \bmod {I}$ in $\cO_N$. Wegen $I
  \cap \Z = q^a \Z$ gilt die Kongruenz $1 \equiv r^n \bmod {q^a}$ auch
  in $\Z$. Indem wir nun $r > 1$ und $a$ so groß wählen, dass $q^a >
  r^n + 1$ wird, erhalten wir den gewünschten Widerspruch.
\end{proof}

Wir bemerken noch:

\begin{enumerate}
\item[1.] Ist $\vartheta_1, \dots, \vartheta_m$ eine 1-Folge in $K$
  und $K \subset L$, dann ist $\vartheta_1, \dots, \vartheta_m$ auch
  eine 1-Folge in $L$; für $k$-Folgen mit $k \ge 2$ gilt dies
  i.A. nicht mehr, weil die Norm der Ideale $(\vartheta_i - \vartheta_j)$
  von $L$ abhängt, falls $\vartheta_i - \vartheta_j$ keine Einheit ist.
\item[2.] Wie Lenstra\index[dN]{Lenstra}  \cite[S. 239]{Len77a} im Falle $k=1$ bereits bemerkt
  hat, kann man durch eine kleine Modifikation des Beweises von
  1.13. eine obere Schranke für die euklidischen Minima $M^k(K)$
  erhalten; es gilt nämlich unter der Voraussetzung von (1.16)
  $M^k(K) \le M_{r,s}/\mu_k$.
\item[3.] Auch die Definition (1.17) von Lenstra \cite[S. 241]{Len77a} lässt
  sich auf $k$-Folgen mit $k > 1$ verallgemeinern: man betrachtet hier
  Folgen $\vartheta_1, \dots, \vartheta_m$, sodass unter $l+1$
  Folgengliedern mindestens zwei sind, deren Differenz in $E_k$ liegt
  (hierbei ist $l$ eine natürliche Zahl; für $l=1$ erhält man die
  gewöhnlichen $k$-Folgen zurück). Bezeichnet man mit $\mu_{k,l}$ die
  maximale Länge solcher verallgemeinerter $k$-Folgen und ist
  $\mu_{k,l} > l \cdot M_{r,s}$, so kann man ähnlich wie in (1.16)
  schließen, dass $K$ $k$-euklidisch ist.
\item[4.] Unter der Annahme der erweiterten Riemannschen Vermutung hat
  Lenstra gezeigt \cite[S. 242f]{Len77a}, dass man mit (1.15) nur endlich
  viele normeuklidische Zahlkörper finden kann; dabei wurde wesentlich
  die Ungleichung $\lambda_1 \le 2^n$, $n=(K:\Q)$, benutzt, die man
  direkt aus der Definition von $\lambda_1$ mit $I = (2)$ erhält. Die
  Abschätzung (1.18) für $\lambda_2$ lässt einen ähnlichen Schluss für
  $k$-euklidische Ringe mit $k \ge 2$ nicht zu.
\end{enumerate}

\subsection*{Die Gaußsche Generalmensur}
Wenn wir zeigen wollen, dass ein Zahlkörper $K$ normeuklidisch ist,
müssen wir zu jedem $x \in K$ ein $y \in R$ mit
$|N_{K/\Q}(x-y)|<1$ finden. Angenommen, wir kennen eine
Funktion $\cM_K : K \to \R$ mit der Eigenschaft
$|N_{K/\Q}(x)| \le \cM_K(x)$ für alle $x \in K$, so
genügt es offenbar, zu jedem $x \in K$ ein $y \in R$ mit
$\cM_K(x-y)<1$ zu finden. Wegen $N_{K/\Q}(x) = \prod
\sigma(x)$, wo das Produkt über alle Einbettungen $\sigma$ von $K$ in
$\C$ läuft, kommt z.B.
\[
\cM_K(x) = \cM_{K/\Q}(x) = \frac{1}{n} \Bigl(
\sum |\sigma(x)|^2 \Bigr)
\]
als solche Funktion in Frage; in der Tat gilt\label{dGMen}

\begin{quote}
  {\bf (1.19)} {\em Sei $K$ ein Zahlkörper mit $(K:\Q) = n$; dann gilt
    für alle $x \in K$}
  \[ |N_{K/\Q}(x)| \le \cM_K(x)^{n/2}. \]
\end{quote}

\begin{proof}
  Aus der Ungleichung zwischen geometrischem und arithmetischem Mittel
  erhält man
  \[ |N_{K/\Q}(x)|^{2/n} = \Bigl\{ \prod |\sigma(x)|^2
     \Bigr\}^{1/n} \le \Bigl( \sum |\sigma(x)|^2 \Bigr)^{1/n} =
     \cM_K(x). \]
\end{proof}

Die Funktion $\cM_K$ ist für Kreisteilungskörper (abgesehen von dem
Faktor $1/n$) bereits von Gauß eingeführt (Werke II, S. 395) und von
Cassels\index[dN]{Cassels} 1969 ausführlich untersucht worden.
Lenstra\index[dN]{Lenstra} ist es dann 1974 gelungen, mit Hilfe von
$\cM_K$ zu zeigen, dass die Kreisteilungskörper $\Q(\zeta_m)$ für $m =
1, 3, 4, 5, 7, 8, 9, 11, 12, 15, 20$ normeuklidisch
sind. Ojala\index[dN]{Ojala} \cite{Oja77} hat 1977 auch
$\Q(\zeta_{16})$ als normeuklidisch nachgewiesen und dabei (neben
einem Computer) ebenfalls $\cM_K$ benutzt.\footnote{Unlängst hat
  R. McKenzie mit einem Computer gezeigt, dass auch $\Q(\zeta_{13})$
  normeuklidisch ist; s. dazu Leutbecher\index[dN]{Leutbecher} \&
  Niklasch\index[dN]{Niklasch} 1987.} Der Vorteil von $\cM_K$ liegt darin, dass
man $\cM_K(x)$ i.A. viel einfacher berechnen kann als $N_{K/\Q}(x)$.

Während aber Norm und Spur in Körpertürmen transitiv sind, können wir
dies für $\cM_K$ nicht erwarten: für $x \in K$ ist nämlich
nicht notwendig $\cM_K(x) \in \Q$. Ist nämlich
$\sigma$ eine Einbettung von $K$ in $\C$, so ist
$|\sigma(x)|^2 = \sigma(x) \ov{\sigma(x)}$ eine Zahl aus dem
größten reellen Teilkörper des normalen Abschlusses von $K/\Q$
(mit $\ov{\phantom{x}}$ haben wir wie üblich die komplexe
Konjugation bezeichnet). Etwas einfacher liegen die Dinge allerdings,
wenn $K/\Q$ abelsch ist: dann ist die komplexe Konjugation ein
Element der Galoisgruppe $\mathrm{Gal}(K/\Q)$ und daher mit
jedem Automorphismus von $K/\Q$ vertauschbar. Damit folgt
leicht $(\cM_K(x))^\sigma = \cM_K(x)$ für jedes
$\sigma \in \mathrm{Gal}(K/\Q)$, und wir haben in der Tat
$\cM_K(x) \in \Q$ (wesentlich in diesem Beweis ist die
Kommutativität von komplexer Konjugation und $\sigma$; deswegen gilt
$\cM_K(x) \in \Q$ sogar für alle CM-Körper,
sh. Washington\index[dN]{Washington} \cite{Was82}).

Direkt aus der Definition folgt für abelsche $K$ (bzw. für CM-Körper
$K$) die Beziehung $\cM_K(x) = \frac1n T_{K/\Q}(x \ov{x})$, ebenfalls
wieder unter Beachtung der Vertauschbarkeit von $\sigma$ mit den
anderen Automorphismen. Mit Hilfe dieser Beziehung hat Lenstra
\cite{Len75a,Len75b} die beiden nächsten Ergebnisse für
Kreisteilungskörper bewiesen und bemerkt, dass sie auch für CM-Körper
gelten. Verzichtet man aber auf diese Formel, so erhält man etwas
allgemeiner:

\begin{quote}
{\bf (1.20)} {\em Seien $K \subseteq L$ Zahlkörper mit $n = (L:K)$;
  dann gilt für alle $x \in L$ und alle $y \in K$ die Beziehung}
  \[ \cM_L(x) - \cM_L(x-y) = \cM_K\Bigl(\frac{1}{n}
     T_{L/K}(x)\Bigr) - \cM_K\Bigl(\frac{1}{n}
     T_{L/K}(x) - y\Bigr). \]
\end{quote}

Dies wird unser Ersatz für die fehlende Transitivität von
$\cM$ sein; mit (1.20) werden wir nämlich Aussagen über den EA
in $L$ machen können, sobald wir nur gewisse Eigenschaften von $K$ gut
genug kennen.

Bevor wir dies beweisen, erinnern wir an einige bekannte Tatsachen aus
der Körpertheorie (sh. z.B. Marcus\index[dN]{Marcus} 1977): Sei $\Q
\subset K \subset L \subset N$, $N$ der normale Abschluss von $L/\Q$,
$(K:\Q)=k$, $(L:K)=n$, also $(L:\Q) = kn$. Dann gibt es $k$
Einbettungen $\sigma_1, \dots, \sigma_k$ von $K$ in $\C$, sowie $n$
Einbettungen $\tau_1, \dots, \tau_n$ von $L$ in $\C$, welche $K$
elementweise fest lassen. Die $\sigma_i$ und $\tau_j$ können wir zu
Automorphismen von $N$ fortsetzen (i.A. auf viele verschiedene
Weisen), und wir wollen diese (ein für allemal gewählten)
Fortsetzungen ebenfalls mit $\sigma_i$ bzw. $\tau_j$ bezeichnen. Damit
erhalten wir jede Einbettung $\rho$ von $L$ in $\C$ als Einschränkung
der $(L:\Q)$ Automorphismen $\sigma_i \tau_j$ von $N$ auf $L$. Damit
haben wir in (1.20)
\begin{align*}
  (L:\Q) & \{\cM_L(x) - \cM_L(x-y)\} =
  \sum |x^\rho|^2 - \sum |(x-y)^\rho|^2 \\
  & = \sum (x^\rho \ov{x^\rho} - (x-y)^\rho \ov{(x-y)^\rho}) \\
  & = \sum (x^\rho \ov{y^\rho} +
    \ov{x^\rho} y^\rho - y^\rho \ov{y^\rho}).
\end{align*}
Jetzt beachten wir $\sigma(y) = y$ (denn $\sigma$ lässt $K$ punktweise
fest):
\[
= \sum_{\sigma, \tau} \{\sigma \tau (x) \ov{\sigma(y)} +
\ov{\sigma \tau (x)} \sigma(y) - \sigma(y)
\ov{\sigma(y)}\}.
\]
Indem wir die Summation über $\tau$ nach innen ziehen und $T_{L/K}(x)
= \sum \tau(x)$ beachten, erhalten wir weiter
\begin{align*}
& = \sum_{\sigma} \{\sigma(T_{L/K}(x)) \ov{\sigma(y)} +
    \ov{\sigma(T_{L/K}(x))} \sigma(y) - n \cdot \sigma(y)
    \ov{\sigma(y)}\} \\
& = n \cdot \Bigl\{ \sum_{\sigma} \sigma\Bigl(\frac{1}{n} T_{L/K}(x)\Bigr)
    \ov{\sigma\Bigl(\frac{1}{n} T_{L/K}(x)\Bigr)} \\
& \qquad - \sigma\Bigl(\frac{1}{n} T_{L/K}(x) - y\Bigr)
   \ov{\sigma\Bigl(\frac{1}{n} T_{L/K}(x) - y\Bigr)} \Bigr\} \\
& = (L:\Q) \Bigl\{ \cM_K\Bigl(\frac{1}{n} T_{L/K}(x)\Bigr)
   - \cM_K\Bigl(\frac{1}{n} T_{L/K}(x) - y\Bigr) \Bigr\}
\end{align*}

Nun behaupten wir

\begin{quote}
  {\bf (1.21)} {\em Sei $L = K(\zeta_m)$, wo $\zeta_m$ eine $m$-te
    Einheitswurzel ist. Dann gilt für alle $x \in L$ die Beziehung}
  $$ (L:K)\cM_L(x) = \frac{1}{m} \sum_{j=1}^{m}
     \cM_K(T_{L/K}(x \zeta_m^j)). $$
\end{quote}

Für den Fall, dass $K$ Kreisteilungskörper ist, stammt (1.21) in
dieser Form von Lenstra\index[dN]{Lenstra}  \cite{Len75a,Len75b}, geht aber im
wesentlichen auf Cassels\index[dN]{Cassels}  \cite{Cas69} zurück.

\begin{proof}
  In den folgenden Rechnungen mögen $\tau$ und $\tau'$ alle $(L:K)$
  Einbettungen von $L$ in $\C$ durchlaufen, die $K$ elementweise fest
  lassen, $\sigma$ dagegen die Einbettungen von $K$ in $\C$. Dann ist
  \begin{align*}
    (K:\Q) & \cdot \sum_{j=1}^{m} \cM_K(T_{L/K}(x \zeta_m^j)) =
    (K:\Q) \sum_{j=1}^{m} \cM_K\Bigl( \sum_{\tau} \tau(x
      \zeta_m^j) \Bigr) \\
     & = \sum_{j=1}^{m} \sum_{\sigma} \sigma\Bigl( \sum_{\tau} \tau(x
      \zeta_m^j) \Bigr) \ov{\sigma\Bigl( \sum_{\tau} \tau(x \zeta_m^j) \Bigr)} \\
     & = \sum_{\sigma} \sum_{\tau,\tau'} \sum_{j=1}^{m} \sigma(\tau(x))
         \ov{\sigma(\tau'(x))} \sigma(\zeta_m^j)
         \ov{\sigma(\zeta_m^j)},
  \end{align*}  
  wobei wir $\sigma, \tau, \tau'$ wieder als Automorphismen von $N/\Q$
  auffassen. Das Bild einer Einheitswurzel unter einem solchen
  Automorphismus ist wieder eine Einheitswurzel, und wir finden
  \[
  \ov{\sigma(\zeta_m^j)} = \sigma(\zeta_m^{-j}) =
  \sigma(\zeta_m^j)^{-1}.
  \]
  Also ist  
  \[  
  (K:\Q) \sum_{j=1}^{m} \cM_K(T_{L/K}(x \zeta_m^j)) =
  \sum_{\sigma} \sum_{\tau,\tau'} \sum_{j=1}^{m} \sigma(\tau(x))
  \ov{\sigma(\tau'(x))} \sigma(\zeta_m^j)
  \ov{\sigma(\zeta_m^j)}.
  \]
  Jetzt setzen wir $\zeta_{\tau,\tau'} = \sigma(\zeta_m)
  \ov{\sigma(\zeta_m')}$ und finden
  \[
  \zeta_{\tau,\tau'} = 1 \Leftrightarrow \tau(\zeta_m) =
  \tau'(\zeta_m) \Leftrightarrow \tau = \tau', \text{ sowie}
  \]
  \[
  \sum_{j=1}^{m} \zeta_{\tau,\tau'}^j =
  \begin{cases}
    m, & \text{falls } \tau = \tau'; \\
    0, & \text{sonst}.
  \end{cases}
  \]
  Indem wir die Summation über $j$ nach innen ziehen, folgt
  \begin{align*}
    (K:\Q) \sum \cM_K(T_{L/K}(x \zeta_m^j))
    & = m \cdot \sum_{\sigma} \sum_{\tau} \sigma(\tau(x))
    \ov{\sigma(\tau(x))}\\
    & = m \cdot \sum_{\sigma} |x^\sigma|^2 = m \cdot (L:\Q)
    \cdot \cM_L(x).
  \end{align*}
  Durch Division mit $m \cdot (K:\Q)$ erhalten wir die Behauptung.
  \end{proof}

Lenstras Idee war nun, die Existenz des EA in $L=K(\zeta_m)$ auf
Eigenschaften von $K$ zurückzuführen. Dazu definiert man
\begin{align*}
  F & = F_K = \{x \in K: \cM_K(x) \le \cM_K(x-y)
        \text{ für alle } y \in R\}, \quad \text{sowie} \\
  c(K) & = \sup \{ \cM_K(x) : x \in F_K \}.
\end{align*}
Direkt aus der Definition folgt dann, dass es zu jedem $x \in K$
ein $y \in R$ gibt mit $\cM_K(x-y) \le c(K)$. Wegen (1.19) ist
$R$ sicher dann normeuklidisch, wenn $c(K)<1$ ist. Tatsächlich genügt
oft schon $c(K) \le 1$, wie wir nun zeigen wollen.

\begin{quote}
  {\bf (1.22)} {\em Sei $K$ ein Zahlkörper, $\sigma$ eine Einbettung
    von $K$ in $\C$, und es gelte $|\sigma(x)| = |\sigma(x-u)| = 1$
    für ein $x \in K$ und eine Einheitswurzel $u \in R$; dann ist
    $x \in R$.}
\end{quote}

\begin{proof}
  Wegen $x^\sigma \ov{x^\sigma} = 1$ und $u^\sigma \ov{u^\sigma} = 1$
  gilt
  \[ 1 = |\sigma(x-u)| = (x^\sigma - u^\sigma)(\ov{x^\sigma} -
         \ov{u^\sigma}) = 2 - x^\sigma \ov{u^\sigma} -
         \ov{x^\sigma} u^\sigma. \]
  Für $y = \sigma(-x/u)$ gilt also $y \ov{y} = 1$ und
  \[ y + \ov{y} = -x^\sigma / u^\sigma - \ov{x^\sigma} /
   \ov{u^\sigma} = -(x^\sigma \ov{u^\sigma} +
   \ov{x^\sigma} u^\sigma)/u^\sigma \ov{u^\sigma} = -1. \]
  Damit ist $y$ Nullstelle von $z^2 + z + 1 = 0$, d.h. $y$ ist eine
  dritte Einheitswurzel. Insbesondere ist $y$ ganz, und dies
  impliziert, dass auch $x/u$ und $x$ ganz sind; dies war zu zeigen.
\end{proof}

Ein reelles $c' \in \R$ heißt eine \textbf{brauchbare Schranke}
\index[dS]{brauchbare Schranke}
(``usable bound'' bei Lenstra) für $K$, wenn $c' \ge c(K)$
gilt und es für alle $x \in F_K$ mit $\cM_K(x)=c'$ eine
Einheitswurzel $u \in R$ gibt mit $\cM_K(x-u)=c'$. Man
beachte, dass jedes $c' > c(K)$ brauchbar ist.

\begin{quote}
  {\bf (1.23)} {\em Ist $c'=1$ brauchbar für $K$, dann ist $K$
    normeuklidisch.}
\end{quote}

\begin{proof}
  Sei ein $x \in K$ gegeben; wir suchen ein $y \in R$ mit
  $|N_{K/\Q}(x-y)| < 1$. OBdA dürfen wir $x \in F_K$ annehmen; dann
  ist $\cM_K(x) \le 1$, wegen (1.19) also
  $|N_{K/\Q}(x)| \le \cM_K(x) \le 1$. Also genügt $y=0$, falls
  nicht gerade $|N_{K/\Q}(x)| = 1$ ist. Tritt dies ein, so muss auch
  $\cM_K(x) = 1$ sein, und da $c=1$ eine brauchbare Schranke
  für $K$ ist, gibt es eine Einheitswurzel $u \in R$ mit
  $\cM_K(x) = \cM_K(x-u) = 1$. Ist $|N_{K/\Q}(x-u)| < 1$,
  so können wir $y=u$ wählen; andernfalls haben wir
  $|N_{K/\Q}(x)| = |N_{K/\Q}(x-u)| = \cM_K(x) = \cM_K(x-u) = 1$.
  Aus der Gleichheit in (1.19) folgt dann
  $|\sigma(x)| = |\sigma(x-u)| = 1$ für alle Einbettungen $\sigma$, und
  nach (1.22) ist dann $x \in R$, sodass hier $y=x$ genügt.
\end{proof}

Zentrales Ergebnis ist nun

\begin{quote}
  {\bf (1.24)} {\em Sei $\zeta_m$ eine $m$-te Einheitswurzel und
    $L = K(\zeta_m)$; dann gilt
    $$ c(L) \le (L:K)\cdot c(K). $$
    Ist darüberhinaus $c'$ eine brauchbare Schranke für $K$, dann ist
    $(L:K) \cdot c'$ eine solche für $L$.}
\end{quote}

\begin{proof}
  Sei $x \in F_L$; wir müssen zeigen, dass $\cM_L(x) \le
  (L:K)c(K)$ gilt. Wegen $x \in F_L$ ist für alle $y \in S$ und damit
  erst recht für alle $y \in R$, $\cM_L(x-y) \ge
  \cM_L(x)$. Also garantiert (1.20) $\frac{1}{n}T_{L/K}(x) \in
  F_K$ für $n = (L:K)$. Nun ist $|\sigma(\zeta_m)| = 1$ für jede
  Einbettung von $L$ in $\C$ und jedes $j \in \N$; also ist
  $\cM_L(x\zeta_m^j) = \cM_L(x)$, und wie oben
  erhalten wir nun $\frac{1}{n} \cdot T_{L/K}(x\zeta_m^j) \in F_K$ für
  $j=1, \dots, m$. Daher ist $\cM_K(T_{L/K}(x\zeta_m^j)) = n^2
  \cdot \cM_K(\frac{1}{n}T_{L/K}(x\zeta_m^j)) \le n^2 \cdot
  c(K)$. Mit (1.21) erhalten wir hieraus
  $$ mn \cdot \cM_L(x) = \sum_{j=1}^{m}
    \cM_K(T_{L/K}(x\zeta_m^j)) \le \sum_{j=1}^{m} n^2 \cdot c(K) =
    mn^2 \cdot c(K), $$
  sodass wir schließlich $\cM_L(x) \le n \cdot c(K)$ haben.
\end{proof}

Sei nun $c'$ eine brauchbare Schranke für $K$ und $\cM_L(x) =
nc'$ für ein $x \in F_L$. Wegen $\cM_L(x) \le n \cdot c(K)$
und $c \ge c(K)$ muss $c'=c(K)$ sein, und aus obigem Beweis folgt
$\cM_K(\frac{1}{n}T_{L/K}(x\zeta_m^j)) = c(K) = c'$ für alle
$j \in \Z$. Indem wir $j=0$ wählen, sehen wir, dass $\alpha =
\frac{1}{n}T_{L/K}(x)$ ein Element von $F_K$ ist mit
$\cM_K(\alpha) = c'$. Da $c'$ eine brauchbare Schranke für $K$
ist, gibt es eine Einheitswurzel $u \in R$ mit $\cM_K(\alpha -
u) = c'$. Mit $y=u$ erhalten wir aus (1.20) $\cM_L(x-u) =
\cM_L(x) = nc'$, und dies zeigt, dass $nc'$ eine brauchbare
Schranke für $L$ ist.

Für $K=\Q$ ist $\cM_\Q(x) = |x|^2$,
$F_\Q = [-0.5, +0.5]$ und $c(\Q) = \frac14$; um zu zeigen,
dass $c' = \frac14$ eine brauchbare Schranke für $\Q$ ist, müssen
wir zu jedem $x \in F_\Q$ mit $\cM_\Q(x) =
\frac14$ eine Einheitswurzel $u \in \Z$ finden mit
$\cM_\Q(x-u) = \frac14$. Da aber $c'$ auf $F_\Q$
nur in den Punkten $0.5$ und $-0.5$ angenommen wird, genügt $u=1$.

Ist also $\zeta_m$ eine primitive $m$-te Einheitswurzel und
$L=\Q(\zeta_m)$, so folgt aus (1.24), dass $c(L) \le
(L:\Q)/4$ gilt und $(L:\Q)/4$ eine brauchbare Schranke
ist. Also sind alle Kreisteilungskörper $L=\Q(\zeta_m)$ mit
$(L:\Q) \le 4$ normeuklidisch: dies sind $\Q(\zeta_m)$
für $m = 3, 4, 5, 8, 12$. Für Kreisteilungskörper
$K=\Q(\zeta_p)$, $p$ ungerade Primzahl, hat Lenstra\index[dN]{Lenstra}  \cite{Len74} die
Konstanten $c(K)$ auf elementare Weise bestimmt:

\begin{quote}
  {\bf (1.25)} {\em Sei $K = \Q(\zeta_p)$, $p \equiv 1 \bmod {2}$
    prim; dann ist $c(K) = \frac{p+1}{12}$, und $c'=c(K)$ ist
    brauchbare Schranke für $K$.}
\end{quote}

Daraus folgt sofort, dass $\Q(\zeta_7)$ und
$\Q(\zeta_{11})$ ebenfalls normeuklidisch sind; man erhält
jedoch weiter für

\begin{itemize}
\item $L = \Q(\zeta_9) = K(\zeta_3)$, $K = \Q(\zeta_3)$,
  $c(K) = \frac13$, also $c(L) \le (L:K) \cdot c(K) = 1$;
\item $L = \Q(\zeta_{15}) = K(\zeta_5)$ für $K = \Q(\zeta_5)$,
  $c(K) = 1/2$, $c(L) \le 2 \cdot 1/2 = 1$;
\item $L = \Q(\zeta_{20}) = K(\zeta_4)$ für $K = \Q(\zeta_5)$,
  $c(K) = 1/2$, $c(L) \le 2 \cdot 1/2 = 1$.
\end{itemize}

Also sind auch $\Q(\zeta_9)$, $\Q(\zeta_{15})$ und $\Q(\zeta_{20})$
normeuklidisch.

Ist $K$ imaginärquadratischer Zahlkörper, so gilt $\cM_K(x) =
N_{K/\Q}(x)$ und damit $c(K) = M(K)$; (0.18) zeigt, dass die
einzigen imaginärquadratischen Körper mit $c(K) \le 1/2$ die beiden
Kreisteilungskörper $\Q(\zeta_3)$ und $\Q(\zeta_4)$
sind.

Ist dagegen $K = \Q(\sqrt{m})$ reellquadratisch, so gilt für
$x=a+b\sqrt{m}$
\[
\cM_K(x) = \frac{(a+b\sqrt{m})^2 + (a-b\sqrt{m})^2}{2} = a^2 +
m b^2,
\]
und der Beweis von (0.21) zeigt

\begin{quote}
  {\bf (1.26)} {\em Sei $K = \Q(\sqrt{m})$ reellquadratischer
    Zahlkörper; dann gilt}
   \[ c(K) =
   \begin{cases}
     \frac{1+m}{4      } & \text{ {\em für }} m \equiv 2,3 \bmod {4} \\
     \frac{(1+m)^2}{16}  & \text{ {\em für }} m \equiv 1 \bmod {4}.
   \end{cases}
   \]
\end{quote}

Für $K = \Q(\sqrt{5})$ ist also $c(K) = 9/20$, und (1.24)
liefert die normeuklidischen Zahlkörper $\Q(\sqrt{5},
\sqrt{-1})$, $\Q(\sqrt{5}, \sqrt{-3})$ und
$\Q(\zeta_5)$.

Um eine weitere Anwendung dieser Methode zu geben, betrachten wir eine
quadratische Erweiterung von $K = \Q(\sqrt{-3})$: eine solche
lässt sich in der Form $L = K(\sqrt{\mu})$ mit quadratfreiem $\mu \in
\Z[\sqrt{-3}]$ schreiben. Wir setzen $\mu = a+b\sqrt{-3}$, $\mu' =
a-b\sqrt{-3}$, $\mu\mu' = a^2+3b^2 = p$, und finden für ein $x =
r+s\rho +t\sqrt{\mu} +u\rho\sqrt{\mu} \in L$ ($\rho =
(-1+\sqrt{-3})/2$) nach etwas Rechnung
\[
\cM_{L/\Q}(x) = r^2 - rs + s^2 + \sqrt{p}(t^2 - tu + u^2).
\]
Im Falle $\mu \equiv 1 \bmod {4}$ ist nun $\beta = \frac{1+\sqrt{\mu}}2$
ganz; indem wir $x$ in der Form $x = e+f\rho + g\beta + h\rho\beta$
schreiben, erhalten wir in obiger Notation $r=e+g/2$, $s=f+h/2$,
$t=g/2$, $u=h/2$, also
\[
\cM_{L/\Q}(x) = (e+g/2)^2 - (e+g/2)(f+h/2) + (f+h/2)^2
+ \sqrt{p}(g^2 - gh + h^2)/4.
\]
Nun wissen wir aber, dass $c(K) = \frac13$ ist; folglich können wir $g$ und
$h \bmod {\Z}$ so wählen, dass $g^2 - gh + h^2 =
\cM_{K/\Q}(g + h\rho) \le \frac13$ wird. Jetzt wählen wir
$e+g/2$ und $f+h/2 \bmod \Z$ so, dass auch $(e+g/2)^2 -
(e+g/2)(f+h/2) + (f+h/2)^2 \le \frac13$ wird, und haben damit gezeigt: Für
alle $x \in L$ können wir ein $y \in S$ finden, sodass
\[
\cM_{L/\Q}(x-y) \le \frac13 + \frac{\sqrt{p}}{12}
              = \frac{4 + \sqrt{p}}{12}
\]
wird.

Insbesondere ist also $c(L) \le (4 + \sqrt{p})/12$, und $L$ ist
normeuklidisch, falls nur $\sqrt{p} \le 8$ (und natürlich
$\mu \equiv 1 \bmod {4}$) ist. Damit erhalten wir die folgenden
norm\-eu\-kli\-di\-schen Körper: $K(\sqrt{\mu})$ für
$$ \mu = -1+2\sqrt{-3}, \ 3+2\sqrt{-3}, \ 5, \ -5+2\sqrt{-3}, \
         -7, \ -3+4\sqrt{-3}, \ 7-2\sqrt{-3}; $$
die Diskriminanten dieser Körper sind übrigens $\disc L = 117$, $189$,
$225$, $333$, $441$, $513$, $549$. Dass diese Körper normeuklidisch
sind, hat schon Lakein\index[dN]{Lakein} \cite{Lak72} bewiesen,
allerdings mit einer völlig anderen Methode.

\begin{quote}
  {\bf (1.27)} {\em Sei $L = \Q(\sqrt{a+b\sqrt{-3}})$, $p=a^2+3b^2$;
    dann ist $c(L) \le \frac{4+\sqrt{p}}{12}$.}
\end{quote}

Es sei noch bemerkt, dass für die Funktion $\cM_K$ ein Analogon
zu (1.14) existiert: setzt man nämlich $\cM_K$ auf den
$\R^n$ fort durch
\[
\cM(x) = \frac{x_1^2 + \dots + x_r^2 + 2(x_{r+1}^2 +
  x_{r+2}^2) + \dots + 2(x_{n-1}^2 + x_n^2)}{n}
\]
für ein $x = (x_1, \dots, x_n) \in \R^n$ und setzt $U = \{x
\in \R^n : \cM(x) < \frac14\}$, so ist $\cM(u-v) <
1$ für alle $u,v \in U$. Eine einfache Koordinatentransformation und
eine klassische Volumenberechnung (\glqq{}Volumen einer Hyperkugel im
$\R^n$\grqq{}) zeigt dann
$$ \mathrm{vol}(U) = 2^{-s} \Big(\frac{\pi}4\Big)^{n/2}
      \Gamma\Big(1 + \frac{n}2\Big), $$
wo $\Gamma$ die Gammafunktion bezeichnet. Mit einem Packungsargument
(sh. dazu Rogers\index[dN]{Rogers} \cite{Rog64} und
Lenstra\index[dN]{Lenstra} \cite{Len77a}) erhält man dann

\begin{quote}
  {\bf (1.28)} {\em Sei $K$ ein Zahlkörper mit $(K:\Q)=n$ und
    $\delta < 1$; dann gibt es für alle $x \in K$ ein $y \in R$ und
    eine Einheit $u \in R^\times$ mit $\cM_{K/\Q}(u x - y) < \delta$,
    wobei
    \[ \delta = \sigma_n  \Gamma\Big(1 + \frac{n}2\Big)
           \Big(\frac{4}{n\pi}\Big)^{n/2}
           \frac{\sqrt{|\disc K|}}{\mu_1} \]
    ist.}
\end{quote}

Hierbei sind die $\sigma_n$ die z.B. bei Lenstra\index[dN]{Lenstra}  \cite{Len77a} angegebenen,
nur von $n$ abhängigen Größen, und $\mu_1$ gibt (wie auf S. 33) die
Länge einer maximalen 1-Folge in $K$ an. Mit (1.19) folgt nun sofort

\begin{quote}
  {\bf (1.29)} {\em Unter den Voraussetzungen von (1.28) gilt: für
    alle $x \in K$ gibt es ein $y \in R$ mit $|N_{K/\Q}(x-y)| <
    \delta^{n/2}$. Insbesondere ist $K$ dann normeuklidisch.}
\end{quote}

Dies ist wie (1.14) ein Ergebnis, das auf Lenstra\index[dN]{Lenstra}  \cite{Len77a} zurückgeht;
die folgende Tabelle zeigt, dass (1.29) in den Fällen $n=2$, $s=1$; $4
\le n \le 7$, $s \ge 2$; $n=8$, $s \ge 3$ bessere Schranken als (1.15)
und (1.16) liefert:

$$ \begin{array}{ccc}
  \toprule n & \sigma_n \Gamma(1 + \frac n2) / \pi^{n/2} &
               \sigma_n \Gamma(1 + \frac n2) (4/n\pi)^{n/2} \\
  \midrule
  1 & 0.5 & 1 \\
  2 & \sqrt{3}/6 & \sqrt{3}/3 \\
  3 & 0.18613 & 0.286566  \\
  4 & 0.13128 & 0.131280 \\
  5 & 0.09988 & 0.037175 \\
  6 & 0.08113 & 0.024039 \\
  7 & 0.06982 & 0.009848 \\
  8 & 0.06327 & 0.003955 \\ \bottomrule
\end{array} $$

Ersetzt man $\mu_1$ in (1.28) durch $\mu_k$, so erhält man natürlich
entsprechende Ergebnisse für $k$-euklidische Ringe. Da $\cM_K$
nicht multiplikativ ist, erhält man aus (1.28) keine Schranken für
$c(K)$; geht man dagegen anders vor und setzt $F(K) = \{x \in K :
\cM(x) \le \cM(u x - y) \text{ für alle } u \in R^\times
\text{ und alle } y \in R\}$, sowie $c(K) = \sup \{ \cM(x) : x
\in F \}$, so liefert (1.28) zwar Abschätzungen für $c(K)$, jedoch hat
man kein Ergebnis wie (1.20) oder (1.24) für die $c(K)$ vorliegen.

\section*{{\sc Anmerkungen zu} \S\ 1}
\addcontentsline{toc}{section}{{\sc Anmerkungen zu} \S\ 1}
                  
Die Existenz rein verzweigter Primideale lässt sich auch zum Auffinden
von Faktoren der Klassenzahl verwenden; das bekannteste Ergebnis
dieser Form ist wohl die Geschlechtertheorie von Gauß, der (in unserer
Sprache) gezeigt hat, dass die Klassenzahl eines reell-quadratischen
Zahlkörpers mit $t$ verzweigten Primidealen durch $2^{t-2}$
bzw. $2^{t-1}$ teilbar ist, je nachdem die Norm der Grundeinheit $+1$
oder $-1$ ist. Für ähnliche Ergebnisse in Zahlkörpern höheren Grades
sh. Ishida\index[dN]{Ishida}  \cite{Ish76}.

Die Verwendung rein verzweigter Ideale im Sinne von (1.5) findet sich
andeutungsweise bereits bei Behrbohm\index[dN]{Behrbohm} \&
Rédei\index[dN]{Redei@R\'edei} \cite[p. 198]{BR36}; deren Methode ist
dann von Erdős\index[dN]{Erdos@Erd\"os} und Ko\index[dN]{Ko} \cite{EK38}
modifiziert und von Heilbronn\index[dN]{Heilbronn}
\cite{Hei38,Hei50,Hei51}, Cioffari\index[dN]{Cioffari} \cite{Cio79} und
Egami\index[dN]{Egami} \cite{Ega79} weiterentwickelt
worden\footnote{Beim letzten Abfassen der Arbeit habe ich bemerkt,
  dass sich zur Abschätzung von $M(K)$ auch Ideale verwenden lassen,
  die nicht rein verzweigt sind, sh. dazu § 4.}.

Die Betrachtung der Ideale $I = (u-1)$ für Einheiten $u$ zur
Abschätzung von $M(K)$ hat erstmals Rédei\index[dN]{Redei@R\'edei}
\cite{Red41b} vorgeschlagen; er konnte damit z.B. zeigen, dass
$\Q(\sqrt{61})$ und $\Q(\sqrt{109})$ nicht normeuklidisch sind.

Ist $0,1,\vartheta_2,\dots,\vartheta_m$ eine 1-Folge, so muss
jedenfalls $\vartheta_i-1$ für $i \ge 3$ eine Einheit sein. Einheiten
$u$ mit der Eigenschaft, dass auch $u-1$ Einheit ist, heißen
\textbf{Ausnahmeeinheiten}\index[dS]{Ausnahmeeinheit} (exceptional
units). Ausnahmeeinheiten sind schon untersucht worden, bevor Lenstra
ihre Bedeutung für die Existenz eines EA erkannt hat. Julia Robinson
hat vermutet, dass es in einem algebraischen Zahlkörper nur endlich
viele Ausnahmeeinheiten gibt; diese Vermutung haben dann
S.~Lang\index[dN]{Lang} \cite{Lan60}, S. Chowla\index[dN]{Chowla}
\cite{Cho61} und T.~Nagell\index[dN]{Nagell} \cite{Nag64} unabhängig
voneinander bewiesen. Nagell hat sich in einer Reihe von Arbeiten mit
der Bestimmung aller Ausnahmeeinheiten in Zahlkörpern mit kleinem
Einheitenrang befasst.

Schließlich haben Chudnovsky u. Chudnovsky\index[dN]{Chudnovsky}
\cite{Chu86} gezeigt, dass die Existenz ``schneller'' Algorithmen zur
Multiplikation von Polynomen über $\Z$ von der Existenz von
Zahlkörpern mit vielen Ausnahmeeinheiten abhängt.

\chapter*{\S\ 2 Die Bestimmung euklidischer Minima}
\setcounter{chapter}{2}
\addcontentsline{toc}{chapter}{\S\ 2 Die Bestimmung euklidischer Minima}
\markboth{Euklidische Ringe}{\S\ 2 Die Bestimmung euklidischer Minima}

Am Beginn des \S\ 1 haben wir gesehen, wie man

 untere Schranken für das
euklidische Minimum $M(K)$ findet; hier nun wollen wir Methoden
vorstellen, die es gestatten, $M(K)$ nach oben abzuschätzen. Wenn
beide Schranken übereinstimmen, hat man $M(K)$ bestimmt.

Sei dazu $\{\alpha_1, \dots, \alpha_n\}$ eine $\Q$-Basis eines
algebraischen Zahlkörpers $K$ (i.\,A. werden wir hier eine GHB wählen;
jedoch empfiehlt es sich manchmal, aus Symmetriegründen etwa
$\{1,\sqrt{m}\}$ für quadratische Zahlkörper $\Q(\sqrt{m})$
auch im Falle $m \equiv 1 \bmod {4}$ zu nehmen -- die Gründe hierfür
werden noch klar werden). Dann definieren wir durch
$$ \varphi : K \to \R^n =: \uK\label{dpuK}:
   (r_1\alpha_1 + \dots + r_n\alpha_n) = (r_1, \dots, r_n) $$
eine Einbettung von $K$ in $\uK = \R^n$. Wie schon in \S\ 1 werden wir
$\R$ und $\varphi(\R)$ i.\,A. identifizieren. Durch
\[
|\cdot|_i : \uK \to \R : x = (x_1, \dots, x_n) \mapsto
\Big| \sum_{j=1}^{n} x_j \tau_i(\alpha_j)\Big|^\delta,
\]
wo $\delta = 1$ für reelle und $\delta = 2$ für nichtreelle
Einbettungen von $K$ in $\C$ ist, haben wir die Bewertungen von $K$
auf $\uK$ \glqq{}fortgesetzt\grqq{} in dem Sinne, dass für alle $x \in
K$ $|x|_1 = |\varphi(x)|_1$ gilt; man beachte jedoch, dass $|\cdot|_i$
auf $\uK$ natürlich keine Bewertungen\label{dpBew2} mehr sind, schon
allein deshalb nicht, weil $\uK$ kein Körper ist (wir fassen $\uK$ nur
als $\R$-Vektorraum auf). Weiter folgt aus $|x|_1 = 0$ i.\,A. nicht,
dass $x = 0$ ist! Durch $N(x) = |x|_1 \cdots |x|_{r+s}$ können wir
auch die Norm $|N_{K/\Q}|$ von $K$ nach $\uK$ fortsetzen und haben
$|N_{K/\Q}(x)| = N(\varphi(x))$ für alle $x \in K$. Auch hier müssen
wir bemerken, dass $N$ auf $\uK$ keine Norm im topologischen Sinne
ist, weil aus $N(x)=0$ nicht $x=0$ folgt.

\medskip\noindent
\textbf{Beispiel:}
sei $K=\Q(\sqrt{7})$; wir wählen die $\Q$-Basis $\{1, \sqrt{7}\}$;
dann ist $\varphi(a+b\sqrt{7}) = (a,b)$ für alle $a,b \in \Q$.  Also
ist $\varphi(K) = \Q^2$, während $\uK = \R^2$ ist. Weiter haben wir
$|x|_1 = |a+b\sqrt{7}|$, $|x|_2 = |a-b\sqrt{7}|$ und $N(x) =
|a^2-7b^2|$ für ein $x = (a,b) \in \uK$. Insbesondere ist $|x|_2 = 0$
und damit auch $N(x) = 0$ für $x = (\sqrt{7},1) \in \uK \setminus
\{0\}$. Man beachte, dass durch $\{x \in \uK: |x|_1 = 0\}$ eine Gerade
im $\R^2$ definiert ist.
\medskip

Indem wir $\uK$ mit der gewöhnlichen euklidischen Metrik
versehen, werden ${|\cdot|_1}$, $\ldots$, $|\cdot|_{r+s}$ und $N$ zu stetigen
Funktionen auf $\uK$, und für jedes $\eps \in
\R^+$ ist $U_\eps = \{x \in \uK: |x|_i < \eps \text{ für } 1 \leq i \leq r+s\}$
eine beschränkte Umgebung von $0$ (die Beschränktheit sieht man bei
Verwendung von (1.10) leicht ein). Ist $\{\alpha_1, \dots, \alpha_n\}$
eine GHB von $K$, so heißt
$$ F = \{(x_1, \dots, x_n) \in \uK : - \frac12 \leq x_i < \frac12
     \text{ für } i = 1, \dots, n\} $$
ein \textbf{Fundamentalbereich}\index[dS]{Fundamentalbereich} von
$K$. Mit $\uF$ bezeichnen wir den kompakten Abschluss
$$ \uF = \{(x_1, \dots, x_n) \in \uK: - \frac12 \leq x_i \leq \frac12
  \text{ für } i = 1, \dots, n\} $$
von $F$.

Um nun zu zeigen, dass für alle $x \in K$ ein $y \in \cO_K$
existiert mit $|N_{K/\Q}|(x-y) < k$, genügt es, zu jedem $x \in \uK$
ein $y \in \cO_K$ zu finden mit $N(x-y) < k$. Damit wird folgende
Definition gerechtfertigt:
\[
M(\uK) := \inf \{k \in \R: \forall x \in \uK \ \exists y
\in \cO_K: N(x-y) < k\};
\]
entsprechend definieren wir auch die weiteren Minima $M_2(\uK)$,
$M_3(\uK)$ usw., soweit sie existieren. Wegen $K \subset \uK$ ist
die Ungleichung $M(K) \leq M(\uK)$ sofort klar. Ob aus der
Isolation von $M(K)$ auch diejenige von $M(\uK)$ folgt, weiß
man nicht; mit $M(\uK)$ ist aber sicher auch $M(K)$ isoliert,
und wir haben $M_2(K) \leq M_2(\uK)$. Entsprechendes gilt für
die höheren Minima.

Setzt man $M(F) = \inf \{k \in \R: \forall x \in F \ \exists y
\in \cO_K \text{ mit } N(x-y) < k\}$, so ist offenbar
$M(F) = M(\uF) = M(\uK)$. Da $N$ als stetige Funktion auf der
kompakten Menge $\uF$ beschränkt ist, ist
$M(\uF) \le \max \{N(x): x \in \uF\}$, und insbesondere ist
$M(\uK)$ (und damit auch $M(K)$) endlich.

Die Methode zur Bestimmung von $M(K)$, die wir nun vorstellen wollen,
geht auf Barnes und Swinnerton-Dyer zurück und verläuft immer
folgendermaßen: man zeigt $M(K) \ge k$ für ein $k \in \R$ mit
den Kriterien aus \S\  1, sowie $M(\uK) \le k$ für dasselbe $k$;
damit ist $M(K) \ge M(\uK) \le k$, also $M(K) = M(\uK) =
k$. Diese Bemerkung zeigt bereits, dass diese Methode nur dann
erfolgreich sein kann, wenn $M(K) = M(\uK)$ gilt; man vermutet
allerdings, dass dies in Zahlkörpern immer der Fall ist. Wir wollen
diese Methode nun am Beispiel $K=\Q(\sqrt{14})$
beschreiben. Dazu wählen wir die GHB $\{1,\sqrt{14}\}$ und setzen
$(a,b) \times (c,d) = \{t+s\sqrt{14}: a \le t \le b, c \le s \le d\}$;
damit ist z.B. $F = (-0.5, 0.5) \times (-0.5, 0.5)$. Um zu zeigen, dass
$M(K) \ge \frac{3}{4}$ ist, wählen wir ein $k \in \R$, das
etwas kleiner als $\frac{3}{4}$ ist, z.B. $k=1.2$ (der genaue Wert von
$k$ spielt hier keine Rolle; um $M(K)$ zu bestimmen, sollte er nur
zwischen $M_0(K)$ und $M(K)$ liegen). Wir nennen dann eine Menge
$(a,b) \times (c,d) \subset \uK$ bedeckt,\index[dS]{bedeckt} wenn es
für alle $x \in (a,b) \times (c,d)$ ein $y \in \cO_K$ gibt mit
$N(x-y) < k$ (ob eine solche Menge bedeckt ist, hängt also von $k$ ab).
Ein $x \in \uK$ mit $N(x-y) \ge k$ für alle $y \in \cO_K$
nennen wir einen \textbf{Ausnahmepunkt}.\index[dS]{Ausnahmepunkt}

Die Programme, die wir in \S\  11 genauer beschreiben werden, erlauben es
uns, ganz $F$ bis auf die folgenden Gebiete zu bedecken:
\begin{align*}
  S_1 & = (0.499, 0.5) \times (0.499975, 0.5), \\
  S_2 & = (0.499, 0.5) \times (-0.5, -0.499975), \\
  S_3 & = (-0.5, -0.499) \times (0.499975, 0.5), \\
  S_4 & = (-0.5, -0.499) \times (-0.5, -0.499975).
\end{align*}
Man sieht jetzt, dass es etwas ökonomischer ist, statt $F$ die Menge
$F' = (0,1) \times (0,1)$ zu betrachten, weil dann nur $S = (0.499,
0.501) \times (0.499975, 0.500025)$ unbedeckt bleibt und man es nur
mit einer Menge statt mit vier zu tun hat. Sei nun $x \in F'$ ein
Ausnahmepunkt (falls es einen solchen gibt; in der Tat wissen wir aus
\S\  1, dass $M(K) = \frac{3}{4}$ für $x = (1+\sqrt{14})$ gilt); dann ist
sicher $x \in S$. Mit $x$ ist aber auch $xu$ ein Ausnahmepunkt, wo $u
= 15+4\sqrt{14}$ die FE von $K$ ist, d.h. es gibt ein $a \in
\R$ mit $xu-a \in S$. Diese Überlegung führt uns dazu, die
Menge $uS = \{ux: x \in S\}$ zu betrachten; eine einfache Rechnung
zeigt
\[
uS = (35.4836, 35.5164) \times (9.495625, 9.504375), \text{ also}
\]
\[
uS - (35 + 9\sqrt{14}) = (0.4836, 0.5164) \times (0.495625, 0.504375).
\]
Mit $a = 35 + 9\sqrt{14}$ gilt daher: für jedes $x \in S$ liegt $ux-a$
entweder in bedecktem Gebiet oder wieder in $S$ (hätten wir weniger
genau gerechnet und z.B. $F'$ nur bis auf die Menge $S' = (0.47, 0.53)
\times (0.47, 0.53)$ bedecken können, so wäre $uS' = (33.37, 37.63)
\times (8.93, 10.07)$ gewesen und es gäbe kein $a \in \R$ mit
obiger Eigenschaft, weil z.B. $uS' - (35 + 9\sqrt{14}) = (-1.63, 1.63)
\times (-0.07, 1.07)$ neben Punkten aus $S'$ auch solche aus $S'-1$
oder $S'+1$ enthält).

Wenn wir aber eine Situation wie oben vorliegen haben, können wir die
Lage etwaiger Ausnahmepunkte recht genau beschreiben; es gilt nämlich

\begin{quote}
  {\bf (2.1)} {\em Sei $K$ Zahlkörper, $u$ eine Einheit in $R$, $S$
    eine Teilmenge von $F$ und $a \in R$ ein Punkt mit der
    Eigenschaft, dass für alle $x \in S$ $ux-a$ in bedecktem Gebiet
    oder wieder in $S$ liegt. Dann gilt für jeden Ausnahmepunkt $x_0$
    und für jede Bewertung $| \cdot |_1$ mit $1 < |u|_1$ : $|x_0 -
    z|_1 = 0$; hierbei ist $z = \frac{a}{u-1}$ der Fixpunkt der Abbildung
    $x \mapsto ux-a$. Gilt darüberhinaus $|u|_1 \neq 1$ für alle
    Bewertungen $| \cdot |_1$ von $K$, so wird durch $x_{i+1} = ux_i -
    a$ eine Folge $x_0, x_1, x_2, \dots$ von Ausnahmepunkten definiert
    mit $x_i \in S$ und $\lim x_i = z$.}
\end{quote}

\begin{figure}[ht!]
\centering  
\begin{tikzpicture}[scale=1.5]
\draw[thin, gray!50] (-2.2,-3) grid [xstep=1, ystep=2] (3.5,2.8);
\draw[thick, ->] (-2.6,0) -- (4,0);
\draw[thick, ->] (0,-2.5) -- (0,3.2);
\draw[line width=1pt] (0.8,1.8) rectangle (1.2,2.2);
\draw[dashed] (0.6,1.6) rectangle (1.4,2.4);
\draw[line width=1pt] (-1,1.8) rectangle (-0.8,2);
\draw (-0.8,1.6) node {$S_3$};
\draw[line width=1pt] (-1,-1.8) rectangle (-0.8,-2);
\draw (-0.8,-1.6) node {$S_4$};
\draw[line width=1pt] ( 1,-1.8) rectangle (1.2,-2);
\draw (1.2,-1.6) node {$S_2$};
\draw (1,2) -- (1.5,3);
\draw (1.5,3.15) node {$z = \frac{1+\sqrt{14}}2$};
\draw (0.1,0.2) node {$0$};
\draw (1.1,0.2)    node {$\frac12$};
\draw (2.1,0.2) node {$1$};
\begin{scope}
  \clip (-2.7,-0.3) rectangle (4,3);
  \draw (-3,-2/7) -- ( 4,26/7);
  \draw ( 5,-2/7) -- (-2,26/7);
  \fill (1,2) circle (1pt);
\end{scope}
\draw (-1.5,0.75) node[rotate= 29] {$|x-z|_1 = 0$};
\draw ( 3.5,0.75) node[rotate=-29] {$|x-z|_2 = 0$};
\end{tikzpicture}
\caption{Ausnahmemengen in $\Q(\sqrt{14}\,)$}
\end{figure}

In unserem Beispiel $K = \Q(\sqrt{14})$ ist $u = 15 + 4\sqrt{14}$,
$a = 35 + 9\sqrt{14}$ und damit
$z = a/(u-1) = \frac{35 + 9\sqrt{14}}{14 + 4\sqrt{14}} =  1 + \sqrt{14}$;
die Menge aller $x \in K$ mit $|x-z|_1 = 0$ ist hier eine Gerade durch
den Punkt $z$. Als einziger Ausnahmepunkt (für $k=1,2$) wird sich hier
$x_0 = z$ herausstellen; damit ist $x_1 = ux_0 - a = x_2 = z$, sodass
die Folge $x_0, x_1, x_2, \dots$ trivialerweise gegen $z$ konvergiert.

Im Falle reellquadratischer Körper stammt (2.1) von Barnes und
Swinnerton-Dyer \cite[Theorem C]{BS52a}; der Beweis, den sie für ihr
Theorem C geben, geht auf Cassels (unveröffentlicht; sh. dazu Bambah
\cite{Bam51} zurück, und ist viel undurchsichtiger als der folgende:

\begin{proof}[Bew. von (2.1)]
Sei $| \cdot |_1$ eine Bewertung von $K$ mit $1 < |u|_1$ und $x_0 \in
F$ ein Ausnahmepunkt. Dann muss auch $x_1 = ux_0 - a$ wieder ein
Ausnahmepunkt sein (sonst läge mit $x_1$ auch $x_0 = (x_1 - a)u^{-1}$
in bedecktem Gebiet), und nach Voraussetzung ist dann $x_1 \in
S$. Durch die Rekursion $x_{i+1} = ux_i - a$ erhalten wir also eine
Folge $x_0, x_1, x_2, \dots$ von in $S$ liegenden Ausnahmepunkten. Nun
ist aber $x_{i+1} - z = u(x_i - z)$ für alle $i \in N_0$; da
alle $x_i$ in $S$ liegen, muss $|x_{i+1} - z|_1$ beschränkt sein für
alle Bewertungen $| \cdot |_1$ von $K$. Im Falle $|u|_1 > 1$ geht dies
aber nur, wenn $|x_i - z|_1 = 0$ ist für alle $i \in N_0$.

Ist dagegen $|u|_1 < 1$, so ist $|x_i - z|_1$ eine Nullfolge; falls
also $|u|_1 \neq 1$ für alle Bewertungen $| \cdot |_1$ von $K$ ist,
gilt $\lim |x_i - z|_1 = 0$ für alle $j$ mit $1 \le j \le r_0$; dies
ist aber nur möglich, wenn $\lim x_i = z$ gilt.
\end{proof}

In $K=\Q(\sqrt{14})$ liegt also jeder Ausnahmepunkt von $S$
auf der Geraden $|x-z|_1 = 0$ für diejenige Bewertung $|\,\cdot\,|_1$, für
welche $|15 + 4\sqrt{14}|_1 > 1$ ist. Um weiterzukommen, betrachten
wir statt $u$ die Einheit $u' = u^{-1} = 15 - 4\sqrt{14}$ und finden wie oben
$u'S = -25 + 5\sqrt{14} + (0.4836, 0.5164) \times (0.495625,0.504375)$.
Wendet man hierauf (2.1) an, so folgt wegen $1<|u'|_2$,
dass jeder Ausnahmepunkt von $S$ auf der Geraden $|x-z|_2=0$ durch den
Punkt $z = \frac{-21 + 5\sqrt{14}}{14-4\sqrt{14}} = \frac{1+\sqrt{14}}2$
liegt. Da sich die beiden Geraden $|x-z|_1=0$ und $|x-z|_2=0$ genau im
Punkt $z$ schneiden, ist $z = \frac{1+\sqrt{14}}2$ der einzige
Ausnahmepunkt von $S$. Bereits in \S\  1 haben wir gesehen, dass
$M(z)=5/4$ ist; wir haben also

\begin{quote}
  {\bf (2.2)} \em{Sei $K=\Q(\sqrt{14})$; dann gilt
  $M(K) = M(\uK) = 5/4$; dieses euklidische Minimum wird
  auf $F = [0,1] \times [0,1]$ genau im Punkt
  $z = (1 + \sqrt{14})/2$ angenommen; für alle $x \in F \setminus \{z\}$
  ist $M(K,x) \le 1.2$, d.h. das erste euklidische Minimum ist isoliert.}
\end{quote}

In vielen Fällen wird man -- wie oben für $\Q(\sqrt{14})$ --
sowohl $u$ als auch $u^{-1}$ zur Bestimmung von $M(K)$ verwenden
können; in diesem Fall nimmt (2.1) die folgende Form an:

\begin{quote}
{\bf (2.3)} {\em Mit den Bezeichnungen von (2.1) möge gelten:
  \begin{enumerate}
  \item es gibt ein $a \in \R$, sodass $ux-a$ für alle $x \in S$
    in bedecktem Gebiet oder wieder in $S$ liegt;
  \item für alle $x \in S$ gibt es ein $b \in \R$, sodass
    $u^{-1}x-b$ in bedecktem Gebiet oder wieder in $S$ liegt;
  \item es ist $|u|_j \ne 1$ für jede Bewertung $||_j$ von $K$;
  \end{enumerate}
  dann ist $z = \frac{a}{u-1}$ der einzig mögliche Ausnahmepunkt von $S$. }
\end{quote}

Bem.: Für reellquadratische $K$ stammt dies von Cassels und ist
erstmals von Bambah \cite{Bam51} veröffentlicht worden; für den
Originalbeweis gilt dieselbe Bemerkung wie für denjenigen von 2.1. Wie
man leicht einsieht, lässt sich in Körpern mit Einheitenrang 1 die
Bedingung 3. durch die Forderung \glqq{}$u$ ist keine Einheitswurzel\grqq{}
ersetzen.

\begin{proof}
  Sei $x \in S$ ein Ausnahmepunkt; für diejenigen Bewertungen mit
  $|u|_j > 1$ ist $|x-z|_j = 0$ wegen (2.4) und der Bedingung 1. Für
  alle andern Bewertungen ist $|u|_j < 1$ wegen 3. sodass dann $1 <
  |u^{-1}|_j$ gilt. Wir stellen nun fest, dass mit $x$ auch $x' =
  u^{-1}x - b \in S$ ein Ausnahmepunkt sein muss; wegen $x' \in S$ und
  1. ist $ux' - a \in S$, also $x - ub - a \in S$. Da $S$ eine
  Teilmenge von $F$ und $a, ub \in \R$ sind, muss $a = -ub$ sein: also
  ist $b = -au^{-1}$ unabhängig von $x$. Das Argument von (2.4) zeigt
  nun, dass $x$ den Gleichungen $|x-z'|_j = 0$ genügt, wobei
  $$ z' = \frac{b}{u^{-1}-1} = \frac{bu}{1-u} = \frac{a}{u-1} = z $$
  gilt. Daher ist $|x-z| = 0$ für alle $j$ mit $1 \le j \le r+s$, und
  folglich ist $x-z$ der einzig mögliche Ausnahmepunkt.
\end{proof}

Mit etwas mehr Mühe können wir auch noch das zweite euklidische
Minimum von $\Q(\sqrt{14})$ bestimmen (Bedocchi \cite{Bed85} zeigte
1985, dass $M_2(K) \le 1$ gilt; van der Linden behauptet auf
\cite[S. 25]{Lin85}, dass $M_2(K) < 1$ ist, gibt dafür aber keinen
Beweis). Wir wählen dazu $k=0.96$ und können ganz $F_1 =
(-0.1,0.9)\times (-0.1,0.9)$ bedecken bis auf die folgenden Gebiete:
\begin{align*}
   S^+ & = (-0.0005, 0.0005) \times (0.37475, 0.37525), \\
   S^- & = (-0.0005, 0.0005) \times (-0.37525, -0.37475), \\
   S & = (0.499, 0.501) \times (0.499975, 0.500025);
\end{align*}
mit $a = 35 + 9\sqrt{14}$ und $b = -21 + 5\sqrt{14}$ liegen $uS-a$ und
$u^{-1}S-b$ in bedecktem Gebiet oder wieder in $S$ (insbesondere haben
beide Mengen mit $S^+$ und $S^-$ keinen Punkt gemeinsam); nach (2.3)
ist daher $z = (1+\sqrt{14})/2$ der einzige Ausnahmepunkt von $S$. Nun
finden wir
\[
uS^+ = 21 + 6\sqrt{14} + (-0.0215,0.0215) \times (-0.38075, -0.36925)
\text{ und}
\]
\[
uS^- = -21 - 6\sqrt{14} + (-0.0215,0.0215) \times (0.36925, 0.38075).
\]
Wir sehen also: für $a=21+6\sqrt{14}$ und $x\in S^+$ liegt $ux-a$
entweder in bedecktem Gebiet oder in $S^-$, während für $x\in S^-$ die
Punkte $ux-a$ in bedecktem Gebiet oder in $S^+$ liegen. Dies führt uns
auf folgende Verallgemeinerung von (2.1):

\begin{quote}
  {\bf (2.4)} {\em Sei $K$ Zahlkörper, $u$ Einheit, $S_1, \dots, S_t$
    Teilmengen von $F$ und $a_1, \dots, a_t \in R$ Punkte mit der
    Eigenschaft: für alle $x\in S_1$ liegt $ux-a_1$ in bedecktem
    Gebiet oder in $S_{i+1}$ (hierbei setzen wir $S_{t+1}=S_1$). Ist
    dann $x_1\in S_1$ ein Ausnahmepunkt, so gilt $|x_1-z|_i = 0$ für
    jede Bewertung von $K$ mit $|u_i| > 1$, wobei $z = a/(u^t-1)$ mit
    $a = u^{t-1}a_1 + u^{t-2}a_2 + \dots + u a_{t-1} + a_t$ Fixpunkt
    der Abbildung $x \to ux - a$ ist. }
\end{quote}

\begin{proof}
  Sei $x_1\in S_1$ ein Ausnahmepunkt; dann ist auch $x_2 = ux_1 - a_1
  \in S_2$ ein solcher, und wir finden nacheinander $x_3 = ux_2 - a_2
  \in S_3, \dots, x_{t+1} = ux_t - a_t \in S_{t+1} = S_1$. Indem wir
  von hier aus zurückrechnen, erhalten wir $x_{t+1} = u^t x_1 - a$ für
  $a = u^{t-1}a_1 + u^{t-2}a_2 + \dots + u a_{t-1} + a_t$; damit haben
  wir ein $a \in R$ gefunden mit der Eigenschaft, dass für jeden
  Ausnahmepunkt $x_1 \in S_1$ auch $u^t x_1 - a \in S_1$ ist für $u^t
  = u$. Mit (2.1) folgt jetzt, dass jeder Ausnahmepunkt $x \in S_1$ die
  Gleichung $|x-z|_i = 0$ erfüllt, wobei $z = a/(u^t-1)$ Fixpunkt der
  Abbildung $x \to u^t x - a$ und $|u^t|_i$ eine Bewertung mit $1 <
  |u^t|_i$ ist. Die Behauptung folgt.
\end{proof}

Auch dieses Ergebnis stammt im reellquadratischen Fall von Barnes und
Swinnerton-Dyer, eine entsprechende Verallgemeinerung von (2.3) von
Cassels. In $\Q(\sqrt{14})$ ist $a_1 = 21 + 6\sqrt{14}$, $a_2 = -a_1$,
also $a = u a_1 - a_1 = a_1(u-1)$ und damit $z = a/(u_2-1) =
a_1/(u_1-1) = 3\sqrt{14}/8 = (0, 3/8)$ (wir wollen in Zukunft statt
$a+b\sqrt{m}$ einfach $(a,b)$ schreiben). Durch Verwenden von $u^{-1}$
erhält man nun leicht, dass dieses $z$ der einzig mögliche
Ausnahmepunkt von $S_1$ ist; damit ist $z = u z - a_1 = -z$ der einzig
mögliche Ausnahmepunkt von $S_2$. Jetzt benutzen wir (1.8) mit $t=2$,
$x_1=z$, $x_2 = z$ und finden $M(K,z) = M(K,z') = 31/32 = 0.96875$;
damit haben wir folgende Verschärfung von (2.2):

\begin{quote}
  {\bf (2.2')} {\em Sei $K = \Q(\sqrt{14})$: dann ist
    $M(K) = M_1(K) = \frac54$ das erste euklidische Minimum von $K$;
    dieses ist isoliert und wird $\bmod\ R$ genau im Punkt
    $(\frac12, \frac12)$ angenommen. Weiter ist
    $M_2(K) = M_2(\uK) = \frac{31}{32}$ das zweite
    euklidische Minimum von $K$; auch dieses ist isoliert und wird
    $\bmod\ R$ genau in den Punkten $(0, \frac38)$ angenommen. }
\end{quote}

Nicht immer kann man die euklidischen Minima so einfach bestimmen wie
in $\Q(\sqrt{14})$: ein etwas komplizierteres Beispiel ist
$K = \Q(\sqrt{13})$. Hier wählen wir die $\Q$-Basis $\{1, \sqrt{13}\}$
und den \glqq{}Fundamentalbereich\grqq{}
$F = (0, \frac12) \times (0, \frac12)$ (wir dürfen uns hierauf
beschränken, weil es zu jedem $x \in \uK$ ein $y \in R$ gibt mit
$x - y \in F \cup -F \cup F^\sigma \cup -F^\sigma$;
hier steht $\sigma$ für den nichttrivialen Automorphismus
$\sigma: \sqrt{13} \mapsto -\sqrt{13}$). Mit $k=0.333$ können wir ganz
$F$ bis auf die folgenden Gebiete
\begin{align*}
  S_1 & = (-0.001, 0.001) \times (0.2128, 0.2129) &
  S_4 & = (0.161, 0.163) \times (0.1678, 0.1682) \\
  S_2 & = (0.116, 0.117) \times (0.1806, 0.1808) &
  S_5 & = (0.165, 0.166) \times (0.1668, 0.1672) \\
  S_3 & = (0.151, 0.152) \times (0.1707, 0.1711) &
  S_6 & = (0.166, 0.167) \times (0.1665, 0.1669),
\end{align*}
sowie die daraus durch Multiplikation mit $-1$ und Anwenden des
Automorphismus $\sigma$ entstehenden Gebiete bedecken. Mit $u =
\frac{3+\sqrt{13}}{2}$ finden wir nun $-uS_1 = -\frac{3+\sqrt{13}}{2}
\times (0.11465, 0.1168) \times (0.18015, 0.1808)$, d.h.
$-uS_1 + u$ liegt in bedecktem Gebiet oder in $S_2$; weiter liegt
$-uS_2 + u$ in bedecktem Gebiet oder in $S_3$, $-uS_3 + u$ in bedecktem
Gebiet oder in $S_4$ usw., und $-uS_6 + u$ in bedecktem Gebiet
oder in $S_6$. Mit $S = S_1 \cup S_2 \cup \dots \cup S_6$ gilt also:
$-uS + u$ liegt in bedecktem Gebiet oder wieder in $S$. Nach (2.1)
liegt jeder Ausnahmepunkt $x$ von $S$ auf der Geraden $|x-z|_1 = 0$
durch den Punkt $z = \frac{u}{u+1} = (\frac16, \frac16)$.

Jetzt benutzen wir $u' = \frac{3-\sqrt{13}}{2}$ und setzen $T_1=S_1$, sowie
\begin{align*}
  T_2 & = (-0.117, -0.116) \times (0.1806, 0.1808) \\
  T_6 & = (-0.167, -0.166) \times (0.1665, 0.1669)
\end{align*}
und $T = T_1 \cup T_2 \cup \dots \cup T_6$. Wie oben findet man nun,
dass die Punkte $-u'T - u'$ in bedecktem Gebiet oder wieder in $T$
liegen. Also liegt jeder Ausnahmepunkt $x$ von $T$ auf der Geraden
$|x-z'|_2 = 0$ durch den Punkt $z' = -\frac{u'}{u'+1} = (- \frac16, \frac16)$.

Nun ist $T \cap S = S_1 = T_1$, sodass jeder Ausnahmepunkt von $S_1$
auf den Geraden $|x-z|_1 = 0$ und $|x-z'|_2 = 0$ liegen muss. Diese
beiden Geraden schneiden sich im Punkt $x_1 = (0, \frac16 + \frac16\sqrt{13})
\in K \setminus \Q$, sodass $x_1$ der einzig mögliche
Ausnahmepunkt von $S_1$ ist. Weiter sieht man sofort, dass
$x_2 = -u x_1 +u$ der einzig mögliche Ausnahmepunkt von $S_2$ ist (denn
$uS_2$ liegt $\bmod R$ in bedecktem Gebiet oder in $S_1$). Die
Rekursion $x_{i+1} = u x_i - a$ für $u = -\frac{3+\sqrt{13}}{2}$
und $a=-u$ liefert nun die Folge von möglichen Ausnahmepunkten
\[
x_{k+1} = \Bigl(\frac{1}{6} - \frac{1}{6} e^k,
               \frac{1}{6} + \frac{e^k}{6\sqrt{13}}\Bigr)
\]
für $k = 0, 1, 2, \dots$ und $e = (-3+\sqrt{13})/2$.

Um zu zeigen, dass dadurch alle möglichen Ausnahmepunkte von $S$
gegeben sind, nehmen wir an, $x$ sei ein solcher. Da $x$ auf der
Geraden $|x-z|_2 = 0$ liegt, gibt es ein $k \in \N$ mit
$|x_{k+1} - z|_2 < |x - z|_2 < |x_k - z|_2$ (d.h. $x$ liegt
\glqq{}zwischen\grqq{} zwei Gliedern der Folge $\{x_i\}$). Dann liegt aber
$(-u)^{-k-1} x \bmod R$ zwischen $x_1$ und $x_2$, was offenbar nicht
der Fall ist.

Schließlich zeigen wir noch, dass $M(x_1) = \frac{1}{3}$ ist und dieses
Minimum nicht erreicht wird. Dazu beachten wir, dass erstens $M(x_1) =
M(x_2) = \dots = M(x_k) = \dots$ ist, sodass aus $\lim N(x_i) = N(\lim
x_i) = N(z) = |N_{K/\Q}(z)| = \frac{1}{3}$ sicher $M(x_k) \ge
\frac{1}{3}$ folgt. Gäbe es nun ein $y \in R$ mit $N(x_k - y) <
\frac{1}{3}$, so folgte aus (1.8) die Existenz eines $z = x_k \bmod R$
mit $z = r + s\sqrt{m}$, $|s| < 0.18957$ und $N(z) <
\frac{1}{3}$. Man sieht aber leicht ein, dass es ein solches $z$ nicht
gibt. Also ist $M(x_k) = \frac{1}{3}$; weil jedoch $N(x_k - y)$ für
alle $y \in R$ irrational ist, kann nicht $N(x_k - y) = \frac{1}{3}$
sein, und dies zeigt, dass $M(x_k)$ nicht erreicht wird.

Wir geben nun noch die numerischen Werte für die Punkte $x_1, \dots,
x_7$ an, um einen Vergleich zwischen den Ausnahmepunkten und den
Mengen $S_1, \dots, S_6$ zu ermöglichen:
\begin{align*}
x_1 &= (0, 0.212891683\dots) \\
x_2 &= (0.116204060\dots, 0.180662475\dots) \\
x_3 &= (0.151387818\dots, 0.170904256\dots) \\
x_4 &= (0.162040603\dots, 0.167949705\dots) \\
x_5 &= (0.165266007\dots, 0.167055139\dots) \\
x_6 &= (0.166242581\dots, 0.166784286\dots) \\
x_7 &= (0.166538263\dots, 0.166702279\dots) \quad \text{usw.}
\end{align*}

Noch größere Schwierigkeiten bereitet die Bestimmung eines
euklidischen Minimums, wenn dieses nicht isoliert ist. Dass es
reell-quadratische Zahlkörper gibt, deren zweites Minimum $M_2(K)$
schon nicht mehr isoliert ist, haben Barnes und Swinnerton-Dyer
vermutet; Godwin hat diese Vermutung dann 1963 für $K =
\Q(\sqrt{23})$ bewiesen. Wahrscheinlich ist es mit den hier
vorgestellten Methoden möglich zu zeigen, dass $M_2(K)$ für $K =
\Q(\sqrt{m})$, $m=(2n+1)^2-2$, $n \ge 2$, ein nicht isoliertes
zweites Minimum hat. Barnes und Swinnerton-Dyer haben gezeigt, dass in
diesem Fall gilt: es ist $M_1(K) = (8n^3 + 6n^2 - 6n + 1)/2m$; dieses
Minimum ist isoliert und wird $\bmod R$ genau in den Punkten
$\bigl(\frac{1}{2}, \frac{1}{2} - \frac{1}{m}\bigr)$ angenommen.
Entsprechend lässt sich Godwins Ergebnis verallgemeinern, indem man
zeigt, dass
$$M _2(K) \ge \frac{2n(2n+1)\sqrt{m} - (2n^2 + m)}{2m} $$
gilt, und dass im Falle von Gleichheit dieses Minimum nicht isoliert
ist (im Falle $m=7$ erhält man übrigens $M_3(K) \ge \frac{6\sqrt{7}-9}{14}$).

In der nachstehenden Tafel geben wir die euklidischen Minima $M_1(K)$
und $M_2(K)$ für $K = \Q(\sqrt{m})$, $m \le 102$, an, soweit
sie bekannt sind. Die erste Tabelle dieser Art findet sich bei
Heinhold \cite{Hei39}; diese wurde dann von Barnes und Swinnerton-Dyer
ergänzt. Die dort noch fehlenden ersten Minima wurden schließlich von
Godwin \cite{God55} bestimmt.

Die zweiten Minima, die in dieser drei Arbeiten noch nicht errechnet
wurden, sind in der Tabelle kursiv gesetzt. Die kleinsten $m$, für die
$M_2(K)$ noch unbekannt ist, sind jetzt $m = 19$ und $m=22$ für $m =
2,3 \bmod 4$, sowie $m = 41$ für $m = 1 \bmod 4$; in diesen Fällen sind
in der Tafel untere Schranken für $M_2(K)$ angegeben.

Es sei noch bemerkt, dass die Tafeln von Barnes und Swinnerton-Dyer
einige kleinere Druckfehler enthalten; so ist das erste Minimum von
$\Q(\sqrt{19})$ fälschlich mit $\frac38$ angegeben; diesen Fehler
haben Barnes und Swinnerton-Dyer im vierten Teil ihrer Arbeit selbst
berichtigt. Außerdem muss dann \glqq{}$C_1 = (\frac12, \frac38)$\grqq{}
zu \glqq{}$C_1 = (0, 5/9)$\grqq{} geändert werden.

Weiter sind die Mengen $C_2$, an denen $M_2$ angenommen wird, für die
Werte $m = 17, 35, 37, 65, 101$ nicht richtig. Für die Werte
$m = 17$, $37$, $65$, $101$ ist jeweils ein Vorzeichen beim zweiten
Punkt in $C_2$ zu ändern, während für $m = 15$ (bzw. $m = 35$)
\glqq{}$C_2 = (0, \frac12)$\grqq{} statt
\glqq{}$C_2 = (\frac12, \frac18)$\grqq{} (bzw.
\glqq{}$C_2 = (0, \frac12)$\grqq{} statt
\glqq{}$C_2 = (\frac12, \frac12)$\grqq{}) stehen muss.

\begin{table}[h]
\centering
\begin{tabular}{llll}
\toprule
$m$ & $M_1$ & $M_2$ & Bem. \\
\midrule
5 & $1/4$ & $1/5$ & sh. Davenport 1946 \\
13 & $1/3$ & $4/13$ & \\
17 & $1/2$ & $8/17$ & \\
21 & $5/7$ & $12/23$ & \\
29 & $4/5$ & $23/29$ & \\
33 & $29/44$ & $6/11$ & \\
37 & $3/4$ & $27/37$ & \\
41 & $23/32$ & $\ge 2\,2819/4100$ & \\
53 & $9/7$ & $68/53$ & \\
57 & $14/19$ & $219/304$ & $M_1$: Godwin 1955 \\
61 & $1611/1525$ & $41/39$ & \\
65 & $1$ & $64/65$ & \\
69 & $25/23$ & & \\
73 & $1541/2136$ & irrational; sh. Godwin 1955 & \\
77 & $19/11$ & $16/11$ & \\
85 & $16/9$ & $151/85$ & \\
89 & $1004287/1000004$ & & \\
93 & $44/31$ & & \\
97 & $33679354/31404817$ & & \\
101 & $5/4$ & $125/101$ & \\
\bottomrule
\end{tabular}
\caption{Euklidische Minima reellquadratischer Zahlkörper
  ($m \equiv 1 \bmod 4$)}\label{dApp17}
\end{table}

\begin{table}[h]
\centering
\begin{tabular}{llll}
\toprule
$m$ & $M_1$ & $M_2$ & Bemerkungen \\
\midrule
2 & $1/2$ & $1/4$ & sh. Varnavides \cite{Var48b} \\
3 & $1/2$ & $1/3$ & \\
6 & $3/4$ & $1/2$ & $M_3(K) \ge (6\sqrt{7}-9)/14$ \\
7 & $9/14$ & $1/2$ & \\
10 & $3/2$ & $39/40$ & \\
11 & $19/22$ & $3125/3971$ & \\
14 & $5/4$ & $31/32$ & \\
15 & $3/2$ & $7/5$ & \\
19 & $170/171$ & $\ge 10579/12844$ & \\
22 & $27/22$ & $\ge 175903/155236$ & \\
23 & $77/46$ & $(20\sqrt{23}-31)/46$ & Godwin 1963 \\
26 & $5/2$ & $207/104$ & \\
30 & $3/2$ & $29/20$ & \\
31 & $45/31$ & $3775/3038$ & \\
34 & $9/4$ & $135/68$ & $M_3 = 137/72$ \\
35 & $5/2$ & $17/7$ & \\
38 & $11/4$& $173/72$ & \\
39 & $5/2$ & $107/48$ & \\
42 & $7/4$ & $41/24$  & \\
43 & $11829/6962$ & $5902/3483$ & Godwin 1955 \\
46 & $79877/48668$ & & \\
47 & $253/94$ & $\ge (42\sqrt{47}-65)/94$ & \\
\bottomrule
\end{tabular}
\caption{$m \equiv 2,3 \bmod {4}$}
\end{table}

\begin{table}[h]
\centering
\begin{tabular}{llll}
\toprule
$m$ & $M_1$ & $M_2$ & Bemerkungen \\
\midrule
51 & $287/102$ & & \\
55 & $9/4$ & $351/176$ & \\
58 & $3/2$ & $27477/19604$ & \\
59 & $125/59$ & & $M_3 = 367/128$ \\
62 & $13/4$ & $367/124$ & \\
66 & $15/4$ & $431/128$ & \\
67 & $341/162$ & & Godwin 1955 \\
70 & $891/500$ & & Godwin 1955 \\
71 & $7393/3479$ & & Godwin 1955 \\
74 & $5/2$ & & \\
78 & $7/2$ & & \\
79 & $585/158$ & $\ge (72\sqrt{79}-111)/158$ & \\
82 & $9/2$ & $1311/328$ & \\
83 & $631/166$ & & Godwin 1955 \\
86 & $10030/5203$ & & Godwin 1955 \\
87 & $169/58$ & & \\
91 & $5/2$ & & \\
94 & $4708623/2143294$ & & Godwin 1955 \\
95 & $7/2$ & & \\
102 & $19/4$ & $869/200$ & \\
\bottomrule
\end{tabular}
\caption{$m \equiv 2,3 \bmod {4}$}
\end{table}

Die Abschätzung $M(K) \le \sqrt{d}/4$ für die euklidischen Minima
$M(K)$ reellquadratischer Zahlkörper stammt von Minkowski
(sh. Ges. Abh. 1, 320--356 oder auch Hardy \& Wright, \S\  24) und ist
bestmöglichst: Letzteres folgt aus

\begin{quote}
  {\bf (2.5)} {\em Sei $n \in \N$ ungerade und $m=n^2+1$
    quadratfrei. Dann ist für den reellquadratischen Zahlkörper
    $K = \Q(\sqrt{m})$ $M(K) = \frac{n}{2}$, und dieses Minimum
    wird $\bmod R$ genau in $z = \frac{\sqrt{m}}{2}$ angenommen. }
\end{quote}

Wegen $\disc K = 4m$ in diesem Fall ist
$\frac{\sqrt{d}}4 = \frac{\sqrt{n^2+1}}2$, d.h. man kann
$\frac{\sqrt{d}}4 - M(K)$ beliebig klein machen, indem man $n$
genügend groß wählt.

Zum Beweis von (2.5) wählen wir die GHB $\{1, \vartheta\}$ mit
$\vartheta = n+\sqrt{m}$. Dann müssen wir zu jedem $z \in F$ ein $a
\in R$ finden mit $|N(z-a)| \le \frac{n}{2}$ und zeigen, dass
Gleichheit nur für $z = \sqrt{m}/2$ auftreten kann. Wir behaupten, dass
(bei dieser geschickten Wahl von $\vartheta$) sogar $a=0$ genügt. Dazu
setzen wir $z = x+y\vartheta$, fassen die Norm $N$ als Funktion von
$x$ und $y$ auf und fragen uns, wo $N$ auf $F$ ein Maximum oder ein
Minimum annehmen kann. Im Innern von $F$ ist dies nicht möglich, denn
wegen
\[
N(x,y) = N(x+y\vartheta) = x^2 + 2nx - y^2,
\]
\[
\frac{\partial N}{\partial x} = 2x + 2ny = 0 \quad \Rightarrow \quad x
= -ny
\]
und
\[
\frac{\partial N}{\partial y} = 2nx - 2y = 0 \quad \Rightarrow \quad y
= nx
\]
ist $z = 0$ der einzige Punkt, in dem beide partiellen Ableitungen
verschwinden (in $z$ wird weder ein Maximum, noch ein Minimum
angenommen: $z$ ist nämlich Sattelpunkt von $N$). Folglich wird das
Minimum und Maximum von $N$ auf dem Rand von $F$ angenommen. Sei nun
$y=y_0 - \frac{1}{2}$; dann ist $N(x, \frac{1}{2}\vartheta) = x^2 + nx
- \frac{1}{4}$, $\frac{\partial N}{\partial x} = 2x+n$, und folglich
wird $N$ auf der Geraden $y=y_0$ nur in $x_0 = -\frac{n}{2}$
extremal. Für $n > 2$ liegt $x_0$ aber nicht in $F$, für $n=1$ ist
$x_0 = -\frac{1}{2}$, und wir sehen (nach ähnlichen Rechnungen für
$y_0 = -\frac{1}{2}$, bzw. $x_0 = \frac{1}{2}$), dass das Maximum von
$|N|$ auf $F$ höchstens in den Eckpunkten $z_e = (\pm\frac{1}{2},
\pm\frac{1}{2})$ angenommen wird. Dort gilt aber $N(z_e) = \frac{n}{2}$,
und die Behauptung folgt wegen $(1+\vartheta)/2 \equiv \sqrt{m}/2 \bmod R$.

Zu zeigen bleibt noch $M(K, x) = \frac{n}{2}$ für $x =
\sqrt{m}/2$. Dazu verwenden wir (1.8) und nehmen an, es sei $M(K, x)
\le \frac{n}{2}$. Dann gibt es ein $z = x \bmod R$ mit $z =
r+s\sqrt{m}$, $|s| < \sqrt{k(n+1)/2m} < \frac{n}{2}$ und $|N(z)| <
k$. Die Bedingungen an $s$ lassen nur $s = \pm\frac{1}{2}$ zu, und aus
Symmetriegründen dürfen wir $s = \frac{1}{2}$ annehmen. Damit ist aber
$N(r+s\sqrt{m}) = |r^2 - \frac{n}{4}|$ und dies wird offenbar minimal,
wenn $|r| \le \frac{n}{2}$ ist. Einsetzen von $r_0 = \frac{n-1}{2}$
und $r_1 = \frac{n+1}{2}$ liefert $N(r_0 + s\sqrt{m}) = N(r_1 +
s\sqrt{m}) = \frac{n}{2}$, und das war zu zeigen.

Ist $n$ gerade und $m=n^2+1$ quadratfrei, so gilt (2.5) in dieser Form
nicht mehr, weil $\{1, \vartheta\}$ keine GHB mehr ist. Dagegen bleibt
(2.5) auch in diesem Fall richtig, wenn man den Ring ganzer Zahlen
$D(m)$ ersetzt durch den Ring $\Z[\sqrt{m}]$; insbesondere ist
$M(f)=1$ für $\Z[\sqrt{5}]$ und $f=|N|$ (sh. die Beispiele auf
S. \pageref{pBei}).

Mit einer etwas andern Methode lässt sich (2.5) noch verschärfen: ist
nämlich $z \in F$ ein Ausnahmepunkt für
$$ k = \frac{n-1}{2} - \frac{1}{4m}, $$
so ist notwendig
\[
|N(x+y\delta)| = |x^2 + 2nxy - y^2| \ge k.
\]
Wegen $|x|, |y| \le \frac{1}{2}$ ist aber $|x^2 - y^2| \le
\frac{1}{4}$ und $|2nxy| \le n|y|$, also
\[
|N(x+y\delta)| \le |x^2 - y^2| + |2nxy| \le n|y| + \frac{1}{4}.
\]
Beide Ungleichungen zusammen ergeben $n|y| \ge k - \frac{1}{4}$ für
jeden Ausnahmepunkt $z \in F$, also
\[
|y| \le \frac{1}{2} - \delta, \quad \delta = \frac{3}{4n} + \frac{1}{4mn}.
\]
Dies ist nur für $-\frac{1}{2} \le y \le -\frac{1}{2} + \delta$ oder
$\frac{1}{2} - \delta \le y \le \frac{1}{2}$ möglich.

Indem wir $F = (-\frac{1}{2}, \frac{1}{2}) \times (-\frac{1}{2},
\frac{1}{2})$ (wie in $\Q(\sqrt{14})$) durch $F' = (0,1) \times
(0,1)$ ersetzen, erkennen wir, dass jeder Ausnahmepunkt in dem Streifen
$|y - \frac{1}{2}| \le \delta$ liegen muss. Entsprechend folgt nun aus
\[
|N(z)| \le n|x| + \frac{1}{4},
\]
dass für einen Ausnahmepunkt $z = x + y\delta \in F'$ auch $|x -
\frac{1}{2}| \le \delta$ gelten muss. Also liegt jeder Ausnahmepunkt $x
+ y\delta$ in dem durch
\[
\Big| x - \frac{1}{2} \Big| \le \delta, \quad
\Big| y - \frac{1}{2} \Big| \le \delta
\]
definierten Rechteck um den Punkt $(1+\delta)/2$, an dem, wie wir oben
gesehen haben, $M_1(K)$ angenommen wird. Da nun $x$ und $y$ beide in
der Nähe von $\frac{1}{2}$ liegen, können wir $|x^2 - y^2|$ etwas
besser nach oben abschätzen als durch $\frac{1}{4}$; für $x,y \in F$
ist nämlich $|x - y| \le \delta$, $|x + y| \le 1$ oder $|x - y| \le
1$, $|x + y| \le \delta$. Also ist
\[
|x^2 - y^2| = |x - y||x + y| \le \delta,
\]
und wie oben folgt nun $n|x| \ge k - \delta$ und $n|y| \ge k - \delta$.
Weiter ersehen wir hieraus, dass jeder Ausnahmepunkt $z \in F'$ in dem durch
\[
\Big|x - \frac{1}{2}\Big| \le \delta_1, \quad
\Big|y - \frac{1}{2}\Big| \le \delta_1
\]
definierten Rechteck $S$ liegt, wobei
\[
\delta_1 = \frac{1}{2n} + \frac{1}{n^2}
\]
ist.

Auch wenn wir diesen Schritt unendlich oft wiederholen, können wir
$\delta_1$ nicht mehr wesentlich verbessern. Wir werden daher nur die
Einheiten $u = \vartheta$ und $u' = \vartheta' = 2n - \vartheta$ verwenden, um
weiterzukommen. Wir finden
\[
z \cdot u = (x + y\vartheta)\vartheta = y \cdot (x + 2ny)\vartheta
\]
und
\[
-z \cdot u' = -(x + y\vartheta)\vartheta' = y - 2nx + x\vartheta.
\]
Da mit $z$ auch $z\vartheta$ Ausnahmepunkt sein muss, existiert ein
$a = r + s\vartheta \in R$ mit $z\vartheta - a \in S$, d.h. mit
\[ \Big|x + 2ny - s - \frac{1}{2}\Big| \le \delta_1; \]
wegen $|x - \frac{1}{2}| \le \delta_1$ muss $|2ny - s| \le 2\delta_1$
sein. Andererseits folgt aus $|y - \frac{1}{2}| \le \delta_1$, dass wir
\[ |2ny - n| \le 2n\delta_1 = 1 + \frac{2}{n} \]
haben, und folglich ist $s \in \{n-1, n, n+1\}$. Dividiert man
$|2ny - s| \le 2\delta_1$ durch $2n$, so erhält man
\[ \Big|y - \frac{1}{2n}\Big| \le \frac{\delta_1}{n} = \eps, \]
und eine entsprechende Rechnung mit $\delta'$ anstelle von $\delta$
zeigt, dass auch $|x - \frac{1}{2n}| \le \eps$ gilt, wo wieder
$r \in \{n-1, n, n+1\}$ ist.

Durch Kombination aller Möglichkeiten erhält man jetzt die folgenden
neun möglichen Ausnahmemengen:
\[ S_{i,k}: \quad \Big|x - \frac{n+i}{2n}\Big| \le \eps, \quad
                 \Big|y - \frac{n-k}{2n}\Big| \le \eps, \]
wo $i$ und $k$ unabhängig voneinander die Zahlen $-1, 0, +1$
durchlaufen. Von diesen neun Mengen enthalten aber $S_{i,k}$ für
$(i,k) = (-1,-1), (-1,+1), (+1,-1), (+1,+1)$ keine Ausnahmepunkte,
falls nur $n \ge 5$ ist: ist z.B. $|x + \frac{n-1}{2n}| \le \eps$
und $|y - \frac{1}{2n}| \le \eps$, so findet man $|x^2 - y^2| \le
(n+1)/n^3$ und $2n|x|y| \le (n-1) \cdot 2\delta_1^2/2n$, und für $n
\ge 5$ folgt daraus $|N(z)| < k$. Damit sind noch fünf Ausnahmemengen
verblieben, nämlich $S_{i,k}$ für
$(i,k) = (-1,0), (+1,0), (0,0), (0,-1), (0,+1)$, das sind
$$ \begin{array}{lll} \smallskip 
  S_1 := S_{-1,0} : & \Big|x - \frac{n-1}{2n}\Big| \le \eps, &
                           \Big|y - \frac{1}{2}\Big| \le \eps, \\ \smallskip 
  S_2 := S_{+1,0} : & \Big|x - \frac{n+1}{2n}\Big| \le \eps, &
                           \Big|y - \frac{1}{2}\Big| \le \eps, \\ \smallskip 
  S_3 := S_{0,0} :  & \Big|x - \frac{1}{2}\Big| \le \eps, &
                           \Big|y - \frac{1}{2}\Big| \le \eps, \\ \smallskip 
  S_4 := S_{0,-1} : & \Big|x - \frac{1}{2}\Big| \le \eps, &
                           \Big|y - \frac{n+1}{2n}\Big| \le \eps, \\
  S_5 := S_{0,+1} : & \Big|x - \frac{1}{2}\Big| \le \eps, &
                           \Big|y - \frac{n-1}{2n}\Big| \le \eps.
\end{array} $$
Nun sieht man sofort, dass $-\vartheta' S_1 \bmod R$ in bedecktem Gebiet
oder in $S_4$ liegt (für $x+y\vartheta \in S_1$ ist ja  $s=x$ mit
$r+s\vartheta = -(x+y\vartheta)\vartheta'$); dies impliziert
\[ \Big|y-2nx-(n-1)-\frac{1}{2}\Big| \le \eps, \]
wegen $|y-\frac{1}{2}|\le \eps$ folgt
\[ |2nx-(n-1)|\le \frac{2\eps}{n^2} \]
und somit
\[ \Big|x-\frac{n-1}{2n}\Big| \le \frac{\delta}{n^2}. \]
Ganz analoge Überlegungen führen zu
$|x-\frac{n-1}{2n}|\le \frac{\delta}{n^2}$ für $z\in S_2$, bzw. zu
$|y-\frac{n-1}{2n}|\le \frac{\delta}{n^2}$ für $z\in S_4, S_5$.

Jetzt betrachten wir $\vartheta S_1$: dieses Gebiet liegt $\bmod R$
entweder in bedecktem Gebiet oder in einer der drei Mengen $S_3, S_4,
S_5$. Liegt $(x+y\vartheta) \bmod R$ in $S_3$, so muss
\[ \Big|x+2ny-n-\frac{1}{2}\Big| \le \eps \]
sein; wegen $|x-\frac{n-1}{2n}|\le \eps$ folgt
$$ \Big|2ny-n-\frac{1}{2n}\Big|\le 2\eps, \quad \text{also} \quad 
   \Big|y-\frac{1}{2}-\frac{1}{4n^2}\Big |\le \frac{\eps}{n}. $$
Entsprechend führt die Möglichkeit, dass $(x+y\vartheta)\vartheta \bmod R$
in $S_4$ liegt auf die Bedingung $|y-\frac{1}{2}|\le \eta$, und
$(x+y\vartheta)\vartheta \bmod R \in S_5$ liefert schließlich
\[ \Big|y-\frac{1}{2}-\frac{1}{2n^2}\Big|\le \eta. \]
Wir behaupten nun, dass die Mengen
$$ \Big|x-\frac{n-1}{2n}\Big|\le \eta, \quad
   \Big|y-\frac{1}{2}\Big|\le \eta \quad \text{ sowie }
   \Big|x-\frac{n-1}{2n}\Big|\le \eta, \quad
   \Big|y-\frac{1}{2}-\frac{1}{4n^2}\Big| \le \eta $$
keinen Ausnahmepunkt enthalten können. Wenn wir dies gezeigt haben,
folgt, dass $\vartheta S_1 - n\vartheta$ nur in bedecktem Gebiet oder
in $S_5$ liegen kann. Weiter liegt $\vartheta S_3 - (n+1)\vartheta$ in
bedecktem Gebiet oder in $S_2$. Nun ist $S_2 = 1+\vartheta - S_1$,
folglich liegt $S_5 - 1 - (n+2)\vartheta$ in bedecktem Gebiet oder in
$-S_1$. Damit können wir (2.4) anwenden mit $t=4$, $a_1 = n\vartheta$,
$a_2 = 1+(n+2)\vartheta$, $a_3 = -a_1$, $a_4 = -a_2$, $u=\vartheta$
und finden
\[ z_1 = \frac{u^3 a_1 + u^2 a_2 - u a_1 - a_2}{u^4 - 1}
       = \frac{u a_1 + a_2}{u^2 + 1} = \frac{m - n + (m + 1)\vartheta}{2m} \]
als einzig möglichen Ausnahmepunkt von $S_1$. Weiter sind nun
\[ z_2 = \frac{m + n + (m - 1)\vartheta}{2m}, \quad
   z_4 = \frac{m - 1 + (m - n)\vartheta}{2m}, \quad
   z_5 = \frac{m + 1 + (m + n)\vartheta}{2m} \]
die einzigen Ausnahmepunkte von $S_2, S_4$ und $S_5$.

Jetzt folgt sofort, dass $\vartheta S_3 \bmod R$ nur in bedecktem Gebiet
oder wieder in $S_3$ liegen kann, und (2.4) liefert
\[
z_3 = \frac{1+\vartheta}{2}
\]
als den einzig möglichen Ausnahmepunkt von $S_3$.

Wir zeigen nun, dass es keinen Ausnahmepunkt $z=x+y\vartheta$ mit
$|x-\frac{n-1}{2n}|\le \eta$ und $|y-\frac{1}{2}|\le \eta$ gibt. Dazu
nehmen wir oBdA $- \frac{1}{2} - \eta \le y \le -\frac{1}{2}$
an und finden
\begin{align*}
  N(z) & = x^2 + 2nxy - y^2 \le \Big(\frac{n-1}{2n} - \eta\Big)
          + n \Big(\frac{n-1}{2n} + \eta\Big) - \frac{1}{4} \\
     & = \frac{n-1}{2} - \frac{1}{2n} + \frac{1}{4n^2} +
         \eta \cdot \frac{n^3 - n^2 - n + \delta}{n^2}.
\end{align*}
Für $n \ge 5$ ist dies aber sicher $< k = \frac{n-1}{2} \cdot
\frac{1}{4m}$, sodass das obige Gebiet dann keinen Ausnahmepunkt
enthält (für \glqq{}große\grqq{} $n$ erhält man dies fast ohne
Rechnung, weil $N(z) \le \frac{n-1}{2} - \frac{1}{4n} + O(n^{-2})$ und
$k = \frac{n-1}{2} \cdot O(n^{-2})$ ist; hierbei haben wir die
bekannte \glqq{}Groß-O-Notation\grqq{} verwendet). Ganz analog rechnet
man nach, dass es auch keinen Ausnahmepunkt mit $|x-\frac{n-1}{2n}|\le
\eta$ und $|y-\frac{1}{2}-\frac{1}{4n^2}|\le \eta$ gibt (hier ist
$N(z) \le \frac{n-1}{2} - \frac{1}{4n} + O(n^{-2})$).

Jetzt müssen wir noch $M(z_i)$ für die möglichen Ausnahmepunkte $z_i$
bestimmen. Dass $M(z_3) = \frac{R}{2}$ ist, hatten wir bereits
gesehen. Wir benutzen nun (1.8) und stellen dazu die Punkte $z_i$ in
der Form $r+s\sqrt{m}$ dar; wir finden
$$ \begin{cases}
  z_1 & = \frac{m+1}{2m} \sqrt{m}, \quad
  z_2 = \frac{m-1)}{2m} \sqrt{m}, \\
  z_4 & = \frac12 + \frac{m-n}{2m} \sqrt{m}, \quad
  z_5 = \frac12 + \frac{m+n}{2m} \sqrt{m} \end{cases}
  \quad (n \equiv 1 \bmod 2), $$
$$ \begin{cases}
    z_1 & = \frac12 + \frac{m+1}{2m} \sqrt{m}, \quad
    z_2 = \frac12 + \frac{m-1}{2m} \sqrt{m}, \\
    z_4 & = \frac{m-n}{2m} \sqrt{m}, \quad
    z_5 = \frac{m+n)}{2m} \sqrt{m} \end{cases}
  \quad (n \equiv 0 \bmod 2). $$
(1.8) liefert nun wegen $a=n$ die Schranke
$\mu_2 = \sqrt{\frac{k(n+1)}{2m}} < 0.5$. Nun ist im Falle
$n \equiv 1 \bmod 2$
\[\min_{x \in \Z} |N(x + z_2)|
      = \Big|N\Big(\frac{n-1}{2} + z_2\Big)\Big|
      = \frac{n-1}{2} + \frac{1}{4m} \]
und
\[ \min_{x \in \Z} |N(x + z_4)|
   = \Big|N\Big(\frac{n-2}{2} + z_4\Big)\Big|
    = \frac{n-1}{2} - \frac{1}{4m}, \]
und eine entsprechende Rechnung für den Fall $n \equiv 0 \bmod 2$ zeigt

\begin{quote}
  {\bf (2.5')} {\em Sei $n \in \N$, $m=n^2+1$,
    $R \equiv \Z[\sqrt{m}]$ und $l$ der Absolutbetrag der
    Norm in $\Q(\sqrt{m})$. Dann ist $M_1(f) = \frac n2$ und im
    Falle $n \ge 2$ $M_2(f) = \frac{n-1}{2} - \frac{1}{4m}$.
    Diese Minima werden $\bmod R$ nur in $C_1 = \{\sqrt{m}/2\}$ und
    \begin{align*}
      C_2 & = \left\{ \frac{m \pm 1}{2m} \sqrt{m},
                  \frac{1}{2}, \frac{(m \pm n)}{2m} \sqrt{m} \right\}
             \quad \text{für } n \equiv 1 \bmod 2, \text{ und} \\
      C_2 & = \left\{ \frac{1}{2} + \frac{m \pm 1}{2m} \sqrt{m},
                     \frac{m \pm n}{2m} \sqrt{m} \right\}
             \quad \text{für } n \equiv 0 \bmod 2
    \end{align*}
    angenommen. Ist außerdem $n$ ungerade und $m$ quadratfrei, so ist
    $R = D[m]$, sowie $M_1(l) = M_1(K)$, $M_2(l) = M_2(K)$. }
\end{quote}

Die ersten Minima $M_1(l)$ der Ringe $\Z[\sqrt{m}]$ wurden
schon von Heinhold \cite{Hei39} gefunden. Analoge Sätze für $m=n^2+r$,
$r=1, 2, 4$, stammen von Barnes und Swinnerton-Dyer und von
Varnavides. Beide haben jedoch zum Beweis recht komplizierte Lemmata
benutzt (z.B. die Theoreme H, J, K von BSD), die sich in unserem
Beweis als überflüssig herausgestellt haben.

Es überrascht nun nicht, dass man mit denselben Methoden auch die Fälle
$r \in \{1, 2, 4, n\}$ behandeln kann (quadratische Zahlkörper
$\Q(\sqrt{m})$ mit $m=n^2+r$, $r|4n$, $-n \le r \le n$ heißen
Körper vom R-D-Typ nach C. Richaud und G. Degert, die die FE in diesen
Körpern explizit angegeben haben), auch wenn man die ersten beiden
Minima nicht immer so einfach bestimmen kann (wir haben bereits die
Vermutung geäußert, dass im Falle $m = n^2-2$, $n \ge 5$, $n \equiv 1
\bmod 2$, das zweite Minimum nicht isoliert ist). So gilt z.B. die
folgende Tabelle:

\section*{(2.6) Euklidische Minima für reellquadratische Körper vom R-D-Typ}

Sei $R = \Z[\sqrt{m}]$ (für $I = 1$) und $R = \Z[\frac{1+\sqrt{m}}2]$
(für $I = 2$), und $f$ der Absolutbetrag der Norm. Dann wird das erste
euklidische Minimum durch die folgende Tabelle gegeben:

$$ \begin{array}{llccc}
  \rsp   m    &  n  & I & M(K) & C(K) \\ \hline
  \rsp n^2+1 & n \ge 1 & 1 & \frac n2 & (0, \frac12) \\
  \rsp       & n \equiv 0\ (2), n \ge 1 & 2 & \frac n8 & \\
  \hline
  \rsp n^2-1 & n \equiv 0\ (2), n \ge 2 & 1 & \frac{n-1}2 & (\frac12, \frac12)
             \\ \hline
  \rsp n^2+2 & n \equiv 0\ (2), n \ge 2 & 1 & \frac{2n-1}4 & \\
  \rsp       & n \equiv 1 \bmod 2, n \ge 3 & 1 & \frac{2n^3-3n^2+6n-7}{4m}
             & \\ \hline 
  \rsp n^2-2 &  n \equiv 0\ (2), n \ge 4 & 1 & \frac{2n-1}4
             & (\frac12, \frac12) \\
  \rsp       & n \equiv 1 (2), n \ge 3 & 1 & \frac{2n^3-3n^2-6n+9}{4m}
             & \\ \hline 
  \rsp n^2+4 & n \equiv 1 (2), n \ge 3 & 2 &  \frac{n^2 - 2n + 1}{4n}
             & \\         \hline
  \rsp n^2-4 & n \equiv 1 (2), n \ge 3 & 2 & \frac{n^2-5}{4n+8} & \\ \hline
  \rsp n^2+n & n \equiv 0\ (2), n \ge 2 & 1 & \frac{n+1}4 & (0, \frac12) \\
  \rsp       & n \equiv 1\ (2), n \ge 1 & 1 & \frac{n+1}4 & (\frac12, 0) 
\end{array} $$

Darüberhinaus ist das zweite Minimum in folgenden Fällen bekannt:

$$ \begin{array}{llccc}
  \rsp   m    &  n  & I & M_2(K) & C_2(K) \\ \hline
  \rsp n^2+1 & n \equiv 0\ (2) & 1 & \frac{n-1}2 - \frac1{4m}
             & (0, \frac{m \pm 1}{2m}), (\frac12, \frac{m \pm n}{2m}) \\
  \rsp       & n \equiv 0\ (2) & 1 & \frac{n}8 \big(1 - \frac1{m}\big) & \\
  \rsp n^2-1 & n \equiv 0\ (2) & 1 & \frac{m-1}{2n+2}
             & (0, \pm \frac{n}{2n+2}) \\
  \rsp n^2+2 &  n \equiv 0\ (2) & 1 & \frac{2n^3-3n^2+4n-2}{4n} & \\
  \rsp n^2+4 &  n \equiv 1\ (2) & 2 & \frac{n^3-2n^2+5n-8}{4m}
\end{array} $$

Die Arbeit an (2.6) ist, wie man sieht, noch nicht abgeschlossen. Für
den Beweis hilfreich ist jedenfalls

\begin{quote}
  {\bf (2.7)} {\em Sei $f$ der Absolutbetrag der Norm in
    $\Q(\sqrt{m})$; dann ist }
  \[ M(f) \le \frac{1}{2} \begin{cases} \max \{p, 2n+1-p\}, & \text{wo }
    R=\Z[\sqrt{m}], \\
             & \quad m=n^2+p, \quad 1 \le p \le 2n \\
    \max \{p, 2n+2-p\}, & \text{wo } R=\Z[\beta_m], \\
             &  m=(2n+1)^2+4p, \quad 1 \le p \le 2n+1
              \end{cases} \]
   Hierbei ist $\beta_m = (1+\sqrt{m})/2$. Weiter wird bei Gleichheit das
   Minimum $\bmod R$ höchstens in den Punkten $z_1 = \sqrt{m}/2$ oder
   $z_2 = (1+\sqrt{m})/2$ (falls $R = \Z[\sqrt{m}]$), bzw. in
   $z_1 = \beta_m/2$ oder $z_2 = (1+\beta_m)/2$ (falls
   $R = \Z[\beta_m]$) angenommen.
\end{quote}

Das Resultat (2.7) stammt von Heinhold \cite{Hei39} (sh. dazu auch
Lekkerkerker \cite[p. 422]{Lek69}); problemlos erhält man daraus die
Minkowskische Ungleichung für quadratische Zahlkörper:

\begin{quote}
  {\bf (2.8)} {\em Sei $K$ quadratischer Zahlkörper mit Diskriminante $d$;
    dann ist $M(K) \le \sqrt{d}/4$. }
\end{quote}

Wählt man die Basis $\{1, \beta\}$ mit $\beta = n+\sqrt{m}$, wo
$m=n^2+p$, $1 \le p \le n$, ist, so findet man schnell, dass
$\max\{N(z):z \in F'\}$ gleich den in (2.7) angegebenen Schranken von
Heinhold ist; hierbei ist $F' = (0, 0.5) \times (0, 0.5)$. Leider
genügt das nicht, um (2.7) zu beweisen, weil
$F' \cup -F' \cup {F'}^\sigma \cup -{F'}^{\sigma} $hier kein
Fundamentalbereich ist (wie das noch im Falle der Basis
$\{1, \sqrt{m}\}$ der Fall gewesen ist).

Mit (2.7) können wir einige $M(K)$ aus (2.6) nach oben abschätzen: ist
z.B. $m = n^2-1$, so schreiben wir $m = (n-1)^2+2n-2$ und erhalten aus
(2.7) die Schranke $M(K) \le \frac{2n-2}{4} = \frac{n-1}2$.

Zur Abschätzung von $M(K)$ nach unten verwenden wir (1.8). Es sei
jedoch darauf hingewiesen, dass die Schranken $\mu_2$ manchmal von $n$
abhängen (z.B. im Falle $m = n^2+2$). Der Beweis, dass $M(K)$ in (2.6)
dann wirklich exakt ist, ist in diesen Fällen recht mühsam (sh. dazu
Barnes und Swinnerton-Dyer).

Barnes und Swinnerton-Dyer haben in ihrer Arbeit einige interessante
Sätze über die euklidischen Minima in reellquadratischen Zahlkörpern
bewiesen; am Ende des zweiten Teils ihrer Arbeit haben sie darauf
hingewiesen, dass sich diese Sätze auf Körper mit beliebigem
Einheitenrang $>1$ verallgemeinern lassen. Insbesondere für ihr
Theorem M (unser 2.12) ist diese Verallgemeinerung bis heute nur für
Körper mit einer Grundeinheit gelungen (sh. van der Linden 1985,
Theorem 8.7). Wir wollen daher die Schwierigkeiten, die sich hierbei
ergeben, genauer beschreiben. Wir zeigen dazu

\begin{quote}
  {\bf (2.9)} {\em Sei $K$ ein Zahlkörper und $\{u_1, \dots, u_t\}$ ein
    System unabhängiger Einheiten. Dann gibt es eine endliche Menge $Z
    \subset R$ mit der Eigenschaft: für alle $x \in K$ gilt
    \[ M(K, x) = \inf_{u} \min_{z \in Z} N(xu + z), \]
    wobei das Infimum über alle Einheiten $u$ der von den $u_i$
    erzeugten Untergruppe von $R^\times$ gebildet wird, und wo $x_u \in F$
    durch die Kongruenz $x_u \equiv xu \bmod {R}$ bestimmt ist. }
\end{quote}

\begin{proof}
  Sei $\{\alpha_1, \dots, \alpha_t\}$ eine Ganzheitsbasis und
  $Z = \{z = \sum r_i\alpha_i : |r_i| < \mu_i + 0.5\}$ für die
  $\mu_i$ aus (1.12). Wir behaupten, dass dann
  $$ M(K, x) = \inf_{u} \min_{z \in Z} N(xu + z) =: k $$
  gilt. Ansonsten wäre nämlich $M(K, x) < k$, und es gibt ein $y \in
  R$ mit $N(x-y) < k$. Wie in (1.12) folgt nun die Existenz einer
  Einheit $u$ mit $(x-y)u = \sum r_i\alpha_i$ und $|r_i| <
  \mu_i$. Jetzt gilt $(x-y)u = xu \equiv x_u \bmod {R}$ für ein $x_u
  \in F$; also ist $z = (x-y)u - xu \in R$. Mit $z = \sum s_i\alpha_i$
  und $xu = \sum t_i\alpha_i$ ist dann $s_i = r_i - t_i$, also $|s_i|
  \le |r_i| + |t_i| \le \mu_i + 0.5$ und damit $z \in Z$ im
  Widerspruch zur Annahme $M(K, x) < k$.
\end{proof}

Ist $x \in K$, so gibt es offenbar nur endlich viele $x_u \in F$ (denn
mit $x = a/b$ ist $x e = x \bmod {R}$ für jede Einheit $e \in R^\times$).
Ist dagegen $x \notin K$, so kann niemals $x_u = x \bmod {R}$
sein (denn aus $x_u = -x = x(e-1) \in R$ folgt ja $x \in K$). In
diesem Fall sind dann die $x_u \in F$ paarweise verschieden, und wegen
der Kompaktheit von $F$ gibt es einen Häufungspunkt $w$ der $x_u$. Zur
Abschätzung von $M(K, w)$ verwenden wir

\begin{quote}
  {\bf (2.10)} {\em Mit den Bezeichnungen von (2.9) sei $w$ ein
    Häufungspunkt der $x_u$. Dann gilt $M(K, x) \le M(K, w)$.}
\end{quote}

\begin{proof}
  Ist $w$ ein Häufungspunkt der $x_u$, so gibt es eine gegen $w$
  konvergierende Teilfolge $(x_{u_e})$ der $(x_u)$, und da die Norm
  stetig ist, erhalten wir
  $$ \inf_{u} N(x_u + z) \le \lim N(x_e + z)
     = N(\lim x_e + z) = N(w + z) \quad
     \text{für alle } z \in Z. $$
  Nun ist aber mit $w$ auch jedes $w_u$ ein Häufungspunkt der $x_u$,
  weil mit $x_e \to w$ auch $x_{eu} \to w_u$ gilt (hier ist
  natürlich $x_{eu} \in F$ durch $xeu \equiv x_{eu} \bmod {R}$ für
  Einheiten $e, u \in R^\times$ definiert). Also können wir in
  obiger Überlegung $w$ durch $w_u$ ersetzen und haben für alle $u$
  und für alle $z \in Z$ die Ungleichung
  $$ \inf_{u} N(x_u + z) \le N(w_u + z), $$
  und die Bildung des Minimums über alle $z \in Z$ liefert
  $$ M(\uK, x) = \inf_{u} \min_{z \in Z} N(x_u + z) \le \min_{z \in Z}
  N(w_u + z). $$
  Dabei durften wir inf und min vertauschen, weil $Z$ endlich
  ist. Weil schließlich $M(\uK, w) = \inf_{u} \min_{z \in Z} N(w_u + z)$
  ist, folgt $M(\uK, x) \le M(\uK, w)$ wie behauptet.
\end{proof}

Jetzt zeigen wir die Verallgemeinerung von Theorem L (Barnes und
Swinnerton-Dyer):

\begin{quote}
  {\bf (2.11)} {\em Sei $K$ ein Zahlkörper
  \begin{enumerate}
  \item[(i)] Zu jedem $x \in \uK$ gibt es ein $w \in \uK$ mit
    $M(\uK,x) = M(\uK,w)$, und $M(\uK,w)$ wird erreicht.
    \item[(ii)] die Menge $\{M(\uK,x) : x \in \uK\} \subset \R$
      ist abgeschlossen.
  \end{enumerate} }
  \end{quote}

Zum besseren Verständnis von (2.11) halte man sich das Beispiel D(13)
vor Augen: hier hatte man eine Folge $x_1, x_2, \dots$ von Punkten mit
$M(K,x_i) = \frac{1}{3}$, wobei dieses Minimum bei keinem $x_i$
erreicht wird. Dagegen ist $M(K,w) = \frac{1}{3}$ im Häufungspunkt $w
= (\frac{1}{6}, \frac{1}{6})$ der $x_i$, und in $w$ wird das Minimum
erreicht.

\begin{proof}[Bew. von (i)]
  Sei $(x_e)$ eine Teilfolge der $x_u$ mit
  $$ M(\uK,x) = \lim_{e \to \infty} \min_{z \in Z} N(x_e + z). $$
  OBdA dürfen wir annehmen, dass die $x_e$ paarweise verschieden sind:
  Wird $M(\uK,x)$ nicht erreicht, ist dies nämlich sicher möglich, im
  andern Fall ist nichts zu zeigen. Sei nun $w \in F$ ein
  Häufungspunkt der $x_e$; dann gibt es eine gegen $w$ konvergente
  Teilfolge der $x_e$, die wir ebenfalls wieder mit $(x_e)$ bezeichnen
  wollen. Es gilt dann
  $$ M(\uK,x) = \lim_{e \to \infty} \min_{z \in Z} N(x_e + z) =
     \min_{z \in Z} N(w + z) \ge M(\uK,w), $$
  während aus (2.10) $M(\uK,x) \le M(\uK,w)$ folgt. Also ist
  $M(\uK,x) =  M(\uK,w)$, und wegen
  $M(\uK,w) = \min_{z \in Z} N(w+z)$ wird $M(\uK,w)$
  tatsächlich erreicht.
  (ii) Sei $k \in \R$ ein Häufungspunkt der Menge
  $\{M(K,x) : x \in K\}$ und $(x_i)_{i \in \N}$ eine Folge von
  Punkten in $F$ mit $\lim_{i \to \infty} M(\uK,x_i) = k$. Wir müssen dann
  ein $w \in F$ finden mit $k = M(\uK,w)$.
  Nach i) dürfen wir annehmen, dass $M(\uK,x_i)$ für jedes $i \in
  \N$ erreicht wird. Also existiert für jedes $i \in
  \N$ eine Einheit $u = u(i) \in R^\times$ mit der Eigenschaft
  $$ M(\uK,x_i) = \min_{u'} \min_{z \in Z} N(y_{u'} + z) $$
  (wobei $u'$ nur die Einheiten $1, u, u^2, \dots, u^{m-1}$ zu
  durchlaufen braucht). Sei nun $w$ ein
  Häufungspunkt der $x_{iu}$ und $(x_e)$ eine gegen $w$ konvergierende
  Teilfolge der $x_{iu}$. Dann ist
  \begin{align*}
  k & = \lim_{i \to \infty} M(\uK,x_i)
      = \lim_{i \to \infty} \min_{z \in Z} N(x_{iu} + z) \\
    & = \lim_{e \to \infty} \min_{z \in Z} N(x_e + z)
      = \min_{z \in Z} N(w + z) \ge M(K,w)
  \end{align*}
  Nach (2.10) gilt aber auch $M(\uK,x) \ge \inf_i M(\uK,x_{iu}) = k$,
  folglich ist $M(\uK,x) = k$.
\end{proof}

Ein wichtiger Spezialfall von (2.11.ii) ist: es gibt ein $x \in \uK$ mit
$M(\uK,x) = M(\uK)$; diese Beobachtung stammt von Heinhold \cite{Hei39}. Barnes
und Swinnerton-Dyer haben vermutet, dass es sogar ein $x \in K$ gibt
mit $M(\uK,x) = M(\uK)$, und dass insbesondere $M(\uK) = M(K)$ gilt. Wenn
man sich an das Beispiel D(13) erinnert, so scheint diese Vermutung
auf den ersten Blick fast trivial: man suche ein $x \in K$ mit
$M(\uK,x) = M(\uK)$ und wähle dann $w$ als den Häufungspunkt der $x_u$. Der
Haken hierbei ist natürlich, dass $w$ nicht in $K$ zu liegen braucht
(häufen die $x_u$ nur endlich viele Häufungspunkte wie in D(13), so
gäbe es eine Einheit $u$ mit $w \equiv w_u \bmod {R}$, und dies
implizierte in der Tat $w \in K$; man scheint aber die Möglichkeit,
dass die $x_u$ unendlich viele Häufungspunkte haben, nicht so einfach
ausschließen zu können).

Selbst die weniger weitgehende Vermutung $M(K) = M(\uK)$ ist
bisher nur für Körper mit Einheitenrang 1 bewiesen worden
(Theorem M von Barnes und Swinnerton-Dyer (für $n = 2$),
Prop. 5.2. von van der Linden (für Einheitenrang 1)):

\begin{quote}
  {\bf (2.12)} {\em Für alle $\eps > 0$ existiert ein $y \in K$
    mit $M(\uK, y) > M(\uK) - \eps$; Insbesondere ist $M(K) = M(\uK)$.}
\end{quote}

\begin{proof}
  Sei $u$ eine Einheit von $K$ mit $|u|_1 > 1$; aus dem Beweis des
  Dirichletschen Einheitensatzes folgt, dass die Potenzen von $u$
  \glqq{}von der $1$ wegbleiben\grqq{} in dem Sinne, dass es ein $c > 0$
  gibt mit $|u^m - 1|_1 > c$ und $|u^m - 1|_2 > c$ für alle $m \in \Z$. Sei
  weiter $x \in K$ ein Punkt mit $M(\uK, x) = M(\uK)$ (ein solches $x$
  gibt es nach 2.11.ii). ObdA dürfen wir annehmen, dass $x$
  Häufungspunkt der $x_u$ in $F$ ist (sonst ersetzen wir $x$ durch
  einen solchen Häufungspunkt $x'$, denn nach (2.9) ist
  $M(\uK, x') \ge M(\uK, x)$, und wegen $M(\uK, x) = M(\uK)$ und
  $M(\uK, x') \le M(\uK)$ muss auch $M(\uK, x') = M(\uK)$ sein).

  Also gibt es eine Potenz $u^m$ von $u$, sodass für ein beliebig
  vorgegebenes $\delta > 0$ $u^m x \equiv x + z \bmod {R}$ und
  $|z|_j < \delta c$ für $j = 1,2$ ist. Daher existiert ein $a \in R$ mit
  $u^m x = a + x + z$, d.h. mit $(u^m - 1)x = a + z$. Jetzt setzen wir
  $y = a/(u^m - 1)$ und haben $y \in K$, sowie $x - y = z/(u^m -
  1)$. Damit haben wir ein $y \in K$ gefunden, das nicht nur nahe bei
  $x$ liegt, sondern für das auch die Potenzen $uy, u^2 y, \dots, u^m
  y$ nahe bei $ux, u^2 x, \dots, u^m x$ liegen; genauer: für $0 \le k
  \le m$ gilt
  $$ |u^k x - u^k y|_1 = |u^k (x - y)|_1 \le |u^m (x - y)|_1
   = \Big|\frac{z}{1 - u^{-m}}\Big|_1 < \frac{\delta c}{c} = \delta $$
   wegen $|u|_1^k \le |u|_1^m$ und $|1 - u^{-m}|_1 < \frac1c$. Analog
   erhält man
   $$ |u^k x - u^k y|_2 = |u^k (x - y)|_2 \le |(x - y)|_2
   = \Big|\frac{z}{u^m - 1}\Big|_2 < \frac{\delta c}{c} = \delta $$
  wegen $|u|_2 < 1$ und $1/|u^m - 1|_2 < \frac1c$. Diese beiden
  Ungleichungen lassen sich auch in der Form $|x_e - y_e|_1 < \delta$,
  $|x_e - y_e|_2 < \delta$ für $e = u^k$ schreiben. Unsere Aufgabe ist
  es nun, $M(\uK, x)$ und $M(\uK, y)$ zu vergleichen. Nach (2.9) ist
  $$ M(\uK, x) = \inf_{u} \min_{z \in Z} N(x_u + z), $$
  während wegen $u^m y \equiv y \bmod {R} $
  $$ M(\uK, y) = \min_{u'} \min_{z \in Z} N(y_{u'} + z) $$
  gilt, wobei $u'$ nur die Einheiten $1, u, u^2, \dots, u^{m-1}$ zu
  durchlaufen braucht. Wegen $x_u, y_u \in F$ und $z \in Z$ liegen die
  Punkte $x_u + z$ und $y_u + z$ in der beschränkten und abgeschlossenen
  Menge $Z' = \{ \sum r_i \alpha_i : |r_i| \le \mu_i + 1\}$.
  Folglich gibt es ein $\delta>0$ mit der Eigenschaft, dass für alle
  $x,y,z'$ gilt: aus $|x-y|_1 < \delta$ und $|x-y|_2 < \delta$ folgt
  $|N(x)-N(y)| < \eps$ für das vorgegebene $\eps > 0$.

  Durchläuft $u'$ die Einheiten $1, u, \dots, u^{m-1}$, so gilt für
  alle $z \in Z$ und $j=1,2$ $|(x_u + z)-(y_u + z)|_j < \delta$; daher
  ist auch $|N(x_u + z)-N(y_u + z)| < \eps$, insbesondere also
  $N(x_u + z) < N(y_u + z) + \eps$. Jetzt folgt sofort
  $$ \min_{u'} \min_{z \in Z} N(x_{u'} + z) <
     \min_{u'} \min_{z \in Z} N(y_{u'} + z) + \eps = M(K, y) + \eps,
     \text{ aber wegen} $$
   $$ M(\uK, x) = \inf_{u} \min_{z \in Z} N(x_u + z) \le
      \min_{u'} \min_{z \in Z} N(x_{u'} + z) $$
   ist das bereits die Behauptung.
\end{proof}

Die Schwierigkeit bei einer Verallgemeinerung dieses Satzes auf Körper
mit beliebigem Einheitenrang $\ge 1$ liegt einzig und allein in der
Konstruktion eines $y \in K$ mit der Eigenschaft, dass die endlich
vielen $y_u$ nahe bei den entsprechenden $x_u$ liegen. Eine solche
Konstruktion ist mir bisher allerdings nicht gelungen.

Wir wollen der Vollständigkeit halber noch die Verallgemeinerung des
\glqq{}Theorem G\grqq{} von Barnes und Swinnerton-Dyer auf Körper mit
Einheitenrang $\ge 1$ angeben; sie lautet

\begin{quote}
  {\bf (2.13)} {\em Gilt $M(\uK, x) < k$ für alle $x \in \uF$ bis auf eine
    endliche Ausnahmemenge $G$, dann gibt es ein $k' < k$, sodass
    $M(\uK, x) \le k'$ für alle $x \in \uF \setminus G$ gilt. }
\end{quote}

Hieraus folgt insbesondere, dass $M_1(K)$ isoliert ist, wenn es nur an
endlich vielen Stellen in $F$ angenommen wird. Der Beweis von Barnes
und Swinnerton-Dyer funktioniert auch im allgemeinen Fall, wenn man
die Einheit $u$ so wählt, dass $|u|_j \ne 1$ für alle $j$ mit $1 \le j
\le r+s$ wird und (2.3) verwendet.

Wir wollen noch ein paar Worte zur Bestimmung von $M^2(K)$ sagen: die
Begriffe und auch die Sätze aus \S\  2 haben direkte Analoga, wenn man
$R$ durch die Menge der Kettenbrüche der Länge $\le k$
ersetzt. Dasselbe gilt auch für (1.8) und (1.12); allerdings sind die
Analoga von (1.8) und (1.12) wahrscheinlich wertlos, da es
z.B. unendlich viele Kettenbrüche $r+s\sqrt{m}$ der Länge 2 mit
beschränktem $r,s \in \Q$ gibt und daher die Bedingungen
(1.8.a,b,c) unendlich viele mögliche $z \in K$ zulassen.

Um also $M^2(K)$ nach unten abzuschätzen, können wir nicht auf (1.8)
oder (1.12) zurückgreifen; außer der Abschätzung (0.13) hat sich in
Ringen mit Klassenzahl $>1$ folgende Methode bewährt, die wir am
Beispiel $K = \Q(\sqrt{10})$ vorführen wollen (damit wird auch
die Frage beantwortet, die Cooke \cite[S. 75 unten]{Coo77}) gestellt
hat). Wir wählen eine Menge von Kettenbrüchen der Länge 2, z.B.
\[
R_2 = \left\{0, \frac{1}{2}, \frac{1+\sqrt{10}}{2}, \frac{4+\sqrt{10}}{6}, \frac{1}{3}, \frac{\sqrt{10}}{3}\right\}
\]
ihre Nenner haben die Normen 1, 4, 4, 6, 9, 9. Dann versuchen wir, zu
jedem Punkt $x$ aus $F = (0, 0.5) \times (0, 0.5)$ ein $y$ zu finden,
das $\bmod R$ einem Kettenbruch $a_2/b_2$ aus $R_2$ kongruent ist und
$|N_{K/\Q}(x-y)| < k \cdot |N_{K/\Q}(b_2)|$ erfüllt.

Wenn wir $k = 0.99$ wählen, können wir so ganz $F$ bis auf die
Ausnahmemenge $S = (0, 0.001) \times (0.499, 0.5)$ bedecken. Wie im
normeuklidischen Fall folgt nun leicht, dass $z = \sqrt{10}/2$ der
einzig mögliche Ausnahmepunkt ist. Damit bleibt nur noch $M(K,z)$ zu
bestimmen. Da $y = a_2/b_2 = \frac{3}{2} = 1 + \frac{1}{2}$ ein
Kettenbruch der Länge 2 mit $N_{K/\Q}(b_2) = 4$ ist und
$|N_{K/\Q}(z-y)| = \frac{1}{4}$ ist, gilt sicher $M^2(K,z) \le
1$; die dazugehörige Teilerkette ist $\sqrt{10} = 2q_1 + r_1$, $2 =
r_1 q_2 + r_2$ mit $q_1 = 1$, $q_2 = 2$, $r_1 = -2 + \sqrt{10}$,
$r_2 = 6 - 2\sqrt{10}$ und $|N_{K/\Q}(r_2)| = 4 = N_{K/\Q}(2)$.

Nun gibt es aber keine Teilerkette der Länge 2 mit
$|N_{K/\Q}(r_2)| \le 4$, denn ist $P = (\sqrt{10}, 2)$ das
Primideal der Norm 2 über $(2)$, so gilt $r_1, r_2 \in P$, und
insbesondere ist $N_{K/\Q}(r_2) = 0 \bmod 2$. Es müsste daher
$|N_{K/\Q}(r_2)| = 0$ oder $|N_{K/\Q}(r_2)| = 2$ sein;
letzteres ist unmöglich, da $P$ kein Hauptideal ist, und aus demselben
Grund kann auch nicht $r_2=0$ sein, da sonst $P = (r_1)$ folgte.

Also ist $M^2(K) = M^2(\uK) = M^2(K,z) = 1$, und $K$ hat
euklidische Tiefe 1 (denn genau wie oben zeigt man, dass es keine von
$(\sqrt{10},2)$ ausgehende Teilerkette mit $|N_{K/\Q}(r_k)|
\le 4$ gibt). Definiert man die Mengen\label{dDefB}
$B_j$ ($j = 1, 2, \dots, \infty$) durch $B_j = (F \cap K) \setminus F_j$
($F$ ist ein Fundamentalbereich; die $F_j \subset K$ sind auf
S. \pageref{pEi} definiert), so
gilt hier $B_1 = B_2 = \dots = B_\infty = \{\pm \sqrt{10}/2\}$.

Auf dieselbe Art und Weise ist Tabelle~\ref{dTab6} zustande gekommen:

\begin{table}[h]
\centering
\begin{tabular}{llll}
\toprule
$m$ & $M^2(K)$ & $B_1$ & $B_2 = B_3 = \dots$ \\
\midrule
6 & $1/4$ & $\{\}$ & $\{\}$ \\
10 & $1$ & $\{(0, 1/2)\}$ & $\{(0, 1/2)\}$ \\
14 & $1/4$ & $\{(1/2, 1/2)\}$ & $\{\}$ \\
15 & $1$ & $?$ & $\{(1/2, 1/2)\}$ \\
26 & $1$ & $?$ & $\{(0, 1/2)\}$ \\
30 & $3/2$ & $?$ & $\{(0, 1/2)\}$ \\
34 & $1$ & $?$ & $\{(1/2, 1/2) + (1/2, 1/2)\}$ \\
35 & $7/5$ & $?$ & $\{(0, 1/2), (1/2, 1/2)\}$ \\
39 & $5/2$ & $?$ & $\{(1/2, 1/2)\}$ \\
65 & $1$ & $\{(1/2, 1/2)\}$ & $\{(1/2, 1/2)\}$ \\
85 & $1$ & $?$ & $\{(1/2, 1/2) + (1/2, 1/2)\}$ \\
\bottomrule
\end{tabular}
\caption{Euklidische Minima}\label{dTab6}
\end{table}

Aus dieser Tabelle kann man entnehmen, dass z.B. die euklidische Tiefe
von $D(26)$ und $D(30)$ gleich 2 ist. Bereits Cooke hat versucht, die
euklidische Tiefe dieser beiden Ringe zu bestimmen. Für $D(30)$ konnte
er zeigen, dass die euklidische Tiefe $\le 3$ ist, und er hat vermutet,
dass sie gleich 2 ist. Im Falle $D(26)$ hat er gezeigt, dass die
euklidische Tiefe endlich ist; dass er diese nicht bestimmen konnte,
lag an einem Rechenfehler. Er schreibt auf \cite[S. 81]{Coo77}, dass er mit
Kettenbrüchen der Länge 2 ganz $F$ bis auf die folgenden Punkte
bedecken kann: $(0, \frac{1}{2})$, sowie
$(\pm \frac{2}{5}, \pm \frac{2}{5})$. Allerdings ist $\sqrt{26}$ ein
Kettenbruch der Länge 2 und z.B.
$$ \Big| N_{K/\Q}\Big(z - \frac{\sqrt{26}}{3}\Big) \Big| = \frac{2}{45} 
     = \frac{2}{5} \cdot N_{K/\Q}(3) < N_{K/\Q}(3) $$
für $z = \frac{2 + 2 \sqrt{26}}{5}$.

Damit folgt auch aus seinen Rechnungen, dass nur der Punkt
$\sqrt{26}/2$ unbedeckt bleibt und die euklidische Tiefe gleich 2 ist.

Weiter bleibt noch anzumerken, dass $D(10)$ euklidische Tiefe 1 hat
(und nicht 2, wie Cooke schreibt), da hier $B_1 = B_\infty =
\{\sqrt{10}/2\}$ ist; Entsprechendes gilt für $D(65)$.

Mit dieser Methode lassen sich auch die Ergebnisse von Cooke \cite{Coo77}
bestätigen, wonach die Ringe $D(m)$ für
$$ m = 14, 22, 23, 31, 38, 43, 46, 53, 61, 69, 77, 89, 93, 97, 113, 129,
       133, 137, 181, 253 $$
2-stufig normeuklidisch sind; darüberhinaus erhält man für
$$ m = 47, 59, 62, 67, 71, 101, 109, 149, 157, 161, 173, 177, 193,
   197, 201, 213 $$
neue 2-euklidische Ringe.

Man muss hierbei allerdings beachten, dass die für $M^2(K)$ erzielten
Schranken sowie die möglichen Ausnahmepunkte von der Anzahl der
verwendeten Kettenbrüche der Länge 2 abhängt; so kann man in $D(14)$
mit Kettenbrüchen der Länge 2, deren Nenner eine Norm $\le 11$ haben,
ganz $F$ bedecken bis auf die beiden Punkte $x_1 = (\frac{1}{2},
\frac{9}{14})$ und $x_2 = (\frac{1}{2}, \frac{1}{2})$. Wegen (0.13)
ist $M^2(x_2) = \frac{1}{4}$; wegen $|N_{K/\Q}(x_1 -
\frac{2+\sqrt{14}}{7 \cdot 2})| = \frac{1}{28}$ und weil
$\frac{2+\sqrt{14}}{14}$ ein Kettenbruch der Länge 2 mit Nenner
$7\cdot 2/14$ ist, gilt sicher auch $M^2(x_1) \le \frac{1}{4}$ (dies
genügt bereits, um $M^2(K) = \frac{1}{4}$ zu beweisen).

Verwendet man aber auch Kettenbrüche, deren Nenner Norm 13 haben, und
beachtet, dass $z = (6+6\sqrt{14})/13$ ein solcher ist, so findet man,
dass $|N_{K/\Q}(x_1 - z)| < \frac{1}{32}$ ist; dies zeigt, dass
$x_2$ der einzige Ausnahmepunkt für $k = 1/4$ ist.

Wie Cooke bemerkt hat, stehen die Punkte aus $B_\infty$ in einer
Beziehung zur Klassenzahl von $K$; um dies einzusehen, definieren wir
eine Abbildung $\varphi : K \to \mathrm{Cl}(K)$ durch
$\varphi(\frac{\alpha}{\beta}) = [(\alpha, \beta)]$ (d.h. einem $x =
\frac{\alpha}{\beta}$ wird die von dem Ideal $(\alpha, \beta)$
erzeugte Idealklasse zugeordnet). $\varphi$ ist wohldefiniert: ist
nämlich $\alpha/\beta = \alpha'/\beta'$, so folgt $\alpha\beta' =
\alpha'\beta$ und damit $(\alpha', \beta) \sim (\alpha')(\alpha,
\beta) = (\alpha\alpha', \alpha'\beta) = (\alpha\alpha', \alpha\beta')
\sim (\alpha', \beta')$, und wir haben $[(\alpha, \beta)] = [(\alpha',
  \beta')]$.

Ganz einfach folgt nun

\begin{quote}
  {\bf (2.14)} {\em $\varphi(B_\infty)$ enthält alle nichttrivialen
    Idealklassen.}
\end{quote}

\begin{proof}
  Sei $I$ ein Ideal in $R$, aber kein Hauptideal. Wir wählen
  $\beta \in I \setminus \{0\}$ so, dass $|N_{K/\Q}(\beta)|$ minimal
  wird, und dann ein $\alpha \in I$ mit $I = (\alpha, \beta)$. Indem
  wir von $\alpha$ ein geeignetes Vielfaches von $\beta$ subtrahieren,
  können wir $\frac{\alpha}{\beta} \in F$ erreichen, ohne die
  Eigenschaft $I = (\alpha, \beta)$ zu zerstören.
  Wir behaupten nun, dass $\frac{\alpha}{\beta} \in B_\infty$ ist. Für
  jede von $(\alpha, \beta)$ ausgehende Teilerkette ist nämlich $r_k
  \in (\alpha, \beta) \setminus \{0\}$ (wäre $r_k = 0$, so wäre $I =
  (r_{k-1})$ Hauptideal), und nach der Wahl von $\beta$ muss
  $|N_{K/\Q}(r_k)| \ge |N_{K/\Q}(\beta)|$ sein. Also ist $K$ in $x =
  \frac{\alpha}{\beta}$ nicht quasi-euklidisch und folglich
  $\frac{\alpha}{\beta} \in B_\infty$.
\end{proof}

Damit kann man aus der Kenntnis von $B_\infty$ bereits die Klassenzahl
von $K$ bestimmen: für $D(10)$, $D(15)$, $D(26)$ usw. enthält $B_2 =
B_\infty$ nur einen einzigen Punkt, also ist $h(K) \mid 2$ in diesen
Fällen. Andererseits ist hier $\varphi(B_\infty)$ nie die
Hauptidealklasse, und wir können $h=2$ in diesen Fällen schließen.

Etwas interessanter ist $D(35)$: hier enthält $B_\infty$ drei
verschiedene Punkte, und wir haben zuerst nur $h(K) \mid
4$. Allerdings liegen die Ideale $5_1 = (2+\sqrt{35}, 5)$ und $2_1 =
(1+\sqrt{35}, 2)$ in derselben Idealklasse, sodass sich auch hier $h=2$
ergibt.

Dies wirft die Frage auf, wann $\varphi(x) = \varphi(x')$ für zwei
Punkte $x,x'\in K$ gilt. Dazu sei $SL_2(R)$ die Gruppe der
$2\times2$-Matrizen mit Einträgen aus $R$ und Determinante $+1$;
$SL_2(R)$ operiert auf $K$ mittels
\[
\begin{pmatrix} a & b \\ c & d \end{pmatrix}
\begin{pmatrix} \alpha \\ \beta \end{pmatrix} =
\frac{a\alpha + b\beta}{c\alpha + d\beta}
\]
(genaugenommen muss man $K$ durch $K \cup \{\infty\}$ ersetzen, will
man nicht von vornherein ausschließen, dass $c\alpha + d\beta = 0$
wird; schränkt man den Definitionsbereich von $SL_2(R)$ jedoch auf
$B_\infty$ ein, kann dies nie eintreten). Zwei Punkte $x =
\alpha/\beta$ und $x' = \alpha'/\beta'$ heißen nun äquivalent ($x
\equiv x'$), wenn es ein $A\in SL_2(R)$ gibt mit $A(x) = x'$. Damit
gilt

\begin{quote}
  {\bf (2.15)} {\em Sei $K$ ein algebraischer Zahlkörper mit
    unendlicher Einheitengruppe. Dann enthält $B_\infty$ genau $h-1$
    Äquivalenzklassen von Punkten.}
\end{quote}

Dazu muss man erstens zeigen, dass $\varphi(x) = \varphi(x')$ genau dann
gilt, wenn $x \equiv x'$ gilt, und dass zweitens $B_\infty$ kein $x\in
K$ enthält, sodass $\varphi(x)$ die Hauptidealklasse ist. Das letztere
Ergebnis scheint alles andere als elementar zu sein und wurde von
Vaserstein \cite{Vas72} bewiesen; es erscheint mir wünschenswert, dessen
Beweis in die Sprache der algebraischen Zahlentheorie zu
übersetzen. Beim Beweis der ersten Tatsache benutzt Cooke die folgende
Beobachtung von Hurwitz \cite{Hur95a}:

\begin{quote}
  {\bf (2.16)} {\em Die Ideale $(\alpha, \beta)$ und $(\alpha', \beta')$
    sind genau dann gleich, wenn es ein $A\in SL_2(R)$ gibt
    mit $A(\frac{\alpha}{\beta}) = \frac{\alpha'}{\beta'}$. }
\end{quote}

Der Beweis, den Cooke für (2.15) gibt, ist ein reiner Existenzbeweis
und kann daher nicht zufriedenstellen. Der Originalbeweis von Hurwitz
dagegen gibt einem zumindest einen Hinweis darauf, wie man $A$
konstruieren kann, außerdem ist er einfacher als Cookes Beweis:

\begin{proof}[Bew. von (2.15)]
  Sei $A \in \text{SL}_2(\R)$ und $A(\alpha/\beta) =
  \alpha'/\beta'$. Dann gibt es ein $\eta \in \uK$ mit $\alpha'\eta =
  a\alpha + b\beta$, $\beta'\eta = c\alpha + d\beta$. Dies impliziert
  $(\alpha, \beta) = (\eta)(\alpha', \beta')$. $(\alpha', \beta')$
  also $\varphi(\frac{\alpha}{\beta}) = \frac{\alpha'}{\beta'}$.

  Sei nun umgekehrt
  $\varphi(\frac{\alpha}{\beta}) = \frac{\alpha'}{\beta'}$; dann ist
  $(\alpha, \beta) = (\alpha', \beta')$, und nach (2.17) gibt es
  $\gamma, \delta \in \R$ mit $(\gamma, \delta) = (\alpha, \beta)$
  und $\frac{\gamma}{\delta} = \frac{\alpha}{\beta}$. Mit (2.16)
  erhalten wir ein $A \in \text{SL}_2(\R)$ mit
  $A(\frac{\alpha}{\beta}) = (\frac{\alpha'}{\beta'})$, und
  dasselbe $A$ erfüllt dann auch
  $A(\frac{\alpha}{\beta}) = (\frac{\alpha'}{\beta'})$.
\end{proof}

Als Beispiel wollen wir zeigen, dass die Punkte $2\sqrt{35}/5$ und
$(1+\sqrt{35})/2$ äquivalent sind. Dazu schreiben wir
$(1+\sqrt{35})/2$ mit (2.16) in der Form $(1+\sqrt{35})/2 =
\alpha/\beta$ mit $\alpha' = 20+3\sqrt{35}$ und $\beta' =
5-\sqrt{35}$; damit ist $(2\sqrt{35}, 5) = (\alpha', \beta')$. Da
$D(35)$ Klassenzahl 2 hat und $(2\sqrt{35}, 5)^2 = (5)$ ist, können
wir wie in (2.16) mit $\alpha = 2\sqrt{35}$ und $\beta = 5$ schreiben:
$5 = \alpha\gamma + \beta\delta = \alpha'\gamma' + \beta'\delta'$ mit
$\gamma, \delta, \gamma', \delta' \in (2\sqrt{35}, 5)$ für $\gamma =
-\sqrt{35}$, $\delta = 15$, $\gamma' = \sqrt{35}$, $\delta' = -20$ und
erhalten die Matrix
\[
A = \begin{bmatrix} -41-4\sqrt{35} & 60+17\sqrt{35} \\
    -7-2\sqrt{35} & 29+3\sqrt{35} \end{bmatrix},
\]
die in der Tat unimodular ist und $A(\alpha/\beta) = \alpha'/\beta'$
erfüllt.

Mit Hilfe der oben definierten Abbildung $\varphi: \uK \to
\text{Cl}(\uK)$ lässt sich auch ein einfaches Kriterium für
semi-euklidische Ringe angeben:

\begin{quote}
  \textbf{(2.18)} {\em Ein Zahlring $R$ ist genau dann semi-euklidisch,
    wenn für jedes $x \in \mathbb{B}_1$ gilt: $\varphi(x) \neq [(1)]$, d.h.
    mit $x = \frac{\alpha}{\beta}$ ist $(\alpha, \beta)$ kein Hauptideal.}    
\end{quote}

\begin{proof}
  Sei $\varphi(x) \neq [(1)]$ für alle $x \in B_1$, und seien
  $\alpha, \beta \in R$ mit $(\alpha, \beta) = R$. OBdA dürfen wir
  $\alpha \bmod \beta$ so wählen, dass $\frac{\alpha}{\beta} \in \uF$
  wird. Wegen $\varphi(x) = [(1)]$ für $x = \alpha/\beta$ liegt $x$
  nicht in $B_1$; dies bedeutet, dass es ein $y \in \R$ mit
  $|N_{\uK/\Q}(x-y)| < 1$ gibt. Also ist $\R$ semi-euklidisch.

  Sei andererseits $R$ semi-euklidisch und
  $x = \frac{\alpha}{\beta} \in \uK$. Ist $(\alpha, \beta)$ ein
  Hauptideal (also $\varphi(x) = [(1)]$), so gibt es, weil $R$
  semi-euklidisch ist, ein $y \in \R$ mit $|N_{\uK/\Q}(x-y)| < 1$.
  Also liegt $x$ nicht in $B_1$.
\end{proof}

Ein analoges Ergebnis erhält man für $k$-stufig semi-euklidische Ringe,
wenn man $B_1$ durch $B_k$ ersetzt. So ist z.B. $D(10)$
semi-euklidisch, weil $B_1 = \{\sqrt{10}/2\}$ und $(\sqrt{10}, 2)$
kein Hauptideal ist.

\section*{{\sc Anmerkungen zu} \S\ 2}
\addcontentsline{toc}{section}{{\sc Anmerkungen zu} \S\ 2}

Der Begriff \glqq{}euklidisches Minimum\grqq{} in Zahlkörpern ist in
der mathematischen Literatur unbekannt; stattdessen wurde bisher der
Ausdruck \glqq{}inhomogenes Minimum der zugehörigen Normform\grqq{}
verwendet. Das erste Resultat über euklidische Minima stammt von
Minkowski, der für reellquadratische Zahlkörper mit Diskriminante $d$
gezeigt hat, dass $M(K) = M_2(K) \le \sqrt{d}/4$ gilt. Minkowski hat
weiter vermutet, dass für totalreelle Zahlkörper $n$-ten Grades die
Ungleichung $M(K) \le 2^{-n}\sqrt{d}$ gilt; dies konnte jedoch bisher
nur in Spezialfällen ($n = 2, 3, 4, 5$; sh. dazu Remak \cite{Rem34},
Dyson \cite{Dys48} und Skubenko \cite{Sku72}) bewiesen werden. Nach
Tschebotarev gilt immerhin $M(K) \le 2^{-n/2}\sqrt{d}$ (sh. dazu auch
Hardy und Wright, \cite[S. 456--458]{HW}.  Für Zahlkörper beliebiger
Signatur ist noch nicht einmal eine analoge Vermutung bekannt.

In der entgegengesetzten Richtung ist das Ergebnis von Davenport
\cite{Dav51} zu erwähnen, wonach für reell-quadratische Körper
$M(K) \ge \frac{\sqrt{d}}{128}$ gilt. Am Ende dieser Arbeit bemerkt
Davenport, dass Prasad die Konstante $\frac{1}{128}$ auf
$\frac{1}{36}$ verbessert hat. Wie Ennola \cite{Enn58b} schreibt, hat
Prasad in seiner Arbeit (vor deren Veröffentlichung) jedoch einen
Fehler entdeckt. Die \glqq{}Davenport-Konstante\grqq{}
$\sup \frac{M(K)}{\sqrt{d}} \ge \frac{1}{128}$ wurde dann von
Cassels \cite{Cas52} auf $\frac{1}{36}$
verbessert, und schließlich gelang es Ennola,
$$M(K) \ge \frac{\sqrt{d}}{16+6\sqrt{6}}
          \approx \frac{\sqrt{d}}{30.69}$$
zu zeigen; dies ist bis heute das beste Ergebnis geblieben.

Die explizite Bestimmung euklidischer Minima in reellquadratischen
Zahlkörpern haben Heinhold \cite{Hei39}, Davenport
\cite{Dav46}, Varnavides \cite{Var48a,Var48b},
Bambah \cite{Bam50,Bam51} und Inkeri \cite{Ink49} vorgenommen. Nach den
grundlegenden Arbeiten von Barnes und Swinnerton-Dyer haben sich nur
noch Godwin \cite{God62,God65b} und Varnavides \cite{Var70} mit der
Bestimmung euklidischer Minima in quadratischen Zahlkörpern
befaßt. Die ersten Minima kubischer Zahlkörper hat Davenport
\cite{Dav47} berechnet; weitere Ergebnisse in
dieser Richtung stammen von Prasad \cite{Pra49},
Clarke \cite{Cla51a}, Samet \cite{Sam54},
Swinnerton-Dyer \cite{Swi54}, Godwin \cite{God55},
Smith \cite{Smi69,Smi71} und Taylor \cite{Tay75,Tay76}
(sh. hierzu \S\ 4).

Für Zahlkörper mit Körpergrad $\ge 4$ sind bisher keine euklidischen
Minima bestimmt worden; lediglich Cohn und Deutsch \cite{CD86} haben
vermutet, dass für $K = \Q(\sqrt{2+\sqrt{2}})$ $M(K) =
\frac{1}{2}$ und $M_2(K) = \frac{1}{4}$ gilt, während sie für $K =
\Q(\sqrt{3+\sqrt{2}})$ vermutet haben, dass $M(K) =
\frac{1}{2}$ und $M_2(K) = \frac{7}{16}$ ist.

\chapter*{\S\ 3 Quadratische Zahlkörper}
\setcounter{chapter}{3}
\addcontentsline{toc}{chapter}{\S\ 3 Quadratische Zahlkörper}
\markboth{Euklidische Ringe}{\S\ 3 Quadratische Zahlkörper}

Daß quadratische Zahlkörper $K = \Q(\sqrt{m})$ genau für die Werte
$$ d = -11, -8, -7, -4, -3, 5, 8, 12, 13, 17, 21, 24, 28, 29, 33,
       37, 41, 44, 57, 73, 76 $$
normeuklidisch sind (wo $d = \disc K$ die Diskriminante von $K$
bezeichnet), ist ein bereits klassisches Ergebnis. Wir werden sehen,
dass sich der Beweis hierfür recht übersichtlich und verhältnismäßig
kurz darstellen lässt, wenn man (1.5) und (1.8) konsequent
anwendet. Außerdem ermöglichen die Ideen, die hier vorgestellt werden,
zumindest teilweise die Klassifikation von normeuklidischen Ringen
höheren Grades (sh. insbesondere \S\ 5).

Wir stellen nun einige bekannte Eigenschaften quadratischer Zahlkörper
zusammen: so ist jeder quadratische Zahlkörper von der Form
$K = \Q(\sqrt{m})$ für ein quadratfreies $m \in \Z$. Weiter
gilt im Falle

\begin{align*}
  m \equiv 1 \bmod{4} : & \left\{1, \frac{1+\sqrt{m}}{2}\right\}
                       \text{ ist GHB, } \disc K = m, \\
  m \equiv 2, 3 \bmod{4} & : \{1, \sqrt{m}\} \text{ ist GHB, }
                       \disc K = 4m.
\end{align*}

Bezeichnet $(\frac{\cdot}{p})$ das Kroneckersymbol, so lässt
sich das Zerlegungsgesetz (darunter verstehen wir eine Beschreibung
des Verhaltens von Primidealen bei Körpererweiterung) wie folgt
beschreiben:

\begin{enumerate}
\item[(I.)] $(\frac{d}{p}) = +1$: $(p) = P_1P_2$, \quad
  $\|P_1\| = \|P_2\| = p$ \quad ($p$ heißt \textbf{zerlegt})
\item[(II.)] $(\frac{d}{p}) = -1$: $(p) = P$, \quad $\|P\|
  = p^2$ \quad ($p$ heißt \textbf{träge})
\item[(III.)] $p \mid d$ : $(p) = P^2$, \quad $\|P\| = p$ \quad ($p$
  heißt \textbf{rein verzweigt}).
\end{enumerate}

Die Einheitengruppe $\cO_K^*$ von quadratischen Zahlringen
hat eine recht einfache Struktur. Ist $K$ imaginärquadratisch, so
ist $\cO_K^* = \{\pm 1\}$ außer für $K =
\Q(\sqrt{-1})$ und $K = \Q(\sqrt{-3})$, wo die 4. bzw.
6. Einheitswurzeln ganz $\cO_K^*$ bilden. In
reell-quadratischen Zahlkörpern lässt sich nach dem Dirichlet'schen
Einheitensatz jede Einheit $e \in \cO_K^*$ in der Form $e =
\pm u^k$ für ein $k \in \Z$ schreiben, wobei $u$ die
Fundamentaleinheit von $K$ bezeichnet (diese wird i.A. durch die
Forderung $u>1$ eindeutig festgelegt).

Auch die Parität der Klassenzahl $h(K)$ kann man leicht angeben. Dazu
sei $m \in \N$; wir unterscheiden dann

$K = \Q(\sqrt{-m})$: $h(K)$ ist genau dann ungerade, wenn $m
\equiv 1, 2$ oder $m \equiv 3 \bmod{4}$ prim ist, d.h. genau dann,
wenn $\disc K$ eine Primzahlpotenz ist.

$K = \Q(\sqrt{m})$: $h(K)$ ist genau dann ungerade, wenn gilt:

\begin{enumerate}
\item[a)] $m = p$ ist prim;
\item[b)] $m = pq$, $p$ und $q$ prim, $p=2$ oder $p\equiv 3 \bmod{4}$,
  $q\equiv 3 \bmod{4}$.
\end{enumerate}

Für imaginärquadratische $K$ steht ein elementarer Beweis bei
Connell\index[dN]{Connell} \cite{Con62}, für reellquadratische bei
Redei\index[dN]{Redei@R\'edei} \cite{Red60}. Einen anderen Zugang zu diesen
Ergebnissen ermöglicht die Geschlechtertheorie von Gauß
(sh. Zagier\index[dN]{Zagier} \cite{Zag81}). Schließlich kann man aus der
analytischen Klassenzahlformel einen weiteren Beweis gewinnen (ein
Teil eines solchen Beweises steht bei Hasse\index[dN]{Hasse}
\cite[Va,Vb]{Has64}).

Da euklidische Ringe Klassenzahl 1 haben, brauchen wir nur die oben
angegebenen Möglichkeiten in Betracht zu ziehen. Als der einfachste
Fall hat sich derjenige herausgestellt, wo $\disc K \equiv 0 \bmod{4}$
ist. Dass $K$ dann höchstens für die Werte $m = 2, 3, 6, 7, 11, 19$
normeuklidisch sein kann, haben unabhängig voneinander
J. Fox\index[dN]{Fox} \cite{Fox35} und E. Berg\index[dN]{Berg} \cite{Ber35}
gezeigt (die Arbeit von J. Fox taucht in den Literaturangaben fast aller
Autoren unter dem Namen Fox Keston auf; wir halten uns hier an
Bull. Am. Math. Soc. 41 (1935), p. 186, wo ihr Name mit Jeanette Fox
angegeben wird). Die einfachsten Beweise im Falle $\disc K \equiv 0
\bmod{4}$ stammen aber wohl von Behrbohm\index[dN]{Behrbohm} und
R\'edei\index[dN]{Redei@R\'edei} \cite{BR36}, und diese wollen wir nun auch
vorstellen.

Zuerst zeigen wir

\begin{quote}
  \textbf{(3.1)} {\em Ist $a \in \N$ und $( \frac{m}{a}) = -1$, so
    ist weder $a$ noch $-a$ Norm aus $\cO_K$.}
\end{quote}

\begin{proof}
  Wegen $( \frac{m}{a} ) = -1$ enthält $a$ einen Primfaktor
  $p$, der in $a$ ungerade oft aufgeht (d.h. es ist $a = p^e b$, $p
  \nmid b$, $e \equiv 1 \bmod{2}$) und für den $( \frac{m}{p}) = -1$
  gilt. Mit $a$ wäre wegen $p \nmid b$ auch $p^e$ Idealnorm aus
  $\cO_K$; weil $e$ ungerade und $\|(p)\| = p^2$ ist, müsste
  es daher ein Ideal der Norm $p$ in $\cO_K$ geben, was aber
  dem Zerlegungsgesetz widerspricht. Also ist $a$ keine Idealnorm,
  somit $a$ keine Norm einer Zahl aus $\cO_K$.
\end{proof}

Unser nächstes Ergebnis wird auch bei der Untersuchung von Zahlkörpern
der Form $\Q(\sqrt[k]{m})$, $k=2^l$, eine große Rolle spielen:

\begin{quote}
  \textbf{(3.2)} {\em Sei $q \equiv 3 \bmod{4}$ prim und $q \geq 11$;
    dann gibt es $a,b \in \N$ mit $2q = a + b$,
    $a \equiv 5 \bmod{8}$ und $(\frac{a}{q}) = \pm 1$.}
\end{quote}

\begin{proof} Wir unterscheiden
  \begin{enumerate}
    \item[1.] $q \equiv 3 \bmod{8}$: für $q=11$ wählen wir $a=5$; also
      dürfen wir $q \geq 19$ annehmen. Es gibt dann ein $c \in
      \N$ mit $2q < 16c < 3q$, und wir haben $q < 16c - q <
      2q$, sowie $0 < 8c - q < q/2$. Also sind $8c-q$ und $16c-q$
      natürliche Zahlen aus dem Intervall $(0, 2q)$, und beide sind
      $\equiv -q \equiv 5 \bmod{8}$. Wegen $(\frac{2}{q}) = -1$ haben
      beide Zahlen aber verschiedenen Restcharakter mod $q$
      (beachte $16c-q \equiv 2(8c-q) \bmod{q}$), d.h. genau eine
      der beiden Zahlen ist quadratischer Rest, und wir dürfen
      $a = 8c-q$ oder $a = 16c-q$ wählen.
    \item[2.] $q \equiv 7 \bmod{8}$: dann ist $q \geq 23$, und wir
      wählen ein $c \in \N$ mit $5q < 16c < 6q$. Wie oben
      liegen dann die beiden Zahlen $16c-5q$ und $16c-3q$ im Intervall
      $(0, 2q)$, sind $\equiv 5 \bmod{8}$ und haben entgegengesetzten
      quadratischen Restcharakter mod $q$. Wie in 1. folgt nun die
      Behauptung.
  \end{enumerate}
\end{proof}

Jetzt folgt sofort

\begin{quote}
  \textbf{(3.3)} {\em Sei $m \equiv 2 \bmod{4}$, $m \in \N$
    quadratfrei. Dann ist $D(m)$ genau für $m=2$ und $m=6$
    normeuklidisch.}
\end{quote}

\begin{proof}
  Dass $D(2)$ und $D(6)$ normeuklidisch sind, wissen wir aus \S\ 2.
  Weiter haben wir schon in \S\ 1 gesehen, dass $D(14)$ nicht
  normeuklidisch ist. Sei also $m \ge 22$ und $D(m)$
  normeuklidisch. Wegen $h(K)=1$ dürfen wir $m=2q$, $q \equiv 3 \bmod
  4$ prim annehmen. Jetzt verwenden wir (1.6) mit $f=2q$ und den
  Werten von $a$, $b$ aus (3.2) und zeigen, dass weder $a$ noch $-b$
  Normen aus $D(m)$ sind. Wegen $(2/a) = -1$ und $(2/b) = 1$ (man
  beachte dazu $a \equiv 5 \bmod{8}$, $b = 2q-a \equiv 1 \bmod{8}$)
  erhalten wir
  $$ (\tfrac{2q}{a}) = -(\tfrac{q}{a}) = -(\tfrac aq) = -1, \text{ sowie }
     (\tfrac{2q}{b}) = (\tfrac qb) = (\tfrac bq) = (\tfrac{2q-a}{q}) =
     (\tfrac{-a}{q}) = -1. $$
  Mit (3.1) folgt nun, dass weder $a$ noch $-b$ Normen aus $D(m)$ sind,
  und (1.6) liefert die Behauptung.
\end{proof}
Genauso elementar lässt sich auch der Fall $m \equiv 3 \bmod{4}$ behandeln:

\begin{quote}
  \textbf{(3.4)} {\em Sei $q \equiv 3 \bmod{4}$ prim und $q \ge 23$;
    dann gibt es $a,b \in \N$ mit $q \equiv a \pm b$,
    $(a/q)=+1$ und $a \equiv 5-q \bmod{8}$.}
\end{quote}

\begin{proof}
  Für die primen $q$ mit $19 < q < 64$ geben wir die Paare $(q,a)$
  direkt an: (23,6), (31,14), (43,10), (47,6), (59, 26). Jetzt dürfen
  wir $q \ge 67$ annehmen; wieder unterscheiden wir:
  1. $q \equiv 3 \bmod{8}$: wir wählen dann ein $c \in \N$ mit
  $2q-16 < 64c < 3q-8$. Damit liegen die beiden Zahlen $8c+2$ und
  $8(8c+2)-2q$ im Intervall $(0,q)$, beide sind $\equiv 2 = 5-q \bmod
  8$, und genau eine von ihnen ist quadratischer Rest $\bmod q$.
  2. $q \equiv 7 \bmod{8}$: sei $c \in \N$ mit $q+16 < c <
  2q+16$; wir sehen, dass sowohl $8c-2$, als auch $2q-8(8c-2)$ zwischen
  0 und $q$ liegen, $\equiv -2 = 5-q \bmod{8}$ sind und
  entgegengesetzten quadratischen Restcharakter $\bmod q$ haben.
\end{proof}

\begin{quote}
  \textbf{(3.5)} {\em Sei $m \equiv 3 \bmod{4}$ quadratfrei; dann ist
    $D(m)$ genau für $m \equiv 3, 7, 11, 19$ normeuklidisch.}
\end{quote}

\begin{proof}
  Für $m \equiv 3, 7, 11, 19$ ist $M(K)<1$, folglich $D(m)$
  normeuklidisch. Sei daher $m=q$, $q \equiv 3 \bmod{4}$ prim und $q \ge
  23$; nach (3.4) gibt es $a,b \in \N$ mit $q \equiv a+b$, $a
  \equiv 5-q \bmod{8}$ und $(a/q)=+1$. Mit diesen Werten für $a$ und
  $-b$ verwenden wir nun (1.6) mit $f=q$ und müssen noch zeigen, dass
  $a$ und $-b$ keine Normen aus $D(m)$ sind. Da mit $a$ auch $c = a/2
  \equiv 1 \bmod{2}$ eine Norm aus $D(m)$ wäre, folgt dies aber sofort
  aus $(q/c) = (q/b) = -1$ und (3.1).
\end{proof}

Zum Teil lässt sich auch der Fall $m \equiv 1 \bmod{4}$ auf diese Art
behandeln; so gilt z.B.

\begin{quote}
  \textbf{(3.6)} (Hofreiter\index[dN]{Hofreiter} \cite{Hof34})
    {\em Sei $m \equiv 3q$, $q \equiv 7 \bmod{8}$; dann ist $D(m)$ genau 
    für $m = 21$ normeuklidisch.}
\end{quote}

\begin{proof}
  Wegen $M(K) = \frac57 < 1$ für $K = \Q(\sqrt{21})$ dürfen wir $q
  \ge 23$ und $q \equiv 7 \bmod{8}$ prim annehmen; dann verwenden wir
  (1.6) mit den folgenden Werten für $a$ und $b$:
  \begin{center}
    \begin{tabular}{|c|c|c|}
      \hline
      $f=q$ & $a$ & $b$ \\
      \hline
      $q \equiv 7 \bmod{24}$ & 18 & $q-18$ \\
      $q \equiv -1 \bmod{24}$, $q \equiv 5 \bmod{9}$ & 2 & $q-2$ \\
      $q \equiv -1 \bmod{24}$, $q \equiv 2 \bmod{9}$ & 8 & $q-8$ \\
      $q \equiv -1 \bmod{24}$, $q \equiv 8 \bmod{9}$ & 32 & $q-32$ \\
      \hline
    \end{tabular}
  \end{center}
  Wegen $q \equiv 7 \bmod{8}$ ist 2 und damit auch $a$ quadratischer
  Rest mod $q$; andererseits wäre mit $a$ auch 2 eine Norm aus $D(m)$:
  dies ist aber wegen $m \equiv 5 \bmod{8}$ nicht der Fall. Weiter ist
  im Falle $q \equiv 7 \bmod{24}$ wegen $b = q-18 \equiv 13 \bmod{24}$:
  \[ \left(\frac{m}{q}\right) =
  \left(\frac{3}{q}\right)\left(\frac{b}{q}\right) =
  \left(\frac{b}{3}\right)\left(\frac{b}{q}\right) =
  \left(\frac{-18}{q}\right) = \left(\frac{-2}{q}\right) = -1
  \]
  und damit $-b$ keine Norm aus $D(m)$. Schließlich ist $a$ im Falle
  $q \equiv -1 \bmod{24}$ so gewählt, dass $b \equiv 3 \bmod{9}$ ist; mit
  $b$ wäre also auch $c = b/3$ Norm aus $D(m)$, aber
  \[
  \left(\frac{m}{c}\right)
    = \left(\frac{c}{m}\right) = \left(\frac{c}{q}\right)
    =  \left(\frac{3}{q}\right) \left(\frac{3c}{q}\right)
    = \left(\frac{b}{q}\right) = -1
  \]
  zeigt, dass dies nicht der Fall ist.
\end{proof}

\begin{quote}
  \textbf{(3.7)} {\em Sei $m=3q$, $q \equiv 3 \bmod{8}$; dann ist $D(m)$
    genau für $m=33$ und $m=57$ normeuklidisch.}
\end{quote}

\begin{proof}
  Der folgende Beweis, dass nur die Werte $m=33$ und $m=57$ in Frage
  kommen, geht auf Schuster\index[dN]{Schuster} \cite{Sch38} zurück. Gibt es
  ein $s \in \N$ mit $0 < 3s < q$, $s \equiv 1 \bmod{6}$ und $(s/q) =
  -1$, so ist $D(m)$ nicht normeuklidisch: setzt man nämlich $f=q$ und
  $a = 3s$, $b = q-3s$ für $q \equiv 1 \bmod{3}$, sowie $a = q-3s$, $b
  = 3s$ für $q \equiv 2 \bmod{3}$, dann gilt $(a/q) = +1$; weiter ist
  z.B. im ersten Falle mit $a$ auch $s = a/3$ Norm, aber
  $(\frac{q}{s}) = (\frac{3}{s})(\frac{9}{s}) = -1$ und auch
  $(\frac{q}{s}) = (\frac{3}{s})(\frac{9}{s}) = -1$. Wir zeigen nun
  die Existenz eines solchen $s$ für alle primen $q \ge 43$; für
  $q<108$ geben wir die Paare $(q,s)$ direkt an: $(43, 7)$, $(59,
  13)$, $(67, 7)$, $(83, 13)$, $(107, 7)$. Für die andern $q$
  unterscheiden wir
  \begin{enumerate}
  \item[(i)] $q \equiv 11 \bmod{24}$: wir wählen $c \in \N$ mit
    $21q < 108c < 22q$; dann liegen die beiden Zahlen $6c-q$ und $36c-7q$
    im Intervall $(0, q/3)$ und haben entgegengesetzten Restcharakter
    mod $q$.
  \item[(ii)] $q \equiv 19 \bmod{24}$: sei $15q < 108c < 16q$; hier genügt
    eine der beiden Zahlen $q-6c$ und $36-5q$ unseren Anforderungen.
  \end{enumerate}
\end{proof}

Damit bleiben im Falle von zusammengesetztem $m \equiv 1 \bmod{4}$ nur
noch solche zu untersuchen, für die $m=pq$, $p \equiv q \equiv 3
\bmod{4}$ prim, $7 < p < q$ (insbesondere $pq \ge 77$) zu
untersuchen. Hofreiters\index[dN]{Hofreiter} Beweis \cite{Hof35},
dass $D(77)$ nicht normeuklidisch ist, enthält jedoch einen Fehler,
der bisher anscheinend unbemerkt blieb (sh. z.B. van der
Linden\index[dN]{Linden@van der Linden} \cite[S. 15 unten]{Lin85}):
er behauptet, die Gleichung $77x^2 - Y^2 = -92$ sei in $\Z$ nicht
lösbar wegen $(77/23) = -1$; allerdings ist $(77/23) = +1$, und
$x=2, Y=20$ eine Lösung dieser Gleichung. Dass $D(77)$ in der Tat
nicht normeuklidisch ist, zeigen wir mit

\begin{quote}
  \textbf{(3.8)} {\em Sei $m=pq$, $p \equiv q \equiv 3 \bmod{4}$ prim;
    gibt es dann ein $r \in \N$ mit
    $(\frac{r}{p}) = -(\frac{r}{q}) = (\frac{p-r}{q})$
    und $1 \le r \le p$, so ist $D(m)$ nicht normeuklidisch.}
\end{quote}

\begin{proof}
  Wir verwenden (1.6) mit $f=p$ und $a=r$, $b=p-r$, falls
  $(\frac{r}{p}) = +1$, sowie $a=p-r$, $b=r$, falls
  $(\frac{r}{p}) = -1$ ist. Damit haben wir $f = a+b$,
  $(\frac{a}{p}) = +1$, und weiter
  $(\frac{m}{a}) = (\frac{pq}{r}) =
  (\frac{r}{p}) (\frac{r}{q}) = -1$,
  $(\frac{m}{b}) = (\frac{b}{m}) =
  (\frac{p-r}{q})(\frac{-r}{p}) = -1$, d.h. $a$
  und $-b$ sind nach (3.1) keine Normen aus $D(m)$.
\end{proof}

Als Korollar hat man sofort

\begin{quote}
  \textbf{(3.9)} {\em Ist $m=7q$, $q \equiv 11, 19 \bmod{24}$, so ist
    $D(m)$ nicht normeuklidisch. }
\end{quote}

\begin{proof}
  Setze $p=7$, $r=2$ in (3.8).
\end{proof}

Insbesondere ist also $D(77)$ nicht normeuklidisch. Um nun die
Klassifikation der normeuklidischen quadratischen Zahlkörper
abzuschließen, verwenden wir das Ergebnis von
Cassels\index[dN]{Cassels} \cite{Cas52}, nach dem quadratische
Zahlkörper $\Q(\sqrt{m})$ mit $m \in \N$ nur dann normeuklidisch sein
können, wenn $\disc K \le 2577$ gilt.  Wir lassen dann einen Computer
für die 92 Werte von $m=pq$ mit $p \equiv q \equiv 3 \bmod{4}$ prim,
$q \ge 7$, $p < q$, eine Lösung von (3.8) suchen und erhalten

\begin{quote}
  \textbf{(3.10)} {\em Sei $m=pq$, $p \equiv q \equiv 3 \bmod{4}$;
    dann ist $D(m)$ genau für $m=21, 33, 57$ normeuklidisch. }
\end{quote}

Zu guter Letzt brauchen wir noch ein Kriterium für prime
$m \equiv 1 \bmod{4}$; ein solches haben in der folgenden Form
Erdös\index[dN]{Erdos@Erd\"os} und Ko\index[dN]{Ko} \cite{EK38}
bereitgestellt:

\begin{quote}
  \textbf{(3.11)} {\em Sei $m = p \equiv 1 \bmod{4}$ prim; gibt es dann
    $r,s,t,u \in \N$ mit $p = rs + tu$, $(r,s) = (t,u) = 1$,
    $(\frac{r}{p}) = (\frac{s}{p}) = (\frac{t}{p}) = (\frac{u}{p}) = -1$,
    so ist $D(m)$ nicht normeuklidisch.}
\end{quote}

\begin{proof}
  (1.6) mit $a=rs$, $b=tu$.
\end{proof}

Wieder lässt man nun einen Computer suchen, und man findet, dass für
alle primen $m \equiv 1 \bmod{4}$ mit $m \ge 2577$ bis auf die
folgenden eine solche Darstellung existiert:
$$ m = 5, 13, 17, 29, 37,
41, 61, 73, 89, 97, 109, 113, 137, 193, 241, 313, 337, 457, 601. $$

Im Hinblick auf spätere Anwendungen in \S\ 5 wollen wir den Fall
$m \equiv 5 \bmod{24}$ ohne Benutzung der Schranke von Cassels behandeln:

\begin{quote}
  \textbf{(3.12)} {\em Ist $m \equiv p \equiv 5 \bmod{24}$ prim, so
    existiert für alle $p \ge 29$ eine Darstellung $p = rs + tu$ wie
    in (3.11). Insbesondere ist $D(m)$ für solche $m$ genau dann
    normeuklidisch, wenn $m=5$ oder $m=29$ ist.}
\end{quote}

\begin{proof}
  Wir suchen ein $s \in \N$ mit $0 < 3s < p$, $(s,3p)=1$, $s
  \equiv 1 \bmod{4}$ und $(\frac{s}{p}) = -1$. Haben wir
  ein solches $s$ gefunden, so ist $p-3s \equiv 2 \bmod{4}$, d.h. es
  ist $p = 3s+2u$ für ein $u \in \N$, $u \equiv 1 \bmod{2}$;
  und wir sind fertig.
  Ist $p \equiv 2 \bmod{5}$, so können wir $s=5$ nehmen; für die
  andern $p < 432$ geben wir die Paare $(p,s)$ direkt an: $(101, 29),
  (149, 13), (269, 29), (389, 29)$. Ist $p \ge 437$, so gibt es ein $c
  \in \N$ mit $2p-36 < 432c < 3p-36$. Eine der beiden Zahlen
  $12c+1$ und $p-144c-12$ erfüllt dann obige Bedingungen an $s$.
\end{proof}

Offenbar beruht der Beweis von (3.12) auf der Kenntnis von zwei
``kleinen'' quadratischen Nichtresten mod $p$, nämlich $t=2$ und $r=3$
für $p \equiv 5 \bmod{24}$. Im Falle $m \equiv 13 \bmod{24}$ kennt man
nur den Nichtrest $t=2$, sodass ein Analogon zu (3.12) in diesem Fall
viel schwieriger zu beweisen ist; Brauer\index[dN]{Brauer} \cite{Bra40}
hat 1940 die Existenz einer Darstellung (3.11) für prime
$p \equiv 13 \bmod{24}$, $p > 109$, auf kompliziertem, aber elementarem
Weg bewiesen.

Da $D(m)$ für $m = 5, 13, 17, 29, 37, 41, 73$ normeuklidisch ist, für
$m = 61$, $89$, $97$, $109$, $113$, $137$ aber nicht (sh. \S\ 1, Tabelle
\ref{dApp15} und \S\ 2, Tabelle \ref{dApp17} etc.), bleiben noch die Werte
$m = 193, 241, 337, 457$ und $601$ zu untersuchen. Dass der EA in
$D(m)$ für diese $m$ nicht gilt, zeigte zuerst Inkeri\index[dN]{Inkeri}
\cite{Ink47}, und zwar mit einer Methode, die mit (1.8) verwandt ist und
auf Redei\index[dN]{Redei@R\'edei} zurückgeht. Unabhängig davon haben dann
Chatland\index[dN]{Chatland} und Davenport\index[dN]{Davenport} \cite{Cha50}
dafür Beweise gegeben, und zwar kürzere als Inkeri; dabei benutzten
sie Davenports Methode, mit der dieser die Schranke $\disc K < 2^{14}$
für alle reellquadratischen normeuklidischen Zahlkörper bewiesen
hatte. Die Rechnungen von Chatland und Davenport lassen sich mit einem
Computer leicht nachkontrollieren (im Falle $m=601$ hat sich ein
Fehler eingeschlichen; sh. dazu Ennola\index[dN]{Ennola} \cite{Enn58b}).

Damit sind alle normeuklidischen quadratischen Zahlkörper bestimmt. Es
wäre allerdings wünschenswert, auch für die Werte $m = 193, \dots,
601$ ``einfache'' Beweise zu haben (z.B. nicht-euklidische Ideale
kleiner Norm in diesen Körpern). Dass sich zumindest einige dieser $m$
mit (1.8) nicht ausschließen lassen, liegt an den recht großen Werten
der jeweiligen Fundamentaleinheiten; so ist z.B. $u =
139\,468\,303\,679\,532 + 5\,689\,030\,769\,845\,601$ die FE in
$D(601)$.

Wir wenden uns nun der Anwendung von (1.15) und (1.16) in
reellquadratischen Zahlkörpern zu. Mit (1.15) lassen sich nur die
normeuklidischen Ringe $D(m)$ für $m = 2, 3, 5, 13$ finden, und zwar
mit der trivialen 1-Folge $0,1$ (für diese Ringe ist $M_{r,\beta} \le
\sqrt{3/2} < 2$). Zu bemerken ist noch, dass in allen $D(m)$ mit
positivem $m$ $\mu_1 = 2$ ist (d.h. $0, 1$ ist maximale 1-Folge) außer
für $D(5)$, wo $\mu_1 = 4$ ist.

Beim Auffinden von 2-euklidischen Zahlringen ist (1.16) etwas
erfolgreicher als im 1-euklidischen Fall. Ein Beispiel für eine
verhältnismäßig lange 2-Folge ist
\begin{align*}
  & 0, 1, 2, 3, 4, \sqrt{2}, 1+\sqrt{2}, 2+\sqrt{2}, 3+\sqrt{2},
  4+\sqrt{2}, 2\sqrt{2}, 1+2\sqrt{2}, 2+2\sqrt{2}, 3+2\sqrt{2}, \\
  & 4+2\sqrt{2}, 3\sqrt{2}, 1+3\sqrt{2}, 2+3\sqrt{2}, 3+3\sqrt{2},
  4+3\sqrt{2}, 5+4\sqrt{2}, 6+4\sqrt{2}, 7+4\sqrt{2}
  \end{align*}
in $D(2)$; diese 2-Folge zeigt $\mu_2 \ge 23$, während wegen (1.17)
$\lambda_2 \le 25$ ist (denn das Ideal $(5)$ besitzt kein PERS).

In der folgenden Tabelle sind 2-Folgen angegeben, die nachweisen, dass
die Ringe $D(m)$ für
\begin{align*}
  m & = 14, 23, 31, 43, 53, 61, 69, 77, 89, 93, 97, 113, 129, 133, 137, \\
    & \qquad 157, 161, 173, 193, 201, 213
\end{align*}
2-euklidisch sind.

\begin{table}[ht!]
  \centering
  \begin{tabular}{|r|l|} \hline
$m$ & $\theta_1, \dots, \theta_k$  \\ \hline
14 & $0, 1, 4+\sqrt{14}, 5+\sqrt{14}$ \\
23 & $0, 1, 2, 3+\sqrt{23}, 4+\sqrt{23}, 5+\sqrt{23}, 7+2\sqrt{23},
      8+2\sqrt{23}$ \\
31 & $0, 1, 2, 5+\sqrt{31}, 6+\sqrt{31}, 7+\sqrt{31}, 11+2\sqrt{31}$ \\
43 & $0, 1, 2, 8+\sqrt{43}, 14+2\sqrt{43}, 15+2\sqrt{43}, 16+2\sqrt{43}$ \\
53 & $0, 1, 2, 3, 4$ \\
61 & $0, 1, 2, 3, 4, 4+\beta, 6+\beta, 10+\beta$ \\
69 & $0, 1, 2, 2+\beta, 3+\beta, 4+2\beta, 5+2\beta, 6+2\beta$ \\
77 & $0, 1, 2, \beta, 1+\beta, 2+\beta$ \\
89 & $0, 1, 2, 3, 4$ \\
93 & $0, 1, 2, 3+\beta, 4+\beta, 5+2\beta, 6+2\beta$ \\
97 & $0, 1, 2, 3, 7+\beta, 11+2\beta, 12+2\beta, 13+2\beta, 14+2\beta$ \\
113 & $0, 1, 2, 3, 4, 13+2\beta$ \\
129 & $0, 1, 2, 11+2\beta, 12+2\beta, 13+2\beta, 16+3\beta$ \\
133 & $0, 1, 2, 6+\beta, 7+\beta, 8+\beta$ \\
137 & $0, 1, 2, 3, 4, 16+3\beta, 17+3\beta$ \\
157 & $0, 2, 3, 5+\beta, 8+\beta, 13+2\beta, 18+3\beta$ \\
161 & $0, 1, 5+\beta, 6+\beta, 10+2\beta, 11+2\beta, 12+2\beta$ \\
173 & $0, 2, 3, 6+\beta, 12+\beta, 18+2\beta$ \\
193 & $0, 2, 4, 5+\beta, 7+\beta, 9+\beta, 27+4\beta$ \\
201 & $0, 1, 2, 3, 6+\beta, 7+\beta, 8+\beta, 13+2\beta, 14+2\beta$ \\
213 & $0, 1, 2, 3, 5+\beta, 6+\beta, 7+\beta, 8+\beta, 10+2\beta, 11+2\beta$ \\
    & $12+2\beta, 13+2\beta, 17+3\beta, 18+3\beta, 19+3\beta, 20+3\beta$ \\
\hline \end{tabular}
  \caption{2-Folgen für 2-euklidische Ringe $D(m)$; \qquad
    ($\beta = (1+\sqrt{m})/2$)}
\end{table}

Manche hier aufgeführten 2-Folgen sind länger, als man zum Nachweis,
dass $D(m)$ 2-euklidisch ist, benötigt, andere lassen sich noch
verlängern. Alle diese Folgen wurden von Hand errechnet;
möglicherweise lassen sich auch noch andere quadratische Körper mit
(1.16) als 2-euklidisch nachweisen.

\section*{{\sc Anmerkungen zu} \S\ 3}

Um die Klassifikation aller normeuklidischen quadratischen Zahlkörper
hat sich eine große Anzahl von Mathematikern verdient gemacht (zwecks
Pluralbildung sh. Mentz\index[dN]{Mentz} \cite{Men88}); die nachfolgende
Tabelle soll einen Überblick geben:

\begin{table}[htbp]
\centering
\begin{tabular}{p{1cm} p{10cm}}
\toprule
\textbf{Jahr} & \textbf{Ereignis} \\
\midrule
1832 & Gauß\index[dN]{Gauss@Gauß} zeigt bei seiner Untersuchung der
       biquadratischen Reste,
       dass $D(-1)$ normeuklidisch ist. In seinem Nachlaß findet sich
       ein ähnlicher Beweis für $D(-3)$. \\
1847 & Wantzel\index[dN]{Wantzel} veröffentlicht den ersten Beweis
       für $D(-3)$. \\
1886 & Legendre\index[dN]{Legendre} gibt die kleinste Lösung von
       $x^2 - Ny^2 = \pm 1$, $N \le 1003$. \\
1893 & Dedekind\index[dN]{Dedekind} gibt einen Beweis für $D(-1)$ mit 
       dem Hinweis, dass sich die Fälle $D(m)$, 
       $m = -11$, $-7$, $-3$, $-2$, $2$, $3$, $5$, $13$
       ähnlich behandeln lassen. \\
1907 & Sommer\index[dN]{Sommer} gibt FE und Klassenzahl von
       $\Q(\sqrt{m})$, $-97 \le m \le 101$. \\
1927 & Dickson\index[dN]{Dickson} zeigt, dass die Ringe $D(m)$,
       $m = -11, -7, -3, -2, -1, 2, 3, 5, 13$ normeuklidisch sind und
       behauptet, dass es keine anderen normeuklidischen quadratischen
       Ringe gibt. \\
1928 & Schaffstein\index[dN]{Schaffstein} berechnet die Klassenzahl
       reellquadratischer Zahlkörper mit Primzahldiskriminante
       $\le 12000$. \\
1933 & Perron\index[dN]{Perron} weist den EA in den Ringen $D(m)$, 
       $m = 6$, $7$, $11$, $17$, $21$, $29$ nach und deckt damit den
       Irrtum Dicksons auf. Dabei deutet er an, dass möglicherweise
       jeder quadratische Zahlkörper mit Klassenzahl 1 normeuklidisch ist.
       I. Schur\index[dN]{Schur} teilt Perron daraufhin mit, dass
       $D(47)$ nicht normeuklidisch ist; dies scheint das erste Ergebnis
       in dieser Richtung zu sein. \\
1934 & In einem Brief an Perron zeigt Oppenheim,\index[dN]{Oppenheim} dass
       der EA außer in den bereits bekannten Fällen auch in $D(33)$, 
       $D(37)$ und $D(41)$ gilt, während $D(23)$, $D(31)$ und $D(53)$
       nicht norm-euklidisch sind. Unabhängig von Oppenheim weist
       Remak\index[dN]{Remak} den EA in $D(33)$, $D(37)$ und $D(41)$ nach. \\
1935 & Fox\index[dN]{Fox} und Berg\index[dN]{Berg} zeigen unabhängig 
       voneinander, dass $D(m)$ im Falle $m \equiv 2, 3 \bmod{4}$
       höchstens für $m \equiv 2, 3, 6, 7, 11, 19$
       normeuklidisch sein kann; Berg zeigt darüberhinaus, dass $D(19)$
       tatsächlich normeuklidisch ist. Hofreiter beweist, dass $D(57)$
       normeuklidisch ist und findet alle norm-euklidischen $D(m)$
       mit $m \equiv 21 \bmod{24}$. \\
1936 & Behrbohm\index[dN]{Behrbohm} und Redei\index[dN]{Redei@R\'edei}
       finden alle normeuklidischen $D(m)$ mit $m \equiv 5 \bmod{24}$. \\
1938 & Schuster\index[dN]{Schuster} findet alle normeuklidischen $D(m)$ mit
       $m \equiv 9 \bmod{24}$; Erdös\index[dN]{Erdos@Erd\"os} und
       Ko\index[dN]{Ko} zeigen, dass es nur endlich
       viele prime $p$ mit $p \equiv 1 \bmod{8}$ gibt, sodass $D(p)$
       normeuklidisch ist. Heilbronn\index[dN]{Heilbronn} zeigt dies auch für
       zusammengesetzte $m \equiv 1 \bmod{8}$. \\ 
1940 & Brauer\index[dN]{Brauer} zeigt, dass der EA in
       $D(m)$, $m \equiv 13 \bmod{24}$, nur bestehen kann, wenn
       $m \le 109$ gilt. \\
\bottomrule
\end{tabular}
\caption{Historische Entwicklung der Klassifikation normeuklidischer
  quadratischer Zahlkörper}
\end{table}

\begin{table}[htbp]
\centering
\begin{tabular}{p{1cm} p{10cm}}
\toprule
\textbf{Jahr} & \textbf{Ereignis} \\
\midrule
1942 & Redei\index[dN]{Redei@R\'edei} zeigt, dass $D(73)$ normeuklidisch 
       ist, findet alle normeuklidischen $D(m)$ mit $m \equiv 17 \bmod{24}$
       und schließt den EA in $D(m)$, $m = 61, 89, 109, 113, 137$ aus.
       Damit bleibt nur noch die Restklasse $m \equiv 1 \bmod{24}$
       unerledigt. \\
1944 & Hua\index[dN]{Hua} zeigt, dass für normeuklidische Körper
       $\Q(\sqrt{p})$, $p \equiv 1 \bmod{4}$ prim, die Ungleichung
       $\disc K < e^{250}$ besteht. Zusammen mit Min\index[dN]{Min}
       zeigt er, dass jedes
       prime $p \equiv 17 \bmod{24}$ mit $p > 137$ eine Darstellung
       $p = rs + tu$ besitzt, wobei $(r,s) = (t,u) = 1$ und
       $(r/p) = (s/p) = (t/p) = -1$ gilt. \\
1947 & Inkeri\index[dN]{Inkeri} zeigt, dass es keine normeuklidischen
       $\Q(\sqrt{m})$ mit $\disc K < 5000$ außer den bereits bekannten gibt. \\
1948 & Davenport\index[dN]{Davenport} zeigt, dass $\disc K \le 16384 = 2^{14}$
       für reellquadratische, normeuklidische Körper gilt. \\
1949 & Chatland\index[dN]{Chatland} findet für jedes prime $p \equiv 1 \bmod{4}$,
       $601 < p \le 16384$, eine Darstellung (3.11). Zusammen mit
       dem Ergebnis von Inkeri\index[dN]{Inkeri} (1947) heißt das, dass
       alle normeuklidischen quadratischen Zahlkörper bekannt sind. \\
1950 & Chatland und Davenport\index[dN]{Davenport} behandeln die Körper mit
       $193 \le \disc K \le 601$ (anders als Inkeri und offenbar in
       Unkenntnis seiner Arbeit). \\
1952 & Barnes\index[dN]{Barnes} und Swinnerton-Dyer\index[dN]{Swinnerton-Dyer}
       zeigen, dass $D(97)$ -- entgegen einer
       Behauptung von Redei aus dem Jahre 1942 -- nicht normeuklidisch ist.
       Varnavides\index[dN]{Varnavides} gibt einen einheitlichen Beweis
       dafür, dass die Körper $\Q(\sqrt{m})$ mit
       $m = 2$, $3$, $5$, $6$, $7$, $11$, $13$, $17$, $19$, $21$,
       $29$, $33$, $37$, $41$, $57$, $73$  normeuklidisch sind. \\
1958 & Ennola\index[dN]{Ennola} zeigt noch einmal, dass $\Q(\sqrt{m})$
       genau für die  oben angegebenen Werte normeuklidisch ist. \\
1977 & Cooke\index[dN]{Cooke} findet 20 2-stufig normeuklidische
       quadratische Zahlkörper. \\
1985 & Johnson,\index[dN]{Johnson} Queen\index[dN]{Queen} und
       Sevilla\index[dN]{Sevilla} zeigen, dass $D(10)$ und $D(65)$
       semi-euklidisch sind. \\
\bottomrule
\end{tabular}
\caption{Historische Entwicklung der Klassifikation normeuklidischer
  quadratischer Zahlkörper (Fortsetzung)}
\end{table}

\chapter*{\S\ 4 Kubische Zahlkörper}
\setcounter{chapter}{4}
\addcontentsline{toc}{chapter}{\S\ 4 Kubische Zahlkörper}
\markboth{Euklidische Ringe}{\S\ 4 Kubische Zahlkörper}

Während ein quadratischer Zahlkörper durch seine Diskriminante
eindeutig bestimmt ist, ist dies bei kubischen Zahlkörpern nicht mehr
der Fall: so haben z.B. die beiden Körper $\Q(\sqrt[3]{6})$
und $\Q(\sqrt[3]{2})$ dieselbe Diskriminante $d=972$, sind
aber nicht isomorph. Dennoch kann man allein an der Diskriminante
einige Eigenschaften des jeweiligen Körpers ablesen. Da für das
Vorzeichen der Diskriminante eines Körpers mit $r$ reellen und $2s$
nichtreellen Einbettungen von $K$ in $\C$ die Beziehung
$\operatorname{sign}(\disc K) = (-1)^s$ gilt und in
kubischen Körpern nur die beiden Möglichkeiten $r=3$, $s=0$ und $r=1$,
$s=1$ auftreten, ist $\disc K$ genau dann positiv, wenn
$K$ total reell ist und Einheitenrang $2$ hat. Dagegen ist
$\disc K$ negativ, wenn $r=s=1$ ist und $K$
Einheitenrang $1$ hat.

Ist $K$ ein Zahlkörper mit GHB $\{\alpha_1, \dots, \alpha_n\}$, und
bezeichnen $\sigma_1, \dots, \sigma_n$ die Einbettungen von $K$ in
$\C$, so gilt bekanntlich $\disc K = \det
(\sigma_i(\alpha_j))^2$. Also ist $\disc K$ das Quadrat
einer Zahl, die im Galoisabschluss $L$ von $K/\Q$ liegt. Dies
zeigt: ist $K/\Q$ galoissch, so muss $\sqrt{\disc
  K} \in R$ sein, wo $R$ wie üblich den Ring ganzer Zahlen in $K$
bezeichnet.

Beispiel: Ist $K = \Q(\zeta_p)$ der Körper der $p$-ten
Einheitswurzeln für ein primes $p\in\N$, so ist
$\disc K = (-1)^{(p-1)/2}p^{p-2}$, setzt man nun $p^* = (-1)^{(p-1)/2}p$,
so hat man, weil $K/\Q$ galoissch ist,
$\Q(\sqrt{p^*}) = \Q(\sqrt{\disc K}) \subset \Q(\zeta_p)$.

Ist nun $K/\Q$ galoissch und $(K:\Q)$ ungerade, so kann
$K$ keinen quadratischen Teilkörper besitzen, und $\sqrt{\disc K} \in
R$ impliziert dann $\sqrt{\disc K} \in \Z$. Im Falle
$(K:\Q) = 3$ gilt sogar die Umkehrung:

\begin{quote}
  {\bf (4.1)} {\em Ein kubischer Zahlkörper $K$ ist genau dann
    galoissch, wenn $d = \disc K$ ein Quadrat in
    $\Z$ ist. Ist $K$ nicht galoissch, so ist $L =
    \Q(\sqrt{d})$ der Galoisabschluss von $K/\Q$.}
\end{quote}

\begin{proof}
  Es genügt zu zeigen, dass $K(\sqrt{d})$ galoissch ist. Mit $K =
  \Q(\alpha)$ ist aber (wegen $\disc
  K/\Q(\alpha) = k^2 d$ für ein $k\in\Z$) $k\sqrt{d} =
  (\alpha-\alpha')(\alpha' - \alpha'')(\alpha'' - \alpha')$, wo
  $\alpha'$, $\alpha''$, $\alpha'''$ die Konjugierten von $\alpha$
  bezeichnen. Ist $f(x) = x^3+a_2x^2+a_1x+a_0 \in \Z[x]$ das
  Minimalpolynom von $\alpha$, so rechnet man leicht nach, dass
  $\alpha' - \alpha'' = \pm k\sqrt{d}/f'(\alpha) = \pm
  k\sqrt{d}/(3\alpha^2 + 2a_2\alpha + a_1)$, sowie $\alpha' + \alpha''
  = a_2 - \alpha$ gilt. Damit ist $\alpha'$, $\alpha'' \in
  K(\sqrt{d})$, also $L$ normal über $\Q$.
\end{proof}

Wir wollen zuerst den nicht-galoisschen Fall untersuchen: sei dazu
$k_3$ kubischer Körper mit Diskriminante $d$, $d$ kein Quadrat in
$\Z$, $k_2 = \Q(\sqrt{d})$ ein quadratischer
Zahlkörper und $K=k_2k_3$ der normale Abschluss von
$k_3/\Q$. Dann ist $K/\Q$ galoissch mit Galoisgruppe
$S_3$ (die symmetrische Gruppe der Ordnung 6). Wir stellen jetzt die
Hilbertsche Untergruppenreihe für $K/\Q$ auf; aus dieser
Tafel und dem Zerlegungsgesetz in quadratischen Zahlkörpern lesen wir
sofort ab:

\begin{quote}
{\bf (4.2)} {\em Sei $k_3$ kubischer Zahlkörper mit Diskriminante $d$;
  dann gilt für alle primen $p \in \N$:}
\begin{align*}
  \Big(\frac dp\Big) = +1 &
      \Leftrightarrow p \text{ ist voll zerlegt oder träge in } k_3 \\  
  \Big(\frac dp\Big) = -1 &
      \Leftrightarrow p \text{ ist Produkt zweier Primideale in } k_3.
\end{align*}
\end{quote}

Für abelsche Körper ist dies trivialerweise richtig, weil dann
einerseits $d$ ein Quadrat und somit nie $(d/p) = -1$ ist, und weil
andererseits $p$ nie das Produkt zweier Primideale wird.

\textbf{Hilbertsche Untergruppenreihe für $K = k_2 k_3$}
Hierbei ist
$k_3$ ein nichtabelscher kubischer Zahlkörper mit Diskriminante $d$
$d = d_2 f^2$, $d_2 = \disc k_2$, $k_2 =
\Q(\sqrt{d}) = \Q(\sqrt{d_2})$
$K = k_2 k_3$ der normale Abschluss von $k_3$.

\begin{table}[htbp]
\centering
\begin{tabular}{c|c|c|c|c|c|c|c|c|c|c}
\toprule
$(e,f,g)$ & $p$ & $k_2$ & $k_3$ & $\Q$ & $K_{\Z}$
       & $K_T$ & $K_1$ & $K_2$ & $K_3$ & $K$ \\
\midrule
$(1,1,6)$ & $(d/p)=+1$ & $(1,1)$ & $(1,1,1)$ & $\Q$
       & $K$ & $K$ & $K$ & $K$ & $K$ & $K$ \\
$(1,3,2)$ & $(d/p)=+1$ & $(1,1)$ & $(3)$ & $\Q$ & $k_2$
       & $K$ & $K$ & $K$ & $K$ & $K$ \\
$(1,2,3)$ & $(d/p)=-1$ & $(2)$ & $(1,2)$ & $\Q$ & $k_3$
       & $K$ & $K$ & $K$ & $K$ & $K$ \\
$(2,1,3)$ & $p=2$, $p+1$, $p\|d_2$ & $(1^2)$ & $(1,1^2)$ & $\Q$
       & $k_3$ & $k_3$ & $K$ & $K$ & $K$ & $K$ \\
& $p=2$, $2+1$, $4\|d_2$ & & & $\Q$ & $k_3$
       & $k_3$ & $k_3$ & $K$ & $K$ & $K$ \\
& $p=2$, $2+1$, $8\|d_2$ & & & $\Q$ & $k_3$
       & $k_3$ & $k_3$ & $k_3$ & $K$ & $K$ \\
$(3,1,2)$ & $p=3$, $p\|f$, $(d_2/p)=+1$ & $(1,1)$ & $(1^3)$
       & $\Q$ & $k_2$ & $k_2$ & $K$ & $K$ & $K$ & $K$ \\
& $p=3$, $9\|f$, $(d_2/3)=+1$ & & & $\Q$ & $k_2$
       & $k_2$ & $k_2$ & $K$ & $K$ & $K$ \\
$(3,2,1)$ & $p=3$, $p\|f$, $(d_2/p)=-1$ & $(2)$ & $(1^3)$
      & $\Q$ & $\Q$ & $k_2$ & $K$ & $K$ & $K$ & $K$ \\
& $p=3$, $9\|f$, $(d_2/3)=-1$ & & & $\Q$
        & $\Q$ & $k_2$ & $k_2$ & $K$ & $K$ & $K$ \\
$(6,1,1)$ & $p=3$, $3\|f$, $3\|d_2$ & $(1^2)$ & $(1^3)$
     & $\Q$ & $\Q$ & $\Q$ & $k_2$ & $K$ & $K$ & $K$ \\
& $p=3$, $9\|f$, $3\|d_2$ & & & $\Q$ & $\Q$
      & $\Q$ & $k_2$ & $k_2$ & $k_2$ & $K$ \\
\bottomrule
\end{tabular}
\caption{Hilbertsche Untergruppenreihe}
\end{table}

Hierbei bedeutet z.B. $(1,2)$, dass $p$ in $k_3$ Produkt zweier
Primideale mit Trägheitsgrad 1 bzw. 2 ist, während $(1,1^2)$ die
Zerlegung $(p) = P_1 P_2^2$ symbolisiert.

Eine genauere Untersuchung der Faktorisierung des Ideals $(2)$ zeigt

\begin{quote}
  {\bf (4.3)} {\em Sei $k_3$ kubischer Zahlkörper mit Diskriminante $d$;
    dann bestehen nur die folgenden Möglichkeiten:}
\[
\begin{array}{ll}
d \equiv 8 \bmod {16} : & (2) = 2_1 \ftw_2^2 \quad (\text{in } k_3) \\
d \equiv 4 \bmod {16} : & (2) = 2_1^3 \\
d \equiv -4 \bmod {16} : & (2) = 2_1 \ftw_2^2 \\
d \equiv 1 \bmod {8} : & (2) = 2_1 \ftw_2 2_3 \quad \text{{\em oder }} (2) = (2) \\
d \equiv 5 \bmod {8} : & (2) = 2_1 \ftw_2
\end{array}
\]
\end{quote}

Der Beweis beruht auf der Berechnung der verschiedenen
Relativdiskriminanten unter Berücksichtigung obiger Tafel.

Wenn wir nun einen kubischen Körper $k_3$ mit Diskriminante $d_3 = d_2
\cdot f^2$ gegeben haben und uns fragen, ob er normeuklidisch ist, so
können wir nur dann mit (1.5) arbeiten, wenn es in $k_3$ rein
verzweigte Primideale gibt. Da genau die Primteiler von $f$ rein
verzweigen, ist (1.5) für Körper mit $f=1$ nutzlos. Wir wollen daher
folgendermaßen vorgehen: wir betrachten bei festem $d_2$ alle
kubischen Körper $k_3$ mit $d_3 = d_2 f^2$; dabei beginnen wir mit dem
einfachsten Fall $d_2 = -3$.

Kubische Zahlkörper mit Diskriminante $d_3 = -3f^2$ sind genau die
``reinen'' Zahlkörper $k_3 = \Q(\sqrt[3]{m})$ (sh. z.B. Delone und
Faddeev\index[dN]{Faddeev} \cite{Del64} oder Cohn\index[dN]{Cohn}
\cite[Übung 18.4]{Coh78}). Wir notieren einige einfache Eigenschaften:

Sei $m = ab^2$ für gewisse $a,b \in \Z$ mit $(a,b)=1$, $a$ und
$b$ quadratfrei, $\theta = \sqrt[3]{ab^2}$, $\theta' =
\sqrt[3]{a^2b}$, $K = \Q(\theta) = \Q(\theta')$;
dann gilt im Falle
\[
m \equiv 1 \bmod {9} : \{1, \theta, \theta'\} \text{ ist GHB, }
\disc K = -27a^2b^2;
\]
\[
ab^2 \equiv a^2b \equiv 1 \bmod {9} : \{1, \theta,
(1+\theta+\theta')/3\} \text{ ist GHB, } \disc K =
-3a^2b^2.
\]

Das Zerlegungsgesetz in $k_3$ lautet
\[\begin{array}{llll}
\text{a)} & p \nmid \disc K : & (p) = P_1 P_2 \quad
& \Leftrightarrow \quad p \equiv 2 \bmod {3} \\
& & (p) = P_1 P_2 P_3
\quad & \Leftrightarrow \quad p \equiv 1 \bmod {3},\ x^3 \equiv m \text{
  lösbar in } \Z; \\
& & (p) = (p) \quad & \Leftrightarrow \quad
p \equiv 1 \bmod {3},\ x^3 \equiv m \text{ nicht lösbar in } \Z; \\
\text{b)} & p \mid \disc K : & (p) = P^3 \quad
& \Leftrightarrow \quad p \neq 3 \text{ oder } p=3,\ m \not\equiv \pm 1
\bmod {9} \\
& & (p) = P_1 P_2^2 \quad & \Leftrightarrow \quad p=3 \text{
  und } m \equiv \pm 1 \bmod {9}
\end{array}
\]

Schließlich ist $N_{K/\Q}(1+\theta+\theta' +
\theta^2) = r^3 + m s^3 + m^2 t^3 - 3 m r s t$.

Da Cassels\index[dN]{Cassels} \cite{Cas52} gezeigt hat, dass ein kubischer
Zahlkörper mit negativer Diskriminante nur dann normeuklidisch sein
kann, wenn $|\disc K| < 420^2$ ist, brauchen wir nur solche zu
betrachten. Die Tafeln von Nakamula\index[dN]{Nakamula} \cite{Nak88}
enthalten alle solchen Körper; aus ihnen liest man ab, dass unter
ihnen genau die folgenden $\Q(\sqrt[3]{m})$ Klassenzahl 1 haben:\label{dpKh1}
$\Q(\sqrt[3]{m})$ für
\begin{align*}
  m & = 2, 3, 5, 6, 10, 12, 17, 23, 29, 33, 41, 44, 45, 46, 53, 55, 59,
        69, 71, 82, 99, \\
    & \quad 107, 116, 145, 179, 188, 197, 226, 332, 404, 575.
\end{align*}

Es wird sich zeigen, dass unter diesen Körpern genau diejenigen mit $m
= 2, 3, 10$ normeuklidisch sind; dieses Ergebnis stammt von
Cioffari\index[dN]{Cioffari} \cite{Cio79}, dessen Arbeit wir bei unserem
Beweis folgen werden. Der Teil von Cioffaris Arbeit, der sich mit
kubischen Zahlkörpern befasst, enthält folgende Druckfehler:
\begin{enumerate}
\item Seite 392, Tabelle zum Korollar von Prop.\ 3: in der Zeile mit
  $d = 107$ muss in der Spalte $f-e$ statt $91$ die Zahl $93$ stehen;
\item Seite 392, Prop.\ 5: es muss ``$N(u) \equiv (-25)^3 \equiv +10
  \bmod {53}$'' heißen statt ``$\dots \equiv -10 \bmod {53}$'';
  entsprechend müssen die Vorzeichen in den beiden darauffolgenden
  Zeilen geändert werden;
\item Seite 393, Prop.\ 7: statt ``belongs to $\theta (2)$; hence
  $a \equiv \theta \bmod {2}$'' muss es ``belongs to $\theta^2(2)$;
  hence $a \equiv \theta^2 \bmod {2}$'' heißen;
\item Seite 395, Prop.\ 10: in der Zeile $d=44$ der Tabelle muss $c =
  5+2\cdot \theta + \theta^2/2$ stehen statt
  $c = 5+2\theta+\varphi$.
\item Seite 396, Prop.\ 12: statt $(a_1-b_1c_1)(a_2-b_1c_1)$ sollte
  $(a_1 - \sqrt{b_ic_j})(a_2 - \sqrt{b_ic_j})$ stehen, und
  Entsprechendes gilt für die beiden darauffolgenden Zeilen.
\end{enumerate}

Nun gilt

\begin{quote}
  {\bf (4.4)} {\em Sei $(K:\Q) = n\equiv 1 \bmod {2}$, und seien die
    rationalen Primzahlen $p_1, \dots, p_t$ paarweise verschieden und
    rein verzweigt in $K$; außerdem sei $(p_1-1,n) = \dots = (p_t-1,n) = 1$.

    Gibt es dann ein $e \in \N$ mit $1 < e < f = p_1\cdots p_t$,
    sodass weder $e$ noch $f-e$ Norm eines Elements aus $R$ sind,
    dann ist $R$ nicht normeuklidisch.}
\end{quote}

\begin{proof}
  Die Voraussetzungen garantieren $(\varphi(f),n) = 1$; somit ist die
  Potenzierung mit $n$ ein Automorphismus auf $(\Z/f\Z)^*$ und folglich
  $e$ ein $n$-ter Potenzrest mod $f$. Die Behauptung folgt nun mit $a = e$,
  $b=f-e$ aus (1.6), denn wegen $(K:\Q) \equiv 1 \bmod {2}$
  ist $b$ genau dann Norm, wenn $-b$ Norm ist.
\end{proof}

\begin{quote}
  {\bf (4.5)} {\em Für die folgenden Werte von $m$ ist
    $\Q(\sqrt[3]{m})$ nicht normeuklidisch:} $m = 23, 29, 33,
    41, 46, 59, 69, 71, 82, 107, 188, 197, 226, 332, 404, 575$.
\end{quote}

\begin{proof}
  Wir verwenden (4.4) mit $t=1$ und $t=2$; vgl. die Tabelle.
\end{proof}

\begin{table}[htbp]
\centering
\begin{tabular}{cccc|ccccc}
\toprule
$m$ & $p_1$ & $e$ & $f-e$ & $m$ & $p_1$ & $p_2$ & $e$ & $f-e$ \\
\midrule
59 & 59 & 7 & 52 & 23 & 3 & 23 & 13 & 56 \\
71 & 71 & 19 & 52 & 29 & 3 & 29 & 26 & 61 \\
82 & 41 & 13 & 28 & 33 & 3 & 11 & 7 & 26 \\
107 & 107 & 14 & 93 & 41 & 3 & 41 & 19 & 104 \\
179 & 179 & 7 & 172 & 46 & 2 & 23 & 7 & 39 \\
197 & 197 & 39 & 158 & 69 & 3 & 23 & 26 & 43 \\
226 & 113 & 37 & 76 & 116 & 2 & 29 & 21 & 37 \\
332 & 83 & 7 & 76 & 145 & 5 & 29 & 26 & 119 \\
404 & 101 & 28 & 73 & 188 & 2 & 47 & 37 & 57 \\
 & & & & 575 & 5 & 23 & 37 & 78 \\
\bottomrule
\end{tabular}
\end{table}

Beispielsweise ist $ a = 14$ für $m = 107 = 14 + 93$ keine Norm aus $R$,
weil es kein Ideal der Norm $7$ gibt: die Kongruenz $x^3 \equiv 107 \bmod 7$
ist nämlich nicht lösbar, und nach dem Zerlegungsgesetz bleibt (7)
träge in R. Ebenso ist b = 93 keine Norm, weil es in R kein Ideal der
Norm 31 gibt.

Die jetzt noch verbliebenen Möglichkeiten $m = 2, 3, 5, 6, 10, 12, 17,
44, 45, 53, 55, 99$ lassen sich mit (4.4) allein nicht
entscheiden. Wir werden daher in diesen Körpern nach Idealen $I$ mit
$M(K,I) > 1$ suchen. Besonders geeignet sind hierfür Teiler des Ideals
$(u-1)$, $u$ Einheit in R. Diesbezüglich gilt

\begin{quote}
  {\bf (4.6)} {\em Sei $K = \Q(\sqrt[3]{m})$, $m \equiv 1 \bmod 2$
    und $3 \nmid h$, wo $h = h(K)$ die Klassenzahl von $K$ ist. Gibt es
    dann eine Primzahl $p \neq m$, die in $K$ rein verzweigt, so ist
    $u \equiv 1 \bmod 2$ für jede Einheit $u$ in $R$.}
\end{quote}

Bem.: Wegen $(2) = \ftw_1 \ftw_2$, $\|\ftw_2\| = 4$, kommt der Teiler
$\ftw_2$ von $(u-1)$ als Kandidat für $I$ in Frage.

\begin{proof}
  Sei $(p) = P^3$ in R; wegen $(h,3) = 1$ ist $P$ ein Hauptideal in R,
  z.B. $P = (\pi)$ für ein $\pi \in R$. Damit wird $\pi^3 = \pm p u^k$
  für ein $k \in \Z$, wo $u$ die FE von R ist. Indem man mit
  einer geeigneten Potenz von $u$ multipliziert, kann man $k \in \{0,
  1\}$ erreichen. Wäre $k=0$, so folgte $\pi = \sqrt[3]{p} \in K$, was
  wegen $p \neq m$ nicht der Fall ist. Indem wir notfalls $u$ durch
  $u^{-1}$ ersetzen, dürfen wir schließlich $k=1$ annehmen. Nun ist
  $(2) = \ftw_1 \ftw_2$ in $R$, d.h. $\Phi(2) = \Phi(\ftw_1)\Phi(\ftw_2) = 3$,
  also $\pi^3 \equiv 1 \bmod 2$ wegen $(\pi, 2) = 1$. Damit ist dann $u
  \equiv u^k p \equiv \pi^3 \equiv 1 \bmod 2$.
\end{proof}

Sei nun $m \equiv 1 \bmod 2$, $K = \Q(\sqrt[3]{m})$; dann ist
$(2) = \ftw_1 \ftw_2$ mit $\|\ftw_1\| = 2$, $\|\ftw_2\| = 4$. Setzt man $\theta
= \sqrt[3]{m}$, so ist $\{1, \theta, 1+\theta\}$ ein primes
Restsystem mod $\ftw_2$, und es gilt $\theta^2 \equiv 1 + \theta
\bmod \ftw_2$. Also ist $\ftw_2$ genau dann normeuklidisch, wenn die
Restklassen $\theta, 1+\theta \bmod \ftw_2$ Elemente der Norm $< 4$
enthalten. Wegen (4.6) enthalten diese Restklassen sicher keine
Einheiten, sodass hierfür nur Elemente der Norm 2 oder 3 in Frage
kommen.
\begin{itemize}
  \item[] $m = 5$: hier ist $\ftw_1 = (3-\theta^2)$, $(3) = Q^3$ mit $Q =
    (2-\theta)$, folglich gilt für alle Elemente $\alpha$ der Norm 2:
    $\alpha \equiv 3 - \theta^2 \equiv 1 + \theta^2 \equiv \theta
    \bmod \ftw_2$, für alle $\beta$ mit $N_{K/\Q}(\beta)=3$: $\beta
    \equiv 2 - \theta \equiv \theta \bmod \ftw_2$. Also enthält die
    Restklasse $1+\theta \bmod \ftw_2$ nur Elemente der Norm $\ge 5$
    (sogar $\ge 6$ wegen $\theta_1 = (\theta)$).
\item[] $m = 45$: mit $\alpha = \theta^2/\theta$ ist $\{1, \theta,
  \alpha\}$ eine GHB, und man findet $\ftw_1 = (22 - 5\theta - \alpha)$,
  $\ftw_2 = (409 + 155\theta + 97\alpha)$, $\theta_1 = (12 - \theta -
  2\alpha)$. Auch hier enthält die Restklasse $1+\theta \bmod \ftw_2$
  keine Elemente der Norm 2 oder 3 (man beachte $\alpha \equiv
  \theta^2 \equiv 1+\theta \bmod \ftw_2$).
\item[] $m = 55$: mit $\alpha = (1+\theta+\theta^2)/3$ ist $\{1,
  \theta, \alpha\}$ GHB, $\alpha \equiv 0 \bmod \ftw_2$ und $\ftw_1 =
  (341+71+76\alpha)$, $(3) = 3_1 3_2^2$, $3_1 = (5-3\theta + \alpha)$,
  $3_2 = (5+\theta + \alpha)$, sodass die Restklasse $\theta \bmod
  \ftw_2$ keine Elemente der Norm $<4$ enthält.
\item[] $m = 99$: Sei $\alpha = \theta^2/3$; dann wird $\{1, \theta,
  \alpha\}$ GHB, $\ftw_1^2 = (16-5\theta + \alpha)$, $(3) = 3_1^3$ und
  $16-5\theta + \alpha \equiv 1 \bmod \ftw_2$. Ist daher $\pi$ ein
  Element der Norm $2$ in $R$, so ist $\pi^2 e = 16 - 5\theta + \alpha
  \equiv 1 \bmod \ftw_2$ für eine Einheit $e$. Wegen $e \equiv 1 \bmod 2$
  und $\Phi(\ftw_2)=3$ ist damit schon $\pi \equiv 1 \bmod \ftw_2$. Egal, in
  welcher Restklasse die Erzeugenden von $3_1$ liegen, kann $\ftw_2$
  nicht normeuklidisch sein.
\item[] $m = 6$: hier ist $u = 1-6\theta + 3\theta^2$ eine FE von $K$,
  die Ideale $(2)$ und $(3)$ sind rein verzweigt, es ist $\ftw_1 =
  (2-\theta)$, $3_1 = (3+2\theta + \theta^2)$, und die Restklasse
  $1+\theta \bmod \ftw_1$ enthält keine Elemente der Norm $<4$.
\item[] $m = 53$: mit $\alpha = (1-\theta + \theta^2)/3$ ist $\{1,
  \theta, \alpha\}$ eine GHB, und wir haben $(2)=\ftw_1 \ftw_2$, $\|\ftw_1\| =
  2$, $\|\ftw_2\| = 4$, und $(53) = P^3$. Wir behaupten, dass die
  Restklasse $-25 \bmod \ftw_1 P$ keine Elemente der Norm $<106$
  enthält. Ist nämlich $\pi \in R$, $\pi \equiv -25 \bmod P$, so
  folgt, da $P$ rein verzweigt ist, $N_{K/\Q}(\pi) = (-25)^3 = 10
  \bmod 53$. Also ist $N_{K/\Q}(\pi) \in \{-96, -43, 10, 63\}$. Nun
  sind die Zahlen $96$ und $10$ durch eine ungerade Potenz von $2$
  teilbar, und dies impliziert $\pi \equiv 0 \bmod \ftw_1$ im Widerspruch
  zu $\pi \equiv -25 \equiv 1 \bmod \ftw_1$. Da weiter $(7)$ und $(43)$
  in $K$ träge sind, sind auch $63$ und $-43$ keine Normen aus
  $R$. Dies war zu zeigen.
\end{itemize}
Damit sind von der Liste von S. \pageref{pKh1} nur noch die Werte
$m = 2, 3, 10, 12, 17, 44$ übrig; für $m = 12, 17, 44$ können wir den
EA mit Hilfe von (1.10) ausschließen: dazu geben wir für diese $m$ eine
(fundamentale) Einheit $u$, die aus ihr resultierenden Schranken
$\mu_1, \mu_2, \mu_3$ aus (1.10) (zunächst mit $k=1$), ein $x \in K$
und schließlich $M(K,x)$:

\begin{table}[h]
\centering
\begin{tabular}{c|ccccc}
\toprule
$m$ & $u^{-1}$ & $|u|$ & $\mu_1$ & $\mu_2$ & $\mu_3$ \\
\midrule
12 & $1 + 3\theta - 3\theta^2$ & $\approx 165$ & 5.5 & 2.4 & 1.1 \\
17 & $18 - 7\theta$ & $\approx 972$ & 9.91 & 3.9 & 1.5 \\
44 & $(113 - 29 - 17\theta^2)/3$ & $\approx 4007$ & 15.89 & 4.5 & 1.28 \\
\bottomrule
\end{tabular}
\end{table}

\begin{table}[h]
\centering
\begin{tabular}{c|cc}
\toprule
$m$ & $x$ & $M(K,x)$ \\
\midrule
12 & $(12 + 15\theta + 8\theta')/18$ & $169/162$ \\
17 & $(-376 + 466\theta - 19\phi)/1028$ & $1115/1028$ \\
44 & $(8 + 12\theta + 36\theta')/59$ & $81/59$ \\
\bottomrule
\end{tabular}
\end{table}

Hierbei ist $\alpha = (1 - \theta + \theta^2)/3$; außerdem gilt für
das bei $m = 44$ angegebene $x$ die Kongruenz
$x \equiv 12/(5 + 2\cdot \theta + \theta') \bmod R$
(sh. Berichtigung 4. vor (4.4)). Damit haben wir

\begin{quote}
  {\bf (4.7)} {\em $K = \Q(\sqrt[3]{m})$ ist genau für
    $m = 2, 3, 10$ normeuklidisch.}
\end{quote}

Dass diese Körper tatsächlich normeuklidisch sind, haben Godwin ($m =
2$), Taylor ($m = 3, 10$) und auch Cioffari ($m = 2, 3, 10$) gezeigt;
dies kann man mit einem Computer leicht nachprüfen. Darüberhinaus kann
man für $K = \Q(\sqrt[3]{2})$ das euklidische Minimum $M(K) = \frac12$
bestimmen.

Ich habe mich (mit Hilfe der Klassenkörpertheorie) davon überzeugt,
dass von den Körpern der Form $d = -4f^2$ höchstens diejenigen mit $f =
9, 11$ und $83$ normeuklidisch sind; für die ersten beiden hat Taylor
den EA nachgewiesen. Für den Körper mit $f=83$ ist es ohne Kenntnis
der erzeugenden Gleichung recht schwer, weitere Aussagen zu
machen. Ich habe deswegen auch auf eine genaue Darstellung der nötigen
Rechnungen verzichtet.

In den allermeisten Fällen wird man zur Entscheidung, ob ein kubischer
Körper normeuklidisch ist oder nicht, wohl (1.10) benutzen; ohne eine
Tafel kubischer Körper, die neben Diskriminante auch Grundeinheit und
Klassenzahl angibt, kommt man daher nicht weiter.

Ich habe versucht, für kubische Körper mit kleiner Diskriminante
($|d|<1300$) eine solche Tafel aufzustellen (sh. IV); diese enthält
neben der erzeugenden Gleichung die Klassenzahl, sowie eine
Einheit. Es muss dabei bemerkt werden, dass die Tafel keinen Anspruch
auf Vollständigkeit erhebt, und dass die angegebenen Einheiten nicht
notwendig fundamental sind. Weiter habe ich für die Körper in der
untenstehenden Tabelle die euklidischen Minima bestimmt.

Der Tafel kann man übrigens entnehmen, dass der kubische Zahlkörper mit
$\disc K = -283$ semi-euklidisch ist und euklidische Tiefe $1$ hat.
Das euklidische Minimum des Körpers mit Diskriminante $-87$ konnte ich
bisher nicht bestimmen; wahrscheinlich (?) ist jedoch $M(K) = 1/3$,
wobei dieses Minimum außer an den beiden oben angegebenen rationalen
Punkten noch an unendlich vielen irrationalen Punkten angenommen wird
(wie bei $\Q(\sqrt{13})$).

Der Ring $R = \Z[\theta]$, wo $\theta$ eine Nullstelle des Polynoms
$f(x) = x^3 + x^2 + 3x -1$ ist, hat Diskriminante $-176$ und ist im
Körper mit der Diskriminante $-44$ enthalten. $R$ ist ein weiteres
Beispiel für einen nicht ganz abgeschlossenen Ring mit $M(f) = 1$,
wobei $M(f)$ mod $R$ im Punkt $(1+\theta^2)/2$ angenommen wird.

\begin{table}[h]
\centering
$$ \begin{array}{c|c|c}
\toprule
\disc K & M(K) & C_1 \\
\midrule
\rsp  -23 & \frac{1}{5} & (1+\theta+2\theta^2)/5,\ (2+2\theta-\theta^2)/5 \\
\rsp  -31 & \frac{1}{3} & (1-\theta-\theta^2)/3 \\
\rsp  -44 & \frac{1}{2} & (1+\theta^2)/2 \\
\rsp  -59 & \frac{1}{2} & (1+\theta+\theta^2)/2 \\
\rsp  -76 & \frac{1}{2} & (1+\theta^2)/2 \\
\rsp  -83 & \frac{1}{2} & (1+\theta+\theta^2)/2 \\
\rsp  -87 & \frac{1}{3} & (1-\theta^2)/3,\ (1+\theta+\theta^2)/3 \\
\rsp -104 & \frac{1}{2} & (\theta+\theta^2)/2 \\
\rsp -107 & \frac{1}{2} & (3+\theta-3\theta^2)/8,\ (1+\theta+\theta^2)/2 \\
\rsp -108 & \frac{1}{2} & \sqrt[3]{4}/2 \\
\rsp -116 & \frac{1}{2} & (\theta+\theta^2)/2,\ (1+\theta^2)/2 \\
\rsp -135 & \frac{3}{5} & (2+2\theta-2\theta^2)/5 \\
\rsp -139 & \frac{1}{2} & (1+\theta+\theta^2)/2 \\
\rsp -140 & \frac{1}{2} & (3+2\theta-3\theta^2)/10,\ (1+\theta^2)/2 \\
\rsp -152 & \frac{1}{2} & (\theta+\theta^2)/2,\ \theta^2/2 \\
\rsp -172 & \frac{3}{4} & (\theta+\theta^2)/2 \\
\rsp -175 & \frac{3}{5} & (2-\theta+2\theta^2)/5 \\
\rsp -199 & 1 & (3+\theta-3\theta^2)/7 \\
\rsp -283 & \frac{3}{2} & (1+\theta+\theta^2)/2 \\
\rsp -307 & \frac{9}{8} & (1+\theta^2)/2 \\
\rsp -528 & \frac{5}{2} & (1+\theta^2)/2 \\
\rsp -891 & \frac{7}{2} & (1+\theta+\theta^2)/2 \\
\bottomrule
\end{array} $$
\end{table}

Wir haben bereits darauf hingewiesen, dass sich die in (1.10)
errechneten Schranken $\mu_1$ im kubischen Fall noch verbessern
lassen; dies wollen wir nun tun. Dazu betrachten wir zuerst den rein
kubischen Fall $K = \Q(\theta)$, $\theta^3 = m$, und
wählen die $\Q$-Basis $\{1, \theta, \theta^2\}$. Dann
gilt mit den Bezeichnungen von (1.10) und $v:=|z|$, $w:=|z'|$
$N_{K/\Q}(|z|) = v w^2 \leq k$, außerdem $v \leq V$ und $w \leq W$
(wo wir $V$ und $W$ noch geeignet bestimmen werden; in (1.10)
war $V = W = \sqrt[3]{ku}$).

Um $f(v,w) = v+2w$ auf dem Gebiet $0 \leq v \leq V$, $0 \leq w \leq W$
nach oben abzuschätzen, bemerken wir, dass $f$ ihr Maximum höchstens
auf dem Rand annimmt. Ist aber $v = V$, so folgt $w \leq k/v = k/V$
und damit $v+2w \leq V + 2\sqrt{k/V}$. Im Falle $w = W$ dagegen folgt
entsprechend $v \leq k/W$ und $v+2w \leq 2W + k/W^2$.

Also haben wir insgesamt
$$ v+2w \leq \max \{V + 2\sqrt{k/V}, 2W + k/W^2\}. $$
Will man nun erreichen, dass $V = 2W$ wird, so muss man nur
$\sqrt[3]{4k/u^2} \leq |z| \leq \sqrt[3]{4ku}$ wählen und erhält dann
$|z'| \leq \sqrt[3]{ku/2}$, also $V = \sqrt[3]{4ku} = 2W$ und
$$ v + 2w \leq \sqrt[3]{4ku} + \max \{2\sqrt{k/V}, k/W^2\}
   = \sqrt[3]{4ku} + \sqrt[3]{4k/\sqrt{u}} =: S $$
wegen $2\sqrt{k/V} = \sqrt[3]{4k/\sqrt{u}}$ und $k/W = \sqrt{4k/u}$. Jetzt
folgen die Schranken $\mu_1 = S/3$, $\mu_2 = S/3\theta$,
$\mu_3 = S/3\theta^2$ (mit $\theta^3 = m$).

Die Asymmetrie im Beweis lässt mich vermuten, dass diese Schranken noch
nicht bestmöglich sind. Allerdings liefern diese Überlegungen schon
weit bessere Schranken als (1.10); so findet man z.B. für $m = 12, 17,
44$.

\begin{table}[h]
\centering
\begin{tabular}{c|ccc}
\toprule $m$ & $\mu_1$ & $\mu_2$ & $\mu_3$ \\ \midrule 12 & 2.92 &
1.27 & 0.59 \\ 17 & 5.25 & 2.05 & 0.80 \\ 44 & 8.41 & 2.38 & 0.68
\\ \bottomrule
\end{tabular}
\end{table}

Beachtet man, dass die Rechenzeit in etwa proportional zu dem Produkt
$\mu_1\mu_2\mu_3$ ist und vergleicht diese Schranken mit den etwas
weiter oben erhaltenen, so wird einem der Fortschritt gegenüber (1.10)
schnell klar.

Im allgemeinen Fall gehen wir einen andern Weg als in (1.10): wir
benutzen die $\Q$-Basis $\{1, \alpha, \beta\}$; anstatt aber
jetzt die Dualbasis zu benutzen, multiplizieren wir $z$ mit einem
Element der ``Spur'' 0 und bilden dann die ``Spur'' (``Spur'' steht
dabei in Anführungszeichen, weil z.B. $\alpha'-\alpha''$ gar nicht in
$K$ liegt und man folglich
$(\alpha'-\alpha'')+(\alpha''-\alpha)+(\alpha-\alpha'')$ nicht als
Spur von $(\alpha'-\alpha'')$ bezeichnen darf). So gilt z.B.
\[
\sum z(\alpha' - \alpha'') = r_3 m, \quad
\sum z(\beta' - \beta'') = r_2 m, \quad
\sum z(\alpha'\beta'' - \alpha''\beta') = r_1 m,
\]
wobei wir die Summanden als Elemente des normalen Abschlusses $L =
K(\sqrt{d})$ von $K$ auffassen und über die Automorphismen von
$L/\Q(\sqrt{d})$ summieren (d.h. z.B. $\sum z' = z' + z'' +
z$). Dabei ist
\[
m = \alpha'\beta + \alpha''\beta' + \alpha\beta'' - \alpha\beta' - \alpha'\beta'' - \alpha''\beta = \begin{vmatrix}
1 & 1 & 1 \\
\alpha & \alpha' & \alpha'' \\
\beta & \beta' & \beta''
\end{vmatrix} = \pm \sqrt{D},
\]
wo bekanntlich $D = \disc_{K/\Q}(1, \alpha, \beta)$ ist.
Mit Hilfe der Dreiecksungleichung erhält man nun
\begin{align*}
  \|r_3\| & = \Big|\sum z(\alpha' - \alpha'') \Big|
           \le |z||\alpha' - \alpha''| + |z'||\alpha' - \alpha'|
             + |z''||\alpha - \alpha'| \\
          & \le \sqrt[3]{ku}(|\alpha' - \alpha''| + 2|\alpha - \alpha'|),
             \intertext{ und entsprechend } \\
  \|r_2\| & \le \sqrt[3]{ku}(|\beta' - \beta''| + 2|\beta - \beta'|), \quad
             \text{sowie} \\
  \|r_1\| & \le \sqrt[3]{ku}(|\alpha'\beta'' - \alpha''\beta'| +
             2|\alpha\beta' - \alpha'\beta'|).
\end{align*}
Division durch $|m| = \sqrt{-D}$ liefert nun die gewünschten Schranken
für $r_1, r_2, r_3$. Wir halten fest:
\begin{quote}
  {\bf (4.8)} {\em Sei $K$ ein kubischer Körper mit Einheitenrang 1,
    $u$ eine Einheit mit $u>1$, und $x_1, \dots, x_t \in K$ seien
    Punkte, die von $u$ permutiert werden. Weiter sei $\{1, \alpha,
    \beta\}$ eine $\Q$-Basis von $K$. Ist dann $M(K, x_i) <
    k$, so gibt es ein $z = r_1 + r_2\alpha + r_3\beta \in K$ mit den
    Eigenschaften}
    \begin{enumerate}
    \item[(a)] $z \equiv x_i \bmod {R}$ für ein $j \in \{1, \dots, t\}$;
    \item[(b)] $N_{K/\Q}(z) < k$;
    \item[(c)] $|r_i| < \mu_i$ {\em für $i=1,2,3$ und}
      \begin{align*}
        \mu_1 & = \frac{\sqrt[3]{ku}(|\alpha'\beta'' - \alpha''\beta'|
          + 2|\alpha\beta' - \alpha'\beta'|)}{\sqrt{-D}}, \\
        \mu_2 & = \frac{\sqrt[3]{ku}(|\beta' - \beta''|
          + 2|\beta - \beta'|)}{\sqrt{-D}}, \\
        \mu_3 & = \frac{\sqrt[3]{ku}(|\alpha' - \alpha''|
          + 2|\alpha - \alpha'|)}{\sqrt{-D}},
      \end{align*}
      {\em wobei $D = \disc_{K/\Q}(1, \alpha, \beta)$ ist. }
    \end{enumerate}
\end{quote}
Mit Hilfe von (4.8) lassen sich einige der schon von Taylor
untersuchten kubischen Körper als nicht normeuklidisch nachweisen: man
lässt $y$ durch die Restklassen mod $(u-1)$ laufen und bestimmt
$M(K,x)$ für $x = y+1$. Es scheint erfolgversprechend zu sein, mit
dieser Methode diejenigen kubischen Körper zu untersuchen, von denen
man bisher nicht weiß, ob sie normeuklidisch sind oder nicht.
In seiner Arbeit aus dem Jahre 1954 hat
Swinnerton-Dyer\index[dN]{Swinnerton-Dyer} \cite{Swi54} unter anderem
folgendes Ergebnis (ohne Beweis) angegeben:
\begin{quote}
  {\em ist $\{1, \theta, \theta^2\}$ GHB von $\cO(\theta)$,
    $\theta^3 + 2a\theta -1 = 0$, und $a\in\N$ hinreichend
    groß, dann ist $M(K) = (a^2 - a +1)/2$, und dieses Minimum wird mod
    $R$ genau im Punkt $P = (\frac12, \frac12, \frac12)$ angenommen.}
\end{quote}
 
Wir werden nun zeigen, dass sogar gilt:
\begin{quote}
  {\bf (4.9)} {\em Sei $\theta^3 + 2a\theta -1 = 0$,
    $R = \Z[\theta]$ und $t$ der Absolutbetrag der
    Norm. Dann gilt für alle $a\in\N$: $M(K) = M(K) = (a^2 - a +1)/2$,
    und dieses Minimum wird mod $R$ nur in $(\frac12, \frac12, \frac12)$
    angenommen.}
\end{quote}
In diesem Fall ist übrigens $\disc_{K/\Q}(1,
\theta, \theta^2) = -(32a^3 + 27)$ und folglich $M(K) \approx
|d|^{2/3}/16^{3/2}$. Für $a = 1, 2, 3$ erhält man aus (4.9)
$\disc(\theta) = -59, -283, -891$; in diesen Fällen
haben wir das erste Minimum bereits per Computer bestimmt und
$M(K) = \frac12$, $\frac32$, $\frac72$ gefunden in Übereinstimmung
mit (4.9). Beim Beweis dürfen wir uns daher auf $a \geq 4$ beschränken.
Die Begleitmatrix von $z = r + s\theta + t\theta^2$ bezüglich
der Basis $\{1, \theta, \theta^2\}$ bestimmt sich zu
\[
\begin{vmatrix}
r & s & t \\
t & r - 2a t & s \\
s & t - 2a s & r - 2a t
\end{vmatrix}
\]
Hierbei stehen in der zweiten und dritten Zeile jeweils die
Koordinaten von $z\theta$ und $z\theta^2$. Bildung der
Determinante liefert die Norm
\[
N(z) = r^3 + s^3 + t^3 - 3rst + 2a r s^2 - 2a s^2 t - 4a t^2 + 4a^2 r t^2.
\]
Da wir $z$ mod $R$ so verschieben können, dass
$-\frac12 \leq r, s, t \leq \frac12$ wird, finden wir sofort
\[
|N(z)| \le \frac18 + \frac18 + \frac18 + \frac38
         + \frac a4 + \frac a4 + \frac a2 + 2a^2 t^2.
\]
Für jeden Ausnahmepunkt $z$ mit $|N(z)| \ge k = (a^2 - a +1)/2$ gilt
also $2a^2 t^2 \ge (2a^2 - 6a -1)/4$, d.h. $t^2 \ge (2a^2 - 6a -1)/(8a^2)$;
für alle $a \ge 3$ hat man die Abschätzung
$(2a^2 - 6a -1)/(8a^2) > (\frac12 - \frac4{5a})^2$, und wir sehen
$|t| \ge \frac12 - \frac4{5a}$.
Der Koeffizient von $\theta^2$ in $z\theta -1$, $z$ und
$z\theta$ ist bzw. $r$, $t$ und $s$; folglich gelten für $a \ge 3$
auch die Abschätzungen $|r| \le \frac12 - \frac58$ und
$|s| \le \frac12 - \frac58$. Indem wir wieder den Fundamentalsbereich
$F = (-0.5,0.5)\times(-0.5, 0.5)\times(-0.5, 0.5)$ ersetzen durch z.B.
$F' = (0, 1) \times(-1, 0) \times(0, 1)$ (diese Wahl von $F'$ wird
von der Überlegung geleitet, dass $|N(\frac12, \frac12, \frac12)| = k$,
aber $|N(\frac12, \frac12, \frac12)| = \frac{a^2}2 > k$ ist), bekommen
wir die einzige Ausnahmemenge
$ S = (\frac12 - \delta, \frac12 + \delta) \times
     (-\frac12 - \delta, \frac12 + \delta) \times
     (-\frac12 - \delta, \frac12 + \delta)$ mit $\delta = \frac4{5a}$.
Sei nun $z \in S$; wir wollen $|N(z)|$ noch einmal nach oben
abschätzen in der Hoffnung, wegen $z \in S$ bessere Schranken zu
erhalten. Indem wir $z$ notfalls mit $-1$ multiplizieren und eine
entsprechende Rechnung durchführen, dürfen wir $0 \le r \le \frac12$, $0
\le s \le \frac12$, $0 \le t \le \frac12$ annehmen und finden nach etwas
Rechnung $|N(z)| \le 3a\delta + \delta^2 + \frac12 - \frac{a}2 + 2a^2 t^2$.
Wie oben erhält man hieraus die Abschätzung
$t^2 \ge \frac14 - \frac6{5a^2} - \frac2{5a^3}$, und für alle
$a \ge 4$ liefert dies $|t - \frac12| \le \frac1{3a}$. Wie oben folgt
nun auch  $|r - \frac12| \le \frac1{3a}$ und $|s - \frac12| \le \frac1{3a}$.

Jetzt beachten wir $z\theta = t + (r-2at)\theta + s\theta^2$;
da mit $z$ auch $z\theta$ Ausnahmepunkt ist, muss
$|r-2at + b + \frac12| \le \frac1{3a}$ für ein $b \in \Z$ sein.
Etwas Rechnung zeigt $b = a-1$, und es folgt nach (2.1), dass
$x = - \frac{(a-1)\theta + \theta^2}{\theta-1} = \frac{1- \theta + \theta^2}2$
der einzig mögliche Ausnahmepunkt von $S$ ist.

Es bleibt noch zu zeigen, dass $M(K,x) = \frac{a^2-a+1}2$ gilt. Dies macht
man wie in der Arbeit von Swinnerton-Dyer; man muss lediglich dessen
asymptotische Abschätzungen durch explizite ersetzen, was nicht schwer
ist.

\subsection*{2-stufig normeuklidische Zahlkörper}

Schließlich geben wir noch einige Beispiele für 2-stufig
normeuklidische kubische Körper:
\begin{enumerate}
\item[1.] $\disc K = -199$: dieser Körper wird durch eine Nullstelle
  $\alpha$ des Polynoms $f(x) = x^3+4x^2+x+1$ erzeugt, $\{1, \alpha,
  \alpha^2\}$ ist GHB, wir haben die Faktorisierungen $(2) = (2)$,
  $(3) = \fth_1 \fth_2$, $(5) = (5)$ und $(7) = \fsv_1 \fsv_2 \fsv_3$.
  Hierbei ist $\fth_1 = (\alpha+1)$, $\fsv_1 = (\alpha-1)$,
  $\fsv_2 = (\alpha+2)$, $\fsv_3 = (\alpha+3)$. Da die Einheit $\alpha$
  fundamental ist, folgt leicht $M(K,\fsv_1) = 1$, da die Restklassen
  $2, 3 \bmod \fsv_1$ nur Elemente der Norm $\ge 7$ enthalten. Nach
  Taylor gilt hier sogar $M(K) = 1$; jedenfalls ist $K$ nicht
  normeuklidisch. Um zu zeigen, dass $K$  wenigstens 2-stufig
  normeuklidisch ist, müssen wir eine 2-Folge der Länge $\ge 4$ finden.
  Da die Ideale $(2)$, $7\fsv_2$ und $\fsv_3$ ein PERS
  besitzen und $|N_{K/\Q}(\alpha+4)| = 3$ ist, ist
  $0$, $1$, $2$, $\alpha+2$, $\alpha+3$, $\alpha+4$ eine 2-Folge der
  Länge $6$, und es folgt, dass $K$ 2-euklidisch ist.
\item[2.] $\disc K = -351$: $K$ wird von einer Nullstelle
$\alpha$ des Polynoms $f(x) = x^3+3x+3$ erzeugt; wir haben die
  Faktorisierungen $(2) = (2)$, $(3) = P^3$, $(5) = (5)$,
  $(7) = fsv_1fsv_2$, $(11) = \fel_1 \fel_2 11_3$, und mit $N = |N_{K/\Q}|$ gilt
$N(\alpha) = 3$, $N(\alpha+1) = 1$, $N(\alpha+2) = 11$, $N(\alpha-1) =
7$, $N(\alpha-2) = 17$, $N(\alpha^2-\alpha+1) = 19$,
$N(\alpha^2-\alpha) = 21$. Man kann dann nachprüfen, dass $0, 1, 2,
\alpha$, $\alpha+1$, $\alpha^2+1$ eine 2-Folge ist, und dies genügt um
zu zeigen, dass $K$ 2-euklidisch ist. Dass $K$ nicht normeuklidisch ist,
sieht man leicht, wenn man sich die Elemente in der Restklasse
$5 \bmod \fel_1$ betrachtet, wo $\fel_1 = (\alpha+2)$ ist.
\end{enumerate}

Wir wollen nun noch etwas zur Situation in totalreellen kubischen
Zahlkörpern sagen. Hier weiß man nur, dass die Anzahl der zyklischen
normeuklidischen Körper endlich ist (Heilbronn\index[dN]{Heilbronn}
\cite{Hei50}). Smith\index[dN]{Smith} \cite{Smi69} hat 1969 die Körper
mit Diskriminante $< 10^8$ mit einem Computer untersucht und dabei
festgestellt, dass der EA höchstens dann existiert, wenn
$$ \sqrt{d} = 7, 9, 13, 19, 31, 37, 43, 61, 67, 73, 103, 109, 127, 157$$
ist. Darüberhinaus hat er für die
Werte $\sqrt{d} \le 43$ und $\sqrt{d} = 73$ das euklidische Minimum
$M(K)$ bestimmt und gezeigt, dass $K$ für $\sqrt{d} \le 67$
normeuklidisch ist, für $\sqrt{d} = 73$ dagegen nicht. Die Körper mit
$\sqrt{d}= 103, 109, 127, 157$ sind nicht untersucht worden.

Das Kriterium für die Nichtexistenz des EA, das sowohl Heilbronn, als
auch Smith\index[dN]{Smith} benutzt haben, ist das folgende:
\begin{quote}
  {\bf (4.10)} {\em Sei $p \equiv 1 \bmod 6$ prim und $\chi$ ein
    nichttrivialer kubischer Charakter mod $p$. Gibt es dann
    $r, s, t, u \in \N$ mit $p = rs-tu$, $(r,s)=(t,u)=1$,
    $\chi(rs) = \chi(tu) = 1$, $\chi(r) \neq 1$, $\chi(s) \neq 1$,
    $\chi(t) \neq 1$, dann ist $K$ nicht normeuklidisch.}
\end{quote}
\begin{proof}
  Wegen $\chi(rs) = 1$ ist $a = rs$ kubischer Rest mod $p$, während
  die andern an $r, s, t, u$ gestellten Bedingungen garantieren, dass
  weder $rs$ noch $-tu$ Norm aus $R$ sind. (1.6) liefert dann die
  Behauptung.
\end{proof}
\begin{table}[h]
\centering
\begin{tabular}{c|ccccc}
\toprule
 & $p$ & $r$ & $s$ & $t$ & $u$ \\
\midrule
 & 79 & 2 & 5 & 3 & 23 \\
 & 97 & 2 & 11 & 3 & 25 \\
 & 139 & 2 & 3 & 7 & 19 \\
 & 151 & 2 & 43 & 5 & 13 \\
 & 163 & 2 & 11 & 3 & 47 \\
 & 181 & 2 & 11 & 3 & 53 \\
 & 193 & 2 & 37 & 7 & 17 \\
 & 199 & 2 & 41 & 9 & 13 \\
\bottomrule
\end{tabular}
\caption{Beispiele}
\end{table}
Zur Berechnung des kubischen Charakters $\chi$ hat Smith\index[dN]{Smith}
eine Primitivwurzel herangezogen. Benutzt man stattdessen das kubische
Reziprozitätsgesetz (sh. z.B. Ireland u. Rosen \cite{IR82}), lässt sich die
Rechenzeit wesentlich verkürzen, und man erhält z.B.

\begin{quote}
  {\bf (4.11)} {\em Sei $K$ zyklischer kubischer Zahlkörper mit
    Diskriminante $d$. Ist dann $197^2 \le d \le 2.5 \cdot 10^{11}$,
    so ist $K$ nicht normeuklidisch.}
\end{quote}

Da es unwahrscheinlich ist, dass für einen Körper mit
$\disc K > 2.5 \cdot 10^{11}$ eine Darstellung wie in
(4.10) nicht existiert, dürfte die oben gegebene Liste mit zyklischen
kubischen Körpern, die möglicherweise normeuklidisch sind, vollständig
sein.

1971 hat Smith\index[dN]{Smith} \cite{Smi71} euklidische Minima einiger
reeller kubischer Körper bestimmt und außerdem den EA in zahlreichen
derartigen Körpern nachgewiesen. In die Liste mit Diskriminanten
normeuklidischer Körper hat sich ein Druckfehler eingeschlichen: statt
disc $K = 1994$ muss disc $K = 1944$ stehen.

Die Ergebnisse von Smith lassen sich relativ problemlos ausdehnen: so
kann man z.B. zeigen, dass die kubischen Körper mit Diskriminante
$\disc K = 2021$, $2024$, $2037$, $2101$, $2213$ normeuklidisch sind.

An dieser Stelle wollen wir noch zeigen, wie man auch nicht rein
verzweigte Ideale zur Abschätzung von $M(K)$ verwenden kann. Es gilt
nämlich

\begin{quote}
  {\bf (4.12)} {\em Sei $K$ Zahlkörper vom Grad $n$, $R$ der Ring
    ganzer Zahlen in $K$, und es sei $pR =
    \fp^{n-1}\fq$ die Faktorisierung des Ideals
    $(p)$ in $R$. Ist dann $\pi \equiv a \bmod \fp$ und
    $\pi \equiv b \bmod \fq$ für $a, b \in \Z$, so gilt
    $N_{K/\Q}(\pi) = a^{n-1}b \bmod p$.}
\end{quote}

\begin{proof}
  Man geht vor wie in (1.1), (1.2) und (1.3); man beachte nur, dass das
  Minimalpolynom $f$ eines $\pi \in P$ die Form $f(x) = x^n +
  a_{n-1}x^{n-1} + \dots + a_0$ hat, wobei für $i = 0, \dots, n-2$ die
  Kongruenz $a_i \equiv 0 \bmod p$ gilt und $a_{n-1} \equiv \pi \bmod \fq$
  ist.
\end{proof}

Man kann hier auch ein Analogon zum Eisensteinschen
Irreduzibilitätskriterium beweisen: ist $f(x) = x^n + a_{n-1}x^{n-1} +
\dots + a_0 \in \Z[x]$, $a_i \equiv 0 \bmod p$ für $0 \le i \le n-2$,
$a_0 \not\equiv 0 \bmod p^2$, und hat $f$ keinen Linearfaktor $x-a$
(mit $a \equiv a_{n-1} \bmod p$), dann ist $f$ irreduzibel.

Benutzt man lokale Methoden (sh. Ishida\index[dN]{Ishida} \cite{Ish76},
so kann man (4.12) ohne Schwierigkeiten verallgemeinern zu

\begin{quote}
  {\bf (4.13)} {\em Sei $pR = \fp_1^{s_1} \dots \fp_g^{s_g}$
    die Faktorisierung von $(p)$ in einem algebraischen Zahlkörper. Dann
    gilt für alle $\alpha \in R$:
    $N_{K/\Q}(\alpha) = a_1^{s_1} \dots a_g^{s_g}$, wobei
    $a_i \in \Z$ ganze rationale Zahlen sind.}
\end{quote}

Wir geben nun einige Beispiele: Sei $K$ kubischer Zahlkörper mit $5
\parallel$ disc $K$; dann ist $(5) = 5_1 5_2^2$. Ist jetzt $\alpha \in
R$, $\alpha \equiv \pm 2 \bmod 5_1$ und $\alpha \equiv b \bmod 5_2$, so
folgt aus (4.12), dass $N_{K/\Q}(\alpha) = \pm 2 b^2 \equiv 0,
\pm 2 \bmod 5$ ist. Wegen $\alpha \neq 0$ ist $5_1$ kein euklidisches
Ideal, wenn es in $R$ keine Elemente der Norm 2 oder 3 gibt. Man
erhält so
\begin{table}[h]
\centering
\begin{tabular}{c|c}
\toprule
$\disc K$ & $M(K, 5_1)$ \\
\midrule
985 & 1 \\
1345 & $7/5$ \\
3305 & $\ge 1$ \\
4345 & $\ge 1$ \\
6185 & $\ge 1$ \\
\bottomrule
\end{tabular}
\end{table}

\noindent
\textbf{Bem.}: Ich bin nicht mehr dazu gekommen, $M(K, 5_1)$ in allen
Fällen zu bestimmen oder zu untersuchen, ob sich diese Überlegung auch
mit Primidealen über andern $p \equiv 1 \bmod 4$ erfolgreich
durchführen lässt.

\section*{{\sc Anmerkungen zu} \S\ 4}
\addcontentsline{toc}{section}{{\sc Anmerkungen zu} \S\ 4}

\subsection*{1. Kubische Zahlkörper mit negativer Diskriminante}

\begin{tabular}{p{1cm} p{10cm}}
\toprule
\textbf{Jahr} & \textbf{Ereignis} \\
\midrule
1892 & Markov\index[dN]{Markov} gibt Einheiten von $\Q(\sqrt[3]{m})$,
  $m \le 70$ \\
1940 & Delone\index[dN]{Delone} und Faddeev\index[dN]{Faddeev}
  veröffentlichen eine Tafel kubischer Körper mit $-999 \le \disc K \le -23$ \\
1949 & Es ist $M(K) = 1/5$ für den Körper mit
  $\disc K = -23$ (Prasad)\index[dN]{Prasad} \\
1950 & Es gibt nur endlich viele normeuklidische kubische Körper
  mit negativer Diskriminante (Davenport)\index[dN]{Davenport} \\
1952 & Cassels\index[dN]{Cassels} verbessert die Davenport'sche Schranke: es gilt
  $|d| < 420^2$ \\
1954 & Swinnerton-Dyer\index[dN]{Swinnerton-Dyer} beweist für komplexe
  kubische Körper mit $d \le -1237$ die Ungleichung
  $M(K) \le |d|^{2/3}/16^{3/2}$ und zeigt, dass diese bestmöglich ist. \\
1957 & Godwin\index[dN]{Godwin} gibt eine Tafel komplexer kubischer Körper \\
1967 & Die Körper mit $-23 \le \disc K \le -152$
  sind normeuklidisch (Godwin) \\
1973 & Angell\index[dN]{Angell} tabelliert die kubischen Körper mit $-23 \le
  \disc K \le -20000$ und gibt jeweils FE und Klassenzahl \\
1976 & Taylor\index[dN]{Taylor} bestimmt alle normeuklidischen Körper mit
  $|\disc K| < 680$ \\
1979 & Cioffari\index[dN]{Cioffari} findet alle rein kubischen
  Zahlkörper mit EA \\
1988 & Nakamula\index[dN]{Nakamula} veröffentlicht eine Tafel rein
  kubischer Körper und gibt FE samt Klassenzahl \\
\bottomrule
\end{tabular}

\subsection*{2. Kubische Zahlkörper mit positiver Diskriminante}

\begin{tabular}{p{1cm} p{10cm}}
\toprule
\textbf{Jahr} & \textbf{Ereignis} \\
\midrule
1923 & Remak\index[dN]{Remak} zeigt $M(K) \le \sqrt{d}/8$;
  dies impliziert insbesondere, dass der kubische Körper mit $\disc K = 49$
  normeuklidisch ist \\
1947 & Für die Körper mit $\disc K = 49, 81$ wird $M(K)$ bestimmt
  (Davenport)\index[dN]{Davenport} \\
1950 & Es gibt nur endlich viele zyklische kubische Körper mit
  EA; eine obere Schranke wird nicht gegeben (Heilbronn)\index[dN]{Heilbronn} \\
1951 & Für den Körper mit $\disc K = 148$ wird
  $M(K)$ bestimmt (Clarke) \index[dN]{Clarke}\\
1954 & Samet\index[dN]{Samet} bestimmt $M(K)$ für eine Reihe kubischer Körper \\
1956 & Billevic\index[dN]{Billevic} veröffentlicht eine Tafel kubischer
  Körper mit $\disc K \le 1296$ und gibt zwei Fundamentaleinheiten \\
1959 & Godwin\index[dN]{Godwin} und Samet\index[dN]{Samet} tabellieren
  reelle kubische Körper mit $\disc K \le 20.000$ \\
1969 & Smith\index[dN]{Smith} untersucht die zyklischen Körper mit
  $\disc K \le 10^8$ \\
1971 & Smith bestimmt die euklidischen Minima zahlreicher
  kubischer Körper und findet alle normeuklidischen Körper mit $\disc
  K \le 1957$ \\
1975 & Gras\index[dN]{Gras} gibt eine Tafel zyklischer Körper samt Einheiten und
  Klassenzahl \\
1985 & Ennola\index[dN]{Ennola} und Turunen\index[dN]{Turunen} bestimmen
  die reellen kubischen Körper mit $\disc K < 5 \cdot 10^5$ \\
1987 & Cusick\index[dN]{Cusick} und Schoenfeld\index[dN]{Schoenfeld}
  tabellieren die reellen kubischen Körper mit $\disc K \le 6885$ und
  geben Klassenzahl und Grundeinheiten an \\
1988 & Llorente\index[dN]{Llorente} und Quer\index[dN]{Quer} berechnen alle
  kubischen Körper mit $0 < \disc K < 10^7$. \\
\bottomrule
\end{tabular}

\chapter*{\S\ 5 Reine Zahlkörper von Zweierpotenzgrad}
\setcounter{chapter}{5}
\addcontentsline{toc}{chapter}{\S\ 5 Reine Zahlkörper von Zweierpotenzgrad}
\markboth{Euklidische Ringe}{\S\ 5 Reine Zahlkörper von Zweierpotenzgrad}

Wir beginnen mit rein biquadratischen Zahlkörpern und verwenden dazu
folgende Notation:
\[
m = ab^2 c^2 \text{ mit } (a,b) = (b,c) = (c,a) = 1, \quad a,b,c
\text{ quadratfrei,}
\]
\[
K_1 = \Q(\sqrt[4]{m}), \quad K_2 = \Q((1+i)\sqrt{m}) =
\Q(\sqrt{-4m}), \quad K_3 = \Q(i,\sqrt{ac}),
\]
\[
k_1 = \Q(\sqrt{ac}), \quad k_2 = \Q(\sqrt{-ac}), \quad
k_3 = \Q(i), \quad L = K_1(i) = K_2(i) = \Q(i,
\sqrt[4]{m}).
\]

Zuerst wollen wir uns fragen, wann das Polynom $f(x) = x^4 - m$ über
$\Q$ irreduzibel ist; da nach Voraussetzung $m$ keine
4. Potenz ist, hat $f$ keine Nullstellen in $\Q$, sodass
$f$ höchstens in das Produkt zweier quadratischer Faktoren zerfallen
kann. Setzt man
\[
f(x) = (x^2 + rx + s)(x^2 + tx + u) \text{ für } r,s,t,u \in
\Z, \text{ so folgt}
\]
\begin{enumerate}
\item[(a)] $r=t=0, s=-u$, also $m=s^2$: $m$ ist Quadrat in $\Z$; oder
\item[(b)] $s=u, r=-t, r^2=2s, r^4=4s^2=-4m$: $-4m$ ist 4. Potenz in $\Z$.
\end{enumerate}

Also ist $f$ genau dann zerlegbar über $\Q$, wenn $ac=1$
oder $ac=-1$, $b \equiv 0 \mod 2$ gilt. Diese beiden Fälle wollen
wir in Zukunft ausschließen.

Damit entsteht $K_1$ durch Adjunktion einer Wurzel $\alpha$ des
Polynoms $f$; die anderen Wurzeln von $f$ sind $-\alpha$ und
$\pm i\alpha$. Falls $K_1/\Q$ galoissch ist, muss $i \in
K_1$ sein. Dann zerfällt $f$ über $\Q(i)$ in das Produkt
zweier quadratischer Polynome mit Koeffizienten aus $\Z[i]$,
und ein Koeffizientenvergleich liefert $m=-b^2$ für ein $b \in
\Z$. Tatsächlich ist dann $f(x) = (x^2 + ib)(x^2 - ib)$,
sowie $K_1=\Q(i,\sqrt{2b})$. Ist also $\pm m$ kein Quadrat
in $\Z$, dann ist $f$ irreduzibel, $L$ der
Zerfällungskörper von $f$ und $[L:K_j] = 2$ für $j=1,2,3$.

Als Galois-Gruppe von $f$ bzw. $L/\Q$ entpuppt sich die
Diedergruppe
\[
D_4 = \{1, s, s^2, s^3, t, ts, ts^2, ts^3\}
\]
mit den Automorphismen $s: i \mapsto i$, $\alpha \mapsto i\alpha$
und $t: i \mapsto -i$, $\alpha \mapsto \alpha$.

$D_4$ hat die folgenden Untergruppen:
\begin{align*}
  G_1 & = \{1, t\}, \quad G_1' = s^{-1}G_1 s = \{1, ts^2\}, \\
  G_2 & = \{1, ts\}, \quad G_2' = s^{-1}G_2 s = \{1, ts^3\}, \\
  G_3 & = \{1, s^2\}, \quad G_3' = s^{-1}G_3 s = \{1, ts^2\}, \\
  G_4 & = \{1, ts\}, \quad G_4' = s^{-1}G_4 s = \{1, ts^3\},
\intertext{sowie das Zentrum $Z(D_4)$:}
  g_1 & = \{1, s^2, t, ts^2\}, \quad g_2 = \{1, s^3, ts, ts^3\},
  \quad g_3 = \{1, s^2, s^3\}.
\end{align*}

Dabei ist $K_j$ Fixkörper von $G_j$, $k_j$ Fixkörper von
$g_j$, jeweils für $j=1,2,3$. Weiter ist
\begin{align*}
  K_1' & = \Q(i\sqrt[4]{m}) \text{ der Fixkörper von $G_1'$ und} \\
  K_2' & = \Q(i\sqrt[4]{-4m}) = \Q((1-i)\sqrt[4]{m})
             \text{ derjenige von $G_2'$.}
\end{align*}

Wir haben also das folgende Körperdiagramm:

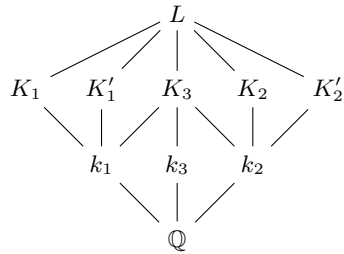
\begin{center}
  \begin{tikzpicture}
    \node (Q) at (0,0) {$\Q$};
    \node (k1) at (-1, 1) {$k_1$};
    \node (k3) at ( 0, 1) {$k_3$};
    \node (k2) at ( 1, 1) {$k_2$};
    \node (K1) at (-2, 2) {$K_1$};
    \node (K1') at (-1, 2) {$K_1'$};
    \node (K3) at ( 0, 2) {$K_3$};
    \node (K2) at ( 1, 2) {$K_2$};
    \node (K2') at ( 2, 2) {$K_2'$};
    \node (L) at ( 0, 3) {$L$};
    \draw (Q) -- (k1) -- (K1) -- (L) -- (K2) -- (k2) -- (Q);
    \draw (Q) -- (k3) -- (K3) -- (L);
    \draw (k1) -- (K1') -- (L);
    \draw (k2) -- (K2') -- (L);
    \draw (k1) -- (K3) -- (k2);
  \end{tikzpicture}
\end{center}

\section*{Hilbertsche Untergruppenreihe}

$$ \footnotesize \begin{array}{ll|ccc|cccccccccccc} \hline
  p & & e & f & g & K_1 & K_2 & Z & T &
       V_1 & V_2 & V_3 & V_4 & V_5 & V_6 & V_7 & V_8 \\ \hline
  1 \ (4) & (\frac mp) = +1, (\frac mp)_4 = +1 &
    1 & 1 & 8 & (1,1,1,1) & (1,1,1,1)
       & 1 & 1 & 1 & 1 & 1 & 1 & 1 & 1 & 1 & 1 \\
                   & (\frac mp) = +1, (\frac mp)_4 = -1 &
    1 & 2 & 4 & (2,2) & (2,2)
       & G_3 & 1 & 1 & 1 & 1 & 1 & 1 & 1 & 1 & 1 \\
                   & (\frac mp) = -1 &
    1 & 4 & 2 & (4) & (4)
       & g_3 & 1 & 1 & 1 & 1 & 1 & 1 & 1 & 1 & 1 \\
    & (\frac{ac}{p}) = +1, p \mod b & 2 & 1 & 4 & (1^2,1^2) & (1^2,1^2)
       & g_3 & G_3 & & 1 & 1 & 1 & 1 & 1 & 1 & 1 \\
    & (\frac{ac}{p}) = -1, p \mod b & 2 & 2 & 2 & (2^2) & (2^2)
       & g_3 & G_3 & & 1 & 1 & 1 & 1 & 1 & 1 & 1 \\
    & p \mid ac & 4 & 1 & 2 & (1^4) & (1^4)
       & g_3 & g_3 & & 1 & 1 & 1 & 1 & 1 & 1 & 1 \\ \hline
    3 \ (4) & (\frac mp) = +1 & 1 & 2 & 4 & (1,1,2) & (2,2)
         & G_1 & 1 & 1 & 1 & 1 & 1 & 1 & 1 & 1 & 1 \\
       & (\frac mp) = +1 & 1 & 2 & 4 & (2,2) & (1,1,2)
         & G_2 & 1 & 1 & 1 & 1 & 1 & 1 & 1 & 1 & 1 \\
       & (\frac{ac}{p}) = +1, p \mod b & 2 & 2 & 2 & (1^2,1^2) & (2^2)
         & g_1 & G_3 & 1 & 1 & 1 & 1 & 1 & 1 & 1 & 1 \\
       & (\frac{ac}{p}) = -1, p \mod b & 2 & 2 & 2 & (1^2) & (1^2,1^2)
         & g_2 & G_3 & 1 & 1 & 1 & 1 & 1 & 1 & 1 & 1 \\
      & p \mid ac & 4 & 1 & 2 & (1^4) & (1^4)
         & G & g_3 & 1 & 1 & 1 & 1 & 1 & 1 & 1 & 1 \\ \hline
    p=2 & m \equiv 1 \bmod 16 & 2 & 1 & 4 & (1,1,1^2) & (1^2,1^2)
         & G_1 & G_1 & G_1 & 1 & 1 & 1 & 1 & 1 & 1 & 1\\
       & m \equiv 9 \mod 16 & 2 & 2 & 2 & (2,1^2) & (1^2,1^2)
         & g_1 & G_1 & G_1 & 1 & 1 & 1 & 1 & 1 & 1 & 1 \\
       & \frac m4 \equiv -1 \mod 16, 2| b & 2 & 1 & 4 & (1^2,1^2) & (1,1,1^2)
         & G_2 & G_2 & G_2 & 1 & 1 & 1 & 1 & 1 & 1 & 1 \\
       & \frac m4 \equiv -9 \mod 16, 2| b & 2 & 2 & 2 & (1^2,1^2) & (2,1^2)
         & g_2 & G_2 & G_2 & 1 & 1 & 1 & 1 & 1 & 1 & 1\\
       & ac \equiv 1 \mod 8, 2| b & 4 & 1 & 2 & (1^2,1^2) & (1^4)
         & g_1 & g_1 & g_1 & G_3 & G_3 & & 1 & 1 & 1 & 1\\
       & ac \equiv 5 \mod 8, 2 \nmid b & 4 & 2 & 1 & (1^2) & (1^4)
         & G & g_1 & g_1 & 1 & 1 & 1 & 1 & 1 & 1 & 1\\
       & ac \equiv 5 \mod 8, 2 \mid b & 4 & 2 & 1 & (1^2) & (1^4)
         & G & g_1 & g_1 & G_3 & G_3 & 1 & 1 & 1 & 1 & 1\\
      & ac \equiv 7 \mod 8, 2\nmid b & 2 & 1 & (1^2) & (1^4)
          & g_2 & g_2 & g_2 & G_3 & G_3 & 1 & 1 & 1 & 1 & 1 & 1\\
      & ac \equiv 3 \mod 8, 2 \nmid b & 4 & 2 & 1 & (1^4) & (2^2)
        & G & g_2 & g_2 & G_3 & G_3 & 1 & 1 & 1 & 1 & 1\\
      & ac \equiv 3 \mod 8, 2 \mid b & 4 & 2 & 1 & (1^4) & (2^2)
        & G & g_2 & g_2 & 1 & 1 & 1 & 1 & 1 & 1 & 1\\
      & 2 \mid ac & 8 & 1 & 1 & (1^4) & (1^4)
        & G & G & G & g_3 & g_3 & G_3 & G_3 & G_3 & G_3 & 1
  \end{array} $$

\section*{Ganzheitsbasis von $\Q(\sqrt[4]{m})$}

\begin{table}[h]
\centering
\begin{tabular}{llll}
\toprule
\textbf{m} & \textbf{disc $K_1$} & \textbf{disc $L$} & \textbf{GHB}  \\
\midrule
$2 \nmid ac$ & $2^8 a^3 b^2 c^3$ & $2^{18} a^6 b^4 c^6$ & $\{1, \alpha_1, \alpha_2, \alpha_3\}$   \\
$ac \equiv 3 \mod 4$, $2 \nmid b$ & $2^8 a^3 b^2 c^3$ & $2^{18} a^6 b^4 c^6$ & $\{1, \alpha_1, \alpha_2, \alpha_3\}$  \\
$ac \equiv 3 \mod 8$, $2 \nmid b$ & $2^4 a^3 b^2 c^3$ & $2^8 a^6 b^4 c^6$ & $\{1, \alpha_1, \beta_2, \beta_3\}$  \\
$ac \equiv 7 \mod 8$, $2 \nmid b$ & $2^2 a^3 b^2 c^3$ & $2^4 a^6 b^4 c^6$ & $\{1, \alpha_1, \beta_2, \beta_3\}$ \\
$ac \equiv 1 \mod 8$, $2 \nmid b$ & $2^2 a^3 b^2 c^3$ & $2^8 a^6 b^4 c^6$ & $\{1, \alpha_1, \beta_2, \beta_3\}$ \\
$ac \equiv 5 \mod 8$, $2 \nmid b$ & $2^4 a^3 b^2 c^3$ & $2^{12} a^6 b^4 c^6$ & $\{1, \alpha_1, \beta_2, \beta_3\}$ \\
$ac \equiv 1 \mod 4$, $2 \nmid b$ & $2^4 a^3 b^2 c^3$ & $2^{12} a^6 b^4 c^6$ & $\{1, \alpha_1, \beta_2, \beta_3\}$  \\
\bottomrule
\end{tabular}

$$ \begin{array}{lll}
\toprule
 m & \beta_2 & \beta_3 \\
ac \equiv 3 \mod 8, 2 \nmid b &
  (1 + \alpha_1 + \alpha_2)/2 &  (\alpha_1 + \alpha_3)/2 \\
ac \equiv 7 \mod 8, 2 \nmid b & 
  (1 + \alpha_1 + \alpha_2)/2 & (\alpha_1 + 2\alpha_2 + \alpha_3)/2 \\
ac \equiv 1 \mod 8, 2 \nmid b &
  (1 + \alpha_2)/2 & (1 + b\alpha_1 + ab\alpha_2 + \alpha_3)/4 \\
ac \equiv 5 \mod 8, 2 \nmid b &
  (1 + \alpha_2)/2 & (\alpha_1 + \alpha_3)/2 \\
ac \equiv 1 \mod 4, 2 \nmid b &
  (1 + \alpha_2)/2 & (\alpha_1 + \alpha_3)/2 
\end{array} $$
\caption{Ganzheitsbasis von $\Q(\sqrt[4]{m})$}
\end{table}

\medskip

Die Galoisgruppe $D_4$ ist 1987 von C.~E. van der Ploeg untersucht
worden; dessen Arbeit enthält jedoch einige Fehler: so ist es
z. B. nicht richtig, dass jeder nichtnormale Zahlkörper 4. Grades einen
quadratischen Zahlkörper enthält, oder dass als Galoisgruppen von
Polynomen 4. Grades nur $Z_4$, $V_4$ oder $D_4$ in Frage
kommen (es fehlen $A_4$ und $S_4$; Zahlkörper 4. Grades, deren
Galoisabschluss diese Gruppen als Galoisgruppen besitzt, haben
z. B. keine quadratischen Zahlkörper); außerdem fehlt in den dort
aufgeführten Körperdiagrammen der Körper, den wir mit $k_3$ bezeichnet
haben.

Um nun eine GHB, die Diskriminante und das Zerlegungsgesetz von
$\Q(\sqrt[4]{m})$ zu finden, gehen wir wie folgt vor: wir
stellen die Hilbertsche Untergruppenreihe für alle Primideale $P$ in
$L$ auf, lesen daraus das Zerlegungsgesetz in $L$, $K_1$, $K_2$ ab,
bestimmen die Beiträge der einzelnen $P$ zur Differente, und errechnen
schließlich die Diskriminante und die Ganzheitsbasis von $K_1$. Die
Ergebnisse dieser Rechnungen sind in den oben stehenden Tafeln
festgehalten; dabei bedeuten
\begin{itemize}
\item $p \in \N$ eine rationale Primzahl
\item $P$ ein Primideal in $L$ über $(p)$
\item $Z = Z(P|p)$ die Zerlegungsgruppe
\item $T = T(P|p)$ die Trägheitsgruppe
\item $V_j = V_j(P|p)$ die $j$-te Verzweigungsgruppe, $V_0 = T$
\item $e = e(P|p)$ der Verzweigungsindex
\item $f = f(P|p)$ der Trägheitsgrad
\item $\alpha_1 = \alpha = \sqrt[4]{m}$, \quad $m = ab^2c^2$, \quad $\alpha_2 = \sqrt{ac}$, \quad $\alpha_3 = \sqrt[4]{a^3b^2c^3}$
\item $( \cdot / p )$ das Legendre-Symbol
\item $( \cdot / p )_4$ das biquadratische Restsymbol (falls $( \cdot / p ) = +1$ ist).
\end{itemize}

Weiter bedeutet z. B. $(1,1,2)$, dass $(p)$ in dem entsprechenden Körper
in das Produkt dreier Primideale zerfällt, von denen eines den
Trägheitsgrad 2, die beiden anderen Trägheitsgrad 1 besitzen,
d. h. $(p) = P_1 P_2 P_3$, $\|P_1\| = \|P_2\| = p$, $\|P_3\| =
p^2$. Entsprechend soll dann $(1,1,1^2)$ die Zerlegung $(p) = P_1 P_2
P_3^2$ symbolisieren usw.

Eine GHB für $\Q(\sqrt[4]{m})$ ist erstmals von
Ljunggren\index[dN]{Ljunggren} \cite{Lju36} angegeben worden,
allerdings ohne Beweis. Auf einem anderen Weg als dem hier
beschriebenen hat Funakura\index[dN]{Funakura} \cite{Fun84} eine GHB
für $\Q(\sqrt[4]{m})$ errechnet. Die Parität der Klassenzahlen von
$K_1$, $K_2$ und $L$ hat Parry\index[dN]{Parry} \cite{Par75a}
bestimmt; die in seinem Beweis verwendeten Zerlegungen des Ideals
$(2)$ in den Zwischenkörpern von $L/\Q$ sind in den beiden Fällen
$m \equiv 1 \mod 8$ und $\sqrt[4]{m} \equiv -1 \mod 8$ nicht
richtig, weil diese -- wie ein Vergleich mit unserer Tabelle zeigt --
von den Restklassen mod $16$ abhängen. Man kann sich aber leicht davon
überzeugen, dass dieser Fehler Parrys Resultate nicht beeinflusst;
diese lauten:

\begin{enumerate}
\item Sei $K = \Q(\sqrt[4]{m})$, $m \in \N$; dann ist
  $h(K)$ genau dann ungerade, wenn $m$ die folgende Gestalt hat:
  \begin{align*}
  m & = 2,\ 2p^2\ (p \equiv 3 \bmod {8}), \\
  m & = p \ (p \equiv \pm 3 \bmod {8}), \\
  m & = 4p\ (p \equiv \pm 3, 7 \bmod {8}), \\
  m & = 2p\ (p \equiv 3 \bmod {8}),\ \\
  m & = 8p\ (p \equiv 3 \bmod {8}).
\end{align*}
\item Sei $K = \Q(\sqrt[4]{-m})$, $m \in \N$; dann ist
$h(K)$ genau dann ungerade, wenn $m$ die folgende Gestalt hat:
$$ m = 2,\ p (p \equiv 3 \bmod {4}),   \
   m = 4p    (p \equiv 3 \bmod {8} $$
\item Sei $L = \Q(i, \sqrt[4]{m})$; dann ist $h(L)$ genau dann
ungerade, wenn $m$ die folgende Gestalt hat:
$$ m = 2,\ p (p \equiv 3 \bmod {8}\,),\
   m = 4p    (p \equiv 3 \bmod {4}. $$
\end{enumerate}
Hierbei bedeutet $p$ immer eine rationale, positive Primzahl.

Berücksichtigt man diese Ergebnisse von Parry, so folgt aus Theorem B
von Cioffari\index[dN]{Cioffari} \cite{Cio79}: ist $m$ kein Quadrat
und $K = \Q(\sqrt[4]{-m})$ normeuklidisch, so ist $m \in \{2, 3, 7,
12, 44, 67\}$.  Wir werden hier zeigen, dass genau die
$\Q(\sqrt[4]{-m})$ mit $m=2, 3, 7, 12$ normeuklidisch sind; dabei
werden wir einen etwas anderen Weg einschlagen als Cioffari; dies wird
es uns erlauben, auf das tiefliegende Ergebnis von Stark (nämlich die
Bestimmung aller imaginärquadratischen Zahlkörper mit Klassenzahl 1)
zu verzichten.

Der zweite Teil von Cioffaris Arbeit enthält noch einmal zwei kleinere
Druckfehler:
\begin{enumerate}
\item In Prop. 14 muss das erste $K$ durch $k$ ersetzt werden;
  \item In der Bemerkung zu Prop. 16 schreibt Cioffari: ``Cassels
    proved that $K$ cannot be Euclidean if $D_{K/\Q} > 5300^2$;
    hence $\Q(\sqrt[4]{-163})$ and
    $\Q(\sqrt[4]{652})$ are not Euclidean.''
\end{enumerate}

Dabei sollte es selbstverständlich $\Q(\sqrt[4]{-652})$
heißen. Zweitens hat van der Linden 1983 bemerkt, dass sich Cassels in
seiner Arbeit verrechnet hat und die richtige Schranke
$D_{K/\Q} > 15170^2$ lautet. Dies führt dazu, dass sich
$\Q(\sqrt[4]{-163})$ mit Hilfe der Schranke von Cassels nicht
mehr ausschließen lässt.

Im Folgenden werden wir Körper der Form $K = \Q(\sqrt[k]{-m})$,
$k = 2^l$, $l \ge 1$ betrachten; hier haben
wir

\begin{quote}
  \textbf{(5.1)} {\em Seien $m \in \Z$ und $\pm m$ kein Quadrat;
    mit $K = \Q(\sqrt[k]{m})$ und $L = \Q(\sqrt[2k]{m})$ für ein
    $k \in \N$ gilt dann $h(K) \mid h(L)$. }
\end{quote}

Bevor wir (5.1) beweisen, erinnern wir an einige Tatsachen aus der
Klassenkörpertheorie: ist $K$ ein algebraischer Zahlkörper mit
Klassenzahl $h$ (im weiteren Sinne), dann existiert ein über $K$
abelscher Körper $L$ mit den Eigenschaften
\begin{itemize}
\item[] (CF--1) $L/K$ ist an allen Primstellen (einschließlich der
  unendlichen) unverzweigt;  
\item[] (CF--2) Die Galoisgruppe $\mathrm{Gal}(L/K)$ ist isomorph zur
  Idealklassengruppe $\mathrm{Cl}(K)$; insbesondere ist $(L:K) = h =
  |\mathrm{Cl}(K)|$.
\end{itemize}

$L$ ist durch $K$ eindeutig bestimmt, heißt der Hilbertklassenkörper
von $K$ und wird mit $\mathrm{CF}(K)$ bezeichnet. Ist umgekehrt $L/K$
abelsch und überall unverzweigt, so gilt $(L:K) \mid h(K)$.

\begin{proof}[Bew. von 5.1]
  Wegen der Inklusion $K \subset L \cap \mathrm{CF}(K) \subset L$
  und $(L:K)=2$ ist entweder $K = L \cap \mathrm{CF}(K)$ oder $L
  \subset \mathrm{CF}(K)$; in letzterem Falle wäre aber $L/K$
  unverzweigt, und dies ist nicht der Fall: denn da $m$ kein Quadrat
  ist, gibt es eine rationale Primzahl $p$, die in $m$ ungerade
  oft aufgeht. Diese ist dann in $L/\Q$ und insbesondere auch in
  $L/K$ rein verzweigt.  Also ist $L \cap \mathrm{CF}(K) = K$,
  somit $L\cdot\mathrm{CF}(K)/L$ eine abelsche, unverzweigte
  Erweiterung vom Grade $(L\cdot\mathrm{CF}(K):L) = (\mathrm{CF}(K):L
  \cap \mathrm{CF}(K)) = (\mathrm{CF}(K):K) = h(K)$, und dies zeigt
  $h(K) \mid h(L)$.
\end{proof}

Im Falle $k=2$ stammt diese Aussage wohl von
Chevalley\index[dN]{Chevalley} \cite{Che31}. Sei nun $m=qs^2$ für ein
quadratfreies $q>1$: mit $L = \Q(\sqrt[k]{-m})$, $k = 2^l$, ist
$K = \Q(\sqrt{-q})$ der quadratische Teilkörper von $L$.  Soll
$L$ euklidisch sein, muss $h(K)=1$ sein wegen (5.1), und
bekanntlich ist dies höchstens dann der Fall, wenn $q=2$ oder $q
\equiv 3 \bmod {4}$ prim ist. Ist aber $q \ge 19$ und $h(K)=1$,
so müssen die Ideale $(2)$ und $(3)$ beim Übergang $K/\Q$ prim
bleiben, weil es in $K$ keine ganzen Zahlen der Norm $2$ oder
$3$ geben kann. Nach dem Zerlegungsgesetz in quadratischen
Zahlkörpern ist somit $(-q/2) = (-q/3) = -1$, nach dem QRG also
$(2/q) = (3/q) = -1$. Weil $(q)$ in $L/\Q$ rein verzweigt,
können wir (1.6) anwenden mit $f=q$, $a=6$, $b=q-6$ und
behaupten, dass $a=6$ ein $2^l$-ter Potenzrest mod $q$
ist. Wegen $q \equiv 3 \bmod {4}$ und $(a/q) = (2/q)(3/q) = +1$
ist dies aber der Fall (es genügt ja, dass $a$ quadratischer Rest
mod $q$ ist). Jetzt müssen wir nur noch zeigen, dass $a$ und
$-b$ keine Normen aus $S$ sind. Wegen $K \subset L$ und der
Normschachtelungsformel brauchen wir dazu nur festzustellen, dass
$a$ und $-b$ keine Normen aus $\Q(\sqrt{-q})$ sind, und dies ist
klar.

\begin{quote}
  \textbf{(5.2)} {\em Sei $L = \Q(\sqrt[k]{-m})$ mit $k = 2^l$,
    $l > 1$ und $m = qs^2$ für ein quadratfreies $q \in \N$,
    $q > 1$. Dann ist $L$ höchstens dann normeuklidisch, wenn $q \in
    \{2, 3, 7, 11\}$ ist. Ist sogar $k \ge 4$, so muss im Falle $q =
    11$ auch $s \equiv \pm 2 \bmod {3}$ gelten.}
\end{quote}

Zu beweisen ist nur noch die Aussage über $q=11$ für $k \ge 4$. Dazu
setzen wir $i=11$, $a=5$, $b=6$ in (1.6); wegen $(\frac{5}{11}) = +1$
und $11 \equiv 3 \bmod {4}$ ist 5 ein $2^1$-ter Potenzrest $\bmod
11$. Ist dann $s \equiv +1 \bmod {5}$, dann gilt
\[
\left(\frac{-m}{5}\right)^4 = \left(-\frac{11}{5}\right)^4
\left(\frac{5}{11}\right)^4 = -1,
\]
und nach dem Zerlegungsgesetz in $K = \Q(\sqrt[4]{-m})$ gibt
es in $K$ (und damit erst recht in $L$) keine Ideale der Norm
$5$. Wegen $(\frac{19}{7}) = -1$ ist 14 nicht einmal
Norm aus $\Q(\sqrt{19})$, und mit (1.6) folgt, dass $L$ auch
für $p=19$ nicht normeuklidisch sein kann. Man sieht nun leicht ein,
dass auch $\Q(\sqrt[4]{11})$ nicht normeuklidisch ist: da es
in $D(11)$ und damit erst recht in $\cO_L$ keine Ideale der
Norm 3 gibt, genügt es zu zeigen, dass die Grundeinheiten von $L$
beide $\equiv 1 \bmod {(3+\sqrt{11})}$ sind, da das Ideal $I =
(3+\sqrt{11})$ wegen $\Phi(11) = 2$ dann nicht euklidisch ist. Nun
sind aber $u_1 = 10 + 3\sqrt{11}$ und
$u_2 = 881 + 477 \theta + 264 \theta^2 + 147 \theta^3$
unabhängig, keine Quadrate und beide $\equiv 1 \bmod I$, sodass die
Behauptung folgt.

\begin{quote}
  \textbf{(5.3)} {\em Sei $L = \Q(\sqrt[4]{-m})$, $m \in \N$
    kein Quadrat und frei von 4. Potenzen; dann ist $L$ genau für
    $m=2,3,7,12$ normeuklidisch.}
\end{quote}

Den Beweis, dass diese Körper normeuklidisch sind, erbringt man wie in
den §§ 2, 3, 4 mit einem Computer (sh. § 11); für $m=3$ hat schon
Lakein\index[dN]{Lakein} \cite{Lak72} den EA nachgewiesen, für $m = 2,
7$ hat dies Cioffari\index[dN]{Cioffari} \cite{Cio79} getan. Es bleibt
noch zu zeigen, dass $L = \Q(\sqrt[4]{-44})$ nicht normeuklidisch
ist. Dazu verwenden wir (1.4) mit $K = \Q(\sqrt{-11})$, $B = (2)$
und $e = (1-\sqrt{-11})/2$ (wegen $\Phi_K(2) = 3$ sind in $K$
alle Zahlen quadratische Reste $\bmod 2$). Um $M(L) \ge
\frac{5}{4}$ zu zeigen, müssen wir die Elemente $r \in D(-11)$
betrachten, die $r \equiv e \bmod {2}$ und $N_{K/\Q}(r) <
\frac{5}{4}$ erfüllen. Das einzige solche $r$ ist aber $r = e$,
und es genügt zu zeigen, dass $e$ keine Norm aus $L$
ist. Bezeichnet $[\,\frac{\cdot}{\cdot}\,]$ das quadratische
Restsymbol in $D(-11)$, so folgt dies sofort aus
\[
[\sqrt{-44}/e] = [2\sqrt{-11}/e] = [-1/e] = (-1/3) = -1
\]
und dem Zerlegungsgesetz in relativquadratischen Zahlkörpern.

Man kann dies aber auch aus dem Zerlegungsgesetz in $L$ folgern: es
gilt nämlich $(3) = 3_1 3_2 3_3$ mit $\|3_1\| = \|3_2\| = 3$,
$\|3_3\| = 9$; also bleibt genau eines
der beiden Ideale $P_1 = \frac{1+\sqrt{-11}}2$ und $P_2 = (e)$ beim
Übergang nach $L$ prim, und zwar ist dies $P_2$ wegen
$P_1 \subset L = 3_1 3_2$ für
$3_1 = (\frac{3 - \sqrt[4]{-44} + \sqrt{-11}}2)$ und
$3_2 = (\frac{3 + \sqrt[4]{-44} + \sqrt{-11}}2)$.
Mit (1.4) folgt nun $M(L) \ge \frac{5}{4}$, und (5.3) ist bewiesen.

Jetzt wenden wir uns den Körpern $\Q(\sqrt[k]{-m})$ mit $k=2^l$,
$l>1$, $m \in \N$ zu. Die einzigen bekannten Ergebnisse für
$l>2$, die den EA in diesen Körpern betreffen, stammen meines
Wissens von Egami\index[dN]{Egami} \cite{Ega79}, der unter Verwendung
der Ergebnisse von Parry und tiefliegender Abschätzungen von
Charaktersummen gezeigt hat:
\begin{quote}
  {\em die Anzahl der normeuklidischen Zahlkörper der Form
    $\Q(\sqrt[4]{m})$, wo $m$ frei von 4. Potenzen und $m \ne
    2p^2$ für prime $p \equiv 3 \bmod {8}$ ist, ist endlich.}
\end{quote}

Wir werden nun zeigen, dass wir uns unter Beschränkung auf rein
elementare Methoden erstens von der Einschränkung $m \ne 2p^2$ frei
machen können, und dass wir zweitens diese endliche Liste von Körpern
$\Q(\sqrt[4]{m})$, die möglicherweise normeuklidisch sind,
sogar explizit angeben können.

\begin{quote}
  \textbf{(5.4)} {\em Sei $L = \Q(\sqrt[k]{2^m p})$ für
    $k=2^l$, $l \ge 1$, $m \equiv 1 \bmod {2}$ und $p \equiv 3
    \bmod {4}$ prim. Dann ist $L$ höchstens für $p=3$ normeuklidisch.}
\end{quote}

\begin{proof}
  Für $l=1$ wissen wir dies bereits aus § 3. Sei also $l \ge 2$; nach
  Parry dürfen wir dann $p \equiv 3 \bmod {8}$, $p \ge 11$
  annehmen. Nach (3.2) gibt es dann $a,b \in \N$ mit $2p = a + b$,
  $a \equiv 5 \bmod {8}$ und $(\frac{a}{p}) = +1$. Wegen $p \equiv 3
  \bmod {4}$ und $a \equiv 1 \bmod {2}$ ist $a$ damit $2^l$-ter
  Potenzrest $\bmod {2p}$, sodass wir (1.6) mit $f=2p$ anwenden können
  (denn die Ideale $(2)$ und $(p)$ sind in $L/\Q$ rein verzweigt). Mit
  $K=\Q(\sqrt{2p})$ sind aber $a$ und $-b$ keine Normen aus
  $\mathcal{O}_K$, wegen $K \subset L$ also erst recht keine Normen
  aus $\mathcal{O}_L$. Also ist $L$ dann nicht normeuklidisch.
\end{proof}

Für $k=4$ bleiben also noch $L=\Q(\sqrt[4]{6})$ und
$L=\Q(\sqrt[4]{24})$ zu untersuchen; beide Körper sind jedoch
nicht normeuklidisch: dazu müssen wir nur zeigen, dass $M(L,I) = 9/8$
für das Ideal $I$ der Norm $8$ für beide Körper $L$ gilt. Ist aber
$\alpha = 1+\sqrt{6} \bmod I$ in $\cO_L$, so folgt schnell, dass
$N_{L/K}(\alpha) = 7 + 2\sqrt{6} \equiv 3 \bmod {4 + 2 \sqrt{6}}$ und
$N_{L/\Q}(\alpha) \equiv 1 \bmod {8}$ gilt, wobei $K = \Q(\sqrt{6})$
der quadratische Teilkörper von $L$ ist. Die Kongruenz in $\cO_K$
zeigt, dass $\alpha$ keine Einheit ist (denn die FE von $\cO_K$ ist
$= 5 + 2 \sqrt{6} \equiv 1\bmod {4 + 2 \sqrt{6}}$). Da $(7)$ in $K$
prim bleibt, ist also $|N_{L/\Q}(\alpha)| \ge 9$, und wegen
$|N_{L/\Q}(3 + \sqrt{6})| = 9$ folgt in der Tat $M(L,I) = 9/8$.

\begin{quote}
  \textbf{(5.5)} {\em Sei $L = \Q(\sqrt[4]{2^m p})$, $k=2^l$,
    $l \ge 2$, $m \equiv 0 \bmod {2}$, $p \equiv 3 \bmod {8}$ prim. Dann
    ist $L$ höchstens für $p \le 19$ normeuklidisch.}
\end{quote}

\begin{proof}
  Sei $p \ge 43$; nach (3.4) gibt es $a,b \in \N$ mit $p = a + b$,
  $a \equiv 2 \bmod {8}$ und $(a/p) = +1$. Wegen $p \equiv 3
  \bmod {4}$ ist $a$ damit $k$-ter Potenzrest $\bmod {p}$. Aus dem
  Beweis von (3.5) folgt, dass weder $a$ noch $-b$ Normen aus
  $\cO_K$ für $K = \Q(\sqrt{p})$ sind. Wegen $K \subset
  L$ ist daher $L$ nicht normeuklidisch.
  Im Falle $p=19$, $m=0 \bmod {4}$, setzen wir $a=5$, $b=14$; dann ist
  $5 \equiv 3^4 \bmod {19}$, und $(\frac{19}{5})_4 = -1$
  zeigt, dass es in $\Q(\sqrt[4]{p})$ keine Ideale der Norm 5
  gibt. Wegen $(\frac{19}{7}) = -1$ ist 14 nicht einmal
  Norm aus $\Q(\sqrt{19})$, und mit (1.6) folgt, dass $L$ auch
  für $p=19$ nicht normeuklidisch sein kann. Man sieht nun leicht ein,
  dass auch $\Q(\sqrt[4]{11})$ nicht normeuklidisch ist: da es
  in $D(11)$ und damit erst recht in $\cO_L$ keine Ideale der
  Norm 3 gibt, genügt es zu zeigen, dass die Grundeinheiten von $L$
  beide $\equiv 1 \bmod {(3+\sqrt{11})}$ sind, da das Ideal $I =
  (3+\sqrt{11})$ wegen $\Phi(11) = 2$ dann nicht euklidisch ist. Nun
  sind aber $u_1 = 10 + 3\sqrt{11}$ und
  $u_2 = 881 + 477 \theta + 264 \theta^2 + 147 \theta^3$
  unabhängig, keine Quadrate und beide $\equiv 1 \bmod I$, sodass die
  Behauptung folgt.
\end{proof}

Ob $\Q(\sqrt[4]{3})$ normeuklidisch ist oder nicht, lässt sich
wohl nicht so einfach entscheiden. Jedenfalls ist $M(L) \ge
\frac{11}{12}$, wie man durch Betrachten des Ideals $I$ der Norm 12
leicht feststellt.

Die Fälle $m=2 \mod 4$, $p=11$ und $p=19$, lassen sich für $k \ge 4$
leicht entscheiden: man benutze einfach (1.6) mit $f = 11$, $a=5$,
$b=6$ bzw. $f = 38$, $a=17$, $b=21$ und beachte $(\frac{44}{3})_4 =
-1$, $(\frac{11}{3})_4 = -1$ resp. $(\frac{76}{17})_4 = -1$,
$(\frac{12}{7})_4 = -1$.

\begin{quote}
  \textbf{(5.6)} {\em Sei $L = \Q(\sqrt[k]{p})$, $k=2^l$,
    $l \ge 2$, $p \equiv 5 \bmod {8}$ prim. Dann ist $L$ höchstens für
    die Werte $p = 5, 13, 37, 61$ normeuklidisch.}
\end{quote}

\begin{proof}
  Aus den Arbeiten von Behrbohm\index[dN]{Behrbohm}
  u. R\'edei\index[dN]{Redei@R\'edei} \cite{BR36}; s. auch § 3) und
  Brauer\index[dN]{Brauer} \cite{Bra40} folgt, dass für prime $p \equiv 5
  \bmod {8}$, $p > 109$, eine Darstellung $p = rs+tu$ existiert mit
  $r,s,t,u \in \N$, $(r,s) = (t,u) = 1$ und $(\frac{r}{p}) =
  (\frac{s}{p}) = (\frac{t}{p}) = -1$ (sh. 3.12.). Da $-1$
  biquadratischer Nichtrest mod $p$ ist, ist entweder $rs$ oder $tu$
  biquadratischer Rest mod $p$ (beachte $rs = -tu \mod p$). Sei also
  oBdA $rs$ biquadratischer Rest mod $p$. Da weder $a=rs$, noch
  $-b=-tu$ Normen aus $\Q(\sqrt[4]{p})$ sind, ist $L$ nach (1.6) nicht
  normeuklidisch.
  \begin{enumerate}
    \item[] $p=29$: hier ist $24 \equiv 4^4 \bmod {p}$, $29 = 24 +
      5$. Da (3) in $\Q(\sqrt[4]{29})$ prim bleibt, ist $24$ sicher
      keine Norm aus $\mathcal{O}_L$. Weiter ist $5$ keine Idealnorm
      wegen $(\frac{29}{5})_4 = -1$;
  \item[] $p=53$: $53 = 15 + 38$, $8^4 = 15 \bmod {53}$;
  \item[] $p=109$: $109 = 104 + 5$, $(\frac{104}{109})_4 =
    (\frac{5}{109})_4 = +1$.
  \end{enumerate}
\end{proof}

Im Falle $l=2$ können wir auch $p=37$ ausschließen: dazu beachten wir,
dass das Ideal $I = (\frac{5 - \sqrt{37}}2)$ der Norm $3$ in $L$ prim
bleibt wegen $[\frac{\sqrt{37}}{I}] = [\frac 5I] = (\frac{5}{3}) =
-1$.  Dann setzen wir $B = (2)$, $e = (\frac{5 - \sqrt{37}}2)$, $k =
\frac74$ in (1.4) und betrachten alle $r \in \mathcal{O}_K$,
$K=\Q(\sqrt[4]{37})$, mit $r \equiv e \bmod {2}$ und $|N_{K/\Q}(r)| <
7$. Da $u = 6+\sqrt[4]{37}$ die FE von $K$ ist und $u \equiv 1 \bmod
{2}$ gilt, ist sicher $|N_{K/\Q}(r)| \in \{3, 5\}$, und weil (5) in
$K$ prim ist, muss sogar $|N_{K/\Q}(r)| = 3$ sein. Also ist $r =
u^{e'}$ für $e' = (5+\sqrt[4]{37})/2$. Letzteres widerspricht der
Kongruenz $r \equiv e \bmod {2}$, während die Möglichkeit $(r) = (e)$
daran scheitert, dass $I$ in $L$ prim bleibt.

\begin{quote}
  \textbf{(5.7)} {\em Sei $L = \Q(\sqrt[k]{2^mp})$,
    $k=2^l$, $m \equiv 0 \bmod {2}$, $p \equiv 5 \bmod {8}$ prim. Dann
    ist $L$ höchstens für $p = 5, 13, 29, 37, 61, 109$
    normeuklidisch.}
\end{quote}

\begin{proof}
  Wie (5.6).
\end{proof}

Im Falle $k=4$ lassen sich die Körper $\Q(\sqrt[4]{4 \cdot 37})$
und $\Q(\sqrt[4]{4 \cdot 61})$ leicht ausschließen: im
ersten Fall benutzt man den Zwischenkörper $K = \Q(\sqrt{37})$
und schließt genau wie oben im Falle
$\Q(\sqrt[4]{37})$. Verwendung von (1.6) mit $f=61$, $a=56$
und $b=5$ erledigt den zweiten Fall.

Schließlich müssen wir noch den Fall behandeln, wo keine rationale
Primzahl außer $p=2$ in $L/\Q$ rein verzweigt, in dem sich
also das Kriterium (1.6) nicht anwenden lässt. Überraschend einfach
sieht man jedoch

\begin{quote}
  \textbf{(5.8)} {\em Sei $L = \Q(\sqrt[4]{2p^2})$,
    $p \equiv 3 \bmod {8}$ prim. Dann ist $L$ nicht normeuklidisch.}
\end{quote}

\begin{proof}
  $L$ hat zwei Fundamentaleinheiten, nämlich $u=1+\sqrt{2}$ und $v$;
  wir zeigen, dass $v \equiv 1 \bmod {\sqrt{2}}$ ist, und wegen
  $\Phi_L(1) = 2$ und $\|1\| = 4$ ist das Ideal $I = (\sqrt{2})$ nicht
  euklidisch (da es keine Elemente der Norm 3 gibt).
  Da das Ideal $(p)$ in $L/K$ verzweigt, gibt es ein $\pi \in \cO_L$ mit
  $p = \pm \pi^2 u^m v^n$. Wäre $n \equiv 0 \bmod {2}$, so folgte
  $\sqrt{\pm p} \in L$ oder $\sqrt{\pm u} \in L$,
  was aber nicht der Fall ist. Nun ist aber $p \equiv 1 \bmod {2}$ und
  $u \equiv 1 \bmod {\sqrt{2}}$, sodass sich $v^n \equiv 1
  \bmod {\sqrt{2}}$ ergibt. Wegen $\Phi_L(1) = 2$ ist $v^2 \equiv 1
  \bmod {I}$, und jetzt folgt sofort $v \equiv 1 \bmod {I}$.
\end{proof}

Zusammenfassend können wir nun feststellen, dass $L =
\Q(\sqrt[4]{m})$ höchstens für die Werte $m = 2, 3, 5, 12, 13,
20, 28, 37, 52, 61, 116, 436$ normeuklidisch ist. Es ist anzunehmen,
dass sich hierunter noch einige nicht normeuklidische Körper
befinden. Weiter kann man mit Hilfe eines Computers bestätigen, dass
die beiden Körper $\Q(\sqrt[4]{2})$ und
$\Q(\sqrt[4]{5})$ normeuklidisch sind.

\section*{{\sc Anmerkungen zu} \S\ 5}
\addcontentsline{toc}{section}{{\sc Anmerkungen zu} \S\ 5}
\medskip

\begin{tabular}{p{1cm} p{10cm}}
\toprule
\textbf{Jahr} & \textbf{Ereignis} \\
\midrule
1936 & Ljunggren gibt Ganzheitsbasen und untersucht Einheiten
  rein biquadratischer Zahlkörper \\
1975 & Parry gibt alle rein biquadratischen Zahlkörper mit ungerader
  Klassenzahl \\
1979 & Egami zeigt, dass es nur endlich viele normeuklidische
  Körper der Form $\Q(\sqrt[4]{m})$ gibt, falls nicht $m =
  2p^2$, $p \equiv 3 \bmod {8}$ prim ist. Der Fall $m = 2p^2$ bleibt
  unerledigt. \\
     & Cioffari zeigt $m \in \{2, 3, 7, 12, 44, 67\}$ oder $m = 2p^2$ für
  normeuklidische Körper $\Q(\sqrt[4]{m})$. Der Fall $m =
  2p^2$ lässt sich aber mit den Ergebnissen von Parry (die Cioffari
  offenbar nicht bekannt waren) erledigen. \\
1980 & Parry entwickelt eine Geschlechtertheorie für rein
  biquadratische Körper \\
1984 & Funakura gibt eine GHB für rein biquadratische Körper \\
1987 & Buchmann berechnet eine Grundeinheit von Ringen der Form
  $\Z[\sqrt[4]{-m}]$ \\
\bottomrule
\end{tabular}

\chapter*{\S\ 6 Bizyklische biquadratische Zahlkörper}
\setcounter{chapter}{6}
\addcontentsline{toc}{chapter}{\S\ 6 Bizyklische biquadratische Zahlkörper}
\markboth{Euklidische Ringe}{\S\ 6 Bizyklische biquadratische Zahlkörper}

Seien $m$ und $n$ quadratfreie ganze Zahlen und $K = \Q(\sqrt{m},\sqrt{n})$.
Dann ist $K/\Q$ galoissche
Erweiterung von $\Q$ vom Grad 4 mit abelscher Galoisgruppe $G
= \mathrm{Gal}(K/\Q) = V_4$ (Kleinsche Vierergruppe). $K$
enthält genau drei nichttriviale Zwischenkörper, nämlich $k_1 =
\Q(\sqrt{m})$, $k_2 = \Q(\sqrt{n})$, und $k_3 =
\Q(\sqrt{mn})$.

Ist $l = (m,n)$, so schreiben wir $m = lm'$, $n = ln'$ und haben
$k_3 = \Q(\sqrt{m'n'})$ mit quadratfreiem $m'n'$. Jedes $\alpha \in K$
können wir in der Form $\alpha = r + s\sqrt{m} + t\sqrt{n} +
u\sqrt{mn}$, $r,s,t,u \in \Q$ darstellen; bezeichnet man die
Relativnorm von $K$ nach $k_i$ mit $N_i$, so findet man leicht
\[
N_1(\alpha) = a + b\sqrt{m}; \quad N_2(\alpha) = c + d\sqrt{n}; \quad
N_3(\alpha) = e + f\sqrt{mn} \quad \text{mit}
\]
\[
a = r^2 + m s^2 - n t^2 - m n u^2, \quad b = 2(rs - n t u),
\]
\[
c = r^2 - m s^2 + n t^2 - m n u^2, \quad d = 2(r t - m s u),
\]
\[
e = r^2 - m s^2 - n t^2 + m n u^2, \quad f = 2(r u - s t).
\]

Die Normschachtelungsformel liefert dann sofort
\[
N_{K/\Q}(\alpha) = a^2 - m b^2 = c^2 - n d^2 = e^2 - m n f^2.
\]

Bezeichnet man mit $T_i$ die Relativspur von $K$ nach $k_i$, so ist
ein $\alpha \in K$ genau dann ganz, wenn $N_i(\alpha)$ und
$T_i(\alpha)$ für ein $i \in \{1,2,3\}$ ganz sind. Nach einer etwas
langwierigen, aber elementaren Rechnung (sh. z. B. Williams \cite{Wil70})
ergibt sich dann: $\alpha = (r + s\sqrt{m} + t\sqrt{n} + u\sqrt{mn})/g
\in K$ ist genau dann ganz, wenn die $r, s, t, u, g \in \Z$
den Bedingungen (x) genügen:

\begin{center}
\begin{tabular}{|c|c|c|c|}
\hline
$m \bmod {4}$ & $n \bmod {4}$ & $d$ & (x) \\
\hline
1 & 1 & 4 & $r \equiv s \equiv t \equiv u \bmod {2}$,
$r+s+t+u \equiv 0 \bmod {4}$ \\
1 & 2 & 2 & $r \equiv s \bmod {2}$, $t \equiv u \bmod {2}$ \\
1 & 3 & 2 & $r \equiv u \bmod {2}$, $s \equiv t \bmod {2}$ \\
2 & 3 & 2 & $r \equiv t \equiv 0 \bmod {2}$, $s \equiv u \bmod {2}$ \\
\hline
\end{tabular}
\end{center}

Dabei sind $m$ und $n$ nur mod 4 angegeben; man macht sich schnell
klar, dass wegen $\Q(\sqrt{m},\sqrt{n}) =
\Q(\sqrt{m},\sqrt{mn}) = \Q(\sqrt{n},\sqrt{m})$
usw. alle möglichen Fälle aufgeführt sind.

Setzt man $d = \mathrm{disc}\, K$, $d_i = \mathrm{disc}\, k_i$
($i=1,2,3$), so rechnet man leicht nach, dass in jedem Fall die
Beziehung $d = d_1 d_2 d_3$ gilt (ein Beweis, der ohne Kenntnis der
Ganzheitsbasis auskommt, verläuft via der
\glqq{}Führer-Diskriminanten-Formel\grqq{}). Das Zerlegungsgesetz in
$K$ lautet dann
\begin{align*}
(p) & = \fp_1 \fp_2 \fp_3 \fp_4,
\quad \text{falls } \big(\frac{d_1}{p}\big) =
\big(\frac{d_2}{p}\big) = \big(\frac{d_3}{p}\big) = +1 \text{ ist}; \\
(p) & = \fp_1 \fp_2, \quad \text{falls }
\left(\frac{d_1}{p}\right) = \left(\frac{d_2}{p}\right) = +1,
\left(\frac{d_3}{p}\right) = -1 \text{ ist}; \\
(p) & = \fp_1^2 \fp_2^2, \quad \text{falls } p \nmid d_1,
   p \nmid d_2 \text{ und } \left(\frac{d_3}{p}\right) = +1 \text{ ist}; \\
   (p) & = \fp^2, \quad \text{falls } p \nmid d_1, p \nmid d_2 \text{ und }
   \left(\frac{d_3}{p}\right) = -1 \text{ ist}, \\
(p) & = \fp^4, \quad \text{falls } p=2, \, 2\mid d_1, \, 2\mid
d_2, \, 2\mid d_3 \text{ ist}.
\end{align*}

Hierbei bezeichnet $(\cdot / p)$ das Kroneckersymbol, das für ungerade
$p$ mit dem Legendresymbol übereinstimmt und für $p=2$ und $d \equiv 1
\bmod {4}$ durch $(d/2) = (2/d)$ definiert ist.

Die bei Cohn\index[dN]{Cohn} \cite[S. 202]{Coh78} gegebene Tafel der höheren
Verzweigungsgruppen für $K$ enthält zwei Fehler; um diese zu
korrigieren, übernehmen wir die dortige Notation und unterscheiden:

\begin{enumerate}
\item[(I)] $(d_1/p) = (d_2/p) = (d_3/p) = +1$
\item[(II)] $(d_1/p) = (d_2/p) = +1$, $(d_3/p) = -1$
\item[(III)] $p \equiv 1 \bmod {2}$, $p \nmid d_1$, $p \nmid d_2$, $(d_3/p) = +1$
\item[(IV)] $p \equiv 1 \bmod {2}$, $p \nmid d_1$, $p \nmid d_2$, $(d_3/p) = -1$
\item[(V)] $p=2$, $d_1 \equiv d_2 \equiv 12 \bmod {16}$, $d_3 \equiv 1 \bmod {8}$
\item[(VI)] $p=2$, $d_1 \equiv d_2 \equiv 12 \bmod {16}$, $d_3 \equiv 5 \bmod {8}$
\item[(VII)] $p=2$, $d_1 \equiv d_2 \equiv 8 \bmod {16}$, $d_3 \equiv 1 \bmod {8}$
\item[(VIII)] $p=2$, $d_1 \equiv d_2 \equiv 8 \bmod {16}$, $d_3 \equiv 5 \bmod {8}$
\item[(IX)] $p=2$, $d_1 \equiv d_2 \equiv 8 \bmod {16}$, $d_3 \equiv 12 \bmod {8}$
\end{enumerate}

Cohn gibt nun für die Fälle V. und VII. (bzw. VI. und VIII.) dieselben
Untergruppenreihen; dies würde jedoch implizieren, dass $d =
\mathrm{disc}\, K$ in den Fällen V. und VII. (bzw. VI. und VIII.)
durch dieselbe Zweierpotenz teilbar wäre, was aber nicht der Fall ist
(beachte $d = d_1 d_2 d_3$).

Die richtigen Untergruppenreihen gibt folgende Tafel:

\begin{center}
$$ \begin{array}{c|cccccccc}
\text{Typ} & \Q & K_{\mathfrak{z}} & K_{\mathrm{T}} & K_1 & K_2 & K_3 & K_4 & K \\
\midrule
\text{I}   & \Q & K & K & K & K & K & K & K \\
\text{II}  & \Q & k_3 & K & K & K & K & K & K \\
\text{III} & \Q & k_3 & k_3 & K & K & K & K & K \\
\text{IV}  & \Q & \Q & k_3 & K & K & K & K & K \\
\text{V}   & \Q & k_3 & k_3 & k_3 & K & K & K & K \\
\text{VI}  & \Q & \Q & k_3 & k_3 & K & K & K & K \\
\text{VII} & \Q & k_3 & k_3 & k_3 & k_3 & K & K & K \\
\text{VIII}& \Q & \Q & k_3 & k_3 & k_3 & k_3 & K & K \\
\text{IX}  & \Q & \Q & \Q & \Q & k_3 & k_3 & K & K \\
\end{array} $$
\end{center}

Seien ab jetzt $m$ und $n$ negativ; dann ist $K$ totalimaginär, und
seine Einheitengruppe hat Rang 1. Ist $u$ die FE von $k_3$ (das ist
der größte reelle Teilkörper von $K$), dann gibt es zwei
Möglichkeiten:

\begin{enumerate}
\item[1.] $u$ ist auch FE von $K$; wir setzen dann $q(K)=1$;
\item[2.] $u$ ist nicht FE von $K$; dann gibt es eine Einheitswurzel
  $\zeta_K$, sodass $u = \zeta \cdot e^2$ wird. In diesem Fall ist $e$
  FE von $K$, und wir schreiben $q(K)=2$.
\end{enumerate}

Die Konstante $q(K)$ heißt der \textbf{Einheitenindex von $K$} und
lässt sich verhältnismäßig leicht bestimmen:

\begin{enumerate}
\item[(a)] $K$ enthält $\Q(i)$: dann ist $q(K)=2$ genau dann, wenn das
  Ideal $(2)$ in $k_3$ Quadrat eines Hauptideals ist: schreibt man
  $(2) = (\alpha)^2$, so ist $e = \alpha^{-1}$ FE von $K$.
\item[(b)] $K$ enthält $\Q(i)$ nicht: dann ist $q(K)=2$ genau dann,
  wenn das Ideal $(m_1)$ in $k_3$ Quadrat eines Hauptideals ist: mit
  $(m_1) = (\alpha)^2$ wird $u = \alpha / \sqrt{m_1}$ FE von $K$.
\end{enumerate}

\noindent
\textbf{Bem.:} Teil (a) steht auch in Cohn\index[dN]{Cohn} \cite{Coh78}
als Teil von Theorem 19.8; jedoch fehlt dort die Bedingung, dass $(2)$
zum Quadrat eines Hauptideals wird, obwohl dies im Beweis verwendet
wird.

Aus der analytischen Klassenzahlformel folgt nun recht leicht die
elegante Beziehung $H = q(K) h_1 h_2 h_3$, wo $H, h_1, h_2, h_3$
jeweils die Klassenzahl von $K, k_1, k_2, k_3$ bedeutet und $K \ne
\Q(\sqrt{-1},\sqrt{2})$ ist. Hieraus wiederum erhält man $h_3 \mid H$
für die Klassenzahl $h_3$ des reellquadratischen Teilkörpers (dies ist
ein Spezialfall der allgemeineren Relation $h' \cdot h$, die für
CM-Körper $K$ mit größtem reellem Teilkörper $K^*$ gilt; sh. dazu
Washington\index[dN]{Washington} 1982). Ist $K$ euklidisch, so muss
insbesondere $H=1$ und damit auch $h_3 =1$ sein.

Beweise für diese und andere Ergebnisse finden sich bei
Kuroda\index[dN]{Kuroda} \cite{Kur43,Kur50},
Kubota\index[dN]{Kubota} \cite{Kub56},
Wada\index[dN]{Wada} \cite{Wad66}, 
Fröhlich\index[dN]{Frohlich@Fröhlich} \cite{Fro83}, sowie in dem
klassischen Werk von Hasse\index[dN]{Hasse} \cite{Has85}.

In ihrer Arbeit \cite{Sau73} aus dem Jahre 1972/73 schreibt
J.~Sauvageot:\index[dN]{Sauvageot}

\begin{quote}
  {\em On sait quels corps quadratiques sont euclidiens. La question
    reste ouverte pour les corps biquadratiques. Elle devrait être
    bientôt (?) résolue pour ceux d'entre eux qui sont bicyclique
    imaginaires, \dots''.}
\end{quote}

In der Tat gilt (den Ring ganzer Zahlen in
$\Q(\sqrt{m},\sqrt{n})$ wollen wir zukünftig mit $D(m,n)$
bezeichnen):

\begin{quote}
  \textbf{(6.1)} {\em Sei $m$ negativ; dann sind genau die folgenden
    Ringe $D(m,n)$ normeuklidisch:
    \begin{itemize}
    \item[] $D(-1,n)$ für $n = 2, 3, 5, 7$;
    \item[] $D(-2,n)$ für $n = -3, 5$;
    \item[] $D(-3,n)$ für $n = 2, 5, -7, -11, 17, -19$, und
    \item[] $D(-7,5)$.
      \end{itemize} }
\end{quote}

Für Körper $\Q(\sqrt{m}, \sqrt{n})$ mit $m \equiv n \equiv 1
\bmod {2}$, $m < 0$, hat bereits J. Sauvageot eine Klassifikation
versucht. Trotz eines Fehlers (im Falle $m = n \equiv 1 \bmod {4}$)
haben sich ihre Vermutungen, die sie am Ende ihrer Arbeit äußert,
bestätigt; sie schreibt dort:

\begin{quote}
  {\em \dots les seuls survivants de ces éliminatoires sont
    \begin{itemize}
      \item $\Q(i\sqrt{3}, i\sqrt{7}); \Q(i\sqrt{3}, i\sqrt{11});
        \Q(i\sqrt{3}, \sqrt{5}) $ dont LAKEIN a montré qu'ils le
        sont.
      \item[] $\Q(i\sqrt{3}, i\sqrt{19})$; $\Q(i\sqrt{3}, \sqrt{17})$
        dont j'espère montrer qu'ils le sont.
      \item[] $\Q(i\sqrt{3}, i\sqrt{43})$, que je crois non-euclidien
        et $\Q(i\sqrt{7}, \sqrt{5})$ dont je ne sais rien.
    \end{itemize}
    Les mots \glqq{}j'espère\grqq{} und \glqq{}je crois\grqq{} sont
    conséquence d'explorations du problème sur ordinateurs que j'exposerai
    \dots si elles aboutissent.}
\end{quote}

Bevor wir (6.1) beweisen, geben wir eine Tafel für die ersten Minima
gewisser $D(m,n)$ zusammen mit einer Menge $C_1$ von Punkten, an denen
diese angenommen werden.

\begin{table}[ht!]
$$ \begin{array}{cr|c|l}
\rsp m & n & M_1(L) & C_1 \\
\midrule
\multirow{10}{*}{$-1$}
\rsp & 2 & \frac{1}{2} & (\frac{1}{2}, \frac{1}{2}, \frac{1}{2}, 0) \\
\rsp & 3 & \frac{1}{4} & (\frac{1}{4}, \frac{1}{4}, \frac{1}{4}, \frac{1}{4}) \\
\rsp & 5 & \frac{5}{16} & (\frac{1}{4}, \frac{1}{4}, \frac{1}{4}, \frac{1}{4}), (\frac{1}{2}, \frac{1}{2}, \frac{1}{2}, 0), (\frac{1}{2}, \frac{1}{2}, \frac{1}{2}, \frac{1}{2}) \\
\rsp & 6 & \frac{3}{2} & (\frac{1}{2}, \frac{1}{2}, \frac{1}{2}, 0) \\
\rsp & 7 & \frac{1}{2} & (\frac{1}{4}, \frac{1}{4}, \frac{1}{4}, \frac{1}{4}) \\
\rsp & 10 & \frac{5}{2} & (\frac{1}{2}, \frac{1}{2}, \frac{1}{2}, 0) \\
\rsp & 11 & \frac{5}{4} & (\frac{1}{4}, \frac{1}{4}, \frac{1}{4}, \frac{1}{4}) \\
\rsp & 13 & \ge 1 & (0, 0, \frac{1}{13}, \frac{1}{13}) \\
\rsp & 14 & \frac{9}{2} & (\frac{1}{2}, \frac{1}{2}, \frac{1}{2}, 0) \\
\rsp & 15 & 1 & (\frac{1}{4}, \frac{1}{4}, \frac{1}{4}, \frac{1}{4}) \\
\midrule
\multirow{4}{*}{$-2$}
\rsp & -3 & \frac{1}{3} & (\frac{1}{2}, \frac{1}{2}, \frac{1}{2}, \frac{1}{2}) \\
\rsp & 5 & \frac{11}{16}& (\frac{1}{4}, \frac{1}{4}, \frac{1}{4}, \frac{1}{4}), (\frac{1}{4}, \frac{1}{4}, \frac{1}{4}, 0), (\frac{1}{2}, \frac{1}{2}, \frac{1}{2}, \frac{1}{2}) \\
\rsp & -7 & \frac{9}{8} & (\frac{1}{4}, \frac{1}{4}, \frac{1}{4}, \frac{1}{4}), M_2 < 0.999 \\
\rsp & -11& \ge \frac{6323}{5808} & (0, 0, \frac{1}{13}, \frac{1}{13}) \\
\midrule
\multirow{7}{*}{$-3$}
\rsp & 2 & \ge \frac{1}{4} & (0, 0, \frac{1}{2}, \frac{1}{2}) \\
\rsp & 5 & \frac{1}{4} & (0, \frac{1}{4}, \frac{1}{4}, 0), (\frac{3}{8}, \frac{1}{8}, \frac{1}{8}, \frac{1}{8}), (\frac{1}{2}, \frac{1}{2}, \frac{1}{2}, 0) \\
\rsp & -7 & \frac{4}{9} & (\frac{1}{4}, \frac{1}{12}, \frac{1}{4}, \frac{1}{12}) \\
\rsp & -11& <0.46 & \\
\rsp & 13 & 1 & (\frac{1}{4}, \frac{1}{12}, \frac{1}{4}, \frac{1}{12}) \\
\rsp & 17 & \frac{13}{16} & (0, \frac{1}{4}, \frac{1}{4}, 0), (\frac{3}{8}, \frac{1}{8}, \frac{1}{8}, \frac{1}{8}), (\frac{1}{2}, \frac{1}{4}, \frac{1}{4}, 0) \\
\rsp & -19& <0.95 & \\
\midrule
\rsp -7 & 5 & \frac{9}{16} & (0, \frac{1}{4}, \frac{1}{4}, 0), (\frac{3}{8}, \frac{1}{8}, \frac{1}{8}, \frac{1}{8}), (\frac{1}{2}, \frac{1}{4}, \frac{1}{4}, 0) \\
\bottomrule
\end{array} $$
  \caption{Euklidische Minima für Ringe $D(m,n)$}
\end{table}

Hierbei sind die Punkte unter $C_1$ bezüglich der Basis $\{1,
\sqrt{m}, \sqrt{n}, \sqrt{mn}\}$ angegeben. Weiter erhält man alle
Punkte, an denen $M(L)$ angenommen wird, durch Multiplikation mit $-1$
und Anwenden der drei nichttrivialen Automorphismen von
$L/\Q$. Man beachte auch, dass die Ringe $D(-1,15)$ und
$D(-3,13)$ semi-euklidisch sind; möglicherweise sind auch $D(-1,13)$
und $D(-1,17)$ semi-euklidisch.

Um (6.1) zu beweisen, gehen wir wie folgt vor: Wir betrachten eine
Restklasse $\alpha \bmod {2}$ in $D(m,n)$. Ist $\sigma$ derjenige
Automorphismus von $\Q(\sqrt{m},\sqrt{n})$, der genau
$\Q(\sqrt{m})$ elementweise fest lässt, so ist
$N_1(\beta) = \beta \beta^\sigma$. Aus $\beta \equiv \alpha \bmod {2}$
folgt nun aber $\beta^\sigma \equiv \alpha^\sigma \bmod {2}$ (wegen
$(2)^\sigma = (2)$), also $N_1(\alpha) \equiv N_1(\beta) \bmod {2}$
für alle $\alpha \equiv \beta \bmod {2}$. Entsprechendes gilt natürlich
auch für die beiden anderen Relativnormen $N_2$ und $N_3$.

Diese Beobachtung wird es uns ermöglichen, in den beiden
imaginärquadratischen Teilkörpern von $\Q(\sqrt{m},\sqrt{n})$
die Existenz ganzer Zahlen mit kleiner Norm nachzuweisen. Hat man
z. B. gezeigt, dass $D(m)$, $m$ negativ, ein Element der Norm 2 enthält,
so folgt leicht $m \in \{-1, -2, -7\}$ (man braucht ja nur die
Lösbarkeit der diophantischen Gleichung $x^2 - m y^2 = 2$ für
$m \equiv 2,3 \bmod {4}$, bzw. von $x^2 - m y^2 = 8$,
$x \equiv y \bmod {2}$, für $m \equiv 1 \bmod {4}$ zu untersuchen).

Den Beweis von (6.1) unterteilen wir in zwei Fälle:

\subsubsection*{I. Es existiert ein Ideal der Norm 2 in $D(m,n)$.}

Da $D(m,n)$ euklidisch ist, ist dieses Ideal Hauptideal, und
Relativnormbildung zeigt, dass beide imaginärquadratischen Teilkörper
Elemente der Norm 2 enthalten. Also sind zwei der drei Körper
$\Q(\sqrt{-1})$, $\Q(\sqrt{-2})$,
$\Q(\sqrt{-7})$ in $\Q(\sqrt{m},\sqrt{n})$ enthalten,
und es verbleiben die Möglichkeiten $D(-1,2)$, $D(-1,7)$ und
$D(-2,-7)$. Die obige Tabelle zeigt, dass hiervon genau die Ringe
$D(-1,2)$ und $D(-1,7)$ normeuklidisch sind.

\subsubsection*{II. Es existiert kein Ideal der Norm 2 in $D(m,n)$.}

Dann kommen nur folgende Möglichkeiten in Betracht:
\begin{enumerate}
\item[(a)] $m \equiv 2, 3 \bmod {4}$, $n \equiv 5 \bmod {8}$
\item[(b)] $m \equiv 1 \bmod {8}$, $n \equiv 5 \bmod {8}$
\end{enumerate}

Wir betrachten zuerst den Fall (a) und unterscheiden, ob $n > 0$ ist
oder nicht:

\begin{enumerate}
\item[(a1)] $m \equiv 2, 3 \bmod {4}$, $n \equiv 5 \bmod {8}$,
            $-m, n \in \N$.
  Hier verzweigt $(2)$ in beiden imaginärquadratischen Teilkörpern
  $\Q(\sqrt{-m})$ und $\Q(\sqrt{-mn})$; wären die
  Primideale über $(2)$ in beiden Körpern keine Hauptideale, so könnte
  nach der Klassenzahlformel nicht $H=1$ sein. Also enthält mindestens
  einer der beiden imaginärquadratischen Körper ein Element der Norm
  2, und wegen $m \equiv 2, 3 \bmod {4}$ ist $m \in \{-1, -2\}$. Ist
  $D(-1,n)$ normeuklidisch, so gibt es ein $\alpha \in D(-1,n)$ mit
  $\alpha \equiv (1+\sqrt{-n})/(1+i) \bmod {2}$ und
  $N_{K/\Q}(\alpha)<16 = N_{K/\Q}(2)$. Bezeichnet
  $N_2$ die Norm von $K$ nach $\Q(\sqrt{-n})$, so findet man
  $N_2(\alpha) \equiv \sqrt{-n} \bmod {2}$. Also enthält $k_2 =
  \Q(\sqrt{-n})$ ein Element $\beta = N_2(\alpha) \equiv
  \sqrt{-n} \bmod {2}$ mit $N(\beta) < 16$ ($N$ steht hier natürlich für
  die Norm von $k_2$ nach $\Q$). Wegen $N(\beta) \ge N(\sqrt{-n}) = n$
  muss also $n<16$ sein, und wegen $n \equiv 5 \bmod {8}$ bleiben nur
  die beiden Möglichkeiten $n=5$ und $n=13$.

Dass $D(-1,13)$ nicht normeuklidisch ist, entnimmt man entweder obiger
Tafel oder man verwendet (1.4) mit $L = \Q(\sqrt{-13})$, $K =
\Q(\sqrt{-7})$, $B = (3+2i)$, $x=2$ und $e=4$. Wir brauchen
dann nur noch zu zeigen, dass es kein $r \in \Z[i]$ gibt mit
$N(r)<13$ und $r \equiv 4 \bmod {3+2i}$, sodass $r$ Norm aus $D(-1,13)$
ist. Die beiden ersten Bedingungen werden aber nur von $r = 1-2i$ und
$r = -1+i$ erfüllt; wegen $\left(\frac{13}{2}\right) =
\left(\frac{13}{3}\right) = -1$ bleiben aber beide $r$ in $L$ prim.

Ist $D(-2,n)$ normeuklidisch, so betrachten wir entsprechend die
Restklasse $\alpha \equiv \sqrt{-2n} \cdot (1 + \sqrt{n})/2 \bmod {2}$
und finden wie oben $\beta = N_2(\alpha) \equiv 1 + \sqrt{-2n}
\bmod {2}$ und $N(\beta) < 16$. Aus $n \equiv 5 \bmod {8}$ und $1 + 2n <
16$ folgt aber $n=5$.

\item[(a2)] $m \equiv 2, 3 \bmod {4}$, $n \equiv 5 \bmod {8}$, $-m, -n
  \in \N$
Hier zeigt (1.4) mit $K = \Q(\sqrt{n})$, $L = K(\sqrt{m})$, $B
= \fp_1$ und $\epsilon = (1 + \sqrt{n})/2$ ($\epsilon$ ist
quadratischer Rest $\bmod {\fp_1}$ wegen
$\Phi_K(\fp_1) = 3$), dass $\Q(\sqrt{n})$ ein Element
ungerader Norm $< 4$ enthält; also ist $n=-3$ oder $n=-11$.

Im Falle $n=-11$ muss außerdem noch $(1 + \sqrt{-11})/2$ Norm aus
$D(m,-11)$ sein, folglich enthält $D(m)$ ein Element der Norm 3, und
dies impliziert $m = -2$.

Ist $D(m,-3)$ normeuklidisch, so betrachten wir die Restklasse
$\alpha \equiv \rho + \sqrt{m} \bmod {2}$, wo $\rho$ primitive 3.
Einheitswurzel ist. Damit ist $N_1(\alpha) \equiv m+1 + \sqrt{m} \bmod {2}$.
Im Falle $m \equiv 2 \bmod {4}$ enthält $\Q(\sqrt{m})$ also ein
$\beta \equiv 1 + \sqrt{m} \bmod {2}$ mit $N(\beta) < 16$, und dies
liefert $m \in \{-2, -6, -10, -14\}$. Unter diesen $D(m,-3)$ haben
nur $D(-2,-3)$ und $D(-2,-3) = D(6,-3)$ Klassenzahl 1.

Entsprechend erhält man im Falle $m \equiv 3 \bmod {4}$ nur die
Möglichkeiten $m \in \{-1, -3, -13\}$, wobei aber nur $D(-1,-3)$
Klassenzahl 1 hat.

\item[(b)] $m \equiv 1 \bmod {8}$, $n \equiv 5 \bmod {8}$.
\end{enumerate}

Sei $\alpha \equiv (1 + \sqrt{m})/2 \bmod {2}$; dann ist $N_1(\alpha)
\equiv (\pm 1 + \sqrt{m})/2 \bmod {2}$. Soll also die Restklasse
$\alpha \bmod {2}$ ein Element der Norm $< 16$ enthalten, so enthält
auch $D(m)$ ein solches; dies liefert $m < 64$. Also enthält jeder der
beiden imaginärquadratischen Zahlkörper ein Element der Norm $< 64$,
und von den nur noch verbliebenen Ringen haben nur die folgenden
Klassenzahl 1: $D(-3,5)$, $D(-3,-7)$, $D(-3,-11)$, $D(-3,17)$,
$D(-3,-19)$, $D(-3,-43)$, $D(-7,5)$, $D(-7,-11)$, $D(-7,-19)$,
$D(-7,-43)$ und $D(-11,-19)$. J. Sauvageot\index[dN]{Sauvageot}
hat 1972/73 versucht, durch Betrachten der Restklasse $\bmod {2}$ noch
einige dieser Ringe auszuschließen, hat dabei jedoch einen Fehler
gemacht. In der Tat kommen wir hier nur mit (1.4) (oder 1.12) weiter:

$$ \begin{array}{cc|c|c|c|c|c}
    \toprule
    m & n & x & K & B & r \bmod {B} & M(K)\ge \\ \midrule
    -3 & -43 & \frac{1 + \sqrt{-43}}2 & \Q(\sqrt{-43}) & (3)
         & \sqrt{-43} & \frac{13}{9} \\
    -7 & -11 & 3 & \Q(\sqrt{-7}) & 2 + \sqrt{-7} & -2 & \frac{16}{11} \\
    -7 & -19 & 2 & \Q(\sqrt{-19}) & \frac{3 + \sqrt{-19}}2 & -3 & 1 \\
    -7 & -43 & \frac{3 + \sqrt{-43}}2 & \Q(\sqrt{-43}) & (7)
         & \frac{-3 + 3\sqrt{-43}}2 & \frac{99}{49} \\
   -11 & -19 & 4 & \Q(\sqrt{-19}) & \frac{5 + \sqrt{-19}}2
         & 5 & 1 \\ \bottomrule
\end{array} $$

Damit ist (6.1) bewiesen.

Ein Teil der unter (6.1) aufgeführten Körper lässt sich leicht als
normeuklidisch erkennen. Dazu beachten wir Folgendes: ist $R$ ein
Zahlring, $\alpha, \beta \in R$ und ist auch $(\alpha + \beta)/2$
ganz, so können wir $x\alpha + y\beta$ (mit $x, y \in \Q$)
$\bmod {R}$ so verschieben, dass $|x|, |y| \le \frac{1}{2}$ und
$|x \cdot y| \le \frac{1}{4}$ wird. Offensichtlich können wir ja
$|x|, |y| \le \frac{1}{2}$ erreichen; ist dann z. B.
$x \cdot y \ge \frac{1}{2}$ (und damit $x, y \ge 0$), so setzen wir
$x' = x - \frac{1}{2}$, $y' = y - \frac{1}{2}$ und haben
$|x'|, |y'| \le \frac{1}{2}$ und $|x' \cdot y'| \le \frac{1}{4}$.

Weiter beachten wir, dass $K$ totalkomplex und somit
$N_{K/\Q}(\alpha) \ge 0$ für alle $\alpha \in K$ ist. Um
$|N_{K/\Q}(\alpha)| < 1$ zu zeigen, brauchen wir also nur
$N_{K/\Q}(\alpha) < 1$ zu zeigen. Dies werden wir in der Regel
dadurch tun, dass wir $N_{K/\Q}(\alpha) = X^2 - n Y^2$
schreiben und dann $|X| < 1$ nachweisen; dann ist nämlich
$N_{K/\Q}(\alpha) = X^2 - n Y^2 \le X^2 < 1$ (bessere
Abschätzungen erhält man natürlich, wenn man auch $|Y| \ge c > 0$ für
gewisse $\alpha$ zeigen kann, weil dann $N_{K/\Q}(\alpha) =
X^2 - n c^2$ ist).

Diese Überlegungen werden wir im Folgenden ständig benutzen; so wählen
wir z. B. im Falle $D(-1,2)$ die Basis $\{1, i, \sqrt{2}, \sqrt{-2}\}$,
beachten, dass $\frac{\sqrt{2}+\sqrt{-2}}2$ ganz ist, und verschieben
$x + y i + z \sqrt{2} + w \sqrt{-2} \bmod {R}$ so, dass
$|x|, |y|, |z|, |w| \le \frac{1}{2}$ und $|z| + |w| \le \frac{1}{2}$ wird.
Die letzte Ungleichung liefert $z^2 + w^2 \le \frac{1}{4}$ und damit
$x^2 + y^2 + 2 z^2 + 2 w^2 \le 1$ (wobei genau dann Gleichheit herrscht,
wenn $|x| = |y| = \frac{1}{2}$ und $|z| = \frac{1}{2}$, $w = 0$ oder $z = 0$,
$|w| = \frac{1}{2}$ ist).

Nun ist aber $N_{K/\Q}(x + y i + z \sqrt{2} + w \sqrt{-2}) =
X^2 - 2 Y^2$ mit $X = x^2 + y^2 + 2 z^2 + 2 w^2$ und $Y = 2(x z + y
w)$, folglich $N_{K/\Q}(x + y i + z \sqrt{2} + w \sqrt{-2})
\le X^2 \le 1$. Da Gleichheit höchstens in den beiden oben angegebenen
Fällen eintreten kann, dort aber $|Y| = 2 |x z| = \frac{1}{2}$ und
somit $X^2 - 2 Y^2 = \frac{1}{2}$ ist, haben wir in jedem Fall
$N_{K/\Q}(x + y i + z \sqrt{2} + w \sqrt{-2}) < 1$ erreicht.

In $D(-1,n)$, $n \in \{-3, 5\}$, gehen wir etwas anders vor: hier
wählen wir die Basis $\{1, i, \sqrt{n}, \sqrt{-n}\}$ und beachten, dass
$(1+\sqrt{n})/2$ und $(i+\sqrt{-n})/2$ ganz sind. Nun können wir $x +
y i + z \sqrt{n} + w \sqrt{-n} \bmod {R}$ so verschieben, dass im Falle
\begin{align*}
n & = -3: & x^2 + 3 w^2 &\le \frac{1}{3}, & y^2 + 3 z^2 &\le \frac{1}{3} \\
n & = 5: & x^2 + 5 z^2 &\le \frac{9}{20}, & y^2 + 5 w^2 &\le \frac{9}{20}
\end{align*}
wird (dies folgt aus dem Beweis von (0.21)). Jetzt folgt wieder
\begin{align*}
  N_{K/\Q}& (x + y i + z \sqrt{n} + w \sqrt{-n})
        = X^2 - |n| Y^2 \le X^2 \quad \text{mit} \\
  |X| & = x^2 + y^2 + 3 z^2 + 3 w^2 \le \frac{4}{3} \quad
         \text{im Falle } D(-1,-3) \quad \text{und} \\
  |X| &= x^2 + y^2 + 5 z^2 + 5 w^2 \le \frac{9}{10} \quad
         \text{im Falle } D(-1,5).
\end{align*}

Ganz analog kann man nun zeigen, dass die Ringe $D(-3,2)$, $D(-3,-2)$,
$D(-3,5)$, und $D(-3,-7)$ normeuklidisch sind. Dass $D(-1,2)$ und
$D(-1,3)$ normeuklidisch sind, wusste schon Eisenstein (1850), dessen
Beweis dem obigen recht ähnlich ist. Masley\index[dN]{Masley}
\cite{Mas72} hat 1972 (in Unkenntnis der Eisensteinschen Arbeit)
weitere Beweise hierfür gegeben, die jedoch weit komplizierter als die
oben angegebenen sind. Lakein \cite{Lak72} wies dann (ebenfalls 1972)
nach, dass der EA auch in den Ringen $D(-1,5)$, $D(-1,7)$, $D(-3,2)$,
$D(-3,-2)$, $D(-3,5)$, $D(-3,-7)$, und $D(-3,-11)$ gilt. Dies lässt
sich per Computer nachprüfen; außerdem erhält man die neuen
normeuklidischen Ringe $D(-2,5)$, $D(-3,17)$, $D(-3,-19)$ und
$D(-7,5)$.

Schließlich kann man wie im reellquadratischen Fall die ersten Minima
$M(K)$ für eine ganze Klasse von Zahlkörpern angeben:

\begin{quote}
  \textbf{(6.2)} {\em Sei $m = n^2 + 1$, $n \equiv 1 \bmod {2}$, $R =
    \Z[i, \sqrt{m}, (\sqrt{m} + \sqrt{-m})/2]$, und $f$ der
    Absolutbetrag der Norm. Dann ist $M(f) = \frac{m}{4}$, und dieses
    Minimum wird mod $R$ nur in dem Punkt $(1 + i + \sqrt{m})/2$
    angenommen. Ist $m$ quadratfrei, so ist $M(f) = M(K)$.}
\end{quote}

Wir gehen vor wie in $\Q(\sqrt{m})$ (sh. § 2) und setzen
$\theta = n + \sqrt{m}$. Damit ist
$\{1, i, \theta, \frac{1+i+\theta+i\theta}2 \}$ eine GHB, und mit
$x = a+bi$, $y = c+di$ gilt $N_0(x+y\theta) = x^2 + 2nxy - y^2$, sowie
$N_{0/\Q}(x+y\theta) = |x^2 + 2nxy - y^2|^2$ (hier fassen wir
$x^2 + 2nxy - y^2$ als Element von $\mathbb{C}$ auf; $|\cdot|$ ist der
gewöhnliche Betrag auf $\mathbb{C}$). Für jeden Ausnahmepunkt
$z = x+y\theta$ gilt daher $\frac{m}{4} \le N_{0/\Q}(x+y\theta)$
und somit $\frac{n}{2} < |x^2 + 2nxy - y^2|$.

Nun können wir $z \bmod {R}$ so wählen, dass zuerst $|c|, |d| \le
\frac{1}{2}$ und $|c|+|d| \le \frac{1}{2}$, und dann $|a|,|b| \le
\frac{1}{2}$ wird. Damit ist dann $|x|^2 = a^2 + b^2 \le \frac{1}{2}$
und $|y|^2 = c^2 + d^2 \le \frac{1}{4}$, also
\[
\frac{n}{2} \le |N_0(z)| = |x^2 + 2nxy - y^2|
            \le |x|^2 + 2n |xy| + |y|^2 \le n |x| \cdot \frac{3}{4}.
\]
Dies impliziert $|x| \ge \frac{1}{2} - \frac{3}{4n}$ für jeden
Ausnahmepunkt $z = x+y\theta$. Indem wir $z \bmod {R}$ aber so wählen,
dass zuerst $|a|, |b| \le \frac{1}{2}$, $|a|+|b| \le \frac{1}{2}$ und
dann $|c|,|d| \le \frac{1}{2}$ wird, folgt entsprechend
$|y| \ge \frac{1}{2} - \frac{3}{4n}$. Wenn wir den Fundamentalbereich $F$ nun
geeignet wählen, liegt jeder Ausnahmepunkt $z = x+y\theta$ in der
durch die Ungleichungen
$\bigl||x| - \frac{1}{2}\bigr| \le \frac{3}{4n}$,
$\bigl||y| - \frac{1}{2}\bigr| \le \frac{3}{4n}$
definierten Menge.

Sei jetzt wieder $|c|,|d| \le \frac{1}{2}$ und $|c|+|d| \le
\frac{1}{2}$. Zusammen mit $|y| \ge \frac{1}{2} - \frac{3}{4n}$
folgt wegen
$c^2 + d^2 = |y|^2 \ge \frac{1}{4} - \frac{3}{4n} + (\frac{3}{4n})^2$
und $c^2 \le (\frac{1}{2} - |d|)^2$, dass
\[
\frac{1}{4} - \frac{3}{4n} \cdot \Big(\frac{3}{4n}\Big)^2
   \le c^2 + d^2 \le 2d^2 - |d| + \frac{1}{4}
\]
gilt. Jetzt finden wir für $n \ge 7$
\[
\Big(|d| - \frac{1}{4}\Big)^2
   \ge \frac{1}{16} - \frac{3}{8n} + \frac{9}{32n^2}
   \ge \Big(\frac{1}{4} - \frac{1}{n}\Big)^2,
\]
und wenn wir berücksichtigen, dass $|x|$ und $|y|$ nahe bei
$\frac{1}{2}$ liegen, können wir dies auch für $n \ge 5$ nachweisen.
Also ist $|d| < \frac{1}{n}$, $|c| - \frac{1}{2} < \frac{1}{n}$ oder
$|c| < \frac{1}{n}$, $|d| - \frac{1}{2} < \frac{1}{n}$. Ganz analog
folgen die entsprechenden Ungleichungen für $a$ und $b$. Damit bleiben
die folgenden möglichen Ausnahmemengen ($\delta = \frac{1}{n}$):
\begin{align*}
  S_1 & = \Big(0, \frac{1}{2}, \frac{1}{2}, 0\Big) +
          (-\delta, \delta) \times (-\delta, \delta) \times
          (-\delta, \delta) \times (-\delta, \delta), \\
   \quad \text{und} \\
   S_2 & = \Big(\frac{1}{2}, 0, \frac{1}{2}, 0\Big) +
          (-\delta, \delta) \times (-\delta, \delta) \times
          (-\delta, \delta) \times (-\delta, \delta), \\
\intertext{ da z.B. $S_1$ und }
S_3 & = \Big(\frac{1}{2}, 0, 0, \frac{1}{2}\Big) +
          (-\delta, \delta) \times (-\delta, \delta) \times
          (-\delta, \delta) \times (-\delta, \delta),
\end{align*}
mod $R$ kongruent sind. Jetzt kann man leicht nachrechnen, dass $S_2$
keinen Ausnahmepunkt enthält und somit $S_1$ die einzige Ausnahmemenge
ist.

Eine nochmalige Abschätzung der Norm auf $S$ zeigt nun, dass man
$\delta = \frac{1}{n}$ durch $\epsilon = \frac{1}{3n}$ ersetzen
darf. Dann beachtet man, dass $\theta$ eine Einheit ist und
$(a+bi+c\theta + di\theta)\theta = c+di+(a+2nc)\theta +
(b+2nd)i\theta$ gilt, und stellt fest, dass die Schranken die Anwendung
von (2.3) erlauben. Also ist $z = (i+\theta)/2 \bmod {R}$ der einzig
mögliche Ausnahmepunkt, und wegen $N_{0/\Q}(z) =
\frac{m}{4}$ folgt $M(K) \le \frac{m}{4}$.

Um $M(K)$ nach unten abzuschätzen, gehen wir vor wie folgt: es sei
z. B. $a = d = 0$, $b \equiv c \equiv \frac{1}{2} \bmod {1}$. Dann gelten mit
$x=a+bi$, $y=c+di$ die Kongruenzen $2xy = 2(ac-bd) + 2(ad+bc)i =
\frac{1}{2} \bmod {\Z[i]}$ und $x^2 - y^2 = a^2 - b^2 - c^2 +
d^2 + 2(ab-cd)i = \frac{1}{2} \bmod {\Z[i]}$. Also ist
$|\operatorname{Re}(x^2 + 2nxy - y^2)| \ge \frac{1}{2}$ und
$|\operatorname{Im}(x^2 + 2nxy - y^2)| \ge \frac{1}{2}$, folglich
$N_{0/\Q}(z-a) \ge (1+n^2)/4 = \frac{m}{4}$ für alle $a \in
R$. Damit sind alle Behauptungen bewiesen.

\section*{{\sc Anmerkungen zu} \S\ 6}
\addcontentsline{toc}{section}{{\sc Anmerkungen zu} \S\ 6}
\medskip
         
\begin{tabular}{p{1cm} p{10cm}}
\toprule
\textbf{Jahr} & \textbf{Ereignis} \\
\midrule
1850 & Eisenstein\index[dN]{Eisenstein} zeigt im Zuge seiner Untersuchungen
 über achte Potenzreste, dass $D(-1,2)$ normeuklidisch ist und bemerkt,
  dass sich  $D(-1,3)$ ganz analog behandeln lässt. \\
1894 & Hilbert\index[dN]{Hilbert} schreibt eine Abhandlung über Dirichletsche
  Zahlkörper, in welcher er insbesondere das quadratische
  Reziprozitätsgesetz in diesen Körpern untersucht. \\
1943 & Kuroda\index[dN]{Kuroda} beginnt, die Dirichletschen Zahlkörper zu
  studieren; diese Untersuchungen werden von Kubota\index[dN]{Kubota}
  \cite{Kub56} und Wada\index[dN]{Wada} \cite{Wad66} fortgeführt. \\ 
1958 & Hasse\index[dN]{Hasse} \cite{Has85} bestimmt alle imaginären bizyklischen
  Körper  mit ungerader Klassenzahl. \\
1970 & Williams\index[dN]{Williams} veröffentlicht eine elementare
  Berechnung der GHB spezieller Dirichletscher Zahlkörper. \\
1972 & Lakein\index[dN]{Lakein} bestimmt einige normeuklidische Dirichletsche
  Zahlkörper; unabhängig davon zeigt Masley, dass $D(-1,2)$ und
  $D(-1,3)$ normeuklidisch sind. Sauvageot gibt erstmals einige
  Dirichletsche Zahlkörper an, die zwar Klassenzahl 1 haben, aber
  nicht normeuklidisch sind.\\
1974 & Brown\index[dN]{Brown} und Parry\index[dN]{Parry} bestimmen alle
  imaginären bizyklischen Zahlkörper mit Klassenzahl 1. Dies wurde
  durch die Lösung des ``Klassenzahl-2''-Problems in imaginärquadratischen
  Körpern ermöglicht. \\
\bottomrule
\end{tabular}

\chapter*{\S\ 7 Dirichletsche Zahlkörper}
\setcounter{chapter}{7}
\addcontentsline{toc}{chapter}{\S\ 7 Dirichletsche Zahlkörper}
\markboth{Euklidische Ringe}{\S\ 7 Dirichletsche Zahlkörper}

Sei $K = \Q(i)$ und $R = \Z[i]$; $R$ hat mit dem Ring $\Z$ viele
Eigenschaften gemein, von denen wir hier einige aufzählen
wollen. Dabei werden wir Elemente aus $R$ mit kleinen ($a$, $b$, $p$,
$q$), ganze rationale Zahlen dagegen mit großen ($A$, $B$, $P$, $Q$ usw.)
Buchstaben bezeichnen. Unter dem Begriff \glqq{}rationale Primzahl\grqq{}
wollen wir immer eine natürliche Zahl verstehen. Weiter sei $N$ immer
die Norm von $K$ nach $\Q$.

\begin{enumerate}
\item Die Einheitengruppe $R^\times$ ist endlich:
  $R^\times = \{ \pm 1, \pm i \}$;
\item $R$ ist ZPE-Ring;
\item jedes prime $p \in R$ ist zu einer der folgenden Zahlen assoziiert:
  \begin{enumerate}
  \item $q = 1+i$, $N(q) = 2$,
  \item $p = A+Bi$, $A-1 \equiv B \equiv 0 \bmod {2}$,
    $N(p) = A^2+B^2=P$, wo $P \equiv 1 \bmod {4}$ rationale Primzahl ist,
  \item $q = Q$, wo $Q \equiv 3 \bmod {4}$ rationale Primzahl ist;
  \end{enumerate}
\item Ist $a \in R$, so gelten die Kongruenzen
  $$ \begin{aligned}
  a^2 &\equiv 0 \bmod {4} &&\Leftrightarrow && a \equiv 0 \bmod {2} \\
  a^2 &\equiv 2i \bmod {4} &&\Leftrightarrow && a \equiv 1+i \bmod {2} \\
  a^2 &\equiv 1 \bmod {4} &&\Leftrightarrow && a \equiv 1 \bmod {2} \\
  a^2 &\equiv 3 \bmod {4} &&\Leftrightarrow && a \equiv i \bmod {2};
  \end{aligned} $$
\item Ist $p \in R$ prim und $(a,p)=1$ für ein $a \in R$, dann gilt
  $a^{Np-1} \equiv 1 \bmod {p}$.
\end{enumerate}

Der \glqq{}kleine Fermat\grqq{} in 5. gibt uns die Möglichkeit, in $R$ ein
quadratisches Restsymbol zu definieren: dazu seien $a,p \in R$, $p$
prim, $Np \equiv 1 \bmod {2}$ und $(a,p)=1$. Dann setzen wir $[a/p] =
\pm 1$, wobei das Vorzeichen mit demjenigen in $a^{(Np-1)/2} \equiv
\pm 1 \bmod {p}$ übereinstimmt. Dieses Restsymbol lässt sich in Analogie
zum Jacobi-Symbol auch auf zusammengesetzte Nenner verallgemeinern;
damit gilt dann

\begin{enumerate}
\item[6.] Das quadratische Reziprozitätsgesetz in $\Z[i]$: Für $a,b \in R$
  mit $a \equiv b \equiv 1 \bmod {2}$ gilt $[a/b] = [b/a]$. Weiter gelten
  mit $a = A+Bi$ und $Q = A+B$ die Ergänzungssätze $[\frac ia] = (-1)^{B/2}$
  und $[\frac{1+i}a] = (\frac{2}{Q})$. Ist außerdem $b = C+Di$ prim und
  $D \neq 0$, so lässt sich mittels $[\frac ab] = (\frac{AC + BD}p)$,
  $P = N(b) = C^2+D^2$ das Restsymbol $[\frac ab]$ auf ein Legendre-Symbol
  zurückführen; diese Identität ermöglicht es auch, das QRG in $R$
  aus demjenigen in $\Z$ herzuleiten.
\item[7.] Sei $L = K(\sqrt{m})$, $m \in R$ quadratfrei, $S = \cO_L$ der
  Ring ganzer Zahlen in $L$; dann hat $S$ die Ganzheitsbasis
  $\{1, \beta\}$ über $R$, sowie die Relativdiskriminante $d$:
  \begin{align*}
    \text{I.} \quad & m \equiv 1 \bmod {4} : \quad
      \beta = \frac{1+\sqrt{m}}{2}, \quad d = m; \\
    \text{II.} \quad & m \equiv \pm 1 + 2i \bmod {4} : \quad
      \beta = \frac{1+\sqrt{m}}{1+i}, \quad d = 2im; \\
    \text{III.} \quad & m \equiv i \bmod {2} \text{ oder }
      m \equiv 0 \bmod {(1+i)} : \quad
      \beta = \sqrt{m}, \quad d = 4m.
  \end{align*}
\end{enumerate}

Ist $\alpha = r + s\sqrt{m}$ für $r,s \in K$, so heißt
$\alpha' = r - s\sqrt{m}$ die Konjugierte von $\alpha$. Damit ist die
Relativdiskriminante durch $d = (\alpha - \alpha')^2$ gegeben, wobei
(auch im Folgenden) zu beachten ist, dass $d$ nur bis auf einen Faktor
$\pm 1$ definiert ist. Allerdings hat dieser Faktor keinen Einfluss auf
das Restsymbol $[\frac{d}{p}]$ oder die Erweiterung
$K(\sqrt{d})$, weil $\pm 1$ Quadrat ist. Schließlich errechnet sich
die Absolutdiskriminante von $L$ aus $d$ durch
$\disc (L/\Q) = 16N(d)$.

Nach 7. gilt für die Relativdiskriminante $d$ quadratischer
Erweiterungen von $K$ die Kongruenz $d \equiv 0 \bmod {2}$ oder $d
\equiv 1 \bmod {4}$; in der Tat gilt dies für Erweiterungen beliebigen
Grades: wäre z. B. für die Relativdiskriminante $d$ einer Erweiterung
$L/K$ $d \equiv \pm 1 \bmod {4}$, so wäre $L/K$ sicher nicht galoissch
(weil $d$ kein Quadrat in $K$ ist). Weiter ist $(1+i)$ in $L/K$
unverzweigt, und folglich gilt dasselbe auch für den normalen Abschluss
$M/K$. Dieser normale Abschluss $M$ enthält aber den quadratischen
Teilkörper $L(\sqrt{d})$, und wir sehen, dass $(1+i)$ auch in
$L(\sqrt{d})$ unverzweigt ist. Nach dem Zerlegungsgesetz in
relativquadratischen Erweiterungen (sh. Hilbert) muss dann die
Kongruenz $x^2 \equiv d \bmod {4}$ in $K$ lösbar sein, was wegen $d
\equiv \pm 1 \bmod {4}$ aber nicht der Fall ist.

Derselbe Beweis funktioniert offensichtlich auch für den Fall
$d \equiv \pm 1 + 2i \bmod {4}$. Ist schließlich $d \equiv 0 \bmod {(1+i)}$,
so ist $(1+i)$ in $L/K$ verzweigt; bezeichnet $e = e(L/K)$ den
Verzweigungsindex, so ist entweder $e=2$ oder $e \ge 3$, und in beiden
Fällen folgt nach bekannten Sätzen (Stichwort: wilde Verzweigung)
$(1+i)^2 \mid d$.

Um das Zerlegungsgesetz in $L$ durch das Restsymbol
$[\frac{d}{p}]$ beschreiben zu können, erweitern wir es
analog zum Kronecker-Symbol und setzen für $p=1+i$
\[
\left( \frac{d}{p} \right) =
\begin{cases}
+1, & \text{falls } d \equiv 1 \bmod {p^5} \\
-1, & \text{falls } d \equiv 5 \bmod {p^5}
\end{cases}
\]

Damit gilt dann

\begin{enumerate}
\item[8.] Sei $p \in R$ prim, $L/K$ quadratische Erweiterung von $K$
  mit Relativdiskriminante $d$; dann ist im Falle
  $$ \begin{array}{rclrcl}
    {} [\tfrac{d}{p}] & = & +1: & (p) & = & P_1 P_2
         \text{ für zwei verschiedene Primideale } P_1, P_2 \text{ in } S; \\
    {} [\tfrac{d}{p}] & = & -1: & (p) & = & (p) \text{ bleibt prim;} \\
    {} [\tfrac{d}{p}] & = &  0: & (p) & = & P^2 \text{ verzweigt in } L/K.
  \end{array} $$
\end{enumerate}

Aus der Tatsache, dass $L$ totalkomplex ist, folgt erstens
$N_{L/\Q}(\alpha) \ge 0$ für alle $\alpha \in L$, und zweitens nach
dem Dirichletschen Einheitensatz
\begin{enumerate}
\item[9.] Jede Einheit $e \in S^\times$ lässt sich in der Form $e =
  \zeta u^k$ schreiben, wo $\zeta$ eine Einheitswurzel und $u$ die
  Fundamentaleinheit von $L$ ist. Dabei ist $\zeta = \zeta_8$, falls
  $L = K(\sqrt{i}) = \Q(i,\sqrt{2})$ ist; $\zeta = \zeta_{12}$, falls
  $L = K(\sqrt{-3})$ ist; $\zeta = \zeta_4 = i$ sonst.
\end{enumerate}

Während man sich im Falle reellquadratischer Zahlkörper dafür
interessiert, ob die Norm der Fundamentaleinheit gleich $\pm 1$ oder
gleich $-1$ ist, lautet die entsprechende Frage hier, ob $N_{L/K}(u) =
\pm 1$ oder gleich $\pm i$ ist, und auch die Antworten sind weitgehend
analog.

Die bisher aufgezählten Eigenschaften Dirichletscher Zahlkörper
(damit werden die quadratischen Erweiterungen von $K$ bezeichnet) sind
klassisch und finden sich bereits bei Dirichlet und Hilbert. Mit ganz
elementaren Methoden lassen sich aber auch hier weitergehende
Ergebnisse erzielen.

Unser erstes Resultat betrifft die Parität der Klassenzahl von $L$;
damit hat sich bisher nur Popović\index[dN]{Popović} \cite{Pop38} befasst;
dieser zeigte, dass die Klassenzahl $h=h(L)$ nur dann ungerade sein
kann, wenn die Relativdiskriminante $d$ prim ist oder $d \equiv \pm 1
\bmod {4}$ gilt und $d$ genau zwei Primteiler besitzt (übrigens sind
die in der Besprechung dieser Arbeit (Reviews in number theory
1940--1972, R 16-33) angegebenen Behauptungen nicht richtig). Dieses
Ergebnis lässt sich jedoch noch verschärfen, und dazu brauchen wir

\begin{quote}
  \textbf{(7.1)} {\em Sei $p \equiv \pm 1 \bmod {4}$ prim in $R$, dann
    gibt es $a, b \in R$ mit $p = a^2 + i b^2$, $a \equiv 1 \bmod {(1+i)}$. }
\end{quote}

\begin{proof}
  Sei $L = K(\sqrt{i}) = \Q(i,\sqrt{2})$; $L$ ist nach \S\ 1 (oder auch
  \S\ 5) normeuklidisch und hat damit Klassenzahl 1. Wegen
  $[\frac{d}{p}] = +1$ ist $p$ in $S$ zerlegt,
  folglich gibt es ein $\pi \in S$ mit $N_{L/K}(\pi) = \zeta p$, wo
  $\zeta$ eine Einheit in $R$ ist (also $\zeta = \pm i, \pm 1$). Weil
  $\{1, \sqrt{i}\}$ eine Ganzheitsbasis von $S$ über $R$ ist,
  dürfen wir $\pi = a + b\sqrt{i}$ für gewisse $a, b \in R$ schreiben
  und finden $\zeta p = a^2 + i b^2$. Division durch $\zeta$ liefert
  nun die Behauptung.
\end{proof}

\textbf{Bem.:} Das Analogon in $\Z$ ist die Darstellbarkeit von primen
$P \equiv 1 \bmod {4}$ als Summe zweier Quadrate; auch im Beweis von
(7.2), der im reellquadratischen Fall auf Rédei\index[dN]{Redei@R\'edei}
\cite{Red60} zurückgeht, übernimmt $p = a^2 + i b^2$ die Rolle von
$P = a^2 + b^2$ in Redeis Beweis.

\begin{quote}
  \textbf{(7.2)} {\em Sei $L = K(\sqrt{m})$, $m \in R$ quadratfrei. Ist
    dann $h(L) \equiv 1 \bmod {2}$, dann kommen nur folgende
    Möglichkeiten in Betracht:}
  \begin{enumerate}
  \item $m = i$
  \item $m = p$, $p \in R$ {\em prim}, $p \not\equiv \pm i \bmod {4}$
  \item $m = pq \equiv \pm 1 \bmod {4}$ {\em für prime} $p, q \in R$,
    {\em wobei entweder $p, q \equiv \pm 1 + 2i \bmod {4}$ oder
      $p \equiv 1 \pm i$, $q \equiv \pm 1 + 2i \bmod {4}$ ist.}
  \end{enumerate}
\end{quote}

\begin{proof}
  Sei $d$ die Relativdiskriminante von $L/K$, wie sie in 7. definiert
  wurde. Da $d$ (nach 7.) keine Einheit ist, gibt es prime $p \in R$
  mit $p \mid d$. Nach dem Zerlegungsgesetz gilt $(p) =
  \mathfrak{p}^2$ für jedes solche $\mathfrak{p}$, wobei
  $\mathfrak{p}$ ein Primideal in $S$ über $p$ ist. Da $L$ nach
  Voraussetzung ungerade Klassenzahl hat, muss mit $\mathfrak{p}^2$
  auch $\mathfrak{p}$ ein Hauptideal sein, und es gibt ein $\pi \in
  S$, sowie eine Einheit $e \in S^\times$ mit $p e = \pi^2$, d.h. mit
  $\sqrt{p e} \in S$ (diese Einheit $e$ hängt natürlich von dem
  jeweiligen $p$ ab, das man gerade betrachtet).

  Indem wir $e$ gegebenenfalls mit einer geeigneten Potenz der
  Fundamentaleinheit $u$ multiplizieren, dürfen wir uns auf die
  Möglichkeiten $e=\zeta$ und $e=\zeta u$ beschränken, wobei $\zeta$
  eine 4. Einheitswurzel ist.
  Ist nun $e=\zeta$ für irgendein $p$ mit $p \mid d$, so folgt
  $\sqrt{\zeta p} \in S \setminus R$, und es muss $L = K(\sqrt{e})$,
  also $m = e p$, sein: somit liegt Fall i) oder Fall ii) vor. Sei
  also $m \ne i$ und $e = \zeta u$ für alle $p \mid d$ (wobei $\zeta$,
  wie schon bemerkt, von $p$ abhängt und nicht für jedes $p$ gleich zu
  sein braucht). Dann folgt
  \begin{enumerate}
  \item[I.] $N_{L/K}(u) = 1$ (und nicht $= i$): mit $p e$ ist nämlich
    auch $pe'$ ein Quadrat in $S$, d.h. es gilt $\sqrt{p e} \cdot
    \sqrt{p e'} = p \sqrt{e e'} \in S$; wegen $m \ne i$ liegt
    $\sqrt{i}$ nicht in $S$, folglich muss $e e' = \pm u u' = \pm 1$
    sein.
  \item[II.] $d$ hat maximal zwei Primteiler: ist nämlich $p \mid d$
    und $q \mid d$ für prime, nicht assoziierte $p,q \in R$, so folgt
    mit $e_1 = \zeta_1 u$ und $e_2 = \zeta_2 u$: $\sqrt{p e_1} \cdot
    \sqrt{q e_2} = u \sqrt{p q \zeta_1 \zeta_2} \in S$, und wir sehen
    $m \sim p q$. Sind außerdem $p$ und $q$ von $1+i$ verschieden, so
    muss die Kongruenz $m \equiv \pm 1 \bmod {4}$ gelten, da sonst
    außer $p$ und $q$ wegen 7. auch $1+i$ ein Teiler von $d$ wäre.
    Jetzt müssen wir noch folgende Möglichkeiten ausschließen (an
    dieser Stelle gehen wir über Popovićs Ergebnisse hinaus):
  \begin{enumerate}
  \item[(a)] $m = p q$, $p = 1 \pm i$, $q \equiv \pm 1, \pm i \bmod {4}$:
    (im Falle $q \equiv 2 \pm i \bmod {4}$ ersetzen wir $p = 1 \pm i$
    durch $p = 1-i$, sodass wir diesen Fall nicht extra behandeln
    müssen). Nach (7.1) existieren $a,b \in R$ mit $m = p q = a^2 + i
    b^2$ und $a \equiv b \equiv 1 \bmod {1+i}$. Wir setzen nun $\pi = a
    \cdot \sqrt{m}$ und $\rho = (b, \pi)$; dabei ist $\pi$ natürlich
    nur bis auf eine Einheit in $S$ bestimmt. Wir sehen nun
    \begin{itemize}
    \item $(\pi, \pi') = 1$: ein gemeinsamer Teiler $\sigma$ von $\pi$
      und $\pi'$ würde nämlich auch $\pi \cdot \pi' = 2a$ und $\pi -
      \pi' = 2\sqrt{m}$ teilen; wegen $(a,m)=1$ ist aber erst recht
      $(a,\sqrt{m})=1$, sodass $\sigma$ Teiler von $(2)$ sein
      müsste. Daraus folgt aber $(1+i) N_{L/K}(\pi) = a^2 - m = i b^2$
      im Widerspruch zu $b \equiv 1 \bmod {1+i}$.
    \item $\rho^2 / \pi$ ist eine Einheit; dies folgt aus
      $(\rho^2) = (b,\pi)^2 = (b^2, b\pi, \pi^2) = (\pi \pi', b\pi,
      \pi^2) = \pi (\pi', b, \pi) = (\pi)$.
    \end{itemize}
    Nun ist $N_{L/K}(\rho) = \zeta b$ für eine 4. Einheitswurzel
    $\zeta$, also $N_{L/K}(\rho)^2 = \pm b^2$, und weiter haben wir
    oben schon gesehen, dass $N_{L/K}(\pi) = i b^2$ gilt. Also hat die
    Einheit $\rho^2 / \pi$ die Norm $N_{L/K}(\rho^2 / \pi) = \pm b^2 /
    i b^2 = \pm i$ im Widerspruch zu I.
  \item[(b)] $m = p q \equiv 1 \bmod {4}$, $p \equiv q \equiv 1 \bmod {4}$:
    wie in (a) ist auch hier $m = a^2 + i b^2$; jedoch gilt nun $a-1
    \equiv b \equiv 0 \bmod {2}$. Wir setzen $\pi = (a \cdot
    \sqrt{m})/2$, $\rho = (\pi', b/2)$ und schließen wie oben, dass
    $\rho^2 / \pi$ Einheit mit Norm $\pm i$ ist.
  \item[(c)] $m = p$, $p \equiv \pm i \bmod {4}$: wir werden zeigen, dass
    $M = L(\sqrt{i})$ eine unverzweigte, quadratische (und damit abelsche)
    Erweiterung von $L$ ist; nach der Klassenkörpertheorie (oder
    Hilberts Satz 94) muss dann die Klassenzahl von $L$ gerade sein.
    Da $L$ totalkomplex ist, genügt es, endliche Primstellen zu
    betrachten. Weil $i$ Einheit ist, können in $M/L$ höchstens
    Primideale über $(2)$ verzweigen. Wäre dies der Fall, müsste $(2)$
    in $M$ rein verzweigt sein. Mit $F = \Q(\zeta_8)$ lässt sich aber
    $M = F(\sqrt{m})$ schreiben, und nach dem Zerlegungsgesetz für
    relativquadratische Erweiterungen (Hilbert) sind die Primideale
    über $(2)$ in $M/F$ unverzweigt (man beachte die Kongruenz $m
    \equiv i \equiv (\sqrt{i})^2 \bmod {4}$ in $F$).
  \end{enumerate}  
  \end{enumerate}
  
\end{proof}

\textbf{Bem.:} Es wäre wünschenswert, auch für den Fall c) einen
Beweis zu haben, der ohne Klassenkörpertheorie auskommt.

Es ist interessant festzustellen, dass sich die Analogie der
Dirichletschen zu quadratischen Zahlkörpern weiter fortsetzt: ganz
genauso wie Cohn\index[dN]{Cohn} \cite{Coh62} lässt sich nämlich nun
zeigen:

\begin{itemize}
\item Sei $m = p \equiv \pm 1 \bmod {4}$ prim und $L = K(\sqrt{m})$;
  dann gibt es eine Einheit $u \in S$ mit $N_{L/K}(u) = i$, und $L$
  hat ungerade Klassenzahl.
\item Seien $p,q \equiv \pm 1 + 2i \bmod {4}$ prim und nicht assoziiert,
  $m = pq$, $L = K(\sqrt{m})$; schreibt man $pS = P^2$, $qS = Q^2$ für
  Primideale $P,Q$ in $S$, dann sind $P$ und $Q$ Hauptideale. Weiter
  gibt es $x,y \in R$ mit $4\zeta = p x^2 + q y^2$, wo $\zeta$ eine
  4. Einheitswurzel ist.
\end{itemize}

Genau wie bei Cohn kann man aus diesen beiden Tatsachen einen Beweis
des QRG in $R$ herleiten.

Wir wenden uns nun der Frage zu, welche Dirichletschen Zahlkörper
normeuklidisch sind. Diese Frage lässt sich wie im reellquadratischen
Fall dann am einfachsten behandeln, wenn die Relativdiskriminante
durch eine hohe Potenz von 2 teilbar ist:

\begin{quote}
  \textbf{(7.3)} {\em Ist die Relativdiskriminante $d$ von $L/K$ durch
    $4$ teilbar, dann ist $L$ genau dann normeuklidisch, wenn $m=i$ oder
    $m=1+i$ gilt.}
\end{quote}

Zum Beweis verwenden wir

\begin{quote}
  \textbf{(7.4)} {\em Ist die Relativdiskriminante $d$ von $L/K$ durch
    $4$ teilbar, die Klassenzahl $h$ von $L$ ungerade und
    $m \not\equiv i, 1+i$, so gilt für die Fundamentaleinheit $u$ die
    Kongruenz $u \equiv 1 \bmod {(1+i)}$.}
\end{quote}

\begin{proof}
  Da das Ideal $(1+i)$ in $L/K$ verzweigt ist, gibt es ein Primideal
  $\mathfrak{p}_1$ der Absolutnorm $2$ in $L$. Da $L$ ungerade
  Klassenzahl haben soll, ist $\mathfrak{p}_1$ ein Hauptideal, d.h. es
  gibt ein $\pi \in S$ mit $(\pi)^2 = (1+i)$. Wegen $d \equiv 0
  \bmod {4}$ ist $\{1, \sqrt{m}\}$ eine Ganzheitsbasis über $R$, sodass
  wir $\pi = a + b\sqrt{m}$ für gewisse $a, b \in R$ schreiben
  können. Also ist $(a + b\sqrt{m})^2 = (1+i)e$ für eine Einheit $e
  \in S$. Nach dem Einheitensatz ist $e = \zeta u^k$ ($\zeta$ ist eine
  4. Einheitswurzel, da $L$ wegen $m \ne i$ keine 8., wegen $m \ne -3$
  keine 12. Einheitswurzeln enthält). Weiter ist $k \equiv 1 \bmod {2}$,
  da sonst $\sqrt{1 \pm i} \in L$ folgte. Indem wir $a + b\sqrt{m}$
  notfalls mit einer geeigneten Potenz von $u$ multiplizieren, dürfen
  wir $k=1$ annehmen.  Nun ist $\pi \equiv \pi' \bmod {2}$, folglich
  $$ \zeta u = \frac{\pi^2}{1+i} = \frac{\pi \pi'}{1+i}
      = \frac{1 \pm i}{1+i} \equiv 1 \bmod (1+i), $$
  und da $\zeta$ eine 4. Einheitswurzel ist, gilt auch
  $\zeta \equiv 1 \bmod {(1+i)}$.
\end{proof}

\begin{proof}[Bew. von (7.3)]
  Wir zeigen $M(L,I) \ge 5/4$ für das Ideal $I = (1+i)$; wegen
  $\Phi_L(I) = 2$ gibt es zwei prime Restklassen mod $I$. Dabei
  enthält die von $1 \bmod {I}$ verschiedene Restklasse keine Einheiten
  wegen (7.2); weiter enthält sie keine Elemente der Norm $3$ (weil
  $(3)$ in $K$ prim ist). Die Behauptung folgt.
\end{proof}

Wir geben noch eine kurze Tabelle einiger euklidischer Minima:

\begin{center}
$$ \begin{array}{|l|c|c|c|}
\hline
\rsp & \disc K & M(K) & C_1 \\ \hline
\rsp i    & 144 & \frac{1}{2} & \text{sh. \S\, 6, } D(-1,2) \\
\rsp 1+i  & 512 & \frac{1}{2} & \\ 
\rsp 2+i  & 1280 & \frac{5}{4} & M_2 \le 0.999 \\
\rsp 3i   & 2304 & \frac{5}{2} & \\ 
\rsp 3+ i & 2560 & \frac{5}{4} & M_2 \le 0.999 \\
\rsp 1+3i & 2560 & \frac{5}{4} & M_2 \le 0.999 \\ \hline
\end{array} $$
\end{center}

Der Fall $L = K(\sqrt{p})$, $p \equiv \pm 1+2i \bmod {4}$, ist schon
etwas schwieriger; hier gilt aber

\begin{quote}
  \textbf{(7.5)} {\em Sei $L = K(\sqrt{p})$, $p \equiv \pm 1+2i \bmod 4$
    prim; gibt es dann $s, t \in R$ mit $s \equiv t \equiv 1 \bmod 2$,
    $(\frac{s}{p}) = (\frac{t}{p}) = -1$, $(s,t)=1$ und
    $|st| \le k \cdot |p|$ für $k = \sqrt{2}/(1+\sqrt{2})$,
    dann ist $L$ nicht normeuklidisch.}
\end{quote}

Hierbei ist $|\cdot|$ der gewöhnliche Betrag auf $\mathbb{C}$.

\begin{proof}
   Wir benutzen (1.4) mit $B = (1+i)p$, $e=st$ und $k=1$; zu zeigen
   ist, dass es in $R$ kein $r \equiv e \bmod {(1+i)p}$ gibt mit $N(r) <
   2P$ und $r = N_{L/K}(\alpha)$ für ein $\alpha \in S$ (hier ist $P =
   N(p) = |p|^2$). Nun gilt $N(r) < 2P$, also $|r| < \sqrt{2}|p|$. Wir
   setzen $r = st + l(1+i)p$ für ein $l \in R$ und finden wegen $|st|
   \le k|p|$ und $|r| < \sqrt{2}|p|$:
   \[ |l| \cdot \sqrt{2} |p| = |r - st| \le |r| + |st|
        \le (\sqrt{2} + k)|p| = 2|p|.\]
   Also ist $|l| < \sqrt{2}$ und damit $l \in \{0, \pm 1, \pm i\}$.
   Wäre $l \equiv 1 \bmod (1+i)$, so folgte $r \equiv i \bmod 2$;
   ein solches $r$ kann aber keine Norm aus $S$ sein wegen (7.6).
   Also ist $l=0$; dann ist $r=st$ aber keine Norm aus $S$ wegen
   $(s,t)=1$ und $(\frac{s}{p}) = (\frac{t}{p}) = -1$ (dies folgt wie
   im reellquadratischen Fall).
\end{proof}

\begin{quote}
  \textbf{(7.6)} {\em Sei $L = K(\sqrt{m})$,
    $m \equiv \pm 1 \pm 2i \bmod {4}$ und $\alpha \in S$. Dann ist}
    $N_{L/K}(\alpha) \equiv \pm i \bmod {2}$.
\end{quote}

\begin{proof}
  Wir müssen noch (7.6) beweisen: sei dazu
  $\alpha = \frac{a + b\sqrt{m}}{1+i}$, $a \equiv b \bmod {1+i}$, und
  $N_{L/K}(\alpha) = c + di = \beta$. Dann ist $a^2 - m b^2 = 2i \beta$,
  und indem man diese Gleichung mod $m$ und mod $\beta$
  betrachtet, folgt $[ \frac{\beta}{m} ] = [ \frac{m}{\beta}] = \pm 1$.
  Wäre nun $\beta \equiv \pm i \bmod {2}$, so folgte aus dem QRG in $R$
  \[ \left[ \frac{\beta}{m} \right] = \left[ \frac{i}{m} \right]
     \left[ \frac{m}{\beta} \right] = \left[ \frac{i}{m} \right]
     \left[ \frac{m}{i} \right] = \left[ \frac{i}{m} \right]
     \left[ \frac{m}{i} \right], \]
  was aber wegen $[ \frac{i}{m}] = -1$ im Widerspruch zu
  obigen Beobachtungen steht.
\end{proof}

\begin{quote}
  \textbf{(7.7)} {\em Sei $L = K(\sqrt{p})$, $p \equiv \pm 1 + 2i \bmod 4$;
    dann ist $L$ genau dann normeuklidisch, wenn $p \in \{1+2i, 3+2i,
    5+2i, 1+6i, 7+2i\}$ ist.}
\end{quote}

\begin{proof}
  Dass $L$ für $p = 1+2i$ und $p = 3+2i$ normeuklidisch ist, hat schon
  Lakein\index[dN]{Lakein} \cite{Lak72} gezeigt. Ist $L$ normeuklidisch,
  so hat van der Linden\index[dN]{Linden@van der Linden}
  \cite[S. 163]{Lin85} gezeigt, dass $\disc L < 16 \cdot 899\,225$
  ist. Da in unserem Fall $\disc L = 64 \cdot N(p)$ ist,
  dürfen wir $N(p) < 224\,805$ annehmen. Für diese $p$ findet man mit
  einem Computer eine Darstellung wie in (7.4), es sei denn, es ist
  $$ N(p) \in \{5, 13, 29, 37, 53, 61, 101, 157, 181, 229, 349, 541\}.$$
  Beispielsweise ist

\begin{center}
\begin{tabular}{|c|c|c|c||c|c|c|c|}
\hline
$p$ & $p$ & $s$ & $t$ & $P$ & $P$ & $s$ & $t$ \\
\hline
109 & $3+10i$ & $1-2i$ & $1+2i$ & 293 & $17+2i$ & $1-2i$ & 3 \\
149 & $7+10i$ & $1-2i$ & $1+2i$ & 317 & $11+14i$ & $1-2i$ & 3 \\
173 & $13+2i$ & $1+2i$ & 3 & 373 & $7+18i$ & $1+2i$ & $3-2i$ \\
197 & $1+14i$ & $1-2i$ & 3 & 389 & $17+10i$ & $1-2i$ & $1+2i$ \\
269 & $13+10i$ & $1-2i$ & $1+2i$ & 397 & $19+6i$ & $1-2i$ & $3-2i$ \\
277 & $9+14i$ & $1+2i$ & $1+4i$ & 421 & $15+14i$ & $1-2i$ & $1+2i$ \\
\hline
\end{tabular}
\end{center}

Die $p \ge 101$ werden wir nun direkt mit (1.4) ausschließen:

\begin{center}
\begin{tabular}{|c|c|c|c|}
\hline
$P$ & $p$ & $s$ & $B$ \\
\hline
61 & $5+6i$ & 5 & $(1+i)p$ \\
101 & $1+10i$ & $4-i$ & $(p)$ \\
157 & $11+6i$ & $-4+5i$ & $(p)$ \\
181 & $9+10i$ & $-4+5i$ & $(p)$ \\
229 & $15+2i$ & $4+i$ & $(p)$ \\
349 & $5+18i$ & $7i$ & $(p)$ \\
541 & $21+10i$& $-5+14i$& $(1+i)p$ \\
\hline
\end{tabular}
\end{center}
\end{proof}

Die für (7.7) nötige Rechenzeit lässt sich etwas verkürzen, wenn man
das folgende Korollar von (7.6) benutzt:

\begin{quote}
  \textbf{(7.8)} {\em Sei $p \equiv \pm 1, 2i \bmod {4}$ prim, $N(p)
    \ge 109$, und $a \equiv \pm 2$, $b \equiv 0 \bmod {5}$ oder $a
    \equiv 0$, $b \equiv \pm 1 \bmod {5}$. Dann ist $L=K(\sqrt{p})$
    nicht normeuklidisch.}
\end{quote}

\begin{proof}
  (7.5) mit $s=1-2i$, $t=1+2i$.
\end{proof}

Auch hier geben wir einige euklidische Minima an:

\begin{center}
\begin{tabular}{|l|c|c|c|}
\hline
$m$ & $\operatorname{disc} K$ & $M(K)$ & $C_1$ \\
\hline
$1+2i$ & 320 & $\frac{1}{2}$ & \\
$3+2i$ & 832 & $\frac{1}{2}$ & \\
$5+2i$ & 1856 & & \\
$1+6i$ & 2368 & & \\
$3+6i$ & 2880 & & \\
$7+2i$ & 3392 & $\ge 50/53$ & \\
\hline
\end{tabular}
\end{center}

\textbf{Bem.:} Der nächstschwierigere Fall ist nun $m \equiv 5
\bmod {(1+i)^6}$; hier gilt

\begin{quote}
  \textbf{(7.9)} {\em Sei $m \equiv 5 \bmod {(1+i)^6}$, $m \ne \pm 3$,
    und $u$ die Fundamentaleinheit von $L$. Hat dann die Form $u =
    \frac{x + y\sqrt{m}}{1+i}$ mit $x \equiv y \equiv 1 \bmod {1+i}$,
    dann ist $L$ nicht normeuklidisch.}
\end{quote}

\begin{proof}
   Wir zeigen, dass das Ideal $(1+i)$ nicht euklidisch ist. Wegen
   $\Phi_L(1+i) = 3$, $\lVert(1+i)\rVert = 4$ und weil $L$ keine
   Ideale der Norm 2 oder 3 enthält, genügt es zu zeigen, dass $u
   \equiv 1 \bmod {1+i}$ gilt. Dies folgt aber sofort aus
   $u-1 = \frac{(x-1) + (y-1)\sqrt{m}}{1+i}$, weil
   $\frac{(x-1) + (y-1)\sqrt{m}}{2}$ ganz ist.
\end{proof}

Mit diesem Kriterium lassen sich jedoch nur wenige $m$ ausschließen;
immerhin gilt

\begin{quote}
  \textbf{(7.10)} {\em Sei $m \equiv \pm 3 \bmod {(1+i)^6}$ und $m = x^2
    \pm 2$ für ein $x \in R$. Dann ist
    $\frac{2 \pm i + \sqrt{m}}{1-i} = (\frac{1 + \sqrt{m}}{2})^3$ für
    $m = 1 + 4i$, während für alle anderen $m$ die Fundamentaleinheit
    $u \equiv \frac{x + \sqrt{m}}{1 + i}$ ist.}
\end{quote}

Damit folgt z. B., dass $L$ für die folgenden Werte von $m$ nicht
normeuklidisch sein kann:

\begin{align*}
m &= -11 &&= (3i)^2 - 2 &&\quad (\text{sh. auch §\,6}) \\
m &= 7 + 12i &&= (3 + 2i)^2 + 2, &&\quad P = 193 \\
m &= 13 + 8i &&= (4 + i)^2 - 2, &&\quad P = 233 \\
m &= 5 + 24i &&= (4 + 3i)^2 - 2, &&\quad P = 601 \\
m &= 23 + 20i &&= (5 + 2i)^2 + 2, &&\quad P = 929.
\end{align*}

Ein etwas nützlicheres Kriterium ist

\begin{quote}
  \textbf{(7.11)} {\em Sei $m \equiv \pm 3 \bmod {(1+i)^6}$,
    $s \equiv t \equiv 1 \bmod {2}$, $(s,t) = 1$,
    $[ \frac{s}{m} ] = [ \frac{t}{m} ] = -1$ und
    $|st| \le (\sqrt{2}-1)|m|$, $|st+m| \ge |m|$. Ist dann
    $|st-m| \ge |m|$ oder $|st-m|$ keine Norm aus $S$, dann
    ist $L = K(\sqrt{m})$ nicht normeuklidisch.}
\end{quote}

\begin{proof}
  Wir verwenden (1.4) mit $e = s \cdot t$ und beachten, dass aus
  $r = e + l m$, $l \in R$,
  $|l m| = |r - e| \le |r| + |e| < |m| + (\sqrt{2}-1)|m| = \sqrt{2}|m|$,
  also $|l| < \sqrt{2}$ und damit $l \in \{0, \pm 1, \pm i\}$ folgt.
  Im Falle $l=0$ ist $r=e$ keine Norm; wäre $l=\pm i$, so folgte
  $r = e + l m \equiv 1 + i \bmod {2}$, und da $(1+i)$ in $L/K$
  träge ist, kann $r$ dann keine Norm aus $S$ sein. Damit bleiben nur
  noch die Möglichkeiten $l = \pm 1$, und die Voraussetzungen garantieren,
  dass dann $|r| \ge |m|$ oder $r$ keine Norm aus $S$ ist.
\end{proof}

\chapter*{\S\ 8 Sonstige Zahlkörper 4. Grades}
\setcounter{chapter}{8}
\addcontentsline{toc}{chapter}{\S\ 8 Sonstige Zahlkörper 4. Grades}
\markboth{Euklidische Ringe}{\S\ 8 Sonstige Zahlkörper 4. Grades}

Wir betrachten zuerst total komplexe Körper $K$ 4. Grades; diese
können wir klassifizieren nach der Galoisgruppe ihres normalen
Abschlusses $N$:

\begin{enumerate}
\item $\Gal(N/\Q) = \Gal(K/\Q) = Z_4$ (zyklische
  Gruppe der Ordnung 4);
\item $\Gal(N/\Q) = \Gal(K/\Q) = V_4$ (Kleinsche
  Vierergruppe);
\item $\Gal(N/\Q) = D_4$ (Diedergruppe der Ordnung 8);
\item $\Gal(N/\Q) = A_4$ (alternierende Gruppe der Ordnung
  12);
\item $\Gal(N/\Q) = S_4$ (symmetrische Gruppe der Ordnung 24).
\end{enumerate}

Hat $K$ einen quadratischen Teilkörper, so kommen nur die ersten drei
Möglichkeiten in Betracht; im 1. und 3. Fall hat $K$ genau einen, im
2. Fall aber drei solcher Teilkörper.

Im 1. Fall gibt es nach van der Linden genau zwei normeuklidische
Körper, nämlich den Körper $\Q(\zeta_5)$ der 5. Einheitswurzeln und
den Teilkörper 4. Grades von $\Q(\zeta_{13})$.

Im 2. Fall haben wir alle normeuklidischen Körper in § 6 bestimmt. Sei
also $\Gal(N/\Q) = D_4$ die Diedergruppe. Dann hat $K$ genau
einen quadratischen Teilkörper $k = \Q(\sqrt{m})$, und wir
unterscheiden:

\begin{enumerate}
  \item[3.1] $\Q(\sqrt{m})$ ist imaginärquadratisch. Damit $K$
    euklidisch sein kann, muss $K$ Klassenzahl 1 haben; dies
    impliziert, dass der quadratische Teilkörper $\Q(\sqrt{m})$
    Klassenzahl 1 oder 2 hat; dies ist eigentlich ein elementares
    Ergebnis (sh. z. B. Narkiewicz \cite{Nar74}), wird jedoch gewöhnlich
    mit Hilberts Satz 94 oder der Klassenkörpertheorie begründet. Hätte
    $k$ Klassenzahl 2, so wäre $K$ der Hilbertklassenkörper von $k$,
    also $K/\Q$ galoissch; dies ist aber nicht der Fall.
    
    Also hat $k = \Q(\sqrt{m})$ Klassenzahl 1, und wir dürfen uns auf
    die Möglichkeiten $m = -1, -2, -3, -7, -11, -43, -67, -163$
    beschränken. Ich erwarte nicht, dass sich bei der Untersuchung
    dieser $k$ andere Probleme auftun als in § 7 (wo $m = -1$ war);
    vielmehr wird man genau wie im Falle der Dirichlet'schen Körper
    vorgehen können. Die beiden schwierigsten Fälle sind naturgemäß
    $m=-1$ und $m=-3$ wegen der Existenz nichttrivialer Einheiten; die
    meisten der hier betrachteten normeuklidischen Körper werden wohl
    auch einen dieser beiden Körper enthalten.
\item[3.2] $\Q(\sqrt{m})$ ist reellquadratisch. Genau wie in
  3.1. findet man, dass auch hier $\Q(\sqrt{m})$ Klassenzahl 1 haben
  muss. Will man hier das Kriterium (1.4) anwenden, so stellt man
  recht schnell fest, dass sich die Anwesenheit von unendlich vielen
  Einheiten in $\Q(\sqrt{m})$ recht störend auswirkt. Dieser Nachteil
  lässt sich aber durch folgende Überlegungen wettmachen: ist nämlich
  $u > 1$ die FE von $k$, so ist entweder $u$ ebenfalls FE von $K$,
  oder es ist $K = k(\sqrt{-u})$, und $\sqrt{-u}$ ist FE von $K$ (denn
  $+1$ und $-1$ sind die beiden einzigen Einheitswurzeln in $K$).
  Damit in letzterem Falle $K$ totalkomplex wird, muss $u$
  total positiv sein (d.h. $u>0$ und $u'>0$), und dann ist
  $N_{K/\Q}(u) = 1$. Dies aber impliziert (sh. z.B. Kubota\index[dN]{Kubota}
  \cite{Kub56}) $\Gal(k(\sqrt{-u})/\Q) = V_4$.  Also ist $u$ auch
  FE von $K$, sodass wir bei vorgegebenem $k$ die Einheitengruppe von $K$
  schon kennen. Wie man das ausnützen kann, wollen wir am Beispiel
  $k = \Q(\sqrt{2})$ demonstrieren:

  Wir unterscheiden die möglichen Faktorisierungen von $(\sqrt{2})$:
  \begin{enumerate}
  \item[(a)] $(\sqrt{2}) = P^2$: dann ist $\Phi(P^2) = 2$,
    $u=1+\sqrt{2} \equiv 1 \pmod{P^2}$ und folglich jede Einheit
    $\equiv 1 \pmod{P^2}$ in $K$. Wäre $K$ normeuklidisch, müsste die
    prime Restklasse $1 + \alpha \bmod P^2$ (wo $\alpha$ ein
    Primelement für $P$ ist) ein Element der Norm $< \|P^2\| = 4$
    enthalten. Dieses Element kann keine Einheit sein, folglich müsste
    es Norm 3 haben. Das kann aber nicht sein, weil (3) in $k$ träge
    ist. Da auch (5) in $k$ träge ist, folgt $M(K) \ge \frac74$ für alle
    solchen $K$.
  \item[(b)] $(\sqrt{2}) = P$ ist träge: dann ist $\Phi(P) = 3$, und
    wie oben folgt, dass das Ideal $P$ (und damit auch $K$) nicht
    normeuklidisch ist.
  \item[(c)] $(\sqrt{2}) = 2_1 2_2$ ist zerlegt: nach dem
    Zerlegungsgesetz für relativ-quadratische Erweiterungen
    (sh. Hilbert) ist hier $K = k(\sqrt{\mu})$ für ein
    $\mu \equiv 1, 3+2\sqrt{2} \bmod 4$, d.h. mit $\mu = a+b\sqrt{2}$ wird
    $a-1 \equiv b \equiv 0 \bmod 2$. Sei nun $P = (p)$ ein Primideal in
    $R = \Z[\sqrt{2}]$ (d.h. $p \in R$), und $P$ möge in $K/k$
    verzweigen: $(p) = Q^2$. Hat $K$ eine ungerade
    Klassenzahl, so muss $Q$ Hauptideal sein, und es gibt ein
    $\pi \in S$ mit $(\pi^2) = (p)$. Da alle Einheiten in $k$ liegen,
    können wir solche in das $p$ hineinziehen und finden $\pi^2 =
    p$. In einem solchen Fall ist also $\mu = \pi$. Dies zeigt auch
    sofort, dass in $K/k$ höchstens ein Primideal verzweigen darf,
    damit $K$ noch ungerade Klassenzahl hat: in diesem Fall ist also
    $a^2 - 2b^2 = q$ eine rationale Primzahl.
  \end{enumerate}  
\end{enumerate}

Damit $K$ auch Klassenzahl 1 hat, müssen die Ideale $2_1$ und $2_2$
Hauptideale sein; dies bedeutet, dass die Gleichung
$4 \cdot p = \xi^2 - \mu \eta^2$ für ein $p \in R$ der Norm $2$
lösbar sein muss (in $R=\Z[\sqrt{2}]$; der Faktor $4$ taucht auf, weil
$(1+\sqrt{\mu})/2$ oder $(1+\sqrt{2} + \sqrt{\mu})/2$ ganz ist). Da
mit $-\mu$ die ganze rechte Seite der Gleichung totalpositiv ist,
können wir $p = 2 + \sqrt{2}$ annehmen (etwaige Faktoren $u^{2k}$
ziehen wir in $\xi$ und $\eta$ hinein) und haben
$2 + \sqrt{2} = \xi^2 - \mu \eta^2$. Jetzt seien durch
$|a+b\sqrt{2}|_1 := |a+b\sqrt{2}|$
und $|a+b\sqrt{2}|_2 := |a - b\sqrt{2}|$ die beiden archimedischen
Bewertungen von $k$ definiert. Dann ist, weil $\xi^2$ und $-\mu
\eta^2$ beide total positiv sind,
\begin{align*}
  4 \cdot (2 + \sqrt{2}) & = |2+\sqrt{2}|_1
    = |\xi^2 - \mu \eta^2|_1 \ge |\mu \eta^2|_1
     \quad \text{und entsprechend} \\
  4 \cdot (2 - \sqrt{2}) & = |2 - \sqrt{2}|_2
  = |\xi^2 - \mu \eta^2|_2 \ge |\mu \eta^2|_2.
\end{align*}
Multiplikation beider Ungleichungen ergibt wegen
$|\alpha_1| \cdot |\alpha_2| = |N_{k/\Q}(\alpha)|$ für $\alpha \in k$:
$32 \ge |N_{k/\Q}(\mu \eta^2)|$.  Wegen $\eta \ne 0$ ist sicher
$|N_{k/\Q}(\eta)| \ge 1$, und wir finden schließlich
$|N_{k/\Q}(\mu)| \le 32$.

Zusammen mit den bereits oben hergeleiteten Bedingungen an $\mu$ lässt
dies nur die Möglichkeit $\mu = \pm 5 + 2\sqrt{2}$ mit
$N_{k/\Q}(\mu) = 17$ und $\disc K = 8^2 \cdot 17 = 1088$ offen.
Wir haben damit gezeigt:

\begin{quote}
  \textbf{(8.1)} {\em Ist $K$ ein totalkomplexer Körper 4. Grades
    mit $\Q(\sqrt{2}) \subset K$, so sind höchstens die folgenden
    $K$ normeuklidisch:}
    $$ \Q(\sqrt{2}, \sqrt{-1}), \quad
       \Q(\sqrt{2}, \sqrt{-3}) \quad \text{{\em und}} \quad
       \Q\Big( \sqrt{-5 - 2\sqrt{2}}\Big). $$
\end{quote}

Entsprechende Rechnungen lassen sich auch für andere reellquadratische
Körper durchführen. Es sei noch bemerkt, dass obige Überlegungen durch
die von Cohn\index[dN]{Cohn} \cite{Coh58} berechneten Klassenzahlen
imaginärquadratischer Erweiterungen von $\Q(\sqrt{5})$ motiviert
wurden. Auch das Buch von Pólya\index[dN]{Polya} \cite{Pol69} sei an
dieser Stelle einmal ausdrücklich erwähnt.

Im Falle $\Gal(N/\Q) = A_4$ ist bis heute kein einziges
Beispiel eines normeuklidischen Körpers 4. Grades bekannt; dies liegt
wohl daran, dass es keine solchen Körper mit $\disc K <
3136 = 2^6 7^2$ gibt.

Die allermeisten totalkomplexen Körper 4. Grades schließlich haben die
Gruppe $S_4$ als Galoisgruppe ihres normalen Abschlusses. Hier wird
man wohl ohne eine Tafel aller solcher Körper (nebst Einheiten- und
Idealklassengruppe) wenig sagen können.

Über Körper 4. Grades mit $r=2$, $s=1$ oder $r=4$, $s=0$ lässt sich
noch weniger sagen. Vereinzelt habe ich mit einem Computer den EA
nachgewiesen; es werden jedoch eine Menge weiterer Rechnungen nötig
sein, bevor man geeignete Ansatzpunkte findet.

\section*{{\sc Anmerkungen zu} \S\ 8}
\addcontentsline{toc}{section}{{\sc Anmerkungen zu} \S\ 8}

\begin{tabular}{p{1cm} p{10cm}}
\toprule
\textbf{Jahr} & \textbf{Ereignis} \\
\midrule
1937 & Dribin\index[dN]{Dribin} bestimmt die Hilbertsche Untergruppenreihen für
  Zahlkörper mit der Galoisgruppe $S_4$. \\
1956 & Godwin\index[dN]{Godwin} gibt eine Tafel totalreeller Körper 4. Grades mit
  kleiner Diskriminante. \\
1957 & Godwin veröffentlicht entsprechende Tafeln für Zahlkörper
  4. Grades mit $r=2$ und $r=0$. \\
1958 & H. Cohn\index[dN]{Cohn} berechnet die Klassenzahl von
  imaginärquadratischen Erweiterungen von $\Q(\sqrt{5})$ (mit einer
  äußerst ungewöhnlichen Methode: bereits Gauß hat entdeckt, dass die
  Anzahl der Darstellungen eines primen $p \equiv 3 \pmod{8}$ als
  Summe dreier Quadrate eng mit der Klassenzahl von $\Q(\sqrt{-p})$
  zusammenhängt; eine entsprechende Formel mit $\Q(\sqrt{5})$ als
  Grundkörper hat sich Cohn zunutze gemacht). \\
1980 & Edgar\index[dN]{Edgar} und Peterson\index[dN]{Peterson}
  beschreiben zyklische Körper 4. Grades. \\
1985 & Godwin erweitert seine Tafel aus dem Jahre 1957 für Zahlkörper
  4. Grades mit $r=2$, $s=1$. \\
1986 & H. Cohn und J. Deutsch\index[dN]{ Deutsch} zeigen per Computer, dass
  $\Q\big(\sqrt{2 + \sqrt{2}}\big)$ und $\Q\big(\sqrt{3 + \sqrt{2}}\big)$
  normeuklidisch sind. \\
1989 & Buchmann\index[dN]{Buchmann} und Pohst\index[dN]{Pohst} geben alle
   totalreellen Zahlkörper 4. Grades mit $\disc K < 10^6$ nebst Einheiten und
  Klassenzahl. \\
\bottomrule
\end{tabular}

\chapter*{\S\ 9 Zahlkörper höheren Grades}
\setcounter{chapter}{9}
\addcontentsline{toc}{chapter}{\S\ 9 Zahlkörper höheren Grades}
\markboth{Euklidische Ringe}{\S\ 9 Zahlkörper höheren Grades}

Wir beginnen mit einem Zitat von H. Stark\index[dN]{Stark} \cite{Sta88}:

\begin{quote}
  {\em Heilbronn seems to suspect \dots that there may be infinitely
    many totally real Euclidean cubic fields. \dots he goes even
    further and names two families of quartic and sextic fields that
    should be investigated.}
\end{quote}

Dabei bezieht sich Stark auf die beiden folgenden Äußerungen von
Heilbronn:\index[dN]{Heilbronn}
\begin{itemize}
\item \cite{Hei50}: {\em Theorem 1: E.A. holds only in a finite number
  of cyclic cubic fields. The question of E.A. in real non-cyclic
  fields is thus left open. Though I cannot prove any result in the
  opposite direction I should be surprised to learn that the analogue
  of theorem 1 is true in that case.}
\item \cite{Hei51}: {\em Finally I should like to mention two types of
  cyclic fields for which E.A. may possibly hold in an infinity of
  cases.  
  \begin{enumerate}
  \item[(a)] The real quartic fields $\Q(\sqrt{\omega p})$,
    $\omega = \frac{5+\sqrt{5}}{2}$, of discriminant $125p^2$, where
    $p \equiv 3 \bmod 20$ is a prime;
  \item[(b)] The complex sextic field $\Q(\zeta_9 + \zeta_9^{-1},\sqrt{-p})$
    of discriminant $-3^8p^3$, where $p \equiv 3 \bmod 4$ is a prime.
  \end{enumerate} }
\end{itemize}

Der Teil (b) der Heilbronnschen Vermutung aus dem Jahre 1951 erscheint
mir recht mysteriös: $K = \Q(\zeta_9 + \zeta_9^{-1}, \sqrt{-p})$ ist
eine kubische Erweiterung von $k=\Q(\sqrt{-p})$, und zwar ist $K/k$
über Primidealen über (3) verzweigt. Also liegt $K$ nicht im
Hilbertklassenkörper von $k$, und es gilt $h(k)\mid h(K)$. Nun hat aber
Heilbronn selbst bewiesen, dass es nur endlich viele $k$ mit
Klassenzahl 1 gibt; die obige Teilbarkeitsrelation impliziert dann,
dass es auch nur endlich viele $K$ mit Klassenzahl 1 gibt. Fast noch
erstaunlicher ist, dass dies auch Stark übersehen hat (der ja
bekanntlich alle imaginär-quadratischen Zahlkörper mit Klassenzahl 1
gefunden hat).

Tatsächlich ist es nicht schwer zu zeigen, dass $K$ nur für $p=3$ oder
$p=11$ normeuklidisch sein kann (für $p=3$ ist $K=\Q(\zeta_9)$, $K$
also normeuklidisch). Ist nämlich $p>11$, und hat $k = \Q(\sqrt{-p})$
Klassenzahl 1, so muss $p \equiv 1 \bmod 3$ sein (andernfalls müsste
$k$ Elemente der Norm 3 besitzen, was nicht der Fall ist). Also ist
$(1 - \sqrt{-p})^2 = 1 - p - 2\sqrt{-p} \equiv \sqrt{-p} \bmod 3$; das
Element minimaler Norm in der Restklasse $\sqrt{-p} \bmod 3$ ist aber
$\frac{3 - \sqrt{-p}}2$ mit der Norm $\frac{p+9}4$, und für $p>43$ ist
diese Norm $>9$. Nach (1.4) kann $K$ dann nicht normeuklidisch
sein. Im Falle $p=19$ ist $\frac{3 - \sqrt{-p}}2$ keine Norm aus
$O_K$, weil (7) in $\Q(\zeta_9 + \zeta_9^{-1})$ träge ist; also ist
$K$ für $p>11$ nicht normeuklidisch. Ist schließlich $p=7$, so sind
sowohl $\sqrt{-p}$, als auch $\frac{3 - \sqrt{-p}}2$ keine Normen aus
$O_K$, womit wir die obige Behauptung bewiesen haben.

Weitere Klassen von Körpern, die man mit einfachen Mitteln untersuchen
kann, sind z. B. $\Q(\rho, \sqrt[3]{m}\,)$ oder $\Q(i, \sqrt[4]{m})$ für
$m \in \mathbb{Z}$. Hier lassen sich ohne weiteres neue Ergebnisse
erzielen. Ähnliches gilt für Zahlkörper der Form $\Q(\sqrt[n]{m})$ mit
z. B. $n=3, 6, 8$ usw., oder überhaupt für $\Q(\sqrt[n]{m}, \sqrt{-l})$,
$l \in \N$.

Verhältnismäßig wenig dagegen weiß man über den euklidischen
Algorithmus in Kreisteilungskörpern. Es ist zwar inzwischen bekannt,
dass $\Q(\zeta_m)$ für die Werte $m = 1, 3, 4, 5, 7, 8, 9, 11, 12, 13,
15, 16, 20, 24$ normeuklidisch ist, und dass $\Q(\zeta_m)$ genau für
\begin{align*}
  m & = 1, 3, 4, 5, 7, 8, 9, 11, 12, 13, 15, 16, 17, 19, 20, 21, 24, 25, 27,
  28, 32, 33, \\
  & \quad 35, 36, 40, 44, 45, 48, 60, 84
\end{align*}
weiß man seit den Untersuchungen von K. Uchida, H.L. Montgomery und
J.M. Masley (sh. z. B. Masley u. Montgomery \cite{MM76}). Außerdem kann man
ohne große Mühe zeigen, dass die Restklasse $1+\lambda^5 \bmod \lambda^6$
keine Einheiten enthält, wo $\lambda = 1+\zeta_{32}$ ein Element der
Norm 2 in $\Q(\zeta_{32})$ ist (sh. dazu Lenstra \cite[S. 14]{Len79c});
da es in $\mathbb{Z}[\zeta_{32}]$ keine Elemente der Norm 33
oder 65 gibt, enthält obige Restklasse nur Elemente der Norm $\ge 97$,
und es folgt $M(K) \ge \frac{97}{64}$ (man beachte
$N_{K/\Q}(\alpha) \equiv 1 \bmod 32$ für $\alpha \in \mathbb{Z}[\zeta_{32}]$;
weitere für einen vollständigen Beweis nötige Einzelheiten,
z. B. Information über die Einheiten in $\mathbb{Z}[\zeta_{32}]$,
findet man bei Washington \cite{Was82}).

\chapter*{\S\ 10 Offene Fragen}
\setcounter{chapter}{10}
\addcontentsline{toc}{chapter}{\S\ 10 Offene Fragen}
\markboth{Euklidische Ringe}{\S\ 10 Offene Fragen}

\begin{quote}
  {\em There are a lot of things I don't understand myself} (Linus van Pelt)
\end{quote}

Die folgenden Fragen sind durchnumeriert; hierbei bezieht sich die
erste Ziffer auf den Paragraphen, dem ich die jeweilige Frage
zugeordnet habe. Fragen, die ich mit den in den §§ 1--9 vorgestellten
Methoden für lösbar halte, sind als Aufgaben formuliert. Schließlich
sind besonders wichtige Fragen (bzw. solche, die ich für sehr wichtig
halte) fett gedruckt.

\subsection*{1.}

\begin{enumerate}
\item[1.1.] Die in (1.10) und (1.12) gefundenen Schranken hängen sehr
  von der Wahl der GHB ab; es stellt sich also die Frage, ob, und wenn
  ja, wie man die GHB so wählen kann, dass die Schranken $\mu_i$
  möglichst klein werden.
\item[1.2.] Sind die vor (1.17) definierten Lenstra-Konstanten $\mu_i$
  und die dazugehörigen Schranken $\lambda_i$ für $n \ge 3$ noch
  endlich (falls $R$ kein Hauptidealring ist, gilt dies sicherlich;
  dafür sorgt die Bedingung (F-2), S.~\pageref{pF2})).
\item[1.3.] In (1.21) wird für $L=Kk$ eine Aussage bewiesen, wobei
  $k/\Q$ abelsch ist. Gilt u. U. etwas Ähnliches, wenn nur
  verlangt wird, dass $k/\Q(\sqrt{-m})$ abelsch ist (d.h. kann
  man die $m$-ten Einheitswurzeln durch elliptische Einheiten
  ersetzen?)
\item[1.4.] Kann man in (1.24) auf die Voraussetzung ``$n$ ist
  brauchbar'' verzichten, oder lässt sie sich zumindest durch eine
  schwächere ersetzen?
\item[1.5.] Bestimme $c(K)$ auch für andere Klassen von Kreiskörpern
  als für volle Kreisteilungskörper mit Primpotenzdiskriminante.
\item[1.6.] Bestimme $c(K)$ für $K = \Q(\sqrt{a+b\sqrt{-3}})$,
  insbesondere im Fall $a=-1$, $b=2$; hier wissen wir, dass $c(K) \le
  (4 + \sqrt{13})/12 < 0.6338$ ist. Wenn wir $c(K) \le \frac12$ zeigen
  könnten, würde folgen, dass der Körper 8. Grades $L = K(i)$
  normeuklidisch ist. Allerdings lassen jüngste Rechnungen meinerseits
  vermuten, dass hier $c(K) > \frac12$ ist.
\end{enumerate}

\subsection*{2.}

\begin{enumerate}
\item[2.1.] Bestimme die fehlenden zweiten Minima in den Tafeln auf
  den Seiten 47 und 48.
\item[2.2.] Kann man in den Sätzen (2.5) bis (2.9) mit der
  vorgestellten Methode auch die höheren Minima ($M_3$, $M_4$, usw.)
  bestimmen?
\item[2.3.] \textbf{Gilt (2.12) auch in Körpern mit Einheitenrang $\ge 2$?}
\item[2.4.] Sei $C(K)$ die Menge aller $x \in K$ mit $M(x) =
  M(K)$. Sind dann die folgenden Aussagen richtig?
  \begin{enumerate}
  \item[(a)] Die Häufungspunkte von $C(K)$ liegen in $K$;
  \item[(b)] $C(K)$ enthält mindestens ein $x$ aus $K$.
  \end{enumerate}
\end{enumerate}

\subsection*{3.}

\begin{enumerate}
\item[3.1.] Finde einfache Beweise dafür, dass $D(m)$ für die Werte
  $$ m = 193, 241, 337, 457, 601 $$
  nicht normeuklidisch ist.
\item[3.2.] Zeige, dass es für alle primen $p \equiv 1 \pmod{24}$ mit
  $p > 601$ natürliche Zahlen $r, s, t, u$ gibt mit $(r,s) = (t,u) = 1$
  und $(\frac rp) = (\frac sp) = (\frac tp) = -1$.
\item[3.3.] Zeige, dass $D(10)$ und $D(65)$ die einzigen
  semieuklidischen quadratischen Ringe mit von 1 verschiedener
  Klassenzahl sind.
\item[3.4.] \textbf{Ist der Ring $R = D(14)$ euklidisch?} Dass $R$
  nicht normeuklidisch ist, haben wir bereits gesehen. Die Frage, ob
  es auf $R$ andere Funktionen gibt, bezüglich derer $R$ euklidisch
  ist, hat in diesem speziellen Fall erstmals Samuel\index[dN]{Samuel}
  \cite{Sam71} explizit gestellt. Lenstra\index[dN]{Lenstra} hat dann
  in \cite{Len74} vorgeschlagen, in dieser Hinsicht ``gewichtete
  Normen'' zu betrachten: dazu wähle man ein Ideal $P$ in $R$, setze
  $N(P)=c$ für ein positives reelles $c$, sowie $N(Q) = || Q ||$ für
  alle anderen Primideale $Q$. Dann setze man $N$ multiplikativ auf
  alle Ideale in $R$ fort und definiere $N(a) := N(aR)$ für alle $a \in K$.
  $N$ heißt die durch $N(P)=c$ gewichtete Norm.  Wir wollen nun statt
  Abbildungen $f:R \to \mathbb{N}$ auch Abbildungen $f:R \to
  \mathbb{R}$ als euklidische Funktionen zulassen, wenn es zu jedem
  $a \in R$ nur endlich viele Funktionswerte
  $<f(a)$ gibt (damit kann man die Menge $f(R)$ nach
  wachsenden Funktionswerten abzählen und erhält so eine Bijektion $g:
  f(R) \to \mathbb{N}$; damit ist dann $g \circ f: R \to \mathbb{N}$
  eine mögliche euklidische Funktion im bisherigen Sinne).  Wir
  betrachten nun das Ideal $P = (4 + \sqrt{14})$ der Norm 2 in $R$ und
  fragen uns, ob man $N(P) = c$ so wählen kann, dass die dadurch
  gewichtete Norm eine auf $R$ euklidische Funktion wird. Damit das
  Ideal $I = (7 + 2\sqrt{14})$ der Norm 7 euklidisch bezüglich $N$
  wird, muss jede prime Restklasse mod $I$ ein Element $a$
  mit $N(a) < 7 = N(7)$ enthalten. Da $u = 15 + 4\sqrt{14}
  \equiv 1 \mod I$ ist, enthalten nur die Restklassen $\pm 1 \mod I$
  Einheiten. Weiter haben die Elemente $\pm 3 + \sqrt{14}$ die Norm 5
  und sind $\equiv 3 \mod I$. Wegen $P = (4 + \sqrt{14})$ und $4 +
  \sqrt{14} \equiv -3 \mod I$ enthält auch $P$ kein Element $\equiv 2
  \mod I$ der Norm $<7$. Da die einzigen Ideale der Norm $<7$ die
  Ideale $R$, $(\pm 3 + \sqrt{14})$ und möglicherweise Potenzen von
  $P$ sind, gilt $N(a) > N(P)^2$ für jedes $a \equiv 2 \mod I$. Damit
  also $N(a) < 7$ wird, müssen wir $c^2 = N(P)^2 < 7$ wählen.

  Jetzt betrachten wir $I = (2) = P^2$; soll die Restklasse $1 +
  \sqrt{14} \mod I$ ein Element mit Norm $< c^2$ enthalten, so müssen
  wir $\delta < c^2$ machen, weil $J = R$ das einzige Ideal
  mit $I + J = R$ und $N(J) < 5$ ist.

  Also kann die durch $N(4 + \sqrt{14}) = c$ gewichtete Norm höchstens
  dann eine euklidische Funktion auf $R$ sein, wenn $\sqrt{5} < c < \sqrt{7}$
  gilt. Die einfachste Möglichkeit, $c$ zu wählen, ist daher
  wohl $c = \sqrt{6}$. Tatsächlich hat Bedocchi 1985 auf einem ganz
  anderen Weg ebenfalls diese Funktion gefunden und vermutet, dass sie
  euklidisch auf $R$ ist.
\end{enumerate}

\subsection*{4.}

\begin{enumerate}
\item[4.1.] In quadratischen Zahlringen $K$ mit $d = \disc K$ gilt
  bekanntlich
  $$ \frac{\sqrt{d}}{16 + 6\sqrt{6}} \le M(K) \le \frac{\sqrt{d}}4. $$
  In kubischen Zahlkörpern mit Einheitenrang 1 dagegen weiß man
  $$ \frac{\sqrt{d}}{420} \le M(K) \le \frac{|d|^{2/3}}{16^{3/2}}, $$
  wobei Exponent $2/3$ und Faktor $1/16^{3/2}$ bestmöglich sind
  (allerdings ist diese letzte Ungleichung nur für
  $\disc K < -1236$ bewiesen). In Analogie zum
  quadratischen Fall könnte man nun vermuten, dass sich die untere
  Schranke $\sqrt{d}/420$ verbessern lässt zu $k_1 |d|^{2/3}$ für ein
  reelles $k_1$. Van der Linden schreibt bei seinem Beweis der
  ``Casselschen Schranke'' $\sqrt{d}/420$, dass zwischen dem Beweis im
  kubischen und im quadratischen Fall eine gewisse Asymmetrie liegt,
  und er vermutet, dass dies daran liegt, dass einige seiner
  Abschätzungen im kubischen Fall nicht bestmöglich sind.
  Falls sich der Exponent $\frac12$ in der Casselschen Schranke
  tatsächlich auf $\frac23$ verbessern lässt, so hätte dies wohl auch
  Auswirkungen auf den totalreellen kubischen Fall (bzw. auf
  Zahlkörper beliebigen Grades und beliebiger Signatur überhaupt):
  hier weiß man ja, dass $M(K) \le \sqrt{d}/8$ gilt und dass diese
  Schranke bestmöglich ist, sodass man auch in diesem Fall vermuten
  würde, dass es ein reelles $k_2$ mit $k_2 \sqrt{d} \le M(K)$ gibt.
  Schließlich würde man so zu der Vermutung geführt, dass es nur
  endlich viele normeuklidische Zahlkörper mit gegebenem Körpergrad
  $n$ gibt.
\item[4.2.] Zeige (z. B. mit einem Computer), dass die Ungleichung $M(K)
  \le |d|^{2/3}/16^{3/2}$ auch für die Körper mit $-300 \ge
  \disc K \ge -1236$ gilt. Für die Körper mit $-23 \ge
  \disc K \ge -199$ folgt dies aus den Werten für
  $M(K)$, die wir in $\S 4$ bestimmt haben. Auch für die euklidischen
  Körper mit $-200 \ge \disc K \ge -1236$ folgt die
  Ungleichung sofort.
\item[4.3.] Bestimme $M(K)$ für den kubischen Körper mit
  $\disc K = -87$.
\item[4.4.] Sind die kubischen zyklischen Körper mit $\sqrt{d} = 103,
  109, 127$ und $157$ normeuklidisch?
\item[4.5.] Gibt es kubische zyklische Körper mit EA und $\sqrt{d} \ge
  5 \cdot 10^5$?
\end{enumerate}

\subsection*{5.}

\begin{enumerate}
\item[5.1.] Finde eine beste obere Schranke für $M(K)$, wo $K$ total
  komplexer Zahlkörper 4. Grades ist. Nach Davenport und
  Swinnerton-Dyer ist $M(K) \le c \cdot d^{3/4}$ bestmöglich; jedoch
  haben sie hierfür keinen Beweis gegeben.
\item[5.2.] Gilt für beliebige Dirichlet'sche Körper die Abschätzung
  $M(K) \le \sqrt{d}/16$?
\item[5.3.] Gilt für quadratische Erweiterungen imaginärquadratischer
  Zahlkörper $M(K) \le c \sqrt{d}$?
\item[5.4.] Finde eine Klasse total komplexer biquadratischer
  Zahlkörper, für die $M(K) \le c \cdot d^{3/4}$ bestmöglich ist.
\item[5.5.] Beende die Klassifikation der reinen biquadratischen
  Zahlkörper mit EA.
\item[5.6.] Gibt es nur endlich viele reine Zahlkörper von
  Zweierpotenzgrad?
\end{enumerate}

\subsection*{6.}

\begin{enumerate}
\item[6.1.] Bestimme $M(K)$ für $D(-1,13)$, $D(-1,17)$, $D(-3,2)$,
  $D(-3,-11)$ und $D(-3,-19)$.
\item[6.2.] Bestimme $M(K)$ für weitere Klassen bizyklischer
  Zahlkörper, beispielsweise für $D(q,m)$, $q = -1, -2, -3$, $m = n^2
  \pm r$, $r|4n$.
\end{enumerate}

\subsection*{7.}

\begin{enumerate}
\item[7.1.] Beende die Klassifikation der normeuklidischen
  Dirichlet'schen Zahlkörper.
\item[7.2.] Bestimme $M(K)$ für weitere Klassen Dirichlet'scher
  Zahlkörper.
\end{enumerate}

\subsection*{8.}

\begin{enumerate}
\item[8.1.] Finde alle normeuklidischen Zahlkörper 4. Grades, die
  einen quadratischen Teilkörper besitzen.
\item[8.2.] Zeige, dass es nur endlich viele normeuklidische Zahlkörper
  $K$ gibt mit $(K:\mathbb{Q}) = 4$, $K/\mathbb{Q}$ zyklisch,
  $\disc K = 125p^3$ (dies würde eine Vermutung von
  Heilbronn widerlegen).
\end{enumerate}

\subsection*{9.}

\begin{enumerate}
\item[9.1.] Sind die Körper der $p$-ten Einheitswurzeln ($p = 17, 19$)
  normeuklidisch?
\end{enumerate}

\chapter*{\S\ 11  Beschreibung der Computerprogramme}
\setcounter{chapter}{11}
\addcontentsline{toc}{chapter}{\S\ 11  Beschreibung der Computerprogramme}
\markboth{Euklidische Ringe}{\S\ 11  Beschreibung der Computerprogramme}

Gegeben sei ein algebraischer Zahlkörper $K$ vom Grad $n$ mit einer
Ordnung $R$, die über $\mathbb{Z}$ von den Elementen $\alpha_1, \dots,
\alpha_n$ erzeugt wird (ist $R$ gleich der Hauptordnung, so bildet
$\{\alpha_1, \dots, \alpha_n\}$ eine GHB). Weiter seien unabhängige
Einheiten $u_1, \dots, u_l$ gegeben, sowie ein $k \in R$, und gefragt
sei nach den Ausnahmemengen bezüglich $k$.

Dazu geben wir uns durch 
\[
S = \{\alpha \in R : \alpha = \sum a_i \alpha_i \, , \, -a \leq a_i \leq a \}
\]
eine Teilmenge $S$ von $R$ vor (z.B. mit $a=10$, $10^2$, \dots). Dann
unterteilen wir den Fundamentalbereich $F = (-\frac{1}{2},
\frac{1}{2}) \times \dots \times (-\frac{1}{2}, \frac{1}{2})$ (wir
übernehmen die Notation von § 2) in kleine $n$-dimensionale Würfel der
Kantenlänge 0.1 (oder auch 0.25). Für jeden solchen Würfel $W$ gehen
wir nun vor wie folgt:

\begin{enumerate}
\item Prüfe, ob für alle $x \in W$ ein $y \in S$ existiert mit $N(x-y)
  \leq k$; falls ja, so prüfe den nächsten Würfel, andernfalls gehe
  nach 2.
\item Unterteile $W$ in $2^n$ Teilwürfel durch Halbierung der
  Seitenlänge. Ist die Seitenlänge nun kleiner als eine vorgegebene
  Schranke (z.B. $10^{-2}$, $10^{-3}$, usw.), so rechnen wir $W$ zur
  Ausnahmemenge hinzu und gehen mit dem nächsten Würfel zurück nach
  1. (beim letzten Würfel gehen wir nach 3.); andernfalls gehen wir
  mit jedem einzelnen Teilwürfel zurück nach 1.  
\item Zu diesem Zeitpunkt haben wir alle Teilwürfel $W$ gefunden, die
  möglicherweise einen Ausnahmepunkt enthalten (falls kein Würfel
  übrig geblieben ist, gilt $M(K) < k$). Wir prüfen nun, für welche
  Würfel $W$ die Transformierten $u_i W \bmod R$ in bedecktem Gebiet
  liegen (indem wir testen, ob $u_i W$ mit einem der andern $W$ einen
  Punkt gemeinsam hat). Jeden Würfel $W$, für den ein $u_i W$ ganz in
  bedecktem Gebiet liegt, können wir vergessen.
    
  Nun ist es denkbar, daß ein $u_i W_1 \bmod R$ einen Würfel $W_2$
  schneidet, welcher zu einem späteren Zeitpunkt ``vergessen'' wird,
  weil $u_j W_2 \bmod R$ in bedecktem Gebiet liegt. Wir wiederholen
  den Schritt 3. daher so lange, bis sich keine Änderung mehr ergibt.
\end{enumerate}

Jetzt prüfen wir (von Hand), ob sich die Sätze aus § 2 anwenden
lassen. Ist dies nicht der Fall, so gibt es folgende Möglichkeiten:

\begin{enumerate}
\item[a)] Verkleinern der Schranke für die Seitenlänge der Teilwürfel in 2.
\item[b)] Vergrößern von $S$.
\item[c)] Falls die $u_i$ kein Grundeinheitensystem bilden, kann man
  versuchen, ein solches zu finden, die $u_i$ dadurch zu ersetzen und
  Schritt 3. zu wiederholen.
\item[d)] Vergrößern von $k$ zu $k_1$; es ist ja durchaus denkbar,
  daß man im ersten Versuch ein viel zu kleines $k$ gewählt hat.
\end{enumerate}

Wir müssen jetzt noch zeigen, wie man die Norm auf einem gegebenen
Teilwürfel $W$ nach oben abschätzt. Sei dazu $W$ (oder auch ein um ein
$y \in S$ verschobenes $W$) durch die Ungleichungen $r_1 \leq x_1 \leq
s_1, \dots, r_n \leq x_n \leq s_n$ gegeben (hier ist $-\frac{1}{2}
\leq r_i < s_i \leq \frac{1}{2}$). Wir müssen dann eine obere Schranke
für $N(\sum x_i \alpha_i)$ finden. Bekanntlich ist aber $N(z) = |z_1|
\cdot \dots \cdot |z_r| \cdot \dots \cdot |z_r+s|$ (für die Definition
der $| \cdot |_j$ sh. § 2), sodaß es genügt, jeden Faktor $|z_j|$ nach
oben abzuschätzen.

Ist $| \cdot |_j$ durch $|z_j| = |\sum x_i \sigma_j(\alpha_i)|$ definiert und $\sigma_j$ eine reelle Einbettung, so gilt 
\[
\sup \{|z_j| : z \in W\} = \max \{ |\sum x_i \sigma_j(\alpha_i)| : x_i \in \{r_i, s_i\} \},
\]
d.h. man kann für $|z_j|$ mit $2^n$ Auswertungen eine obere Schranke
finden, wenn man nur die reellen Zahlen $\sigma_j(\alpha_i)$
(hinreichend genau) kennt. Ist $\sigma_j$ dagegen eine nicht-reelle
Einbettung, so kann man ganz analog den Real- und Imaginärteil von
$\sum x_i \sigma_j(\alpha_i)$ abschätzen, oder man berechnet die
reellen Koeffizienten von $x_i x_k$ im Produkt $\sum x_i
\sigma_j(\alpha_i) \sum x_i \overline{\sigma_j(\alpha_i)}$ und schätzt
dieses Produkt nach oben ab.

Als ein Beispiel geben wir an, wie dieses Verfahren im Falle $n=3$,
$r=s=1$, aussieht; dazu nehmen wir an, $K$ sei durch ein Polynom $f$
mit Index $1$ gegeben (ist der Index $\neq 1$, so muß man das Folgende
entsprechend abändern): Input für die Funktion NORMAX, die den
Normbetrag auf einem Würfel $W$ nach oben abschätzt, sind hier $r_1,
s_1, t_1$, die ``Schrittweite'' $sw$, die reelle Nullstelle $x$ von
$f$, sowie die komplexen Nullstellen $w_1 \pm i w_2$. Wir bilden dann
$r_2 = r_1 + sw$, $s_2 = s_1 + sw$, $t_2 = t_1 + sw$ und $y = w_1^2 -
w_2^2$ (da $1, x, x^2$ eine Ganzheitsbasis bilden, ist $\sigma(x^2) =
w_1^2 - w_2^2 + 2 w_1 w_2 i$). Bezeichnet schließlich $x_1$ die obere
Schranke für $|z_1|$, $y_2$ die obere Schranke für den Realteil, $y_3$
diejenige für den Imaginärteil von $z'$ (bzw. $z''$), so haben wir in
GFA-Basic:

\begin{verbatim}
Function NORMAX(r1,s1,t1,sw)
    r2=r1+sw
    s2=s1+sw
    t2=t1+sw
    x2=x^2
    y=w1^2 - w2^2
    x1=MAX(ABS(r1*s1*x+t1*x2),ABS(r1*s1*x+t2*x2),ABS(r1*s2*x+t1*x2))
    x1=MAX(ABS(r1*s2*x+t2*x2),ABS(r2*s1*x+t1*x2),ABS(r2*s1*x+t2*x2),x1)
    x1=MAX(ABS(r2*s2*x+t1*x2),ABS(r2*s2*x+t2*x2), x1)
    y2=MAX(ABS(r1*s1*w1+t1*y),ABS(r1*s1*w1+t2*y),ABS(r1*s2*w1+t1*y))
    y2=MAX(ABS(r1*s2*w1+t2*y),ABS(r2*s1*w1+t1*y),ABS(r2*s1*w1+t2*y),y2)
    y2=MAX(ABS(r2*s2*w1+t1*y),ABS(r2*s2*w1+t2*y), y2)
    y3=MAX(ABS(s1+2*t1*w1),ABS(s1+2*t2*w1),ABS(s2+2*t1*w1),ABS(s2+2*t2*w2))
    y3=y3*abs(w2)
    y1=y2^2 + y3^2
    RETURN x1*y1
ENDFUNC
\end{verbatim}

$y_1 = y_2^2 + y_3^2$ ist also obere Schranke für $|z_2|$, und $x_1
y_1$ eine solche für den Betrag der Norm auf $W$.

In der Praxis läuft das Programm so ab: nach Eingabe eines erzeugenden
Polynoms (mit möglichst kleinen Koeffizienten; dies ist für die
Fehlerabschätzung von Bedeutung) berechnet man mit Hilfe des
Newton-Verfahrens die reelle Nullstelle $x$ von $f$ mit einer
Genauigkeit von mindestens 10 Dezimalstellen, und erhält nach
Abdivision dieser Nullstelle ein quadratisches Polynom mit negativer
Diskriminante, aus dem man die Werte $w_1$ und $w_2$ leicht bestimmen
kann. Die Menge $S$ beschreiben wir durch drei Arrays $a(n)$, $b(n)$
und $c(n)$, die die (ganzen) Koordinaten bezüglich $1$, $x$ und $x^2$
bezeichnen; anfangs starten wir mit nur einem $y \in S$, nämlich
$y=0$, d.h. wir setzen $n=0$, $a(0)=b(0)=c(0)=0$. Dann wählen wir
z.B. die Schrittweiten $sw_1=0.1$ und $sw_2=0.05$ und teilen $F=(0,
0.5)\times(-0.5, 0.5)\times(-0.5, 0.5)$ in kleine Teilwürfel der
Kantenlänge $sw_1$ ein. Dann schätzen wir die Norm auf $y+W$ für jeden
Teilwürfel $W$ und alle $y \in S$ mit NORMAX ab; ist die Abschätzung
größer als die eingegebene Schranke, so halbieren wir die Kantenlängen
von $W$ und wiederholen den Schritt für jeden der acht Teilwürfel $V$
von $W$. Ist die Schranke für einen dieser Teilwürfel für jedes $y \in
S$ größer als die Schranke, so versuchen wir, ein $y \in R$ zu finden,
sodaß die Norm auf $y+V$ kleiner als die Schranke wird. Dies machen
wir mit der Prozedur SEARCH, indem wir jedes $y \in R$ ausprobieren,
dessen Koordinaten einen Betrag $\leq 8$ haben:

\begin{verbatim}
PROCEDURE SEARCH(r1,s1,t1,sw)
    FOR e = -8 TO 8
        FOR f = -8 TO 8
            FOR g = -8 TO 8
                norm=NORMAX(r1+e,s1+f,t1+g,sw)
                IF norm<schranke
                    n=n+1
                    a(n)=e
                    b(n)=f
                    c(n)=g
                    e=9
                    f=9
                    g=9
                ENDIF
            NEXT g
        NEXT f
    NEXT e
    RETURN
ENDPROC
\end{verbatim}

\section*{Das Hauptprogramm}

Das Hauptprogramm sieht damit wie folgt aus:
\lstset{
  basicstyle=\ttfamily\small,
  breaklines=true,
  columns=fullflexible
}
\begin{lstlisting}
FOR a1=0 TO 0.49 STEP sw1 
FOR b1=0 TO 0.49 STEP sw1 
FOR c1=0 TO 0.49 STEP sw1 
index=-1 
REPEAT 
  index=index+1 
  norm1=NORMAX(a1+a(index),b1+b(index),c1+c(index),sw1) 
UNTIL index=n or norm1<schranke 
IF norm1>schranke 
  FOR a2=a1 TO a1+0.049 STEP sw2 
  FOR b2=b1 TO b1+0.049 STEP sw2 
  FOR c2=c1 TO c1+0.049 STEP sw2 
  index=-1 
  REPEAT 
    index=index+1 
    norm2=NORMAX(a2+a(index),b2+b(index),c2+c(index),sw2) 
  UNTIL index=n or norm2<schranke 
  IF norm2>schranke 
    SEARCH(a2,b2,c2,sw2) 
    IF NORMAX(a2+a(n),b2+b(n),c2+c(n),sw2)>schranke 
      SCHREIB(a2,b2,c2,sw2) 
    ENDIF 
  ENDIF 
  NEXT c2 
  NEXT b2 
  NEXT a2 
ENDIF 
NEXT c1 
NEXT b1 
NEXT a1
\end{lstlisting}

Hierbei ist \textbf{SCHREIB} eine Prozedur, die den betrachteten
Teilwürfel $V$ nebst den transformierten Würfeln $uV$ für gegebene
Einheiten $u$ (z.B. die FE und ihr Inverses) auf Bildschirm, Drucker
oder Diskette ausgibt. Wie man dieses Programm zu ergänzen hat, wenn
man die Teilwürfel weiter unterteilen will, ist klar.

\chapter*{Tabellen}
\setcounter{chapter}{12}
\addcontentsline{toc}{chapter}{Tabellen}
\markboth{Euklidische Ringe}{Tabellen}

\section*{Kubische Zahlkörper mit negativer Diskriminante}
\addcontentsline{toc}{section}
                     {Kubische Zahlkörper mit negativer Diskriminante}
Es bedeuten:

\begin{enumerate}
\item[] $\disc K$: die Diskriminante des Zahlkörpers
\item[] E: $K$ ist normeuklidisch
\item[] NE: $K$ ist nicht normeuklidisch
\item[] H: $K$ hat von $1$ verschiedene Klassenzahl
\item[] $a_2, a_1, a_0$: die Koeffizienten eines erzeugenden Polynoms
  \[ f(x) = x^3 + a_2 x^2 + a_1 x + a_0 \]
\item[] GHB: gibt eine Ganzheitsbasis, falls nicht
  $\{1, \vartheta, \vartheta^2\}$ eine solche ist,
  wo $\vartheta$ Wurzel von $f$ ist
\item[] $h$: die Klassenzahl des Körpers
\item[] $u$: eine Einheit
\end{enumerate}

$$ 
\begin{array}{rrrrrrrrcc}
\toprule
Nr. & \disc K & & a_2 & a_1 & a_0 & GHB & h & u & M(K) \\
\midrule
1. & -23 & E & 0 & -1 & & 1 & 1 & \vartheta & \frac{1}{5} \\
2. & -31 & E & 0 & 1 & & 1 & 1 & \vartheta & \frac{1}{3} \\
3. & -44 & E & 1 & -1 & & 1 & 1 & \vartheta & \frac{1}{2} \\
4. & -59 & E & 0 & 2 & & 1 & 1 & \vartheta & \frac{1}{2} \\
5. & -76 & E & 3 & 1 & & 1 & 1 & \vartheta & \frac{1}{2} \\
6. & -83 & E & 2 & -2 & & 1 & 1 & \vartheta & \frac{1}{2} \\
7. & -87 & E & 2 & -1 & & 1 & 1 & \vartheta & \geq \frac{1}{3} \\
8. & -104 & E & 0 & -1 & & 2 & 1 & 3+2\vartheta & \frac{1}{2} \\
9. & -107 & E & 4 & 2 & & 1 & 1 & \vartheta & \frac{1}{2} \\
10. & -108 & E & 3 & -3 & & 1 & 1 & \vartheta & \frac{1}{2} \\
\midrule
11. & -116 & E & 1 & 0 & & 2 & 1 & 1 - \vartheta + \vartheta^2 & \frac{1}{2} \\
12. & -135 & E & 3 & 0 & & 1 & 1 & \vartheta & \frac{3}{5} \\
13. & -139 & E & 4 & 6 & & 1 & 1 & \vartheta & \frac{1}{2} \\
14. & -140 & E & 5 & 3 & & 1 & 1 & \vartheta & \frac{1}{2} \\
15. & -152 & E & 1 & -2 & & 2 & 1 & 3 - \vartheta - \vartheta^2 & \frac{1}{2} \\
16. & -172 & E & 2 & 0 & & 2 & 1 & 1 + \vartheta^2 & \frac{3}{4} \\
17. & -175 & E & -2 & 3 & & 1 & 1 & \vartheta & \frac{3}{5} \\
18. & -199 & NE& 4 & 1 & & 1 & 1 & \vartheta & 1 \\
19. & -200 & E & 1 & 2 & & -2& 1 & 1 - \vartheta - \vartheta^2 & <0.6 \\
20. & -204 & E & 1 & 1 & & 3 & 1 & 1 - \vartheta - \vartheta^2 & <0.64 \\
\bottomrule
\end{array} $$

\clearpage

\small
$$ \begin{array}{rrrrrrrrcc}
\toprule
Nr. & \disc K & & a_2 & a_1 & a_0 & GHB & h & u & M(K) \\
\midrule
21. & -211 & E & 3 & 1 & 2  & & 1 & 3+\vartheta & \\
22. & -212 & E & 1 & 4 & 2  & &  1 & 1+2\vartheta & \frac{5}{8} \\
23. & -216 & E & 0 & 3 & 2  & &  1 & 1+\vartheta-\vartheta^2 & \\
24. & -231 & E & -4 & 5 & 1 & &  1 & \vartheta & \geq \frac{7}{9} \\
25. & -239 & E & 0 & -1 & 3 & &  1 & 5+3\vartheta & \frac{8}{9} \\
26. & -243 & E & 0 & 0 & 3 & &  1 & 2-\vartheta^2 & \\
27. & -244 & E & 2 & -3 & 2 &  & 1 & 1-2\vartheta+2\vartheta^2 & \\
\bottomrule
\end{array} $$

\clearpage

$$
\begin{array}{rrrrrrrlrc}
\toprule
Nr. & \disc K & a_2 & a_1 & a_0 & GHB & h & u & M(K) & Typ \\
\midrule
61. & -516 & 4 & 1 & 6 & \frac{\vartheta+\vartheta^2}{2} & 1 & -16+(21\vartheta+7\vartheta^2)/2 & & E \\
62. & -519 & 2 & -3 & 3 &  & 1 & 1-\vartheta+\vartheta^2 & \geq \frac{44712}{45747} & E \\
63. & -524 & -2 & 4 & 2 &  & 1 & 1+2\vartheta-\vartheta^2 & & NE \\
64. & -527 & 0 & 5 & 1 &  & 1 & \vartheta & & NE \\
65. & -543 & -2 & 3 & 3 &  & 1 & 1+\vartheta-\vartheta^2 & \geq \frac{7593}{8343} & E \\
66. & -547 & 1 & -3 & 4 &  & 1 & 9-2\vartheta-2\vartheta^2 & & NE \\
67. & -563 & 2 & 6 & 1 &  & 2 & \vartheta & & H \\
\bottomrule
\end{array} $$

\clearpage

$$ \begin{array}{rrrrrlrrrl}
\toprule
Nr. & \disc K & a_2 & a_1 & a_0 & \text{GHB} & h & u & M(K) & Typ \\
\midrule
101 & -823 & 3 & -2 & 3 & & 1 & 1-\vartheta+\vartheta^2 & & NE \\
102 & -835 & 2 & 0 & 5 & & 1 & 18+4\vartheta-\vartheta^2 & & \\
103 & -839 & 4 & 3 & 5 & & 1 & 3+\vartheta+2\vartheta^2 & & NE \\
104 & -843 & 7 & -7 & 2 & & 1 & 3-11\vartheta+12\vartheta^2 & & \\
105 & -856 & -1 & 1 & 11 & frac{1+\vartheta^2}{2} & 1 & 5+10\vartheta+4\vartheta^2 & & \\
106 & -863 & 7 & -4 & 1 & & 1 & \vartheta & & NE \\
107 & -867 & 7 & 5 & 2 & & 1 & 9+21\vartheta+28\vartheta^2 & & NE \\
\bottomrule
\end{array} $$

\clearpage

$$ \begin{array}{rrlrrrlrrr}
\toprule
\text{Nr.} & \disc K & & a_2 & a_1 & a_0 & \text{GHB} & h & u & M(K) \\
\midrule
141 & -1099 & & 1 & -1 & 6 & & & & \\
142 & -1107 & & 0 & 6 & 3 & & & & \\
143 & -1108 & & 5 & 0 & 2 & & 1 & & \\
144 & -1135 & NE & 5 & 2 & 3 & & 1 & & \\
145 & -1144 & & -1 & 6 & 2 & & 1 & & \\
146 & -1147 & NE & 2 & -2 & 5 & & 1 & & \\
147 & -1164 & NE & -1 & -1 & 7 & & 1 & & \\
148 & -1172 & NE & 5 & -4 & 2 & & 1 & & \\
149 & -1175 & NE & 3 & 8 & 1 & & 1 & & \\
150 & -1176 & & 1 & -2 & 6 & & & & \\
\midrule
151 & -1187 & NE & 1 & 7 & 2 & & 1 & & \\
152 & -1188 & & 6 & -3 & 4 & \frac{\vartheta+\vartheta^2}{2} & 1 & & \\
153 & -1191 & NE & 4 & -5 & 3 & & 1 & & \\
154 & -1192 & & 7 & -8 & 4 & \frac{\vartheta+\vartheta^2}{2} & 1 & & \\
155 & -1203 & & 4 & 2 & 5 & & 1 & & \\
156 & -1207 & NE & 4 & -1 & 3 & & 1 & & \\
157 & -1208 & NE & 4 & 13 & 2 & \frac{\vartheta+\vartheta^2}{2} & 1 & & \\
158 & -1219 & NE & 3 & -5 & 4 & & 1 & & \\
159 & -1231 & NE & 2 & -3 & 5 & & 1 & & \\
160 & -1235 & NE & 2 & -4 & 5 & & 1 & & \\
\midrule
161 & -1236 & & 7 & -3 & 3 & \frac{1+\vartheta^2}{2} & 1 & & \\
162 & -1255 & & 4 & 9 & 1 & & & & \\
163 & -1259 & NE & 5 & -1 & 2 & & 1 & & \\
164 & -1267 & NE & 1 & 7 & 4 & & 1 & & \\
165 & -1272 & & -7 & 13 & 9 & \frac{1+\vartheta^2}{2} & 1 & & \\
166 & -1291 & NE & 3 & 1 & 6 & & 1 & & \\
167 & -1292 & NE & 6 & 2 & 2 & & 1 & & \\
168 & -1295 & NE & -4 & 5 & 5 & & 1 & & \\
169 & -1300 & & 9 & -8 & 2 & & 1 & & \\
170 & -1315 & NE & 5 & -3 & 2 & & 1 & & \\
\midrule
171 & -1316 & NE & 7 & 4 & 2 & & 1 & & \\
172 & -1319 & NE & 11 & -6 & 1 & & 1 & & \\
173 & -1323 & & 12 & 6 & 1 & & & & \\
174 & -1327 & NE & 0 & 1 & 7 & & 1 & & \\
175 & -1336 & & 5 & -2 & 2 & & 1 & & \\
176 & -1347 & & 4 & -2 & 3 & & 1 & & \\
177 & -1351 & NE & 1 & 0 & 7 & & 1 & & \\
178 & -1355 & NE & 0 & 2 & 7 & & 1 & & \\
179 & -1363 & NE & 4 & -4 & 3 & & 1 & & \\
180 & -1371 & & -4 & 8 & 1 & & & & \\
\bottomrule
\end{array} $$

\clearpage

\section*{Totalreelle kubische Zahlkörper}
\addcontentsline{toc}{section}
                {Totalreelle kubische Zahlkörper}

\begin{table}[h]
\centering
\small
\begin{tabular}{rrrrrrccc}
\toprule
Nr. & disc K & $a_2$ & $a_1$ & $a_0$ & GHB & h & M(K) \\
\midrule
1.  & 49   & 1  & -2 & -1 &  & 1 & $1/7$   & E \\
2.  & 81   & 0  & -3 & 1  &  & 1 & $1/3$   & E \\
3.  & 148  & 3  & -1 & -1 &  & 1 & $1/2$   & E \\
4.  & 169  & 1  & -4 & 1  &  & 1 & $5/13$  & E \\
5.  & 229  & 0  & -4 & 1  &  & 1 & $1/2$   & E \\
6.  & 257  & 3  & -2 & -1 &  & 1 & $1/3$   & E \\
7.  & 316  & 1  & -4 & -2 &  & 1 & $1/2$   & E \\
8.  & 321  & 1  & -4 & -1 &  & 1 & $1/3$   & E \\
9.  & 361  & -2 & -5 & -1 &  & 1 & $8/19$  & E \\
10. & 404  & 1  & -5 & 1  &  & 1 & $1/2$   & E \\
\midrule
11. & 469  & 4  & -2 & -1 &  & 1 & $1/2$   & E \\
12. & 473  & 0  & -5 & 1  &  & 1 & $1/3$   & E \\
13. & 564  & -7 & 11 & 1  &  & 1 & $1/2$   & E \\
14. & 568  & 1  & -6 & 2  &  & 1 & $1/2$   & E \\
15. & 621  & -6 & 6  & 1  &  & 1 & $1/2$   & E \\
16. & 697  & 3  & -4 & -1 &  & 1 & $13/31$ & E \\
17. & 733  & -2 & -6 & -1 &  & 1 & $1/2$   & E \\
18. & 756  & 0  & -6 & 2  &  & 1 & $1/2$   & E \\
19. & 761  & 1  & -6 & 1  &  & 1 & $1/2$   & E \\
20. & 785  & -2 & -5 & 1  &  & 1 &         & E \\
\midrule
21. & 788  & 9  & 7  & 1  &  & 1 &         & E \\
22. & 837  & 0  & -6 & 1  &  & 1 & $1/2$   & E \\
23. & 892  & -5 & 0  & 2  &  & 1 &         & E \\
24. & 940  & 3  & -4 & -2 &  & 1 &         & E \\
25. & 961  & -2 & -9 & -2 & $\frac{\theta+\theta^2}{2}$ & 1 & $16/31$ & E \\
26. & 985  & 1  & -6 & -1 &  & 1 & 1       & NE \\
27. & 993  & 1  & -6 & -3 &  & 1 &         & E \\
28. & 1016 & -1 & -6 & 2  &  & 1 &         & E \\
29. & 1076 & -5 & -3 & 1  &  & 1 &         & E \\
30. & 1101 & -5 & -1 & 2  &  & 1 &         & E \\
\bottomrule
\end{tabular}
\end{table}

$$ \begin{array}{rrrrrrrrl}
\toprule
Nr. & \disc K & a_2 & a_1 & a_0 & GHB & h & M(K) & \text{type} \\
\midrule
31. & 1129 & -2 & -6 & 1 &  & 1 &        & E \\
32. & 1229 & -6 & -2 & 1 &  & 1 &        & E \\
33. & 1257 & 2  & -7 & 1 &  & 1 &        & E \\
34. & 1300 & 3  & -7 & 1 &  & 1 &        & E \\
35. & 1304 & 11 & 10 & 2 &  & 1 &        & E \\
36. & 1345 & 0  & -7 & 1 &  & 1 & 7/5    & NE \\
37. & 1369 & 4  & -7 & 1 &  & 1 &        & E \\
\bottomrule
\end{array} $$

$$ \begin{array}{rrrrrrrl}
\toprule
Nr. & disc K & a_2 & a_1 & a_0 & \text{GHB} & h & M(K) \\
\midrule
71. & 2241 & 0  & -9 & 5  &  & 1 & E \\
72. & 2292 & 10 & -10& 2  &  & 1 & E \\
73. & 2296 & -7 & 2  & 2  &  & 1 & E \\
74. & 2300 & -1 & -8 & 2  &  & 1 & E \\
75. & 2349 & -6 & 0  & 3  &  & 1 & E \\
76. & 2429 & -7 & 1  & 2  &  & 1 & E \\
77. & 2505 & -4 & -5 & 3  &  & 1 & E \\
78. & 2557 & 1  & -9 & 2  &  & 1 & E \\
79. & 2589 & 10 & -12& 3  &  & 1 & E \\
80. & 2597 & -8 & -2 & 1  &  & 3 & E \\
\midrule
81. & 2636 & -7 & 0  & 2  &  & 1 & E \\
82. & 2673 & 0  & -9 & 3  &  & 1 & E \\
83. & 2677 & -3 & -7 & 2  &  & 1 & E \\
84. & 2700 & -6 & -3 & 2  &  & 1 & E \\
85. & 2708 & -4 & -6 & 2  &  & 1 & E \\
86. & 2713 & -6 & -1 & 3  &  & 1 & E \\
87. & 2777 & -5 & -6 & 1  &  & 2 & E \\
88. & 2804 & 1  & -9 & 1  &  & 1 & E \\
89. & 2808 & 0  & -9 & 2  &  & 1 & E \\
90. & 2836 & -2 & -8 & 2  &  & 1 & E \\
\bottomrule
\end{array} $$

\clearpage

\section*{Totalkomplexe biquadratische Zahlkörper mit quadratischem Teilkörper}
\addcontentsline{toc}{section}
    {Totalkomplexe biquadratische Zahlkörper mit quadratischem Teilkörper} 

$$ \begin{array}{rrclll}
\toprule
Nr. & \disc K & \disc k & \mu & M(K) & \text{type} \\
\midrule
1.  & 117 & -3     & -1+2\sqrt{-3}          & \geq 1/7 & E \\
2.  & 125 & 5      & -5+2\sqrt{5}           & \geq 1/5 & E \\
3.  & 144 & -3,-4,12 &                        & 1/4      & E \\
4.  & 189 & -3     & 3+2\sqrt{-3}           & \geq 1/3 & E \\
5.  & 225 & -3,-15,5 &                      & 1/4      & E \\
6.  & 256 & -4,-8,8 &                        & 1/2      & E \\
7.  & 272 & -4     & 1+4i                   & \geq 1/4 & E \\
8.  & 320 & -4     & 1+2i                   & 1/2      & E \\
9.  & 333 & -3     & -5+2\sqrt{-3}          &          & E \\
10. & 392 & -7     & (-1-\sqrt{-7})/2       &          & E \\
\midrule
11. & 400 & -4,-20,5 &                      & 5/16     & E \\
12. & 432 & -3     & \sqrt{-3}              &          & E \\
13. & 441 & -3,-7,21 &                      & 4/9      & E \\
14. & 512 & -4     &                        & 1/2      & E \\
15. & 513 & -3     & -3+4\sqrt{-3}          &          & E \\
16. & 549 & -3     & 7-2\sqrt{-3}           &          & E \\
17. & 576 & -3,24,-8 &                      & 1/3      & E \\
18. & 576 & -3,-24,8 &                      &          & E \\
19. & 605 & -11    & (3+\sqrt{-11})/2       &          & E \\
20. & 656 & -4     & 5+4i                   & 1/2      & E \\
\midrule
21. & 657 & -3     & 5-4\sqrt{-3}           & 1/2      & E \\
22. & 784 & -4,-7,28 &                      & 1/2      & E \\
23. & 832 & -4     & 3+2i                   & 1/2      & E \\
24. & 837 & -3     & -9+2\sqrt{-3}          &          & E \\
25. & 873 & -3     & -7+4\sqrt{-3}          &          & E \\
26. & 981 & -3     & -1+6\sqrt{-3}          &          & E \\
27. & 1008& -3     & 2-\sqrt{-3}            &          & E \\
28. & 1008& -3     & -2+\sqrt{-3}           &          & E \\
29. & 1025& 5      & -11+4\sqrt{5}          &          & E \\
30. & 1040& -4     & 7+4i                   &          & E \\
\bottomrule
\end{array} $$

\clearpage

$$ \begin{array}{rrclll}
\toprule
Nr. & \disc K & \disc k & \mu & M(K) & \text{type} \\
\midrule
31. & 1040 & -4     & 1+8i                    &      & E \\
32. & 1088 & -8     & 3+2\sqrt{-2}            &      & E \\
33. & 1088 & 8      & -5+2\sqrt{2}            &      & E \\
34. & 1089 & -3,-11,33 &                      &      & E \\
35. & 1161 & -3     & 9+4\sqrt{-3}            &      & E \\
36. & 1168 & -4     & 3+8i                    &      & E \\
37. & 1197 & -3     & 11+2\sqrt{-3}           &      & E \\
\bottomrule
\end{array} $$

\clearpage

$$ \begin{array}{rrccll}
\toprule
Nr. & \disc K & \disc k & \mu & M(K) & \text{type} \\
\midrule
71. & 2304 & -4,-24,24 &         & 5/2  & NE \\
72. & 2304 & -4        & 3i      & 5/2 & NE \\
73. & 2304 & -8,12,-24 &         &     & E   \\
74. & 2320 & -4        & 1+12i   & & \\
75. & 2320 & -4        & 9+8i    & & \\
76. & 2368 & -4        & 1+6i    & & E \\
77. & 2448 & -4        & 3+12i   & & NE \\
78. & 2457 & -3        & 9+8\sqrt{-3} & & \\
79. & 2493 & -3        & -13+6\sqrt{-3} & & \\
80. & 2560 & -4        & 1+3i    & 5/4 & NE \\
\midrule
81. & 2560 & -4        & 3+i     & 5/4 & NE \\
82. & 2597 & -7        & -5+2\sqrt{-7} & & \\
83. & 2601 & -3,17,-51 &               & 13/16   & E \\
84. & 2624 & -8        & -3+4\sqrt{-2} & & \\
85. & 2624 & 8         & -7+2\sqrt{-2} & & \\
86. & 2704 & -4,13,-52 &               & \geq 1  & NE \\
87. & 2709 & -3        & 1+10\sqrt{-3} & & \\
88. & 2709 & -3        & 17+2\sqrt{-3} & & \\
89. & 2736 & -3        & 4+\sqrt{-3} & & \\
90. & 2744 & -7        & (7+\sqrt{-7})/2 & &\\
\midrule
91. & 2781 & -3        & 3+10\sqrt{-3} & & \\
92. & 2817 & -3        & -11+8\sqrt{-3} & & \\
93. & 2873 & 13        & -15+4\sqrt{13} & & \\
94. & 2880 & -4        & 3+6i    & & NE \\
95. & 2925 & -3        & -5+10\sqrt{-3} & & \\
96. & 2960 & -4        & 13+4i   & & \\
97. & 2960 & -4        & 11+8i   & & \\
98. & 3024 & -3        & -3+2\sqrt{-3} & &  \\
99. & 3033 & -3        & 17+4\sqrt{-3} & & \\
100.& 3072 & -8        & 1+\sqrt{-2} & & \\
\midrule
101.& 3072 & -8        & -1-\sqrt{-2} & & \\
102.& 3088 & -4        & 7+12i   & & NE \\
103.& 3136 & -7,-8,56  &         & 9/8     & NE \\
104.& 3141 & -3        & 7+10\sqrt{-3} & & \\
105.& 3249 & -3,-19,57 &         & & E \\
106.& 3328 & -4        & 2+3i    & & \\
107.& 3357 & -3        & 19+2\sqrt{-3} & & \\
108.& 3392 & -4        & 7+2i    & 50/53 & E \\
109.& 3429 & -3        & -9+10\sqrt{-3} & &  \\
110.& 3573 & -3        & 17+6\sqrt{-3} & & \\
\bottomrule
\end{array} $$

\newpage

$$ \begin{array}{rccll}
\toprule
\text{Nr.} & \disc K & \disc k & \mu & \\
\midrule
111. & 3600 & -4,-15,60 &  & NE \\
112. & 3648 & -8 & 7 + 2\sqrt{-2} & \\
113. & 3648 & -8 & 5 + 4\sqrt{-2} & \\
114. & 3681 & -3 & -19 + 4\sqrt{-3} & \\
115. & 3725 & 5 & -13 - 2\sqrt{5} & NE \\
116. & 3728 & -4 & 13 + 8i & NE \\
117. & 3753 & -3 & -15 + 8\sqrt{-3} & \\
118. & 3773 & -7 & 7 + 2\sqrt{-7} & \\
119. & 3789 & -3 & 11 + 10\sqrt{-3} & \\
120. & 3856 & -4 & 15 + 4i & \\
\midrule
121. & 3897 & -3 & 1 + 12\sqrt{-3} & \\
122. & 3904 & -4 & 5 + 6i & NE \\
123. & 4032 & -3 & 4 + 2\sqrt{-3} & \\
124. & 4032 & -3 & -4 - 2\sqrt{-3} & \\
125. & 4077 & -3 & -21 + 2\sqrt{-3} & \\
126. & 4112 & -4 & 1 + 16i & \\
127. & 4113 & -3 & 5 - 12\sqrt{-3} & \\
128. & 4221 & -3 & -13 + 10\sqrt{-3} & \\
129. & 4221 & -3 & 19 + 6\sqrt{-3} & \\
130. & 4240 & -4 & 3 + 16i & \\
\midrule
131. & 4240 & -4 & 11 + 12i & \\
132. & 4329 & -3 & -7 + 12\sqrt{-3} & \\
133. & 4352 & -4 & 4 + i & NE \\
134. & 4352 & -8 & -3 + 2\sqrt{-2} & \\
135. & 4400 & -5 & -8 + 2\sqrt{5} & \\
136. & 4401 & -3 & 21 + 4\sqrt{-3} & \\
137. & 4464 & -3 & 2 + 3\sqrt{-3} & \\
138. & 4464 & -3 & -2 - 3\sqrt{-3} & \\
139. & 4496 & -4 & 5 + 16i & NE \\
140. & 4508 & -7 & 4 + \sqrt{-7} & \\
\midrule
141. & 4624 & -4,17,-68 &  & NE \\
142. & 4672 & -8 & -1 + 6\sqrt{-2} & \\
143. & 4725 & -3 & 15 + 10\sqrt{-3} & \\
144. & 4761 & -3,-23,69 &  & NE \\
145. & 4869 & -3 & 23 + 2\sqrt{-3} & \\
146. & 4880 & -4 & 17 + 4i & \\
147. & 4880 & -4 & 7 + 16i & \\
148. & 4977 & -3 & -19 + 8\sqrt{-3} & \\
149. & 4977 & -3 & -11 + 12\sqrt{-3} & \\
150. & 5008 & -4 & 13 + 12i & \\
\bottomrule
\end{array} $$

\newpage

$$ \begin{array}{rccll}
\toprule
\text{Nr.} & \disc K & \disc k & \mu & \\
\midrule
151. & 5193 & -3 & -23 + 4\sqrt{-3} & \\
152. & 5301 & -3 & -17 + 10\sqrt{-3} & \\
153. & 5301 & -3 & -1 + 14\sqrt{-3} & \\
154. & 5328 & -3 & 5 + 2\sqrt{-3} & \\
155. & 5341 & -7 & -9 + 2\sqrt{-7} & \\
156. & 5392 & -4 & 9 + 16i & \\
157. & 5440 & -4 & 9 + 2i & \\
158. & 5488 & -7 & 2\sqrt{-7} & \\
159. & 5537 & -7 & 1 + 4\sqrt{-7} & \\
160. & 5648 & -4 & 17 + 8i & \\
\midrule
161. & 5696 & -8 & -9 + 2\sqrt{-2} & NE \\
162. & 5776 & -4,-19,76 & & NE \\
163. & 5929 & -7,-11,77 & & NE \\
164. & 6032 & -4 & 19 + 4i & \\
165. & 6032 & -4 & 11 + 16i & \\
166. & 6144 & -8 & 2 + \sqrt{-2} & \\
167. & 6144 & -8 & -2 + \sqrt{-2} & \\
168. & 6208 & -8 & -5 + 6\sqrt{-2} & \\
169. & 6400 & -4 & 5i & \\
170. & 6400 & -8,-20,40 & & NE \\
\midrule
171. & 6416 & -4 & 1 + 20i & \\
172. & 6464 & -4 & 1 + 10i & \\
173. & 6544 & -4 & 3 + 20i & \\
174. & 6656 & -4 & 5 + i & NE \\
175. & 6656 & -4 & -1 + 5i & NE \\
176. & 6713 & -7 & 5 + 4\sqrt{-7} & \\
177. & 6725 & 5 & -17 - 2\sqrt{5} & \\
178. & 6800 & -4 & 5 + 20i & NE \\
179. & 6928 & -4 & 17 + 12i & \\
180. & 6976 & -4 & 3 + 10i & NE \\
\midrule
181. & 7056 & -4,21,-84 & & NE \\
182. & 7184 & -4 & 7 + 20i & \\
183. & 7312 & -4 & 21 + 4i & \\
184. & 7424 & -4 & 2 + 5i & NE \\
185. & 7488 & -4 & 9 + 6i & NE \\
186. & 7696 & -4 & 9 + 20i & \\
187. & 7696 & -4 & 15 + 16i & NE \\
188. & 7744 & -8,-11,44 & & NE \\
189. & 8448 & -8 & 5 + 2\sqrt{-2} & \\
\bottomrule
\end{array} $$

\clearpage

\begin{center}
\textbf{Reelle biquadratische K\"orper mit quadratischem Teilk\"orper}
\end{center}

\medskip
Bezeichnungen wie oben; weiter ist

\begin{tabbing}
\hspace*{1.5cm}\= \kill
$V_4$\> die Kleinsche Vierergruppe \\
$Z_4$\> die zyklische Gruppe der Ordnung 4 \\
$D_4$\> die Diedergruppe der Ordnung 8 \\
$A_4$\> die alternierende Gruppe der Ordnung 12 \\
$S_4$\> die symmetrische Gruppe der Ordnung 24
\end{tabbing}

\medskip
\[
\begin{array}{r@{\quad}r@{\quad}c@{\quad}l@{\quad}l@{\quad}c}
\hline
\text{Nr.} & \text{disc}\,K &  & \text{disc}\,k & \mu & \mathrm{Gal}(K/\Q) \\
\hline
1.  & 725  & \mathrm{E} & 5           & 7+2\sqrt{5}     & D_4 \\
2.  & 1125 & \mathrm{E} & 5           & 15+6\sqrt{5}    & Z_4 \\
3.  & 1600 & \mathrm{E} & 5,8,40      &                 & V_4 \\
4.  & 2000 &           & 5           & 5+2\sqrt{5}     & Z_4 \\
5.  & 2048 & \mathrm{E} & 8           & 2+\sqrt{2}      & D_4 \\
6.  & 2225 & \mathrm{E} & 5           & 13+4\sqrt{5}    & D_4 \\
7.  & 2304 & \mathrm{E} & 8,12,24     &                 & V_4 \\
\hline
\end{array}
\]

\clearpage

\section*{Übersicht über die Verteilung normeuklidischer Zahlkörper}
\addcontentsline{toc}{section}
    {Übersicht über die Verteilung normeuklidischer Zahlkörper}

Die Anzahl der bekannten normeuklidischen Zahlkörper hat in den
letzten zwanzig Jahren stark zugenommen. Die folgenden Tabellen, die
die Verteilung der bekannten normeuklidischen Zahlringe bezüglich $n$
(Körpergrad) und $r+s$ (Einheitenrang + 1) zeigen, sollen dies
verdeutlichen:

\bigskip

\textbf{I. 1967 (Kummer (1844), Eisenstein (1850), Godwin (1965b, 1967))}

\begin{center}
\begin{tabular}{c|ccccc|c}
\toprule
$r+s$ & 1 & 2 & 3 & 4 & 5 & $\Sigma$ \\
\midrule
1 & 1 & 5 & & & & 6 \\
2 & & 16 & 15 & 3 & & 34 \\
3 & & & 9 & 0 & 0 & 9 \\
4 & & & & 8 & 0 & 8 \\
5 & & & & & 1 & 1 \\
\midrule
$\Sigma$ & 1 & 21 & 24 & 11 & 1 & 58 \\
\bottomrule
\end{tabular}
\end{center}

\textbf{II. 1974 (Lenstra (1974), Smith (1969, 1971), Lakein (1972))}

\begin{center}
\begin{tabular}{c|cccccccccc|c}
\toprule
$r+s$ & 1 & 2 & 3 & 4 & 5 & 6 & 7 & 8 & 9 & 10 & $\Sigma$ \\
\midrule
1 & 1 & 5 & & & & & & & & & 6 \\
2 & & 16 & 15 & 31 & & & & & & & 62 \\
3 & & & 57 & 2 & 1 & 2 & & & & & 62 \\
4 & & & & 9 & 0 & 0 & 0 & 2 & & & 11 \\
5 & & & & & 1 & 1 & 0 & 0 & 0 & 1 & 3 \\
6 & & & & & & 0 & 0 & 0 & 0 & 0 & 0 \\
\midrule
$\Sigma$ & 1 & 21 & 72 & 42 & 2 & 3 & 0 & 2 & 0 & 1 & 144 \\
\bottomrule
\end{tabular}
\end{center}

\clearpage

\textbf{III. Sept. 1976 (Lenstra (1977), Taylor (1976))}

\begin{center}
\begin{tabular}{c|cccccccccc|c}
\toprule
$r+s$ & 1 & 2 & 3 & 4 & 5 & 6 & 7 & 8 & 9 & 10 & $\Sigma$ \\
\midrule
1 & 1 & 5 & & & & & & & & & 6 \\
2 & & 16 & 52 & 32 & & & & & & & 100 \\
3 & & & 57 & 11 & 12 & 28 & & & & & 108 \\
4 & & & & 9 & 10 & 25 & 23 & 24 & & & 91 \\
5 & & & & & 1 & 2 & 0 & 0 & 0 & 1 & 3 \\
6 & & & & & & 2 & 0 & 0 & 0 & 0 & 0 \\
\midrule
$\Sigma$ & 1 & 21 & 109 & 52 & 23 & 57 & 23 & 24 & 0 & 1 & 311 \\
\bottomrule
\end{tabular}
\end{center}

\textbf{VII. August 1989:} Außer den in dieser Arbeit gefundenen
normeuklidischen Körpern ist auch $\mathbb{Q}(\zeta_{13})$ in der
folgenden Tabelle aufgeführt (sh. dazu Leutbecher u. Niklasch 1987)

\begin{center}
\begin{tabular}{c|cccccccccccc|c}
\toprule
$r+s$ & 1 & 2 & 3 & 4 & 5 & 6 & 7 & 8 & 9 & 10 & 11 & 12 & $\Sigma$ \\
\midrule
1 & 1 & 5 & & & & & & & & & & & 6 \\
2 & & 16 & 52 & 50 & & & & & & & & & 118 \\
3 & & & 62 & 13 & 13 & 28 & & & & & & & 116 \\
4 & & & & 15 & 10 & 37 & 39 & 45 & & & & & 146 \\
5 & & & & & 1 & 12 & 26 & 65 & 92 & 50 & & & 246 \\
6 & & & & & & 3 & 3 & 3 & 0 & 0 & 0 & 1 & 10 \\
\midrule
$\Sigma$ & 1 & 21 & 114 & 78 & 24 & 80 & 68 & 113 & 92 & 50 & 0 & 1 & 642 \\
\bottomrule
\end{tabular}
\end{center}

\addtocontents{toc}{\protect\enlargethispage{2\baselineskip}}

\vskip 1cm

Remarks on the Bibliography.

Further literature references are to be found in
\begin{itemize}
\item Lekkerkerker (1969) (probabilistic homogeneous linear forms,
  geometry of numbers etc.)
\item Lenstra (1974) (in particular on works concerning the
  Euclidean algorithm in function fields),
\item van der Linden (1985)
\item Zimmer (1972) (determination of discriminant, integral basis,
  unit and ideal class groups of an algebraic
  number field)
\end{itemize}
That's it.

\printindex[dN]
\printindex[dS]

\end{document}